\documentclass{article}

\usepackage[T1]{fontenc}
\usepackage[paperwidth=7in,paperheight=11in,
hmargin=.25in,vmargin=1in]{geometry} 
\usepackage{epic}
\usepackage{amsmath, amscd, amsfonts, amsthm, amssymb, enumerate, graphicx, cases}
\usepackage[dvipsnames]{xcolor}
\usepackage[shortlabels]{enumitem}
\usepackage{mathtools, verbatim, 
  rotating, comment}
\usepackage{url}
\usepackage{hyperref} 
\usepackage{multirow}
\usepackage[matrix,arrow,curve]{xy}

\newtheorem{thm}{Theorem}[section]
\newtheorem{thmstar}[thm]{Theorem*}
\newtheorem{lem}[thm]{Lemma}

\newtheorem{prop}[thm]{Proposition}
\newtheorem{cor}[thm]{Corollary}

\newtheorem{conj}[thm]{Conjecture}
\newtheorem{prob}[thm]{Problem}
\theoremstyle{definition}
\newtheorem{conv}[thm]{Conventions}
\newtheorem{nota}[thm]{Notation}
\newtheorem{defn}[thm]{Definition}
\newtheorem{examp}[thm]{Example}
\theoremstyle{remark}
\newtheorem{rem}[thm]{Remark}

\renewcommand{\sf}[1]{\mathcal{#1}}
\DeclareMathOperator{\tr}{tr}
\DeclareMathOperator{\id}{id}
\DeclareMathOperator{\ct}{ct}
\DeclareMathOperator{\ev}{ev}

\DeclareMathOperator{\ntc}{ntc}

\DeclareMathOperator{\sgn}{sgn}
\DeclareMathOperator{\ch}{char}

\DeclareMathOperator{\coef}{coef}
\DeclareMathOperator{\disc}{disc}

\DeclareMathOperator{\im}{im}
\DeclareMathOperator{\coker}{coker}

\DeclareMathOperator{\inv}{inv}
\DeclareMathOperator{\loc}{loc}
\DeclareMathOperator{\res}{res}
\DeclareMathOperator{\supp}{supp}
\DeclareMathOperator{\sym}{sym}

\DeclareMathOperator{\ldr}{ldr}

\let\th\undefined \DeclareMathOperator{\th}{th}
\newcommand{\ur}{\mathrm{ur}}
\newcommand{\ram}{\mathrm{ram}}
\newcommand{\nonmax}{\mathrm{nonmax}}
\DeclareMathOperator{\Res}{res}
\DeclareMathOperator{\Cor}{cor}
\DeclareMathOperator{\Hom}{Hom}
\DeclareMathOperator{\Surj}{Surj}
\DeclareMathOperator{\Frac}{Frac}
\DeclareMathOperator{\Disc}{Disc}
\DeclareMathOperator{\Stab}{Stab}

\DeclareMathOperator{\Coord}{Coord}
\DeclareMathOperator{\Kum}{Kum}
\DeclareMathOperator{\Ind}{Ind}

\DeclareMathOperator{\Pic}{Pic}

\DeclareMathOperator{\Div}{Div}
\DeclareMathOperator{\Cl}{Cl}
\DeclareMathOperator{\Sel}{Sel}
\DeclareMathOperator{\Sym}{Sym}
\DeclareMathOperator{\Aut}{Aut}
\DeclareMathOperator{\Gal}{Gal}
\DeclareMathOperator{\End}{End}
\DeclareMathOperator{\Mat}{Mat}
\DeclareMathOperator{\AS}{AS}
\DeclareMathOperator{\charpoly}{char\,poly}
 
\let\I\Im
\let\Im\undefined \DeclareMathOperator{\Im}{Im}

\newcommand{\GL}{\mathrm{GL}}
\newcommand{\GA}{\sf{GA}}
\newcommand{\SL}{\mathrm{SL}}
\newcommand{\PGL}{\mathrm{PGL}}
\newcommand{\PSL}{\mathrm{PSL}}
\newcommand{\GGamma}{\mathrm{G\Gamma}}
\newcommand{\SO}{\mathrm{SO}}
\renewcommand{\AA}{\mathbb{A}}
\newcommand{\CC}{\mathbb{C}}
\newcommand{\FF}{\mathbb{F}}
\newcommand{\GG}{\mathbb{G}}
\newcommand{\NN}{\mathbb{N}}
\newcommand{\QQ}{\mathbb{Q}}
\newcommand{\PP}{\mathbb{P}}
\newcommand{\RR}{\mathbb{R}}

\newcommand{\ZZ}{\mathbb{Z}}
\renewcommand{\aa}{\mathfrak{a}}
\newcommand{\bb}{\mathfrak{b}}
\newcommand{\cc}{\mathfrak{c}}
\newcommand{\dd}{\mathfrak{d}}
\newcommand{\ff}{\mathfrak{f}}
\newcommand{\mm}{\mathfrak{m}}
\newcommand{\pp}{\mathfrak{p}}
\newcommand{\qq}{\mathfrak{q}}
\renewcommand{\ss}{\mathfrak{s}}
\renewcommand{\tt}{\mathfrak{t}}
\newcommand{\C}{\mathcal{C}}
\newcommand{\D}{\mathcal{D}}

\newcommand{\F}{\mathcal{F}}
\newcommand{\G}{\mathcal{G}}
\newcommand{\Idls}{\mathcal{I}}
\renewcommand{\L}{\mathcal{L}}
\newcommand{\M}{\mathcal{M}}
\newcommand{\N}{\mathcal{N}}
\newcommand{\OO}{\mathcal{O}}
\renewcommand{\P}{\mathfrak{P}}

\newcommand{\R}{\mathcal{R}}
\renewcommand{\S}{\mathcal{S}}

\newcommand{\U}{\mathcal{U}}
\newcommand{\V}{\mathcal{V}}

\newcommand{\Sm}{\mathfrak{S}}
\newcommand{\fcr}{\mathfrak{r}}
\newcommand{\X}{\mathfrak{X}}
\newcommand{\Y}{Y}
\newcommand{\sfa}{a}
\newcommand{\g}{g}
\newcommand{\h}{h}
\renewcommand{\t}{t}
\newcommand{\x}{x}
\newcommand{\y}{y}
\newcommand{\spi}{\pi}
\newcommand{\Ell}{\mathrm{Ell}}
\newcommand{\A}{\mathcal{A}}
\newcommand{\B}{\mathcal{B}}
\renewcommand{\1}{\mathbf{1}}
\newcommand{\0}{{\mathord{0}}}
\newcommand{\tee}{{\mathord{\top}}}
\renewcommand{\*}{{\mathord{\star}}}
\newcommand{\ba}{\overline}
\newcommand{\bs}{\backslash}

\newcommand{\cross}{\times}
\newcommand{\tensor}{\otimes}

\newcommand{\toto}{\twoheadrightarrow}
\newcommand{\longto}{\mathop{\longrightarrow}\limits}
\newcommand{\textand}{\quad \text{and} \quad}
\newcommand{\textor}{\quad \text{or} \quad}
\renewcommand{\to}{\mathop{\rightarrow}\limits}
\renewcommand{\Dot}[1]{#1^\triangle}
\newcommand{\size}[1]{\lvert #1 \rvert}
\newcommand{\Size}[1]{\left\lvert #1 \right\rvert}
\newcommand{\floor}[1]{\left\lfloor #1 \right\rfloor}
\newcommand{\ceil}[1]{\left\lceil #1 \right\rceil}
\newcommand{\intsec}{\cap}
\newcommand{\union}{\cup}

\newcommand{\nequiv}{\not\equiv}
\newcommand{\isom}{\cong}
\newcommand{\<}{\left\langle}
\renewcommand{\>}{\right\rangle}
\renewcommand{\(}{\left(}
\renewcommand{\)}{\right)}
\newcommand{\laurent}[1]{(\!(#1)\!)}

\newcommand{\afterhat}{\,\widehat{}}
\newcommand{\ignore}[1]{}
\newcommand{\ds}{\displaystyle}

\renewcommand{\epsilon}{\varepsilon}

\newcommand{\sh}[2]{%
  \put(#1,#2){\rule{\the\unitlength}{\the\unitlength}}
}

\newcommand{\bbq}[8]{
  \begin{minipage}{0.1\linewidth}
    \xymatrix@!0{
      & #5 \ar@{-}[rr]\ar@{-}[dd]
      & & #6 \ar@{-}[dd]
      \\
      #1 \ar@{-}[ur]\ar@{-}[rr]\ar@{-}[dd]
      & & #2 \ar@{-}[ur]\ar@{-}[dd]
      \\
      & #7 \ar@{-}[rr]
      & & #8
      \\
      #3 \ar@{-}[rr]\ar@{-}[ur]
      & & #4 \ar@{-}[ur]
    }
  \end{minipage}
}

\excludecomment{wild}
\excludecomment{todofn}
\excludecomment{old}

\begin{document}
\title{Reflection theorems for number rings}
\author{Evan O'Dorney}
\maketitle
\setcounter{tocdepth}{2}
\tableofcontents

\begin{abstract}
The Ohno-Nakagawa reflection theorem is an unexpectedly simple identity relating the number of $\GL_2 \ZZ$-classes of binary cubic forms (equivalently, cubic rings) of two different discriminants $D$, $-27D$; it generalizes cubic reciprocity and the Scholz reflection theorem. In this paper, we provide a framework for generalizing this theorem using a global and local step. The global step uses Fourier analysis on the adelic cohomology $H^1(\AA_K, M)$ of a finite Galois module, modeled after the celebrated Fourier analysis on $\AA_K$ used in Tate's thesis. The local step is combinatorial, more elementary but much more mysterious. We establish reflection theorems for binary quadratic forms over number fields of class number $1$, and for cubic and quartic rings over arbitrary number fields, as well as binary quartic forms over $\ZZ$; the quartic results are conditional on some computational algebraic identities that are probabilistically true. Along the way, we find elegant new results on Igusa zeta functions of conics and the average value of a quadratic character over a box in a local field.
\end{abstract}

\part{Introduction}
\section{Introduction}

\subsection{Historical background}

In 1932, using the then-new machinery of class field theory, Scholz \cite{ScholzRefl} proved that the class groups of the quadratic fields $\QQ(\sqrt{D})$ and $\QQ(\sqrt{-3D})$, whose discriminants are in the ratio $-3$, have $3$-ranks differing by at most $1$.
This is a remarkable early example of a \emph{reflection theorem.} A
generalization due to Leopoldt \cite{Leopoldt} relates different components of the $p$-torsion
of the class group of a number field containing $\mu_p$ when decomposed under
the Galois group of that field. Applications of such reflection theorems are
far-ranging: for instance, Ellenberg and Venkatesh \cite{EV} use reflection theorems of Scholz type to prove upper bounds on $\ell$-torsion in class groups of number fields, while Mih\u ailescu \cite{MihCat2} uses Leopoldt's
generalization to simplify a step of his monumental proof of the Catalan conjecture that $8$ and $9$ are the only consecutive perfect powers. Through the
years, numerous reflection principles for different generalizations of ideal
class groups have come into print. A very general reflection theorem for Arakelov class
groups is due by Gras \cite{Gras}.

A quite different direction of generalization was discovered by accident in 1997: The following relation was conjectured by Ohno \cite{Ohno} on the basis of numerical data and proved by Nakagawa \cite{Nakagawa}, for which reason we will call it the \emph{Ohno-Nakagawa (O-N) reflection theorem}:
\begin{thm}[Ohno--Nakagawa] \label{thm:O-N}
For a nonzero integer $D$, let $h(D)$ be the number of $\GL_2(\ZZ)$-orbits of \emph{binary cubic forms}
\[
  f(x,y) = ax^3 + bx^2y + cxy^2 + dy^3
\]
of discriminant $D$, each orbit weighted by the reciprocal of its number of symmetries (i.e.~stabilizer in $\GL_2(\ZZ)$). Let $h_3(D)$ be the number of such orbits $f(x,y)$ such that the middle two coefficients $b, c$ are multiples of $3$, weighted in the same way.

Then for every nonzero integer $D$, we have the exact identity
\begin{equation}
\label{eq:O-N_cubic}
h_3(-27 D) = \begin{cases}
3 h(D), & D > 0 \\
h(D), & D < 0.
\end{cases}
\end{equation}
\end{thm}

By the well-known index-form parametrization (see \ref{thm:hcl_cubic_ring} below), $h(D)$ also counts the cubic \emph{rings} of discriminant $D$ over $\ZZ$, weighted by the reciprocal of the order of the automorphism group. It turns out that $h_3(D)$ counts those rings $C$ for which $3 | \tr_{C/\ZZ} \xi$ for every $\xi \in C$. When $D$ is a fundamental discriminant, the corresponding cubic extensions are closely related, via class field theory, to the $3$-class group of $\QQ(\sqrt{D})$ and we get back Scholz's reflection theorem, as Nakagawa points out (\cite{Nakagawa}, Remark 0.9).

Theorem \ref{thm:O-N} was quite unexpected, because $\GL_2(\ZZ)$-orbits of binary cubics have been tabulated since Eisenstein without unearthing any striking patterns. Even the exact normalizations $h(D)$, $h_3(D)$ had been in use for over two decades. They appear in the \emph{Shintani zeta functions}
\begin{align*}
  \zeta^{\pm}(s) &= \sum_{n \geq 1} \frac{h(\pm n)}{n^s} \\
  \hat\zeta^{\pm}(s) &= \sum_{n \geq 1} \frac{h_3(\pm 27 n)}{n^s},
\end{align*}
a family of Dirichlet series which play a prominent role in understanding the distribution of cubic number fields, similar to how the famous Riemann zeta function controls the distribution of primes. As Shintani proved as early as 1972 \cite{Shintani}, the Shintani zeta functions satisfy a matrix functional equation (see Nakagawa \cite{Nakagawa}, eq.~(0.1))
\begin{equation} \label{eq:fnl_eq}
  \begin{bmatrix}
    \zeta^+(1-s) \\ \zeta^-(1-s)
  \end{bmatrix}
  = 2^{-1}3^{3 s - 2} \pi^{-4s} \Gamma\(s - \frac{1}{6}\)\Gamma(s)^2 \Gamma\(s + \frac{1}{6}\)
  \begin{bmatrix}
    \sin 2\pi s & \sin \pi s \\ 3 \sin \pi s & \sin 2\pi s
  \end{bmatrix}
  \begin{bmatrix}
    \hat\zeta^+(s) \\ \hat\zeta^-(s)
  \end{bmatrix}
\end{equation}
The condition that $3$ divide $b$ and $c$ is equivalent to requiring that the cubic form $f$ is \emph{integer-matrix}, that is, its corresponding symmetric trilinear form
\[
\bbq {a}{b/3}{b/3}{c/3}{b/3}{c/3}{c/3}{d}
\]
has integer entries. This condition arose in Shintani's work by taking the dual lattice to $\ZZ^4$ under the pairing
\begin{equation} \label{eq:ip cubic}
  \<(a,b,c,d), (a',b',c',d')\> = ad' - \frac{1}{3} bc' + \frac{1}{3} cb' - da',
\end{equation}
which plays a central role in proving the functional equation. However, as we will find, the pairing \eqref{eq:ip cubic} does not figure in the proof of our reflection theorems, which indeed often relate lattices that are not dual under it.

Using the functional equation, Shintani proved that the $\zeta^{\pm}$ admit meromorphic continuations to the complex plane with simple poles at $1$ and $5/6$, inspiring him to conjecture that the number $N_\pm(X)$ of cubic fields of positive or negative discriminant up to $X$ has the shape
\[
N_\pm(X) = a_\pm X + b_\pm X^{5/6} + o(X^{5/6})
\]
for suitable constants $a_\pm$ and $b_\pm$. This conjecture was proven by Bhargava, Shankar, and Tsimerman \cite{BST-2ndOrd} and independently by Taniguchi and Thorne \cite{TT_rc}. Neither proof needs the Ohno-Nakagawa reflection theorem (Theorem \ref{thm:O-N}), which appears in the notation of Shintani zeta functions in the succinct form
\begin{equation} \label{eq:O-N_Shintani}
  \hat\zeta^+(s) = \zeta^-(s) \quad \text{and} \quad \hat\zeta^-(s) = 3\zeta^+(s).
\end{equation}
\begin{rem}
In the earlier papers, the term ``Ohno-Nakagawa identities'' was used, referring to the pair \eqref{eq:O-N_Shintani}. Our work confirms the intuition that, despite the different scalings, both identities are essentially one theorem.
\end{rem}

\subsection{Methods}

Several proofs of O-N are now in print (\cite{Nakagawa,Marinescu,OOnARemarkable,Gao}), all of which consist of two main steps:
\begin{itemize}
  \item A ``global'' step that uses global class field theory to understand cubic fields, equivalently $\GL_2(\QQ)$-orbits of cubic forms;
  \item A ``local'' step to count the rings in each cubic field, equivalently the $\GL_2(\ZZ)$-orbits in each $\GL_2(\QQ)$-orbit, and put the result in a usable form.
\end{itemize}
In this paper, the distinction between these steps will be formalized and clarified.

For the global step, we take inspiration from Tate's celebrated thesis \cite{Tate_thesis}, which uses Fourier analysis on the adeles to give illuminating new proofs of the functional equations for the Riemann $\zeta$-function and various $L$-functions. Taniguchi and Thorne (see \cite{TT_oexp}) used Fourier analysis on the space of binary cubic forms over $\FF_q$ to get the functional equation for the Shintani zeta function of forms satisfying local conditions at primes. Despite the similarities, their work is essentially independent from ours. We are also inspired by a remark due to Calegari in a paper of Cohen, Rubinstein-Salzedo, and Thorne (\cite{CohON}, Remark 1.6), pointing out that their reflection theorem counting dihedral fields of prime order can also be derived from a theorem of Greenberg and Wiles for the sizes of Selmer groups in Galois cohomology.

We present a notion of \emph{composed variety}, a scheme $\V$ over the ring of integers $\OO_K$ of a number field admitting an action of an algebraic group $\G$ over $\OO_K$. Our guiding example is the scheme $\V$ of binary cubic forms of discriminant $D$ with its action of $\V = \SL_2$. The term ``composed'' refers to the presence of a composition law on the orbits, which relate naturally to a Galois cohomology group $H^1(K,M)$. Our (global) reflection theorems can be stated as saying that two composed varieties $\V^{(1)}$, $\V^{(2)}$ have the same number of $\OO_K$-points, with a suitable weighting. Introducing a new technique of Fourier analysis on the adelic cohomology group $H^1(\AA_K, M) = \prod'_v H^1(K_v,M)$, based on Poitou-Tate duality, we present a generalized reflection engine (Theorems \ref{thm:main_compose} and \ref{thm:main_compose_multi}) that reduces global reflection theorems to \emph{local reflection theorems,} that is, statements involving only the $\OO_{K_v}$-points of $\V^{(1)}$ and $\V^{(2)}$ for a single place $v$ of $K$. A typical case is Theorem \ref{thm:O-N_cubic_local}.

These local reflection theorems are approachable by elementary methods but can be difficult to prove. We present two kinds of proofs. The first is a bijective argument involving Bhargava's self-balanced ideals that is very clean but has only been discovered at the ``tame primes'' ($\pp \nmid 3$ in the cubic case, $\pp \nmid 2$ in the quartic). The second is by explicitly computing the number of orders of given resolvent in a cubic or quartic algebra. We express it as a generating function in a number of variables depending on the splitting type of the resolvent. The generating function is rational, and local reflection can be written as an equality between two rational functions; but these functions are so complicated that the best approximation to a proof of the identity that we can find is a Monte Carlo proof, namely, substituting random values for the variables in some large finite field and verifying that the equality holds. The reader is invited to recheck this verification using the source code in Sage that will be made available with the final version of this paper. 

\subsection{Results}

We are able to prove O-N for binary cubic forms over all number fields $K$, verifying and extending the conjectures of Dioses \cite[Conjecture 1.1]{Dioses}. However, we go further and ask whether every $\SL_2(\OO_K)$-invariant lattice within the space $V(K)$ of binary cubic forms admits an O-N-style reflection theorem. Over $\ZZ$, this question was answered affirmatively for each of the ten invariant lattices by Ohno and Taniguchi \cite{10lat}. Over $\OO_K$, such lattices were classified by Osborne \cite{Osborne}, and they differ from one another only at the primes dividing $2$ and $3$. The lattices at $3$ yield an elegant reflection theorem (Theorem \ref{thm:O-N_traced}) in which the condition $b,c \in \tt$, where $\tt$ is an ideal dividing $3$ in $\OO_K$, reflects to $b,c \in 3\tt^{-1}$, the complementary divisor. At $2$, the corresponding reflection theorems still exist, though they become difficult to write explicitly: see Theorem \ref{thm:invarlat}.

We also find a new reflection theorem (Theorem \ref{thm:O-N_quad}) counting binary \emph{quadratic} forms, not by discriminant, but by a curious invariant: the product $a(b^2 - 4ac)$ of the discriminant and the leading coefficient. Over $\ZZ$, the reflection theorem (Theorem \ref{thm:O-N_quad_Z}) has the potential to be proved simply using quadratic reciprocity, eschewing the machinery of Galois cohomology, though it seems unlikely that the theorem would have ever been discovered without it.

Nakagawa has also conjectured \cite{NakPairs} a reflection theorem for \emph{pairs of ternary quadratic forms,} which parametrize \emph{quartic rings.} The natural invariant to count by is the discriminant, but it is more natural from our perspective to subdivide further and ask for a reflection theorem for rings with fixed \emph{cubic resolvent}, which holds in the known cases \cite[Theorem 1]{NakPairs}. Here our global framework applies without change, but the local enumeration of orders in a quartic field presents formidable combinatorial difficulties, especially in the wildly ramified ($2$-adic) setting, which have been attacked in another work of Nakagawa \cite{NakOrders}. Our methods have the potential to finish this work, but because we count by resolvent rather than discriminant, our answers do not directly match his. 

The process of proving local quartic O-N leads us down some fruitful routes that do not at first sight have any connection to reflection theorems or to the enumeration of quartic rings. These include new cases of the Igusa zeta functions of conics (Lemmas \ref{lem:conic_1} and \ref{lem:conic_pi}) and a result on the average value of a quadratic character on a box in a local field (Theorem \ref{thm:char_box}). If quartic O-N holds true in all cases, it implies that the cubic resolvent ring (in the sense of Bhargava) of a maximal quartic order has a second natural  characterization: it is the ``conductor ring'' for which the Galois-naturally attached extension $K_6/K_3$ is a ring class field (Theorem* \ref{thm*:cond_ring}).

\subsection{Outline of the paper}
In Section \ref{sec:layman's_appendix}, we state and give examples of the main global reflection theorems of the paper over $\ZZ$, in a fashion that requires a minimum of prior knowledge, for the end of further diffusing interest in, and appreciation of, the beauty of number theory.

In Part \ref{part:Gal_coho}, we lay out preliminary matter, much of which is closely related to results that have appeared in the literature but under different guises. It includes a simple characterization (Proposition \ref{prop:H1}) of Galois $H^1$ in terms of \'etale algebras whose Galois group is a semidirect product. It also includes a theorem (Theorem \ref{thm:levels}) on the structure of $H^1(K, M)$ in the case that $K$ is local and $M \isom \C_p$ (with any Galois structure), which will be invaluable in what follows.

In Part \ref{part:composed}, we lay out the framework of composed varieties, on which we perform the novel technique of Fourier analysis of the local and global Tate pairings to get our main local-to-global reflection engine (Theorems \ref{thm:main_compose} and \ref{thm:main_compose_multi}). The remainder of the paper will concern applications of this engine.

In Part \ref{part:first}, we prove two relatively simple reflection theorems: one for quadratic forms (Theorem \ref{thm:O-N_quad}), and a version of the Scholz reflection principle for class groups of quadratic orders (Theorem \ref{thm:Scholz_for_locally_dual_orders}).

In Part \ref{part:cubic}, we prove our extensions of Ohno-Nakagawa for cubic forms and rings.

The quartic case is dealt with in Parts \ref{part:quartic} and \ref{part:quartic_count}: the first part dealing with the bijective methods, and the second with the (long) work of explicitly counting orders in each quartic algebra. The case of partially ramified cubic resolvent (splitting type $1^21$) is still in progress, so we restrict our attention to the four tamely splitting types in the present version. 

We conclude the paper with some unanswered questions engendered by this research.

\subsection{Acknowledgements}
For fruitful discussions, I would like to thank (in no particular order):
Manjul Bhargava, Xiaoheng Jerry Wang, Fabian Gundlach, Levent Alp\"oge, Melanie Matchett Wood, Kiran Kedlaya, Alina Bucur, Benedict Gross, Sameera Vemulapalli, Brandon Alberts, Peter Sarnak, and Jack Thorne.

\section{Examples for the lay reader}
\label{sec:layman's_appendix}
Fortunately for the non-specialist reader, the statements (though not the proofs) of the main results in this thesis can be stated in a way requiring little more than high-school algebra. We here present these statements and some examples to illustrate them.

\subsection{Reflection for quadratic equations}
\begin{defn}\label{defn:quad_inv}
Let $f(x) = ax^2 + bx + c$ be a quadratic polynomial, where the coefficients $a$, $b$, $c$ are integers. The \emph{superdiscriminant} of $f$ is the product
\[
  I = a \cdot (b^2 - 4ac)
\]
of the leading coefficient with the usual discriminant.
\end{defn}

\begin{lem} \label{lem:quad_inv}
If we replace $x$ by $x + t$ in a quadratic polynomial $f$, where $t$ is a fixed integer, then the superdiscriminant does not change.
\end{lem}
\begin{proof}
This can be verified by brute-force calculation, but the following method is more illuminating. The discriminant is classically related to the two roots of $f$,
\[
  x_1 = \frac{-b + \sqrt{b^2 - 4ac}}{2a} \textand
  x_2 = \frac{-b - \sqrt{b^2 - 4ac}}{2a},
\]
through their \emph{difference:}
\begin{align*}
  x_1 - x_2 &= \frac{2\sqrt{b^2 - 4ac}}{2a} = \frac{\sqrt{b^2 - 4ac}}{a} \\
  (x_1 - x_2)^2 &= \frac{b^2 - 4ac}{a^2} \\
  a^3(x_1 - x_2)^2 &= a \cdot (b^2 - 4ac) = I.
\end{align*}
If we replace $x$ by $x + t$, then $a$ does not change, and both roots $x_1, x_2$ are decreased by $t$, so their difference $x_1 - x_2$ is unchanged. Therefore $I$ is unchanged.
\end{proof}

\begin{defn}
  Call two quadratics $f_1$, $f_2$ \emph{equivalent} if they are related by a translation $f_2(x) = f_1(x + t)$. If $I$ is a nonzero integer, let $q(I)$ be the number of quadratics of superdiscriminant $I$, up to equivalence. Let $q_2(I)$, $q^+(I)$, $q_2^+(I)$ be the number of such quadratics that satisfy certain added conditions:
  \begin{itemize}
    \item For $q_2(I)$, we require that the middle coefficient $b$ be even.
    \item For $q^+(I)$, we require that the roots be real, that is, that $b^2 - 4ac > 0$.
    \item For $q_2^+(I)$, we impose both of the last two conditions.
  \end{itemize}
\end{defn}

We are now ready to state a \emph{quadratic reflection theorem,} the main result of this section.
\begin{thm}[\textbf{``Quadratic O-N''}] \label{thm:quadratic_Z}
For every nonzero integer $n$,
\begin{align*}
  q_2^+(4n) &= q(n) \\
  q_2(4n) &= 2 q^+(n).
\end{align*}
\end{thm}
\begin{proof}
The proof is not easy. See Theorem \ref{thm:O-N_quad_Z}.
\end{proof}

It's not hard to compute all quadratics of a fixed superdiscriminant $I$. The leading coefficient $a$ must be a divisor of $I$ (possibly negative), and there are only finitely many of these. Then, by replacing $x$ by $x + t$ where $t$ is an integer nearest to $-b/(2a)$, we can assume that $b$ lies in the window $-|a| < b \leq |a|$. We can try each of the integer values in this window, checking whether
\[
c = \frac{ab^2 - I}{4a^2}
\]
comes out to an integer.
\begin{examp}\label{examp:15}
There are five quadratics of superdiscriminant $15$:
\[
\begin{tabular}{l|cccc}
  $f(x)$ & $q$ & $q^+$ & $q_2$ & $q_2^+$ \\ \hline
  ${-x^2} + x - 4 $ & $\checkmark$ & & &    \\
  $15x^2 + x      $ & $\checkmark$ & $\checkmark$ & & \\
  $15x^2 - x      $ & $\checkmark$ & $\checkmark$ & & \\
  $15x^2 + 11x + 2$ & $\checkmark$ & $\checkmark$ & & \\
  $15x^2 - 11x + 2$ & $\checkmark$ & $\checkmark$ & & 
\end{tabular}
\]
You might think we left out $-x^2 - x - 4$, but it is equivalent to another quadratic on the list:
\[
  -x^2 - x - 4 = -(x + 1)^2 + (x + 1) - 4.
\]
So we get the totals
\[
  q(15) = 5 \textand
  q^+(15) = 4.
\]
There are $18$ quadratics of superdiscriminant $60$:
\[
\begin{tabular}{l|cccc}
  $f(x)$ & $q$ & $q^+$ & $q_2$ & $q_2^+$ \\ \hline
  ${x^2 - 15        }$ & $\checkmark$ & $\checkmark$ & $\checkmark$ &  $\checkmark$ \\
  ${-x^2 - 15       }$ & $\checkmark$ & & $\checkmark$ & \\
  ${-3x^2 + 2x - 2  }$ & $\checkmark$ & & $\checkmark$ & \\
  ${-3x^2 - 2x - 2  }$ & $\checkmark$ & & $\checkmark$ & \\
  ${-4x^2 + x - 1   }$ & $\checkmark$ & & & \\
  ${-4x^2 - x - 1   }$ & $\checkmark$ & & & \\
  ${15x^2 + 2x      }$ & $\checkmark$ & $\checkmark$ & $\checkmark$ &  $\checkmark$ \\
  ${15x^2 - 2x      }$ & $\checkmark$ & $\checkmark$ & $\checkmark$ &  $\checkmark$ \\
  ${15x^2 + 8x + 1  }$ & $\checkmark$ & $\checkmark$ & $\checkmark$ &  $\checkmark$ \\
  ${15x^2 - 8x + 1  }$ & $\checkmark$ & $\checkmark$ & $\checkmark$ &  $\checkmark$ \\
  ${60x^2 + x       }$ & $\checkmark$ & $\checkmark$ & & \\
  ${60x^2 - x       }$ & $\checkmark$ & $\checkmark$ & & \\
  ${60x^2 + 31x + 4 }$ & $\checkmark$ & $\checkmark$ & & \\
  ${60x^2 - 31x + 4 }$ & $\checkmark$ & $\checkmark$ & & \\
  ${60x^2 + 41x + 7 }$ & $\checkmark$ & $\checkmark$ & & \\
  ${60x^2 - 41x + 7 }$ & $\checkmark$ & $\checkmark$ & & \\
  ${60x^2 + 49x + 10}$ & $\checkmark$ & $\checkmark$ & & \\
  ${60x^2 - 49x + 10}$ & $\checkmark$ & $\checkmark$ & &
\end{tabular}
\]
Counting carefully, we get
\[
  q(60) = 18, \quad q_2(60) = 8, \quad q^+(60) = 13, \quad q_2^+(60) = 5.
\]
The equalities
\[
  q_2^+(60) = 5 = q(15) \textand q_2(60) = 8 = 2\cdot 4 = 2q^+(15)
\]
are instances of Theorem \ref{thm:quadratic_Z}. From the same theorem, we derive, without computation, that
\[
  q_2^+(240) = q(60) = 18 \textand q_2(240) = 2q^+(60) = 26.
\]
\end{examp}

This short investigation raises many questions. The superdiscriminant $I = a(b^2 - 4ac)$ does not seem to have been considered before. Is there an explicit formula for $q(I)$? Is there an elementary proof of Theorem \ref{thm:quadratic_Z}? See Example \ref{ex:QR} for a connection to Gauss's celebrated law of quadratic reciprocity.

\subsection{Reflection for cubic equations}
\begin{defn}
For a cubic polynomial
\[
  f(x) = ax^3 + bx^2 + cx + d,
\]
we define the \emph{discriminant} to be
\begin{equation} \label{eq:cubic_disc_rts}
  \disc f = a^4(x_1 - x_2)^2(x_1 - x_3)^2(x_2 - x_3)^2,
\end{equation}
where $x_1, x_2, x_3$ are the roots. Explicitly,
\begin{equation} \label{eq:cubic_disc_coefs}
  \disc f = b^2c^2 - 4ac^3 - 4b^3d - 27a^2d^2 + 18abcd.
\end{equation}
\end{defn}

There are many transformations of a cubic polynomial that don't change the discriminant. One is changing $x$ to $x + t$, where $t$ is a constant. Another is reversing the coefficients,
\[
  f(x) = ax^3 + bx^2 + cx + d \longmapsto x^3 f\(\frac{1}{x}\) = dx^3 + cx^2 + bx + a.
\]
Both of these are special cases of the following construction.

\begin{defn} \label{defn:cubic_eqvt}
Two cubic polynomials $f_1$, $f_2$ with integer coefficients are \emph{equivalent} if there is a matrix
\[
  \begin{bmatrix}
    p & q \\
    r & s
  \end{bmatrix}
\]
whose determinant $ps - qr$ is $\pm 1$ such that
\[
  f_2(x) = (rx + s)^3 \cdot f_1\(\frac{px + q}{rx + s}\).
\]
A matrix that makes $f$ equivalent to \emph{itself,} that is,
\[
f(x) = (rx + s)^3 \cdot f\(\frac{px + q}{rx + s}\),
\]
is called a \emph{symmetry} of $f$. The number of symmetries of $f$ is denoted by $s(f)$.
\end{defn}

\begin{defn}
If $D$ is a nonzero integer, define $h(D)$ to be the number of cubic polynomials
\[
  f(x) = ax^3 + bx^2 + cx + d
\]
of discriminant $D$, up to equivalence, each $f$ counted not once but $1/s(f)$ times, where $s(f)$ is the number of symmetries. Define $h_3(D)$ to be the number of cubics of discriminant $D$ for which the middle two coefficients, $b$ and $c$, are multiples of $3$, up to equivalence, each $f$ counted $1/s(f)$ times as before.
\end{defn}

We can now state the Ohno-Nakagawa reflection theorem that got this research project started:
\begin{thm}[Ohno-Nakagawa; Theorem \ref{thm:O-N}] \label{thm:O-N_Z}
For every nonzero integer $D$,
\[
  h_3(-27 D) = \begin{cases}
    3 h(D), & D > 0 \\
    h(D), & D < 0.
  \end{cases}
\]
\end{thm}
\begin{proof}
Several proofs are in print (see the Introduction). In this paper, we prove this theorem as a special case of Theorem \ref{thm:O-N_traced}.
\end{proof}

\begin{examp}
Take $D = 1$. There is just one cubic with integer coefficients and discriminant $1$, namely
\[
  f(x) = x(x + 1) = x^2 + x.
\]
The reader may balk at considering a quadratic polynomial as a ``cubic'' with leading coefficient $0$, but the polynomial can be replaced by any number of equivalent forms, for instance
\[
  (x-1)^3 \cdot f\(\frac{x}{x-1}\) = x(x - 1)(2x - 1).
\]
We will suppress this detail in subsequent examples.
(A program for computing all cubics of a given discriminant is found in the attached file \verb|cubics.sage|, based on an algorithm of Cremona \cite{Crem_Redn, Crem_Redn_2}). The cubic $f$ has six symmetries, which is related to the fact that three linear factors can be permuted in $3! = 6$ ways. In terms of $f(x) = x(x + 1)$, the symmetries are
\[
  \begin{bmatrix}
    1 & 0 \\
    0 & 1
  \end{bmatrix},
  \begin{bmatrix}
  -1 & -1 \\
  0 & 1
\end{bmatrix},
  \begin{bmatrix}
  0 & 1 \\
  1 & 0
\end{bmatrix},
  \begin{bmatrix}
  -1 & -1 \\
  1 & 0
\end{bmatrix},
  \begin{bmatrix}
  1 & 0 \\
  -1 & -1
\end{bmatrix},
  \begin{bmatrix}
  0 & 1 \\
  -1 & -1
\end{bmatrix}.
\]
So $h(1) = 1/6$.

Correspondingly, we look at cubics of discriminant $-27$. There are two:
\begin{equation*}
  f(x) = x^2 + x + 7 \textand f(x) = x^3 + 1.
\end{equation*}
Each admits two symmetries: the first has
\[
\begin{bmatrix}
  1 & 0 \\
  0 & 1
\end{bmatrix},
\begin{bmatrix}
  -1 & -1 \\
  0 & 1
\end{bmatrix},
\]
and the second has
\[
\begin{bmatrix}
  1 & 0 \\
  0 & 1
\end{bmatrix},
\begin{bmatrix}
  0 & 1 \\
  1 & 0
\end{bmatrix}.
\]
So $h(-27) = 1/2 + 1/2 = 1$ and $h_3(-27) = 1/2$. In particular,
\[
  h_3(-27) = 3 h_3(1),
\]
in conformity with Theorem \ref{thm:O-N_Z}.
\end{examp}

\subsection{Reflection for \texorpdfstring{$2\times n\times n$}{2 x n x n} boxes}
Bhargava \cite{B3} studied $2 \times 3 \times 3$ boxes as a visual representation for quartic rings, as cubic polynomials do for cubic rings. We think that reflection holds not only for $2\times 3\times 3$ boxes but for $2\times 5\times 5$, $2\times 7\times 7$, and so on. We nearly prove the $2\times 3\times 3$ case in this paper. We are quite far from proving it for the larger boxes.

\begin{defn}
A \emph{box} is a pair $(A, B)$ of $n \times n$ integer symmetric matrices. The \emph{resolvent} of a box is the polynomial
\[
  f(x) = \det(Ax - B).
\]
\end{defn}
It is a polynomial in $x$, of degree at most $n$. If $A$ is the identity matrix, the resolvent devolves into the standard \emph{characteristic polynomial.}
\begin{defn}
Two boxes $(A_1, B_1)$ and $(A_2, B_2)$ are \emph{equivalent} if there is an integer $n \times n$ matrix $X$, whose inverse $X^{-1}$ also has integer entries, such that
\[
  A_2 = XA_1X^\top \textand B_2 = XB_1X^\top.
\]
If $(A_2, B_2) = (A_1, B_1) = (A, B)$ are the same pair, then $X$ is called a \emph{symmetry} of $(A, B)$. The number of symmetries of $(A, B)$ will be denoted by $s(A,B)$.
\end{defn}

\begin{conj}[\textbf{``O-N for $2\times n \times n$ boxes''}]\label{conj:O-N_2xnxn_Z}
Let $n$ be a positive odd integer. Let $f$ be a polynomial of degree $n$ with no multiple roots and only one real root. Denote by $h(f)$ the number of $2\times n \times n$ boxes with resolvent $f$, up to equivalence, each box weighted by the reciprocal of its number of symmetries. Denote by $h_2(f)$ the number of such boxes with even numbers along the main diagonals of $A$ and $B$, weighted the same way. Then
\begin{equation}
  h_2(2^{n-1} f) = 2^{\frac{n-1}{2}} \cdot h(f).
\end{equation}
\end{conj}
\begin{rem}
The condition that $f$ have no multiple roots (even complex ones) is needed to ensure that there are only finitely many boxes with $f$ as a resolvent. The condition that $f$ have no more than one real root can be eliminated, but then we must impose conditions on the real behavior of the boxes that are difficult to state succinctly.
\end{rem}
\begin{examp} \label{examp:-23_boxes}
Take as resolvent $f(x) = x^3 - x - 1$, the simplest irreducible cubic. It has one real root $\xi \approx 1.3247$ and discriminant $-23$. There are two boxes with resolvent $f$, up to equivalence:
\[
\left[\left(\begin{array}{rrr}
  0 & 0 & -1 \\
  0 & -1 & 0 \\
  -1 & 0 & -1
\end{array}\right), \left(\begin{array}{rrr}
  0 & -1 & 0 \\
  -1 & 0 & -1 \\
  0 & -1 & -1
\end{array}\right)\right],
\left[\left(\begin{array}{rrr}
  0 & -1 & 0 \\
  -1 & 0 & -1 \\
  0 & -1 & -1
\end{array}\right), \left(\begin{array}{rrr}
  -1 & 0 & -1 \\
  0 & -1 & -1 \\
  -1 & -1 & -1
\end{array}\right)\right].
\]
(These were computed from the balanced pairs $(\OO_R, 1)$ and $(\OO_R, \xi)$ in the number field $R = \ZZ[\xi]/(\xi^3 - \xi - 1)$ corresponding to $f$.) Neither has any symmetries besides the two trivial ones, the identity matrix and its negative, so 
\[
  h(f) = \frac{1}{2} + \frac{1}{2} = 1.
\]
There are many boxes with resolvent $2f$, but just one with even numbers all along the main diagonals of $A$ and $B$, namely
\[
\left[\left(\begin{array}{rrr}
  0 & 0 & 1 \\
  0 & -2 & 0 \\
  1 & 0 & 2
\end{array}\right), \left(\begin{array}{rrr}
  0 & 1 & 0 \\
  1 & 0 & 0 \\
  0 & 0 & -2
\end{array}\right)\right].
\]
(This was computed from the unique quartic ring $\OO = \ZZ \cross \OO_R$ with resolvent $\OO_R$.) It too has only the trivial symmetries, to $h_2(2f) = 1/2$, in accord with Conjecture \ref{conj:O-N_2xnxn_Z}.

\end{examp}

\subsection{Reflection for quartic equations}
There are also reflection theorems that appear when counting \emph{quartic} polynomials.
\begin{defn}
If
\[
  f(x) = a x^4 + b x^3 + c x^2 + d x + e
\]
is a quartic polynomial with integer coefficients, its \emph{resolvent} is
\begin{equation}
  g(y) = y^{3} - c y^{2} + (b d - 4 a e) y + 4 a c e - b^2 e - a d^2;
\end{equation}
equivalently, if
\[
  f(x) = a(x - x_1)(x - x_2)(x - x_3)(x - x_4),
\]
then
\[
  g(y) =
  \big(y - a (x_1 x_2 + x_3 x_4)\big)
  \big(y - a (x_1 x_3 + x_2 x_4)\big)
  \big(y - a (x_1 x_4 + x_2 x_3)\big).
\]
\end{defn}
\begin{rem}
Cubic resolvents of this type have been used since the 16th century as a step in solving quartic equations. For instance, it is well known that if $f(x)$ factors as the product of two quadratics with integer coefficients, then $g(y)$ has a rational root (the converse is not true).
\end{rem}
Analogously to Definition \ref{defn:cubic_eqvt}, we put:
\begin{defn}
  Two quartic polynomials $f_1$, $f_2$ with integer coefficients are \emph{equivalent} if there is a matrix
  \[
  \begin{bmatrix}
    p & q \\
    r & s
  \end{bmatrix}
  \]
  whose determinant $ps - qr$ is $\pm 1$ such that
  \[
  f_2(x) = (rx + s)^4 \cdot f_1\(\frac{px + q}{rx + s}\).
  \]
  A matrix that makes $f$ equivalent to \emph{itself,} that is,
  \[
  f(x) = (rx + s)^4 \cdot f\(\frac{px + q}{rx + s}\),
  \]
  is called a \emph{symmetry} of $f$. The number of symmetries of $f$ is denoted by $s(f)$.
\end{defn}
We have:
\begin{lem}
\begin{enumerate}[$($a$)$]
  \item If two quartics $f_1$, $f_2$ are equivalent, then their resolvents $g_1$, $g_2$ are related by a translation
  \[
    g_2(x) = g_1(x + t)
  \]
  for some integer $t$.
  \item A quartic and its resolvent have the same discriminant
  \[
    \disc f = \disc g = 
    \parbox[t]{0.7\linewidth}{
    $b^2 c^2 d^2 - 4 a c^3 d^2 - 4 b^3 d^3 + 18 a b c d^3 - 27 a^2 d^4 - 4 b^2 c^3 e + 16 a c^4 e + 18 b^3 c d e - 80 a b c^2 d e - 6 a b^2 d^2 e + 144 a^2 c d^2 e - 27 b^4 e^2 + 144 a b^2 c e^2 - 128 a^2 c^2 e^2 - 192 a^2 b d e^2 + 256 a^3 e^3.$
    }
  \]
\end{enumerate}
\end{lem}
\begin{proof}
Exercise.
\end{proof}

As before, our reflection theorem will relate general quartics to quartics satisfying certain divisibility relations. Here the relations are quite peculiar:
\begin{defn}
A quartic polynomial
\[
  f(x) = a x^4 + b x^3 + c x^2 + d x + e
\]
is called \emph{supereven} if $b$, $c$, and $e$ are multiples of $4$ and $d$ is a multiple of $8$.
\end{defn}

Not every quartic equivalent to a super-even quartic is itself supereven. (For instance, $f_1 = x^4 + 4$ and $f_2 = 4x^4 + 1$ are equivalent under the flip
$\big[\begin{smallmatrix}
  0 & 1 \\
  1 & 0
\end{smallmatrix}\big]$
, but $f_2$ is not supereven.) We therefore make the following definition.
\begin{defn}
  Two quartic polynomials $f_1$, $f_2$ with integer coefficients are \emph{evenly equivalent} if there is a matrix
  \[
  \begin{bmatrix}
    p & q \\
    r & s
  \end{bmatrix}
  \]
  whose determinant $ps - qr$ is $\pm 1$, \textbf{and \emph{r} is even,} such that
  \[
  f_2(x) = (rx + s)^4 \cdot f_1\(\frac{px + q}{rx + s}\).
  \]
  Such a matrix that makes $f$ equivalent to \emph{itself,} that is,
  \[
  f(x) = (rx + s)^4 \cdot f\(\frac{px + q}{rx + s}\),
  \]
  is called an \emph{even symmetry} of $f$. The number of even symmetries of $f$ is denoted by $s_2(f)$.
\end{defn}

\begin{thm}[\textbf{``Quartic O-N''}]
  Let $g$ be an integer cubic with leading coefficient $1$, no multiple roots, and odd discriminant. Denote by $h(g)$ the number of quartics whose resolvent is $g(y + t)$ for some $t$, up to equivalence and weighted by the reciprocal of the number of symmetries. Denote by $h_2(g)$ the number of supereven quartics whose resolvent is $g(y + t)$ for some $t$, up to even equivalence and weighted by the reciprocal of the number of even symmetries. Define $g_2$ by
  \[
    g_2(y) = 64 g\(\frac{y}{4}\)
  \]
   Then:
  \begin{itemize}
    \item If $g$ has one real root, then
    \[
    4 h(g) = h_2(g_2).
    \]
    \item If $g$ has three real roots, then we subdivide
    \[
      h(g) = h^+(g) + h^-(g) + h^\pm(g)
    \]
    where the respective terms count only quartic functions that are always positive, always negative, and have four real roots. We subdivide
    \[
      h_2(g) = h^+_2(g) + h^-_2(g) + h^\pm_2(g).
    \]
    Then:
    \begin{align*}
      2 h(g) &= h_2^{\pm}(g_2) \\
      4\big(h^{+}(g) + h^\pm(g)\big) &= h_2^{+}(g_2) + h_2^\pm(g_2) \\
      4\big(h^{-}(g) + h^\pm(g)\big) &= h_2^{-}(g_2) + h_2^\pm(g_2) \\
    \end{align*}
  Also, denote by $k(g)$ the number of integral $3\times 3$ symmetric matrices of characteristic polynomial $g$. Then
  \[
    k(g) = 24\big(h^{\pm}(g) - h^{+}(g) - h^-(g)\big).
  \]
  \end{itemize}
\end{thm}
\begin{proof}
See Theorem \ref{thm:BQ}.
\end{proof}
\begin{rem}
We think that the hypothesis of odd discriminant is removable, but we have not yet finished the proof.
\end{rem}
\begin{examp} \label{ex:bq1}
Let $g(y) = y^3 - y - 1$. By techniques presented in Section \ref{sec:bq}, it is possible to transform the boxes found in example \ref{examp:-23_boxes} into binary quartic forms. We find that there is only one quartic with resolvent $g$, namely
\[
  f(x) = x^3 - x - 1
\]
(which, as before, can be transformed by an equivalence to one with nonzero leading coefficient); and four supereven binary quartics with resolvent $g_2(y) = y^3 - 16y - 64$, namely
\begin{align*}
  f(x) &= 4x^3 + 12x^2 + 8x - 4 = 4\big((x + 1)^3 - (x + 1) - 1\big) \\
  f(x) &= -x^4 + 4x^3 + 12x^2 + 8x \\
  f(x) &= -x^4 + 8x - 4 \\
  f(x) &= -x^4 + 4x^3 - 4.
\end{align*}
\end{examp}
All these have one pair of complex roots (as must occur for a resolvent with negative discriminant) and only the trivial symmetries
$\pm\big[\begin{smallmatrix}
  1 & 0 \\
  0 & 1
\end{smallmatrix}\big]$, so
\[
  h(g) = \frac{1}{2} \textand h_2(g_2) = 2 = 4 \cdot \frac{1}{2},
\]
in accord with the first part of the theorem.
\begin{examp} \label{ex:bq2}
Consider $g(y) = y^3 - 2y^2 - 3y + 6 = (y - 2)(y + \sqrt{3})(y - \sqrt{3})$, a cubic with three real roots. The quartics with resolvent $g$ are
\[
  f(x) = -x(2x - 1)(3x^2 - 1),
\]
which has four real roots, and
\[
  f(x) = (x^2 + x + 1)(x^2 + 1),
\]
which has no real roots and is positive for all real $x$. Each has only the trivial symmetries, so
\[
  h^\pm(g) = \frac{1}{2}, \quad h^+(g) = \frac{1}{2}, \quad h^-(g) = 0.
\]
(Note the discrepancy between $h^+$ and $h^-$.) Correspondingly, there are eight supereven binary quartics with resolvent $g_2(y) = (y - 8)(y + 4\sqrt{3})(y - 4\sqrt{3})$:
\begin{align*}
  f(x) &= -2x^4 - 8x^3 - 4x^2 + 8x            = -2  x  (x + 2)  (x^2 + 2x - 2) \\
  f(x) &= 4x^3 - 4x^2 - 16x - 8           = 4  (x + 1)  (x^2 - 2x - 2) \\
  f(x) &= 8x^4 - 16x^3 + 20x^2 - 12x + 4 = 4  (x^2 - x + 1)  (2x^2 - 2x + 1) \\
  f(x) &= x^4 - 6x^3 + 20x^2 - 32x + 32   = (x^2 - 4x + 8)  (x^2 - 2x + 4) \\
  f(x) &= -x^4 + 8x^2 - 12                   = -(x^2 - 2)(x^2 - 6) \\
  f(x) &= -3x^4 + 8x^2 - 4                  = -(x^2 - 2)(3x^2 - 2) \\
  f(x) &= x^4 + 8x^2 + 12                   = (x^2 + 2)(x^2 + 6) \\
  f(x) &= 3x^4 + 8x^2 + 4                  = (x^2 + 2)(3x^2 + 2).
\end{align*}
Thus
\[
  h_2^\pm(g_2) = 2, \quad h_2^+(g_2) = 2, \quad h_2^-(g_2) = 0.
\]
This is in accord with the theorem, from which we also learn that
\[
  k(g) = 48\big(h^{\pm}(g) - h^{+}(g) - h^-(g)\big) = 0,
\]
so $g$ is \emph{not} the characteristic polynomial of any integer $3 \times 3$ symmetric matrix, despite having three real roots (which is a necessary, but not a sufficient, condition).
\end{examp}
\begin{examp}
Let $g(y) = y^3 - y$. Knowing that $f(x) = x^3 - x$ is the only quartic with cubic resolvent $f$, and it has four symmetries, the powers of 
$\big[\begin{smallmatrix}
  0 & -1 \\
  1 & 0
\end{smallmatrix}\big]$, we get
\[
  k(g) = 24\(h^\pm(g) - h^+(g) - h^-(g)\) = 24\(\frac{1}{4} - 0 - 0\) = 6.
\]
So there are six symmetric matrices with characteristic polynomial $y^3 - y$. Indeed, they are the diagonal matrices with $1$, $0$, and $-1$ along the diagonal in any of the $3! = 6$ possible orders.
\end{examp}

\section{Notation}\label{sec:notation}
The following conventions will be observed in the remainder of the paper.

We denote by $\NN$ and $\NN^+$, respectively, the sets of nonnegative and of positive integers.

If $P$ is a statement, then
\[
  \1_P = \begin{cases}
    1 & \text{$P$ is true} \\
    0 & \text{$P$ is false}.
  \end{cases}
\]
If $S$ is a set, then $\1_S$ denotes the characteristic function $\1_S(x) = \1_{x \in S}$.

An \emph{algebra} will always be commutative and of finite rank over a field, while a \emph{ring} or \emph{order} will be a finite-dimensional, torsion-free ring over a Dedekind domain, containing $1$. An order need not be a domain.

If $a,b \in L$ are elements of a local or global field, a separable closure thereof, or a finite product of the preceding, we write $a|b$ to mean that $b = c a$ for some $c$ in the appropriate ring of integers $\OO_L$. If $a|b$ and $b|a$, we say that $a$ and $b$ are \emph{associates} and write $a \sim b$. Note that $a$ and $b$ may be zero-divisors.

If $S$ is a finite set, we let $\Sym(S)$ denote the set of permutations of $S$; thus $S_n = \Sym(\{1,\ldots,n\})$. If $\size{S} = \size{T}$, and if $\g \in \Sym(S)$, $\h \in \Sym(T)$ are elements, we say that $\g$ and $\h$ are \emph{conjugate} if there is a bijection between $S$ and $T$ under which they correspond. Likewise when we say that two subgroups $G \subseteq \Sym(S)$, $H \subseteq \Sym(T)$ are conjugate.

We will use the semicolon to separate the coordinates of an element of a product of rings. For instance, in $\ZZ \cross \ZZ$, the nontrivial idempotents are $(1;0)$ and $(0;1)$.

If $n$ is a positive integer, then $\zeta_n$ denotes a primitive $n$th root of unity in $\bar\QQ$, while $\bar\zeta_n$ denotes the $n$th root of unity
\[
\bar\zeta_n = \(1; \zeta_n; \zeta_n^2; \ldots; \zeta_n^{n-1}\) \in \bar\QQ^n.
\]

Throughout the proofs of the local reflection theorems, we will fix a local field $K$, its valuation $v = v_K$, its residue field $k_K$ of order $q$, and a uniformizer $\pi = \pi_K$. The letter $e$ will denote the absolute ramification index ($e = v_K(2)$ in the quadratic and quartic cases, $v_K(3)$ in the cubic). We let $\mm_K$ denote the maximal ideal, and likewise $\mm_{\bar K}$ be the maximal ideal of the ring $\OO_{\bar K}$ of algebraic integers over $K$; note that $\mm_{\bar K}$ is not finitely generated. We also allow $v = v_K$ to be applied to elements of $\bar K$, the valuation being scaled so that its restriction to $K$ has value group $\ZZ$. We use the absolute value bars $\size{\bullet}$ for the corresponding metric, whose normalization will be left undetermined.

If $K$ is a local field, an \emph{$m$-pixel} is a subset of an affine or projective space over $\OO_K$ defined by requiring the coordinates to lie in specified congruence classes modulo $\pi^n$. For instance, in $\PP^2(\OO_K)$, a $0$-pixel is the whole space, which is subdivided into $(q^2+q+1)q^{2n - 2}$-many $n$-pixels for each $n \geq 1$.

If $R/K$ is a finite-dimensional, locally free algebra over a ring, we denote by $R^{N=1}$ the subgroup of units of norm $1$. The group operation is implicitly multiplication, so $R^{N=1}[n]$, for instance, denotes the $n$th roots of unity of norm $1$.

\part{Galois cohomology}
\label{part:Gal_coho}

\section{\'Etale algebras and their Galois groups}
\label{sec:etale}

\subsection{\'Etale algebras}
If $K$ is a field, an \emph{\'etale algebra} over $K$ is a finite-dimensional separable commutative algebra over $K$, or equivalently, a finite product of finite separable extension fields of $K$. A treatment of \'etale algebras is found in Milne (\cite{MilneFields}, chapter 8): here we summarize this theory and prove a few auxiliary results that will be of use.

An \'etale algebra $L$ of rank $n$ admits exactly $n$ maps $\iota_1,\ldots, \iota_n$ (of $K$-algebras) to a fixed separable closure $\bar K$ of $K$. We call these the \emph{coordinates} of $L$; the set of them will be called $\Coord(L/K)$ or simply $\Coord(L)$. Together, the coordinates define an embedding of $L$ into $\bar K^n$, which we call the \emph{Minkowski embedding} because it subsumes as a special case the embedding of a degree-$n$ number field into $\CC^n$, which plays a major role in algebraic number theory, as in Delone-Faddeev \cite{DF}.

For any element $\gamma$ of the absolute Galois group $G_K$, the composition $\gamma \circ \iota_i$ with any coordinate is also a coordinate $\iota_j$, so we get a homomorphism $\phi = \phi_L : G_K \to \Sym(\Coord(L))$) such that
\[
  \gamma(\iota(\x)) = (\phi_\gamma\iota)(\x)
\]
for all $\x \in L, \iota \in \Coord(L)$. This gives a functor from \'etale $K$-algebras to $G_K$-sets (sets with a $G_K$-action), which is denoted $\F$ in Milne's terminology. A functor going the other way, which Milne calls $\A$, takes $\phi : G_K \to S_n$ to
\begin{equation} \label{eq:Gset_to_etale}
  L = \{ (x_1,\ldots,x_n) \in \bar K^n \mid \gamma(x_i) = x_{\phi_\gamma(i)} \, \forall \gamma \in G_K, \forall i \}
\end{equation}

\begin{prop}[\cite{MilneFields}, Theorem 7.29] \label{prop:G_sets}
The functors $\F$ and $\A$ establish a bijection between
\begin{itemize}
  \item \'etale extensions $L/K$ of degree $n$, up to isomorphism, and
  \item $G_K$-sets of size $n$ up to isomorphism; that is to say, homomorphisms $\phi : G_K \to S_n$, up to conjugation in $S_n$.
\end{itemize}
\end{prop}

Moreover, the bijection respects base change, in the following way:

\begin{prop}\label{prop:etale_respects_base_chg}
Let $K_1/K$ be a field extension, not necessarily algebraic, and let $L/K$ be an \'etale extension of degree $n$. Then $L_1 = L \tensor_{K} K_1$ is \'etale over $K_1$, and the associated Galois representations $\phi_{L/K}$, $\phi_{L_1/K_1}$ are related by the commutative diagram
\begin{equation}\label{eq:etale_respects_base_chg}
\xymatrix{
  G_{K_1} \ar[r]^{\bullet |_{\bar K}} \ar[d]_{\phi_{L_1/K_1}} & G_K \ar[d]^{\phi_{L/K}} \\
  \Sym(\Coord_{K_1}(L_1)) \ar@{-}[r]^{\sim} & \Sym(\Coord_K(L))
}
\end{equation}
\end{prop}
\begin{proof}
That $L_1/K_1$ is \'etale is standard (see Milne \cite{MilneFields}, Prop.~8.10). For the second claim, consider the natural restriction map $r : \Coord_{K_1}(L_1) \to \Coord_K(L)$. It is injective, since a $K_1$ linear map out of $L_1$ is determined by its values on $L$; and since both sets have the same size, $r$ is surjective and is hence an isomorphism of $G_{K_1}$-sets (the $G_{K_1}$-structure on $\Coord_{K}(L)$ arising by restriction from the $ G_{K} $-structure).
\end{proof}

We will use this proposition most frequently in the case that $K$ is a global field and $K_1 = K_v$ one of its completions. The resulting $L_1$ is then the product $L_v \cong \prod_{w|v} L_w$ of the completions of $L$ at the places dividing $v$. Note the departure from the classical habit of studying the completion $L_w$ at each place individually. The preservation of degrees, $[L_1 : K_1] = [L : K]$ will be important for our applications.

\subsection{The Galois group of an \'etale algebra}

Define the \emph{Galois group} $G(L/K)$ of an \'etale algebra to be the image of its associated Galois representation $\phi : G_K \to \Sym(\Coord(L))$. It transitively permutes the coordinates corresponding to each field factor. For example, if $L$ is a quartic field, then $G(L/K)$ is one of the five (up to conjugacy) transitive subgroups of $\sf S_4$, which (to use the traditional names) are $\sf S_4$, $\sf A_4$, $\sf D_4$, $\sf V_4$, and $\sf C_4$. Galois groups in this sense are used in the tables of cubic and quartic fields in Delone-Faddeev \cite{DF} and the Number Field Database \cite{NFDB}. Note that the Galois group $G(L/K)$ is defined whether or not $L$ is a Galois extension. If it is, then the Galois group is \emph{simply} transitive and coincides with the Galois group in the sense of Galois theory.

Important for us will be two notions pertaining to the Galois group.
\begin{defn}
Let $G \subseteq \S_n$ be a subgroup. A \emph{$G$-extension} of $K$ is a degree-$n$ \'etale algebra $L$ with a choice of subgroup $G' \subseteq \Sym(\Coord(L))$ that is conjugate to $G$ and contains $G(L/K)$, plus a conjugacy class of isomorphisms $G' \cong G$: the conjugacy being in $G$, not in $\S_n$. The added data is called a \emph{$G$-structure} on $L$.
\end{defn}
\begin{prop} \label{prop:G_sets_G_strucs}
$G$-extensions $L/K$ up to isomorphism are in bijection with homomorphisms $\phi : G_K \to G$, up to conjugation in $G$.
\end{prop}
\begin{proof}
Immediate from Proposition \ref{prop:G_sets}.
\end{proof}

\begin{examp}\label{ex:C4-structure}
$L = \QQ(\zeta_5)$ is a $\sf C_4$-extension (taking $\sf C_4 = \<(1234)\> \subseteq \sf S_4$), indeed its Galois group is isomorphic to $\sf C_4$; and $L$ admits two distinct $\sf C_4$-structures, as there are two ways to identify $\sf C_4$ with its image in $\sf S_4$, which are conjugate in $\sf S_4$ but not in $\sf C_4$. Likewise, $L = \QQ \cross \QQ \cross \QQ \cross \QQ$ admits six $\sf C_4$-structures, one for each embedding of $\sf C_4$ into $\sf S_4$, as its Galois group is trivial.
\end{examp}

\subsection{Resolvents}
This will be an important notion.

\begin{defn}
Let $G \subseteq \sf S_n$, $H \subseteq \S_m$ be subgroups and $\rho : G \to H$ be a homomorphism. Then for every $G$-extension $L/K$, the corresponding $\phi_L : G_K \to G$ may be composed with $\rho$ to yield a map $\phi_R : G_K \to H$, which defines an \'etale extension $R/K$ of degree $m$. This $R$ is called the \emph{resolvent} of $L$ under the map $\rho$.
\end{defn}

\begin{examp}\label{ex:rsv43}
Since there is a surjective map $\rho_{4,3} : \sf S_4 \to \sf S_3$, every quartic \'etale algebra $L/K$ has a cubic resolvent $R$. This resolvent appears in Bhargava \cite{B3}, but it is much older than that. It is generated by a formal root of the \emph{resolvent cubic} that appears when a general quartic equation is to be solved by radicals.
\end{examp}

\begin{examp} \label{ex:rsv_n2}
Likewise, the sign map can be viewed as a homomorphism $\sgn : \sf S_n \to \sf S_2$, attaching to every \'etale algebra $L$ a quadratic resolvent $T$. If $L = K[\theta]/f(\theta)$ is generated by a polynomial $f$, and if $\ch K \neq 2$, then it is not hard to see that $T = K[\sqrt{\disc f}]$ where $\disc f$ is the polynomial discriminant. Note that $T$ still exists even if $\ch K = 2$. We have that $T \cong K \cross K$ is split if and only if the Galois group $G(L/K)$ is contained in the alternating group $\sf A_n$.
\end{examp}

\begin{examp}\label{ex:D4}
The dihedral group $\sf D_4$ has an outer automorphism, because rotating a square in the plane by $45^\circ$ does not preserve the square but does preserve every symmetry of the square. This map $\rho : \sf D_4 \to \sf D_4$ associates to each $\sf D_4$-algebra $L$ a new $\sf D_4$-algebra $L'$, not in general isomorphic. This is the classical phenomenon of the \emph{mirror field}. For instance, if $L = \QQ[\sqrt{1 + \sqrt{2}}]$, then
\[
  L' = \QQ\left[\sqrt{1 + \sqrt{2}} + \sqrt{1 - \sqrt{2}}\right] = \QQ\left[\sqrt{2 + 2\sqrt{-1}}\right].
\]
Both $L$ and $L'$ have the same Galois closure, a $\sf D_4$-octic extension of $\QQ$. Likewise, the outer automorphism of $\S_6$ permits the association to each sextic \'etale algebra $L/K$ a mirror sextic \'etale algebra $L'$. 
\end{examp}

\begin{examp}
The \emph{Cayley embedding} is an embedding of any group $G$ into $\Sym(G)$, acting by left multiplication. The Cayley embedding $\rho : \S_n \hookrightarrow \S_{n!}$ attaches to every \'etale algebra $L$ of degree $n$ an algebra $\tilde L$ of degree $n!$ with an $\S_n$-torsor structure. This is none other than the \emph{$\S_n$-closure} of $L$, constructed by Bhargava in a quite different way in \cite[Section 2]{B3}.

More generally, for any $G \subseteq S_n$, the Cayley embedding $G \hookrightarrow \Sym(G)$ allows one to associate to each $G$-extension $L$ a $G$-torsor $T$, which we may call the \emph{$G$-closure} of $L$. The name ``closure'' is justified by the following observation: if $G \subseteq S_n$ is a transitive subgroup, then, since any transitive $G$-set is a quotient of the simply transitive one, we can embed $L$ into $T$ by Proposition \ref{prop:etale_sub} below. More generally, $G$-closures of ring extensions, not necessarily \'etale or even reduced, have been constructed and studied by Biesel \cite{Biesel_thesis, Biesel}.
\end{examp}

If $\rho : G \to H$ is invertible, as in many of the above examples, then the map from $G$-extensions to $H$-extensions is also invertible: we say that the two extensions are \emph{mutual resolvents}.

\subsection{Subextensions and automorphisms}

The Galois group holds the answers to various natural questions about an \'etale algebra. The next two propositions are given without proof, since they follow immediately from the functorial character of the correspondence in Proposition \ref{prop:G_sets}
\begin{prop} \label{prop:etale_sub}
The subextensions $L' \subseteq L$ of an \'etale extension $L/K$, correspond to the equivalence relations $\sim$ on $\Coord(L)$ stable under permutation by $G(L/K)$, under the bijection
\[
  \mathord{\sim} \mapsto L' = \{ \x \in L : \iota(\x) = \iota'(\x) \text{ whenever } \iota \sim \iota' \}.
\]
\end{prop}
\begin{rem}
Note that if $L$ is a Galois field extension, the image of $\phi_L$ is a simply transitive subgroup $\sf \Gamma$, and identifying $\Coord(L)$ with $\sf \Gamma$, the stable equivalence relations are just right congruences modulo subgroups of $\sf \Gamma$: so we recover the Galois correspondence between subgroups and subfields.
\end{rem}

The Galois group is \emph{not} a group of automorphisms of $L$. However, the automorphisms of $L$ as a $K$-algebra can be described in terms of the Galois group readily.

\begin{prop}\label{prop:etale_aut}
Let $L$ be Minkowski-embedded by its coordinates $\iota_1,\ldots,\iota_n$. Then the automorphism group $\Aut(L/K)$ is given by permutations of coordinates,
\[
  \tau_\spi(x_1;\ldots;x_n) = x_{\spi^{-1}(1)} ; \ldots ; x_{\spi^{-1}(n)},
\]
for $\spi$ in the centralizer $C(S_n, G(L/K))$ of the Galois group.
\end{prop}
(For $H \subseteq G$ groups, the \emph{centralizer} $C(G,H)$ of $H$ in $G$ is the subgroup of elements of $G$ that commute with every element of $H$.)

This provides a characterization, in terms of the Galois group, of rings having various kinds of automorphisms.

\begin{itemize}
  \item Since $S_2$ is abelian, any \'etale algebra $L$ of rank $2$ has a unique non-identity automorphism, the conjugation $\bar \x = \tr \x - \x$.
  \item If $L$ has rank $4$, automorphisms $\tau$ of $L$ of order $2$ whose fixed algebra is of rank $2$ are in bijection with $D_4$-structures on $L$. Indeed, the conditions force $\tau$ to correspond to the permutation $\pi = (12)(34)$ or one of its conjugates, and the centralizer of this permutation is $D_4$.
  \item Particularly relevant is the case that $L$ has a complete set of automorphisms that permute the coordinates simply transitively: this is a generalization of a Galois field extension called a \emph{torsor}. This case is sufficiently important to merit its own subsection.
\end{itemize}

\subsection{Torsors}

\begin{defn}
Let $G$ be a finite group. A \emph{$G$-torsor} over $K$ is an \'etale algebra $L$ over $K$ equipped with an action of $G$ by automorphisms $\{\tau_\g\}_{\g \in G}$ that permute the coordinates simply transitively, that is, such that $L \tensor_K \bar K$ is isomorphic to
\[
  \bigoplus_{g \in G} \bar K
\]
with $G$ acting by right multiplication on the indices.
\end{defn}

\begin{prop}
\label{prop:torsor}
Let $G$ be a group of order $n$. An \'etale algebra $L$ is a $G$-torsor if and only if it is a $G$-extension, where $G$ is embedded into $S_n$ by the Cayley embedding ($G$ acting on itself by left multiplication). Moreover, there is a bijection between
\begin{itemize}
  \item $G$-torsor structures on $L$, up to conjugation in $G$, and
  \item $G$-structures on $L$.
\end{itemize}
The bijection is given in the following way: there is a labeling $\{\iota_\g\}$ of the coordinates of $L$ with the elements of $G$ such that the Galois action is by left multiplication
\begin{equation} \label{eq:torsor_left}
  \g(\iota_\h(\x)) = \iota_{\phi_\g \h}(\x)
\end{equation}
while the torsor action is by right multiplication
\begin{equation} \label{eq:torsor_right}
  \iota_\g(\tau_{\h}(\x)) = \iota_{\g \h^{-1}}(\x).
\end{equation}
\end{prop}

\begin{proof}
We first claim that the only elements of $\Sym(G)$ commuting with all right multiplications are left multiplications, and vice versa. If $\pi : G \to G$ is a permutation commuting with left multiplications, then
\[
  \pi(g) = \pi(g \cdot \id_G) = g \cdot \pi(\id_G),
\]
so $\pi$ is a right multiplication. So the embedded images of $G$ in $\Sym(G)$ given by left and right multiplication (which are conjugate under the inversion permutation $\bullet^{-1} \in \Sym(G)$) are centralizers of one another. It is then clear that conjugates $G'$ of $G$ in $\Sym(\Coord(L))$ that \emph{contain} $G(L/K)$ are in bijection with conjugates $G''$ that \emph{commute} with $G(L/K)$. This establishes the first assertion. For the bijection of structures, if an embedding $G \cong G' \subseteq \Sym(\Coord(L))$ is given, then we can label the coordinates with elements of $G$ so that $G$ acts on them by multiplication; then $G''$ gets identified with $G$ by the corresponding right action. The only ambiguity is in which embedding is labeled with the identity element; if this is changed, one computes that the resulting identification of $G''$ with $G$ is merely conjugated, so the map is well defined. The reverse map is constructed in exactly the same way.
\end{proof}

Here is another perspective on torsors.
\begin{prop} \label{prop:fld_fac}
$G$-torsors over a field $K$, up to isomorphism, are determined by their field factor, a Galois extension $L_1/K$ equipped with an embedding $\Gal(L/K) \hookrightarrow G$ up to conjugation in $G$.
\end{prop}
\begin{proof}
If $T$ is a $G$-torsor, then since $G$ permutes the coordinates simply transitively, all the coordinates have the same image; that is, the field factors of $G$ are all isomorphic to a Galois extension $L/K$. The torsor operations fixing one field factor $L_i$ of $T$ realize the Galois group $\Gal(L/K)$ as a subgroup of $G$; changing the field factor $L_i$ and/or the identification $L_i \cong L$ corresponds to conjugating the map $\Gal(L/K) \hookrightarrow G$ by an element of $G$.

Conversely, suppose $L$ and an embedding
\[
  \Gal(L/K) \longto^\sim H \subseteq G
\]
are given. Let $1 = \g_1,\ldots, \g_r$ be coset representatives for $G/H$. Then $\g_2,\ldots, \g_r$ must map any field factor $L_1 \cong L$ isomorphically onto the remaining field factors $L_2,\ldots,L_r$, each $L_i$ occurring once. To finish specifying the $G$-action on $T \cong L_1 \cross \cdots \cross L_r$, it suffices to determine $\g|_{L_i}$ for each $\g \in G$. Factor $\g\g_i = \g_j \h$ for some $j \in \{1,\ldots,r\}$, $\h \in H$. Then for each $\x \in L_1$, $\g(\g_i(\x)) = \g_j(\h(\x))$, and the value of this is known because the $H$-action on $L_1$ is known. It is easy to see that we get one and only one consistent $G$-torsor action in this way.
\end{proof}

Because all field factors of a torsor are isomorphic, we will sometimes speak of ``the'' field factor of a torsor.

\subsubsection{Torsors over \'etale algebras}
On occasion, we will speak of a $G$-torsor over $L$, where $L$ is itself a product $K_1 \cross \cdots \cross K_r$ of fields. By this we simply mean a product $T_1 \cross \cdots \cross T_r$ where each $T_i$ is a $G$-torsor over $K_i$. This case is without conceptual difficulty, and some theorems on torsors will be found to extend readily to it, such as the following variant of the fundamental theorem of Galois theory:
\begin{thm} \label{thm:tor_Gal}
Let $T$ be a $G$-torsor over an \'etale algebra $L$. For each subgroup $H \subseteq G$,
\begin{enumerate}[(a)]
  \item The fixed algebra $T^H$ is uniformly of degree $[G : H]$ over $L$ (that is, of this same degree over each field factor of $L$);
  \item $T$ is an $H$-torsor over $T^H$, under the same action;
  \item If $H$ is normal, then $T^H$ is also a $G/H$-torsor over $L$, under the natural action.
\end{enumerate}
\end{thm}
\begin{proof}
Adapt the relevant results from Galois theory.
\end{proof}


\subsection{A fresh look at Galois cohomology}
Galois cohomology is one of the basic tools in the development of class field theory. It is usually presented in a highly abstract fashion, but certain Galois cohomology groups, specifically $H^1(K, M)$ for finite $M$, have explicit meaning in terms of field extensions of $M$. It seems that this interpretation is well known but has not yet been written down fully, a gap that we fill in here. We begin by describing Galois modules.

\begin{prop}[\textbf{a description of Galois modules}] \label{prop:Gal_mod}
  Let $M$ be a finite abelian group, and let $K$ be a field. Let $M^-$ denote the subset of elements of $M$ of maximal order $m$, the exponent of $M$. The following objects are in bijection:
  \begin{enumerate}[$($a$)$]
      \item \label{it:Gal_mod} Galois module structures on $M$ over $K$, that is, continuous homomorphisms $\phi : G_K \to \Aut M$;
      \item \label{it:Gal_mod_T} $(\Aut M)$-torsors $T/K$;
      \item \label{it:Gal_mod_L0} $(\Aut M)$-extensions $L_0/K$, where $\Aut M \hookrightarrow \Sym M$ in the natural way;
      \item \label{it:Gal_mod_L-} $(\Aut M)$-extensions $L^-/K$, where $\Aut M \hookrightarrow \Sym M^-$ in the natural way.
    \end{enumerate}
\end{prop}
\begin{proof}
For item \ref{it:Gal_mod_L-} to make sense, we need that $M^-$ generates $M$; this follows easily from the classification of finite abelian groups.

The bijections are immediate from Propositions \ref{prop:G_sets_G_strucs} and \ref{prop:torsor}.
\end{proof}

We will denote $M$ with its Galois-module structure coming from these bijections by $M_{\phi}$, $M_{T}$, or $M_{L_0}$. Note that $T$, $L_0$, and $L^-$ are mutual resolvents.

\begin{examp} \label{ex:cubic_Gal_mod}
For example (and we will return to this case frequently), if we let $M = \sf C_3$ be the smallest group with nontrivial automorphism group: $\Aut M \isom \sf C_2$. Then the Galois module structures on $M$ are in natural bijection with $\sf C_2$-torsors over $K$, that is, quadratic \'etale extensions $T/K$. If $\ch K \neq 2$, these can be parametrized by Kummer theory as $T = K[\sqrt{D}]$, $D \in K^\cross/\( K^\cross\) ^2$. The value $D = 1$ corresponds to the split algebra $T = K \cross K$ and to the module $M$ with trivial action. We have an isomorphism
\[
  M_T \cong \{0, \sqrt{D}, -\sqrt{D}\}
\]
of $G_K$-sets, and of Galois modules if the right-hand side is given the appropriate group structure with $0$ as identity.

In particular, the Galois-module structures on $\sf C_3$ form a group $\Hom(G_K, \sf C_2) \isom K^\cross/\(K^\cross\)^2$: the group operation can also be viewed as \emph{tensor product} of one-dimensional $\FF_3$-vector spaces with Galois action.
\end{examp}

\subsubsection{Galois cohomology}

Note that the zeroth cohomology group $H^0(K, M)$ has a ready parametrization:
\begin{prop}\label{prop:H0}
Let $M = M_{L_0}$ be a Galois module. The elements of $H^0(K, M)$ are in bijection with the degree-$1$ field factors of $L_0$.
\end{prop}
\begin{proof}
Proposition \ref{prop:Gal_mod} establishes an isomorphism of $G_K$-sets between the coordinates of $L_0$ and the points of $M$. A degree-$1$ field factor corresponds to an orbit of $G_K$ on $\Coord(L_0)$ of size $1$, which corresponds exactly to a fixed point of $G_K$ on $M$.
\end{proof}

Deeper and more useful is a description of $H^1$. For an abelian group $M$, let $\GA(M) = M \rtimes \Aut M$ be the semidirect product under the natural action of $\Aut M$ on $M$. We can describe $\GA(M)$ more explicitly as the group of \emph{affine-linear transformations} of $M$; that is, maps
\[
  \sfa_{\g,\t}(\x) = \g \x + \t, \quad \g \in \Aut M, \t \in M
\]
composed of an automorphism and a translation, the group operation being composition. In particular, we have an embedding
\[
  \GA(M) \hookrightarrow \Sym(M).
\]
\begin{prop}[\textbf{a description of $H^1$}]\label{prop:H1}
Let $M = M_{\phi} = M_{L_0}$ be a Galois module.
\begin{enumerate}[$($a$)$]
\item \label{it:z1_hom} $Z^1(K,M)$ is in natural bijection with the set of continuous homomorphisms $\psi : G_K \to \GA(M)$ such that the following triangle commutes:
\begin{equation} \label{eq:tri_h1}
\xymatrix{
  G_K \ar[r]^\psi \ar[dr]_\phi & \GA(M) \ar[d] \\
  & \Aut M
}
\end{equation}
\item \label{it:h1_hom} $H^1(K,M)$ is in natural bijection with the set of such $\psi : G_K \to \GA(M)$ up to conjugation by $M \subseteq \GA(M)$.
\item \label{it:h1_ext} $H^1(K,M)$ is also in natural bijection with the set of $\GA(M)$-extensions $L/K$ (with respect to the embedding $\GA(M) \hookrightarrow \Sym(M)$) equipped with an isomorphism from their resolvent $(\Aut M)$-torsor to $T$.
\end{enumerate}
\end{prop}
\begin{proof}
By the standard construction of group cohomology, $Z^1$ is the group of continuous crossed homomorphisms
\[
Z^1(K, M) = \{\sigma : G_K \to M \mid \sigma(\gamma \delta) = \sigma(\gamma) + \phi(\gamma) \sigma(\delta)\}.
\]
Send each $\sigma$ to the map
\begin{align*}
  \psi : G_K &\to \GA(M) \\
  \gamma &\mapsto a_{\phi(\gamma), \sigma(\gamma)}.
\end{align*}
It is easy to see that the conditions for $\psi$ to be a homomorphism are exactly those for $\sigma$ to be a crossed homomorphism, establishing \ref{it:z1_hom}. For \ref{it:h1_hom}, we observe that adding a coboundary $\sigma_a(\gamma) = \gamma(a) - a$ to a crossed homomorphism $\sigma$ is equivalent to post-conjugating the associated map $\psi : G_K \to \GA(M)$ by $a$. As to \ref{it:h1_ext}, a $\GA(M)$-extension carries the same information as a map $\psi$ up to conjugation by \emph{the whole of $\GA(M)$}. Specifying the isomorphism from the resolvent $(\Aut M)$-torsor to $T$ means that the map $\pi \circ \psi = \phi : G_K \to \Aut(M)$ is known exactly, not just up to conjugation. Hence $\psi$ is known up to conjugation by $M$.
\end{proof}
\begin{rem}
The zero cohomology class corresponds to the extension $L_0$, with its structure given by the embedding $\Aut M \hookrightarrow \GA(M)$. This can be seen to be the unique cohomology class whose corresponding $\GA(X)$-extension has a field factor of degree $1$.
\end{rem}

If $K$ is a local field, a cohomology class $\alpha \in H^1(K,M)$ is called \emph{unramified} if it is represented by a cocycle $\alpha : \Gal(\bar K/K) \to M$ that factors through the unramified Galois group $\Gal(K^{\ur}/K)$. The subgroup of unramified coclasses is denoted by $H^1_\ur(K, M)$. If $M$ itself is unramified (and we will never have to think about unramified cohomology in any other case), this is equivalent to the associated \'etale algebra $L$ being unramified.

If $X = M_{T}$ is a Galois module and $\sigma \in Z^1(K,M)$ is the Galois module corresponding to a $\GA(X)$-extension $L/K$, we can also take the $ \GA(X) $-closure of $L$, a $ \GA(X) $-torsor $E$ which fits into the following diagram:
\begin{equation}
\begin{gathered}
\xymatrix@dr{
    E \ar@{-}[r]\ar@{-}[d] & T \ar@{-}[d] \\
    L \ar@{-}[r] & K }
\end{gathered}
\end{equation}
Because of the semidirect product structure of $ \GA(X) $, we have $ E \isom L \tensor_K T $. It is also worth tabulating the permutation representations of finite groups that yield each of the \'etale algebras discussed here:
\begin{equation}
\xymatrix{
  0 \ar[r] & M \ar[r] & \GA(M) \ar[r] \ar@{^{(}->}[ld]_{\text{yields $L$}} \ar@{^{(}->}[d]^{\text{yields $E$}} & \Aut(M) \ar[r] \ar@{^{(}->}[d]^{\text{yields $T$}} & 0 \\
  & \Sym(M) & \Sym(\GA(M)) & \Sym(\Aut(M))
}
\end{equation}

\subsubsection{The Tate dual}

If $M$ is a Galois module and the exponent $m$ of $M$ is not divisible by $\ch K$, then
\[
  M' = \Hom(M, \mu_m)
\]
is also a Galois module, called the \emph{Tate dual} of $M$. The modules $M$ and $M'$ have the same order and are isomorphic as abstract groups, though not canonically; as Galois modules, they are frequently not isomorphic at all.

\begin{examp}\label{ex:cubic_Tate_dual}
If $M = M_{K[\sqrt{D}]}$ is one of the order-$3$ modules studied in Example \ref{ex:cubic_Gal_mod}, then the relevant $\mu_m$ is
\[
  \mu_3 \isom M_{K[\sqrt{-3}]}.
\]
Examining the Galois actions (here it helps to use the theory of $G_K$-sets of size $2$ presented in Knus and Tignol \cite{QuarticExercises}), we see that
\[
  M' = M_{K[\sqrt{-3D}]}.
\]
This explains the $D \mapsto -3D$ pattern in the Scholz reflection theorem and its generalizations, including cubic Ohno-Nakagawa.
\end{examp}

\begin{examp}
A module $M$ of underlying group $\C_2 \cross \C_2$ is always self-dual, regardless of what Galois-module structure is placed on it. This can be proved by noting that $M$ has a unique alternating bilinear form
\begin{align*}
  B : M \cross M &\to \mu_2 \\
  (x,y) &\mapsto \begin{cases}
    1, & x = 0, y = 0, \text{ or } x = y \\
    -1, & \text{otherwise.}
  \end{cases}
\end{align*}
Being unique, it is Galois-stable and induces an isomorphism $M' \isom M$.
\end{examp}

Particularly notable for us are the cases when $\GA(M)$ is the full symmetric group $\Sym(M)$, for then \emph{every} \'etale algebra $L/K$ of degree $\size{M}$ has a (unique) $\GA(M)$-affine structure. It is easy to see that there are only four such cases:
\begin{itemize}
  \item $M = \{1\}$, $\GA(M) \cong \S_1$
  \item $M = \ZZ/2\ZZ$, $\GA(M) \cong \S_2$
  \item $M = \ZZ/3\ZZ$, $\GA(M) \cong \C_3 \rtimes \C_2 \cong \S_3$
  \item $M = \ZZ/2\ZZ \cross \ZZ/2\ZZ$, $\GA(M) \cong (\C_2 \cross \C_2) \rtimes \S_3 \cong \S_4$.
\end{itemize}
For degree exceeding $4$, not every \'etale algebra arises from Galois cohomology, a restriction that plays out in the existing literature on reflection theorems. For instance, Cohen, Rubinstein-Salzedo, and Thorne \cite{CohON} prove a reflection theorem in which one side counts $\sf D_p$-dihedral fields of prime degree $p \geq 3$. From our perspective, these correspond to cohomology classes of an $M = \sf C_p$ whose Galois action is by $\pm 1$. The Tate dual of such an $M$ can have Galois action by the full $(\ZZ/p\ZZ)^\cross$, and indeed they count extensions of Galois group $\GA(\sf C_p)$ on the other side of the reflection theorem. This will appear inevitable in light of the motivations elucidated in Part \ref{part:composed}.


\section{Extensions of Kummer theory to explicitize Galois cohomology}
\label{sec:Kummer}
Now that Galois cohomology groups $ H^1(K,M) $ have been parametrized by \'etale algebras, can invoke parametrizations of \'etale algebras by even more explicit objects. The most familiar instance of this is \emph{Kummer theory,} an isomorphism
\[
  H^1(K, \mu_m) \cong K^\cross/(K^\cross)^m
\]
coming from the long exact sequence associated to the \emph{Kummer sequence}
\[
  0 \longto \mu_m \longto \bar K^\cross \longto^{\bullet^m} \bar K^\cross \longto 0. 
\]
In favorable cases, the cohomology $H^1(K, M)$ of other Galois modules $M$ can be embedded into $R^\cross/(R^\cross)^m$ for some finite extension $R$ of $K$.

We first state the hypothesis we need:
\begin{defn}
Let $M$ be a finite Galois module of exponent $m$ over a field $K$, and let $X$ be a Galois-stable generating set of $M$. We say that $M$ equipped with $X$ is a \emph{good} module if the natural map of Galois modules
\begin{align*}
  \X = \bigoplus_{x \in X} (\ZZ/m\ZZ) &\to M \\
  e_x &\mapsto x
\end{align*}
is split, that is, its kernel admits a Galois-stable complementary direct summand $\tilde M$. Such a direct summand is known as a \emph{good structure} on $M$.
\end{defn}
\begin{prop}\label{prop:good}
The following examples of a Galois module $M$ with generating set $X$ are good:
\begin{enumerate}[$($a$)$]
  \item $M \isom \sf C_p$, with any action, and $X = M \setminus \{0\}$.
  \item $M \isom \sf C_m^{n}$, with any action preserving a basis $X$.
  \item $M \isom \sf C_m^{n - 1}$, $n \geq 2$ with $\gcd(m,n)$, with an action that preserves a \emph{hyperbasis} $X$, that is, a generating set of $n$ elements with sum $0$.
\end{enumerate}
\end{prop}
\begin{proof}
\begin{enumerate}[$($a$)$]
  \item Here the Galois modules are representations of $\FF_p^\cross \isom \sf C_{p-1}$ over $\FF_p$. Since the group and field are of coprime order, complete reducibility holds: any subrepresentation is a direct summand. In fact, $\X$ is the regular representation, $M$ is the tautological representation in which each $\lambda \in \FF_p^\cross$ acts by multiplication by $\lambda$, and $\tilde M$ can be taken (uniquely in general) to be the product of all the other isotypical components of $\X$.
  \item Here the natural map $\X \to M$ is an isomorphism, so $\tilde M = \X$.
  \item Here the natural map $\X \to M$ is the quotient by the one-dimensional space
  \[
    \<\sum_{x \in X} e_x\>.
  \]
  This space has a Galois-stable direct complement, namely the kernel $\tilde M$ of the linear functional
  \begin{align*}
    \X &\to \FF_2 \\
    e_x &\to 1. \qedhere
  \end{align*}
\end{enumerate}
\end{proof}
\begin{prop} \label{prop:Kummer_good}
Let $M$ be a Galois module with a good structure $(X, \tilde M)$, and let $R$ be the resolvent algebra corresponding to the $G_K$-set $X$. For any Galois module $A$ with underlying group $\ZZ/m\ZZ$, there is a natural injection
\[
  H^1(K, M \tensor A) \to H^1(R, A)
\]
as a direct summand. The cokernel is naturally isomorphic to
\[
  H^1(K, (\X/\tilde M) \tensor A).
\]
\end{prop}
\begin{proof}
We use the good structure
\[
  M \isom \tilde M \hookrightarrow \X
\]
to embed
\[
  H^1(K, \tilde M \tensor A) \hookrightarrow H^1(\X \tensor A).
\]
Since $\tilde M$ is a direct summand, this is an injection with cokernel naturally isomorphic to $H^1(K, (\X/\tilde M) \tensor A)$. It remains to construct an isomorphism
\[
  H^1(\X \tensor A) \stackrel{\sim}{\longrightarrow} H^1(R, A).
\]
If $R$ decomposes as a product
\[
  R \isom R_1 \cross \cdots \cross R_s  
\]
of field factors corresponding to the orbits $X = \bigsqcup_i X_i$ of $G_K$ on $X$, then $\X$ has a corresponding decomposition
\[
  \X = \bigoplus_{i = 1}^s \X_i
\]
where $\X_i = \<e_x : x \in X_i\>$ is none other than the induced module $\Ind_K^{R_i} \ZZ/m\ZZ$. Its cohomology is computed by Shapiro's lemma:
\[
  H^1(K, \X \tensor A) \isom \bigoplus_{i=1}^s H^1(K, \X_i \tensor A) = \bigoplus_{i=1}^s H^1(K, \Ind_K^{R_i} A) \isom \bigoplus_{i=1}^s H^1(R_i, A) = H^1(R, A).
\]
This is the desired isomorphism.
\end{proof}



We can harness Kummer theory to parametrize cohomology of other modules as follows. 
\begin{thm}[\textbf{an extension of Kummer theory}]\label{thm:Kummer_new} Let $ M $ be a finite Galois module, and assume that $ m = \exp M $ is not divisible by $ \ch K $. Let $ G_K $ act on the set $M'^-$ of surjective characters $ \chi : M \toto \mu_m $ through its actions on $ M $ and $ \mu_m $, and let $ F $ be the \'etale algebra corresponding to this $ G_K $-set.
\begin{enumerate}[$($a$)$]
  \item\label{it:Kum_exists} There is a natural group homomorphism
  \[
    \Kum : H^1(K, M) \to F^\cross / (F^\cross)^m.
  \]
  \item\label{it:Kum_p-ic} If $ M \cong \C_p $ is cyclic of prime order, then $\Kum$ is injective, $ F $ is naturally a $ (\ZZ/p\ZZ)^\cross $-torsor, and
  \[
    \im(\Kum) = \left\{ \alpha \in F^\cross / (F^\cross)^p : \tau_{c}(\alpha) = \alpha^c \, \forall c \in (\ZZ/p\ZZ)^\cross \right\}.
  \]
  If $ p = 3 $, then the image simplifies to
  \[
    \im(\Kum) = \left\{ \alpha \in F^\cross / (F^\cross)^3 : N_{F/K}(\alpha) = 1\right\},
  \]
  and the $\GA(\C_3) \isom S_3 $-extension $L$ corresponding to a given $\alpha \in F^\cross/ (F^\cross)^3$ of norm $ 1 $ can be described as follows: Define a $K$-linear map
  \begin{align*}
    \kappa : K &\to \bar K^3 \\
    \xi &\mapsto \(\tr_{\bar K^2/K} \xi \omega \sqrt[3]{\delta}\)_\omega,
  \end{align*}
  where $\sqrt[3]{\delta} \in \bar K^2$ is chosen to have norm $1$, and $\omega$ ranges through the set \[
    \{(1;1), (\zeta_3; \zeta_3^2); ( \zeta_3^2; \zeta_3)\}
  \] of cube roots of $1$ in $\bar K^2$ of norm $1$. Then
  \[
    L = K + \kappa(F).
  \]
  \item\label{it:Kum_quartic} If $ M \cong \C_2 \cross \C_2 $, then $ \Kum $ is injective and
  \[
    \im(\Kum) = \left\{ \alpha \in F^\cross / (F^\cross)^2 : N_{F/K}(\alpha) = 1\right\}.
  \]
  Moreover, the $ \GA(M) \isom S_4$-extension corresponding to a given $\alpha \in F^\cross/ (F^\cross)^2$ of norm $ 1 $ can be described as follows: Define a $K$-linear map
  \begin{align*}
    \kappa : K &\to \bar K^4 \\
    \xi &\mapsto \(\tr_{\bar K^3/K} \xi \omega \sqrt{\delta}\)_\omega,
  \end{align*}
  where $\sqrt{\delta} \in \bar K^3$ is chosen to have norm $1$, and $\omega$ ranges through the set \[
  \{(1;1;1), (1;-1;-1); (-1;1;-1); (-1;-1;1)\}
  \] of square roots of $1$ in $\bar K^3$ of norm $1$. Then
  \[
  L = K + \kappa(F).
  \]
\end{enumerate}
\end{thm}
\begin{proof}
If $ \chi : M \toto \mu_m $ is a surjective character, let $ F_\chi $ be the fixed field of the stabilizer of $ \chi $; thus $ F_\chi $ is the field factor of $ F $ corresponding to the $ G_K $-orbit of $ \chi $. If $ \chi_1,\ldots, \chi_\ell $ are orbit representatives, we can map
\[
  H^1(K, M) \mathop{\longrightarrow}\limits^{\prod\Res} \prod_{i} H^1(F_{\chi_i}, M) \mathop{\longrightarrow}\limits^{\prod\chi_{i*}} \prod_i H^1(F_{\chi_i}, \mu_m) \isom \prod_i F_{\chi_i}^\cross/(F_{\chi_i}^\cross)^m \isom F^\cross / (F^\cross)^m.
\]
This yields our map $ \Kum $. Alternatively, note that by Shapiro's lemma,
\[
  \prod_i H^1(F_{\chi_i}, \mu_m) \isom \prod_i H^1(K, \Ind_{F_{\chi_i}}^K \mu_m) \isom H^1(K, I),
\]
where
\[
  I_M = \Ind_{F}^K \mu_m = \bigoplus_{\chi : M \toto \mu_m} \mu_m,
\]
a Galois module under the action
\[
  g\big((a_\chi)_i\big) = \big(g(c_{g^{-1}(\chi)})_\chi\big) = \big(g(c_{\chi(g\bullet)})_\chi\big).
\]
Under this identification, it is not hard to check that $ \Kum = j_* $, where $ j $ is the inclusion $ M \hookrightarrow I_M $ given by
\[
  a \mapsto (\chi(a))_\chi.
\]

Although $j$ is injective (because the characters of maximal order $m$ generate the group of all characters), it is not obvious whether $ j $ induces an injection on cohomology, nor what the image is. What makes the modules $ M $ in parts \ref{it:Kum_p-ic} and \ref{it:Kum_quartic} tractable is that, in these cases, $M \bs \{0\}$ is a \emph{good} generating set for $M$, so $ M $ is a direct summand of $ I_M $. In part \ref{it:Kum_p-ic}, we can identify
\[
  I_M \isom \Ind_{\{1\}}^{(\ZZ/p\ZZ)^\cross} \FF_p \tensor_{\FF_p} \mu_m
\]
as a twist of the regular representation of $ (\ZZ/p\ZZ)^\cross $ over $ \FF_p $. Since $ \FF_p $ has a complete set of $ (p-1) $st roots of unity, this representation splits completely into one-dimensional subrepresentations. The image of $ j $ is the eigenspace generated by $ (c)_{c \in \ZZ/p\ZZ^\cross} $, so $ \Kum $ is injective and its image is the subspace of $ F^\cross / (F^\cross)^p $ cut out by the same relations $ \tau_c(x) = c x $ (where $\tau_c$ is the torsor operation on $F$, resp.~the automorphism of $I_M$, indexed by $c$) that cut out $ j(M) $ in $ I_M $.

As to part \ref{it:Kum_quartic}, since $\ZZ/2\ZZ \cross \ZZ/2\ZZ$ has three surjective characters whose product is $1$, we have $ I_M/M \cong \mu_2 $ with the map $ \eta : I_M \to \mu_2 $ given by multiplying the coordinates. Since $ \mu_2 $ also injects diagonally into $ I_M $, we easily get a direct sum decomposition, which shows that $ \Kum $ is injective. As to the image, it is not hard to show that the diagram
\[
  \xymatrix{
    F^\cross/(F^\cross)^2 \ar[r]^\sim \ar[d]^{N} & H^1(F,\C_2) \ar[r]^\sim \ar[d]^\Cor & H^1(K, I) \ar[dl]^{\eta_*} \\
    K^\cross/(K^\cross)^2 \ar[r]^\sim & H^1(K, C_2)
  }
\]
commutes, establishing the desired norm characterization of $ \im(\Kum) $.

The formulas by radicals for the cubic and quartic algebras corresponding to a Kummer element follow easily by chasing through the Galois actions on the appropriate \'etale algebras. The quartic case is also considered by Knus and Tignol, where a closely related description of $L$ is given (\cite{QuarticExercises}, Proposition 5.13).
\end{proof}

\begin{rem}
    For general $ M $, the map $ H^1(K,M) \to F^\cross/(F^\cross)^m $ may be made by the construction in Theorem \ref{thm:Kummer_new}, but its image is hard to characterize, and it may not even be injective: for instance, when $M \cong \sf C_4$, coclasses correspond to $\sf D_4$-extensions, and $\Kum$ conflates each extension with its mirror extension (compare Example \ref{ex:D4}).
\end{rem}

Though it will not be used in the sequel, it is worth noting that Artin-Schreyer theory is amenable to the same treatment.

\begin{thm}\label{thm:A-S_new} Let Let $ M $ be a finite Galois module with underlying abelian group $ A $ of exponent $ m = p = \ch K $.
    \begin{enumerate}[$($a$)$]
        \item\label{it:A-S_exists} There is a natural map
        \[
        \AS : H^1(K, M) \to F / \wp(F).
        \]
        \item\label{it:A-S_p-ic} If $ A \cong \C_p $, then $\AS$ is injective, $ F $ is naturally a $ (\ZZ/p\ZZ)^\cross $-torsor, and
        \[
        \im(\AS) = \left\{ \alpha \in F / \wp(F) : \tau_{c}(\alpha) = c \alpha \; \forall c \in (\ZZ/p\ZZ)^\cross \right\}.
        \]
        \item\label{it:A-S_quartic} If $ p = 2 $ and $ A \cong \C_2 \cross \C_2 $, then $ \AS $ is injective and
        \[
        \im(\AS) = \left\{ \alpha \in F / \wp(F) : \tr_{F/K}(\alpha) = 0\right\}.
        \]
    \end{enumerate}
\end{thm}

\subsection{The Tate pairing and the Hilbert symbol}
\label{sec:Tate_pairing}
Assume now that $ K $ is a local field. Our next step will be to understand the (local) Tate pairing, which is given by a cup product
\[
  \langle \, , \, \rangle_T : H^1(K, M) \cross H^1(K, M') \to H^2(K, \mu_m) \isom \mu_m.
\]
As we were able to parametrize the cohomology groups $H^1(K, M)$ in favorable cases, it should not come as a surprise that we can often describe the Tate pairing with similar explicitness.

Recall the definitions of the Artin and Hilbert symbols. If $M \cong \ZZ/m\ZZ $ has \emph{trivial} $ G_K $-action, then $ M' \cong \mu_m $, and we have a Tate pairing
\[
  \langle \, , \, \rangle_T : H^1(K, \ZZ/m\ZZ) \cross H^1(K, \mu_m) \to \mu_m
\]
Now $H^1(K, \ZZ/m\ZZ) \cong \Hom(K, \ZZ/m\ZZ)$ parametrizes $\ZZ/m\ZZ$-torsors, while by Kummer theory, $H^1(K, \mu_m) \cong K^\cross / (K^\cross)^m$. The Tate pairing in this case is none other than the \emph{Artin symbol} (or \emph{norm-residue symbol}) $\phi_L(x)$ which attaches to a cyclic extension $L$, of degree dividing $m$, a mapping $\phi_L : K^\cross \to \Gal(L/K) \to \mu_m$ whose kernel is the norm group $N_{L/K}(L^\cross)$ (see Neukirch \cite{NeukirchCoho}, Prop.~7.2.13). If, in addition, $\mu_m \subseteq K$, then $H^1(K, \ZZ/m\ZZ)$ is also isomorphic to $K^\cross/\( K^\cross\) ^m$, and the Tate pairing is an alternating pairing
\[
  \langle \, , \, \rangle : K^\cross/\( K^\cross\) ^m \cross K^\cross/\( K^\cross\) ^m \to \mu_m
\]
classically called the \emph{Hilbert symbol} (or \emph{Hilbert pairing}). It is defined in terms of the Artin symbol by
\begin{equation}
  \<a, b\> = \phi_{K[\sqrt[m]{b}]}(a).
\end{equation}
In particular, $\<a, b\> = 1$ if and only if $a$ is the norm of an element of $K[\sqrt[m]{b}]$. This can also be described in terms of the splitting of an appropriate Severi-Brauer variety; for instance, if $m = 2$, we have $ \<a,b\> = 1 $ exactly when the conic
\[
  ax^2 + by^2 = z^2
\]
has a $K$-rational point. See also Serre (\cite{SerreLF}, \textsection\textsection XIV.1--2). (All identifications between pairings here are up to sign; the signs are not consistent in the literature and are totally irrelevant for this paper.) Pleasantly, for the types of $ M $ featured in Theorem \ref{thm:Kummer_new}, the Tate pairing can be expressed simply in terms of the Hilbert pairing.

We extend the Hilbert pairing to \'etale algebras in the obvious way: if $L = K_1 \cross \cdots \cross K_s$, then
\[
  \<(a_1;\ldots;a_s), (b_1;\ldots;b_s)\>_L := \<a_1, b_1\>_{K_1} \cdot \cdots \cdot \<a_s, b_s\>_{K_s}.
\]
Note that if $a$ is a norm from $L[\sqrt[m]{b}]$ to $L$, then $\<a, b\>_L = 1$, but the converse no longer holds. We then have the following:

\begin{thm}[\textbf{a formula for the local Tate pairing}]\label{thm:Tate_pairing} Let $ K $ be a local field. For $ M $, $ F $ as in Theorem \ref{thm:Kummer_new}, let $ M' $ be the Tate dual of $ M $, and let $ F' $ be the corresponding \'etale algebra, corresponding to the $G_K$-set $M^-$ of elements of maximal order in $M$, just as $F$ corresponds to $M'^-$. The Tate pairing
    \[
      \<\bullet,\bullet\> : H^1(K, M) \cross H^1(K, M') \to H^2(K, \mu_m) \isom \C_m
    \]
    can be described in terms of the Hilbert pairing in the following cases:
\begin{enumerate}[(a)]
\item\label{it:Tate_p} If $A \isom \C_p$, then both $F$ and $ F' $ embed naturally into $ F'' := F[\mu_p] $, and the Tate pairing is the restriction of the Hilbert pairing on $ F'' $.
\item\label{it:Tate_quartic} If $A \cong (\ZZ/2\ZZ)^2$, then we have natural isomorphisms $ M \isom M' $, $ F \isom F' $, and the Tate pairing is the restriction of the Hilbert pairing on $ F $.
\end{enumerate}
\end{thm}
\begin{proof}
  In case \ref{it:Tate_quartic}, set $F'' = F = F[\mu_2]$. We will do the two cases largely in parallel.
  
    Let $ \Surj(A,B) \subseteq \Hom(A,B) $ denote the set of surjections between two groups $ A,B $. Note that if $ A,B $ are Galois modules, then $ \Surj(A,B) $ is a $G_K$-set. Note that $ F'' $ is the \'etale algebra corresponding to the $ G_K $-set
    \[
      Z = \Surj(M, \mu_p) \cross \Surj(\mu_p, \ZZ/m\ZZ).
    \]
    There is an obvious map $ Z \to \Surj(M, \mu_p) $ given by projection to the first factor, which allows us to recover the identification $ F'' = F[\mu_m] $. There is also a map of $ G_K $-sets
    \[
    \Psi : Z \to \Surj(M', \mu_p)
    \]
    which sends a pair $ (\chi,u) $ (where $ \chi: M \toto \mu_m $, $ u : \mu_m $) to the unique surjective $ \psi : M' \to \mu_p $ satisfying
    \[
      \begin{cases}
        \psi(\chi) = u^{-1}(1), & A \isom \sf C_p \\
        \psi(\chi) = 1, & A \isom \sf C_2 \cross \sf C_2.
      \end{cases}
    \]
    This allows us to embed $ F' $ into $ F'' $. It is worth noting that when $ A \isom \sf C_2 \cross \sf C_2 $, $ u $ carries no information and $ F \isom F' \isom F'' $.
    
    Let $ F_1'',\ldots, F_\ell'' $ be the field factors of $ F'' $; each $ F_i'' $ corresponds to an orbit $ G_K(\chi_i, u_i) $ on $ Z $. Let $ \psi_i = \Psi(\chi_i, u_i). $ Then for $ \sigma \in H^1(K, M) $, $ \tau \in H^1(K, M') $,
    \begin{align*}
      \<\sigma, \tau\>_{\text{Hilb}} &= \<\Kum(\sigma), \Kum(\tau)\>_{\text{Hilb; } F''} \\
      &= \prod_{F''_i} \inv_{F''_i} \big(\chi_{i*} \res^K_{F_i''}\sigma \cup \psi_{i*} \res^{K}_{F_i''}\tau\big).
    \end{align*}
    Since $ \inv_{F''_i} = \inv_K \circ \Cor^{K}_{F''_i}$ (a standard fact), we have
    \begin{align*}
      \<\sigma, \tau\>_{\text{Hilb}} &= \inv_K \sum_{F_i''} \Cor^{K}_{F_i''} \big(\chi_{i*} \res^K_{F_i''}\sigma \cup \psi_{i*} \res^{K}_{F_i''}\tau\big) \\
      &= \inv_K \sum_{F_i''} \Cor^{K}_{F_i''} \big(\ev \circ (\chi_i \tensor \psi_i )\big)_* \res^K_{F_i''}(\sigma \cup \tau)
    \end{align*}
    where $ \ev : M \tensor M' \to \mu_m $ is the evaluation map. We now apply the following lemma, which slightly generalizes results seen in the literature.
    \begin{lem}\label{lem:cor_res}
        Let $H \subseteq G$ be a subgroup of finite index. Let $X$ and $\Y$ be $G$-modules, and let $f : X \to \Y$ be a map that is $H$-linear (but not necessarily $G$-linear). Denote by $\tilde f$ the $G$-linear map
        \[
        \tilde{f}(\x) = \sum_{gH \in G/H} \g f \g^{-1}(\x).
        \]
        Let $\sigma \in H^n(H, \Y)$. Then
        \[
        \Cor_H^G(f_* \Res_H^G \sigma) = \tilde{f}_* \sigma.
        \]
    \end{lem}
    \begin{proof}
        Since we are concerned with the equality of a pair of $\delta$-functors, we can apply dimension shifting to assume that $n = 0$. The proof is now straightforward.
    \end{proof}
    Applying with $ G = G_K $, $ H = G_{F''_i} $, and $ f = \ev \circ (\chi_i \tensor \psi_i ) : M \tensor M' \to \mu_m $, we get
    \begin{align*}
      \<\sigma, \tau\>_{\text{Hilb}} &= \inv_K \sum_{F_i''} \Big(\sum_{g \in G_K/G_{F_i''}} g \circ \ev \circ (\chi_i \circ g^{-1} \tensor \psi_i \circ g^{-1}) \Big) \\
      &= \inv_K \sum_{F_i''} \Big(\sum_{g \in G_K/G_{F_i''}} \big(\ev \circ (\chi_{i,g} \tensor \psi_{i,g} )\big)\Big)_* (\sigma \cup \tau), \\
    \end{align*}
    where $ \chi_{i,g} = g(\chi_i) $ and $ \psi_{i,g} = g(\psi_i) $ are given by the natural action. Now the outer sum runs over all $ G_K $-orbits of $ Z $ while the inner sum runs over the elements of each orbit, so we simply get
    \[
      \<\sigma, \tau\>_{\text{Hilb}} = \inv_K \Big(\sum_{(\chi,u) \in Z} \ev \circ (\chi \tensor \Psi(\chi,u)) \Big)_* (\sigma \cup \tau).
    \]
    Since the Tate pairing is given by
    \[
      \<\sigma, \tau\>_{\text{Tate}} = \inv_K \ev_* (\sigma \cup \tau),
    \]
    it remains to check that
    \[
      \sum_{(\chi,u) \in Z} \ev \circ (\chi \tensor \Psi(\chi,u)) = \ev
    \]
    as maps from $ M \tensor M' $ to $ \mu_m $. In the case $ M \cong \sf C_p $, each term is actually equal to $ \ev $, and there are $ (p - 1)^2 \equiv 1 $ mod $ p $ terms. In the case $ M \cong \sf C_2 \cross \sf C_2 $, a direct verification on a basis of $ M \tensor M' $ is not difficult.
\end{proof}

\section{Rings over a Dedekind domain}
\label{sec:rings}
Thus far, we have been considering \'etale algebras $L$ over a field $K$. We now suppose that $K$ is the fraction field of a Dedekind domain $\OO_K$ (not of characteristic $2$), which for us will usually be a number field or a completion thereof, although there is no need to be so restrictive. Our topic of study will be the subrings of $L$ that are lattices of full rank over $\OO_K$---the \emph{orders,} to use the standard but unfortunately overloaded word.

There is always a unique \emph{maximal order} $\OO_L$, the integral closure of $\OO_K$ in $L$. If $L = L_1 \cross \cdots \cross L_r$ is a product of field factors, we have $\OO_L = \OO_{L_1} \cross \cdots \cross \OO_{L_r}$.

\subsection{Indices of lattices}
There is one piece of notation that we explain here to avoid confusion. If $V$ is an $n$-dimensional vector space over $K$ and $A,B \subseteq V$ are two full-rank lattices, we denote by the \emph{index} $[A:B]$ the unique fractional ideal $\cc$ such that
\[
\Lambda^n A = \cc \Lambda^n B
\]
as $\OO_K$-submodules of the top exterior power $\Lambda^n V$. Alternatively, if $A \supseteq B$, then the classification theorem for finitely generated modules over $\OO_K$ lets us write
\[
A/B \cong \OO_K/\cc_1 \oplus \cdots \oplus \OO_K/\cc_r,
\]
and the index equals
\[
\cc_1\cc_2 \cdots \cc_r.
\]
The index satisfies the following basic properties:
\begin{itemize}
  \item $[A:B][B:C] = [A:C]$;
  \item If $V$ is a vector space over both $K$ and a finite extension $L$, and $A$ and $B$ are two $\OO_L$-sublattices, then $[A:B]_K = N_{L/K}[A:B]_L$;
  \item If $V = L$ is a $K$-algebra and $\alpha \in L^\cross$, then $[A : \alpha A] = N_{L/K}(\alpha)$.
\end{itemize}
Despite the apparent abstractness of its definition, the index $[A:B]$ is \emph{not} hard to compute in particular cases: localizing at a prime ideal, we can assume $\OO_K$ is a PID, and then it is the determinant of the matrix expressing any basis of $B$ in terms of a basis of $A$.

If $L$ is a $K$-algebra, $\OO \subseteq L$ is an order, and $\aa \subseteq L$ is a fractional ideal, the index $[\OO : \aa]$ is called the \emph{norm} of $\aa$ and will be denoted by $N_\OO(\aa)$ or, when the context is clear, by $N(\aa)$. Note the following basic properties:
\begin{itemize}
  \item If $\aa = \alpha\OO$ is principal, then $N_\OO(\aa) = N_{L/K}(\alpha)$.
  \item If $\aa$ and $\bb$ are two $\OO$-ideals and $\aa$ is invertible, then $N(\aa\bb) = N(\aa)N(\bb)$. This is easily derived from the theorem that an invertible ideal is locally principal (Lemma \ref{lem:inv=pri}). It is false for two arbitrary $\OO$-ideals.
  \item If $\OO_K = \ZZ$ or $\ZZ_p$, then for any $L$ and $\OO$, the norm $N_\OO(\aa)$ of an integral ideal is the ideal generated by the absolute norm $\size{\OO/\aa}.$
\end{itemize}

\subsection{Discriminants}
As is standard, we define the \emph{discriminant ideal} of an order $\OO$ in an \'etale algebra $L$ to be the ideal $\dd$ generated by the trace pairing
\begin{equation}
  \begin{aligned}
    \tau : \OO^{2n} &\to \OO_K \\
    (\xi_1,\xi_2,\ldots,\xi_n,\eta_1,\eta_2,\ldots,\eta_n) &\mapsto \det[\tr \xi_i \eta_j]_{i,j = 1}^n.
  \end{aligned}
\end{equation}
The trace pairing is nondegenerate, that is, $\dd \neq 0$ (this is one equivalent definition of \'etale). The primes dividing $\dd$ are those at which $L$ is ramified and/or $\OO$ is nonmaximal. This notion is standard and widely used. However, it does not quite extend the (also standard) notion of the discriminant of a $ \ZZ $-algebra over $ \ZZ $, which has a distinction between positive and negative discriminants. The Ohno-Nakagawa theorem involves this distinction prominently; Dioses \cite{Dioses} and Cohen--Rubinstein-Salzedo--Thorne \cite{CohON} each frame their extensions of O-N in terms of an ad-hoc notion of discriminant that incorporates the splitting data of an order at the infinite primes. Here we explain the variant that we will use.

Since the trace pairing $\tau$ is alternating in the $\xi$'s and also in the $\eta$'s, it can be viewed as a bilinear form on the rank-$1$ lattice $\Lambda^n \OO$. Identifying $\Lambda^n \OO$ with a (fractional) ideal $\cc$ of $\OO_K$ (whose class is often called the \emph{Steinitz class} of $\OO$), we can write
\[
\tau(\xi) = D \xi^2
\]
for some nonzero $D \in \cc^{-2}$. Had we rescaled the identification $\Lambda^n \OO \to \cc$ by $\lambda \in K^\cross$, $D$ would be multiplied by $\lambda^2$. We call the pair $(\cc, D)$, up to the equivalence $(\cc, D) \sim (\lambda\cc, \lambda^{-2}D)$, the \emph{discriminant} of $\OO$ and denote it by $\Disc \OO$.

There is another perspective on the discriminant $\Disc \OO$. Let $\tilde{L}$ be the $S_n$-torsor corresponding to $L$, which comes with $n$ embeddings $\kappa_1,\ldots, \kappa_n : L \to \tilde{L}$ freely permuted by the $S_n$-action (not to be confused with the $n$ coordinates of $L$). Noting that, for any $\alpha \in L$,
\[
\tr(\alpha) = \sum_i \kappa_i(\alpha),
\]
we can factor the trace pairing matrix:
\[
[\tr \xi_i\eta_j]_{i,j} = [\kappa_h(\xi_i)]_{i,h} \cdot [\kappa_h(\eta_j)]_{h,j}.
\]
Define
\[
\tau_0(\xi_1,\ldots,\xi_n) = \det [\kappa_h(\xi_i)]_{i,h},
\]
so that
\[
\tau(\xi_1,\ldots,\xi_n, \eta_1,\ldots, \eta_n) = \tau_0(\xi_1,\ldots, \xi_n) \cdot \tau_0(\eta_1,\ldots,\eta_n).
\]
Now look more carefully at the map $\tau_0$. First, $\tau_0$ is alternating under permutations of the $\xi_i$'s, so it defines a linear map
\[
\tau_0 : \cc \to \tilde L.
\]
Moreover, $\tau_0$ is alternating under postcomposition by the torsor action of $S_n$ on $\tilde{L}$, which permutes the $\kappa_h$ freely. Thus the image of $\tau_0$ lies in the $C_2$-torsor $T_2 = \tilde{L}^{A_n}$, which we call the \emph{discriminant torsor} of $L$, and even more specifically in the $(-1)$-eigenspace of the nontrivial element of $C_2$. By (the simplest case of) Kummer theory, we may write $T_2 = K[\sqrt{D'}]$. If $(\xi_1,\ldots,\xi_n)$ is any $K$-basis of $\OO_L$, so that $\xi_1 \wedge \cdots \wedge \xi_n$ corresponds to some nonzero element $c \in \cc$, then
\[
Dc^2 = \tau(c,c) = \tau_0(c)^2 = \( a \sqrt{D'}\) ^2 = D'a^2.
\]
Thus $T_2 = K[\sqrt{D}]$. We summarize this result in a proposition.

\begin{prop}
  If $L$ is an \'etale algebra over $K$ of discriminant $(\cc, D)$, then $K[\sqrt{D}]$ is the discriminant torsor of $L$; that is, the diagram of Galois structure maps
  \[
  \xymatrix{
    G_K \ar[r]^{\phi_L} \ar[rd]_{\phi_{K[\sqrt{D}]}} & S_n \ar[d]^{\sgn} \\
    & S_2
  }
  \]
  commutes.
\end{prop}

There is notable integral structure on $D$ as well.
\begin{lem}[\textbf{Stickelberger's theorem over Dedekind domains}]\label{lem:Stickelberger}
  If $(\cc, D)$ is the discriminant of an order $\OO$, then $D \equiv t^2$ mod $4\cc^{-2}$ for some $t \in \cc^{-1}$.
\end{lem}
\begin{rem}
  When $\OO_K = \ZZ$, Lemma \ref{lem:Stickelberger} states that the discriminant of an order is congruent to $0$ or $1$ mod $4$: a nontrivial and classical theorem due to Stickelberger. Our proof is a generalization of the most familiar one for Stickelberger's theorem, due to Schur \cite{Schur1929}.
\end{rem}
\begin{proof}
  Since $D \in \cc^{-2}$, the conclusion can be checked locally at each prime dividing $2$ in $\OO_K$. We can thus assume that $\OO_K$ is a DVR and in particular that $\cc = (1)$. Now there is a simple tensor $\xi_1 \wedge \cdots \wedge \xi_n$ that corresponds to the element $1 \in \cc$. By definition,
  \begin{equation} \label{eq:x_Stickelberger}
    \sqrt{D} = \tau_0(\xi_1,\ldots,\xi_n) = \det [\kappa_h(\xi_i)]_{i,h} =
    \sum_{\pi \in S_n} \Big( \sgn(\sigma) \prod_i \kappa_{\pi(i)} (\xi_i) \Big)
    = \rho - \bar\rho,
  \end{equation}
  where
  \[
  \rho = \sum_{\pi \in A_n} \prod_i \kappa_{\pi(i)} (\xi_i)
  \]
  lies in $T_2$ by symmetry and $\bar\rho$ is its conjugate. By construction, $\rho$ is integral over $\OO_K$, that is to say, $\rho + \bar\rho$ and $\rho\bar\rho$ lie in $\OO_K$. Now
  \[
  D = (\rho - \bar\rho)^2 = (\rho + \bar\rho)^2 - 4 \rho\bar\rho
  \]
  is the sum of a square and a multiple of $4$ in $\OO_K$.
\end{proof}
\begin{rem}
  One can write \eqref{eq:x_Stickelberger} in the suggestive form
  \[
  \det [\kappa_h(\xi_i)]_{i,h} = \det \begin{bmatrix}
    \theta_1(1) & \theta_2(1) \\
    \theta_1(\rho) & \theta_2(\rho)
  \end{bmatrix}
  \]
  where $\theta_1, \theta_2$ are the two automorphisms of $T_2$. This equates discriminants of orders in $L$ with those of orders in $T_2$. Equalities of determinants of this sort reappear in Bhargava's parametrizations of quartic and quintic rings and appear to be a common feature of many types of resolvent fields.
\end{rem}

We can now state the notion of discriminant as we would like to use it.
\begin{defn}
  A \emph{discriminant} over $\OO_K$ is an equivalence class of pairs $(\cc,D)$, with $D \in \cc^{-2}$ and $D \equiv t^2$ mod $4\cc^2$ for some $t \in \cc^{-1}$, up to the equivalence relation
  \[
  (\cc, D) \sim (\lambda\cc, \lambda^{-2}D).
  \]
  If $\OO$ is an \'etale order, the \emph{discriminant} $\Disc \OO$ is defined as follows: Pick any representation $\phi : \Lambda^n \OO \to \cc$ of the Steinitz class as an ideal class; then $\Disc \OO$ is the unique pair $(\cc, D)$ such that
  \[
  \det[\tr \xi_i\eta_j]_{i,j = 1}^n = D \cdot \phi(\xi_1 \wedge \cdots \wedge \xi_n) \cdot 
  \phi(\eta_1 \wedge \cdots \wedge \eta_n).
  \]
\end{defn}

Note the following points.
\begin{itemize}
  \item The discriminant recovers the discriminant ideal via $\dd = D \cc^2$.
  \item If $L$ has degree $3$, the discriminant also contains the splitting information of $L$ at the infinite primes. Namely, for each real place $\iota$ of $K$, if $\iota(D) > 0$ then $L_v \cong \RR \cross \RR \cross \RR$, while if $\iota(D) < 0$ then $L_v \cong \RR\cross \CC$.
  \item By a usual abuse of language, if $L$ is an \'etale algebra over a number field $K$, its \emph{discriminant} is the discriminant of the ring of integers $\OO_L$ over $\OO_K$.
  \item The discriminants over $\OO_K$ form a cancellative semigroup under the multiplication law
  \[
  (\cc_1, D_1)(\cc_2, D_2) = (\cc_1\cc_2, D_1D_2).
  \]
  \item If $\OO_K$ is a PID, then we can take $\cc = \OO_K$, and then the discriminants are simply nonzero elements $D \in \OO_K$ congruent to a square mod $4$, up to multiplication by squares of units.
  \item We will often denote a discriminant by a single letter, such as $\D$. When elements or ideals of $\OO_K$ appear in discriminants, they are to be understood as follows:
  \begin{align}
    D \quad (D \in \OO_K) \quad &\text{means} \quad ((1), D) \label{eq:discs_conv_1} \\
    \cc^2 \quad (\cc \subseteq K) \quad &\text{means} \quad (\cc, 1). \label{eq:discs_conv_2}
  \end{align}
  The seemingly counterintuitive convention \eqref{eq:discs_conv_2} is motivated by the fact that, if $\cc = (c)$ is principal, then $(\cc, 1)$ is the same discriminant as $((1), c^2)$.
\end{itemize}

With these remarks in place, the reader should not have difficulty reading and proving the following relation:
\begin{prop}
  If $\OO \supseteq \OO'$ are two orders in an \'etale algebra $L$, then
  \[
  \Disc \OO' = [\OO : \OO']^2 \cdot \Disc \OO.
  \]
\end{prop}

\subsection{Quadratic rings}
We will spend a lot of time investigating the number of rings over $\OO_K$ of given degree $n$ and discriminant $\D = (\cc, D)$. For \emph{quadratic} rings, the problem has a complete answer:
\begin{prop}[\textbf{the parametrization of quadratic rings}]
  Let $\OO_K$ be a Dedekind domain of characteristic not $2$. For every discriminant $\D$, there is a unique quadratic \'etale order $\OO_\D$ having discriminant $\D$.
\end{prop}
\begin{proof}
  Note first that the theorem is true when $\OO_K = K$ is a field: by Kummer theory, quadratic \'etale algebras over $K$ are parametrized by $K^\cross/\( K^\cross\) ^2$, as are discriminants; and it is a simple matter to check that $\Disc K[\sqrt{D}] = D$. We proceed to the general case.
  
  For existence, let $\D = (\cc, D)$ be given. By definition, $D$ is congruent to a square $t^2$ mod $4\cc^2$, $t\in \cc^{-1}$. Consider the lattice
  \[
  \OO = \OO_K \oplus \cc \xi, \quad \xi = \frac{t + \sqrt{D}}{2} \in L = K[\sqrt{D}].
  \]
  To prove that $\OO$ is an order in $L$, it is enough to verify that $(c\xi)(d\xi) \in \OO$ for any $c,d \in \cc$, and this follows from the computation
  \[
  \xi^2 = \xi(t - \bar\xi) = t\xi - \( \frac{t^2 - D}{4}\) 
  \]
  and the conditions $t \in \cc^{-1}, t^2 - D \in 4\cc^{-2}$.
  
  Now suppose that $\OO_1$ and $\OO_2$ are two orders with the same discriminant $\D = (\cc, D)$. Their enclosing $K$-algebras $L_1$, $L_2$ have the same discriminant $D$ over $K$, and hence we can identify $L_1 = L_2 = L$. Now project each $\OO_i$ along $\pi : L \to L/K$ is an $\OO_K$-lattice $\cc_i$ in $L/K$, which is a one-dimensional $K$-vector space: indeed, we naturally have $L/K \cong \Lambda^2 L$, and upon computation, we find that $\Disc \OO_i = (\cc_i, D)$. Consequently $\cc_1 = \cc_2 = \cc$. Now, for each $\beta \in \cc$, the fiber $\pi^{-1}(\cc) \intsec \OO_i$ is of the form $\beta_i + \OO_K$ for some $\beta_i$. The element $\beta_1 - \beta_2$ is integral over $\OO_K$ and lies in $K$, hence in $\OO_K$. Thus $\OO_1 = \OO_2$.
\end{proof}

If $\OO$ is any order in an \'etale algebra $ L/K $ ($ \ch K \neq 2 $), the quadratic order $B = \OO_{\Disc \OO}$ having the same discriminant as $\OO$ is called the \emph{quadratic resolvent ring} of $\OO$. It embeds into the discriminant torsor $T_2$, in two conjugate ways. Indeed, it is not hard to show that $B$ is generated by the elements
\[
\rho(\xi_1,\ldots,\xi_n) = \sum_{\pi \in A_n} \prod_i \kappa_{\pi(i)} (\xi_i) \in T_2
\]
appearing in the proof of Lemma \ref{lem:Stickelberger}.

\begin{rem}
The notion of a quadratic resolvent ring extends to characteristic $2$, being always an order in the quadratic resolvent algebra constructed in Example \ref{ex:rsv_n2}. We omit the details.
\end{rem}

\subsection{Cubic rings}
Cubic and quartic rings have parametrizations, known as \emph{higher composition laws,} linking them to certain forms over $\OO_K$ and also to ideals in resolvent rings. The study of higher composition laws was inaugurated by Bhargava in his celebrated series of papers (\cite{B1,B2,B3,B4}), although the gist of the parametrization of cubic rings goes back to work of F.W.~Levi \cite{Levi}. Later work by Deligne and by Wood \cite{WQuartic, W2xnxn} has extended much of Bhargava's work from $\ZZ$ to an arbitrary base scheme. In a previous paper \cite{ORings}, the author explained how a representative sample of these higher composition laws extend to the case when the base ring $A$ is a Dedekind domain. In the present work, we will need a few more; fortunately, there are no added difficulties, and we will briefly run through the statements and the methods of proof.

\begin{thm}[\textbf{the parametrization of cubic rings}] \label{thm:hcl_cubic_ring} Let $A$ be a Dedekind domain with field of fractions $K$, $\ch A \neq 3$.
  \begin{enumerate}[$($a$)$]
    \item Cubic rings $\OO$ over $A$, up to isomorphism, are in bijection with cubic maps
    \[
    \Phi : M \to \Lambda^2 M
    \]
    between a two-dimensional $A$-lattice $M$ and its own Steinitz class, up to isomorphism, in the obvious sense of a commutative square
    \[
    \xymatrix{
      M_1 \ar[r]^\sim_{i} \ar[d]^{\Phi_1} & M_2 \ar[d]^{\Phi_2} \\
      \Lambda^2 M_1 \ar[r]^\sim_{\det i} & \Lambda^2 M_2.
    }
    \]
    The bijection sends a ring $\OO$ to the index form $\Phi : \OO/A \to \Lambda^2(\OO/A)$ given by
    \[
    \x \mapsto \x \wedge \x^2.
    \]
    \item \label{cubic:field} If $\OO$ is nondegenerate, that is, the corresponding cubic $K$-algebra $L = K \tensor_A \OO$ is \'etale, then the map $\Phi$ is the restriction, under the Minkowski embedding, of the index form of $\bar{K}^3$, which is
    \begin{equation}\label{eq:k3_phi}
      \begin{aligned}
        \Phi : \bar{K}^3 / \bar{K} &\to \bar{K}^2 / \bar{K} \\
        (x;y;z) &\mapsto \big((x - y)(y - z)(z - x), 0\big).
      \end{aligned}
    \end{equation}
    \item \label{cubic:lift} Conversely, let $L$ be a cubic \'etale algebra over $K$. If $\bar{\OO} \subseteq L/K$ is a lattice such that
    $\Phi$ sends $\bar{\OO}$ into $\Lambda^2 \bar{\OO}$, then there is a unique cubic ring $\OO \subseteq L$ such that, under the natural identifications, $\OO/A = \bar{\OO}$.
  \end{enumerate}
  \end{thm}
  \begin{proof}
    \begin{enumerate}[$($a$)$]
      \item 
    The proof is quite elementary, involving merely solving for the coefficients of the unknown multiplication table of $\OO$. The case where $A$ is a PID is due to Gross (\cite{cubquat}, Section 2): the cubic ring having index form
    \[
      f(x\xi + y\eta) = ax^3 + bx^2y + cxy^2 + dy^3
    \]
    has multiplication table
    \begin{equation}\label{eq:mult_table}
    \begin{aligned}
      \xi\eta = -ad,
      \xi^2 = -ac + b\xi - a\eta,
      \eta^2 = -bd + d\xi - c\eta.
    \end{aligned}
    \end{equation}
    For the general Dedekind case, see my \cite{ORings}, Theorem 7.1. It is also subsumed by Deligne's work over an arbitrary base scheme; see Wood \cite{WQuartic} and the references therein.
    \item This follows from the fact that the index form respects base change. The index form of $\bar K^3/\bar K$ is a Vandermonde determinant that can easily be written in the stated form.
    \item We have an integral cubic map $\Phi\big|_{\bar \OO} : \bar \OO \to \Lambda^2 \bar \OO$, which is the index form $\Phi_\OO$ of a unique cubic ring $\OO$ over $\OO_K$. But over $K$, $\Phi_\OO$ is isomorphic to the index form of $L$. Since $L$ (as a cubic ring over $K$) is determined by its index form, we obtain an identification $\OO \tensor_{\OO_K} K \isom L$ for which $\OO/A$, the projection of $\OO$ onto $L/K$, coincides with $\bar{\OO}$. The uniqueness of $\OO$ is obvious, as $\OO$ must lie in the integral closure $\OO_{L}$ of $K$ in $L$.
  \end{enumerate}
  \end{proof}
    In this paper we only deal with nondegenerate rings, that is, those of nonzero discriminant, or equivalently, those that lie in an \'etale $K$-algebra. Consequently, all index forms $\Phi$ that we will see are restrictions of \eqref{eq:k3_phi}. When cubic algebras are parametrized Kummer-theoretically, the resolvent map becomes very explicit and simple:
    \begin{prop}[\textbf{explicit Kummer theory for cubic algebras}]\label{prop:Kummer_resolvent_cubic}
      Let $R$ be a quadratic \'etale algebra over $K$ ($\ch K \neq 3$), and let
      \[
      L = K + \kappa(R)
      \]
      be the cubic algebra of resolvent $R' = R \odot K[\mu_3]$ (the Tate dual of $R$) corresponding to an element $\delta \in K^\cross$ of norm $1$ in Theorem \ref{thm:Kummer_new}\ref{it:Kum_p-ic}, where
      \[
      \kappa(\xi) = \(\tr_{\bar K^2/K} \xi \omega \sqrt[3]{\delta}\)_{\omega \in \(\bar{K}^2\)^{N=1}[3]} \in \bar K^3
      \]
      so $\kappa$ maps $R$ bijectively onto the traceless plane in $L$. Then the index form of $L$ is given explicitly by
      \begin{equation}\label{eq:Kummer_rsv_cubic}
        \begin{aligned}
          \Phi : L/K &\to \Lambda^2(L/K) \\
          \kappa(\xi)& \mapsto 3\sqrt{-3}\delta\xi^3 \wedge 1,
        \end{aligned}
      \end{equation}
      where we identify
      \[
        \Lambda^2 L/K \isom \Lambda^3 L \isom \Lambda^2 R' \isom R'/K \isom \sqrt{-3} \cdot R/K
      \]
      using the fact that $R'$ is the discriminant resolvent of $L$. 
    \end{prop}
    \begin{proof}
      Direct calculation, after reducing to the case $K = \bar{K}$.
    \end{proof}

\begin{thm}[\textbf{self-balanced ideals in the cubic case}]
  \label{thm:hcl_cubic_sbi}
  Let $\OO_K$ be a Dedekind domain, $\ch K \neq 3$, and let $R$ be a quadratic \'etale extension. A \emph{self-balanced triple} in $R$ is a triple $(B, I, \delta)$ consisting of a quadratic order $B \subseteq R$, a fractional ideal $I$ of $B$, and a scalar $\delta \in (KB)^\cross$ satisfying the conditions
  \begin{equation}
    \delta I^3 \subseteq B, \quad N(I) = (t) \text{ is principal}, \textand N(\delta) t^3 = 1,
  \end{equation}
  \begin{enumerate}[$($a$)$]
    \item \label{cubic:idl} 
    Fix $B$ and $\delta \in R^\cross$ with $N(\delta)$ a cube $t^{-3}$. Then the mapping
    \begin{equation}\label{eq:sbi_cubic}
      I \mapsto \OO = \OO_K + \kappa(I)
    \end{equation}
    defines a bijection between
    \begin{itemize}
      \item self-balanced triples of the form $(B, I, \delta)$, and
      \item subrings $\OO \subseteq L$ of the cubic algebra $L = K + \kappa(R)$ corresponding to the Kummer element $\delta$, such that $\OO$ is \emph{$3$-traced}, that is, $\tr(\xi) \in 3\OO_K$ for every $\xi \in L$.
    \end{itemize} 
  \item Under this bijection, we have the discriminant relation
    \begin{equation} \label{eq:cubic zmat disc}
      \disc C = -27 \disc B.
    \end{equation}
  \end{enumerate}
\end{thm}
\begin{proof}
The mapping $\kappa$ defines a bijection between lattices $I \subseteq R$ and $\kappa(I) \subseteq L/K$. The difficult part is showing that $I$ fits into a self-balanced triple $(B, I, \delta)$ if and only if $\kappa(I)$ is the projection of a $3$-traced order $\OO$. Note that if $(B, I, \delta)$ exists, it is unique, as the requirement $[B : I] = (t)$ pins down $B$. 

Rather than establish this equivalence directly, we will show that both conditions are equivalent to the symmetric trilinear form
\begin{align*}
  \beta : I \cross I \cross I &\to \Lambda^2 R \\
  (\alpha_1, \alpha_2, \alpha_3) &\to \delta \alpha_1 \alpha_2 \alpha_3
\end{align*}
taking values in $t^{-1} \cdot \Lambda^2 I$.

In the case of self-balanced ideals, this was done over $\ZZ$ by Bhargava \cite[Theorem 3]{B1}. Over a Dedekind domain, it follows from the parametrization of \emph{balanced} triples of ideals over $B$ \cite[Theorem 5.3]{ORings}, after specializing to the case that all three ideals are identified with one ideal $I$. It also follows from the corresponding results over an arbitrary base in Wood \cite[Theorem 1.4]{W2xnxn}.

In the case of rings, we compute by Proposition \ref{prop:Kummer_resolvent_cubic} that $\beta$ is the trilinear form attached to the index form of $\kappa(I)$. By Theorem \ref{thm:hcl_cubic_ring}\ref{cubic:lift}, the diagonal restriction $\beta(\alpha, \alpha, \alpha)$ takes values in $t^{-1} \Lambda^2 (I)$ if and only if $\kappa(I)$ lifts to a ring $\OO$. We wish to prove that $\beta$ itself takes values in $t^{-1} \Lambda^2(I)$ if and only if $\OO$ is $3$-traced. Note that both conditions are local at the primes dividing $2$ and $3$, so we may assume that $\OO_K$ is a DVR. With respect to a basis $(\xi,\eta)$ of $I$ and a generator of $t^{-1} \Lambda^2(I)$, the index form of $\OO$ has the form
\[
  f(x,y) = ax^3 + bx^2y + cxy^2 + dy^3, \quad a,\ldots, d\in \OO_K.
\]
If this is the diagonal restriction of $\beta$, then $\beta$ itself can be represented as a $3$-dimensional matrix
\[
  \bbq {a}{b/3}{b/3}{c/3}{b/3}{c/3}{c/3}{d},
\]
which is integral exactly when $b, c \in 3\OO_K$. Since the trace ideal of $\OO$ is generated by
\[
  \tr(1) = 3, \quad \tr(\xi) = -b, \quad \tr(\eta) = c
\]
(by reference to the multiplication table \eqref{eq:mult_table}), this is also the condition for $\OO$ to be $3$-traced, establishing the equivalence.

The discriminant relation \eqref{eq:cubic zmat disc} follows easily from the definition of $\kappa$.
\end{proof}

\subsection{Quartic rings and their cubic resolvent rings}

The basic method for parametrizing quartic orders is by means of \emph{cubic resolvent rings,} introduced by Bhargava in \cite{B3} and developed by Wood in \cite{WQuartic} and the author in \cite{ORings}.

\begin{defn}[\cite{ORings}, Definition 8.1; also a special case of \cite{WQuartic}, p.~1069]
  Let $A$ be a Dedekind domain, and let $\OO$ be a quartic algebra over $A$. A \emph{resolvent} for $\OO$ (``numerical resolvent'' in \cite{ORings}) consists of a rank-$2$ $A$-lattice $Y$, an $A$-module isomorphism $\Theta : \Lambda^2 Y \to \Lambda^3 (\OO/A)$, and a quadratic map $\Phi : \OO/A \to Y$ such that there is an identity of biquadratic maps
  \begin{equation}\label{eq:resolvent}
    x \wedge y \wedge xy = \Theta(\Phi(x) \wedge \Phi(y))
  \end{equation}
  from $\OO \cross \OO$ to $\Lambda^3 (\OO/A)$.
\end{defn}

We collect some basic facts about these resolvents.

\begin{thm}[\textbf{the parametrization of quartic rings}]\label{thm:hcl_quartic} The notion of resolvent for quartic rings has the following properties. 
  \begin{enumerate}[$($a$)$]
    \item \label{quartic:to_ring} If $X$ is a rank-$3$ $A$-lattice and $\Theta : \Lambda^2 Y \to \Lambda^3 X$, $\Phi : X \to Y$ satisfy \eqref{eq:resolvent}, then there is a unique (up to isomorphism) quartic ring $\OO$ equipped with an identification $\OO/A \isom X$ making $(\OO, Y, \Theta, \Phi)$ a resolvent.
    \item \label{quartic:cubic_ring} There is a canonical (in particular, base-change-respecting) way to associate to a resolvent $(\OO, Y, \Theta, \Phi)$ a cubic ring $C$ and an identification $C/A \isom X$ with the following property: For any element $x \in \OO$ and any lift $y \in C$ of the element $\Phi(x) \in C/A$, we have the equality
    \[
    x \wedge x^2 \wedge x^3 = \Theta(y \wedge y^2).
    \]
    It satisfies
    \[
    \Disc C = \Disc \OO.
    \]
    (Here the discriminants are to be seen as quadratic resolvent rings, as in \cite{ORings}; this implies the corresponding identity of discriminant ideals.)
    If $\OO$ is nondegenerate, then $C$ is unique.
    \item \label{quartic:exist} Any quartic ring $\OO$ has at least one resolvent.
    \item \label{quartic:max} If $\OO$ is maximal, the resolvent is unique (but need not be maximal).
    \item \label{quartic:count} The number of resolvents of $\OO$ is the sum of the absolute norms of the divisors of the \emph{content} of $\OO$, the smallest ideal $\cc$ such that $\OO = \OO_K + \cc \OO'$ for some order $\OO'$.
    \item \label{quartic:field} Let $(Y, \Theta, \Phi)$ be a resolvent of $\OO$ with associated cubic ring $C$, and let $K = \Frac A$. If the corresponding quartic $K$-algebra $L = K \tensor_A \OO$ is \'etale, then the cubic $K$-algebra $R = C \tensor_A \OO$ is none other than the cubic resolvent of $L$, as defined in Example \ref{ex:rsv43}. The maps $\Theta$ and $\Phi$ are the restrictions, under the Minkowski embedding, of the unique resolvent of $\bar{K}^4$, which is $\bar{K}^3$ with the maps
    \begin{equation}\label{eq:k4_theta}
      \begin{aligned}
        \Theta : \Lambda^2(\bar{K}^3 / \bar{K}) &\to \Lambda^3(\bar{K}^4 / \bar{K}) \\
        (0;1;0) \wedge (0;0;1) &\mapsto (0;1;0;0) \wedge (0;0;1;0) \wedge (0;0;0;1)
      \end{aligned}
    \end{equation}
    and
    \begin{equation}\label{eq:k4_phi}
      \begin{aligned}
        \Phi : \bar{K}^4 / \bar{K} &\to \bar{K}^3 / \bar{K} \\
        (x;y;z;w) &\mapsto (xy + zw; xz + yw; xw + yz).
      \end{aligned}
    \end{equation}
    \item \label{quartic:lift}Conversely, let $L$ be a quartic \'etale algebra over $K$ and $R$ its cubic resolvent. Let
    \[
    \Theta_K : \Lambda^2 R/K \to \Lambda^3 L/K, \quad \Phi_K : L/K \to R/K
    \]
    be the resolvent data of $L$ as a (maximal) quartic ring over $K$. Suppose $\bar{\OO} \subseteq L/K$, $\bar{C} \subseteq R/K$ are lattices such that
    \begin{itemize}
      \item $\Phi_K$ sends $\bar{\OO}$ into $\bar{C}$,
      \item $\Theta_K$ maps $\Lambda^3 \bar{\OO}$ isomorphically onto $\Lambda^2 \bar{C}$.
    \end{itemize}
    Then there are unique quartic rings $\OO \subseteq L$, $C \subseteq R$ such that, under the natural identifications, $\OO/A = \bar{\OO}$, $C/A = \bar{C}$, and $\bar{C}$ is a resolvent with the restrictions of $\Theta_K$ and $\Phi_K$.
  \end{enumerate}
  
\end{thm}
\begin{proof}
  \begin{enumerate}[$($a$)$]
    \item See \cite{ORings}, Theorem 8.3.
    \item See \cite{ORings}, Theorems 8.7 and 8.8.
    \item See \cite{ORings}, Corollary 8.6.
    \item This is a special case of the following part.
    \item See \cite{ORings}, Corollary 8.5.
    \item By base-changing to $K$, we see that $Y \tensor_A K = R/K$ is a resolvent for $L$. Since the resolvent is unique, it suffices to show that the cubic resolvent $R'$ from Example \ref{ex:rsv43} is a resolvent for $L$ also. The maps $\Theta$ and $\Phi$ defined in the theorem statement are seen, by symmetry, to restrict to maps of the appropriate $K$-modules. The verification of \eqref{eq:resolvent} and of the fact that the multiplicative structure on $R'$ is the right one can be checked at the level of $\bar{K}$-algebras.
    \item Letting $X = \bar{\OO}$, $Y = \bar{C}$ in part \ref{quartic:to_ring}, we construct the desired $\OO$ and $C$. By comparison to the situation under base-change to $K$, we see that $\OO$, $C$ naturally inject into $L$, $R$ respectively. Uniqueness is obvious, as $\OO$ must lie in the integral closure $\OO_L$.
  \end{enumerate}
\end{proof}

In this paper we only deal with nondegenerate rings, that is, those of nonzero discriminant, or equivalently, those that lie in an \'etale $K$-algebra. Consequently, all resolvent maps $\Theta$, $\Phi$ that we will see are restrictions of \eqref{eq:k4_theta} and \eqref{eq:k4_phi}. When quartic algebras are parametrized Kummer-theoretically, the resolvent map becomes very explicit and simple:
\begin{prop}[\textbf{explicit Kummer theory for quartic algebras}]\label{prop:Kummer_resolvent_quartic}
  Let $R$ be a cubic \'etale algebra over $K$ ($\ch K \neq 2$), and let
  \[
  L = K + \kappa(R)
  \]
  be the quartic algebra of resolvent $R$ corresponding to an element $\delta \in K^\cross$ of norm $1$ in Theorem \ref{thm:Kummer_new}\ref{it:Kum_quartic}, where
  \[
  \kappa(\xi) = \(\tr_{\bar K^3/K} \xi \omega \sqrt{\delta}\)_{\omega \in \(\bar{K}^3\)^{N=1}[2]} \in \bar K^4
  \]
  so $\kappa$ maps $R$ bijectively onto the traceless hyperplane in $L$. Then the resolvent of $L$ is given explicitly by
  \begin{equation}\label{eq:Kummer_rsv_theta}
    \begin{aligned}
      \Theta : \Lambda^3(R) &\to \Lambda^3(L/K) \\
      \alpha \wedge \beta \wedge \gamma &\mapsto \frac{1}{16 \sqrt{N(\delta)}} \cdot \kappa(\alpha) \wedge \kappa(\beta) \wedge \kappa(\gamma)
    \end{aligned}
  \end{equation}
  \begin{equation}\label{eq:Kummer_rsv_phi}
    \begin{aligned}
      \Phi : L/K &\to R/K \\
      \kappa(\xi)& \mapsto 4\delta\xi^2
    \end{aligned}
  \end{equation}
\end{prop}
\begin{proof}
  Since the resolvent is unique (over a field, any \'etale extension has content $1$), it suffices to prove that \eqref{eq:Kummer_rsv_theta} and \eqref{eq:Kummer_rsv_phi} define a resolvent. This can be done after extension to $\bar{K}$, and then it is enough to prove that \eqref{eq:Kummer_rsv_theta} and \eqref{eq:Kummer_rsv_phi} agree with the standard resolvent on $\bar{K}^4$, given in Theorem \ref{thm:hcl_quartic}\ref{quartic:field}.
  
  As to \eqref{eq:Kummer_rsv_theta}, since both sides are alternating in $\alpha$, $\beta$, and $\gamma$, it suffices to prove it in the case that
  \[
  \alpha = (1;0;0), \quad \beta = (0;1;0) \quad \gamma = (0;0;1)
  \]
  form the standard basis of $R = \bar{K}^3$. Let $\delta = (\delta^{(1)}, \delta^{(2)}, \delta^{(3)})$. Then
  \begin{align*}
    \kappa(\alpha) &= \(\sqrt{\delta^{(1)}}, \sqrt{\delta^{(1)}}, -\sqrt{\delta^{(1)}}, -\sqrt{\delta^{(1)}}\) \\
    \kappa(\beta)  &= \(\sqrt{\delta^{(2)}}, -\sqrt{\delta^{(2)}}, \sqrt{\delta^{(2)}}, -\sqrt{\delta^{(2)}}\) \\
    \kappa(\gamma) &= \(\sqrt{\delta^{(3)}}, -\sqrt{\delta^{(3)}}, -\sqrt{\delta^{(3)}}, \sqrt{\delta^{(3)}}\)
  \end{align*}
  and hence the wedge product of these differs from the standard generator of $\Lambda^3(L/K)$ by a factor of
  \begin{align*}
    & \begin{vmatrix}
      1 & 1 & 1 & 1 \\
      \sqrt{\delta^{(1)}} & \sqrt{\delta^{(1)}} & -\sqrt{\delta^{(1)}} & -\sqrt{\delta^{(1)}} \\
      \sqrt{\delta^{(2)}} & -\sqrt{\delta^{(2)}} & \sqrt{\delta^{(2)}} & -\sqrt{\delta^{(2)}} \\
      \sqrt{\delta^{(3)}} & -\sqrt{\delta^{(3)}} & -\sqrt{\delta^{(3)}} & \sqrt{\delta^{(3)}}
    \end{vmatrix} \\
    &= \sqrt{N(\delta)} \begin{vmatrix}
      1 & 1 & 1 & 1 \\
      1 & 1 & -1 & -1 \\
      1 & -1 & 1 & -1 \\
      1 & -1 & -1 & 1
    \end{vmatrix} \\
    &= 16 \cdot \sqrt{N(\delta)}.
  \end{align*}
  The calculation for \eqref{eq:Kummer_rsv_phi} is even more routine.
\end{proof}

\begin{rem}
  The datum $\Theta$ of a resolvent carries no information, in the following sense. It is unique up to scaling by $c \in A^\cross$, and the resolvent data $(X, Y, \Theta, \Phi)$ and $(X, Y, c\Theta, \Phi)$ are isomorphic under multiplication by $c^{-1}$ on $X$ and by $c^{-2}$ on $Y$. If $A$ is a PID, indeed, neither $X$ nor $Y$ carries any information, and the entire data of the resolvent is encapsulated in $\Phi$, a pair of $3\cross 3$ symmetric matrices over $A$ (with formal factors of $1/2$ off the diagonal) defined up to the natural action of $\GL_3 A \cross \GL_2 A$. This establishes the close kinship with Bhargava's parametrization of quartic rings in \cite{B3}. However, it is useful to keep $\Theta$ around.
\end{rem}

\subsubsection{Traced resolvents}
Just as we found it natural to study not just binary cubic $1111$-forms, but also $1331$-forms and their analogue for each divisor of the ideal $(3)$, so too we study not just quartic rings in general but those satisfying a natural condition at the primes dividing $2$.

\begin{defn}
  Let $A$ be a Dedekind domain, $\ch A \neq 2$, and let $\tt$ be an ideal dividing $(2)$ in $A$. A resolvent $(\OO, Y, \Theta, \Phi)$ over $A$ is called \emph{$\tt$-traced} if, for all $x$ and $y$ in $A$, the associated bilinear form
  \[
  \Phi(x,y) = \frac{\Phi(x+y) - \Phi(x) - \Phi(y)}{2}
  \]
  whose diagonal restriction is $\Phi(x,x) = \Phi(x)$ takes values in $2^{-1}\tt Y$.
  If $A$ is a PID, this is equivalent to saying that the off-diagonal entries in the matrix representation of $\Phi$, which a priori live in $\frac{1}{2} A$, actually belong to $\frac{\tt}{2} A$. We say that $A$ is \emph{$\tt$-traced} if it admits a $\tt$-traced resolvent.
\end{defn}

Here are some facts about traced resolvents:
\begin{prop}\label{prop:traced} Let $\OO$ be a quartic ring over a Dedekind domain $\OO_K$.
  \begin{enumerate}[$(a)$]
    \item\label{traced:conds} $\OO$ is $\tt$-traced if and only if
    \begin{enumerate}[$(i)$]
      \item\label{traced:trace} $\tt^{2}|\tr x$ for all $x\in \OO$;
      \item\label{traced:sq} $x^2 \in A + \tt\OO$ for all $x\in \OO$.
    \end{enumerate}
    \item\label{traced:count} If $\OO$ is not an order in the \emph{trivial algebra} $K[\epsilon_1, \epsilon_2, \epsilon_3]/(\epsilon_i \epsilon_j)_{i,j=1}^3$, the number of $\tt$-traced resolvents of $\OO$ is the sum of the absolute norms of the divisors of its \emph{$\tt$-traced content,} which is the smallest ideal $\cc$ such that $\OO = A + \cc \OO'$ and $\OO'$ is also $\tt$-traced.
    \item\label{traced:reduced} If $(\OO,Y,\Theta,\Phi)$ is a $\tt$-traced resolvent with associated cubic ring $C$, then $\tt^{2}|\ct(C)$, that is, $C = A + \tt^{2}C'$ for some cubic ring $C'$. We call $C'$ a ``reduced resolvent'' of the $\tt$-traced ring $A$. Also, $\tt^8 | \disc A$.
  \end{enumerate}
\end{prop}
\begin{proof}
  \begin{enumerate}[$(a)$]
    \item
    Since both statements are local at the primes dividing $2$, we can assume that $\OO_K$ is a DVR, and thus that $\tt = (t)$ is principal. With respect to bases $(1 = \xi_0, \xi_1,\xi_2,\xi_3)$ for $\OO$ and $(1 = \eta_0, \eta_1, \eta_2)$ for a resolvent $C$, the structure constants $c_{ij}^k$ of the ring $\OO$, defined by
    \[
      \xi_i \xi_j = \sum_k c_{ij}^k \xi_k,
    \]
    are determined by the entries of the resolvent
    \[
      \Phi = \([a_{ij}], [b_{ij}]\)
    \]
    via the determinants
    \[
      \lambda^{ij}_{k\ell} = 2^{\1_{i\neq j} + \1_{k\neq \ell}} \begin{vmatrix}
        a_{ij} & a_{k\ell} \\
        b_{ij} & b_{k\ell}
      \end{vmatrix}
    \]
    and a set of formulas appearing in Bhargava \cite[equation (21)]{B3} and over a Dedekind domain by the author \cite[equation (12)]{ORings}:
    \begin{equation}\label{eq:c-lam}
      \begin{aligned}
        c_{ii}^j &= -\epsilon \lambda^{ii}_{ik} \\
        c_{ij}^k &= \epsilon \lambda^{jj}_{ii} \\
        c_{ij}^j - c_{ik}^k &= \epsilon \lambda^{jk}_{ii} \\
        c_{ii}^i - c_{ij}^j - c_{ik}^k &= \epsilon \lambda^{ij}_{ik},
      \end{aligned}
    \end{equation}
  where $(i,j,k)$ denotes any permutation of $(1,2,3)$ and $\epsilon = \pm 1$ its sign. (Here the nonappearance of some of the individual $c_{ij}^k$ on the left-hand side of \eqref{eq:c-lam} stems from the ambiguity of translating each $\xi_i$ by $\OO_K$, which does not change the matrix of $\Phi$.)
  
  Assume first that $\Phi : \OO/\OO_K \to C/\OO_K$ is $\tt$-traced. Then 
  \begin{equation}\label{eq:traced}
    \lambda^{ij}_{k\ell} \in \tt^{\1_{i\neq j} + \1_{k\neq \ell}}.
  \end{equation}
  We then prove that the conditions \ref{traced:trace} and \ref{traced:sq} must hold:
  \begin{enumerate}[$(i)$]
    \item The trace
    \begin{align*}
      \tr(\xi_1) &= c_{11}^1 + c_{12}^2 + c_{13}^3 \\
      &= \lambda^{12}_{13} + 2\lambda^{23}_{11} + 4 c_{13}^3 \\
      &\equiv 0 \mod \tt^2,
    \end{align*}
    and likewise $\tr(\xi_2), \tr(\xi_3) \in \tt^2$.
    \item The coefficients $c_{11}^i$ of ${\xi_1}^2$ satisfy:
    \begin{align*}
      c_{11}^2 &= \lambda_{13}^{11} \in \tt
    \end{align*}
    and likewise for $c_{11}^3$; and then $c_{11}^1 \in \tt$ also, since the trace $c_{11}^1 + c_{12}^2 + c_{13}^3 = \tr(\xi_1) \in \tt^2 \subseteq \tt$. So the desired relation $\xi^2 \in \OO_K + \tt \OO$ holds when $\xi = \xi_1$, indeed $\xi = a_{1}\xi_{1}$ for any $a_1 \in \OO_K$. The same proof works for $\xi = a_2\xi_2$ or $\xi = a_3\xi_3$. Since the case $\xi = a_0 \in \OO_K$ is trivial and squaring is a $\ZZ$-linear operation modulo $2$, we get the result for all $\xi \in \OO$.
  \end{enumerate}
  Conversely, suppose that \ref{traced:trace} and \ref{traced:sq} hold. We first establish \eqref{eq:traced}. We have
  \begin{itemize}
    \item $\lambda^{11}_{13} = c_{11}^2 \in \tt$
    \item $\lambda^{23}_{11} = c_{12}^2 - c_{13}^3 = \tr \xi_1 - c_{11}^1 - 2c_{13}^3 \in \tt$
    \item $\lambda^{12}_{13} = c_{11}^1 - c_{12}^2 - c_{13}^3 = \tr \xi_1 - 2\lambda^{23}_{11} + 4c_{13}^3 \in \tt^2$.
  \end{itemize}
  Permuting the indices as needed, this accounts for all the $\lambda^{ij}_{k\ell}$ about which \eqref{eq:traced} makes a nontrivial assertion.
  
  Now we work from the $\lambda^{ij}_{k\ell}$ back to the resolvent $(\A,\B)$. We may assume that $C$ is nontrivial (the trivial rings, one for each Steinitz class, are plainly $2$-traced with $(\A,\B) = (0,0)$.) Then, in the proof of \cite{ORings}, Theorem 8.4, the author established that there are vectors $\mu_{ij}$ in a two-dimensional vector space $V$ over $K$, unique up to $\GL_2(V)$, such that
  \[
    \mu_{ij} \wedge \mu_{k\ell} = \lambda^{ij}_{k\ell} \cdot \omega
  \]
  for some fixed generator $\omega \in \Lambda^2 V$. (The proof uses the \emph{Pl\"ucker relations}, which are a consequence of the associative law on $\OO$.) This $V$ is none other than $R/K$, the resolvent module of the quartic algebra $L = \OO \tensor_{\OO_K} K$, which admits the unique resolvent
  \begin{equation}\label{eq:rsv_field}
    \begin{aligned}
      \Phi(a_1\xi_1 + a_2\xi_2 + a_3\xi_3) &= \sum_{i<j} a_i a_j \mu_{ij} \\
      \Theta(\xi_1 \wedge \xi_2 \wedge \xi_3) &= \omega.
    \end{aligned}
  \end{equation}
  The resolvents of $\OO$ were found to be exactly the lattices $M$ containing the span $M_0$ of the six $\mu_{ij}$, with the correct index
  \[
    [M : M_0] = \cc = \(\lambda^{ij}_{k\ell}\)_{i,j,k,\ell} = \(c_{ii}^j, c_{ij}^k, c_{ij}^j - c_{ik}^k, c_{ii}^i - c_{ij}^j - c_{ik}^k : i \neq j \neq k \neq i\),
  \]
  the content ideal of $\OO$. By inspection of \eqref{eq:rsv_field} that $M$ is $\tt$-traced if and only if it actually contains the span $\tilde M_0$ of the six vectors
  \[
    \tilde\mu_{ij} = t^{-\1_{i \neq j}} \mu_{ij}.
  \]
  Condition \eqref{eq:traced} is interpreted as saying that the $\tilde\mu_{ij} \wedge \tilde\mu_{k\ell}$ are still integer multiples of $\omega$. Then the $\tt$-traced resolvents are the lattices $M \supseteq \tilde M_0$. The needed index
  \[
    \tilde \cc = [M : \tilde M_0] = \(\tilde\lambda^{ij}_{k\ell}\)_{i,j,k,\ell}
  \]
  is an integral ideal, so such $M$ exists, finishing the proof of \ref{traced:conds}.
  
  \item It suffices to prove that $\tilde \cc$ is the $\tt$-traced content of $\OO$. To see this, note that if $\OO = \OO_K + a \OO'$ has content divisible by $a$, then the structure coefficients $c_{ij}^k$ of $\OO'$ are obtained from those of $\OO$ by dividing by $a$. This means that the $\lambda^{ij}_{k\ell}$ and $\tilde\lambda^{ij}_{k\ell}$ are divided by $a$, and so remain integral (indicating that $\OO'$ is also $\tt$-traced) exactly when $a \mid \tilde \cc$.
  
  \item We can again reduce to the case that $\OO_K$ is a DVR so $\OO$ has an $\OO_K$-basis. Recall that the index form of the resolvent $C$ is given by
  \[
    f(x,y) = 4 \det (\A x + \B y)
  \]
  (\cite{B3}, Proposition 11; \cite{ORings}, Theorem 8.7).
  If $\A$ and $\B$ have off-diagonal entries in $2^{-1} \tt$, it immediately follows that $f$ is divisible by $t^2$, so $\t^2 \mid \ct(C)$. Consequently $\disc \OO = \disc C$, being quartic in the coefficients of $f$, is divisible by $t^8$. \qedhere
  \end{enumerate}
\end{proof}

Similar to Theorem \ref{thm:hcl_cubic_sbi}, we have the following relation between $2$-traced quartic rings and self-balanced ideals:
\begin{thm}[\textbf{self-balanced ideals in the quartic setting}] \label{thm:hcl_quartic_sbi}
  Let $\OO_K$ be a Dedekind domain, $\ch K \neq 2$, and let $R$ be a cubic \'etale extension. A \emph{self-balanced triple} in $R$ is a triple $(C, I, \delta)$ consisting of a cubic order $C \subseteq R$, a fractional ideal $I$ of $C$, and a scalar $\delta \in (K C)^\cross$ satisfying the conditions
  \begin{equation} \label{eq:quartic_sbi}
    \delta I^2 \subseteq C, \quad N(I) = (t) \text{ is principal}, \textand N(\delta) t^2 = 1,
  \end{equation}
    Fix an order $C\subseteq R$ and a scalar $\delta \in R^\cross$ with $N(\delta)$ a square $t^{-2}$. Then the mapping
    \begin{equation}\label{eq:sbi_quartic}
      I \mapsto \OO = \OO_K + \kappa(I)
    \end{equation}
    defines a bijection between
    \begin{itemize}
      \item self-balanced triples of the form $(C, I, \delta)$, and
      \item subrings $\OO \subseteq L$ of the quartic algebra $L = K + \kappa(R)$ corresponding to the Kummer element $\delta$, such that $\OO$ is $2$-traced with reduced resolvent $C$.
    \end{itemize} 
\end{thm}
\begin{proof}
The proof is very similar to that of \ref{thm:hcl_cubic_sbi}, so we simply summarize the main points. The linear isomorphism $\kappa$ establishes a bijection between lattices $I\subseteq R$ and $\kappa(I) \subseteq L/K$. We wish to prove that $(C, I,\delta)$ is balanced if and only if $\kappa(I)$ is the projection of a $2$-traced order with reduced resolvent $C$.

First note that either of these conditions uniquely specifies
\[
  [C : I] = (t),
\]
the former by the balancing condition $N(I) = (t)$, and the latter by the $\Theta$-condition that $\OO$ have discriminant $256 \disc C$.

Once again, it is difficult to proceed directly, and we instead prove that both conditions are equivalent to the bilinear map
\begin{align*}
  \Phi : I \cross I &\to R/K \\
  (\alpha_1,\alpha_2) &\to 4\delta \alpha_1\alpha_2
\end{align*}
taking values in $4C/K$.
\end{proof}

On the self-balanced ideals side, this follows from the parametrization of \emph{balanced} pairs of ideals by $2\times 3\times 3$ boxes performed over $\ZZ$ by Bhargava \cite[Theorem 2]{B2} and over a general base by Wood \cite[Theorem 1.4]{W2xnxn}.

On the quartic rings side, the diagonal restriction of $\Phi$ is precisely the resolvent of $\kappa(I)$, by Proposition \ref{prop:Kummer_resolvent_quartic}. That $\Phi(\alpha,\alpha) \in 4C/K$ for each $\alpha \in I$ expresses the one condition remaining for $\kappa(I)$ to lift (by Theorem \ref{thm:hcl_quartic}\ref{quartic:lift}) to a quartic ring $\OO$ with resolvent $\OO_K + 4C$. Then, by definition, this resolvent is $2$-traced exactly when $\Phi$ itself has image in $4C/K$.

\section{Cohomology of cyclic modules over a local field}
\label{sec:loc_fld_struc}

Let $ M $ be a Galois module with underlying group $\C_p$ over a local field $ K \supseteq \QQ_p $ (that is, a \emph{wild} local field of characteristic $ 0 $). Denote by $ T $ and $T'$, respectively, the $ (\ZZ/p\ZZ)^\cross $-torsors corresponding to the action of $ G_K $ on $M \bs \{0\}$ and on
\[
\Surj(M, \mu_p) = M' \setminus \{0\},
\]
and denote by $ \tau_c : T' \to T' $ the torsor operation corresponding to $c \in (\ZZ/p\ZZ)^\cross \isom \Aut M$. By Theorem \ref{thm:Kummer_new}, Kummer theory gives an isomorphism
\begin{equation}\label{eq:mul_Kum}
  H^1(K,M) \isom \left\{ \alpha \in T'^\cross / (T'^\cross)^p : \tau_{c}(\alpha) = \alpha^c \, \forall c \in (\ZZ/p\ZZ)^\cross \right\}.
\end{equation}
Our objective in this section is to understand the group on the right: that is, to describe a basis of it (a generalization of the well-known Shafarevich basis for $ T'^\cross $) and understand how the Tate pairing respects it. Much of our work parallels that of Del Corso and Dvornicich \cite{DCD} and Nguyen-Quang-Do \cite{Nguyen}.

If $ \sigma \in H^1(K,M) $, we let $ L = L_\sigma $ be the $ \GA(\C_p) $-extension of degree $ p $ coming from the affine action of $ G_K $ on $ M $, while we let $ E = E_\sigma $ be the associated $ \GA(\C_p) $-torsor. Owing to the semidirect product structure of $ \GA(\C_p) $, we get a natural decomposition
\[
E \isom L \tensor_K T.
\]
Using the division algorithm in $ \ZZ $, we let $ \ell = \ell(L) = \ell(\sigma)$ and $ \theta = \theta(L) = \theta(\sigma) $ the integers such that
\[
v_K(\disc L) = p(e - \ell) + \theta, \quad -1 \leq \theta \leq p - 2.
\]
We call $ \ell $ the \emph{level}, and $ \theta $ the \emph{offset}, of the $ \GA(\C_p) $-extension $ L $ or of the coclass $ \sigma $. Although these definitions appear strange, they allow us to state concisely the following theorem, which will be the main theorem of this section.
\begin{thm}[\textbf{levels and offsets}]\label{thm:levels}
  Let $ M $ be a Galois module with underlying group $\C_p$ over a local field $ K $ with $ \ch K \neq p $.
  \begin{enumerate}[$($a$)$]
    \item\label{lev:offset} The level $ \ell $ of a coclass determines its offset $ \theta $ uniquely in the following way:
    \begin{enumerate}[$ ($i$) $]
      \item If $ \ell = e $, then $ \theta = v_K(\disc T) $.
      \item If $ 0 \leq \ell < e $, then $ \theta \geq 0 $ and 
      \[
      \theta \equiv \ell - v_K(\beta) \mod p - 1,
      \]
      where $ \beta $ is the Kummer element corresponding to the resolvent $ \mu_{p-1} $-torsor of $ M $. 
      \item If $ \ell = -1 $, then $ \theta = -1 $.
    \end{enumerate}
    \item\label{lev:subgp} For all $ i \geq 0 $, the \emph{level space}
    \[
    \L_{i} = \L_{i}(M) = \{ \sigma \in H^1(K, M) : \ell(\sigma) \geq i \}
    \]
    consisting of coclasses of level at least $ i $ is a subgroup of $ H^1(K, M) $.
    \item\label{lev:ur} $ \L_e = H^1_{\ur}(K,M)$.
    \item\label{lev:size_Li} For $ 0 \leq i \leq e $,
    \[
    \size{\L_i} = q^{e - i} \size{H^0(K, M)}.
    \]
    \item\label{lev:size_all} $ \L_{-1} $ is the whole of $ H^1(K,M) $, and
    \[
    \size{H^1(K,M)} = q^e \size{H^0(K, M)} \cdot \size{H^0(K, M')}.
    \]
    \item\label{lev:distance} For $d \leq 1$, a neighborhood
    \[
      \left\{[\alpha] \in H^1(K,M) : \alpha \in T'^\cross, \size{\alpha - 1} \leq d \right\}
    \]
    is a level space $\L_i$ whose index $i$ is given by
    \[
      i = \begin{cases}
        \ds\Biggl\lceil\frac{v_K(\beta)}{p} + \frac{p-1}{p} \ceil{\frac{\log d}{\log \size{\pi_K}} - 1 - \frac{v_K(\beta)}{p-1}}\Biggr\rceil, & d \geq d_{\min} \\
        e + 1, & d < d_{\min},
      \end{cases}
    \]
    where
    \[
      d_{\min} = \size{p}^{p/(p-1)}.
    \]
    \item\label{lev:perp} For $ 0 \leq i \leq e $, with respect to the Tate pairing between $ H^1(K,M) $ and $ H^1(K, M') $,
    \[
    \L_i(M)^\perp = \L_{e-i}(M').
    \]
  \end{enumerate}
\end{thm}
One corollary is sufficiently important that we state it before starting the proof:
\begin{cor}\label{cor:levels}
  For $0 \leq i \leq e$, the characteristic function $L_i$ of the level space $\L_i$ has Fourier transform given by
  \begin{equation}\label{eq:hat 1 Li}
    \widehat{{L_i}} = q^{e - i} {L_{e-i}}.
  \end{equation}
  where $q = \size{k_K}$.
\end{cor}
\begin{proof}
  Immediate from Theorem \ref{thm:levels}, parts \ref{lev:size_Li} and \ref{lev:perp}.
\end{proof}

\subsection{Discriminants of Kummer and affine extensions}

The starting point for our investigation of discriminants is as follows:
\begin{thm} \label{thm:disc_Kummer}
  Let $ K $ be a local field with $\mu_p \subseteq K$, and let $ u \in K^\cross $ be a minimal representative of a class in  $K^\cross/(K^\cross)^p$. The discriminant ideal of the associated Kummer extension $L = K[\sqrt[p]{u}]$ is given by
  \[
  \Disc(L/K) = \begin{cases}
    \dfrac{p^{p} \cdot {\pi_K}^{p-1}}{(u-1)^{p-1}} \OO_K & v_K(u-1) < \dfrac{pe_K}{p-1} \\
    (1) & v_K(u-1) \geq \dfrac{pe_K}{p-1}.
  \end{cases}
  \]
\end{thm}
\begin{proof}
  One can find an explicit basis for $ \OO_L $ and compute the discriminant. For details, see Del Corso and Dvornicich \cite[Lemmas 5, 6, and 7]{DCD}.
\end{proof}

In this section, we will prove the following generalization:

\begin{thm}\label{thm:disc_Kummer_aff}
  Let $ K $ be a local field, and let $ M $ be a $ G_K $-module with underlying group $ \C_p $. Let $ T' $ be the $ (\ZZ/p\ZZ)^\cross $-torsor corresponding to the $ G_K $-set $ \Surj(M, \mu_m) $, and let $ T'_1 $ be the field factor of $ T' $. Let $u \in T'^\cross_1$ be a minimal representative for a class in $ (T'^\cross_1)_\omega $ parametrizing, via Theorem \ref{thm:Kummer_new} a coclass $ \sigma \in H^1(K, M) $, and let $ L $ be the corresponding $\GA(\C_p)$-extension. Then
  \begin{equation} \label{eq:dka_disc}
    \Disc(L/K)\OO_{T'_1} = \begin{cases}
      \dfrac{p^{p} \cdot {\pi_K}^{p-1}}{(u-1)^{p-1}} \OO_{T'_1} & v_{T'_1}(u-1) < \dfrac{pe_{T'_1}}{p-1} \\
      \Disc(T/K) \cdot \OO_{T'_1} & v_{T'_1}(u-1) \geq \dfrac{pe_{T'_1}}{p-1},
    \end{cases}
  \end{equation}
  where $ T $ is the $ (\ZZ/p\ZZ)^\cross $-torsor corresponding to $ M $.
\end{thm}
\begin{rem}
  Note that $ (u-1)^{p-1} \OO_{T'_1} $ is the extension of an ideal of $ K $, since $ e_{T'_1/K} $ divides $ p-1 $.    
\end{rem}

\begin{proof}
  If $ L $ is \emph{not} a field, then the image of $ G_K $ in $ \GA(\C_p) $ lies in a nontransitive subgroup (viewing $ \GA(\C_p) $ as embedded in $ \Sym(\C_p) $). It is not hard to show that every nontransitive subgroup of $ \GA(\C_p) $ has a fixed point. Moving this fixed point to $ 0 $, we get that $ \sigma = 0 $, $ u = 1 $, and $ L \isom K \cross T $, in accord with the second case of the formula.
  
  We may now assume that $ L $ is a field. Although the extension $ L/K $ need not be Galois, we have
  \[
  T[\mu_p] \tensor_K L \isom T[\mu_p, \sqrt[p]{u}],
  \]
  a Kummer extension of $ T[\mu_p] $. Let $ E $ be a field factor of $ T[\mu_p] $ containing $ T'_1 $. Then
  \[
  E \tensor_K L \isom E[\sqrt[p]{u}]
  \]
  as extensions of $ E $. Note that $ [E : K] | (p-1)^2 $; in particular, $ [E : K] $ is prime to $ p $. So $ u $ remains a minimal representative in $ E^\cross $, and since $ E $ and $ L $ must be linearly disjoint, $ E \tensor_K L = EL $ is a field unless $ u = 1 $. So
  \[
  \Disc(EL/E) = \begin{cases}
    \dfrac{p^{p} \cdot {\pi_E}^{p-1}}{(u-1)^{p-1}} \OO_E & v_K(u-1) < \dfrac{pe_K}{p-1} \\
    (1) & v_K(u-1) \geq \dfrac{pe_K}{p-1}.
  \end{cases}
  \]
  We must now relate $ \Disc(EL/E) $ to $ \Disc(L/K) $. If $ v_K(u-1) \geq {pe_K}/(p-1) $, then $ EL/E $ is unramified, so $ p \nmid e_{EL/K} $ and $ L/K $ is unramified as well. In particular, $ L/K $ is Galois, so $ M $ is trivial, $ T $ is totally split, and the formula again holds.
  
  We are left with the case that $ v_K(u-1) < \dfrac{pe_K}{p-1} $. Here $ EL/E $, and hence $ L/K $, are totally ramified. We relate their discriminants by the following trick, which also appears in Del Corso and Dvornicich \cite{DCD}. An $ \OO_K $-basis for $ \OO_L $ is given by
  \begin{equation} \label{eq:basis_L}
    1, \pi_L, \ldots, \pi_L^{p-1}.
  \end{equation}
  The same elements form an $ \OO_E $-basis for an order $ \OO \subseteq \OO_{EL} $, but their $ EL $-valuations are $0, e', \ldots, (p - 1)e'$, where $ e' = e_{E/K} $. Divide each basis element by $ \pi_E $ as many times as possible so that it remains integral. We get a new system of elements
  \begin{equation} \label{eq:basis_EL}
    1, \frac{\pi_L}{\pi_E^{a_1}}, \ldots, \frac{\pi_L^{p-1}}{\pi_E^{a_{p-1}}}.
  \end{equation}
  Since $ e' $ is coprime to $ p $, these elements have $ EL $-valuations $0,1, \ldots, e' - 1$ in some order and thus form an $ \OO_E $-basis for $ \OO_{EL} $. We have
  \begin{align*}
    v_{E}([\OO_{EL} : \OO]) &= a_1 + \cdots + a_{p-1}\\
    &= \frac{[e' + \cdots + (p-1)e'] - [1 + \cdots + (p-1)]}{p} \\
    &= \frac{(p-1)(e'-1)}{2},
  \end{align*}
  and hence
  \begin{align*}
    \Disc(L/K) &= \Disc(\OO/\OO_E) \\
    &= \Disc(\OO_{EL}/\OO_E) \cdot [\OO_{EL} : \OO]^2 \\
    &= \dfrac{p^{p} \cdot {\pi_E}^{p-1}}{(u-1)^{p-1}} \cdot \pi_{E}^{(p-1)(e'-1)} \cdot \OO_E \\
    &= \dfrac{p^{p} \cdot {\pi_E}^{e'(p-1)}}{(u-1)^{p-1}} \cdot \OO_E \\
    &= \dfrac{p^{p} \cdot {\pi_K}^{(p-1)}}{(u-1)^{p-1}} \cdot \OO_E. \qedhere
  \end{align*}
\end{proof}
\begin{rem}
  Along the lines of the preceding argument, we can prove the following more general result on discriminants in extensions of coprime degree:
\end{rem}
\begin{prop} \label{prop:disc_base_change}
  Let $ L $ and $ M $ be two extensions of a local field $ K $ with $\gcd([L:K],[M:K]) = 1$. Then
  \[
  \Disc(L/K) = \Disc(LM/M) \cdot \( \frac{\pi_M}{\pi_K} \)^{f_{L/K} \( e_{L/K} - 1 \)}.
  \]
\end{prop}

\subsection{The Shafarevich basis}

We start with the following exposition of the Shafarevich basis theorem. Although this theorem has appeared many times in the literature (see Del Corso and Dvornicich, \cite{DCD}, Proposition 6), we include a proof here by a method that will establish some important corollaries for us.

Filter $ U = K^\cross $ by the subgroups
\[
U_i = \{ x \in \OO_K^\cross : x \equiv 1 \mod \pi^i \},
\]
and let $ \bar{U}_i $ be the projection of $ U_i $ onto $ \bar{U} := K^\cross/(K^\cross)^p $. Note that $ \bar{U}_i = 0 $ for $ i > pe_K/(p-1) $, as the Taylor series for $ \sqrt[p]{x} $ about $ 1 $ converges for $ x \equiv 1 $ mod $ \pi^{\floor{pe_K/(p-1)}+1} $. So
\begin{equation} \label{eq:comp_series_K}
  \bar{U} \isom \bar{U}/\bar{U}_0 \oplus \bigoplus_{0 \leq i \leq \frac{pe_{T'_1}}{p-1}} \bar{U}_{i}/\bar{U}_{i+1}
\end{equation}
as $ \FF_p $-vector spaces, and we can produce a basis for $ \bar{U} $ by lifting a basis for each of the composition factors on the right-hand side.

\begin{prop}[\textbf{the Shafarevich basis theorem}] \label{prop:Sh_basis}
  Let $K \supseteq \QQ_p $ be a local field, and let $ i \geq 0 $. The structure of $ \bar{U}_{i}/\bar{U}_{i+1} $ is as follows:
  \begin{itemize}
    \item If 
    \[
    0 < i < \frac{p}{p-1} e_K, \quad p \nmid i,
    \]
    then $ \bar{U}_{i}/\bar{U}_{i+1} \isom U_i / U_{i+1} $ has a basis of $ f_K $ units of the form $1 + x_j \pi^i$, where $x_j$ ranges over an $\FF_p$-basis of $k_K$. We call these \emph{generic units.}
    \item if $ i = \frac{p}{p-1} e_K $ and $\mu_p \subseteq K$, then $ \bar{U}_{i}/\bar{U}_{i+1} $ has dimension $ 1 $ 
    and is generated by
    \[
    u = 1 + p(\zeta_p - 1) a,
    \]
    for any $ a \in \OO_K $ with $ p \nmid \tr_{\OO_K/\QQ_p}(a) $. We call such a generator an \emph{intimate unit,} and we let
    \[
      d_{\min} = \size{p(\zeta_p - 1)} = \size{p}^{p/(p-1)},
    \]
    the distance of an intimate unit to $1$.
    \item For all other $ i, $ we have $ \bar{U}_{i}/\bar{U}_{i+1} = 0 $.
  \end{itemize}
\end{prop}
\begin{proof}
  
  \ignore{
    In the tame case, it is clear that only $ i = 0 $ can yield a nonzero group $ \bar{U}_{i}/\bar{U}_{i+1} $, and in this case
    \[
    \bar{U}_0/\bar{U}_1 \isom k_K^\cross / (k_K^\cross)^p
    \]
    is generated by a generator of $k_K^\cross$ if $ p \mid \size{k_K^\cross} $, that is, $ \mu_p \subseteq K$, and is $ 0 $ otherwise.
    
    In the wild case,
  }
  Note that $ \bar{U}_0/\bar{U}_1 = 0$ because $ U_0/U_1 \isom k_K^\cross $ has order prime to $ p $. To compute $ \bar{U}_{i}/\bar{U}_{i+1} $, where $ i \geq 1 $, we must see how many of the congruence classes $ 1 + x \pi^i $ mod $ \pi^{i+1} $ (where $ x \in k_K $) contain a $ p $th power.
  
  Consider a general $ p $th power $ u^p $, $ 1 \neq u \in \OO_K^\cross $. Write $u = 1 + y \pi^j $, $ \pi \nmid y $. By the binomial theorem,
  \[
  u^p \equiv 1 + p y \pi^j + y^p \pi^{p j} \mod \pi^{2j + e_K},
  \]
  so
  \[
  v(u^p - 1) \begin{cases}
    = pj, & j < \frac{e_K}{p-1}, \\
    \geq \frac{p e_K}{p-1}, & j = \frac{e_K}{p-1}, \\
    = j + e_K, & j > \frac{e_K}{p-1}.
  \end{cases}
  \]
  We now perform the needed analysis in each case:
  \begin{itemize}
    \item If $ 0 < i < pe_K/(p-1) $ and $ p \nmid i $, then $ v(u^p - 1) $ can never attain the value $ i $, so the map $ U_i/U_{i+1} \to \bar{U}_i/\bar{U}_{i+1} $ is an isomorphism, and we must include an entire basis of $ f $ elements in our basis for $ \bar{U} $.
    \item If $ 0 < i < pe_K/(p-1) $ and $ p | i $, then the $ p $th powers
    \[
    (1 + y \pi^{i/p})^p \equiv 1 + y^p \pi^i \mod \pi^{i+1}
    \]
    cover all the desired congruence classes, as the map $ y \mapsto y^p $ on $ k_K $ is surjective; so $ \bar{U}_i/\bar{U}_{i+1} = 0 $.
    \item If $ i > pe_K/(p-1) $, then the $ p $th powers of elements of the form $ 1 + x\pi^{i - e_K} $ surject onto the congruence classes, repeating what we knew from the Taylor series.
    \item Finally, if $ i = pe_K/(p-1) $, then we can only use powers $ u^p = (1 + \pi^j)^p $ where $ j = e_K/(p-1) $. We have
    \[
    (1 + y \pi^j) \equiv 1 + p y \pi^j + y^p \pi^{p j} \equiv 1 + (y^p - c y) \pi^{i} \mod \pi^{i+1},
    \]
    where 
    \[
    c = \frac{-p}{\pi^{(p-1)j}} \in \OO_K^\cross.
    \]
    So we must analyze the (clearly linear) map $ \wp_c : k_K \to k_K $ given by $ \wp_c(y) = y^p - c y $.
    
    If $ c $ is \emph{not} a $ (p-1) $st power in $ k_K $ (or in $ \OO_K $, which amounts to the same thing by Hensel's lemma), then $ \wp_c $ is injective and hence surjective, so $ \bar{U}_i/\bar{U}_{i+1} = 0 $. Also, $ K $ has no nontrivial $ p $th roots of unity, as $ u = \zeta_p $ would yield a nontrivial element of $ \ker \wp_c $ (since $ v_p(\zeta_p) = e_K/(p-1) $).
    
    If $ c = b^{p-1} $ is a $ (p-1) $st power in $ \OO_K $, then $ y = b $ is an element of $ \ker \wp_c $. Note that $ u = 1 + b \pi^j $ lifts to a nontrivial $ p $th root of unity $ \zeta_p $, since $ u^p \equiv 1 $ mod $ \pi^{i+1} $ has (by the Taylor series again) a $ p $th root that is $ 1 $ mod $ \pi^{j+1} $. Note that $ \ker \wp_c $ has dimension only $ 1 $, since $ b $ is unique up to $ \mu_{p-1} = \FF_p^\cross $. Consequently $ \coker \wp_c \isom \bar{U}_i/\bar{U}_{i+1} $ has dimension exactly $ 1 $.
    
    This does not tell us how to find a generator for $ \bar{U}_i/\bar{U}_{i+1} $. For this, put $ y = b y' $ so
    \[
    (1 + b y' \pi^j) \equiv 1 + b^p (y^p - y) \pi^i \equiv 1 + b^p \wp(y) \mod \pi^{i+1},
    \]
    where $ \wp(y) = y^p - y $ is the usual Artin-Schreyer map. Since $ y^p $ and $ y $ are Galois conjugates over $ \FF_p $, we have $ \tr_{k_K/\FF_p}(y^p - y) = 0$, so if $ a \in k_K $ is an element with nonzero trace to $ \FF_p $, then $ \bar{U}_i/\bar{U}_{i+1} $ is generated by
    \[
    1 + a b^p \pi^i \equiv 1 + a b c \pi^i \equiv 1 + ap(\zeta_p - 1) \mod \pi^{i+1}. \qedhere
    \]
  \end{itemize}
\end{proof}
We draw two corollaries of the above method.

\begin{cor} \label{cor:Sh_basis}
  $ \bar U $ has dimension
  \[
  [K : \QQ_p] + 1 + \1_{\mu_p \subseteq K}
  \]
  with a basis consisting of $ \pi_K $ and (arbitrary lifts of) the elements in the bases of $ \bar U_i/\bar U_{i+1} $ in Proposition \ref{prop:Sh_basis}.
\end{cor}
We call this basis the \emph{Shafarevich basis} for $ U $.
\begin{rem}
  The dimension of $\bar U$ follows also from the Euler-characteristic computation
  \[
    \frac{\size{H^0(K, \mu_p)} \size{H^2(K, \mu_p)}}{\size{H^1(K, \mu_p)}} = p^{-[K : \QQ_p]}
  \]
\end{rem}
\begin{proof}
  Clearly $ \bar{U}/\bar{U}_0 \isom \C_p $ is generated by $ \pi_K $. The result follows from the composition series \eqref{eq:comp_series_K}.
\end{proof}

\begin{cor} \label{cor:Sh_dist_red}
  An element $ x \in \bar U $ belongs to $ \bar{U}_i $ if and only if, in the expansion of $ x $ in the Shafarevich basis, only the basis elements in $ U_i $ appear (to nonzero exponents).
\end{cor}
\begin{proof}
  Filtering $ \bar{U}_i $ by the $ \bar{U}_j $, for $ j \geq i $, we see that a basis for $ \bar{U}_i $ is given by the portion of the Shafarevich basis coming from $ \bar{U}_j/\bar{U}_{j+1} $ for $ j \geq i $. This is just the basis elements that lie in $ U_i $.
\end{proof}

The following simple result is one I have not seen in the literature before:
\begin{cor} \label{cor:Kum_cyclo}
  The cyclotomic extension $ \QQ_p[\mu_p] $ is isomorphic to $ \QQ_p[\sqrt[p-1]{-p}] $ and has Kummer element $ -p $ as a $ \mu_{p-1} $-torsor.
\end{cor}
\begin{proof}
  Let $ K = \QQ_p[\sqrt[p-1]{-p}] $, with uniformizer $ \pi_{K} = \sqrt[p-1]{-p} $. In the notation of the intimate unit case of Proposition \ref{prop:Sh_basis}, we have $ i = p $, $ j = 1 $,
  \[
  c = \frac{-p}{\pi_K^{p-1}} = 1,
  \]
  and we can take $ b = \sqrt[p-1]{c} = 1 $. Accordingly, $ 1 + \pi_K $ is congruent mod $ \pi_K^2 $ to a unique $ p $th root of unity $ \zeta_p $. Note that $ \zeta_p^i \equiv 1 + i \pi_K $ mod $ \pi_K^2 $. Since $ \mu_{p-1} $ acts on $ \pi_K $ by multiplication, it must act on the powers of $ \zeta_p $ by $ \omega $. Hence $ K \isom \QQ_p[\mu_p] $ as $ \mu_{p-1} $-torsors.
\end{proof}

In the rest of this section we will study how $ \bar{U} = \bar{U}(K) $ behaves under field extension. We use the following notational conventions:

{    
  Elements $u \in K^\cross/(K^\cross)^p$ are classified by their \emph{distance,} by which we mean the closest distance of a representative from $1$:
  \[
  d(u) = d_K(u) = \min_{\y \in K^\cross} \size{u \y^p - 1}.
  \]
  Here the absolute value is the local one on $K$. (We could choose a normalization of this absolute value, but we prefer to express $d(u)$ in terms of an undetermined $\size{\pi_K}$ and $\size{p} = \size{\pi_K}^{e_{K}}$.) Note that $d(u) \leq 1$, since $u$ can always be taken to have nonnegative valuation. Also, it is easy to see that $d(uv) \leq \max\{d(u), d(v)\}$, so $d$ defines a norm on $K^\cross/(K^\cross)^p$.
  
  For ease in stating theorems involving distances, we note that an ideal (or even a fractional ideal) $\aa$ of a local field $K$ is uniquely determined by the largest absolute value of its elements, which we denote by $\size{\aa}$. We have
  \[
  \aa = \{\x \in K : \size{\x} \leq \size{\aa}\}.
  \]
  
  For any real $d > 0$, let $B_{\leq d}$ denote the closed ball of radius $d$ about $1$ in $K^\cross/(K^\cross)^p$:
  \[
  B_{\leq d} = \{u \in K^\cross/(K^\cross)^p : d(u) \leq d\}
  \]
  and likewise for $B_{< d}$. It is easy to prove that these are subgroups. If $\ff \subsetneq \OO_K$ is an ideal, then $B_{\leq \size{\ff}}$ is the projection of $\U_\ff$; but this fails for $\ff = (1)$.
}

Note the following:
\begin{lem} \label{lem:dist_base_change}
  If $L/K$ is an extension of local fields whose degree $n$ is prime to $p$, then the canonical map from $K^\cross/(K^\cross)^p$ to $L^\cross/(L^\cross)^p$ is injective and preserves distance: that is, for every $u \in K^\cross/(K^\cross)^p$,
  \[
  d_K(u) = d_L(u).
  \]
\end{lem}
\begin{proof}
  The injectivity follows from the fact that if $u \in K$, then $N_{L/K} = u^n$ and $n$ is prime to $p$.
  
  It is obvious that $d_L(u) \leq d_K(u)$, so it suffices to prove the opposite inequality. It's easy to see that if $\x \in \OO_L$, then
  \[
  \size{N_{L/K}(\x) - 1} \leq \size{\x - 1}.
  \]
  Let $\y \in L^\cross$ achieve $\size{uy^p - 1} = d_L(u)$. Then
  \[
  d_L(u) = \size{uy^p - 1} \geq \size{N_{L/K}(u \y^p) - 1} = \size{u^{[L:K]} N_{L/K}(\y)^p - 1} \geq d_K\( u^{[L:K]} \).
  \]
  But $u$ is a power of $u^{[L:K]}$ up to $p$th powers, so $d_K(u) \leq d_K\( u^{[L:K]} \)$, completing the proof.
\end{proof}

Assume $ K \supseteq \QQ_p $. We first parametrize $ M $ itself. By (classical) Kummer theory, we have canonical isomorphisms
\[
\Hom(G_K, \Aut(\C_p)) \isom H^1(G_K, (\ZZ/p\ZZ)^\cross) \isom H^1(G_K, \mu_{p-1}) \isom K^\cross/(K^\cross)^{p-1}
\]
where the isomorphism $ (\ZZ/p\ZZ)^\cross \isom \mu_{p-1} $ is given by Teichm\"uller lift. (Note that we do not need to pick a generator of $ (\ZZ/p\ZZ)^\cross $ to do this.) We have the following:
\begin{lem} \label{lem:Kum_Tate_dual}
  If $ M $ is the Galois module with underlying group $ \C_p $ corresponding to the Kummer element $ \beta \in K^\cross/(K^\cross)^{p-1} $, then the Tate dual $ M' $ has Kummer element $ \beta' = -p/\beta $.
\end{lem}
\begin{proof}
  When $ M $ is trivial, the result was proved as Corollary \ref{cor:Kum_cyclo}. The lemma then follows by noting that if $ M_1 $, $ M_2 $ are cyclic Galois modules of order $ p $ with Kummer elements $ \beta_1 $, $ \beta_2 $, respectively, then $ \Hom(M_1,M_2) $ is also cyclic of order $ p $ and has Kummer element $ \beta_2/\beta_1 $.
\end{proof}

If $ T' $ has $ r $ field factors (all necessarily isomorphic to one $ T'_1 $), then, by Corollary \ref{cor:Sh_basis},
\[
\dim_{\FF_p} T'^\cross / (T'^\cross)^p = \begin{cases}
  n + 2r, & \mu_p \subseteq T'_1 \\
  n + r & \text{otherwise}.
\end{cases}
\]
The group $ T'^\cross / (T'^\cross)^p $ is a representation of $ (\ZZ/p\ZZ)^\cross $ over the field $ \FF_p $. Since $\FF_p$ has the $(p-1)$st roots of unity, such a representation splits as a direct sum of $1$-dimensional representations; there are $p-1$ of these, and they are the powers of the standard representation
\[
\omega : (\ZZ/p\ZZ)^\cross \to \GL_1(\FF_p)
\]
given by the obvious isomorphism. By Theorem \ref{thm:Kummer_new}, $ H^1(K,M) $ is parametrized by the $ \omega $-isotypical component of $ T'^\cross/(T'^\cross)^p $, which we denote by $ T'^\cross_{\omega} $ for brevity.

We can reduce the problem from $ T' $ to $ T'_1 $ in the following way:
\begin{lem}\label{lem:torsor_spl}
  Let $T'$ be a $ \mu_t $-torsor, $ t | p-1 $. Let $ T'_1 $ be the field factor of $ T' $, and let $ r = [T'_1 : K] $. Then:
  \begin{enumerate}[$ ( $a$ ) $]
    \item\label{it:spl_stab} The subgroup fixing $ T'_1 $ (as a set) is $ \mu_r \subseteq \mu_{p-1} \isom (\ZZ/p\ZZ)^\cross $, and $ T'_1 $ is a $ \mu_r $-torsor;
    \item\label{it:spl_sp} Projection to $ T'_1 $ defines an isomorphism $T'^{\cross}_\omega \isom (T'^{\cross}_1)_{\omega} $, where
    \begin{align}
      T'^\cross_{\omega} &=\left\{ \alpha \in T'^\cross / (T'^\cross)^p : \tau_{c}(\alpha) = \alpha^c \, \forall c \in \mu_t \right\} \\
      (T'^\cross_1)_{\omega} &= \left\{ \alpha \in T'^\cross_1 / (T'^\cross_1)^p : \tau_{c}(\alpha) = \alpha^c \, \forall c \in \mu_r \right\}.
    \end{align}
    \item\label{it:spl_ge} More generally, for any $s$ with $r|s|t$, the orbit $\mu_s(T'_1)$ consists of $s/n$ field factors whose product $T'_s$ is a $\mu_s$-torsor. Projection onto $T'_s$ and then onto $ T'_1 $ defines isomorphisms
    \[
    T'^\cross_{\omega} \cong (T'^\cross_s)_{\omega} \cong (T'^\cross_1)_{\omega}.
    \]
  \end{enumerate}
\end{lem}
\begin{rem}
  A result with much the same content, but in a slightly different setting, is proved by Del Corso and Dvornicich (\cite{DCD}, Proposition 7).
\end{rem}
\begin{proof}
  The torsor action must permute the field factors transitively; since $\mu_n$ is cyclic, a generator $\zeta_n$ must cyclically permute them, and the stabilizer of $T_1$ (as a set) is $\<\zeta^{n/r}\> = \mu_r$, proving\ref{it:spl_stab}. Since the action simply transitively permutes the coordinates of $T_1$, $T_1$ is a $\mu_r$-torsor. If we know the $T_1$-component $\alpha|_{T_1^\cross}$ of an $\alpha \in \(T^\cross/(T^\cross)^p\)_\omega$, then all the other components are uniquely determined by the eigenvector condition; it is only necessary for $\alpha|_{T_1^\cross}$ to behave properly under $\mu_r$, namely that $\alpha_{T_1^\cross} \in (T_1^\cross)_{\omega^k}$.
  
  This proves \ref{it:spl_sp}. Also, it is clear from our analysis that $\mu_s(T_1)$ is a $\mu_s$-torsor. Applying \ref{it:spl_sp} to this torsor proves \ref{it:spl_ge}.
\end{proof}

We now filter $ T'^\cross_1 $ as above to discover its $ \omega $-component.

\begin{prop}[\textbf{\textbf{the Shafarevich basis for $T'^\cross_{\omega}$}}] \label{prop:Sh_basis_eig}
  As above, let $T' = K[\sqrt[p-1]{\beta}]$ be the $(\ZZ/p\ZZ)^\cross$-torsor corresponding to the Tate dual $ M' $ of a cyclic Galois module $ M $ of order $ p $ with Kummer element $\beta \in K^\cross/(K^\cross)^{p-1}$, and let $ T'_1 $ be the field factor of $ T' $. Filter $ \bar{U} = T'^\cross_1/(T'^\cross_1)^p $ by the subgroups $ \bar U_i $ as in the previous subsection. Since the $ \mu_r $-torsor action on $ T'_1 $ preserves the valuation, each $ U_i $ is a subrepresentation of $ T'^\cross_1 / (T'^\cross_1)^p $. Then:
  
  \begin{enumerate}[$($a$)$]
    \item The $ \omega $-isotypical component of $ \bar U / \bar U_{0} $ has dimension $ 1 $ if $ M' $ is trivial, $ 0 $ otherwise.
    \item The $ \omega $-isotypical component of $ \bar U_i / \bar U_{i+1} $ has dimension
    \begin{itemize}
      \item $ f $ if
      \[
      0 < i < \frac{pe_{T'_1/\QQ_p}}{p-1}, \quad i \nequiv 0 \mod p, \quad i \equiv \frac{e_{T'_1/K} v_K(\beta)}{p-1} \mod e_{T'_1/K};
      \]
      \item $ 1 $ if $ i = \frac{pe_{T'_1/\QQ_p}}{p-1} $ and $ M $ is trivial;
      \item $ 0 $ otherwise.
    \end{itemize}
  \end{enumerate}
\end{prop}
\begin{proof}
  Note that $ U/U_0 = \<\pi_{T'_1}\> $ is a copy of the trivial representation and that $ \omega $ is trivial (as a representation of $ \mu_r $) exactly when $ r = 1 $, that is, $ M' $ is trivial.
  
  Our convention for the Kummer map is such that
  \[
  \tau_{c}(\sqrt[p-1]{\beta}) = \tilde{c} \sqrt[p-1]{\beta}.
  \]
  By a standard result in Kummer theory, the degree $ r $ is the least integer such that $ \beta = \beta_1^{(p-1)/r} $ is a $ (p-1)/r $th power, and $ T'_1 = K[\sqrt[r]{\beta_1}] $ with
  \[
  \tau_{c}(\sqrt[r]{\beta_1}) = \tilde{c} \sqrt[r]{\beta_1}.
  \]
  Let $ f' = \gcd\big(v_K(\beta_1), r\big) $ and $ e' = r/f' $. We claim that $ e' $ and $ f' $ are respectively the ramification and inertia indices of $T'_1$ over $K$. By the Euclidean algorithm, we may choose integers $g$ and $h$ such that
  \begin{equation} \label{eq:Euclid}
    gv(\beta_1) - hr = f'.
  \end{equation}
  Construct the elements
  \[
  \pi' = \frac{\( \sqrt[r]{\beta_1} \)^{g}}{\pi_K^h} \textand u' = \frac{\( \sqrt[r]{\beta_1} \)^{e'}}{\pi_K^{v_K(\beta_1)/f'}}.
  \]
  Note that $ v_K(\pi') = 1/e' $, so 
  \begin{equation} \label{eq:e_geq}
    e_{T'_1/K} \geq e'.
  \end{equation}
  On the other hand, $ v_K(u') = 0 $, and $ u' $ is an $ f' $th root of the unit
  \[
  \beta_2 = \frac{\beta_1}{\pi_K^{v_K(\beta_1)}}.
  \]
  Note that $ \beta_2 $ is not an $ \ell $th power for any prime $ \ell | f' $, as otherwise $ \beta $ would be a $ (p-1)\ell/r $th power, contradicting what we know about $ r $. So the residue class of $ u' $ generates a degree-$ f' $ extension of $ k_K $ inside $ k_{T'_1} $; in particular,
  \begin{equation} \label{eq:f_geq}
    f_{T'_1/K} \geq f'.
  \end{equation}
  Equality must hold in \eqref{eq:e_geq} and \eqref{eq:f_geq}, so $ T'_1 $ has uniformizer $ \pi' $ and residue field generator $ u' $. Note that for $ c \in \FF_p^\cross $,
  \[
  \tau_c(\pi') = \tilde c^g \pi' \textand \tau_c(u') = \tilde c^{e'} u'.
  \]
  
  We now have what we need to compute the Galois action on $ \bar U_i/\bar U_{i+1} $. By Proposition \ref{prop:Sh_basis}, the space $ \bar U_i/\bar U_{i+1} $ is nonzero only for
  \[
  0 < i < \frac{p e_{T'_1}}{p-1}, \quad p \nmid i
  \]
  (the generic units) and possibly for $ i = pe_K/(p-1) $ also (the intimate units).
  
  We begin with the first case. Here $ \bar U_i/\bar U_{i+1} \isom U_i/ U_{i+1} $ has a basis
  \[
  1 + \pi'^i u'^j, \quad 0 \leq j < f'.
  \]
  For $ c \in \mu_r $,
  \begin{align*}
    \tau_c(1 + \pi'^i u'^j) &= 1 + \tilde c^{gi + e'j} \pi'^i u'^j \\
    & \equiv (1 + \pi'^i u'^j)^{c^{gi + e'j}} \mod \pi'^{i+1}.
  \end{align*}
  Thus the basis element $ 1 + \pi'^i u'^j $ generates a $ 1 $-dimensional $ \mu_r $-submodule of $ U_i/U_{i+1} $ isomorphic to $ \omega^{gi + e'j} $. Accordingly, we select the generic units satisfying
  \[
  g i + e'j \equiv 1 \mod r.
  \]
  Since $ r = e'f' $, there are exactly $ f' $ values of $ j $ satisfying this when
  \begin{equation}\label{eq:cong_i}
    e'| gi - 1,
  \end{equation}
  and none otherwise. By \eqref{eq:Euclid}, $ g $ is the multiplicative inverse of $ v(\beta_1)/f' = e' v(\beta)/(p-1) $ mod $ e' $, so we can rewrite \eqref{eq:cong_i} as
  \[
  i \equiv \frac{e' v(\beta)}{p-1} \mod e',
  \]
  as desired.
  
  As to the case that $ i = pe_{T'_1}/(p-1) $ (the intimate units), we simply note that, by Proposition \ref{thm:disc_Kummer_aff}, we have
  \[
  (\bar U_i)_\omega \cong H^1_\ur(K, M),
  \]
  so
  \[
  \size{(\bar U_i)_\omega} = \size{H^1_\ur(K,M)} = \size{H^0(K,M)}.
  \]
  (A direct computation of the torsor action on the intimate units is also possible; it turns out that $ \bar U_i \cong \mu_m $ as $ \mu_{m-1} $-modules.)
\end{proof}

By complete reducibility, we can get a basis for $ \bar U_\omega $ from the bases for its composition factors:
\begin{cor}
  If $M \isom \sf C_p$, then $ H^1(K,M) $ has dimension
  \[
  [K : \QQ_p] + \dim_{\FF_p} H^0(K, M) + \dim_{\FF_p} H^0(K, M'),
  \]
  with a basis consisting of appropriate lifts of the Shafarevich basis elements picked out by Proposition \ref{prop:Sh_basis_eig}.
\end{cor}
In particular, we have proved Theorem \ref{thm:levels}\ref{lev:size_all}.

\subsection{Proof of Theorem \ref{thm:levels}}
It now remains to recast the above results in terms of levels and offsets and prove the remaining parts of Theorem \ref{thm:levels}.

Let $ \alpha \in (T'^\cross_1)_\omega $, $ \alpha $ being a minimal-distance element. We consider the possibilities for the leading factor in $ \alpha $ with respect to the Shafarevich basis; this determines $ \size{\alpha - 1} $ by Corollary \ref{cor:Sh_basis}, and thence $ d := v_K(\Disc(L/K)) $ and hence the level and offset of $ \alpha $.

\begin{itemize}
  \item If $ \alpha $ is led by the uniformizer, then $ M' $ is trivial. From Theorem \ref{thm:disc_Kummer_aff}, we get $ d = pe_K + p - 1 $, so $ \ell(\alpha) = -1 $ and $ \theta(\alpha) = -1 $.
  \item If $ \alpha $ is led by a generic unit, then we have $ v_{T'_1}(\alpha - 1) = i,$ where $ i $ is an integer satisfying
  \[
  0 < i < \frac{pe_{T'_1/\QQ_p}}{p-1}, \quad i \nequiv 0 \mod p, \quad i \equiv \frac{e_{T'_1/K} v_K(\beta)}{p-1} \mod e_{T'_1/K},
  \]
  and each of these values is attained by some $ \alpha $. Using the one-to-one correspondence of Theorem \ref{thm:disc_Kummer_aff}, we get that $ v_{K}(\disc L) $ attains exactly the values $ d $ such that
  \[
  p-1 < d < pe_{K} + p - 1, \quad d \nequiv -1 \mod p, \quad d \equiv e_K - v_K(\beta) \mod p - 1.
  \]
  Thus $ 0 \leq \ell(\alpha) < e_K $, $ 0 \leq \theta(\alpha) \leq p-2 $, and $ \theta $ is determined by $ \ell $ via the condition mod $ p-1. $
  \item If $ \alpha $ is led by an intimate unit or $ \alpha = 1 $, then $ d = v_K(\disc L) = v_K(\disc T') $ was already computed in proving Theorem \ref{thm:disc_Kummer_aff}. Since $ T' $ is a product of tamely ramified extensions, we have $ d < [T':K] = p-1 $, so $ \ell = 0 $.
\end{itemize}
This proves \ref{lev:offset} and \ref{lev:ur}.

In the case that $\alpha$ is led by a generic unit of level $\ell$, $0\leq \ell < e$, there are alternative ways to characterize $\theta$. We have
\[
  \theta = \left\{\frac{\ell - v_K(\beta)}{p-1}\right\} (p-1),
\]
from which
\[
  v_K(\disc L) = p(e - \ell) + \theta = p(e - \ell) + \left\{\frac{\ell - v_K(\beta)}{p-1}\right\} (p - 1)
\]
Since $v_K(\disc L)$ is in bijection with $\size{\alpha - 1}$ by Theorem \ref{thm:disc_Kummer_aff}, we can likewise determine
\[
  v_K(\alpha - 1) = \floor{\frac{p\ell - v_K(\beta)}{p-1}} + 1 + \frac{v_K(\beta)}{p-1}.
\]
Item \ref{lev:distance} demands that we invert this to determine how the level $\ell(\alpha)$ changes as $\alpha$ ranges in a ball
\[
  \size{\alpha - 1} \leq d.
\]
If $d < d_{\min}$, then every $\alpha$ in this ball is a $p$th power (by Proposition \ref{prop:Sh_basis}, or simply by noting that the Taylor series for $p$th root converges on this ball), so the range of $[\alpha]$ is $\{1\} = \L_{e + 1}$. If $d \geq d_{\min}$, then the intimate units are certainly included, so the range of $[\alpha]$ is at least $\L_e$; it also includes all generic units whose levels $\ell$ satisfy
\begin{align*}
  v_K(\alpha - 1) &\geq \frac{\log d}{\log \size{\pi_K}} \\
  \floor{\frac{p\ell - v_K(\beta)}{p-1}} + 1 + \frac{v_K(\beta)}{p-1} &\geq \frac{\log d}{\log \size{\pi_K}} \\
  \floor{\frac{p\ell - v_K(\beta)}{p-1}} &\geq \frac{\log d}{\log \size{\pi_K}} - 1 - \frac{v_K(\beta)}{p-1}.
\end{align*}
Using the exchange
\[
  \floor{x} \geq y \iff \floor{x} \geq \ceil{y} \iff x \geq \ceil{y},
\]
valid for all real numbers $x$ and $y$, we can get this into a form solvable for $\ell$:
\begin{align*}
  \frac{p\ell - v_K(\beta)}{p-1} &\geq \ceil{\frac{\log d}{\log \size{\pi_K}} - 1 - \frac{v_K(\beta)}{p-1}} \\
  \ell &\geq \frac{v_K(\beta)}{p} + \frac{p-1}{p} \ceil{\frac{\log d}{\log \size{\pi_K}} - 1 - \frac{v_K(\beta)}{p-1}}.
\end{align*}
So the range of $[\alpha]$ is $\L_i$, where
\begin{equation}\label{eq:dist_to_lev}
  \Biggl\lceil\frac{v_K(\beta)}{p} + \frac{p-1}{p} \ceil{\frac{\log d}{\log \size{\pi_K}} - 1 - \frac{v_K(\beta)}{p-1}}\Biggr\rceil,
\end{equation}
as claimed in \ref{lev:distance}.

For \ref{lev:subgp}, we note that as $d$ decreases from $1$ to $d_{\min}$, the corresponding $i$ in \eqref{eq:dist_to_lev} hits every value from $0$ to $e$, since the argument to the outer ceiling increases by jumps of $(p - 1)/p < 1$. So each $ \L_i $ is a subgroup. For \ref{lev:size_Li}, we note that $ \size{\L_0} = \size{H^0(K,M)} $, while for $ 1 \leq i \leq e, $
\[
\L_{i+1}/\L_i \cong \bar U_j/\bar{U}_{j+1}
\]
has $ p^f = q $ elements, where $ j $ is the unique value of $ v_{T'_1}(\alpha - 1) $ for values of $ \alpha $ having level $ i $.

Finally, we have claimed a relation \ref{lev:perp} regarding how level spaces interact with the Tate pairing. In the case of the Hilbert pairing, the result we need is as follows:
\begin{lem}[\textbf{an explicit reciprocity law}] \label{lem:Hilb_prod_size}
  Let $ K $ be a local field with $ \mu_p \subseteq K $. If $ \alpha, \beta \in \OO_K^\cross $ satisfy
  \[
  \size{\alpha - 1} \cdot \size{\beta - 1} < d_{\min} = \size{p}^{p/(p-1)},
  \]
  then the Hilbert pairing $ \<\alpha, \beta\>_K $ vanishes.
\end{lem}
\begin{proof}
  This is a consequence of the \emph{conductor-discriminant formula} (see Neukirch \cite{Neukirch}, VII.11.9): For a Galois extension $L/K$,
  \[
    \Disc(L/K) = \prod_\chi \ff(\chi)^{\chi(1)},
  \]
  where $\chi$ ranges over the irreducible characters of $\Gal(L/K)$. Here we apply the formula to $L = K[\sqrt[p]{\alpha}]$. Scale $\alpha$ by $p$th powers to be as close to $1$ as possible. If $\alpha = 1$ or $\alpha$ is an intimate unit, the Hilbert pairing clearly vanishes since $L$ is unramified and $\beta$ is a unit. So we can assume that $L$ is a ramified extension of degree $p$. Then there are $p$-many characters on $\Gal(L/K)$, all of dimension $1$. One is the trivial character, whose conductor is $1$. The others all have the same conductor $\ff$, so
  \[
    \Disc(L/K) = \ff^{p-1}.
  \]
  By Theorem \ref{thm:disc_Kummer_aff}, we have
  \[
    \Disc(L/K) \sim \frac{p^p \cdot \pi_K^{p-1}}{(\alpha-1)^{p-1}},
  \]
  so $\ff$ is generated by any element $f$ with
  \[
    \size{f} = \Size{\frac{p^{p/(p - 1)} \cdot \pi_K}{\alpha - 1}}.
  \]
  Note that $d_{\min} = \size{p}^{p/(p-1)}$ is actually an attainable norm of an element of $K$, namely $(\zeta_p - 1)^p$.
  By the given inequality, $\beta \equiv 1 \mod \ff$ which implies that the Hilbert symbol
  \[
    \<\alpha, \beta\>_K = \phi_{L/K}(\beta)
  \]
  vanishes.
\end{proof}

If $ 0 \leq i \leq e $, $ \alpha \in \L_i(M) $, and $ \beta \in \L_{e-i}(M) $, then it is easy to verify that the hypothesis of Lemma \ref{lem:Hilb_prod_size} holds in each field factor of $ T[\mu_m] $, in which the Hilbert pairing is being computed. Hence
\[
\L_i(M)^\perp \supseteq \L_{e-i}(M').
\]
However, since
\[
\size{\L_i(M)} \cdot \size{\L_{e-i}(M')} = q^{e_K} \cdot \size{H^0(K,M)} \cdot \size{H^0(K,M')} = \size{H^1(K,M)},
\]
equality must hold.
\qed

\subsection{The tame case}
If $ K $ is a tame local field, that is, $ \ch k_K \neq p $, the structure of $ H^1(K, M) $ is well known. We put
\[
e = 0, \quad \L_{-1} = \{0\}, \quad \L_0 = H^1_\ur(K,M), \quad \L_1 = H^1(K,M)
\]
and observe that Theorem \ref{thm:levels}\ref{lev:ur}, \ref{lev:size_Li}, \ref{lev:size_all}, \ref{lev:perp} and Corollary \ref{cor:levels} still hold.

The \emph{wild function field} case $K = \FF_{p^r} \laurent{t}$ admits a similar treatment, but now the number of levels is infinite. We do not address this case here.

\part{Composed varieties}
\label{part:composed}

\section{Composed varieties}
\label{sec:composed}
It has long been noted that orbits of certain algebraic group actions on varieties over a field $K$ parametrize rings of low rank over $K$, which can also be identified with the cohomology of small Galois modules over $K$. The aim of this section is to explain all this in a level of generality suitable for our applications.

\begin{defn}
  Let $K$ be a field and $\bar K$ its separable closure. A \emph{composed variety} over $K$ is a quasi-projective variety $V$ over $K$ with an action of a quasi-projective algebraic group $\Gamma$ over $K$ such that:
  \begin{enumerate}[$($a$)$]
    \item $V$ has a $K$-rational point $x_0$;
    \item the $\bar K$-points of $V$ consist of just one orbit $\Gamma(\bar K) x_0$;
    \item the point stabilizer $M = \Stab_{\Gamma(\bar K)} x_0$ is a \emph{finite abelian} subgroup.
  \end{enumerate}
\end{defn}

The term \emph{composed} is derived from Gauss composition of binary quadratic forms and the ``higher composition laws'' of the work of Bhargava and others, from which we derive many of our examples.

\begin{prop}~
  \begin{enumerate}[$($a$)$]
    \item Once a base orbit $\Gamma(K)x_0$ is fixed, there is a natural injection
    \[
    \psi : \Gamma(K)\backslash V(K) \hookrightarrow H^1(K,M)
    \]
    by which the orbits $\Gamma(K)\backslash V(K)$ parametrize some subset of the Galois cohomology group $H^1(K, M)$.
    \item The $\Gamma(K)$-stabilizer of every $x \in V(K)$ is canonically isomorphic to $H^0(K, M)$.
  \end{enumerate}  
\end{prop}
\begin{proof}
  \begin{enumerate}[$($a$)$]
    \item Let $x \in V(K)$ be given. Since there is only one $\Gamma(\bar K)$-orbit, we can find $\gamma \in \Gamma(\bar K)$ such that $\gamma(x_0) = x$. For any $g \in \Gal(\bar K/K)$, $g(\gamma)$ also takes $x_0$ to $x$ and so differs from $\gamma$ by right-multiplication by an element in $\Stab_{\Gamma(\bar K)} x_0 = M$. Define a cocycle $\sigma_x : \Gal(\bar K/K) \to M$ by
    \[
    \sigma_x(g) = g(\gamma) \cdot \gamma^{-1}.
    \]
    It is routine to verify that
    \begin{itemize}
      \item $\sigma_x$ satisfies the cocycle condition $\sigma_x(gh) = \sigma_x(g) \cdot g(\sigma_x(h))$ and hence defines an element of $H^1(K,M)$;
      \item If a different $\gamma$ is chosen, then $\sigma_x$ changes by a coboundary;
      \item If $x$ is replaced by $\alpha x$ for some $\alpha \in \Gamma_K$, the cocycle $\sigma_x$ is unchanged;
      \item If the basepoint $x_0$ is replaced by $\alpha x_0$ for some $\alpha \in \Gamma(K)$, the cocycle $\sigma_x$ is unchanged, up to identifying $M$ with $\Stab_{\Gamma(\bar K)} (\alpha x_0) = \alpha M \alpha^{-1}$ in the obvious way. (This is why we can fix merely a base \emph{orbit} instead of a basepoint.)
    \end{itemize}
    So we get a map
    \[
    \psi : \Gamma(K)\backslash V(K) \to H^1(K,M).
    \]
    We claim that $\psi$ is injective. Suppose that $x_1, x_2 \in V(K)$ map to equivalent cocycles $\sigma_{x_1}$, $\sigma_{x_2}$. Let $\gamma_i \in \Gamma(\bar K)$ be the associated transformation that maps $x_0$ to $x_i$. By right-multiplying $\gamma_1$ by an element of $M$, as above, we can remove any coboundary discrepancy and assume that $\sigma_{x_1} = \sigma_{x_2}$ on the nose. That is, for every $g \in \Gal(\bar K/K)$,
    \[
    g(\gamma_1) \cdot \gamma_1^{-1} = g(\gamma_2) \cdot \gamma_2^{-1},
    \]
    which can also be written as
    \[
    g(\gamma_2 \gamma_1^{-1}) = \gamma_2 \gamma_1^{-1}.
    \]
    Thus, $\gamma_2 \gamma_1^{-1}$ is Galois stable and hence defined over $K$. It takes $x_1$ to $x_2$, establishing that these points lie in the same $\Gamma(K)$-orbit, as desired.
    \item If $\gamma(x_0) = x$, then the $\Gamma(\bar K)$-stabilizer of $x$ is of course $\gamma M \gamma^{-1}$. We claim that the obvious map
    \begin{align*}
      M \to \gamma M \gamma^{-1} \\
      \mu \mapsto \gamma \mu \gamma^{-1}
    \end{align*}
    is an isomorphism of Galois modules. We compute, for $g \in \Gal(\bar K/K)$,
    \[
    g\(\gamma \mu \gamma^{-1}\) = g(\gamma) g(\mu) g(\gamma)^{-1}
    = \gamma \sigma_x(g) g(\mu) \sigma_x(g)^{-1} \gamma^{-1}
    = \gamma g(\mu) \gamma^{-1},
    \]
    establishing the isomorphism. In particular, the Galois-stable points $\Stab_{\Gamma(K)} x_0 = H^0(K,M)$ are the same at $x$ as at $x_0$. Note the crucial way that we used that $M$ is abelian. By the same token, the identification of stabilizers is independent of $\gamma$ and is thus canonical. \qedhere
  \end{enumerate}
\end{proof}

The base orbit is distinguished only insofar as it corresponds to the zero element $0 \in H^1(K,M)$. Changing base orbits changes the parametrization minimally:
\begin{prop}
  The parametrizations $\psi_{x_0}, \psi_{x_1} : \Gamma(K)\backslash V(K) \to H^1(K,M)$ corresponding to two basepoints $x_0, x_1 \in V(K)$ differ only by translation:
  \[
  \psi_{x_1}(x) = \psi_{x_0}(x) - \psi_{x_0}(x_1),
  \]
  under the isomorphism between the stabilizers $M$ established in the previous proposition.
\end{prop}
\begin{proof}
  Routine calculation.
\end{proof}

While $\psi$ is always injective, it need not be surjective, as we will see by examples in the following section.
\begin{defn}~
  \begin{enumerate}[$($a$)$]
    \item A composed variety is \emph{full} if $\psi$ is surjective, that is, it includes a $\Gamma(K)$-orbit for every cohomology class in $H^1(K,M)$.
    \item If $K$ is a global field, a composed variety is \emph{Hasse} if for every $\alpha \in H^1(K,M)$, if the localization $\alpha_v \in H^1(K_v,M)$ at each place $v$ lies in the image of the local parametrization
    \[
    \psi_v : V(K_v)\backslash \Gamma(K_v) \to H^1(K_v, M),
    \]
    then $\alpha$ also lies in the image of the global parametrization $\psi$.
  \end{enumerate}
\end{defn}
\subsection{Examples}
In this section, $K$ is any field not of one of finitely many bad characteristics for which the exposition does not make sense.

\begin{examp}\label{ex:Kummer}
  The group $\Gamma = \GG_m$ can act on the variety $V = \AA^1 \backslash \{0\}$, the punctured affine line, by 
  \[
  \lambda (x) = \lambda^n \cdot x.
  \]
  There is a unique $\bar K$-orbit. The point stabilizer is $\mu_n$, and the parametrization corresponding to this composed variety (choosing basepoint $x_0 = 1$) is none other than the Kummer map
  \[
  K^\cross / (K^\cross)^n \to H^1(K, \mu_n).
  \]
  That $V$ is full follows from Hilbert's Theorem 90.
\end{examp}
\begin{examp}\label{ex:bin_cubic}
  Let $V$ be the variety of binary cubic forms $f$ over $K$ with fixed discriminant $D_0$. This has an algebraic action of $\SL_2$, which is transitive over $\bar K$ (essentially because $\PSL_2$ carries any three points of $\PP^1$ to any other three), and there is a ready-at-hand basepoint
  \[
  f_0(X,Y) = X^2 Y - \frac{D}{4} Y^3.
  \]
  The point stabilizer $M$ is isomorphic to $\ZZ/3\ZZ$, but twisted by the character of $K(\sqrt{D})$; that is, $M \isom \{0, \sqrt{D}, -\sqrt{D}\}$ as sets with Galois action. Coupled with the appropriate higher composition law (Theorem \ref{thm:hcl_cubic_ring}), this recovers the parametrization of cubic \'etale algebras with fixed quadratic resolvent by $H^1(K,M)$ in Proposition \ref{prop:Gal_mod}. To see that it is the same parametrization, note that a $\gamma \in \Gamma(\bar K/K)$ that takes $f_0$ to $f$ is determined by where it sends the rational root $[1:0]$ of $f_0$, so the three $\gamma$'s are permuted by $\Gal(\bar K/K)$ just like the three roots of $f$. In particular, $V$ is full.
\end{examp}
\begin{examp}
  Continuing with the sequence of known ring parametrizations, we might study the variety $V$ of pairs of ternary quadratic forms with fixed discriminant $D_0$. This has one orbit over $\bar K$ under the action of the group $\Gamma = \SL_2 \cross \SL_3$; unfortunately, the point stabilizer is isomorphic to the alternating group $A_4$, which is not abelian.
  
  So we narrow the group, which widens the ring of invariants and requires us to take a smaller $V$. We let $\Gamma = \SL_3$ alone act on pairs $(A, B)$ of ternary quadratic forms, which preserves the resolvent
  \[
  g(X,Y) = 4 \det\(AX + BY\),
  \]
  a binary cubic form. We let $V$ be the variety of $(A,B)$ for which $g = g_0$ is a fixed separable polynomial. These parametrize quartic \'etale algebras $L$ over $K$ whose cubic resolvent $R$ is fixed. There is a natural base orbit $(A_0, B_0)$ whose associated $L \isom K \cross R$ has a linear factor. The point stabilizer $M \isom \ZZ/2\ZZ \cross \ZZ/2\ZZ$, with the three non-identity elements permuted by $\Gal(\bar K / K)$ in the same manner as the three roots of $g_0$. We have reconstructed the parametrization of quartic \'etale algebras with fixed cubic resolvent by $H^1(K,M)$ in Proposition \ref{prop:Gal_mod}. In particular, $V$ is full.
\end{examp}
\begin{examp}
  Alternatively, we can consider the space $V$ of binary quartic forms whose invariants $I = I_0$, $J = J_0$ are fixed. The orbits of this space have been found useful for parametrizing $2$-Selmer elements of the elliptic curve $E : y^2 = x^3 + xI + J$, because the point stabilizer is $M \isom \ZZ/2\ZZ \cross \ZZ/2\ZZ$ with the Galois-module structure $E[2]$. This space $V$ embeds into the space of the preceding example via a map which we call the \emph{Wood embedding} after its prominent role in Wood's work \cite{WoodBQ}:
  \begin{align*}
    f &\mapsto (A,B) \\
    ax^4 + bx^3y + cx^2 y^2 + dxy^3 + ey^4 &\mapsto
    \(
    \begin{bmatrix}
      & & 1/2 \\
      & -1 & \\
      1/2 & &
    \end{bmatrix},
    \begin{bmatrix}
      a & b/2 & c/3 \\
      b/2 & c/3 & d/2 \\
      c/3 & d/2 & e
    \end{bmatrix}\).
  \end{align*}
  In general, $V$ is \emph{not} full. For instance, over $K = \RR$, if $E$ has full $2$-torsion, there are only three kinds of binary quartics over $\RR$ with positive discriminant (positive definite, negative definite, and those with four real roots) which cover three of the four elements in $H^1(\RR, \ZZ/2\ZZ)$. Two of these three (positive definite, four real roots) form the subgroup of elements whose corresponding $E$-torsor $z^2 = f(x,y)$ is soluble at $\infty$: these are the ones we retain when studying $\Sel_2 E$. The fourth element of $H^1(\RR, \ZZ/2\ZZ)$ yields \'etale algebras whose $(A,B)$ has
  \[
    A = \begin{bmatrix}
      1/2 & & \\
      & 1 & \\
      & & 1/2
    \end{bmatrix},
  \]
  a conic with no real points. However, over global fields, it is possible to show that $V$ is Hasse, using the Hasse-Minkowski theorem for conics.
\end{examp}

\begin{rem}
Because of the extreme flexibility afforded by general varieties, it is reasonable to suppose that any finite $K$-Galois module $M$ appears as the point stabilizer of some full composed variety over $K$. However, we do not pursue this question here.
\end{rem}

\subsection{Integral models; localization of orbit counts}
Let $K$ be a number field and $\OO_K$ its ring of integers. Let $(V, \Gamma)$ be a composed variety, and let $(\V, \G)$ be an \emph{integral model,} that is, a pair of a flat separated scheme and a flat algebraic group over $\OO_K$ acting on it, equipped with an identification of the generic fiber with $(V, \Gamma)$. Then $\G(\OO_K) \hookrightarrow \Gamma(K)$, and the $\Gamma(K)$-orbits on $V(K)$ decompose into $\G(\OO_K)$-orbits.

\begin{lem}[\textbf{localization of global class numbers}] \label{lem:loczn_main}
Let $(\V, \G)$ be an integral model for a composed variety $(V, \Gamma)$. For each place $v$, let
\[
  w_v : \G(\OO_v) \backslash \V(\OO_v) \to \CC
\]
be a function on the local orbits, which we call a \emph{local weighting.} Suppose that:
\begin{enumerate}[$($i$)$]
  \item $(V, \Gamma)$ is Hasse.
  \item $\G$ has \emph{class number one}, that is, the natural localization embedding
  \[
    \G(\OO_K) \backslash \Gamma(K) \hookrightarrow
    \bigoplus_v \G(\OO_v) \backslash \Gamma(K_v)
  \]
  is surjective.
  \item For each place $v$, there are only finitely many orbits of $\G(\OO_v)$ on $\V(\OO_v)$. This ensures that the \emph{weighted local orbit counter}
  \begin{align*}
    g_{v, w_v} : H^1(K_v, M) &\to \CC \\
    \alpha &\mapsto \sum_{\substack{\G(\OO_{K_v})\gamma \in \G(\OO_{v}) \backslash \Gamma(K_v) \\ \text{such that }\gamma x_\alpha \in \V(\OO_v)}} w_v(\gamma x_\alpha)
  \end{align*}
  takes finite values. (Here $x_\alpha$ is a representative of the $\Gamma(K_v)$-orbit corresponding to $\alpha$. If there is no such orbit because $V$ is not full, we take $g_{v, w_v}(\alpha) = 0$.)
  \item\label{loczn:one_orbit} For almost all $v$, $\G(\OO_v)\backslash \V(\OO_v)$ consists of at most one orbit in each $\Gamma(K_v)$-orbit, and $w_v = 1$ identically.
\end{enumerate}
Then the global integral points $\V(\OO_K)$ consist of finitely many $\G(\OO_K)$-orbits, and the \emph{global weighted orbit count} can be expressed in terms of the $g_{v, w_v}$ by
\begin{equation} \label{eq:x_wtd}
  h_{\{w_v\}} \coloneqq \sum_{\G(\OO_K) x \in \G(\OO_K) \backslash \V(\OO_K)} \frac{\prod_v w_v(x)}{\size{\Stab_{\G(\OO_K)} x}} = \frac{1}{\size{H^0(K,M)}} \sum_{\alpha \in H^1(K, M)} \prod_v g_{v, w_v}(\alpha).
\end{equation}
\end{lem}
\begin{proof}
Grouping the $\G(\OO_K)$-orbits into $\Gamma(K)$-orbits, it suffices to prove that for all $\alpha \in H^1(K,M)$,
\begin{equation} \label{eq:x_orbit}
  \sum_{\G(\OO_K)x \subseteq \Gamma(K)x_\alpha } \frac{\prod_v w_v(x)}{\size{\Stab_{\G(\OO_K)} x}} = \frac{1}{\size{H^0(K,M)}} \prod_v g_{v, w_v}(\alpha).
\end{equation}
If there is no $x_\alpha$, the left-hand side is zero by definition, and at least one of the $g_{v, w_v}(\alpha)$ is also zero since $V$ is Hasse. So we fix an $x_\alpha$. The right-hand side of \eqref{eq:x_orbit}, which is finite by hypothesis \ref{loczn:one_orbit} since $\alpha$ is unramified almost everywhere, can be written as
\[
  \frac{1}{\size{H^0(K,M)}} \sum_{\{\G(\OO_v)\gamma_v\}_v} \prod_{v} w_v(\gamma_v x_\alpha),
\]
the sum being over systems of $\gamma_v \in \Gamma(K_v)$
such that $\gamma_v x_\alpha$ is $\OO_v$-integral. Since $\G$ has class number one, each such system glues uniquely to a global orbit $\G(\OO_K)\gamma, \gamma \in \Gamma(K)$, for which $\gamma x_\alpha$ is $\OO_v$-integral for all $v$, that is, $\OO_K$-integral. Thus the right-hand side of \eqref{eq:x_orbit} is now transformed to
\[
  \frac{1}{\size{H^0(K,M)}} \sum_{\substack{\G(\OO_K)\gamma \\ \quad \gamma x_\alpha \in \V(\OO_K)}} \prod_v w_v(\gamma_v x_\alpha).
\]
Now each $\gamma$ corresponds to a term of the left-hand side of \eqref{eq:x_orbit} under the map
\begin{align*}
  \G(\OO_K)\backslash \Gamma(K) &\to \G(\OO_K) \backslash V(K) \\
  \G(\OO_K)\gamma &\mapsto \G(\OO_K) \gamma x_\alpha.
\end{align*}
The fiber of each $\G(\OO_K)x$ has size
\[
  [\Stab_{\Gamma(K)} x : \Stab_{\G(\OO_K)} x] = \frac{\size{H^0(K,M)}}{\size{\Stab_{\G(\OO_K)} x}}.
\]
So we match up one term of the left-hand side, having value
\[
  \prod_v w_v(x)/\size{\Stab_{\G(\OO_K)} x},
\]
with $\size{H^0(K,M)}/\size{\Stab_{\G(\OO_K)} x}$-many elements on the right-hand side. In view of the outlying factor $1/\size{H^0(K,M)}$, this completes the proof.
\end{proof}

\subsection{Fourier analysis of the local and global Tate pairings}
\label{sec:Fourier}

We now introduce the main innovative technique of this thesis: Fourier analysis of local and global Tate duality. In structure we are indebted to Tate's celebrated thesis \cite{Tate_thesis}, in which he
\begin{enumerate}
  \item\label{it:Tate_local} constructs a perfect pairing on the additive group of a local field $K$, taking values in the unit circle $\CC^{N=1}$, and thus furnishing a notion of Fourier transform for $\CC$-valued $L^1$ functions on $K$;
  \item\label{it:Tate_product} derives thereby a pairing and Fourier transform on the adele group $\AA_K$ of a global field $K$;
  \item\label{it:Tate_Poisson} proves that the discrete subgroup $K \subseteq \AA_K$ is a self-dual lattice and that the Poisson summation formula 
  \begin{equation}
    \sum_{x \in K} f(x) = \sum_{x \in K} \hat f(x)
  \end{equation}
  holds for all $f$ satisfying reasonable integrability conditions.
\end{enumerate}

In this paper, we work not with the additive group $K$ but with a Galois cohomology group $H^1(K, M)$. The needed theoretical result is \emph{Poitou-Tate duality}, a nine-term exact sequence of which the middle three terms are of main interest to us:
\[
  (\text{finite kernel}) \to H^1(K, M) \to \sideset{}{'}{\bigoplus}_v H^1(K_v, M) \to H^1(K, M')^\vee \to (\text{finite cokernel}).
\]
This can be interpreted as saying that $H^1(K, M)$ and $H^1(K, M')$ (where $M' = \Hom(M, \mu)$ is the Tate dual) map to dual lattices in the respective adelic cohomology groups
\[
  H^1(\AA_K, M) = \sideset{}{'}{\bigoplus}_v H^1(K_v, M) \textand 
  H^1(\AA_K, M') = \sideset{}{'}{\bigoplus}_v H^1(K_v, M'),
\]
which are mutually dual under the product of the local Tate pairings
\[
  \<\{\alpha_v\}, \{\beta_v\}\> = \prod_v \<\alpha_v, \beta_v\>
  \in \mu.
\]
Here, for $K$ a local field, the local Tate pairing is given by the cup product
\[
  \<\bullet, \bullet\> : H^1(K, M) \cross H^1(K, M') \to H^2(K, \mu) \isom \mu.
\]
It is well known that this pairing is perfect. (The Brauer group $H^2(K, \mu)$ is usually described as being $\QQ/\ZZ$ but, having no need for a Galois action on it, we identify it with $\mu$ to avoid the need to write an exponential in the Fourier transform.) Now, for any sufficiently nice function $f : H^1(\AA_K, M) \to \CC$ (locally constant and compactly supported is more than enough), we have Poisson summation
\[
  \sum_{\alpha \in H^1(K, M)} f(\alpha) = c_M \sum_{\beta \in H^1(K, M')} \hat f(\beta)
\]
for some constant $c_M$ which we think of as the covolume of $H^1(K,M)$ as a lattice in the adelic cohomology. (In fact, by examining the preceding term in the Poitou-Tate sequence, $H^1(K,M)$ need not inject into $H^1(\AA_K, M)$, but maps in with finite kernel; but this subtlety can be absorbed into the constant $c_M$.)

We apply Poisson summation to the local orbit counters $g_v$ defined in the preceding subsection and get a very general reflection theorem.
\begin{defn} \label{defn:dual}
Let $K$ be a local field. Let $(V^{(1)}, \Gamma^{(1)})$ and $(V^{(2)}, \Gamma^{(2)})$ be a pair of composed varieties over $K$ whose associated point stabilizers $M^{(1)}$, $M^{(2)}$ are Tate duals of one another, and let $(\V^{(i)}, \G^{(i)})$ be an integral model of $(V^{(i)}, \Gamma^{(i)})$. Two weightings on orbits
\[
  w^{(i)} : \G(\OO_K) \backslash \V(\OO_K) \to \CC
\]
are called (mutually) \emph{dual} with \emph{duality constant} $c \in \QQ$ if their local orbit counters $g_{w^{(i)}}$ are mutual Fourier transforms:
\begin{equation}\label{eq:general_local_refl}
  g^{(2)} = c \cdot \hat{g}^{(1)}.
\end{equation}
where the Fourier transform is scaled by
\[
  \hat f(\beta) = \frac{1}{H^0(K,M)} \sum_{\alpha \in H^1(K, M)} f(\alpha).
\]
An equation of the form \eqref{eq:general_local_refl} is called a \emph{local reflection theorem.} If the constant weightings $w^{(i)} = 1$ are mutually dual, we say that the two integral models $(\V^{(i)}, \G^{(i)})$ are \emph{naturally dual.}
\end{defn}

\begin{thm}[\textbf{local-to-global reflection engine}]\label{thm:main_compose}
Let $K$ be a number field. Let $(V^{(1)}, \Gamma^{(1)})$ and $(V^{(2)}, \Gamma^{(2)})$ be a pair of composed varieties over $K$ whose associated point stabilizers $M^{(1)}$, $M^{(2)}$ are Tate duals of one another. Let $(\V^{(i)}, \G^{(i)})$ be an integral model for each $(V^{(i)}, \Gamma^{(i)})$, and let
\[
  w_v^{(i)} : \G^{(i)}(\OO_v) \backslash \V^{(i)}(\OO_v) \to \CC
\]
be a local weighting on each integral model. Suppose that each integral model and local weighting satisfies the hypotheses of Lemma \ref{lem:loczn_main}, and suppose that at each place $v$, the two integral models are dual with some duality constant $c_v \in \QQ$. Then the weighted global class numbers are in a simple ratio:
\[
  h_{\left\{w_v^{(2)}\right\}} = \prod_v c_v \cdot h_{\left\{w_v^{(1)}\right\}}.
\]
\end{thm}
\begin{proof}
By Lemma \ref{lem:loczn_main},
\[
  h_{\left\{w_v^{(i)}\right\}} = \frac{1}{\size{H^0(K, M^{(i)})}} \sum_{\alpha \in H^1(K, M^{(i)})} \prod_v g_{v,w_v^{(i)}}(\alpha).
\]
At almost all $v$, each $g_{v,w_v^{(i)}}$ is supported on the unramified cohomology, and must be constant there because otherwise its Fourier transform would not be supported on the unramified cohomology. However, $g_{v,w_v^{(i)}}$ cannot be identically $0$ because of the existance of a global basepoint. So for such $v$,
\[
  g_{v,w_v^{(i)}} = \1_{H^1_\ur(K, M^{(i)})} \textand c_v = 1.
\]
In particular, the product $\prod_v g_{v,w_v^{(i)}}$ is a locally constant, compactly supported function on $H^1(K, \AA_{K})$, which is more than enough for Poisson summation to be valid.

Since the pairing between the adelic cohomology groups $H^1(\AA_K, M^{(i)})$ is made by multiplying the local Tate pairings, a product of local factors has a Fourier transform with a corresponding product expansion:
\[
\widehat{\prod_v g_{v,w_v^{(1)}}} = \prod_v \hat g_{v,w_v^{(1)}} = \prod_v c_v \cdot \prod_v g_{v,w_v^{(2)}}.
\]
We then apply Poisson summation to get a formula for the ratio of the global weighted class numbers:
\[
h_{\left\{w_v^{(2)}\right\}} = \frac{\size{H^0(K, M^{(1)})}}{\size{H^0(K, M^{(2)})}} \cdot c_{M^{(1)}} \prod_v c_v \cdot h_{\left\{w_v^{(1)}\right\}}.
\]
This gives the desired identity, except for determining the scalar $c_M$, which depends only on the Galois module $M = M^{(1)}$. This can be ascertained by applying Poisson summation to just one function $f : H^1(\AA_K, M) \to \CC$ for which either side is nonzero. The easiest such $f$ to think of is the characteristic function of a compact open box
\[
  X = \prod_v X_v,
\]
with $X_v = H^1(K, M)$ for almost all $v$. Such a specification is often called a \emph{Selmer system,} and the sum
\[
  \sum_{\alpha \in H^1(K,M)} \1_{\alpha_v \in X_v \forall v}
\]
is the order of the \emph{Selmer group} $\Sel(X)$ of global cohomology classes obeying the specified local conditions. Poisson summation becomes a formula for the ratio $\size{\Sel(X)}/\size{\Sel(X^\perp)}$ as a product of local factors, commonly known as the \emph{Greenberg-Wiles formula.} By appealing to any of the known proofs of the Greenberg-Wiles formula (see Darmon, Diamond, and Taylor \cite[Theorem 2.19]{DDT} or Jorza \cite[Theorem 3.11]{Jorza}), we pin down the value
\[
  c_M = \frac{\size{H^0(K, M')}}{\size{H^0(K, M)}}. \qedhere
\]
\end{proof}

At certain points in this paper, it will be to our advantage to consider multiple integral models at once. The following theorem has sufficient generality.
\begin{thm}[\textbf{local-to-global reflection engine: general version}]\label{thm:main_compose_multi}
Let $K$ be a number field. Let $(V^{(1)}, \Gamma^{(1)})$ and $(V^{(2)}, \Gamma^{(2)})$ be a pair of composed varieties over $K$ whose associated point stabilizers $M^{(1)}$, $M^{(2)}$ are Tate duals of one another. For each place $v$ of $K$, let
\[
  \left\{(\V^{(i)}_{j_v}, \G^{(i)})_{j_v} : j_v \in J_v^{(i)}\right\}
\]
be a family of integral models for each $(V^{(i)}$ indexed by some finite set $J_v^{(i)}$, and let
\[
w^{(i)}_{j_v} : \G^{(i)}(\OO_v) \backslash \V^{(i)}(\OO_v) \to \CC
\]
be a weighting on the orbits of each integral model. Similarly to Lemma \ref{lem:loczn_main} and Theorem \ref{thm:main_compose}, assume that
\begin{enumerate}[$($i$)$]
  \item $(V, \Gamma)$ is Hasse.
  \item For each combination of indices $j = (j_v)_v$, $j_v \in J_v^{(i)}$, the local integral models $(\V^{(i)}_{j_v}, \G^{(i)}_{j_v})$ glue together to form a global integral model $(\V^{(i)}_{j}, \G^{(i)}_{j})$. (Since the integral models are equipped with embeddings $\V^{(i)}_{j_v} \to V^{(i)}_v$, the gluing is seen to be unique; and its existence will be obvious in all the examples we consider.)
  \item Each such $\G^{(i)}_{j}$ has class number one.
  \item For each $j_v$, there are only finitely many orbits of $\G^{(i)}_{j_v}$ on $\V^{(i)}_{j_v}$, ensuring that the local orbit counter $g_{j_v, w_{j_v}}$ takes finite values.
  \item For almost every $v$, the index set $J_v = \{j_v\}$ has just one element, with the corresponding integral model $\V^{(i)}_{j_v}$ consisting of at most one orbit in each $\Gamma(K_v)$-orbit, and $w^{(i)}_{j_v} = 1$ identically.
  \item At each $v$, we have a local reflection theorem
  \[
    \sum_{j_v \in J_v^{(1)}} \hat g_{j_v, w_{j_v}} =
    \sum_{j_v \in J_v^{(2)}} g_{j_v, w_{j_v}}.
  \]
\end{enumerate}
Then the class numbers of the global integral models $(\V_j^{(i)}, \G_j^{(i)})$ with respect to the weightings $w_j^{(i)} = \prod_v w_{j_v}^{(i)}$ satisfy global reflection:
\[
  \sum_{j \in \prod_v J_v^{(1)}} h\(\V_j^{(1)}, w_j^{(1)}\) =
  \sum_{j \in \prod_v J_v^{(2)}} h\(\V_j^{(2)}, w_j^{(2)}\).
\]
\end{thm}
\begin{proof}
Except for complexities of notation, the proof closely follows the preceding one. The first five hypotheses ensure that each global integral model $\(\V^{(i)}_{j}, \G^{(i)}_{j}\)$ satisfies the hypotheses of Lemma \ref{lem:loczn_main}, so its class number is representable as a sum over the lattice of global points in adelic cohomology:
\[
  h\(\V_j^{(i)}, w_j^{(i)}\) = \frac{1}{\size{H^0(K, M^{(i)})}} \sum_{\alpha \in H^1(K, M^{(i)})} \prod_v g_{j_v,w_{j_v}}(\alpha).
\]
When we sum over all $j$, the contributions of each $\alpha$ factor to give
\[
  \sum_{j \in \prod_v J_v^{(i)}} h\(\V_j^{(i)}, w_j^{(i)}\) = \frac{1}{\size{H^0(K, M^{(i)})}} \sum_{\alpha \in H^1(K, M^{(i)})}\prod_v \sum_{j_v \in J_v^{(i)}} g_{j_v,w_{j_v}}(\alpha).
\]
But by the assumed local reflection identity, we have
\[
  \(\prod_v \sum_{j_v \in J_v^{(1)}} g_{j_v,w_{j_v}}\)^{\ds\widehat{}}
  =
  \prod_v \(\sum_{j_v \in J_v^{(1)}} g_{j_v,w_{j_v}}\)^{\ds\widehat{}}
  =
  \prod_v \sum_{j_v \in J_v^{(2)}} g_{j_v,w_{j_v}}.
\]
So we get the desired identity from Poisson summation. The scale factor was determined in proving the previous theorem.
\end{proof}
\begin{rem}
Unlike in the previous theorem, we have not included duality constants $c_v$, but the same effect can be obtained by taking the appropriate constant for the weighting $w_j^{(i)}$.
\end{rem}

\subsubsection{Examples}

As one might guess, there are many pairs of composed varieties whose point stabilizers $M^{(1)}$, $M^{(2)}$ are Tate duals; and, given any integral models, it is usually possible to concoct weights $w^{(i)}$ that are mutually dual, thereby getting reflection theorems from Theorem \ref{thm:main_compose}. More noteworthy is when a pair of integral models are \emph{naturally} dual at all finite places. Even more significant is if a group $\G$ acts on a large variety $\Lambda$, leaving certain functions $I$ on $\Lambda$ invariant, such that \emph{every} level set of $I$ is an integral model for a composed variety with natural duality. This is the case for O-N.

We have found three families of naturally dual composed varieties of this sort:
\[
\begin{tabular}{ccccc}
  $ \Gamma $ & $ \Lambda $ & $ I $ & Parametrizes & $ M $ \\ \hline
  $ \displaystyle \left\{\left[\begin{array}{cc}
    \lambda & t            \\
    0       & \lambda^{-2}
  \end{array}\right]\right\} \subset \GL_2 $ & \begin{tabular}{@{}c@{}} Quadratic forms, \\ $ \Sym^2(2) $ \end{tabular} & $ a(b^2 - 4ac) $ & ? & $ \C_2 $ \\
  $ \SL_2 $ & \begin{tabular}{@{}c@{}} Cubic forms, \\ $ \Sym^3(2) $ \end{tabular} & Discriminant & \begin{tabular}{@{}c@{}}Cubic rings / \\ $ 3 $-torsion in \\ quadratic rings\end{tabular} &  $ \C_3 $\vspace{1ex} \\
  $ \SL_3 $ & \begin{tabular}{@{}c@{}} Pairs of ternary \\ quadratic forms, \\ $ \Sym^2(3)^{\oplus 2} $ \end{tabular} & Cubic resolvent & \begin{tabular}{@{}c@{}}Quartic rings / \\ $ 2 $-torsion in \\ cubic rings \end{tabular} & $ \C_2 \cross \C_2 $
\end{tabular}
\]
These three representations will be considered in detail in Section \ref{sec:quadratic}, Part \ref{part:cubic}, and Parts \ref{part:quartic}--\ref{part:quartic_count}, respectively. In each case, there is a local reflection that pairs two integral models of $ V $ over $ \OO_K $ which look alike over $ K $.

\begin{rem}
In the latter two cases, the integral models are \emph{dual} under an identification of $ V $ with its dual $ V^* $ (which are isomorphic, up to an outer automorphism of $ \Gamma = \SL_3 $ in the last case). But in the quadratic case, $ V^* $ decomposes into $\bar K$-orbits according to a different invariant $J = a/\Delta^2$, and the integral orbit counts are infinite, so the alignment with duals in the classical sense must be considered at least partly coincidental.
\end{rem}

Closely related to the quartic rings example is the action of $ \SL_2 $ on binary quartic forms $ \Sym^4 (2) $. Here, the orbits are parametrized by a subset of a cohomology group $ H^1(K, M) $ ($ M \isom \C_2 \cross \C_2 $ as a group) cut out by a quadratic relation. Nevertheless, we will state some interesting reflection identities for these spaces in Section \ref{sec:bq}.

More generally, we can consider the space $\Lambda$ of pairs $(A,B)$ of $n$-ary quadratic forms, on which $\Gamma = \SL_n$ acts preserving a binary $n$-ic resolvent
\[
  I = \det(Ax - By).
\]
Although we do not consider it in this paper, preliminary investigations suggest that its integral models are naturally dual to one another for $n$ odd, yielding a corresponding global reflection theorem (Conjecture \ref{conj:O-N_2xnxn_Z}). 
This composed variety figures prominently in the study of Selmer elements of hyperelliptic curves \cite{bhargava2013hyperelliptic}.

On the other hand, the following families of composed varieties do \emph{not} admit natural duality:
\begin{itemize}
  \item The action of $\GG_m$ on the punctured affine line by multiplication by $n$th powers. The orbits do parametrize $H^1(K, \mu_m) = K^\cross / (K^\cross)^{n}$. But over a local or global field, there are infinitely many integral orbits in each rational orbit.
  \item The action of $\SO_2$ (the group of rotations preserving the quadratic form $x^2 + xy + y^2$) on binary cubic forms of the shape
  \[
    f(x,y) = ax^3 + bx^2 y + (-3a + b)x y^2 + a y^3
  \]
  which is symmetric under the threefold shift $x \mapsto y$, $y \mapsto -x-y$. This representation is used by Bhargava and Shnidman \cite{BhSh} to parametrize \emph{cyclic} cubic rings, that is, those with an automorphism of order $3$. The reason for failure of natural duality is quite simple. Within the representation over $\ZZ_p$ for $p \equiv 1 \mod 3$, take the composed variety where the discriminant is $p^2$. The cohomology group $H^1(\QQ_p, \ZZ/3\ZZ)$ is isomorphic to $\ZZ/3\ZZ \cross \ZZ/3\ZZ$, and a function $f$ on it may be written as a matrix
  \[
    \begin{matrix}
      \multicolumn{1}{c|}{f(0)} & f(\alpha) & f(2\alpha) \\ \hline
      f(\beta) & f(\alpha + \beta) & f(2\alpha + \beta) \\
      f(2\beta) & f(\alpha + 2\beta) & f(2\alpha + 2\beta)
    \end{matrix}
  \]
  in which the zero-element and the unramified cohomology $\<\alpha\>$ are marked off by dividers.
  
  The six \emph{ramified} cohomology elements each have one integral orbit, corresponding to the maximal order; the three \emph{unramified} cohomology elements---the zero element for $\QQ_p^3$, and the other two for the degree-$3$ unramified field extension in its two orientations---all have \emph{no} integral orbits, because the three orders
  \[
    \{(x_1, x_2, x_3) \in \ZZ_p^3 : x_i \equiv x_j \mod p\}
  \]
  are all asymmetric under the threefold automorphism of $\ZZ_p^3$. So we get a local orbit counter
  \[
    \begin{matrix}
      \multicolumn{1}{c|}{0} & 0 & 0 \\ \hline
      1 & 1 & 1 \\
      1 & 1 & 1
    \end{matrix}
  \]
  whose Fourier transform
  \[
   \begin{matrix}
     \multicolumn{1}{c|}{2} & -1 & -1 \\ \hline
     0 & 0 & 0 \\
     0 & 0 & 0
   \end{matrix}
  \]
  has mixed signs and thus cannot be the local orbit counter of any composed variety. Similar obstructions to natural duality have obtained in many of the composed varieties parametrizing rings with automorphisms found by Gundlach \cite{Gundlach}.
\end{itemize}

\part{Reflection theorems: first examples}
\label{part:first}

The remainder of this thesis will be devoted to stating and proving explicit reflection theorems for various objects of interest.

\section{Quadratic forms by superdiscriminant}\label{sec:quadratic}

We begin with the simplest Galois module $M \isom \ZZ/2\ZZ$.

There are many full composed varieties whose point stabilizer is of order $2$, and the one we take is, to say the least, one of the more unexpected. The group $\GL_2$ acts on the space
\[
  V = \Sym^2(2) = \{ ax^2 + bxy + cy^2 : a,b,c \in \GG_a \}
\]
of binary quadratic forms in the natural way. Let $\Gamma$ be the algebraic subgroup, defined over $\ZZ$, of elements of a peculiar form:
\[
  \left\{
  \begin{bmatrix}
    u & t \\
    0 & u^{-2}
  \end{bmatrix}
  : u \in \GG_m, t \in \GG_a \right\}.
\]
Abstractly, this group is a certain semidirect product of $\GG_a$ by $\GG_m$. As is not too hard to verify, the restriction of $\Lambda$ to $\Gamma$ has a single polynomial invariant, the \emph{superdiscriminant}
\[
  I := aD = a(b^2 - 4ac).
\]
Because of the asymmetry between $x$ and $y$, there is no harm in writing forms in $V$ inhomogeneously as $f(x) = ax^2 + bx + c$, as was done in Section \ref{sec:layman's_appendix}.

Then the variety
\[
  V(I) = \{f \in V : I(f) = I\}
\]
is full composed. We take the basepoint
\[
  f_0 = I x^2 + \frac{1}{4I}
\]
of discriminant $1$. Then the rational orbits are parametrized by $D = b^2 - 4ac \in K^\cross/(K^\cross)^2$ consistent with the parametrization of their splitting fields via Kummer theory.

\begin{rem}
The group $\Gamma$ is not reductive, that is, does fit into the classical Dynkin-diagram parametrization for Lie groups. Non-reductive groups are decidedly in the minority within the whole context of using orbits to parametrize arithmetic objects, but they have occurred before: Altu\u{g}, Shankar, Varma, and Wilson \cite{ASVW} count $D_4$-fields using orbits of pairs of ternary quadratic forms under a certain nonreductive subgroup of $\GL_2 \cross \SL_3$.
\end{rem}

Now we introduce integral models. Suppose $\OO_K \subseteq K$ is a PID with field of fractions $K$. If $\tau \in \OO_K$ divides $2$, then
\[
  V_\tau = \{ a x^2 + b x + c : a,c \in \OO_K, \tau \mid b \}
\]
is a $\GL_2(\OO_K)$-invariant lattice in $V$. For any $I \in \OO_K$, we can take $(\V, \G) = (V_\tau(I), \Gamma(\OO_K))$ as an integral model for $V(I)$. For it to have any integral points, we must have $\tau^2 \mid I$.

Our first local reflection theorem says that each of these integral models has a natural dual.

\begin{thm}[\textbf{``Local Quadratic O-N''}]\label{thm:O-N_quad_local}
Let $K$ be a non-archimedean local field, $\ch K \neq 2$. For $I, \tau \in \OO_K$ elements dividing $2$, the integral models
\[
  V_{\tau}(I) \textand V_{2\tau^{-1}} \(4\tau^{-4} I\)
\]
are naturally dual with scale factor the absolute norm $N(\tau) = \size{\OO_K / \tau\OO_K}$. In other words, the local orbit counters are related by
\begin{equation} \label{eq:O-N_quad_local}
  \hat g_{V_{\tau}(I)} = N(\tau) \cdot g_{V_{2\tau^{-1}} \(4\tau^{-4} I\)}.
\end{equation}
\end{thm}

\begin{proof}
We prove this result by explicitly computing the local orbit counter $g_{V_{\tau(I)}}$, which sends each $[D] \in K^\cross/(K^\cross)^2 \isom H^1(K, \ZZ/2\ZZ)$ to the number of cosets $[\gamma] \in \Gamma(\OO_K)\bs \Gamma(K)$ such that $\gamma v_0 \in V_\tau(I)(\OO_K)$, where $v_0$ is an arbitrary vector in $V(K)$ with $I(v_0) = I$ and $D(v_0) = D$. Let $t = v(\tau)$ and $e = v(2)$; we have
$e > 0$ exactly when $K$ is $2$-adic, and
\[
  0 \leq t \leq e.
\]

A coset $[\gamma]$ is specified by two pieces of information. First is the valuation $v(u)$ of the diagonal elements; this is equivalent to specifying $v(D)$ and $v(a)$, where, as is natural we set
\[
  \gamma v_0 = a x^2 + b x y + c y^2 \textand D = b^2 - 4 a c.
\]
Second, we specify $t$ modulo the appropriate integral sublattice. If (as we may assume) $v(u) = 0$, then $t$ is defined modulo $1$, which is the same as specifying $b$ modulo $2a$. So the problem of computing $g_{V_{\tau}(I)}$ devolves onto computing how many $b \in \tau\OO_K$, up to translation by $2a$, yield an integral value for
\[
  c = \frac{b^2 - D}{4a};
\]
that is, we must solve the quadratic congruence
\begin{equation} \label{eq:qfc}
  b^2 \equiv D \mod 4a.
\end{equation}
The answer, in general, depends on how close $D$ is to being a square in $K$. So we will express our answer in terms of the \emph{level spaces} introduced in Theorem \ref{thm:levels}. Here the level of a coclass $[\alpha], \alpha \in K^\cross$, is defined in terms of the discriminant of $K[\sqrt{\alpha}]$, which, by Theorem \ref{thm:disc_Kummer}, can be computed from the minimal distance $\size{\alpha - 1}$, over all rescalings of $\alpha$ by squares. The level spaces thus correspond to the natural filtration of $K^\cross/(K^\cross)^2$ by neighborhoods of $1$:
\[
  \L_i = \begin{cases}
    K^\cross/(K^\cross)^2, & i = -1 \\
    \{[\alpha] \in \OO_K^\cross/(\OO_K^\cross)^2 : \alpha \in 1 + \pi^{2e-2i} \OO_K\}, & 0 \leq i \leq e \\
    \{1\}, & i = e+1.
  \end{cases}
\]
Let $L_i$ be the characteristic function of $\L_i$. By Corollary \ref{cor:levels}, the Fourier transform of each $L_i$ is a scalar multiple of $L_{e - i}$.

We now claim that, if we fix $v(I)$ and $v(a)$ (and hence $v(D)$), then the contribution of all solutions of \eqref{eq:qfc} to $g_{\tau,I}$ can be expressed as a linear combination of the $L_i$. The basic idea, which will be a recurring one, is to group the solutions into \emph{families} that have a constant number of solutions over some subset $S \subseteq K^\cross/(K^\cross)^2$. The subset $S$ will be called the \emph{support} of the family, and the number of solutions for each $D \in S$ will be called the \emph{thickness} of the family.

If $v(D) \geq v(4a)$, then \eqref{eq:qfc} simplifies to $4a|b^2$, that is,
\[
  v(b) \geq \ceil{e + \frac{1}{2}v(a)}.
\]
Since we are counting values of $b$ modulo $2a$, the number of solutions is simply
\[
  q^{\( e + v(a)\) - \ceil{e - \frac{1}{2} v(a)}} = q^{\floor{\frac{1}{2}v(a)}}.
\]
We get a family with this thickness, supported on either $\L_0$ or $\L_{-1} \setminus \L_{0}$ according as $v(D)$ is even or odd.

If $v(D) < v(4a)$, then $b^2$ must be actually able to cancel at least the leading term of $D$ to get any solutions. In particular, $v(D)$ must be even. Let $\tilde{D} = D / \pi^{v(D)}$, and let $\tilde{b} = b/\pi^{\frac{1}{2} v(D)}$, so $\tilde{b}$ must be a unit satisfying
\begin{equation} \label{eq:qfc2}
  \tilde b^2 \equiv \tilde D \mod \frac{4a}{D}.
\end{equation}
Let $m = v(4a/D)$. If $m \geq 2e+1$, a unit is a square modulo $\pi^m$ only if it is a square outright, so we get a family supported just on the trivial class $1 \in K^\cross/(K^\cross)^2$. Otherwise, we have $1 \leq m \leq 2e$, and the support is $L_{\ceil{m/2}}$. The corresponding thicknesses are easy to compute. The $\tilde b$ satisfying \eqref{eq:qfc2} form a fiber of the group homomorphism
\[
  \phi = \bullet^2 \colon \( \OO_K / \pi^{v(2a) - \frac{1}{2}v(D)} \)^\cross
  \to \( \OO_K / \pi^{v(4a) - v(D)} \)^\cross,
\]
and the cokernel of this homomorphism has size $[\L_0 : \L_i]$, so the thickness is
\begin{align*}
  \size{\ker \phi} &= [\L_0 : \L_i] \cdot \frac{\Size{\(\OO_K / \pi^{v(2a) - \frac{1}{2}v(D)} \)^\cross}}{\Size{\(\OO_K / \pi^{v(4a) - v(D)} \)^\cross}} \\
    &= [\L_0 : \L_i] \cdot \frac{\( 1 - \frac{1}{q} \) q^{v(2a) - \frac{1}{2}v(D)}}
  {\( 1 - \frac{1}{q} \) q^{v(4a) - v(D)}} \\
  &= [\L_0 : \L_i] \cdot q^{\frac{1}{2} v(D) - e} \\
  &= \begin{cases}
    q^{\frac{1}{2} v(D) - e + \ceil{m/2}} = q^{\floor{v(a)/2}}, & 1 \leq m \leq 2e \\
    2q^{\frac{1}{2} v(D)}, & m \geq 2e + 1.
  \end{cases}
\end{align*}

We have not mentioned the condition $b \in (\tau)$, because it is equivalent to $v(D) \geq 2 v(\tau)$, and eliminates some families, leaving the others intact.

By way of illustration, we tabulate the contributions to $g_{\tau,I}$ in the example where $e = 2$. It is already easy to check many examples of Theorem \ref{thm:O-N_quad_local}.
\[
\newcommand{\rline}[1]{\multicolumn{1}{c|}{#1}}
\begin{tabular}{r|ccccccccc}
  $\downarrow v(a); v(D) \rightarrow$ & $0$ & $1$ & $2$ & $3$ & $4$ & $5$ & $6$ & $7$ & $8$ \\ \hline
  $0$ & ${L_2}$ & & ${L_1}$ & \rline{} & ${L_0}$ & ${L_{-1}} - {L_0}$ & ${L_0}$ & ${L_{-1}} - {L_0}$ & ${L_0}$ \\ 
  \cline{2-2}\cline{6-6}
  $1$ & \rline{${L_3}$} & & ${L_1}$ & & \rline{${L_0}$} & ${L_{-1}} - {L_0}$ & ${L_0}$ & ${L_{-1}} - {L_0}$ & ${L_0}$ \\
  \cline{3-3}\cline{7-7}
  $2$ & ${L_3}$ & \rline{} & $q{L_2}$ & & $q{L_1}$ & \rline{} & $q{L_0}$ & $q\({L_{-1}} - {L_0}\)$ & $q{L_0}$ \\
  \cline{4-4}\cline{8-8}
  $3$ & ${L_3}$ & & \rline{$q{L_3}$} & & $q{L_1}$ & & \rline{$q{L_0}$} & $q\({L_{-1}} - {L_0}\)$ & $q{L_0}$ \\
  \cline{5-5}\cline{9-9}
  $4$ & ${L_3}$ & & $q{L_3}$ & \rline{} & $q^2{L_2}$ & & $q^2{L_1}$ & \rline{} & $q^2{L_0}$ \\
  \cline{6-6}\cline{10-10}
\end{tabular}
\]
We have shown the subdivision of the table into three \emph{zones} given by the inequalities:
\begin{itemize}
  \item Zone I: $v(D) \geq v(4a)$
  \item Zone II: $v(a) < v(D) \leq v(4a)$
  \item Zone III: $v(a) > v(D)$.
\end{itemize}
(A more general definition of a zone will be given later.) In general, the shapes of these zones, together with the needed condition $v(D) \geq 2t$, will look as follows:
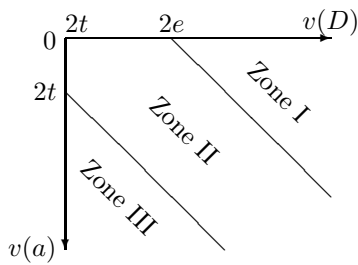
\[
\setlength{\unitlength}{1em}
\begin{picture}(13,9)(-2,-8.2)
  \put(0,0){\vector(1,0){10}}
  \put(0,0){\vector(0,-1){8}}
  \put(0,-2){\line(1,-1){6}}
  \put(4,0){\line(1,-1){6}}
  \put(-2,-9){\makebox(2,2)[r]{$v(a)$~}}
  \put(-2,-3){\makebox(2,2)[r]{$2t$~}}
  \put(-2,-1){\makebox(2,2)[r]{$0$\vphantom{$^2$}~}}
  \put( 0,0){\makebox(2,2)[bl]{$2t$\vphantom{|}}}
  \put( 3,0){\makebox(2,2)[b]{$2e$\vphantom{|}}}
  \put( 9,0){\makebox(2,2)[b]{$v(D)$\vphantom{\big|}}}
  \put(2,-1){\rotatebox{-45}{\makebox(7,0){Zone II}}}
  \put(6,0){\rotatebox{-45}{\makebox(5.7,0){Zone I}}}
  \put(0,-4){\rotatebox{-45}{\makebox(5.7,0){Zone III}}}
\end{picture}
\]
The feature to be noted is that, under the transformation $t \mapsto e - t$, the shape of Zone II is flipped about a diagonal line and Zones I and III are interchanged. This will be the basis for our proof of Theorem \ref{thm:O-N_quad_local}; but there will be irregularities owing to the floor functions in the formulas and the fact that $\L_{-1} \setminus \L_0$, instead of $\L_{-1}$, appears as a support.

There are two ways to finish the proof. One is to establish a bijection of families, as outlined in the previous paragraph, so that $L_i$ and $L_{e-i}$ are interchanged as supports and all the thicknesses correspond appropriately. Such an approach will be used for cubic O-N in Section \ref{sec:wild_bij}. The other is to verify the local reflection computationally, by means of a generating function. We choose the second, admittedly less elegant, method, mainly because it shows, in a context simple enough to be worked by hand, transformations that we will relegate to a computer in the succeeding sections.

Let
\[
  F(Z) = \sum_{n \geq 0} g_{\tau, \pi^n} Z^n,
\]
a formal power series whose coefficients are functions of $D \in K^\cross/(K^\cross)^2$. We write $F = F_{\mathrm{I}} + F_{\mathrm{II}} + F_{\mathrm{III}}$, where $F_{\mathfrak{X}}$ is the contribution coming from Zone $\mathfrak{X}$ in the preceding analysis.

Writing $i = v(a)$ and $d = v(D)$, we proceed to compute
\begin{align*}
  F_{\mathrm{I}} &= \sum_{i \geq 0} \sum_{d \geq 2e+i+1} \begin{cases}
    q^{\floor{i/2}} L_0 Z^{i+d}, & d \text{ even} \\
    q^{\floor{i/2}} (L_{-1} - L_0) Z^{i+d}, & d \text{ odd}
  \end{cases}\\
  &= Z^{2e} \sum_{i \geq 0} \sum_{d \geq i+1} \begin{cases}
    q^{\floor{i/2}} L_0 Z^{i+d}, & d \text{ even} \\
    q^{\floor{i/2}} (L_{-1} - L_0) Z^{i+d}, & d \text{ odd}.
  \end{cases} \\
\intertext{Splitting $i = 2i_f + i_p$, where $0 \leq i_p \leq 1$, and likewise $d = 2d_f + d_p$, we get}
  F_{\mathrm{I}} &= Z^{2e} \sum_{i_p = 0}^1 \sum_{i_f \geq 0} \left(
    \sum_{d \geq 2i_f + i_p + 1} q^{i_f} (-1)^d L_0 Z^{d + 2i_f + i_p} + \sum_{d_f \geq i_f} q^{i_f} L_{-1} Z^{2d_f + 2i_f + i_p + 1} \right) \\
  &= Z^{2e} \sum_{i_p = 0}^1 \sum_{i_f \geq 0} \left(\frac{q^{i_f}(-1)^{i_p} L_0 Z^{4i_f + 2i_p + 1}}{1 + Z} +
    \frac{q^{i_f} Z^{4i_f + i_p + 1}L_{-1}}{1 - Z^2} \right) \\
  &= Z^{2e} \sum_{i_p = 0}^1  \left( \frac{(-1)^{i_p} Z^{2i_p} L_0}{(1+Z)(1 - qZ^4)} + \frac{Z^{i_p + 1} L_{-1}}{(1 - Z^2)(1 - qZ^4)} \right) \\
  &= Z^{2e}\(\frac{(1-Z^2)}{(1+Z)(1 - qZ^4)} L_0 + \frac{Z(1 + Z)}{(1-Z^2)(1 - qZ^4)} L_{-1} \) \\
  &= \frac{Z^{2e}(1-Z)}{1 - qZ^4} L_0 + \frac{Z^{2e+1}}{(1-Z)(1 - qZ^4)} L_{-1}.
\end{align*}
For Zone II, which appears only when $e > 0$, the most sensible way to evaluate the sum
\[
  F_{\mathrm{II}} = \sum_{i\geq 0} \sum_{\substack{i \leq d < i + 2e \\ d \geq 2t \\ d\text{ even}}} q^{e + \floor{\frac{i}{2}} - \frac{d}{2}} L_{e + \floor{\frac{i}{2}} - \frac{d}{2}} Z^{i+d}
\]
is to group terms with the same level $L_j$. We have $j = e + \floor{\frac{i}{2}} - \frac{d}{2}$, so the values of $i$ and $j$ determine $d$. The condition $d \geq 2t$ reduces to $i \geq 2(j + t - e)$; the other condition $i \leq d < i + 2e$ is automatically satisfied if $1 \leq j \leq e-1$, while if $j = 0$ or $j = e$, we must have $i$ odd or $i$ even respectively. For $1 \leq j \leq e-1$, the $L_j$-piece of $F_{\mathrm{II}}$ is therefore
\begin{align*}
&\sum_{i \geq \max\{0, 2(j + t - e)\}} q^{\floor{i/2}} L_j Z^{i + 2\( e + \floor{\frac{i}{2}} - 2j\)} \\
&= \sum_{i_p = 0}^1 \sum_{i_f \geq \max\{0, j + t - e\}} q^{i_f} Z^{4i_f + i_p + 2e - 2j} L_j \\
&= \left( \sum_{i_p = 0}^1 Z^{i_p} \right) \left( \sum_{i_f \geq \max\{0, j + t - e\}} \( qZ^4\)^{i_f} \right) Z^{2e - 2j} L_j \\
&= \frac{(1+Z)\( qZ^4\)^{\max\{0,j+t-e\}}}{1 - qZ^4} \cdot Z^{2e - 2j}L_j.
\end{align*}
For $j = 0$ and $j = e$, since $i_p$ can only take one of its two values, the initial factor $1+Z$ is to be replaced by $Z$ and $1$ respectively. Finally, Zone III presents no particular difficulties:
\begin{align*}
  F_{\mathrm{III}} &= \sum_{\substack{d\geq 2t \\ d\text{ even}}} \sum_{i \geq d + 1} 2q^{d/2} L_{e+1} Z^{i+d} \\
  &= 2 \sum_{\substack{d\geq 2t \\ d\text{ even}}} q^{d/2} \cdot \frac{Z^{2d+1}}{1 - Z} \cdot L_{e+1} \\
  &= \frac{2q^sZ^{4t+1}}{(1 - Z)(1 - qZ^4)} L_{e+1}.
\end{align*}
Summing up, we get for $e \geq 1$ (the case $e = 0$ can be handled similarly)
\begin{align*}
  F &= F_{\mathrm{I}} + F_{\mathrm{II}} + F_{\mathrm{III}} \\
  &= \frac{Z^{2e+1}}{(1-Z)(1 - qZ^4)} L_{-1} + \frac{Z^{2e}(1-Z)}{1 - qZ^4} L_0
  + \frac{Z^{2e+1}}{1 - qZ^4} L_0 + \sum_{1 \leq j \leq e-1} \frac{(1+Z)\( qZ^4\)^{\max\{0,j+t-e\}}}{1 - qZ^4} \cdot Z^{2e - 2j}L_j \\ & \quad {} + \frac{\( qZ^4 \)^t}{1 - qZ^4} L_e + \frac{2q^sZ^{4t+1}}{(1 - Z)(1 - qZ^4)} L_{e+1} \\
  &= \frac{Z^{2e+1}}{(1-Z)(1 - qZ^4)} L_{-1} + \frac{Z^{2e}}{1 - qZ^4} L_0 + \sum_{1 \leq j \leq e-1} \frac{(1+Z)\( qZ^4\)^{\max\{0,j+t-e\}}}{1 - qZ^4} \cdot Z^{2e - 2j}L_j \\ & \quad {} + \frac{\( qZ^4 \)^t}{1 - qZ^4} L_e + \frac{2Z \( qZ^4 \)^t}{(1 - Z)(1 - qZ^4)} L_{e+1}.
\end{align*}
Now the evident symmetry between the coefficients of $L_j$ and $L_{e-j}$, when the transformation $t \mapsto e - t$ is made, establishes the theorem.
\end{proof}

Inserting this into the machinery of Part \ref{part:composed} produces global reflection theorems:
\begin{thm}[\textbf{``Quadratic O-N''}] \label{thm:O-N_quad}
Let $K$ be a number field of class number $1$. Then for any $I, \tau \in \OO_K$ with $\tau \mid 2$,
\[
  \sum_{\substack{f \in \Gamma\(\OO_K\)\bs V_\tau(I)(\OO_K) \\
     \disc f > 0 \text{ at every real place}}}
  \frac{1}{\size{\Stab_{\Gamma(\OO_K) f}}}
  = \frac{\size{N_{K/\QQ}(\tau)}}{2^{r_2(K)}}
  \sum_{f \in \Gamma\(\OO_K\)\bs V_{2\tau^{-1}}(4\tau^{-4}I)(\OO_K)} \frac{1}{\size{\Stab_{\Gamma(\OO_K) f}}}
\]
where $r_2(K)$ is the number of complex places of $K$.
\end{thm}
\begin{proof}
We verify the hypotheses of Lemma \ref{lem:loczn_main} on the integral models $V_\tau(I)$ and $V_{2\tau^{-1}}(4\tau^{-4} I)$:
\begin{enumerate}[$($i$)$]
  \item $V(I)$ is Hasse because it is full, as previously noted.
  \item To check that $\Gamma$ has class number $1$, it suffices to check the factors $\GG_m$ and $\GG_a$ of which $\Gamma$ is a semidirect product. The former of these has the same class number as $K$, explaining the restriction in the theorem statement.
  \item The finiteness of the local orbit counter follows from the formulas for it computed in the previous theorem.
  \item Finally, at almost all places, we plug in $e = t = 0$ to get $F = L_0$, establishing the needed convergence.
\end{enumerate}
Now we need the local reflection itself. We keep track of the constants $c_v$ accrued:
\begin{itemize}
  \item If $v \nmid 2\infty$, the integral models are naturally dual with constant $c_v = 1$.
  \item If $v \mid 2$, the integral models are naturally dual with constant $c_v = [\OO_v : \tau\OO_v]$. Multiplying over all $v \mid 2$ and using that $\tau \mid 2$ yields a factor
  \[
    \prod_{v|2} c_v = [\OO_{K} : \tau\OO_K] = \size{N_{K/\QQ} \tau}.
  \]
  \item If $v$ is real, the integral models are no longer naturally dual at $v$. We place the non-natural weighting
  \[
    w^{(2)} = \1_{0}
  \]
  that picks out $\alpha \in H^1(K, \ZZ/2\ZZ)$ that vanish at $v$, that is, forms with positive discriminant at $v$. This is the Fourier transform of $w^{(1)} = 1$, so $c_v = 1$.
  \item Finally, if $v$ is complex, then the integral models are certainly naturally dual at $v$, because $\size{H^1} = 1$. However, the scaling of the Fourier transform by $1/\size{H^0(\CC, M)} = 1/2$ requires that we take $c_v = 1/2$.
\end{itemize}
Multiplying these constants gives the constant claimed.
\end{proof}

\begin{rem}
The condition that $\OO_K$ be a PID can be dropped, but then $\Gamma$ no longer has class number $1$, and each side of the theorem becomes a sum of orbit counts on $\Cl(\OO_K)$-many global integral models that locally look alike. We do not spell out the details here. We wonder whether such a method works in general to circumvent the class-number-$1$ hypothesis in Theorem \ref{thm:main_compose}.
\end{rem}

We conclude by specializing further to the case $K = \QQ$. We replace $\Gamma(\ZZ)$ by its index-$2$ subgroup, the group $\ZZ$ of translations. This merely doubles all orbit counts, and it acts freely on quadratics with nonzero discriminant, so we can suppress all mention of stabilizers for the following charmingly simple statement, also featured in Section \ref{sec:layman's_appendix}:

\begin{thm}[\textbf{``Quadratic O-N''}] \label{thm:O-N_quad_Z}
If $n$ is a nonzero integer, let $q(n)$ be the number of integer quadratic polynomials $f(x) = ax^2 + bx + c$ with
\[
  a(b^2 - 4ac) = n,
\]
up to the trivial change $x \mapsto x + t$ ($t \in \ZZ$). Let $q_2(n)$, $q^+(n)$, and $q_2^+(n)$, respectively, be the number of these $f$ such that $2|b$ (for $q_2$), such that the roots of $f$ are real (for $q^+$), or which satisfy both conditions (for $q_2^+$). Then for all nonzero integers $n$,
\begin{align*}
  q_2^+(4n) &= q(n) \\
  q_2(4n) &= 2 q^+(n).
\end{align*}
\end{thm}
\begin{examp} \label{ex:QR}
Looking at $n = p_1p_3$, where $p_1 \equiv 1$ (mod $4$) and $p_3 \equiv 3$ (mod $4$) are primes, the counts involve certain Legendre symbols. For instance, the combination $a = p_3$, $b^2 - 4ac = p_1$ is feasible if and only if the congruence
\[
  b^2 \equiv p_1 \mod 4p_3
\]
has a solution, which happens exactly when $\( \frac{p_1}{p_3} \) = 1$. Working out all cases, we find that
\[
  q^+(p_1p_3) = 5 + \( \frac{p_1}{p_3} \) \textand
  q_2(4p_1p_3) = 10 + 2 \( \frac{p_3}{p_1} \).
\]
Thus our reflection theorem recovers the quadratic reciprocity law
\[
  \( \frac{p_1}{p_3} \) = \( \frac{p_3}{p_1} \).
\]
We wonder: does there exist a proof of Theorem \ref{thm:O-N_quad_Z} using no tools more advanced than quadratic reciprocity?
\end{examp}

\section{Class groups: generalizations of the Scholz and Leopoldt reflection theorems}
We now return to the consideration with which we began: reflection theorems for class groups. Scholz \cite{ScholzRefl} proved a relation between the $3$-torsion in the class groups of $\QQ(\sqrt{D})$ and $\QQ(\sqrt{-3D})$. Leopoldt \cite{Leopoldt} significantly generalized this result. We here present a generalization of Leopoldt's result to orders in $\GA(\FF_p)$-extensions, where exact formulas (as opposed to bounds) can often be obtained. We will not use composed varieties; instead, we will use Poisson summation in the form of the Greenberg-Wiles formula to get reflection theorems.

Let $T/K$ be a $(\ZZ/p\ZZ)^\cross$-torsor. If $\OO \subseteq T$ is a Galois-invariant $\OO_K$-order, then $(\ZZ/p\ZZ)^\cross$ acts on the class group $\Cl(\OO)$. The $p$-primary part $\Cl(\OO)_p$ is broken up into eigenspaces, one for each character $\chi : (\ZZ/p\ZZ)^\cross \to \mu_{p-1} \subseteq \ZZ_p^\cross$. There is a distinguished character $\chi = \omega$ lifting the reduction map modulo $p$ (the Teichm\"uller lift). We will concern ourselves with the $\omega$-part $\Cl(\OO)_{p,\omega}$. (The remaining parts are related to the $\omega$-parts of the class groups of other torsors.) We look at the $p$-torsion, or equivalently the $p$-cotorsion:
\[
  \Cl(\OO)[p]_\omega \cong (\Cl(\OO)/p\Cl(\OO))_\omega.
\]

\subsection{Dual orders}

We now develop a condition on two orders $\OO_1 \subseteq T$, $\OO_2 \subseteq T'$ that will suffice to produce a reflection theorem between their class groups. First, a simple lemma:
\begin{lem}\label{lem:ring_cl_fld}
Let $L$ be an \'etale algebra over a number field $K$, let $\OO \subseteq L$ be an order, and let $M/L$ be a $G$-torsor. The following conditions are equivalent:
\begin{enumerate}[$($a$)$]
\item $M$ is a ring class algebra for $\OO$; that is, the global Artin map
\[
  \psi_{M/L} = \prod_{i} \psi_{M_i/L_i} : \Idls(L, \mm) \to G
\]
factors through $\Cl(\OO)$, where $\mm \subseteq \OO_K$ is an admissible modulus for $M/L$, $\Idls(L, \mm)$ is the group of invertible fractional ideals of $L$ prime to $\mm$, and the product runs through all field factors $M_i$ of $M$, with $L_i$ being the corresponding field factor of $L$;
\item For every valuation $v$ of $K$, the local Artin map
\[
  \phi_{M_v/L_v} = \prod_{u|w|v} \phi_{M_u/L_w} : L_v^\cross \to G
\]
vanishes on $\OO_{v}^\cross$.
\end{enumerate}
\end{lem}
\begin{proof}
By local-global compatibility, the global Artin map can be described idelically as the product of the local ones. Indeed, $\Idls(L, \mm)$ embeds into $\AA_L^\cross/\prod_{v\nmid \mm} \OO_{L_v}^\cross$, and
\[
  \psi_{M/L} = \prod_{v} \phi_{M_v/L_v} : \AA_L^\cross \to G.
\]
Now the idele-theoretic description of $\Cl(\OO)$ is
\[
  \Cl(\OO) = \AA_L^\cross / \(L^\cross \cdot \prod_v\OO_v^\cross\).
\]
Since $\psi_{M/L}$ always vanishes on the principal ideles $L^\cross$, it factors through $\Cl(\OO)$ if and only if it vanishes on each $\OO_v^\cross$ ($v$ a place of $K$), where it reduces to the product $\phi_{M_v/L_v}$ of the local Artin maps at the primes dividing $v$, as desired.
\end{proof}
There is an analogue for \emph{narrow} ring class algebras: here $\phi_{M_v/L_v}$ is required to vanish on $\OO_v^\cross$ for $v$ finite only.

This motivates the following definitions.
\begin{defn}\
\begin{enumerate}[(a)]
\item Let $K$ be a local field, $T/K$ a $(\ZZ/p\ZZ)^\cross$-torsor and $T'$ its Tate dual. Two $(\ZZ/p\ZZ)^\cross$-invariant orders $\OO_1 \subseteq T$, $\OO_2 \subseteq T'$ are called \emph{dual} if the $\omega$-parts of the multiplicative groups, $(\OO_1^\cross)_\omega$ and $(\OO_2^\cross)_\omega$, are orthogonal complements under the Hilbert pairing, which as we know is perfect between $(T^\cross)_\omega$ and $(T'^\cross)_\omega$.
\item Let $K$ be a global field, $T/K$ a $(\ZZ/p\ZZ)^\cross$-torsor and $T'$ its Tate dual. Two orders $\OO_1 \subseteq T$, $\OO_2 \subseteq T'$ are called \emph{dual} if the completions $\OO_{1,\qq}$, $\OO_{2,\qq}$ are dual for all primes $\qq$ of $K$.
\end{enumerate}
\end{defn}
A dual pair yields a reflection theorem, as follows.

\begin{thm} \label{thm:Scholz_for_locally_dual_orders}
Let $\OO_1 \subseteq T$, $\OO_2 \subseteq T'$ be dual orders. Then
\begin{equation} \label{eq:Scholz}
  \frac{\size{\Cl^+(\OO_1)[p]_\omega}}{\size{\Cl(\OO_2)[p]_\omega}} =
  \frac{p^{\1_{T \text{ is totally split}}}}{p^{\1_{T' \text{ is totally split}}}}
  \cdot \prod_{v \text{ of } K} \frac{\size{\OO_{2,\ell,\omega}^\cross}}{p^{\1_{T_p \text{ is totally split}}}}
\end{equation}
\end{thm}
\begin{proof}[Proof of Theorem \ref{thm:Scholz_for_locally_dual_orders}]
The maps $\psi : \Cl(\OO_1) \to \ZZ/p\ZZ$ are the Artin maps of ring class algebras $E/T$ of $\OO_1$. Now $(\ZZ/p\ZZ)^\cross$ acts both on maps $\psi$ and algebras $E$, and it is easy to see that the $\psi$ belonging to the $\omega$-component correspond to $E$ that are \emph{symmetric,} that is, are $\GA(\FF_p)$-torsors with resolvent subtorsor $T$. 
So we get an injection of groups
\[
  i : \Hom(\Cl(\OO_1), \ZZ/p\ZZ) \to H^1(K, M_T).
\]
By Lemma \ref{lem:ring_cl_fld}, the image of $i$ is a Selmer group $\Sel_{X}(K, M_T)$, where the local conditions $X_v$ are given by
\begin{align*}
  X_v &= \text{the whole of } H^1(K_v, M_T) \text{ if } v | \infty \\
  X_v &= \OO_{1,v,\omega}^\perp = \OO_{2,v,\omega}, \quad v \text{ finite}
\end{align*}
where the second equality uses the duality of $\OO_1$ and $\OO_2$ and the Kummer parametrization
\[
  H^1(K_v, M_T) \cong T'^\cross_\omega.
\]
By the exact same argument, the dual Selmer system
\begin{align*}
  X_v &= 0 \text{ if } v | \infty \\
  X_v &= \OO_{2,v,\omega}^\perp = \OO_{1,v,\omega}, \quad v \text{ finite}
\end{align*}
has Selmer group naturally identified with $\Hom(\Cl^+(\OO_2), \ZZ/p\ZZ)$. (Of course, the distinction between wide and narrow class groups is only relevant if $p = 2$, a case which we will exclude in the next section.)

To finish, we apply the Greenberg-Wiles formula, as mentioned in the end of the proof of Theorem \ref{thm:main_compose}, and use that $\size{H^0(M_T)}$ is either $p$ or $1$ according as $T$ is totally split.
\end{proof}

\subsection{Dual orders are plentiful for quadratic extensions}
It's not hard to show that the maximal orders $\OO_T$, $\OO_{T'}$ are dual at primes $\ell \nmid p$. At $p$, however, it is not obvious how one might find a pair of dual orders, or whether such orders exist. However, there is a case in which this is manageable, and it specializes to the Scholz reflection principle in the case $p = 3$.

Let $p$ be an odd prime. We will assume that our base field $K$ contains the element
\[
  \rho_p = \zeta_p + \zeta_p^{-1} = 2 \cos \frac{2\pi}{p}.
\]
(Note that $\rho_3 = -1$, so this assumption always holds when $p = 3$.) This entails in particular that $K(\zeta_p)$ is an extension of $K$ of degree at most $2$, being $K(\sqrt{D})$ where
\[
  D = (\zeta_p - \zeta_p^{-1})^2 = \rho_p^2 - 4.
\]
We note that $D$ is a unit locally at all finite primes $\qq$ except those dividing $p$, in which case $\size{D}_\qq = \size{p}_\qq^{1/(p-1)} = d_{\min,\qq}^{1/p}$.

If $Q = K[\sqrt{a}]$ is an \'etale quadratic algebra, we may form the $(\ZZ/p\ZZ)^\cross$-torsor $T = Q^{(p-1)/2}$, with the unique possible torsor action. $Q$ is a $\mu_2$-torsor, and $Q^\cross/(Q^\cross)^p$ is the direct sum of two components: $Q^\cross_{\omega^0} \cong K^\cross/(K^\cross)^p$, and $Q^\cross_{\omega} = Q^{N=1}/(Q^{N=1})^p$ which parametrizes $\GA(M)$-extensions whose resolvent torsor is $T$. Due to the splitting of $T$, these are in fact $D_p$-extensions, where $D_p$ is the dihedral group (the permutation group that the symmetries of a regular $p$-gon induce on its vertices).

The Tate dual $T'$ is $Q'^{(p-1)/2}$, a product of copies of the quadratic algebra $Q' = K[\sqrt{D a}]$.

\subsection{Local dual generalized orders}

Suppose our base field $K$ has a distinguished subring of integers $\OO_K$, a Dedekind domain with field of fractions $K$. If $\OO \subseteq \OO_Q$ is an order over $\OO_K$, denote by $\R(\OO)$ the projection of $\OO^\cross$ onto $(Q^\cross)_{\omega}$, quotienting out by both $p$th powers and the eigenspace corresponding to the trivial character (namely $(\OO_K^\cross)/(\OO_K^\cross)^p$).

If $K$ is local, we call a pair of orders $\OO \subseteq Q$, $\OO' \subseteq Q'$ \emph{dual} if the associated \emph{unit class subgroups} $\R(\OO) \subseteq (Q^\cross)_{\omega}$, $\R(\OO') \subseteq (Q'^\cross)_{\omega}$ are orthogonal complements. For example, it is not hard to prove that if $\ch k_K \neq p$, the maximal orders in $Q$ and $Q'$ are dual to one another. We pose the question of whether any order in $Q$ admits a dual order. The answer is no, because $\R(\OO)$ can be as small as $\{1\}$ but cannot be as big as $(Q^\cross)_{\omega}$, being always contained in $(\OO_Q^\cross)_{\omega}$. This is essentially the only obstruction, and we remedy it by introducing a notion of \emph{generalized order.}

\begin{defn} If $Q$ is a quadratic \'etale algebra over a field $K$, in which a Dedekind domain $\OO_K$ is fixed as a ring of integers, a \emph{generalized order} in $Q$ is a finitely generated $\OO_K$-subalgebra $\OO$ spanning $Q$ over $K$ and closed under the conjugation automorphism of $Q$.
\end{defn}

If $\OO_K$ is local, then as soon as $\OO$ contains an element of $Q$ with negative valuation, even with respect to only one of the valuations on $Q$ (if $Q$ is split), then taking conjugates and powers shows that $\OO$ contains all elements of $Q$. Thus the only generalized orders in this case are that $\OO = Q$ or $\OO$ is an order in the ordinary sense, that is, a subring of $\OO_{Q}$ that spans $Q$. Letting $\OO_Q = \OO_K[\xi]$, these orders have the form $\OO = \OO_K[\pi^i\xi]$ for $i \geq 0$.

In general, a generalized order $\OO$ over a Dedekind domain $\OO_K$ is specified by a collection $(\OO_\qq)_\qq$ of orders in the completions $K_\qq$, almost all maximal; and particular has the form $\OO_1[\qq_1^{-1},\ldots,\qq_r^{-1}]$ where $\OO_1$ is an order in $\OO_Q$ and the $\qq_i$ are finitely many primes of $K$, at which $\OO_1$ can be taken maximal. Class groups of generalized orders over number fields are not hard to study: in the foregoing notation, we have that $\Cl(\OO) \cong \Cl(\OO_1)/\<\qq_1,\ldots,\qq_r\>$ is formed by quotienting out by the classes of the relevant primes.

\begin{lem} \label{lem:ord cond}
If $K$ is a local field and $\OO \subseteq \OO_Q$ is a quadratic generalized order, then $\R(\OO)$ is a level space in $Q^\cross_{\omega}$ (in the sense of Theorem \ref{thm:levels}). Moreover, all level spaces arise in this way.
\end{lem}
\begin{proof}
In the tame case that $\ch k_K \neq p$, there are at most three level spaces, and it is easy to identify the generalized orders to which they correspond:
\begin{align*}
  Q^\cross_{\omega} &= \R(Q) \\
  (\OO_Q)^\cross_{\omega} &= \R(\OO_Q) \\
  \{1\} &= \R(\OO), \text{ any } \OO \subsetneq \OO_Q.
\end{align*}
The last holds because any $\x \in \OO^\cross$ is the product of $x_0 \in K^\cross$, which maps into the $\omega^0$-component, and an $x_1 \equiv 1$ mod $\pi_K\OO_Q$ which is necessarily a $p$th power.

In the wild case we use similar methods. Since $p \neq 2$, we may write $\OO_Q = \QQ_K[\sqrt{\beta_Q}]$, where $v_K(\beta_Q)$ is $0$ or $1$. The Kummer element $\beta$ corresponding to the torsor $T = Q^{(p-1)/2}$ is $\beta = \beta_Q^{(p-1)/2}$.

The generalized order $Q$ has unit class subgroup
\[
  \R(Q) = \L_{-1}.
\]
The remaining orders can be described as
\[
  \OO_j = \OO_K\left[\pi^{j - v(\beta)/2}\sqrt{\beta}\right],
\]
where $j$, the valuation of a generator, ranges over the nonnegative elements of $\ZZ$ (if $\beta \sim 1$) or $\ZZ + 1/2$ (if $\beta \sim \pi$). A unit in such an order is of the form
\[
  u = a\left(1 + b\pi^{j - v(\beta)/2} \sqrt{\beta}\right), \quad a \in \OO_K^\cross, b \in \OO_K
\]
Since the factor $a$ belongs to the $\omega^0$-component, it can be ignored. The range of $[u] \in H^1(K,M)$, by Theorem \ref{thm:levels}\ref{lev:distance}, is $\L_i$, where
\[
  i = \begin{cases}
    \ds\ceil{\frac{(p-1)(j-1)}{p}}, & j \leq \ds\frac{p e}{p-1} + 1 \\
    e + 1, & j > \ds\frac{p e}{p-1} + 1.
  \end{cases}
\]
It is easy to see that all $i$ ($0 \leq i \leq e + 1$) are attained thereby.
\end{proof}

\begin{prop} \label{prop:ord ord}
Every generalized order in a quadratic extension $Q/K$ has a (not necessarily unique) dual order in the reflection extension $Q'$.
\end{prop}
\begin{proof}
Follows immediately from Lemma \ref{lem:ord cond} and Theorem \ref{thm:levels}\ref{lev:perp}.
\end{proof}

In the tame case, we evidently have the dual pairs 
\[
  \OO_Q \longleftrightarrow \OO_{Q'}, \quad Q \longleftrightarrow \OO
\]
for any $\OO \subsetneq \OO_{Q'}$. In the wild case, things are only a bit more involved:
\begin{prop} \label{prop:dual_orders_wild}
Let $\OO_{Q,j}$, $\OO_{Q',j'}$ be the orders in $Q$ and $Q'$ as parametrized in the proof of Lemma \ref{lem:ord cond}. A dual to $Q$ is any $\OO_{Q',j'}$ for which
\[
  j > \frac{p e}{p-1} + 1.
\]
For $0 \leq j \leq \frac{p e}{p - 1}$, a dual to $\OO_{Q,j}$ is $\OO_{Q',j'}$ where
\[
  j' = \frac{pe}{p-1} + 1 - j.
\]
\end{prop}
\begin{proof}
The only slightly nontrivial step is to show that, in the second case, the corresponding level indices
\[
  i = \ceil{\frac{(p-1)(j-1)}{p}} \textand i' = \ceil{\frac{(p-1)(j'-1)}{p}} 
\]
have sum $e$. But after noting that the arguments to the two ceilings have sum $\equiv 1/p \mod 1$, the summation becomes easy.
\end{proof}

The method of proof of Theorem \ref{thm:Scholz_for_locally_dual_orders} applies without change to generalized orders and yields the following.
\begin{thm}\label{thm:Scholz_gen}
Let $p \geq 3$ be a prime, let $K$ be a global field with $\zeta_p + \zeta_p^{-1} \in K$, and let Let $\OO_1 \subseteq Q$, $\OO_2 \subseteq Q'$ be dual generalized quadratic orders. Then
\begin{equation} \label{eq:Scholz_gen}
  \frac{\size{\Cl(\OO_1)[p]_\omega}}{\size{\Cl(\OO_2)[p]_\omega}} =
  \frac{p^{\1_{Q \isom K \cross K}}}{p^{\1_{Q' \isom K \cross K}}}
  \cdot \prod_{v \text{ of } K} \frac{\size{\OO_{2,\ell,\omega}^\cross}}{p^{\1_{Q_p \isom K \cross K}}}
\end{equation}
\end{thm}

\subsection{Relation to the Scholz reflection theorem}
\begin{examp}
Let $p = 3$, $K = \QQ$, $T = \QQ[\sqrt{D}]$, and $T' = \QQ[\sqrt{-3D}]$, where $D$ is a fundamental discriminant with $3 \nmid D$. Construct a pair of dual orders $\OO$, $\OO'$ by specification at each prime $\ell$ of $\ZZ$ as follows:
\begin{itemize}
  \item If $\ell \neq 3, \infty$, we take $\OO$ and $\OO'$ to be maximal at $\ell$, contributing nothing to the product in Theorem \ref{thm:Scholz_gen}.
  \item If $\ell = 3$, using Proposition \ref{prop:dual_orders_wild}, we see that the orders $\QQ_3[\sqrt{D}]$ and $\QQ_3[\sqrt{-27D}]$ are dual, as are $\QQ_3[\sqrt{9D}]$ and $\QQ_3[\sqrt{-3D}]$. The first contributes $1$ to the product, and the second contributes $3$.
\end{itemize}
The prime $\ell = \infty$ does not enter into the construction of the dual orders, but it introduces a factor $\size{H^0(\RR, M_D)}$ that depends on the sign of $D$. Finally, note that all of the class group $\Cl(\OO)$ of a quadratic order belongs to the $\omega$-eigenspace, the $1$-eigenspace being $\Cl(\ZZ) = 0$. So we get an equality, which was also noticed by Nakagawa (\cite{Nakagawa}, Theorem 0.5): 
\begin{cor}\label{cor:Scholz_Z}
If $D \equiv 0, 1$ mod $4$ is an integer, write $\Cl(D)$ for the class group of the quadratic ring over $\ZZ$ having discriminant $D$. Let $D$ be a fundamental discriminant not divisible by $3$. Then
\begin{align}
  {\size{\Cl(-27D)[3]}} &= {\size{\Cl(D)[3]}} \cdot 3^{\1_{D = -3} - \1_{D = 1} + \1_{D > 0}} \label{eq:ScholzZ1} \\
  {\size{\Cl(-3D)[3]}} &= {\size{\Cl(9D)[3]}} \cdot 3^{\1_{D = -3} - \1_{D = 1} + \1_{D > 0} - 1}. \label{eq:ScholzZ2}
\end{align}
\end{cor}
Both equations are generalizations of the Scholz reflection principle, which states that for $D \neq 1, -3$,
\[
  \size{\Cl(-3D)[3]} = {\size{\Cl(D)[3]}} \cdot 3^{\1_{D > 0} - \epsilon}
\]
where $\epsilon \in \{0,1\}$. This theorem shows that $\epsilon$ can be explained by the size of the kernel of either of the maps
\begin{equation} \label{eq:ord maps}
  \Cl(9D)/\Cl(9D)^3 \to \Cl(D)/\Cl(D)^3 \textor \Cl(-27D)/\Cl(-27D)^3 \to \Cl(-3D)/\Cl(-3D)^3.
\end{equation}
It also shows that exactly one of the maps \eqref{eq:ord maps} is an isomorphism, the other having kernel of size $3$---a theorem, perhaps, that has not appeared in the literature yet?
\end{examp}

\part{Reflection theorems: cubic rings}\label{part:cubic}

\section{Cubic Ohno-Nakagawa}
\label{sec:cubic}

The space $V(K)$ of binary cubic forms over a local or global field $K$ can have many integral models. Let $V_{\OO_K}$ be the lattice of binary cubic forms with trivial Steinitz class; these can be written as
\[
V_{\OO_K} = ax^3 + bx^2y + cxy^2 + dy^3 : a,b,c,d \in \OO_K,
\]
and we abbreviate the form $ax^3 + bx^2y + cxy^2 + dy^3$ to $(a,b,c,d)$. A theorem of Osborne classifies all lattices $L \subseteq V_{\OO_K}$ that are $\GL_2(\OO_K)$-invariant and \emph{primitive}, in the sense that $\pp^{-1}L \nsubseteq V(\OO_K)$ for all finite primes $\pp$ of $\OO_K$:
\begin{thm}[Osborne \cite{Osborne}, Theorem 2] \label{thm:Osborne}
  A primitive $\GL_2(\OO_K)$-invariant lattice in $V(\OO_K)$ is determined by any combination of the primitive $\GL_2(\OO_{K,\pp})$-invariant lattices in the completions $V(\OO_{K,\pp})$, which are:
  \begin{enumerate}[(a)]
    \item\label{it:3} If $\pp|3$, the lattices $\Lambda_{\pp,i} = \{(a,b,c,d) : b \equiv c \equiv 0 \bmod \pp^i \}$, for $0 \leq i \leq v_\pp(3)$;
    \item\label{it:2} If $\pp|2$ and $N_{K/\QQ}(\pp) = 2$, the five lattices
    \begin{align*}
      \Lambda_{\pp,1} &= V(\OO_{K,\pp}), \\
      \Lambda_{\pp,2} &= \{(a, b, c, d) \in V(\OO_{K,\pp}) : a + b + d \equiv a + c + d \equiv 0 \mod \pp \} \\
      \Lambda_{\pp,3} &= \{(a, b, c, d) \in V(\OO_{K,\pp}) : a + b + c \equiv b + c + d \equiv 0 \mod \pp \} \\
      \Lambda_{\pp,4} &= \{(a, b, c, d) \in V(\OO_{K,\pp}) : b + c \equiv 0 \mod \pp \} \\
      \Lambda_{\pp,5} &= \{(a, b, c, d) \in V(\OO_{K,\pp}) : a \equiv d \equiv b + c \mod \pp \},
    \end{align*}
    \item For all other $\pp$, the maximal lattice $V(\OO_{K,\pp})$ only.
  \end{enumerate}
\end{thm}
From the perspective of algebraic geometry, if $\pp | 2$, the latter four lattices are not true integral models, because they lose their $\SL_2$-invariance as soon as we extend scalars so that the residue field has more than $2$ elements. By contrast, if $\pp|3$, the $\SL_2$-invariance of the space $L_{\pp_i}$ can be established purely formally. This integral model, which we will call the space of \emph{$\pp^i$-traced} forms, will be the subject of our main reflection theorem in this part.

Although Osborne deals only with the case of $\Gamma(\OO_K)$, his method generalizes easily to the lattice
\[
  V(\OO_K,\aa) = \{ax^3 + bx^2y + cxy^2 + dy^3 : a \in \aa, b \in \OO_K, c \in \aa^{-1}, d \in \aa^{-2}\}
\]
that pops up when considering the maps
\[
  \Phi : M \to \Lambda^2 M
\]
that appear in the higher composition law Theorem \ref{thm:hcl_cubic_ring}. Here the relevant action of
\[
  \Gamma(\OO_K,\aa) = \Aut_{\OO_K}(\OO_K \oplus \alpha) =
  \left\{\begin{bmatrix}
    a_{11} & a_{12} \\
    a_{21} & a_{22}
  \end{bmatrix} \in \GL_2(K): a_{ij} \in \aa^{j-i}\right\}
\]
is nontrivial on both $M$ and $\Lambda^2 M$, thus affecting $V(\OO_K, \aa)$ via a twisted action
\begin{equation}\label{eq:twisted_action}
\left(\begin{bmatrix}
  a_{11} & a_{12} \\
  a_{21} & a_{22}
\end{bmatrix}
\mathop{.\vphantom{I}} \Phi\right)(x,y) = \frac{1}{a_{11} a_{22}-a_{12} a_{21}} \Phi(a_{11} x + a_{21} y, a_{12} x + a_{22} y).
\end{equation}
(Compare \cite{potf}, p.~142 and \cite{WGauss}, Theorem 1.2.) The twist by the determinant does not affect invariance of lattices but renders the action faithful, while otherwise scalar matrices that are cube roots of unity would act trivially. We sidestep this issue entirely by restricting the action to the group $\SL_2$, which preserves the discriminant $D \in \aa^{-2}$ of the form. The corresponding ring has discriminant $(\aa, D)$.

For instance, over $K = \ZZ$ there are ten primitive invariant lattices, comprising five types at $2$ and two types at $3$. The O-N-like reflection theorems relating all the types at $2$ were computed by Ohno and Taniguchi \cite{10lat} and will be considered later in this paper (Section \ref{sec:invarlat2}). While the behavior at $2$ admits only mild generalization, being based on the combinatorics of the finitely many cubic forms over $\FF_2$, the behavior at $3$ is robust. We begin by making some definitions needed to track the behavior of cubic forms and rings at primes dividing $3$.

If $\OO$ is a ring of finite rank over a Dedekind domain $\OO_K$, define its \emph{trace ideal} $\tr(\OO)$ to be the image of the trace map $\tr_{\OO/\OO_K} : \OO \to \OO_K$. Note that $\tr(\OO)$ is an ideal of $\OO_K$ and, since $1 \in \OO$ has trace $n = \deg(\OO/\OO_K)$, it is a divisor of the ideal $(n)$. In particular, if $\OO_K$ is a DVR, this notion is uninteresting unless $\OO_K$ has residue characteristic dividing $n$. Let $\tt$ be an ideal of $\OO_K$ dividing $(n)$. We say that the ring $\OO$ is \emph{$\tt$-traced} if $\tr(\OO) \subseteq \tt$.

By Theorem \ref{thm:hcl_cubic_ring}, we can parametrize cubic orders $\OO$ by their Steinitz class $\aa$ and index form
\[
\Phi(x\xi + y\eta) = (ax^3 + bx^2y + cxy^2 + dy^3)(\xi \wedge \eta)
\]
relative to a decomposition $\OO = \OO_K \oplus \OO_K\xi \oplus \aa\eta$, where $a \in \aa$, $b \in \OO_K$, $c \in \aa^{-1}$, and $d \in \aa^{-2}$. Then a short computation using the multiplication table from Theorem \ref{thm:hcl_cubic_ring} shows that, if $(1, \xi, \eta)$ is a normal basis, then $\tr(\xi) = -b$ and $\tr(\eta) = c$, so $\tr(\OO) = \<3, b, \aa c\>$. Thus the based $\tt$-traced rings over $\OO_K$ are parametrized by the rank-$4$ lattice of cubic forms
\[
\V_{\aa,\tt}(\OO_K) := \{ax^3 + bx^2y + cxy^2 + dy^3 : a \in \aa, b \in \tt, c \in \tt\aa^{-1}, d \in \aa^{-2}\},
\]
on which $\GL(\OO \oplus \aa)$ acts by the twisted action \eqref{eq:twisted_action}. For instance, if $\OO_K = \ZZ$, $\aa = (1)$, and $\tt = (3)$, this is the lattice of integer-matrix cubic forms considered in the introduction. Our goal in this section is to prove a generalization for all number fields $K$ and spaces $V(\OO_K, \aa, \tt)$.

\begin{thm}[\textbf{``Local cubic O-N''}] \label{thm:O-N_cubic_local}
  Let $K$ be a nonarchimedean local field, $\ch K \neq 3$. Let $V(D)$ be the composed variety of binary cubic forms of discriminant $D$, under the action of the group $\Gamma = \SL_2$. If $\alpha \in K^\cross$ and $\tau \mid 3$ in $\OO_K$, let $\V_{\alpha, \tau}(D)$ be the integral model of $V(K)(D)$ consisting of forms of the shape
  \[
  f(x,y) = a\alpha x^3 + b\tau x^2 y + c\alpha^{-1}\tau x y^2 + d\alpha^{-2} y^3,
  \]
  together with its natural action of $\G_{\alpha} = \SL(\OO_K \oplus \alpha\OO_K)$. Then the integral models
  \begin{equation}\label{eq:x_cubic_dual_t}
  \(\V_{1,\tau}(D), \SL_2\OO_K\) \textand \(\V_{1,3\tau^{-1}}(-27\tau^{-6}D), \SL_2\OO_K\),
  \end{equation}
  and consequently
  \begin{equation}\label{eq:x_cubic_dual_a_t}
  \(V_{\alpha,\tau}(D), \G_{\alpha}\) \textand \(V_{\alpha\tau^{-3},3\tau^{-1}}(-27 D), \G_{\alpha\tau^{-3}}\)
  \end{equation}
  are naturally dual with duality constant $N_{K/\QQ}(\tau) = \size{\OO_K/\tau\OO_K}$.
\end{thm}

The two formulations are easily seen to be equivalent. The first one is the one we will prove, but the second one has the needed form of a local reflection theorem to apply at each place to get the following global reflection theorem:
\begin{thm}[\textbf{O-N for traced cubic rings}]\label{thm:O-N_traced}
  Let
  \[
  \V_{\aa,\tt}(\OO_K) := \{f(x,y) = ax^3 + bx^2y + cxy^2 + dy^3 : a \in \aa, b \in \tt, c \in \tt\aa^{-1}, d \in \aa^{-2}\},
  \]
  a representation of
  \[
    \G_\aa \coloneqq \SL(\OO_K \oplus \aa).
  \]
  Note that $\V_{\aa,\tt}$ is the integral model of $\(V(K), \Gamma(K)\)$ parametrizing $\tt$-traced cubic rings over $\OO_K$ with Steinitz class $\aa$. For $D \in \tt^2\aa^{-2}$, define the class number
  \[
    h_{\aa,\tt}(D)
    = \sum_{\substack{\Phi \in \G_{\aa} \backslash \V_{\aa,\tt}(\OO_K) \\ \disc \Phi = D}} \frac{1}{\size{\Stab \Phi}}
    = \sum_{\substack{\Disc \OO = (\aa,D) \\ \tt\text{-traced}}} \frac{1}{\size{\Aut_K \OO}}.
  \]
  Then we have the global reflection theorem
  \begin{equation} \label{eq:O-N_traced}
    h_{\aa,\tt}(D) = \frac{3^{\#\{v|\infty : D \in (K_v^\cross)^2\}}}{N_{K/\QQ}(\tt)} \cdot h_{\aa\tt^{-3}, 3 \tt^{-1}}(-27D).
  \end{equation}
\end{thm}
\begin{proof}
Use Theorem \ref{thm:O-N_cubic_local} at each finite place. At the infinite places, the two integral models are necessarily naturally dual because $H^1(\RR, M_D) = 0$; but the duality constant depends on $H^0$, which depends on the sign of $D$ at each real place, as desired.
\end{proof}

Observe that taking $K = \QQ$, $\aa = 1$, $\tt = 1$ recovers Ohno-Nakagawa (Theorem \ref{thm:O-N}).

This also yields the extra functional equation for the Shintani zeta functions (see Corollary \ref{cor:Shintani} below).

We can rewrite our results in terms of Shintani zeta functions.

\begin{defn}
Let $K$ be a number field. If $\OO/\OO_K$ is a cubic ring of nonzero discriminant, the \emph{signature} $\sigma(\OO)$ of $\OO$ is the Kummer element $\alpha \in (K\tensor_\QQ \RR)^\cross / ((K\tensor_\QQ \RR)^\cross)^2$ corresponding to the quadratic resolvent of $\OO$. That is, it takes the value $\alpha_v = +1$ or $-1$ at each real place $v$ of $K$ according as $\OO_v \isom \RR \cross \RR \cross \RR$ or $\RR \cross \CC$, and $\alpha_v = 1$ at each complex place.
\end{defn}

\begin{defn}\label{defn:Shintani}
Given a number field $K$, a signature $\sigma \in (K\tensor_\QQ \RR)^\cross / ((K\tensor_\QQ \RR)^\cross)^2$, an ideal class $[\aa] \in \Cl(K)$, and an ideal $\tt \mid 3$, we define the \emph{Shintani zeta function}
\[
  \xi_{K, \sigma, [\aa], \tt}(s) = \sum_{\OO} \frac{1}{\size{\Aut_K(\OO)}} N_{K/\QQ}(\disc_K L)^{-s}
\]
where the sum ranges over all cubic orders $\OO$ over $\OO_K$ having signature $\sigma$, Steinitz class $\aa$, and trace ideal contained in $\tt$. We also define the \emph{Shintani zeta function} with unrestricted Steinitz class
\[
  \xi_{K, \sigma, \tt}(s) = \sum_{[\aa] \in \Cl(K)} \xi_{K, \sigma, [\aa], \tt}(s).
\]
\end{defn}
\begin{rem}
Confusingly, it is traditional to denote Shintani zeta functions by the Greek letter xi.
\end{rem}
\begin{rem}
By Minkowski's theorem on the finite count of number fields with bounded degree and discriminant, each term $n^{-s}$ has a finite coefficient, so the Shintani zeta function at least makes sense as a formal Dirichlet series. It generalizes the Shintani zeta functions for rings over $\ZZ$ mentioned in the introduction. Datskovsky and Wright \cite{DW2} study an adelic version of the Shintani zeta function; they show that $\xi_{K, \sigma, (1)}$ and $\xi_{K, \sigma, (3)}$ are entire meromorphic with at most simple poles at $s = 1$ and $s = 5/6$, satisfying an explicit functional equation. We surmise that the same method will prove the same for $\xi_{K, \sigma, [\aa], \tt}$. However, we do not consider the analytic properties here.
\end{rem}

Then we have the following corollary, which generalizes Conjecture 1.1 of Dioses \cite{Dioses}.
\begin{cor}[\textbf{the extra functional equation for Shintani zeta functions}]\label{cor:Shintani}
Let $K$ be a number field, $\sigma \in (K\tensor_\QQ \RR)^\cross / ((K\tensor_\QQ \RR)^\cross)^2$ a signature, $[\aa] \in \Cl(K)$ an ideal class, and $\tt \mid 3$ an ideal. Then the Shintani zeta function $\xi_{K, \sigma, [\aa], \tt}(s)$ satisfies an extra functional equation
\begin{equation}\label{eq:Shintani_a_t}
  \xi_{K, \sigma, [\aa], \tt}(s) = \frac{3^{\#\{v|\infty : \sigma_v = 1\} + 3[K : \QQ] s}}{N_{K/\QQ}(t)^{1 + 6 s}} \xi_{K, -\sigma, [\aa \tt^{-3}], 3\tt^{-1}}(s),
\end{equation}
Hence, summing over all $\aa$,
\begin{equation}\label{eq:Shintani_t}
  \xi_{K, \sigma, \tt}(s) = \frac{3^{\#\{v|\infty : \sigma_v = 1\} + 3[K : \QQ] s}}{N_{K/\QQ}(t)^{1 + 6 s}} \xi_{K, -\sigma, 3\tt^{-1}}(s),
\end{equation}
\end{cor}
\begin{proof}
Fix $\aa$ and $\tt$. Sum Theorem \ref{thm:O-N_traced} over all $D \in \tt^3\aa^{-2}$ of signature $\sigma$, weighting each $D$ by 
\[
  N_{K/\QQ}\(D\aa^2\)^{-s},
\]
the norm of the discriminant of the associated cubic rings. Then the left-hand side of the summed equality matches that of \eqref{eq:Shintani_a_t}. The right-hand side involves rings with discriminant ideal $27D\aa^2\tt^{-6}$, so a compensatory factor of
\[
  \frac{N_{K/\QQ}\(D\aa^2\)^{-s}}{N_{K/\QQ}\(27D\aa^2\tt^{-6}\)^{-s}} = \frac{3^{[K:\QQ]s}}{N_{K/\QQ}(\tt)^{6s}}
\]
must be added to the right-hand side to pull out the desired Shintani zeta function.
\end{proof}

In the succeeding subsections, we present three approaches to the local duality (Theorem \ref{thm:O-N_cubic_local}). First, we present a short conceptual proof in the special case that $\ch k_K \neq 3$, a ``tame'' case. Second, we explicitly compute the local orbit counters for a computational proof. Third, we organize the local orbits into families for a more conceptual general proof.

\subsection{A bijective proof of the tame case}
\label{sec:cubic_tame}
\begin{proof}[Proof of the tame case of Theorem \ref{thm:O-N_cubic_local}.]
Fix $D \in \OO_K$. Let $T = K[\sqrt{D}]$ be the corresponding quadratic algebra, and $T' = K[\sqrt{-3D}]$. For brevity we will write $H^i(T)$ for the cohomology $H^i(K, M_T)$ of the corresponding order-$3$ Galois module, and $H^i(T')$ likewise.

 Denote by $f(\sigma)$, for $\sigma \in H^1(T)$, the number of orders of discriminant $D$ in the corresponding cubic algebra $L_\sigma$; and likewise, denote by $f'(\tau)$, for $\tau \in H^1(T')$, the number of orders of discriminant $-3D$ in $L_\tau$. Our task is to prove that $f' = \hat f$. We note that if $-3$ is a square in $\OO_K^\cross$, then $T = T'$ and $f = f'$.

Note that $f$ is even: $\sigma$ and $-\sigma$ are parametrized by the same cubic algebra with opposite orientations of its resolvent. The Fourier transform of an even, rational-valued function on a $3$-torsion group $H^1(T)$ is again even and rational-valued. So far, so good.

Our method will be first to prove the duality at $0$: that is, that
\begin{align}
  f'(0) &= \hat f(0) \label{eq:dual 0 1} \\
  f(0) &= \widehat{f'}(0). \label{eq:dual 0 2}
\end{align}
Let us explain how this implies that $f' = \hat f$. We compute $\size{H^1(T)}$ using the self-orthogonality of unramified cohomology:
\[
  \size{H^1(T)} = \size{H^1(T)^\ur} \cdot \size{H^1(T')^\ur}
  = \size{H^0(T)} \cdot \size{H^0(T')}.
\]
So there are basically three cases:
\begin{enumerate}[(a)]
  \item If neither $D$ nor $-3D$ is a square in $K_v$, then $H^1(T) \cong H^1(T') \cong 0$, and \eqref{eq:dual 0 1} trivially implies that $f' = \hat f$. 
  \item If one of $D$, $-3D$ is a square, then $H^1(T)$ and $H^1(T')$ are one-dimensional $\FF_3$-vector spaces. The space of even functions on each is $2$-dimensional, and
  \begin{align*}
    \g &\mapsto (\g(0),\hat{\g}(0)) \\
    \g' &\mapsto (\widehat{\g'}(0), \g'(0))
  \end{align*}
  are corresponding systems of coordinates on them. Consequently, the two equations \eqref{eq:dual 0 1} and \eqref{eq:dual 0 2} together imply that $\hat f = f'$.
  \item Finally, if $D$ and $-3D$ are both squares, then $f$ is a function on the two-dimensional $\FF_3$-space $H^1(T) \cong H^1(T')$ which we would like to prove self-dual. Note that $H^1(T)$ has four subspaces $W_1,\ldots, W_4$ of  dimension $1$. Consider the following basis for the five-dimensional space of even functions on $H^1(T)$:
  \[
    f_1 = \1_{W_1}, \ldots, f_4 = \1_{W_4}, f_5 = \1_0.
  \]
  Note that $f_1,\ldots, f_4$ are self-dual (the Tate pairing is alternating, so any one-dimensional subspace is isotropic), while $f_5$ is not: indeed $\hat f_5(0) \neq f_5(0)$. Thus if \eqref{eq:dual 0 1} holds, then $f$ is a linear combination of $f_1,\ldots,f_4$ only and hence $\hat f = f$.
\end{enumerate}

We have now reduced the theorem to a pair of identities, \eqref{eq:dual 0 1} and \eqref{eq:dual 0 2}. By symmetry, it suffices to prove \eqref{eq:dual 0 2}, which may be written 
\begin{equation} \label{eq:dual 0}
  f(0) \stackrel{?}{=} \hat f(0) = \frac{1}{\size{H^0(T')}}\sum_{\tau \in H^1(T')} f(\tau).
\end{equation}
The proof is clean and bijective.


The left-hand side of \eqref{eq:dual 0} counts orders of discriminant $D$ in the split algebra $K \cross T$. These can be straightforwardly parametrized as
\[
  \OO_{K} + 0 \cross \aa,
\]
where $\aa$ is a multiplicatively closed lattice in $T$, that is, an invertible ideal in some quadratic order $\OO \subseteq T$. (Here we use that, in a quadratic algebra, any lattice $\aa$ is an invertible ideal with respect to its endomorphism ring $\End \aa$. This fails for higher-degree algebras, which we will encounter later.)

The sum on the right-hand side of \eqref{eq:dual 0} counts all cubic orders of discriminant $-3D$. Any cubic order $C$ of discriminant $-3D$ can be assigned an ideal in $T$ as follows. Let $L$ be the fraction algebra of $C$. By Theorem \ref{thm:Kummer_new}, we have the description
\[
  L = K + \{\xi \sqrt[3]{\delta} + \bar\xi \sqrt[3]{\bar\delta} \mid \xi \in T'\}
\]
for some $\delta \in T'^{N=1}/\( T'^{N=1}\) ^3$; and so, since $3$ is invertible in $\OO_K$,
\begin{equation} \label{eq:lattice C}
  C = \OO_K + \{\xi \sqrt[3]{\delta} + \bar\xi \sqrt[3]{\bar\delta} \mid \xi \in \cc\}
\end{equation}
for some lattice $\cc$ in $T'$. Now by Theorem \ref{thm:hcl_cubic_sbi}, we get that $(\OO_D,\cc,\delta)$ is a \emph{self-balanced ideal}, that is,
\begin{equation} \label{eq:balanced cubic}
  \delta \cc^3 \subseteq \OO_D, \quad N(\cc) = (t) \text{ is principal}, \textand N(\delta) t^3 = 1,
\end{equation}
Now $\cc$ need not be invertible in $\OO_D$; but let $\OO := \End \cc$. Note that $\OO = \OO_{D/\pi^{2i}}$ for some positive integer $i$. Then form the \emph{shadow}
\[
  \aa = \frac{\delta \cc^3}{\pi^{i}}.
\]
Since the norm is multiplicative  on invertible ideals, the properties of $\cc$ in \eqref{eq:balanced cubic} can be recast as properties of $\aa$:
\begin{equation} \label{eq:balanced aa}
  \aa \subseteq \End \aa \textand N_{\OO_D}(\aa) = 1.
\end{equation}
The first of these says that $\aa$ is multiplicatively closed, and the second that the ring $C = \OO_K + 0 \cross \aa$ corresponding to $\aa$ has discriminant $D$.

It remains to show that there are exactly $\size{H^0(T)}$ cubic orders $C$ corresponding to each shadow $\aa$ satisfying \eqref{eq:balanced aa}, weighting each $C$ by the number of isomorphic copies of $C$ in its fraction algebra $L = KC$.

First, we simply count $C$ up to isomorphism. This is the same as counting the self-balanced ideals $(\OO_D,\cc,\delta)$ up to the equivalence relation
\begin{equation} \label{eq:equiv rel 2}
  \cc \mapsto \lambda \cc, \quad \delta \mapsto \lambda^{-3} \delta \quad (\lambda \in T'^\cross)
\end{equation}
Note the slight subtlety in this step: we would like to define an isomorphism between the associated cubic algebras $L_{\delta}$, $L_{\lambda^{-3}\delta}$ by
\[
  \xi \sqrt[3]{\delta} + \bar\xi \sqrt[3]{\bar\delta} \mapsto
  \lambda\xi \sqrt[3]{\lambda^{-3}\delta} + \bar\lambda \bar\xi \sqrt[3]{\bar\lambda^{-3} \bar\delta},
\]
which works, but only after potentially rescaling the cube roots on the right-hand side by a suitable 3rd root of unity in $K$ so that their product is $N(\lambda) t = N(\lambda) \sqrt[3]{\delta} \sqrt[3]{\bar\delta}$.

Let $\aa = \alpha\OO$ be the given shadow, where $\OO = \OO_{D/\pi^{2i}}$ is its ring of invertibility. Clearly $\cc$ must be invertible with regard to $\OO$. The possible $(\cc,\delta)$ may be found by fixing $\cc = \OO$ and taking $\delta = \pi^i\alpha\epsilon$ where $\epsilon \in \OO^\cross$ is constrained by the requirement that $N(\delta)$ be a cube. If, without loss of generality, we scale $\alpha$ so that $N(\pi^i\alpha)$ is a cube, then the admissible values are $\{\epsilon \in \OO^{\cross} : N(\epsilon) \in (K^\cross)^3\}$. Now in the equivalence relation \ref{eq:equiv rel 2}, the multipliers $\lambda$ preserving $\cc = \OO$ are $\lambda \in \OO^\cross$, so we must consider $\epsilon$ up to $\big(\OO^\cross\big)^3$. Since $N : \OO^\cross / \big(\OO^\cross\big)^3 \to \OO_K^\cross / \big(\OO_K^\cross\big)^3$ is surjective, the number of distinct $\epsilon$, which is the number of nonisomorphic $C$, is simply
\[
  \frac{\Size{\OO^\cross / \big(\OO^\cross\big)^3}}{\Size{\OO_K^\cross / \big(\OO_K^\cross\big)^3}}.
\]
Next, we weight each $C$ by the number of isomorphic copies of $C$ in $L = KC$. The automorphisms of $C$ are given by
\[
  \xi \sqrt[3]{\delta} + \bar\xi \sqrt[3]{\bar\delta} \mapsto \omega \xi \sqrt[3]{\delta} + \bar\omega \bar\xi \sqrt[3]{\bar\delta}
\]
where $\omega \in \OO_{T'}$ satisfies $\omega^3 = \omega \bar\omega = 1$. But choices $\omega \in \OO$ fix $C$, so we must quotient out by those, and the number of isomorphic copies is
\[
  \frac{\Size{\OO_{T'}^{N=1}[3]}}{\Size{\OO^{N=1}[3]}} = \frac{\Size{\OO_{T'}^\cross[3]}}{\Size{\OO^\cross[3]}},
\]
and the total number of cubic orders we seek is the product
\begin{equation}\label{eq:ans 1}
  \frac{\Size{\OO^\cross / \big(\OO^\cross\big)^3}}{\Size{\OO_K^\cross / \big(\OO_K^\cross\big)^3}} \cdot
  \frac{\Size{\OO_{T'}^\cross[3]}}{\Size{\OO^\cross[3]}}.
\end{equation}
To maneuver this into the required form, first note that
\[
  \frac{\Size{\OO^\cross / \big(\OO^\cross\big)^3}}{\Size{\OO^\cross[3]}} = \frac{\Size{\OO_{T'}^\cross / \big(\OO_{T'}^\cross\big)^3}}{\Size{\OO_{T'}^\cross[3]}}
\]
by the Snake Lemma, since $\OO_{T'}^\cross/\OO^\cross$ is finite; so \eqref{eq:ans 1} takes the form
\[
  \frac{\Size{\OO_{T'}^\cross / \big(\OO_{T'}^\cross\big)^3}}{\Size{\OO_K^\cross / \big(\OO_K^\cross\big)^3}},
\]
which we can now compute directly to equal
\[
  \Size{\OO_{T'}^{N=1} / \big(\OO_{T'}^{N=1}\big)^3} = \frac{T^{\prime,N=1}/\big(T^{\prime,N=1}\big)^3}{\size{H^1(T')}} = \frac{\size{H^1(T)}}{\size{H^1(T')}} = \size{H^0(T)}.
\]
\end{proof}

\subsection{A computational proof}

We now turn our attention to the wild case of Theorem \ref{thm:O-N_cubic_local}. We present two proofs, one computational, one more conceptual.

\subsubsection{Trace ideals of maximal orders}
Our first step is to compute $\tr(\OO_L)$ for the maximal orders of all cubic extensions $L$. The answer is delightfully simple.
\begin{prop}\label{prop:tr_idl}
Let $L/K$ be a cubic \'etale algebra over a $3$-adic field. The trace ideal of the maximal order of $L$ is
\[
  \tr(\OO_L) = (3) + \mm_K^{e - \ell(L)}
\]
where $\ell(L)$ is the level.
\end{prop}
\begin{rem}
Hyodo (\cite{Hyodo}, equation (1--4); see also Xia and Zhukov \cite{XiaZhu14}) proves a theorem like this one for an invariant he calls the \emph{depth} of a ramified extension. The depth is in fact closely related to the level and offset.
\end{rem}
\begin{proof}
We first dispose of the case that $L$ is not totally ramified, that is, has splitting type $111$, $12$, $3$, or $1^21$, by noting that in all these cases $\tr(\OO_L) = (1)$ and $\ell(L) = e$.

Now let $L = K[\pi_L]$ be a totally ramified extension. We have $\tr(\OO_L) = (3, \tr \pi_L, \tr \pi_{L}^2)$. There are two cases.

\begin{enumerate}[$\<1\>$]
  \item \label{case:off0}
In this case we assume $v_K(\tr \pi_L) \geq v_K(\tr \pi_L^2)$, including the case that both are infinite. In this case it is possible to adjust $\pi_L$ by a multiple of $\pi_L^2$ so as to make the trace vanish. Therefore we may assume that $\tr \pi_L = 0$, so the Eisenstein minimal polynomial of $\pi_L$ is a depressed cubic,
\[
  \phi(\x) = x^3 + u x + v,
\]
with $v_K(u) \geq 1$ and $v_K(v) = 1$. We have $\tr(x^2) = 2 u$, so $\tr(\OO_K) = (3, u)$.
\[
  \disc(L) = \disc \phi = 4 u^3 - 27 v^2.
\]
Now $v_K(27 v^2) = 3e + 2$ is not a multiple of $3$, so the two terms have unequal valuation. If $v_K(u) \leq e$, then $4u^3$ dominates so $L$ has level $\ell = e - v_K(u)$ and offset $\theta = 0$; and the trace ideal is $(u) = \mm_K^\ell$. If $v_K(u) > e$, then $-27 v^2$ dominates so $v_K(\disc L) = 3e + 2$. This is the case $L = K[\sqrt[3]{\pi_K}]$ of a uniformizer radical extension. The level is $-1$, the offset is $-1$ and the trace ideal is $(3)$.
  \item \label{case:off1}
In this case we assume $v_K(\tr \pi_L) < v_K(\tr \pi_L^2)$. Note that $v_K(\tr \pi_K\pi_L) \leq v_K(\tr \pi_L^2)$, so it is possible to adjust $\pi_L^2$ by a multiple of $\pi_K\pi_L$ to produce an element $\rho$ such that $\tr \rho = 0$ and $v_L(\rho) = 2$. The minimal polynomial of $\rho$ is a depressed cubic
\[
  \phi(\x) = \x^3 + u \x + v,
\]
with $v_K(u) \geq 1$ and $v_K(v) = 2$. Now $\rho$ does not generate all of $\OO_L$; instead, an $\OO_K$-basis of $\OO_L$ is $(1, \rho^2/\pi_K, \rho)$ so $\tr(\OO_L) = (3, u/\pi_K)$.
\[
  \disc(L) = \frac{\disc \phi}{\pi_K^2} = \frac{4 u^3 - 27 v^2}{\pi_K^2}.
\]
Now $v_K(27 v^2) = 3e + 4$ is not a multiple of $3$, so the two terms in the numerator have unequal valuation. If $v_K(u) \leq e + 1$, then $4u^3$ dominates so $L$ has level $\ell = e - v_K(u) + 1$ and offset $\theta = 1$; and the trace ideal is $(u/\pi_K) = \mm_K^\ell$. If $v_K(u) > e + 1$, then $-27 v^2$ dominates and we have a uniformizer radical extension again. \qedhere
\end{enumerate}
\end{proof}
The foregoing proof has a corollary on the structure of totally ramified cubic extensions which will be important to us.
\begin{cor}
Let $L/K$ be a totally ramified extension.
\begin{enumerate}[(a)]
\item If $\theta(L) = 0$, then $L$ has a traceless uniformizer.
\item If $\theta(L) = 1$, then $L$ has a traceless element of valuation $2$.
\end{enumerate}
\end{cor}
\begin{proof}
Simply note that Case \ref{case:off0} occurs only when $\theta = 0$ or $-1$, and Case \ref{case:off1} occurs only when $\theta = 1$ or $-1$.
\end{proof}
\begin{prob}
Find an analogue of Proposition \ref{prop:tr_idl} for $\GA(\ZZ/p\ZZ)$-extensions of $p$-adic fields, $p \geq 5$.
\end{prob}

\subsubsection{The subring zeta function}
If $L/K$ is a cubic $3$-adic algebra and $0 \leq t \leq e$, let
\[
  g(L, t) = \sum_{\substack{\text{orders } \OO \subseteq L \\ \tr(\OO) \subseteq \mm_K^t}} z^{v_K(\disc \OO)/2} \in \ZZ[[z^{1/2}]]
\]
be the generating function of $\mm_K^t$-traced orders in $L$. This is related to the \emph{subring zeta function} (see Section \ref{sec:subring}). The factor of $1/2$ in the exponent is used (quite arbitrarily) to make a factor of $z$ correspond to passing to a subring of index $\mm_K$. Note that $g(L, t) \in z^{1/2}\ZZ[[z]]$ or $\ZZ[[z]]$ according as $T$ is ramified or not.

Note that if $t \leq e - \ell(L)$, then every order in $L$ is automatically $\mm_K^t$-traced, so $g(L, t) = g(L, 0)$ is simply the generating function for all orders that was computed by Datskovsky and Wright.

We will proceed to compute $g(L, t)$ for all $L$ and $t$. We begin by tabulating the possible splitting types for a cubic algebra $L$:
\[
\begin{tabular}{cccc}
$\sigma(L)$ & $\ell(L)$ & $\theta(L)$ & $v_K(\Disc L)$ \\ \hline
$111$ & $e$ & $0$ & $0$ \\
$12$ & $e$ & $0$ & $0$ \\
$3$ & $e$ & $0$ & $0$ \\
$1^21$ & $e$ & $1$ & $1$ \\
$1^3$ & $0 \leq \ell < e$ & $0$ or $1$ & $3(e - \ell) + \theta$ \\
$1^3$ & $-1$ & $-1$ & $3e + 2$
\end{tabular}
\]
\begin{lem}\label{lem:traced_count_overview}
We have
\[
  g(L, t) = z^{v_K(\disc \OO_L)/2} \cdot s^{\sigma,\theta}(n)
\]
for a certain power series $s^{\sigma,\theta}(n) \in \ZZ[[z]]$ that depends only on three parameters: the splitting type $\sigma = \sigma(L)$, the offset $\theta = \theta(L)$ (which is redundant unless $\sigma = 1^3$), and the \emph{trace deficit} $n = t - e + \ell(L)$.
\end{lem}
In due course, we will prove this theorem and determine the series $s^{\sigma,\theta}(n)$ by directly counting the rings involved. For now, we assume it and proceed to deduce Theorem \ref{thm:O-N_cubic_local}, which, in this notation, says:
\begin{thm}[\textbf{Local cubic O-N}]\label{thm:O-N_traced_local}
Considering $g(L, t)$ as a function of $L$, its Fourier dual is given by
\[
  \hat{g}(L, t) = q^t z^{3t - \frac{3e}{2}} g(L, e - t).
\]
\end{thm}
\begin{proof}
The proof proceeds by writing $g(L, t)$ as a linear combination of characteristic functions of level spaces. Let $g(T, \ell, t) = g(L, t)$ for any $L$ of resolvent torsor $T$ and level $\ell$, where
\[
  e \geq \ell \geq \ell_{\min} = \begin{cases}
    -1 & \text{if $T'$ is split} \\
    0 & \text{otherwise.}
  \end{cases}
\]
By Theorem \ref{thm:levels}\ref{lev:size_Li}, such $L$ exists, and by Lemma \ref{lem:traced_count_overview}, the series $g(T, L, t)$ is independent of which $L$ of this level we choose, with one exception: if $T \cong K \cross K$ is split and $\ell = 0$, then $L$ could have splitting type $111$ or $3$. We resolve the ambiguity as follows: give $g(K \cross K, 0, t)$ the value of $g(L, t)$ when $\sigma(L) = 3$, and introduce a symbol $g(K \cross K, -1, t)$ with the value $g(K \cross K \cross K, t)$ for splitting type $111$. Correspondingly, define the level space
\[
  \L_{e + 1} = \{0\},
\]
in spite of the fact that cubic algebras of splitting types $3$ and $111$ both have level $e$. Then in all cases, if we set
\[
  \ell_{\max} = \begin{cases}
    e + 1 & \text{if $T$ is split} \\
    e & \text{otherwise}
  \end{cases}
\]
then we have
\[
  g(L, t) = \sum_{\ell_{\min} \leq \ell \leq \ell_{\max}} \1_{L \in \L_{\ell}\setminus\L_{\ell+1}} g(T, \ell, t) = \sum_{\ell_{\min} \leq \ell \leq \ell_{\max}} \1_{L \in \L_{\ell}} (g(T, \ell, t) - g(T, \ell-1, t))
\]
where
\begin{equation}\label{eq:g=0}
  g(T, \ell_{\min} - 1, t) = 0
\end{equation}
Now $\hat{\1}_{\L_{\ell}} = c_\ell \1_{\L_{e - \ell}}$, where
\[
  c_{\ell} = \frac{\size{\L_\ell}}{\size{H^0(T)}} = \begin{cases}
    \frac{1}{3} & \ell = e+1 = \ell_{\max} \\
    q^{e - \ell} & 0 \leq \ell \leq e \\
    3 q^e & \ell = -1 = \ell_{\min}.
  \end{cases}
\]
Thus it suffices to prove that
\begin{equation} \label{eq:hpart symm}
  c_{\ell} (g(T, \ell, t) - g(T, \ell - 1, t)) = q^t z^{3t - \frac{3e}{2}} (g(T', e - \ell, t) - g(T', e - \ell - 1, t))
\end{equation}
for $\ell_{\min} \leq \ell \leq \ell_{\max}$. It is easy to verify that flipping $T \mapsto T'$, $\ell \mapsto e - \ell$, $t \mapsto e - t$ transforms \eqref{eq:hpart symm} to an equivalent equation, so we assume that $\ell > 0$ and $\ell + t \geq e$, which cuts down the number of cases.

We now enumerate the cases of \eqref{eq:hpart symm}, which by Lemma \ref{lem:traced_count_overview}, depend on the splitting types and offsets of the fields $L$ appearing. Recall that when the level $\ell$ ($0 \leq \ell \leq e$) and resolvent torsor $T$ of a cubic algebra $L$ are known, the offset can be determined by the congruence
\[
  3\ell + \theta = v_K(\disc L) \equiv v_K(\beta_T) \mod 2
\]
where $\beta_T$ is a Kummer element for $T$ (here we simply have $v_K(\beta_T) \equiv v_K(\disc L) \mod 2$). Now since $\beta_T\beta_{T'} = -3$ (up to squares), we have
\[
  v_K(\beta_T) + v_K(\beta_T') \equiv e \mod 2.
\]
From this, we find that each offset in \eqref{eq:hpart symm} determines the other three, even without knowing $T$, and there are only six cases:
\[
\begin{tabular}{cc|cccc}
$\sigma(T)$ & $\ell$ & $\sigma(L_{T, \ell})$ & $\sigma(L_{T, \ell+1})$ & $\sigma(L_{T', e-\ell})$ & $\sigma(L_{T', e-\ell+1})$ \\ \hline
any & $0 < \ell < e$ & $1^3(\theta = 0)$ & $1^3(\theta = 1)$ & $1^3(\theta = 0)$ & $1^3(\theta = 1)$ \\
any & $0 < \ell < e$ & $1^3(\theta = 1)$ & $1^3(\theta = 0)$ & $1^3(\theta = 1)$ & $1^3(\theta = 0)$ \\
$1^2$ & $e$ & $1^21$ & $1^3(\theta = 0)$ & $1^3(\theta = 1)$ & --- \\
$2$ & $e$ & $1 2$ & $1^3(\theta = 1)$ & $1^3(\theta = 0)$ & --- \\
$11$ & $e$ & $3$ & $1^3(\theta = 1)$ & $1^3(\theta = 0)$ & $1^3(\theta = -1)$ \\
$11$ & $e$ & $111$ & $3$ & $1^3(\theta = -1)$ & ---
\end{tabular}
\]
The dashes for $\ell = \ell_{\min}$ indicate that the corresponding term was declared zero in \eqref{eq:g=0}. We now write each term of \eqref{eq:hpart symm} in terms of $s^{\sigma,\theta}(n)$ using Lemma \ref{lem:traced_count_overview}. Thanks to our assumptions that $\ell > 0$ and $\ell + t \geq e$, the terms on the right side involve only $s^{\sigma,\theta}(n)$ for $n < 0$. In this case, we can replace $n$ by $0$ because all orders in $L$ automatically satisfy the trace condition. The corresponding generating function was computed by Datskovsky and Wright and will soon be recovered by us (see \eqref{eq:cub0})
\[
  s_0 = s^{1^3,\theta}(0) = \frac{1}{(1-z)(1-qz^3)},
\]
independent of $\theta$. The theorem is now reduced to the following lemma.
\end{proof}
\begin{lem} \label{lem:cub count}
For $t \geq 0$,
\begin{align}
  s^{1^3,\theta=0}(n) - z^2 s^{1^3,\theta=1}(n - 1) &= q^n z^{3n}(1 - z^2) s_0 \label{eq:cub long} \\ 
  s^{1^3,\theta=1}(n) - z s^{1^3,\theta=0}(n - 1) &= q^n z^{3n}(1 - z) s_0 \label{eq:cub short} \\
  s^{3}(n) - z^2 s^{1^3,\theta=1}(n - 1) &= q^n z^{3n}(1 - z) s_0 \label{eq:cub3} \\
  s^{12}(n) - z^2 s^{1^3,\theta = 1}(n - 1) &= q^n z^{3n} s_0 \label{eq:cub12} \\
  s^{1^21}(n) - z s^{1^3,\theta = 0}(n - 1) &= q^n z^{3n} s_0 \label{eq:cub1^21} \\
  s^{111}(n) - s^{3}(n) &= 3q^n z^{3n+1} s_0. \label{eq:cub111}
\end{align}
\end{lem}
This lemma can be viewed as a set of coupled difference equations for computing the $s^{\sigma,\theta}(t)$.

\subsubsection{Counting traced cubic orders: the proofs of Lemma \ref{lem:traced_count_overview} and Lemma \ref{lem:cub count}}

It now remains to count the $\mm_K^t$-traced orders in each cubic algebra $L$. These are sublattices of $\OO_L$ containing $1$ and satisfying (a) the \emph{ring condition,} that is, closure under multiplication, and (b) the \emph{trace condition} that each of their elements has trace in $\mm_K^t$. We first simplify these two conditions.

If $\OO \subseteq \OO_L$ is a sublattice containing $1$, the quotient $\OO_L/\OO \cong (\OO_L/\OO_K)/(\OO/\OO_K)$ is a finite group generated as an $\OO_K$-module by two elements. We may write
\begin{equation} \label{eq:x ord ij}
  \OO_L/\OO = (\OO_K/\mm_K^i) \xi \oplus (\OO_K/\mm_K^j) \eta,
\end{equation}
where $\xi$ and $\eta$ are generators such that $\OO_L$ and $\OO$ have $\OO_K$-bases $(1, \xi, \eta)$ and $(1, \pi_K^i\xi, \pi_K^j\eta)$. By symmetry we may assume that $i \geq j$. By varying the basis $(1, \xi, \eta)$, we get all lattices $\OO$ of ``index $(i,j)$'' in the sense that \eqref{eq:x ord ij} holds. (This is an example of a \emph{reduced basis}, which we will use more systematically in the quartic case.) Note that
\begin{equation} \label{eq:x ord}
  \OO = \OO_K + \mm_K^i\OO_L + \mm_K^j \eta,
\end{equation}
from which it is easy to see that two bases $(1, \xi, \eta)$ and $(1, \xi', \eta')$ yield the same $\OO$ if and only if
\begin{equation}\label{eq:lattice eqv}
  \eta \equiv u \eta' \mod \mm_K^{j - i} \OO_L, \quad \text{some } u \in (\OO_K/\mm_K^{i - j})^\cross.
\end{equation}
Note that $\xi$ is irrelevant. Also note that when $i = j$ there is a single lattice, the \emph{content ring} $\OO_K + \mm_K^i\OO_L$.

Having constructed all lattices of index $(i,j)$, we test whether they are rings using Theorem \ref{thm:hcl_cubic_ring}\ref{cubic:lift}: $\OO$ is a ring if and only if its index form is integral. Let the index form of $L$ be
\[
  \Phi_L(x\xi + y\eta) = (ax^3 + bx^2y + cxy^2 + dy^3) (\xi \wedge \eta).
\]
Then
\begin{equation} \label{eq:index form of subring}
  \Phi_L(x\pi_K^i\xi + y\pi_K^j\eta) = (a\pi_K^{2i - j} + b\pi_K^i + c\pi_K^j + d\pi_K^{2j - i}) (\pi_K^i\xi \wedge \pi_K^j\eta).
\end{equation}
\begin{figure}[ht]
\[
\setlength{\unitlength}{0.5cm}
\begin{picture}(7,7.2)(-0.5,-0.5)
\put(-0.5,0){\line(1,0){7}}
\put(0,-0.5){\line(0,1){7}}
\drawline(0,0)(6.2,6.2)
\drawline(0,0)(6.2,3.1)
\put(6,0.1){\makebox(0,0)[b]{$i$}}
\put(-0.1,6.5){\makebox(0,0)[r]{$j$}}
\put(4,0.5){\makebox(0,0)[b]{Root zone}}
\put(4,2.7){\rotatebox{35}{Free zone}}
\end{picture}
\]
\caption{Two zones for the indices $(i,j)$ for a candidate subring of $\OO_L$}
\label{fig:ij}
\end{figure}
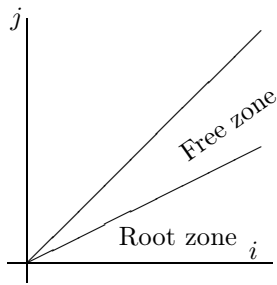

From this we can deduce that there are two kinds of pairs $(i,j)$. If $j \leq i \leq 2j$, then every lattice with indices $(i,j)$ is a ring: we say that $(i,j)$ is in the \emph{free zone.} If $i > 2j$, the ring condition is that $\pi_K^{i - 2j} | d$, that is,
\begin{equation}\label{eq:ring}
    \Phi_{\OO_L}(\eta) \equiv 0 \mod \mm_K^{i - 2j}.
\end{equation}
Because $\eta$ must in this sense be a root of $\Phi$, we call this range of $(i,j)$ values the \emph{root zone.} Note that if $L$ has splitting type $(3)$, then $\Phi$ has no roots even mod $\mm_K$ and hence no orders in the root zone. 
For the remaining splitting types, the root-zone orders can be subdivided according to which root of $\Phi$---which ``$1$'' in the splitting type---$\eta$ reduces to mod $\mm_K$.

Finally we must test our orders for the trace condition $\tr \OO \subseteq \mm_K^{e - \ell + n}$, where $\ell = e - v_K(\tr \OO_L) = \max\{\ell(L), 0\}$. We may assume that $0 \leq n \leq \ell$. The trace ideal of an order \eqref{eq:x ord} can be computed by
\[
  \tr \OO = \tr (\OO_K + \mm_K^i\OO_L + \mm_K^j \eta) = (3) + \mm_K^{\ell + i} + \mm_K^j \tr \eta.
\]
Now if $n > i$, the trace condition is impossible; if $n \leq j$, it is automatic; and if $j < n \leq i$, it is contingent on $\eta$: the condition is
\begin{equation}\label{eq:trace}
\frac{\tr \eta}{\pi_K^\ell} \equiv 0 \mod \mm_K^{n - j}.
\end{equation}
The ring condition \eqref{eq:ring} is cubic in $\eta$, while the trace condition \eqref{eq:trace} is linear. We will presently see how they interact.

\subsubsection{Traced orders in the free zone}
In this subsubsection we evaluate
\[
  s^{\mathrm{FZ}}(n) = \sum_{\substack{\OO \subset \OO_L \text{ free-zone orders,} \\ \tr(\OO) \subseteq \mm_K^n \tr(\OO_L)}} z^{v_K([\OO_L : \OO])}.
\]
In the free zone, we know that the lattices of index $(i,j)$ are parametrized by elements $\eta \in \OO_L/\OO_K$ not divisible by $\mm_K$ up to the equivalence relation \eqref{eq:lattice eqv}, and all are orders. We get $q^{i-j-1}(q+1)$ orders if $i > j$, just $1$ if $i = j$; and all of these satisfy the trace condition if $n \leq j$. It remains to test them on the trace condition \eqref{eq:trace} under the hypothesis that $j < n \leq i$.

It is not hard to see that $\OO_L$ has a basis $(1, \xi_0, \eta_0)$ where $\tr(\xi_0) = 0$. Then $\eta_0$ can be taken to generate $\tr(\OO_L) = \mm_K^{e - \ell}$. Write $\eta = x\xi_0 + y\eta_0$ (the $\OO_K$-component being irrelevant). The trace condition becomes $\mm_K^{n - j} | \x$. In particular, $\x$ is not a unit so $\y$ must be, and we can scale so that $\y = 1$; then there are $q^{i - n}$ choices for $\x$ (which is defined modulo $\mm_K^{i - j}$).

Summing, the generating function for the free zone is
\begin{equation} \label{eq:cubFZ}
  s^{\mathrm{FZ}}(n) = \sum_{j < n \leq i \leq 2j} q^{i - n} z^{i + j} + \sum_{n \leq j < i \leq 2j} (q+1) q^{i-j-1} z^{i+j} + \sum_{n \leq j = i} z^{i+j}.
\end{equation}

\subsubsection{Traced orders in the root zone}
In this subsubsection we fix a root $\bar\eta_0$ of the index form $\Phi$ of $\OO_L$ modulo $\mm_K$. We evaluate the generating function of $\bar\eta_0$-orders, that is, orders whose corresponding basis element $\eta$ reduces to $\bar\eta_0$ modulo $\mm_K$:
\[
  s^{\mathrm{RZ},\bar\eta_0}(n) = \sum_{\substack{\OO \subset \OO_L \text{ $\bar{\eta_0}$-orders,} \\ \tr(\OO) \subseteq \mm_K^n \tr(\OO_L)}} z^{v_K([\OO_L : \OO])}
\]
In each of the various cases that we will encounter, we will choose a basis $(1, \xi_0, \eta_0)$ such that $\eta_0$ reduces to $\bar\eta_0$ mod $\OO_K + \mm_K\OO_L$. Then the lattices belonging to this root are parametrized by elements $\eta = x\xi_0 + y\eta_0$ where $\mm_K | \x$; so $\mm_K \nmid \y$ and we can scale so that $\y = 1$. Consequently we can take $\eta = \x'\pi_K\xi_0 + \eta_0$ where $\x'$ runs over the residue classes mod $\mm_K^{i - j - 1}$.

The first case is that of the \emph{simple} root, the ``1'' with no exponent that appears in the splitting types $111$, $12$, and $1^21$. In this case Hensel's lemma tells us that $\bar\eta_0$ can be lifted to an element $\eta_0$ with $\Phi(\eta_0) = 0$, and we can identify it explicitly: $L = K \cross T$ splits and $\eta_0 = (1; 0)$. Note that $\tr(\eta_0) = 1$ and $\ell = e$. We can complete to a basis $(1,\xi_0,\eta_0)$ with $\tr(\xi_0) = 0$. The index form of $L$ has the form
\[
  \Phi(x\xi_0 + y\eta_0) = (ax^3 + bx^2y + cxy^2)(\xi_0 \wedge \eta_0)
\]
where $\mm_K \nmid c$, since the root is simple. When plugging in a value $\eta = \x'\pi_K\xi_0 + \eta_0$ with $\y = 1$, the $cxy^2$ term will dominate so the ring condition is $\mm_K^{i-2j-1}|\x$: we get $q^j$ rings. But $\tr(\eta) = 1$ so the trace condition cannot be fulfilled unless $n \leq j$, in which case it is vacuous. Thus we get a root subring generating function
\begin{equation} \label{eq:cubSR}
  s^{\mathrm{SR}}(n) = \sum_{\substack{j \geq n \\ i \geq 2j}} q^j z^{i+j}.
\end{equation}

Now assume that $\bar\eta_0$ is the \emph{multiple} root of one of the splitting types $1^21$ and $1^3$: that is, in a basis $(1, \xi_0, \eta_0)$ with $\eta_0$ lifting $\bar\eta_0$, the index form is
\[
  \Phi_L(x\xi_0 + y\eta_0) = (ax^3 + bx^2y + cxy^2 + dy^3) (\xi_0 \wedge \eta_0)
\]
with $\mm_K | c$ and $\mm_K | d$. A now-standard trick shows that $\mm_K^2 \nmid d$, that is, $v_K(d) = 1$: if not, then applying the formula \eqref{eq:index form of subring} with $i = 0$ and $j = -1$ would show that $(1, \xi_0, (\eta_0 + u)/\pi_K)$ is a basis of a ring for some $u \in \OO_K$: that is, $\OO_L$ would not be the maximal order in $L$. So $\mm_K^2 \nmid d$, and the root $\bar\eta_0$ mod $\mm_K$ has \emph{no} lift to mod $\mm_K^2$. This shows that the only rings belonging to this root occur for $i = 2j + 1$, the very edge of the root zone, where the ring condition \eqref{eq:ring} is mod $\mm_K$ only (and is therefore automatically satisfied). There are $q^{i - j - 1} = q^j$ rings.

The trace condition can be expressed in terms of the index form using the fact, previously mentioned, that
\[
  \tr(\xi_0) = -b \textand \tr(\eta_0) = c.
\]
Thus $\tr(\eta) = \tr(\x \pi_K \xi_0 + \eta_0) = c - \pi_K bx$. We may assume that $\ell = e - \min\{v_K(b), v_K(c)\} > 0$, as when $\ell = 0$ we have $n = 0$ and no trace condition. There are two cases.
\begin{enumerate}[$\<1\>$]
\item \label{case:lean}
If $v_K(b) \geq v_K(c)$, then $v_K(\tr(\eta)) = v_K(c) = e - \ell$ no matter what $\x$ we pick, and so the trace condition is unsatisfiable for all $n > j$.
\item \label{case:fat}
But if $v_K(b) < v_K(c)$, then the equation $c - \pi_K bx = 0$ has a solution $x_0$, and the solutions to the trace condition
\[
  \frac{c - \pi_K bx}{\pi_K^{e - \ell}} \equiv 0 \mod \mm_K^{n - j}
\]
are $\x = x_0 + u$ for $v_K(u) \geq n - j - 1$ (since $v_K(b) = e - \ell$). There are $q^{i - n} = q^{2j + 1 - n}$ solutions.
\end{enumerate}

We must now determine for which algebras $L$ Cases \ref{case:lean} and \ref{case:fat} occur. For splitting type $1^3$, if we take our basis $(1, \xi_0, \eta_0) = (1, \pi_L, \pi_L^2)$, they match up exactly with Cases \ref{case:off0} and \ref{case:off1} in the proof of Proposition \ref{prop:tr_idl} and therefore correspond to the offsets $\theta = 0$ and $\theta = 1$, respectively. For splitting type $1^21$, we have $\mm_K \nmid b$ (since otherwise $\Phi$ would have a triple root) and again Case \ref{case:fat} occurs. Since $\theta = 1$ in this case too, we can divide up the generating functions by $\theta$-value rather than splitting type:
\begin{align}
  s^{\mathrm{MR},\theta=0}(n) &= \sum_{j \geq n} q^j z^{3j + 1} \label{eq:cubMR lean} \\
  s^{\mathrm{MR},\theta=1}(n) &= \sum_{j \geq n} q^j z^{3j + 1} + \sum_{\frac{n - 1}{2} \leq j < n} q^{2j + 1 - n} z^{3j + 1}. \label{eq:cubMR fat}
\end{align}
We have left out the case $\theta = -1$ of the uniformizer radical extensions. For these, $n = 0$, so the trace condition is vacuous and either \eqref{eq:cubMR lean} or \eqref{eq:cubMR fat} yields the correct generating function
\begin{equation} \label{eq:cub0}
  s^{\mathrm{MR}}(0) = \sum_{j \geq 0} q^j z^{3j + 1}.
\end{equation}

\subsubsection{Putting it together}
The subring generating functions $s^{\sigma,\theta}(n)$ are now derived by summing the free-zone and root-zone contributions:
\begin{equation}
\begin{aligned}
  s^{111}(n) &= s^{\mathrm{FZ}}(n) + 3 s^{\mathrm{SR}}(n) \\
  s^{12}(n) &= s^{\mathrm{FZ}}(n) + s^{\mathrm{SR}}(n) \\
  s^{3}(n) &= s^{\mathrm{FZ}}(n) \\
  s^{1^21}(n) &= s^{\mathrm{FZ}}(n) + s^{\mathrm{SR}}(n) + s^{\mathrm{MR},\theta=1}(n) \\
  s^{1^3, \theta}(n) &= s^{\mathrm{FZ}}(n) + s^{\mathrm{MR},\theta}(n).
\end{aligned}
\end{equation}
The proof of Lemma \ref{lem:traced_count_overview} is now complete.

\begin{proof}[Proof of Lemma \ref{lem:cub count}]
It is possible to prove Lemma \ref{lem:cub count} in an automated fashion by summing the doubly geometric series \eqref{eq:cubFZ}, \eqref{eq:cubSR}, \eqref{eq:cubMR lean}, and \eqref{eq:cubMR fat}, solving the linear recurrences \eqref{eq:cub long}--\eqref{eq:cub111}, and checking that the resulting rational functions agree for both even and odd $n$. We here present a more illuminating method, which does not attempt to sum all the series but simply manipulates their terms---that is, it is very nearly a bijective proof.

We begin with a simplification of the generating functions $s^{1^3,\theta = 0}$ and $s^{1^3, \theta = 1}$.
\begin{lem}
\begin{align}
  s^{1^3,\theta = 0}(n) &= \sum_{\substack{a \geq 0 \\ b \geq \max \{3a, \frac{3a + 3n}{2}\}}} q^a z^b \label{eq:region0} \\
  s^{1^3,\theta = 1}(n) &= \sum_{\substack{a \geq 0 \\ b \geq \max \{3a, \frac{3a + 3n - 1}{2}\}}} q^a z^b \label{eq:region1}
\end{align}
\end{lem}
\begin{proof} We have
\begin{align*}
  s^{1^3, \theta = 0}(n) &= s^{\mathrm{FZ}}(n) + s^{\mathrm{MR},\theta = 0}(n) \\
  &= \sum_{j < n \leq i \leq 2j} q^{i - n} z^{i + j} + \sum_{n \leq j < i \leq 2j} (q+1) q^{i-j-1} z^{i+j} + \sum_{n \leq j = i} z^{i+j}
    + \sum_{j \geq n} q^j z^{3j + 1}.
\intertext{Splitting the two summands in the $(q + 1)$ factor and combining the former with the third sum,}
  s^{1^3, \theta = 0}(n) &= \sum_{j < n \leq i \leq 2j} q^{i - n} z^{i + j} + \sum_{n \leq j < i \leq 2j} q^{i - j} z^{i+j} + \sum_{n \leq j < i \leq 2j} q^{i - j - 1} z^{i+j} + \sum_{j \geq n} q^j z^{3j + 1}.
\intertext{Reindexing the third sum by $i \mapsto i+1$,}
  s^{1^3, \theta = 0}(n) &= \sum_{j < n \leq i \leq 2j} q^{i - n} z^{i + j} + \sum_{n \leq j < i \leq 2j} q^{i - j} z^{i+j} + \sum_{n \leq j \leq i \leq 2j + 1} q^{i - j} z^{i + j + 1} + \sum_{j \geq n} q^j z^{3j + 1} \\
  &= \sum_{j < n \leq i \leq 2j} q^{i - n} z^{i + j} + \sum_{n \leq j < i \leq 2j + 1} q^{i - j} z^{i+j} + \sum_{n \leq j \leq i \leq 2j + 1} q^{i - j} z^{i + j + 1}.
\end{align*}
It is now not hard to check that the last two sums together include every term $q^a z^b$ with $a \geq 0$ and $b \geq \max\{3a,a + 2n\}$ exactly once. As for the first sum, it includes every term with $a \geq 0$ and $\frac{3a + 3n}{2} \leq b \leq a + 2n$ exactly once. So we have \eqref{eq:region0}.

To prove \eqref{eq:region1}, it suffices to tack on the additional root-zone terms
\[
  s^{\mathrm{MR},\theta = 1}(n) - s^{\mathrm{MR},\theta = 0}(n)
  = \sum_{\frac{n - 1}{2} \leq j < n} q^{2j + 1 - n} z^{3j + 1}
  = \sum_{\substack{0 \leq a < n \\ b = \frac{3a + n - 1}{2}}} q^a z^b. \qedhere
\]
\end{proof}
In other words, $s^{1^3,\theta}(n)$ can be viewed as a sum of terms $q^a z^b$ for integer points $(a,b)$ in a certain region $\R$:
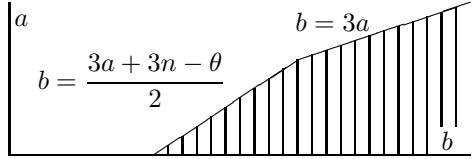
\begin{figure}[ht]
\setlength{\unitlength}{0.25in}
\[
\begin{picture}(9.6,3.2)
\put(0,0){\line(1,0){9.6}}
\put(0,0){\line(0,1){3.2}}
\put(3,0){\line(3,2){3}}
\put(6,2){\line(3,1){3.6}}
\put(0.1,2.7){$a$}
\put(9,0.1){$b$}
\put(2,1){\makebox(2.5,1)[rb]{$b = \dfrac{3a + 3n - \theta}{2}$}}
\put(5,2.5){\makebox(2.5,1)[rb]{$b = 3a\vphantom{q}$}}
\put(3.3,0){\line(0,1){0.2}}
\put(3.6,0){\line(0,1){0.4}}
\put(3.9,0){\line(0,1){0.6}}
\put(4.2,0){\line(0,1){0.8}}
\put(4.5,0){\line(0,1){1.0}}
\put(4.8,0){\line(0,1){1.2}}
\put(5.1,0){\line(0,1){1.4}}
\put(5.4,0){\line(0,1){1.6}}
\put(5.7,0){\line(0,1){1.8}}
\put(6.0,0){\line(0,1){2.0}}
\put(6.3,0){\line(0,1){2.1}}
\put(6.6,0){\line(0,1){2.2}}
\put(6.9,0){\line(0,1){2.3}}
\put(7.2,0){\line(0,1){2.4}}
\put(7.5,0){\line(0,1){2.5}}
\put(7.8,0){\line(0,1){2.6}}
\put(8.1,0){\line(0,1){2.7}}
\put(8.4,0){\line(0,1){2.8}}
\put(8.7,0){\line(0,1){2.9}}
\put(9.0,0.6){\line(0,1){2.4}}
\put(9.3,0.6){\line(0,1){2.5}}
\end{picture}
\]
\caption{Terms $q^a z^b$ appearing in the count of traced orders of a totally ramified cubic field}
\label{fig:region}
\end{figure}

In proving \eqref{eq:cub long} and \eqref{eq:cub short}, we observe that changing $n$ and $\theta$ moves only the left side of the region $\R$, while multiplication by powers of $z$ moves the whole region to the right. The differences appearing in \eqref{eq:cub long} and \eqref{eq:cub short} have the pleasant property that the composite of these two transformations causes the left sides to coincide, causing cancellation everywhere but along the top side $b = 3a$. In symbols:
\begin{align*}
  s^{1^3,\theta=0}(n) - z^2 s^{1^3,\theta=1}(n - 1)
  &= \sum_{\substack{a \geq 0 \\ b \geq \max\{3a, \frac{3a + 3n}{2}\}}} q^a z^b
  - \sum_{\substack{a \geq 0 \\ b \geq \max\{3a + 2, \frac{3a + 3n}{2}\}}} q^a z^b \\
  &= \sum_{\substack{a \geq n \\ 3a \leq b \leq 3a + 1}} q^a z^b \\
  &= (1 + z) \sum_{a \geq n} q^a z^{3a} \\
  &= \frac{(1 + z) q^n z^{3n}}{1 - qz^3} \\
  &= q^n z^{3n} (1 - z^2) s^{1^3}(0).
\end{align*}
This proves \eqref{eq:cub long}, and the proof of \eqref{eq:cub short} is analogous.

As for equations \eqref{eq:cub3}--\eqref{eq:cub111}, the linear combination $\eqref{eq:cub long} - \eqref{eq:cub3}$ reduces to a formula for $s^{\mathrm{MR}, \theta = 0}$,
\[
  s^{\mathrm{MR},\theta = 0}(n) = q^n z^{3n+1}(1 - z) s^{1^3}(0),
\]
while the linear combinations $\eqref{eq:cub12} - \eqref{eq:cub3}$, $\eqref{eq:cub short} - \eqref{eq:cub1^21}$, and \eqref{eq:cub111} all reduce to a formula for $s^{\mathrm{SR}}$:
\[
  s^{\mathrm{SR}}(n) = q^n z^{3n+1} s^{1^3}(0).
\]
These formulas are easily derived from the series formulations \eqref{eq:cubMR lean} and \eqref{eq:cubSR}.
\end{proof}

Also easily derived from the foregoing is a formula for the number of traced orders in a given cubic algebra, in other words, for the coefficient of a single power $z^b$ in one of the series $g(L,t)$:

\begin{thm}[\textbf{the traced subring zeta function}]\label{thm:traced_count}
Let $L$ be a cubic algebra over a $3$-adic field $K$ with discriminant $\disc L = \mm_K^{d_0}$. Let $d$ and $t$ be integers satisfying the necessary restrictions
\[
  d \geq d_0, \quad d \equiv d_0 \mod 2, \quad 0 \leq t \leq e_K.
\]
Then the number $g(L, d, t)$ of $\mm_K^t$-traced orders of discriminant $\mm_K^d$ in $L$ is a linear combination of the three functions
\begin{align*}
  g^{1^3}(d_0, d, t) &= \frac{q^r - 1}{q - 1}, \quad r = \begin{cases}
    0, & d < 3t \\
    \floor{\frac{d}{3}} - t + 1, & 3t \leq d \leq 6t - d_0 \\
    \floor{\frac{d - d_0}{6}} + 1, & d \geq 6t - d_0
  \end{cases} \\
  g^{3}(d_0 = 0, d, t) &= \frac{q^r - 1}{q - 1}, \quad r = \begin{cases}
    0, & d < 3t \\
    \floor{\frac{d}{3}} - t + 1, & 3t \leq d \leq 6t \\
    \frac{d}{2} - 2 \ceil{\frac{d}{6}} + 1, & d \geq 6t
  \end{cases} \\
  g^{\mathrm{SR}}(d_0 = 0 \text{ or } 1, d, t) &= \frac{q^r - q^t}{q - 1}, \quad r = \begin{cases}
    t, & d < 6t \\
    \floor{\frac{d}{6}} - t + 1, & d \geq 6t
  \end{cases}
\end{align*}
in a manner dependent on the splitting type of $L$:
\begin{itemize}
  \item $\sigma(L) = 1^3$: $g = g^{1^3}$
  \item $\sigma(L) = 1^21$: $g = g^{1^3} + g^{\mathrm{SR}}$
  \item $\sigma(L) = 3$: $g = g^3$
  \item $\sigma(L) = 12$: $g = g^3 + g^{\mathrm{SR}}$
  \item $\sigma(L) = 111$: $g = g^3 + g^{3\mathrm{SR}}.$
\end{itemize}
\end{thm}
\begin{rem}
When $t = 0$ (and here we need no longer assume that $\ch k_K = 3$), we recover the formulas for orders in a cubic field computed by Datskovsky and Wright and written more explicitly by Nakagawa and the author.
\end{rem}

\subsection{A bijective proof of the wild case}\label{sec:wild_bij}
Ideally we would desire a bijective proof of global reflection identities such as Theorem \ref{thm:O-N_traced}. This is still beyond reach. However, using the machinery of Galois cohomology and Poitou-Tate duality, we reduced this theorem to a local result, Theorem \ref{thm:O-N_traced_local}. In this section we will prove this theorem using a bijective method.

It is not obvious what ``bijective'' means when trying to prove that two local weightings are Fourier transforms of each other. Recall that in Theorem \ref{thm:levels}, we constructed \emph{level spaces} $H^1(T) \supseteq \L_0 \supsetneq \L_1 \supsetneq \cdots \supsetneq \L_e$ with the property that, for $i$,
\begin{equation}\label{eq:level space dual}
  \L_i^\perp = \L_{e-i}.
\end{equation}
As in the preceding proof, we extend this notation slightly. If $\L_e \neq \{0\}$, which happens exactly when $T \cong K \cross K$ is split, then let $\L_{e+1} = \{0\}$; and dually, if $\L_0 \neq H^1(T)$, which happens exactly when $T'$ is split, then let $\L_{-1} = H^1(T)$. Then $\L_e \setminus \L_{e+1}$ consists of the unramified field extension of $K$, and $\L_{-1} \setminus \L_0$ consists of the \emph{uniformizer radical extensions (UREs)} $K[\sqrt[3]{\pi}]$, the cubic extensions of maximal discriminant-valuation $3e + 2$. We have
\[
  \widehat{\1_{\L_{e+1}}} = \frac{1}{3} \1_{\L_{-1}} \textand \widehat{\1_{\L_{-1}}} = 3q^e \1_{\L_{e+1}}.
\]
We also let
\begin{align*}
  \ell_{\min} &= \begin{cases}
    -1 & T' \cong K \cross K \\
    0 & \text{otherwise},
  \end{cases} &
  \L_{\min} &= \L_{\ell_{\min}} = H^1(T), \\
  \ell_{\max} &= \begin{cases}
    e+1 & T \cong K \cross K \\
    e & \text{otherwise},
  \end{cases} &
  \L_{\max} &= \L_{\ell_{\max}} = \{0\}.
\end{align*}
The level spaces $\L_i$ and their associated characteristic functions $L_i$ will be central to our proof. Our strategy is as follows: we group all cubic rings whose resolvent torsor is $T$ into \emph{families} $\F$ with the following properties:
\begin{itemize}
  \item All rings in a family have the same discriminant and trace ideal.
  \item All rings in a family $\F$ are contained in \'etale algebras $L$ belonging to some level space $\L_i$; this $\L_i$ is called the \emph{support} $\supp(\F)$ of the family.
  \item Each $L \in \L_i$ has the same number of orders in the family; this number is called the \emph{thickness} $\th(\F)$ of the family.
\end{itemize}
The proof of Theorem \ref{thm:O-N_traced_local} will then consist in exhibiting an involution between the families of support $\L_i$ and $\L_{e-i}$ which affects their thicknesses, discriminants, and trace ideals in such a manner that they contribute equally to both sides of the theorem. The remainder of this section will be spent in carrying this out.

\begin{lem} \label{lem:cubic families}
Let $T/K$ be a quadratic torsor and let $n_T = v_K(\disc T) \in \{0,1\}$.

Then the cubic orders whose resolvent torsor is $T$ can be partitioned into families $\F_{n,k}$ indexed by the pairs of integers $(n,k)$ satisfying the conditions
\[
  0 \leq k \leq \floor{\frac{n}{3}}, \quad n \equiv n_T \mod 2,
\]
with the following properties:
\begin{enumerate}
  \item The rings in $\F_{n,k}$ have discriminant ideal $(\pi^n)$ and trace ideal $(\pi^{\min\{k,e\}})$.
  \item The support and thickness of each $\F_{n,k}$ depend on which of three zones the pair $(n,k)$ belongs to, as follows:
  \begin{equation}
    \begin{tabular}{cccc}
      Zone & $(n,k)$ & $\supp(\F_{n,k})$ & $\th(\F_{n,k})$ \\ \hline
      I & $0 \leq k < \dfrac{n}{6}$ & $\L_{\max}$ & $\size{H^0(T)}$ \\
      II & $\dfrac{n\vphantom{h}}{6} \leq k \leq \floor{\dfrac{n}{3}} - \dfrac{n}{6} + \dfrac{e}{2}$ & $\L_{e - 2k + \floor{\frac{n}{3}}}$ & $q^{\floor{\frac{n}{3}} - k}$ \\
      III & $\floor{\dfrac{n}{3}} - \dfrac{n\vphantom{h}}{6} + \dfrac{e}{2} < k \leq \floor{\dfrac{n}{3}}$ & $\L_{\min}$ & $q^{\floor{\frac{n}{3}} - k}$
    \end{tabular}
  \end{equation}
\end{enumerate}
\end{lem}
The relative positions of these zones follow a pattern like that in quadratic O-N, up to $O(1)$ discrepancies in the indices:
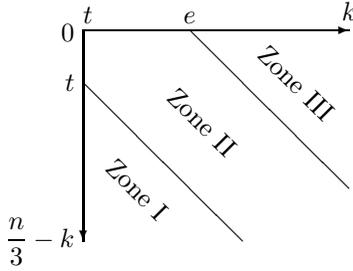
\[
\setlength{\unitlength}{1em}
\begin{picture}(13,9)(-2,-8.2)
  \put(0,0){\vector(1,0){10}}
  \put(0,0){\vector(0,-1){8}}
  \put(0,-2){\line(1,-1){6}}
  \put(4,0){\line(1,-1){6}}
  \put(-2,-9){\makebox(2,2)[r]{$\dfrac{n}{3} - k$~}}
  \put(-2,-3){\makebox(2,2)[r]{$t$~}}
  \put(-2,-1){\makebox(2,2)[r]{$0$\vphantom{$^2$}~}}
  \put( 0,0){\makebox(2,2)[bl]{$t$\vphantom{|}}}
  \put( 3,0){\makebox(2,2)[b]{$e$\vphantom{|}}}
  \put( 9,0){\makebox(2,2)[b]{$k$\vphantom{\big|}}}
  \put(2,-1){\rotatebox{-45}{\makebox(7,0){Zone II}}}
  \put(6,0){\rotatebox{-45}{\makebox(5.7,0){Zone III}}}
  \put(0,-4){\rotatebox{-45}{\makebox(5.7,0){Zone I}}}
\end{picture}
\]
\begin{proof}[Proof of Theorem \ref{thm:O-N_traced_local}]
Once Lemma \ref{lem:cubic families} is proved, we can prove Theorem \ref{thm:O-N_traced_local} quite simply by sending the family $\F = \F_{n,k}$ to $\F' = \F_{n',k'}$, where
\begin{align*}
  n' = n + 3e - 6t \\
  k' = \floor{\frac{n}{3}} - k + e - t.
\end{align*}
If the original $\F$ satisfied the bounds $t \leq k \leq \floor{\frac{n}{3}}$, then it is easy to see that $t' \leq k' \leq \floor{\frac{n}{3}}$ where $t' = e - t$, and likewise $n \equiv n_T$ mod $2$ implies $n' \equiv n_{T'}$. Thus $\F'$ is a family of rings of resolvent torsor $T'$ whose trace ideal is contained in $(\pi^{e-t})$. It is not hard to see that $\F'$ lies in zone III, II, or I according as $\F$ lies in zone I, II, or III. We leave it to the reader to check the needed identities
\begin{align*}
  \supp(\F') &= \supp(\F)^\perp \\
  \th(\F') &= \frac{\size{\supp(\F)}}{\size{H^0(T)}} \cdot \th(\F).
\end{align*}
\end{proof}
\begin{rem}
Zone I, which is supported on $\L_{\max} = \{0\}$, consists precisely of those rings whose structure uses in an essential way that $L$ is \emph{split,} that is, has more than one field factor. Zone II has the approximate shape of a band of constant width,
\[
  \frac{n}{6} \leq k \leq \frac{n}{6} + \frac{e}{2};
\]
but the dependency on the value of $n$ modulo $3$ attests to a waviness of the boundary between zones II and III that cannot be avoided.
\end{rem}
\begin{proof}[Proof of Lemma \ref{lem:cubic families}]
We now begin to enumerate all the orders in every cubic $K$-algebra $L$ and arranging them into families. Let $\disc \OO_L = (\pi^{n_0})$ and $\tr(\OO_L) = (\pi^{k_0})$, and let $\theta = n_0 - 3k_0$. As our investigations of the structure of cubic fields have found, we have $\theta \in \{0,1,2\}$, the case $\theta = 2$ corresponding to the URE.

Let $[1, \xi_0, \eta_0]$ be a basis for $\OO_L$. We can arrange so that $\eta_0$ is traceless and $\tr(\xi_0)$ is a generator for the trace ideal $\tr(\OO_K)$. Any order $C \subseteq L$ then has a unique basis of the form
\begin{equation} \label{eq:trace basis}
  [1, \quad \xi = \pi^i \xi_0 + u \eta_0, \quad \eta = \pi^j \eta_0]
\end{equation}
where $i$ and $j$ are nonnegative integers and $u$ ranges over a system of coset representatives in $\OO_K/\pi^j\OO_K$. (Note the departure from the reduced basis used in the preceding subsection.) Such a $C$ has discriminant valuation
\[
  n = v(\disc C) = n_0 + 2i + 2j
\]
and trace ideal
\[
  (3, \tr(\pi^i \xi_0)) = \pi^{\min\{e, k_0 + i\}}.
\]
Let $k = k_0 + i$. With one exception, namely when $L$ is a URE (which case we will handle later), we will place such a ring $C$ into the family $\F_{n,k}$. We must now compute the sizes of the families we have thus constructed.

Whether or not a lattice $C$ with a basis \eqref{eq:trace basis} is actually a ring is determined by the integrality of its \emph{index form.} The following lemma reduces the number of coefficients to be checked from four to two.
\begin{lem}
Let $\xi, \eta \in \OO_L$ be integral elements of a nondegenerate cubic algebra $L$ over the field of fractions $K$ of a Dedekind domain $\OO_K$ such that the sublattice $C = \OO_K\<1, \xi, \eta\>$ is of full rank. If the outer coefficients
\[
  \Phi_C(\xi) = \frac{\Phi_{\OO_L}(\xi)}{\pi^{[\OO_L : C]}} \textand \Phi_C(\eta) = \frac{\Phi_{\OO_L}(\eta)}{[\OO_L : C]}
\]
of the index form of $C$ are integral, then the entire index form of $C$ is integral and $C$ is a ring.
\end{lem}
\begin{proof}
If the whole index form of $C$ is integral, then there is a ring $C' = \<1, \xi', \eta'\>$ with the same index form (with respect to its basis) as $C$. Using the identity of index forms for $L$ over $K$, we can embed $C'$ into $L$ with $\xi' = \xi + u$, $\eta' = \eta + v$ for some $u, v \in K$. But since $\xi', \xi, \eta', \eta$ are integral elements and $\OO_K$ is integrally closed, we have $u, v \in \OO_K$ so $C' = C$.

So it suffices to prove that the index form is integral. This is a local statement, so we may assume that $\OO_K$ is a DVR. Passing to a finite extension, we may assume that $L \isom K \cross K \cross K$ is totally split. Let
\[
  \xi = (a_1; a_2; a_3) \textand \eta = (b_1; b_2; b_3).
\]
Then
\[
  D = [\OO_L : C] = \det \begin{bmatrix}
    1 & 1 & 1 \\
    a_1 & a_2 & a_3 \\
    b_1 & b_2 & b_3
  \end{bmatrix}
\]
We are given that the outer coefficients
\[
  c_0 = \frac{1}{D}(a_1 - a_2)(a_2 - a_3)(a_3 - a_1) \textand
  c_3 = \frac{1}{D}(b_1 - b_2)(b_2 - b_3)(b_3 - b_1)
\]
of the index form of $C$ are integral. We wish to prove that the same applies to the two middle coefficients. By symmetry, we can consider just the $x^2y$-coefficient
\[
  c_1 = \frac{1}{D}\left[\sum_{i=1}^3(a_i - a_{i+1})(a_{i+1} - a_{i+2})(b_{i+2} - b_i)\right].
\]
(Here, and for the rest of the proof, indices are modulo $3$.) It is easy to verify that
\[
  c_1 = -a_1 + 2a_2 - a_3 + \frac{3}{D}(a_1 - a_2)(a_2 - a_3)(b_3 - b_1).
\]
So it is enough to show that
\[
  d_1 = \frac{1}{D}(a_1 - a_2)(a_2 - a_3)(b_3 - b_1)
\]
is integral, or more generally any of the three
\[
  d_i = \frac{1}{D}(a_i - a_{i+1})(a_{i+1} - a_{i+2})(b_{i+2} - b_i).
\]
But we see that
\[
  a_0^2 a_3 = d_1 d_2 d_3.
\]
Since $a_0$ and $a_3$ have nonnegative valuation, the three $d_i$ cannot all have negative valuation, completing the proof.
\end{proof}
\begin{rem}
The hypothesis that $L$ be nondegenerate is likely nonessential.
\end{rem}

We can now resume the proof of Lemma \ref{lem:cubic families}. Let the index form of $\OO_L$ be
\[
  \Phi_{\OO_L}(x\xi_0 + y\eta_0) = ax^3 + bx^2y + cxy^2 + dy^3.
\]
Of the coefficients of the index form of $C$, we focus on the outer coefficients,
\[
  \Phi_C (\xi)
  = \frac{a \pi^{3i} + b \pi^{2i} u + c \pi^{i} u^2 + d u^3}{\pi^{i+j}}
  \textand \Phi_C (\eta) = d \pi^{2 j - i}.
\]
The latter coefficient is the simpler one, depending only on $i$ and $j$. Due to the tracelessness of $\eta_0$, we have $3|c$. We must then have $\pi^2 \nmid d$, or else $\OO_L$ would be nonmaximal, and the discriminant analysis of Proposition \ref{prop:tr_idl} shows that $v_K(d) = \theta$. Thus the condition $\Phi_C(\eta) \in \OO_K$ comes out to $2j - i + \theta \geq 0$, which simplifies to $n \geq 3k$.

(Incidentally, when $k \leq e$, the relation $n \geq 3k$ expresses an important relation between the trace ideal of a ring and its discriminant, generalizing the observation that an integer-matrix cubic form has discriminant divisible by $27$.)

There thus remains the $\xi$-condition $\Phi_C (\xi) \in \OO_K$, which informally states that $\xi$ is a root of $\Phi_{\OO_L}$ modulo $\pi^{i + j}$. Now $\Phi_{\OO_L}$ is a homogeneous binary form, and it is natural to consider its roots on the projective lines $\PP^1(\OO_K/\pi^m)$; the factorization over the field $\OO_K/\mm$, for instance, gives the splitting type $\sigma(\OO_L)$. However, in our situation there is a distinguished point on this projective line, at least for $m < e - k_0$, namely the traceless point $\eta_0$: and the line is thereby subdivided into an affine line and a portion at infinity. The point at infinity mod $\mm$ is a root of $\Phi_{\OO_L}$ if and only if $\theta > 0$. This is the motivation for the calculations to be undertaken now.

Suppose first that we are in zone I, that is, $n > 6k$, which translates into $j > 2i + \theta$. Then the $\Phi_C(\xi)$ condition
\[
  \pi^{j - 2i} \mid a + b u \pi^{-i} + c u^2 \pi^{-2i} + d u^3 \pi^{-3i}
\]
is clearly dominated by a non-integral last term unless $u$ is of the form $\pi^i u'$, in which case it simplifies to
\[
  \pi^{j - 2i} \mid a + b u' + c u'^2 + d u'^3.
\]
In other words, $\xi_0 + u' \eta_0$ must be a root of $\Phi_{\OO_L}$ modulo $\mm^{j - 2i}$. The condition $j - 2i > \theta$ rules out any contribution from a multiple root modulo $\mm$, which never lifts to mod $\mm^2$ (or else $\OO_L$ would be nonmaximal). So the only roots that contribute are the \emph{simple} roots occurring if $L$ has splitting type $111$, $12$, or $1^21$. There are $\size{H^0(T)}$ simple roots, and none of them are traceless (to be explicit, they are at $(1;0)$ for each decomposition $L \cong K \cross T$ into a linear and a quadratic algebra). By Hensel's lemma, each simple root has a unique lift to any modulus. So the solutions $u'$ form a union of $\size{H^0(T)}$ congruence classes modulo $\mm^{j - 2i}$. Since $u' = u/\pi^i$ is defined modulo $\mm^{j - i}$, there are
\[
  \size{H^0(T)} \cdot q^i = \size{H^0(T)} \cdot q^k
\]
rings for each pair $(i,j)$. This completes the construction of the families $\F_{n,k}$ in zone I.

Now suppose that $(n,k)$ is in zone II or III, still assuming that $L$ is \emph{not} a URE: we have
\[
  3k \leq n \leq 6k
\]
or, and the $(i,j)$ coordinates,
\begin{equation} \label{eq:zones II-III}
  i \leq 2j + \theta \textand j \leq 2i + \frac{n_0}{2} - \theta.
\end{equation}
We claim that there are rings for this pair $(i,j)$ if and only if $L \in \L_{e - 2k + \floor{\frac{n}{3}}}$, which may be also written as $k_0 \leq 2k - \frac{n - 2}{3}$ or in $(i,j)$ coordinates as
\begin{equation} \label{eq:zone supp}
  j \leq 2i + 1 - \theta.
\end{equation}
Assume first that $k_0 > 0$, that is, $L$ has splitting type $1^3$. Then the structure of $L$, and the fact that $\xi$ is a \emph{depth element} (a generator of the trace ideal) imply that $v_L(\xi) = 2 - \theta$, and hence that $v_K(a) = 1 - \theta$. Now it is easy to show that
\[
  v(\Phi_{\OO_L}(\xi)) = v(a\pi^{3i} + bu \pi^{2i} + cu^2 \pi^i + du^3) = \min\{3i + 1 - \theta, 3 v(u) + \theta \}
\]
because the sum is dominated by its first term if $\pi^{i+1-\theta} | u$ and by its last term otherwise. So a necessary condition for there to be rings is that $i + j \leq 3i + 1 - \theta$, which is equivalent to \eqref{eq:zone supp}. If this condition holds, then the $\xi$-condition simply becomes
\begin{equation} \label{eq:vu}
  v(u) \geq \frac{i + j - \theta}{3}
\end{equation}
and we get
\[
  q^{j - \ceil{\frac{i+j-\theta}{3}}} = q^{\floor{\frac{2j + \theta - i}{3}}} = q^{\floor{\frac{n}{3} - k}}
\]
solutions.

If $L$ has one of the other splitting types, then our task is simplified by the facts that $k_0 = 0$ and $n_0 = \theta \in \{0,1\}$. Note that the second inequality of \eqref{eq:zones II-III} implies \eqref{eq:zone supp}, so we are only trying to prove that there \emph{are} solutions in this case. If \eqref{eq:vu} does \emph{not} hold, then the $du^3$ term dominates in $\Phi_{\OO_L}(\xi)$ and we do not get a solution. If \eqref{eq:vu} holds, we leave it to the reader to check the inequalities that imply (even without knowing anything about $a$, $b$, and $c$) that $\pi^{i+j}$ divides each term of $\Phi_{\OO_L}(\xi)$. So we get the same number of solutions as in the preceding case.

Lastly, we must address the exceptional case that $L = K[\sqrt[3]{\pi}]$ is a URE. We take $\xi = \sqrt[3]{\pi}$ and $\eta = (\sqrt[3]{\pi})^2$, which are both traceless; and we have the explicit index form
\[
  \Phi_{\OO_L}(x\xi + y\eta) = x^3 - \pi y^3.
\]
This resembles the index form for a ramified $L$ with $\theta = 1$, and analogously to that case, we compute that rings appear only for the pairs $(i,j)$ with
\begin{equation} \label{eq:URE zone}
  i \leq 2j + 1 \textand j \leq 2i,
\end{equation}
each such $(i,j)$ yielding $q^{\floor{\frac{2j + \theta - i}{3}}}$ solutions. Now we come to the least satisfying part of the bijection. In the absence of a distinguishing $k$ (since all these rings are $(3)$-traced), we we simply have to place these rings into the families $\F_{n,k}$ of zone III so that the discriminant valuations and thicknesses match up. There is a unique choice:
\begin{align*}
  n &= n_0 + 2i + 2j = 3e + 2i + 2j + 2 \\
  k &= \floor{\frac{n}{3}} - \floor{\frac{2j + 1 - i}{3}}.
\end{align*}
We leave it to the reader that this establishes a bijection between the pairs $(i,j)$ in the region \eqref{eq:URE zone} with the pairs $(n,k)$ in zone III. The discriminant and thickness are correct by construction, completing the proof.
\end{proof}

\begin{rem}
This lemma also yields a second proof of the number of $t$-traced rings of discriminant $\mm_K^n$ in a cubic algebra (Theorem \ref{thm:traced_count}).
\end{rem}

\begin{rem}
The method of the above proof can also be adapted to the tame case.
\end{rem}

\section{Non-natural weightings}
\label{sec:non-natural}

Now that we know that the integral models $V_\tt$, $V_{3\tt^{-1}}$ of binary cubic forms are naturally dual, we can further look for duals for \emph{non-natural} weightings. This has applications to counting cubic rings satisfying local conditions. We restrict ourselves to primes not dividing $3\infty$.

For simplicity we work over $\ZZ$, though the techniques extend. Denote by $M_D$ the group $\ZZ/3\ZZ$ with Galois action given by the quadratic character corresponding to $\QQ(\sqrt{D})$. Denote by $V(R)$ the space of binary cubic forms over a ring $R$.

As in Section \ref{sec:composed}, if $W : V(\OO_{\AA_\QQ}) \to \CC$ is a locally constant weighting invariant under $\GL_2(\OO_{\AA_\QQ})$, we denote by $h(D, W)$ the number of $\GL_2 \ZZ$-classes of binary cubic forms over $\ZZ$ of discriminant $D$, each form $\Phi$ weighted by
\[
\frac{W(\Phi)}{\lvert \Stab_{\GL_2\ZZ} \Phi \rvert}.
\]
If $W = \prod_p W_p$ is a product of local weightings, then our local-to-global reflection engine (Theorems \ref{thm:main_compose} and \ref{thm:main_compose_multi}) produces identities relating different $h(D, W)$, if we can find a dual for each $W_p$.

\subsection{Local weightings given by splitting types}

Let $\sigma \in \{111,12,3,1^21, 1^3, 0\}$ be one of the six \emph{splitting types} a binary cubic form can have at a prime. Let
\[
  T(\sigma) = T_p(\sigma) : V(\ZZ_p) \to \{0, 1\}
\]
be the selector that takes the value $1$ on binary cubic forms of splitting type $\sigma$. Then the associated weighted local orbit counter
\[
  g_D\(T_p(\sigma)\) : H^1\(\QQ_p, M_D\) \to \NN
\]
attaches to each cubic algebra $L$ of discriminant $K(\sqrt{D})$ its number of orders of discriminant $D$ and splitting type $\sigma$.

There is another construction of interest to us. If $a \in \QQ$, then the varieties
\begin{equation}\label{eq:vars_for_Z}
  V_\ZZ(D = D_0) \textand V_\ZZ(D = a^2 D_0)
\end{equation}
do not in general look alike. However, their base-changes to $\QQ$ are isomorphic, being related by any $g \in \GL_2(\QQ)$ of determinant $a$. Hence the two varieties \eqref{eq:vars_for_Z} can be viewed as two integral models for $V_\QQ(D = D_0)$. Coupled with Theorem \ref{thm:main_compose_multi}, this viewpoint is very flexible. We denote by $Z_p$ the transformation that applies
\[
  \begin{bmatrix}
    1/p & \\
      & 1
  \end{bmatrix}
\]
to the vectors of an integral model of $V_\QQ$ and conjugates $\G$ accordingly. Observe that
\[
  g_D\left(T_p(\sigma) Z_p^n\right) = g_{D p^{-2n}}\left(T_p(\sigma)\right),
\]
and similarly for global class numbers. It is not hard to see that $h(D, W)$ is meaningful for any $W$ in the $\QQ$-algebra generated by the $T_p(\sigma)$'s and the $Z_p$'s for all $p$.

Over $\ZZ_p$, we still have $Z_p$, and we sometimes omit the subscript, as every $Z_\ell$ with $\ell \neq p$ has no effect on the integral model. We define $Z = Z_\pi$ for integral models over a general local field similarly.

\begin{lem} \label{lem:disc red} Let $K$ be a local field, $\ch k_K \neq 3$, and let $D \in \OO_K\bs \{0\}$. Then the weightings
\[
  T(1^3)Z \textand 2 \cdot T(111) - T(3)
\]
are dual with duality constant $1$; that is, the associated local orbit counters satisfy
\begin{equation} \label{eq:disc red}
  \hat{g}_{-3 \pi^2 D}(1^3) = 2 g_{D}(111) - g_D(3).
\end{equation}
\end{lem}
\begin{proof}
The right-hand side of \eqref{eq:disc red} can be written as
\[
  \size{H^0(M_D)} \cdot \1_{0} - \1_{H^1_\ur},
\]
so it suffices to show that the left-hand side is the Fourier transform of this, namely
\[
  1 - \1_{H^1_\ur} = \1_{H^1_{\ram}}.
\]
Look at binary cubic forms $f(x,y)$ of splitting type $1^3$ and discriminant $-3 \pi^2 D$. Changing coordinates, we can assume
\[
  f(x,y) = a x^3 + b x^2 y + c x y^2 + d y^3 \equiv x^3 \mod \pi.
\]
Then note that $\disc f \equiv - 27a^2d^2 \equiv -27d^2$ mod $\pi^3$, so the only way that $f$ can have discriminant $-3 \pi^2 D$ is if $\pi^2 \nmid d$. Then $f$ is an Eisenstein polynomial, the index form of a maximal order in a totally ramified extension $L$. Hence the weighting counting such $f$ is $\1_{H^1_{\ram}}$, as desired.
\end{proof}

Plugging this, together with the natural duality of Theorem \ref{thm:O-N_cubic_local} at the other primes, into Theorem \ref{thm:main_compose} yields results such as the following:
\begin{thm}\label{thm:1^3}
    Let $p \in \ZZ$ be a prime, $p \neq 3$. For all integers $D$ such that $p\nmid D$,
    \begin{align}
    \frac{1}{c_\infty} h_3(-27 p^2 D, T_p(1^3)) &= 2 h(D, T_p(111)) - h(D, T_p(3)) \label{eq:1^3 1} \\
    c_\infty h(p^2 D, T_p(1^3)) &= 2 h_3(-27 D, T_p(111)) - h_3(-27 D, T_p(3)) \label{eq:1^3 2}
    \end{align}
    where $c_\infty = 3$ for $D > 0$, $c_\infty = 1$ for $D < 0$.
\end{thm}

\subsection{Discriminant reduction}
\label{sec:disc_red}
This can be used to improve a step that often occurs in arithmetic statistics, namely the production of \emph{discriminant-reducing} identities that express the number of forms with certain non-squarefree discriminant in terms of lower discriminants.

For $N$ a positive integer, let $ h(D, R_N) $ be the number of classes of binary cubic forms of discriminant $D$, each weighted not only by the reciprocal of its number of automorphisms but also by its \emph{number of roots} in $\PP^1(\ZZ/N\ZZ)$. Equivalently, consider the natural congruence subgroup
\[
\GGamma^0(N) = \left\{\begin{bmatrix}
  a & b \\
  c & d
\end{bmatrix} \in \GL_2(\ZZ) : b \equiv 0 \mod N\right\},
\]
and let $h(D, R_N)$ be the number of $\GGamma^0(N)$-orbits of cubic $111N$-forms (integral forms with a marked root) of discriminant $D$ over $\ZZ$, each weighted by the reciprocal of its stabilizer in $\GGamma^0(N)$. (If $N > 1$, it is easy to prove that these stabilizers are trivial.) If $3 \nmid N$, denote by $h_3(D, R_p)$ the analogous weighted count of $133N$-forms. For a prime $N = p$, we have
\[
h(D, R_p) = 3 h\(D, T_p(111)\) + h\(D, T_p(12)\) + 2h\(D, T_p(1^21)\) + h\(D, T_p(111)\) + (p + 1)h\(D, T_p(0)\).
\]
This enables us to state succinctly the following theorem.
\begin{thm}[\textbf{discriminant reduction}]
Let $p \neq 3$ be a prime, and $D$ an integer divisible by $p^2$. Then
\[
  h(D) = h\(\frac{D}{p^2}, R_p\) + h\(\frac{D}{p^4}\) - h\(\frac{D}{p^4}, R_p\)
    + \frac{1}{c_\infty}\(2 h_3\(\frac{-27 D}{p^2}, T_p(111)\) - h_3\(\frac{-27 D}{p^2}, T_p(3)\)\),
\]
where $c_\infty = 3$ for $D > 0$, $c_\infty = 1$ for $D < 0$.
\end{thm}
\begin{proof}
The cubic rings $C$ counted by the left-hand side can be divided into maximal and nonmaximal at $p$. If $C$ is maximal at $p$, then $C \tensor_\ZZ \ZZ_p$ is the ring of integers of a totally tamely ramified cubic extension of $\QQ_p$ and $p^2 \parallel D$. By Theorem \ref{thm:1^3}, such rings are counted by the last term.

If $C$ is nonmaximal at $p$, then $C$ sits with $p$-power index inside an overring $C'$. By considering $C' + p C$, we can take an inclusion $C \subset C'$ of one of the following forms:
\begin{itemize}
  \item $C$ has index $p$ in a $C' = C_1$ of discriminant $D/p^2$. Here the index form of $C_1$ must have a marked root modulo $p$ so that the transformation
  \[
    \Phi_{C_1}(x,y) = a x^3 + b x^2 y + c x y^2 + d y^3
    \longmapsto \Phi_C(x,y) = p^2 a x^3 + p b x^2 y + c x y^2 + \frac{d}{p} y^3
  \]
  keeps the form integral. This accounts for the term $h(D/p^2, R_p)$.
  \item $C$ has index $p^2$ in a $C' = C_2$ of discriminant $D/p^4$ with $C_2/C \isom (\ZZ/p\ZZ)^2$. This requires that the index form of $C$ have \emph{content} divisible by $p$; we have
  \[
    \Phi_{C} = \frac{1}{p} C_2.
  \]
  This accounts for the term $h(D/p^4)$.
\end{itemize}
Observe that $C_2$ is unique if it exists. A choice of $C_1$ corresponds to a choice of multiple root of $\Phi_C$, which is unique if $\Phi_C$ is nonzero modulo $p$. Thus, the only chance of overcounting occurs when a $C$ admits both a $C_2$ and one or more $C_1$'s. The $C_1$'s are all the subrings of index $p$ in $C_2$ and thus correspond to the roots of $\Phi_{C_2}$ modulo $p$. So we subtract $1$ (more precisely, $1/\size{\Aut C_2}$) for each root of a form $\Phi_{C_2}$ counted in the term $h(D/p^4)$. That is, we subtract $h(D/p^4, R_p)$, yielding the claimed total.
\end{proof}

More generally, we can reduce at multiple primes at once. Let $T_p^{\max} : V(\OO_{\AA_\QQ}) \to \ZZ$ be the selector for rings maximal at $p$.

\begin{thm}[\textbf{discriminant reduction}]
Let $q = p_1\cdots p_r$ be a squarefree integer, $3 \nmid q$. If $D$ is a nonzero integer divisible by $q^2$, then for any $t < q$,
\begin{align*}
  h(D) &= \sum_{\substack{q = q_1q_2q_3 \\ q_1 \leq t}} h\(\frac{D}{q_2^2 q_3^4}, \prod_{p|q_1} T_p^{\max} \prod_{p|q_2} R_p \prod_{p|q_3} (1 - R_p)\) + {} \\
  & \quad + \frac{1}{c_\infty} \sum_{\substack{q = q_1q_2q_3 \\ q_1 > t}} h_3 \(\frac{-27 D}{q_1^2 q_3^2}, \prod_{p|q_1} \1_{p^2\parallel D}(R_p - 1) \prod_{p|q_3} \1_{p^2 \parallel D}(1 - R_p) \) \\
  \intertext{and}
  h_3(-27 D) &= \sum_{\substack{q = q_1q_2q_3 \\ q_1 \leq t}} h_3\(\frac{-27 D}{q_2^2 q_3^4}, \prod_{p|q_1} T_p^{\max} \prod_{p|q_2} R_p \prod_{p|q_3} (1 - R_p)\) + {} \\
  & \quad + {c_\infty} \sum_{\substack{q = q_1q_2q_3 \\ q_1 > t}} h_3 \(\frac{D}{q_1^2 q_3^2}, \prod_{p|q_1} \1_{p^2\parallel D}(R_p - 1) \prod_{p|q_3} \1_{p^2 \parallel D}(1 - R_p) \).
\end{align*}
\end{thm}
\begin{rem}
If we take $t = \sqrt{q}$, we find that all discriminants appearing are at most $-27D / q$.
\end{rem}
\begin{proof}
Since a ring of discriminant $D$ is maximal or nonmaximal at each of the primes dividing $q$, we have
\[
  h(D) = \sum_{q_1 q_2' = q} h\(D, \prod_{p|q_1} T_p^{\max}(1^3) \prod_{p|q_2'} T_p^{\nonmax}\),
\]
where $T_p^{\max}(1^3) = T_p^{\max} \cdot T_p(1^3)$ and $T_p^{\nonmax} = 1 - T_p^{\max}$ (as is natural).
We transform each term in one of two ways, depending on whether $q_1 \leq t$.

If $q_1 \leq t$, we simply replace each $T_p^{\nonmax}$ by $Z_p R_p + Z_p^2 (1 - R_p)$ by the method of the preceding theorem, which respects local conditions at other primes. We get a sum
\begin{align*}
  & h\(D, \prod_{p|q_1} T_p^{\max}(1^3) \prod_{p|q_2'} T_p^{\nonmax}\) \\
  &= h\(D, \prod_{p|q_1} T_p^{\max}(1^3) \prod_{p|q_2'} (Z_p R_p + Z_p^2 (1 - R_p))\) \\
  &= \sum_{q_2' = q_2q_3}  h\(\frac{D}{q_2^2 q_3^4}, \prod_{p|q_1} T_p^{\max} \prod_{p|q_2} R_p \prod_{p|q_3} (1 - R_p)\).
\end{align*}

If $q_1 > t$, we reflect. A dual of $T_p^{\max}$, when restricted to discriminants $D$ that are divisible by $p^2$, is $\1_{p^2 \parallel D} Z_p (R_p - 1)$ by Lemma \ref{lem:disc red}. Hence a dual of $T_p^{\nonmax} = 1 - T_p^{\max}$ on the same discriminants is $1 + \1_{p^2 \parallel D} Z_p (1 - R_p)$. Applying the reflection theorem,
\begin{align*}
  & h\(D, \prod_{p|q_1} T_p^{\max}(1^3) \prod_{p|q_2'} T_p^{\nonmax}\) \\
  &= \frac{1}{c_\infty} h_3\(-27D, \prod_{p|q_1} \1_{p^2\parallel D} Z_p(R_p - 1) \cdot \prod_{p|q_2'} \(1 + \1_{p^2 \parallel D} Z_p (1 - R_p)\)\) \\
  &= \frac{1}{c_\infty} \sum_{\substack{q_2' = q_2q_3}} h_3 \(\frac{-27 D}{q_1^2 q_3^2}, \prod_{p|q_1} \1_{p^2\parallel D}(R_p - 1) \prod_{p|q_3} \1_{p^2 \parallel D}(1 - R_p) \).
\end{align*}
Summing over $q_1$ yields the first identity. The second is proved in the same way.
\end{proof}

\subsection{Subring zeta functions}\label{sec:subring}
Fix a local field $K$. Instead of restricting ourselves to local weightings taking values in $\NN$ or $\CC$, consider the following (generalized) local weighting:
\begin{align*}
  \eta_D : H^1(K, M_D) &\to \ZZ\laurent{Z} \\
  L &\mapsto \sum_{k \in \ZZ} g_{\pi^{2k}D}(L) \cdot Z^k.
\end{align*}
Note that the sum is a Laurent series since the discriminant $\pi^{2k}D$ is nonintegral for $k$ sufficiently negative. This $\eta_D$ is, up to renormalizing, (a local factor of) the \emph{subring zeta function} that plays a central role in the study of Shintani zeta functions in works such as Datskovsky and Wright \cite{DW2} and Nakagawa \cite{Nakagawa}. In like manner we can define 
\[
  \eta_D(\sigma, L) = \sum_{k \in \ZZ} g_{\pi^{2k}D}(\sigma, L) \cdot Z^k \in \ZZ\laurent{Z},
\]
a partial subring zeta function that picks out the subrings of splitting type $\sigma$. We can even define
\[
  \eta_D(\sigma Z^n) = \eta_{\pi^{-2n} D}(\sigma) = Z^n \cdot \eta_D(\sigma),
\]
which explains the use of the same symbol $Z$ for the discriminant-shift operator and the formal variable in the subring zeta function.

A formula for the subring zeta function $\eta_{\disc L}(L)$, without splitting-type selector, is computed by Datskovsky and Wright \cite{DW2} and put into a more explicit form by Nakagawa:
\begin{thm}[Datskovsky--Wright; Nakagawa] \label{thm:subring_zeta}
The subring zeta function $\eta_{\disc L}(L)(Z = p^{-s})$ is given by
\[
  \eta_{\disc L}(L) = \frac{F}{(1-Z)(1 - q Z^3)}
\]
where $F$ is a polynomial depending on the splitting type $\sigma(\OO_L)$ as follows.
\begin{equation}
\begin{tabular}{cc}
$\sigma_\chi$ & F \\ \hline
$111$ & $(1 + Z)^2$ \\
$12$ & $1 + Z^2$ \\
$3$ & $1 - Z + Z^2$ \\
$1^21$ & $1 + Z$ \\
$1^3$ & $1$
\end{tabular}
\end{equation}
\end{thm}
\begin{proof}
Although Datskovsky--Wright \cite{DW2} and Nakagawa (\cite{Nakagawa}, Lemma 3.2; see Lemma 3.6 for the simplification method) work only with cubic extensions of $\QQ$, their method applies to this case. Alternatively, it follows immediately from Theorem \ref{thm:traced_count}.
\end{proof}
\begin{rem}
When working over the field $K = \QQ$, the local subring zeta functions at each prime $p$ form an Euler product expansion of the Dirichlet series
\[
\frac{\zeta_{L}(s)}{\zeta_{L}(2s)}\zeta(2s)\zeta(3s-1),
\]
where $\zeta$ and $\zeta_{L}$ are the Riemann and Dedekind zeta functions, respectively. In principle, Theorem \ref{thm:traced_count} allow us to write each Shintani zeta function $\xi_{K, \sigma, \tt}$ in Definition \ref{defn:Shintani} as an infinite sum of Euler products, one factor for each cubic \'etale algebra $L/K$, as was done for $\tt = (1)$ by Datskovsky and Wright (\cite{DW2}; see also \cite{NoEuler}). 
\end{rem}

We turn to the computation of $\eta_{\disc L}(\sigma, L)$, in which subrings are filtered by splitting type. Happily the answers are not too hard to deduce from Theorem \ref{thm:subring_zeta}. The easiest cases are $\sigma = 111$, $12$, and $3$, which only occur in maximal orders: thus for these three values of $\sigma$,
\[
  \eta_{\disc L}(\sigma, L) = \begin{cases}
    1 & \text{if } \sigma(\OO_L) = \sigma \\
    0 & \text{otherwise.}
  \end{cases}
\]
Next we compute $\eta_{\disc L}(1^21, L)$.
\begin{lem} \label{lem:orders of type 1^21}
\[
  \eta_{\disc L}(1^21, L) = \frac{G}{1-Z}
\]
where $G$ is a polynomial depending on the splitting type $\sigma_\chi$ as follows.
\begin{equation}
\begin{tabular}{cc}
$\sigma_\chi$ & G \\ \hline
$111$ & $3 Z$ \\
$12$ & $Z$ \\
$3$ & $0$ \\
$1^21$ & $1$ \\
$1^3$ & $0$
\end{tabular}
\end{equation}
\end{lem}
\begin{proof}
A ring $C$ of splitting type $1^21$ has its corresponding cubic form (in suitable coordinates) congruent to a multiple of $xy^2$ modulo $p$. Write
\[
\Phi_C(x,y) = pax^3 + pbx^2y + cxy^2 + pdy^3
\]
where $p\nmid c$. The ring $C$ is non-maximal iff $p|a$, in which it is of index $p$ in a unique overring $C'$ whose cubic form
\[
\Phi_{C'}(x,y) = \frac{a}{p}x^3 + bx^2y + p c x y^2 + p^2dy^3
\]
has a distinguished \emph{simple} root, here at $[0:1]$. Conversely, a ring with a distinguished simple root has one subring of index $p$ and splitting type $1^21$. Now the only non-maximal rings whose cubic forms have a simple root are themselves of splitting type $1^21$: thus the number of these is constant at index $p, p^2, \ldots.$ The value of the constant is the number of simple roots of the maximal order. Thus we have the desired claim.
\end{proof}

As for the subrings of splitting type $0$, they are $\OO_K + \pi C$ for each subring $C$, and so
\[
  \eta_D(0, L) = \eta_{D/\pi^2}(L) = Z \eta_{D}(L)
\]
for all $D$ and $L$. Finally, $\eta_D(1^3, L)$ can be computed by subtracting off all the other splitting types from $\eta_D(L)$. We do not need the explicit value in this case.

Before computing the Fourier transforms of the $\eta_D(\sigma, L)$, it is helpful to compute them for some simpler local weightings.

\begin{lem} \label{lem:td}
Let $p \neq 3$ be a prime.
For $\sigma \in \{111, 12, 3, 1^21, 1^3\}$ and $D \in K^\cross$, let $t_D(\sigma) : H^1(M_D) \to \NN$ be the local weighting given by
\[
  t_D(L) := \begin{cases}
    1 & \sigma(L) = \sigma \\
    0 & \text{otherwise.}
  \end{cases}
\]
Let $T_D$ be the $\QQ$-linear span of the $t_D(\sigma)$. Then the Fourier transform yields an isomorphism
\begin{equation}
\afterhat : T_{-3D} \to T_{D}
\end{equation}
given explicitly by
\begin{enumerate}[$($a$)$]
  \item \label{td:ur} $(t_{-3D}(111) + t_{-3D}(12) + t_{-3D}(3))\afterhat = (t_{-3D}(111) + t_{-3D}(12) + t_{-3D}(3))$
  \item \label{td:121} $t_{-3D}(1^21)\afterhat = t_D(1^21)$
  \item \label{td:2a} $(2t_{-3D}(111) - t_{-3D}(3))\afterhat = t_D(1^3)$
  \item \label{td:2b} $t_{-3D}(1^3)\afterhat = 2t_{D}(111) - t_{D}(3)$
  \item \label{td:cases} $\displaystyle (t_{-3D}(12))\afterhat = \begin{cases}
    t_{D}(12) & q \equiv 1 \mod 3 \\
    t_{D}(111) + t_{D}(3) & q \equiv 2 \mod 3
  \end{cases}$
\end{enumerate}
where $q = \size{k_K}$.
\end{lem}
\begin{proof}
We may assume that $D = D_0 \in \OO_K$ is a fundamental discriminant. Parts \ref{td:ur} and \ref{td:121} reduce to the self-orthogonality of $H^1_{D,\ur}$ for $D \in \OO_K^\cross$ and $D \in \pi \OO_K^\cross$, respectively. Parts \ref{td:2a} and \ref{td:2b} were proved after showing that both sides of Lemma \ref{lem:disc red} concern only $D \in \OO_K^\cross$. Finally, part \ref{td:cases} reduces to part \ref{td:ur}, using the fact that $3 \in \OO_K^\cross$ is a square if and only if $q \equiv 1$ mod $3$.
\end{proof}

\begin{thm}\label{thm:custom}
The weightings $\eta_D$, $\eta_D(111)$, $\eta_D(12)$, $\eta_D(3)$, $\eta_D(1^21)$ span $T_D\laurent{Z} := T_D \tensor_\QQ \QQ\laurent{Z}$. Their Fourier transforms are given by
\begin{enumerate}[$($a$)$]
    \item $\hat\eta_{-3D} = \eta_{D}$;
    \item If $q \equiv 1$ mod $3$, then
    \[
      (\eta_{-3D}(111) + \eta_{-3D}(3))\afterhat = \eta_{D}(111) + \eta_{D}(3) \textand \widehat{\eta}_{-3D}(12) = \eta_{D}(12);
    \]
    \item If $q \equiv 2$ mod $3$, then
    \[
      (\eta_{-3D}(111) + \eta_{-3D}(3))\afterhat = \eta_{D}(12) \textand \widehat{\eta}_{-3D}(12) = \eta_{D}(111) + \eta_{D}(3);
    \]
    \item $(2\eta_{-3D}(111) - \eta_{-3D}(3))\afterhat =
    \frac{1}{Z}\big[(1 - Z)(1 - qZ^3)\eta_{D} + (Z - 1)(\eta_{D}(111) + \eta_{D}(12) + \eta_{D}(3)) + Z^2(2 \eta_{D}(111) - \eta_{D}(3)) - (1+Z)^2 \eta_{D}(1^21) \big]$;
    \item $\hat \eta_{-3D}(1^21) = \frac{1}{1 + Z}\big[(1 - Z)(1 - qZ^3)\eta_{D} + (Z - 1)(\eta_{D}(111) + \eta_{D}(12) + \eta_{D}(3)) + Z(Z-1)(2 \eta_{D}(111) - \eta_{D}(3)) - Z(1+Z) \eta_{D}(1^21)\big]$.
\end{enumerate}
\end{thm}
\begin{proof}
We may assume that $D = D_0 \in \OO_K$ is a fundamental discriminant. We need the following result, which is useful in its own right:
\begin{lem} \label{lem:cubic_discs}
Let $D = \disc \OO_L$ be the discriminant of a maximal cubic order over a local field $K$, $\ch k_K \neq 3$. Write $D = \pi^{2k} D_0$ where $D_0$ is the associated fundamental discriminant. Then
\[
  k = \begin{cases}
    1 & \sigma(\OO_L) = 1^3 \\
    0 & \text{otherwise.}
  \end{cases}
\]
\end{lem}
\begin{proof}
If $L$ is unramified, then $\pi \nmid D$ and $k = 0$. If $\sigma(L) = 1^21$, then
\[
  D = \disc L = \disc (K \cross K(\sqrt{D_0})) = D_0.
\]
So we are left with the case that $\sigma(L) = 1^3$. Note that $L/K$ is tamely ramified so $D = \pi^2 D'$, $\pi \nmid D'$. We need to prove that $D'$ is a discriminant, which is only nontrivial when $\ch k_K = 2$. Let $f(x) = a x^3 + \pi bx^2 + \pi cx + \pi d$ be an Eisenstein polynomial for $L$. Then
\[
  D' = \pi^2 b^2 c^2 - 4 \pi a c^3 - 4 \pi^2 b^3 d - 27 a^2 d^2 + 18 \pi a b c d \equiv (a d - \pi b c)^2 \mod 4,
\]
so $D' = D_0$.
\end{proof}

This allows us to change the scaling in Lemma \ref{lem:orders of type 1^21} (on $\eta_D(1^21)$) and Theorem \ref{thm:subring_zeta} (on $\eta_D$) from $\eta_{\disc L}$ to $\eta_{D}$, multiplying the weighting by $Z$ in the case that $\sigma(\OO_L) = 1^3$. Also, trivially $\eta_D(\sigma) = t_D(\sigma)$ for unramified $\sigma$. So we have five equations
\begin{align*}
  \eta_D(111) &= t_D(111) \\
  \eta_D(12) &= t_D(12) \\
  \eta_D(3) &= t_D(3) \\
  \eta_D(1^21) &= \frac{3t_D(111) + t_D(12) + t_D(1^21)}{1 - Z} \\
  \eta_D &= \frac{(1 + Z)^2t_D(111) + (1 + Z^2)t_D(12) + (1 - Z + Z^2)t_D(3) + (1 + Z)t_D(1^21) + Zt_D(1^3)}{(1 - Z)(1 - q Z^3)}.
\end{align*}
We see that there is an invertible transition matrix between the sets
\[
  \{t_D(111), t_D(12), t_D(3), t_D(1^21), t_D\} \textand \{\eta_D(111), \eta_D(12), \eta_D(3), \eta_D(1^21), \eta_D\}
\]
is invertible (we cannot call these ``bases,'' because $t_D(\sigma)$ can vanish for certain $D$). Accordingly, we can rewrite Lemma \ref{lem:td} in terms of the $\eta$'s and get the identities claimed in the theorem.
\end{proof}

\begin{prob}
Can the rich structure found in this section be carried out, to some extent, when $\ch k_K = 3$? For instance, $\hat V(\FF_3)$, the space of integral $1331$-forms modulo $3993$-forms, has six $\GL_2(\FF_3)$-orbits, the analogues of splitting types, and it is natural to wonder whether the Fourier transform relates them to the six splitting types on $V(\FF_3)$.
\end{prob}

\subsection{Invariant lattices at \texorpdfstring{$2$}{2}} \label{sec:invarlat2}

The foregoing investigation also allows us to extend the work of Ohno and Taniguchi \cite{10lat} on extending O-N to counting binary cubic forms satisfying certain congruence conditions at $2$. Assume that $q = 2$. Recall from Section \ref{sec:cubic} the five primitive invariant lattices $\Lambda_i \subseteq V(\OO_K)$ of binary cubic forms:
\begin{align*}
  \Lambda_1 &= \{f(x,y) = ax^3 + bx^2y + cxy^2 + dy^3 : a,b,c,d \in \OO_K\} \\
  \Lambda_2 &= \{f \in \Lambda_1 : \sigma(f) \in \{0,12\} \} \\
  \Lambda_3 &= \{f \in \Lambda_1 : \sigma(f) \in \{0,111,3\} \} \\
  \Lambda_4 &= \{f \in \Lambda_1 : \sigma(f) \in \{0,111,12,1^3\} \} \\
  \Lambda_5 &= \{f \in \Lambda_1 : \sigma(f) \in \{0,111\} \}
\end{align*}
Let
\begin{align*}
  \lambda_{i}(D) : H^1(K, M_D) \to \ZZ\laurent{Z} \\
  L \mapsto \sum_{k \in \ZZ} g_{\pi^{2k}D}(\Lambda_i, L) \cdot Z^k
\end{align*}
be the analogue of the subring zeta function counting only those rings whose corresponding form is in $\Lambda_i$. Each $\lambda_{i}(D)$ is a linear combination of the appropriate $\eta_{\sigma,D}$, and the transition matrix is again invertible. (No reason is known, beyond pure coincidence, why the number of invariant lattices over $\ZZ_2$ should equal $5$, the number of independent $\eta_{\sigma,D}$.) So we get:
\begin{cor}[\textbf{Local O-N for general invariant lattices;} cf.~\cite{10lat}, Theorems 1.2, 1.3 and 1.4] \label{cor:10lat_local}
Let $K$ be an unramified extension of $\QQ_2$.
The Fourier dual of each $\lambda_{i}(-27D)$ lies again in the span of the $\lambda_{i}(D)$, explicitly:
\begin{align*}
  \hat \lambda_1(-27D) &= \lambda_{1}(D) \\
  \hat \lambda_2(-27D) &= \lambda_{3}( D) \\
  \hat \lambda_3(-27D) &= \lambda_{2}( D) \\
  \hat \lambda_4(-27D) &= \frac{1}{3Z}\big[(-8Z^3 + 6Z^2 - Z)\lambda_1(D) + (2Z - 1)\lambda_2(D) + (-8Z^2 + 4Z)\lambda_3(D) \\
  &\qquad {}+ (-4Z^2 + 1)\lambda_4(D) + (16Z^2 - 1)\lambda_5(D)\big] \\
  \hat \lambda_5(-27D) &= \frac{1}{3Z}\big[(-2Z^3 + 3Z^2 - Z)\lambda_1(D) + (2Z - 1)\lambda_2(D) + (-2Z^2 + Z)\lambda_3(D)\\
  &\qquad {}+ (-Z^2 + 1)\lambda_4(D) + (4Z^2 - 1)\lambda_5(D)\big]
\end{align*}
\end{cor}
\begin{proof}
The $\lambda_i(D)$ are simply the $\eta_D(\sigma)$ in disguise: for instance,
\[
  \lambda_2(D) = \eta_D(0) + \eta_D(12),
\]
and so on. Rewriting the results of Theorem \ref{thm:custom} in terms of the $\lambda_i$ and plugging in $q = 2$ proves the theorem.
\end{proof}
The five equalities in this corollary can be viewed as local reflection theorems in the sense of Theorem \ref{thm:main_compose_multi}, relating different integral models of the same composed variety on each side.

We can now get global results in great generality.
\begin{thm}[\textbf{O-N for general invariant lattices}]\label{thm:invarlat}
Let $K$ be a number field. Let $\aa$ be an ideal of $K$, and let $\Lambda$ be a $\G_a = \SL(\OO_K \oplus \aa)$-invariant lattice of full rank in the space $V(K)$ of binary cubic forms. Let $\tt$ be the trace ideal of $\Lambda$, that is, the unique $\tt \mid 3$ such that $\Lambda_\pp \isom \V_{(1, \tt)}(\pp)$ at every prime $\pp \mid 3$ (after identifying $\G_\aa$ with $\SL_2$, which we can do after localization). For $D \in K^\cross$, denote by $h_\Lambda(D)$ the number of $\G_a$-orbits of binary cubic forms of discriminant $D$, each orbit weighted by the reciprocal of its stabilizer. 

Then:
\begin{enumerate}
  \item If $\Lambda$ is of type $\Lambda_1$, $\Lambda_2$, or $\Lambda_3$ at every prime whose residue field is $\FF_2$, then there is a lattice $\Lambda^*$ invariant under $\G_{\aa\tt^{-3}}$ such that we have the global reflection theorem: for each $D \in K^\cross$,
  \begin{equation}
    h_{\Lambda}(D) = \frac{3^{\#\{v|\infty : D \in (K_v^\cross)^2\}}}{N_{\OO_K/\ZZ}(\tt)} h_{\Lambda^*}(-27D).
  \end{equation}
\item In general, there is a family $\Lambda^*_1, \ldots, \Lambda^*_m$ of lattices, all equal away from $2$, each $\Lambda^*_j$ invariant under a group $\G_{\aa\ss_j\tt^{-3}}, \ldots$ where $\ss_j$ is an ideal having nonzero valuation only at primes dividing $2$, and a global reflection theorem of the shape
\[
  h_{\Lambda}(D) = \frac{3^{\#\{v|\infty : D \in (K_v^\cross)^2\}}}{N_{\OO_K/\ZZ}(\tt)}
  \sum_{j = 1}^m c_j h_{\Lambda^*_j}(-27D),
\]
where the $c_i \in \QQ$ depend only on $\Lambda$.
\end{enumerate}
\end{thm}
\begin{proof}
Theorem \ref{thm:Osborne} limits the lattices we must consider. If $\Lambda$ is of type $\Lambda_1$ at each prime dividing $2$, we have $\Lambda = \cc \V_{\aa,\tt}(\OO_K)$. But since
\[
  \cc \oplus \cc \isom 1 \oplus \cc^2,
\]
there is an isomorphism of integral models
\[
  \(\cc \V_{\aa,\tt}, \G_\aa\) \isom \(\V_{\aa\cc^{-2}, \tt}, \G_{\aa\cc^{-2}}\).
\]
So we can take
\[
  \Lambda^* = \cc \V_{\aa \tt^{-3}, 3\tt^{-1}},
\]
and the desired reflection theorem follows from Theorem \ref{thm:O-N_traced}.

The remaining cases can be solved with a bit of fiddling at $2$. Let $\bar\Lambda \supseteq \Lambda$ be the lattice that sits over $\Lambda$ at each $\pp \mid 2$ as $\Lambda_1$ sits over the relevant $\Lambda_i$, $1\leq i \leq 5$, and let $\bar\Lambda^*$ be the corresponding reflection lattice. All $\Lambda^*_i$ will look like $\Lambda^*$ away from $2$.

Locally at each $\pp \mid 2$, we construct a collection of $\Lambda^*_{j, \pp}$, $1\leq j \leq m_\pp$ as follows: for each term $Z^k \lambda_i(D)$ in Corollary \ref{cor:10lat_local}, take the lattice
\[
  \Lambda^*_{j, \pp} =
  \begin{bmatrix}
    \pi^{k} &   \\
      & 1
  \end{bmatrix}\Lambda_i,
\]
for which $\G_{\pi^{-1}}$-orbits of discriminant $D$ correspond to $\SL_2$-orbits of discriminant $D/\pi^{2k}$ in $\Lambda_i$. Then Corollary \ref{cor:10lat_local} appears as a local reflection theorem
\[
  \hat g_{\Lambda_{\pp}(D)} = \sum_{j=1}^{m_\pp} c_{j,\pp} \cdot g_{\Lambda^*_{j, \pp}(D)}.
\]
for the integral models $\Lambda_{\pp}$ and $\Lambda^*_{j, \pp}$. (Strictly speaking, these are not truly integral models, inasmuch as the $\SL_2$-invariance of $\Lambda_2,\ldots,\Lambda_5$ is not given by an algebraic integrality; but we get the same local and global class numbers by going up to a lattice of type $\Lambda_1$ and imposing non-natural weights to pick out the appropriate splitting types.)
  
We then apply Theorem \ref{thm:main_compose_multi} to these local reflection theorems and, as usual, the ones relating $\bar\Lambda$ and $\bar\Lambda^*$ away from $2$. The resulting reflection theorem involves the global integral models given by gluing the $\Lambda^*_{j, \pp}$ at each $\pp\mid 2$ in all possible ways. When $\Lambda$ is of type $\Lambda_1, \Lambda_2, \Lambda_3$ at each $\pp \mid 2$, there is only one $\Lambda^*$, and $\ss$, which arises from the $Z$-operators in Corollary \ref{cor:10lat_local}, disappears.
\end{proof}

\subsection{Binary cubic forms over \texorpdfstring{$\ZZ[1/N]$}{Z[1/N]}}
For $N$ a squarefree integer, it is natural to ask what happens if we invert finitely many primes and count binary cubic forms of discriminant $D \neq 0$ over $\ZZ[1/N]$, up to the action of the relevant group $\SL_2(\ZZ[1/N])$. There are still only finitely many for each degree, owing to Hermite's theorem on the finiteness of the number of number fields with prescribed degree and set of ramified primes.

Note that $D$ matters only up to multiplication by the squares in $\ZZ[1/N]^\cross$; hence we can restrict our attention to $D\in \ZZ$ that are \emph{fundamental} at each prime $p \mid N$. (If $p \neq 2$, this means that $p^2 \nmid D$. If $p = 2$, this means that $D \equiv 1 \mod 4$ or $D \equiv 8, 12 \mod 16$. However, we allow $D$ to be non-fundamental at primes not dividing $N$.)

We do not have O-N for forms over $\ZZ[1/N]$ in the same formulation as over $\ZZ$. Nevertheless, the other side of the reflection theorem is noticeably not too complicated.

\begin{thm}
Let $N$ be a squarefree integer.

For $0 \neq D \in \ZZ[1/N]$, let $h_{\ZZ[1/N]}(D)$ be the number of $\SL_2(\ZZ[1/N])$-orbits of integral binary cubic forms over $\ZZ[1/N]$, each weighted by the reciprocal of its stabilizer in $\SL_2(\ZZ[1/N])$. If $3 \nmid N$, define $h_{3, \ZZ[1/N]}(D)$ to be the same count, counting only $1331$-forms (that is, forms whose middle two coefficients belong to the ideal $3\ZZ[1/N] \subsetneq \ZZ[1/N])$.

Now let $D$ be a discriminant that is fundamental at all primes dividing $N$. For each $p \mid N\infty$, let
\[
  c_{D,p} = \begin{cases}
    3 & \text{if } D \in \(\ZZ_p^\cross\)^2 \\
    1 & \text{if } D \text{ is a non-square unit modulo $p$}, \\
    1/2 & \text{if } p \mid D.
  \end{cases}
\]
and let
\[
  c_{D, N\infty} = \prod_{p \mid N\infty} c_{D, p}.
\]
Then:
\begin{enumerate}[$($a$)$]
  \item If $3 \nmid N$, then
  \begin{align*}
    h_{\ZZ[1/N]}(D) &= \frac{c_{-27D,N\infty}}{3} \cdot h_{3}(-27 D, R_N) \\
    h_{3, \ZZ[1/N]}(-27 D) &= c_{D, N\infty} \cdot h(D, R_N).
  \end{align*}
  \item If $3 \mid N$, then
  \[
    h_{\ZZ[1/N]}(-27 D) = 3 c_{D, N\infty} \cdot h_{3}(D, R_N).
  \]
\end{enumerate}
\end{thm}

\begin{proof}
We take the same composed variety $(V, \Gamma)$ of binary cubic forms as before. However, we take integral models $(\V^{(i)}, G^{(i)})$ that are not even over the same ring of integers $\OO_\QQ^{(i)}$!

On the left-hand side, we take the scheme $\V^{(1)}$ of binary cubic $111N$- or $133N$-forms of discriminant $D$ or $-27 D$. This does not admit an algebraic action of $\SL_2\ZZ$, but it \emph{does} admit an algebraic action of $G^{(1)} = \GGamma^0(N)$.

On the right-hand side, we take the scheme $\V^{(2)}$ of binary cubic forms over $\ZZ[1/N]$ of discriminant $D$, or $1331$-forms as appropriate, with the natural action of $G^{(2)} = \SL_2$ over $\ZZ[1/N]$.

It is evident that the global class numbers of these integral models match the quantities studied in the theorem. The checking of most of the conditions of Theorem \ref{thm:main_compose} is routine, so we content ourselves with checking the local duality.

When $p \nmid N$, the integral model is identical to that used for O-N, so we already have the needed duality with an appropriate duality constant $c_{D, p}$ or $c_{-27D, p}$. This includes the infinite prime, at which the duality constant $c_{D, \infty}$ tracks the sign of $D$ as in O-N.

When $p \mid N$, the computation of the local class numbers is not difficult:
\begin{itemize}
  \item As to $\V^{(1)}$, we look for forms of discriminant $D$ with a marked root modulo $p$. We first observe that forms corresponding to nonzero cohomology classes are not counted, because they either have splitting type
  \begin{itemize}
    \item $ (3) $, and have no roots modulo $p$, or
    \item $ (1^3) $, and have discriminant non-fundamental at $p$, by Lemma \ref{lem:cubic_discs}.
  \end{itemize}
  So $g^{(1)} : H^1(\ZZ_p, M') \to \NN$ is a scalar multiple of $\1_{0}$, nonzero because the split ring $\ZZ_p \cross \ZZ_p[(D + \sqrt{D})/2]$ has an index form with a root.
  \item As to $\V^{(2)}$, since the completion of $\ZZ[1/N]$ at $p$ is $\QQ_p$, the local orbit counter counts cosets in $\SL_2(\QQ_p) \bs \SL_2(\QQ_p)$ that keep a certain form $f$ ``integral'' over $\QQ_p$. There is obviously only one such coset, regardless of the cohomology class of $f$, so $g^{(2)} : H^1(\ZZ_p, M) \to \NN$ is identically $1$.
\end{itemize}
  It remains only to compute the duality constant.
  \begin{itemize}
    \item If $\QQ_p[\sqrt{D}]$ is split, then there are three roots of $f(x,y) = xy(x+y)$ to mark, but they all wind up equivalent. So $g^{(1)} = \1_{0}$, but because of the $\size{H^0} = 3$ in the scaling of the Fourier transform, we need to insert a factor of $c_{D,p} = 3$.
    \item If $\QQ_p[\sqrt{D}]$ is inert, there is only one root to mark, and $\size{H^0} = 1$, so $c_{D,p} = 1$.
    \item If $\QQ_p[\sqrt{D}]$ is ramified, then we can mark either the single or the double root modulo $p$. These are non-equivalent $\GGamma^0(p)(\ZZ_p)$-orbits inside the same $\SL_2(\QQ_p)$-orbit, so $g^{(1)} = 2 \cdot \1_{0}$ and $\size{H^0} = 1$, so $c_{D,p} = 1/2$. (If we modified the theorem by counting $111N$-forms whose third coefficient $c$ is coprime to $N$, another family stable under $\GGamma^0(N)$, then we would be forced to mark the simple root, and this factor of $1/2$ would disappear.)
  \end{itemize}
Multiplying the duality constants obtained completes the proof. 
\end{proof}

\begin{prob}
Does the integral model of $111N$-forms of non-fundamental discriminant $D$ have a natural dual? The first step in answering this is to check whether the Fourier transform of its local orbit counter takes nonnegative values. 
\end{prob}

\part{Reflection theorems: quartic rings and related objects}\label{part:quartic}

\section{Reflection for 2-adic quartic orders, and applications}
Analogously to the cubic case, the reflection theorem that we state and prove is going to swap $\tt$-traced and $2\tt^{-1}$-traced orders. It is not the most general reflection theorem that one can try to state: see Section \ref{sec:doubly_traced} below.

Fix a nondegenerate cubic ring $C$ over a Dedekind domain $\OO_K$. We can also fix a basis
\[
  C = \OO_K \oplus \OO_K \xi \oplus \aa\eta,
\]
making the index form $\Phi_C(x\xi + y\eta)$ a cubic in the $\V_{\aa,1}$ of Theorem \ref{thm:O-N_traced}. However, none of our work will depend on this basis.

We can then look at the scheme $\V_{\tt, C}$ of pairs of ternary quadratic forms
\[
(\A,\B) : \OO_K \cross \OO_K \cross \aa \rightrightarrows \OO_K \cross \aa,
\]
that is, pairs of $3\times 3$ symmetric matrices with entries in the ideals
\[
  \(\begin{bmatrix}
    (1) & 2^{-1}\tt & 2^{-1}\tt\aa^{-1} \\
       & (1) & 2^{-1}\tt\aa^{-1} \\
           & & \aa^{-2}
  \end{bmatrix},
  \begin{bmatrix}
    \aa & 2^{-1}\tt\aa & 2^{-1}\tt \\
    & \aa & 2^{-1}\tt \\
     & & \aa^{-1}
  \end{bmatrix}
\)
\]
satisfying the four equations
\[
\det(\A x + \B y) = \Phi_C(x,y)
\]
asserting that $(\A,\B)$ parametrizes a quartic ring $L$ that is $\tt$-traced with reduced resolvent $C$. This $\V_{\tt,C}$, together with the natural action of the group $\G = \GL(\OO_K \cross \OO_K \cross \aa)$, is an integral form of the composed variety $(V, \GL_3)$ of pairs of ternary quadratic forms over $K$. We assert that the integral models
\[
\V_{\tt,C} \textand \V_{2\tt^{-1},C}
\]
are naturally dual at all finite places.

\begin{nota}
  Here and in the sequel, we use the label ``Theorem*'' to denote a theorem proved with the following caveats:
  \begin{enumerate}[$($1$)$]
    \item Resolvents that are wildly ramified at a $2$-adic place are excluded.
    \item In general, the theorem depends on a Monte Carlo verification of a rational algebraic identity (as we will explain). However:
    \item The results for $K/\QQ_2$ unramified (e.g.~the cases over $\ZZ$) are known unconditionally.
    \item The results for the reduced resolvent being maximal $(C_\tt = \OO_R)$ are also known unconditionally.
  \end{enumerate}
\end{nota}

\begin{thmstar}[\textbf{``Local Quartic O-N''}] \label{thm*:quartic_local}
Let $K$ be a nonarchimedean local field and $C$ an order in an \'etale algebra $R$ that is not wildly ramified over $K$. For $\tau$ a divisor of $2$, let $\V_{\tau,C}$ be the integral model parametrizing $(\tau)$-traced orders with reduced resolvent $C$. Then $\V_{\tau,C}$ and $\V_{2\tau^{-1}, C}$ are naturally dual with duality constant $q^{2v_K(\tau)}$; in order words, the associated local orbit counters
\[
  g_{\tt, C} : H^1(K, M_R) \to \NN
\]
satisfy the local reflection theorem 
\[
  \hat{g}_{\tt,C} = \Size{\OO_K/\tau\OO_K}^2 \cdot  g_{2\tt^{-1},C}.
\]
\end{thmstar}

At the infinite places, we no longer have natural duality. (We did not have this problem in the cubic case because $H^1(\RR, M)$ is trivial for $\size{M}$ odd.)

Therefore, call a quartic algebra $L/K$ over a number field \emph{nowhere totally complexified (ntc)} if there is no real place $\pp$ of $K$ such that $L_\pp \isom \CC \cross \CC$. This is equivalent to the cohomology element $\sigma_L \in H^1(K, (\ZZ/2\ZZ)^2)$ being trivial at all infinite places. Then the local specifications $\1_{\{0\}}$ for ntc quartic algebras and $1$ for all quartic algebras are mutually dual, provided that one inserts the correct scale factor.

\begin{thmstar}[\textbf{``Quartic O-N''}] \label{thm*:O-N quartic}
  Let $K$ be a number field. Let $C$ be an order in a cubic $K$-algebra $R$, and let $\tt \subseteq \OO_K$ be an ideal such that $\tt \mid (2)$. Let $h(C, \tt)$ count the number of $\tt$-traced quartic rings with reduced resolvent $C$, respectively, each weighted by the reciprocal of its number of resolvent-preserving automorphisms. Let $h^{\ntc}(C, \tt)$ count the subset of the foregoing that are ntc, weighted in the same way. Then
  \[
  h(C, \tt) = \frac{N(\tt)^2}{2^{r_\infty}} \cdot h^{\ntc}(C, 2\tt^{-1}),
  \]
  where $r_\infty$ is the number of real places of $K$ over which $R$ is not totally real plus twice the number of complex places of $K$.
\end{thmstar}
\begin{proof}
  We apply Theorem \ref{thm:main_compose} to the composed varieties $\V^{(1)} = \V_{\tt,C}$ and $\V^{(2)} = \V_{2\tt^{-1},C}$ just defined, with the following local weightings $w_v^{(i)}$:
  \begin{itemize}
    \item At finite $v \nmid 2$, we take $w_v^{(i)} = 1$, which are mutually dual with duality constant $1$ by Theorem \ref{thm:O-N_quartic_local_tame}.
    \item At $v = \pp|2$, we take $w_v^{(i)} = 1$, which are mutually dual with duality constant
    \[
      N(\pp)^{2v_\pp(\tt)}
    \]
    by Theorem* \ref{thm*:quartic_local}.
    \item At complex $v$, we take $w_v^{(i)} = 1$, which are mutually dual with duality constant
    \[
    \hat f_v = \frac{1}{\size{H^0(K_v, M_R)}} = \frac{1}{4}.
    \]
    \item At real $v$ for which $R_v \isom \RR \cross \CC$, we take $w_v^{(i)} = 1$, which are mutually dual with duality constant
    \[
    \hat f_v = \frac{1}{\size{H^0(K_v, M_R)}} = \frac{1}{2}.
    \]
    \item At real $v$ for which $R_v \isom \RR \cross \RR \cross \RR$, we take $w_v^{(1)} = 1$ and $w_v^{(2)} = \1_0$, the selector for rings that are not totally complex at $v$. The corresponding duality constant is
    \[
    \frac{\Size{H^1(K_v, M_R)} }{\size{H^0(K_v, M_R)}} = \frac{4}{4} = 1.
    \]
  \end{itemize}
  The product of all duality constants is thus
  \[
  \prod_{\pp|2} N(\pp)^{2v_\pp(\tt)} \cdot \prod_{K_v \isom \CC} \frac{1}{4} \cdot \prod_{R_v \isom \RR \cross \CC} \frac{1}{2}
  = \frac{N(\tt)^2}{2^{r_\infty}},
  \]
  as desired.
\end{proof}

Although we have been counting quartic rings by resolvent, the corresponding result where we count by discriminant follows quickly. We present the reflection theorem in two forms, one dealing with $2\times 3\times 3$ symmetric boxes, the other with quartic rings (which are in bijection only in the case of rings of content $1$):

\begin{thm}\label{thm:O-N_2x3x3}
Let $K$ be a number field. Denote by $\V_{\tt,\aa}(\OO_K)$ the space of pairs of ternary quadratic forms
\[
(\A,\B) : \OO_K \cross \OO_K \cross \aa \rightrightarrows \OO_K \cross \aa,
\]
that are $\tt$-traced in the sense that the entries belong to the ideals
\[
\(\begin{bmatrix}
  (1) & 2^{-1}\tt & 2^{-1}\tt\aa^{-1} \\
  & (1) & 2^{-1}\tt\aa^{-1} \\
  & & \aa^{-1}
\end{bmatrix},
\begin{bmatrix}
  \aa & 2^{-1}\tt\aa & 2^{-1}\tt \\
  & \aa & 2^{-1}\tt \\
  & & (1)
\end{bmatrix}
\).
\]
It has a natural action of the group
\[
  \G_{\aa} = \SL(\OO_K \oplus \OO_K \oplus \aa) \cross \SL(\OO_K \oplus \aa)
\]
that preserves discriminant.
Denote by $h_\tt(\D)$ the number of $\G_{\tt,\aa}$ orbits of pairs of ternary quadratic forms having discriminant $\D = (\aa, D)$, each orbit weighted by the reciprocal of the order of its stabilizer in $\G_{\aa}$. Denote by $h^{\ntc}_\tt(\D)$ the number of such orbits (weighted in the same way) which are ntc, in the sense that at each real place of $K$, the conics $\A$ and $\B$ have a common point in $\RR\PP^2$. Then for all discriminants $\D$ prime to $2$, 
\[
h_\tt(\tt^8\D) = \frac{N(\tt)^2}{2^{r_\infty}} \cdot h^{\ntc}_{2\tt^{-1}}(256\tt^{-8}\D),
\]
where $r_\infty$ is the number of real places of $K$ at which $D < 0$ plus twice the number of complex places of $K$.
\end{thm}
\begin{proof}
We sum the preceding theorem over all cubic rings $C$ of discriminant $\D$, weighting each $C$ by the reciprocal of the number of \emph{orientation-preserving} automorphisms of $C$, which is the stabilizer of the corresponding form in $\SL(\OO_K \oplus \aa)$. It is easy to see that each orbit is counted the number of times it appears in the theorem. Because the reduced discriminant is prime to $2$ (a needed condition to avoid involving wildly ramified resolvents in the sum), we can state the theorem unconditionally.
\end{proof}

\begin{thm} \label{thm:O-N_quartic_by_disc}
Let $K$ be a number field, and let $\D = (\aa, D)$ be a discriminant.
Denote by $h^\circ_\tt(\D)$ the number of $\tt$-traced quartic rings $\OO$ over $\OO_K$ having discriminant $\D$, each $\OO$ weighted by $1/\size{\Aut_K(\OO)}$. Denote by $h_\tt^{\circ,\ntc}(\D)$ the number of such that are ntc, weighted in the same way. Then for all discriminants $\D$ prime to $2$, 
\begin{equation}\label{eq:O-N_by_disc}
h^\circ_\tt(\tt^8\D) = \frac{N(\tt)^2}{2^{r_\infty}} \cdot h^{\circ,\ntc}_{2\tt^{-1}}(256\tt^{-8}\D),
\end{equation}
where $r_\infty$ is the number of real places of $K$ at which $D < 0$ plus twice the number of complex places of $K$.
\end{thm}

\begin{proof}
In the previous theorem, we studied $h_\tt(\tt^8\D)$, which can be interpreted as the number of quartic rings $\OO$ equipped with a resolvent $C$ and an \emph{orientation,} that is, an identification $\Lambda^4 \OO \isom \aa$ for which the discriminant is $\tt^8\D$. Every quartic ring admits two orientations (there are $\size{\OO_K^\cross}$-many identifications $\Lambda^4 \OO \isom \aa$, but all but one and its negative yield a $D$ scaled by a different square of a unit). So $\frac{1}{2}h_\tt(\tt^8\D\D)$ is the number of resolvents $(\OO, C, \Theta, \Phi)$ of discriminant $\tt^8\D$, up to isomorphism, each weighted by the reciprocal of its number of automorphisms.

Let $h^{1}_\tt(\tt^8\D)$ be the number of quartic rings $\OO$ of discriminant $(\aa,D)$ with $\tt$-traced content $1$, weighted by $1/\lvert \Aut \OO\rvert$. This is related to  If $\OO$ is an ntc quartic ring of discriminant $\tt^8\D$ having some $\tt$-traced content $\cc$, then $\OO = \OO_K + \cc \OO'$, where $\OO'$ has $\tt$-traced content $1$, discriminant $(\aa\cc^{-3}\tt^8,D) = \cc^{-6}\tt^8\D$, and the same automorphism group as $\OO$. Thus
\[
  \frac{1}{2}h_\tt(\tt^8\D) = \sum_{\cc^3 \mid D\aa^2} h_\tt^{1}(\tt^8\D\cc^{-6}).
\]
On the other hand, the number of resolvents of $\OO$ depends on the $\tt$-traced content $\cc$ (Proposition \ref{prop:traced}\ref{traced:count}): it is
\[
  \sigma_1(\cc) = \sum_{\dd \mid \cc} N_{K/\QQ}(\dd).
\]
These are resolvents \emph{as maps out of $\OO$} (as pointed out in \cite{ORings}, end of Section 8), which is the correct manner of counting to make
\[
  h_\tt^\circ(\tt^8\D) = \sum_{\cc^3 \mid D\aa^2} \sigma_1(\cc) h_\tt^{1}(\cc^{-6}\tt^8\D).
\]
We can now write $h_\tt^\circ$ in terms of $h_\tt$:
\begin{align*}
  h_\tt^\circ(\tt^8\D) &= \sum_{\cc^3 \mid D\aa^2} \sigma_1(\cc) h_\tt^{1}(\cc^{-6}\tt^8\D) \\
  &= \sum_{\cc^3 \mid D\aa^2} \sum_{\bb \mid \cc} N(\bb) h_\tt^{1}(\cc^{-6}\tt^8\D) \\
  &= \sum_{\bb^3 \mid D \aa^2} \sum_{\cc'^3 \mid D \aa^2\bb^{-3}} N(\bb) h_\tt^{1}(\bb^{-3}\cc'^{-6}\tt^8\D) \\
  &= \frac{1}{2}\sum_{\bb^3 \mid D \aa^2} N(\bb) h_\tt (\bb^{-3}\tt^8\D).
\end{align*}
Transforming both sides of \eqref{eq:O-N_by_disc} in this manner reduces it to Theorem \ref{thm:O-N_2x3x3}.
\end{proof}

\subsection{Results on binary quartic forms} \label{sec:bq}
We can also derive a reflection theorem about binary quartic forms, which correspond (via a completely general construction for binary $n$-ic forms) to a certain subclass of quartic rings. This subclass was identified explicitly by Wood:
\begin{thm}[\cite{WoodBQ}, Theorem 1.1]\label{thm:WoodBQ}
  There is a natural, discriminant preserving bijection between the set of
  $\GL_2(\ZZ)$-equivalence classes of binary quartic forms and the set of isomorphism classes of
  pairs $(Q, C)$ where $Q$ is a quartic ring and $C$ is a monogenized cubic resolvent of Q (where
  isomorphisms are required to preserve the generator of $C$ modulo $\ZZ$).
\end{thm}
\begin{proof}
  Regarding pairs $(Q, C)$ of quartic rings as pairs $(A,B)$ of $3 \times 3$ symmetric matrices via Bhargava's parametrization, we send a form $\Phi(x,y) = ax^4 + bx^3y + cx^2y^2 + dxy^3 + ey^4$ to
  \[
  (A_0, B) = \(
  \begin{bmatrix}
    & & 1/2 \\
    & -1 & \\
    1/2 & &
  \end{bmatrix},
  \begin{bmatrix}
    a & b/2 & \\
    b/2 & c & d/2 \\
    & d/2 & e
  \end{bmatrix}\).
  \]
  The distinguished generator arises because the resolvent form $g(x,y) = 4 \det (A_0 x - B y)$ is monic, since $\det A_0 = 1/4$. Further details will be found in \cite{WoodBQ}.
\end{proof}
To apply this theorem, we need to know the number of automorphisms of the quartic ring corresponding to a given form:
\begin{lem} \label{lem:bq aut}
  In this bijection, the group of resolvent-preserving automorphisms of a quartic ring is in natural isomorphism with the stabilizer (in $\PGL_2(\ZZ)$) of the corresponding form.
\end{lem}
\begin{proof}
  The conclusion follows easily from the method of proof of the preceding theorem. By  \cite{WoodBQ}, Theorem 2.5, we can choose bases for $Q$ and $C$ so that the corresponding pair of ternary quadratic forms, has the form $(A_0, B)$ above. A resolvent-preserving automorphism is a change of variables $h \in \SL_3(\ZZ)$ that preserves both $A_0$ and $B$. By \cite{WoodBQ}, Lemma 3.2, $h = \epsilon_{A_0}(\tilde h)$ lies in the image of the map
  \begin{align*}
    \epsilon_{A_0} : \GL_2(\ZZ) &\to \SL_3(\ZZ) \\
    \begin{bmatrix}
      a & b \\
      c & d
    \end{bmatrix}
    &\mapsto
    \frac{1}{ad - bc}
    \begin{bmatrix}
      a^2 & ab      & b^2 \\
      2ac & ad + bc & 2bd \\
      c^2 & cd      & d^2
    \end{bmatrix}
  \end{align*}
  By \cite{WoodBQ}, Theorem 3.1, $h$ preserves $B$ if and only if $\tilde h$ preserves the binary quartic form $f$. Moreover, it is easy to see that $\ker \epsilon_{A_0} = \pm 1$. This constructs the desired isomorphism.
\end{proof}

We fix a monic binary cubic form $g(x,y)$ and let $C = \ZZ[\xi]$ be the corresponding monogenized cubic ring, with generator $\xi$. Then in the notation of Theorem* \ref{thm*:O-N quartic},
\[
h\big(C,(1)\big) = \sum_{\substack{\text{quartic rings } Q \\ \text{with resolvent } C, \\ \text{up to $C$-isom}}} \frac{1}{\Size{\Aut_C Q}}
= \sum_{\substack{\text{binary quartics $f(x,y)$} \\ \text{with resolvent $g$,} \\ \text{up to $\PGL_2 \ZZ$}}} \frac{1}{\Size{\Stab_{\PGL_2 \ZZ}} f}.
\]

The quantity $h^{\ntc}\big(C,(2)\big)$ appearing on the opposite side of Theorem* \ref{thm*:O-N quartic} is not so straightforward to interpret. Here we are counting pairs $(A, B)$ of integer symmetric matrices with $\det (Ax - By) = g$ (with a certain condition at $\infty$), so we need to classify integer symmetric matrices $A$ with $\det A = 1$. There are, up to similarity, two:
\begin{lem}\label{lem:sym mat}
  Every integer symmetric matrix $A$ with $\det A = 1$ is similar to
  \[
  A_1 = \begin{bmatrix}
    & & 1 \\
    & -1 & \\
    1 & &
  \end{bmatrix}
  \textor
  I = \begin{bmatrix}
    1 & & \\
    & 1 & \\
    & & 1
  \end{bmatrix}.
  \]
\end{lem}
\begin{proof}
  Let $A$ be an integral symmetric matrix of determinant $1$. Look at the corresponding conic $\C$ defined by $x^\top A x = 0$. Note that for each rational prime $p \not\in \{2, \infty\}$, we have $p \nmid 4 = \det \C$, so $\C$ has good reduction to $\FF_p$: by the Chevalley-Warning theorem, $\C$ has an $\FF_p$-point and hence a $\QQ_p$-point. Then, by Hilbert reciprocity, there are only two possibilities for the isomorphism type of $\C$ over $\QQ$:
  \begin{itemize}
    \item If $\C$ has an $\RR$-point, then $\C$ also has a $\QQ_2$-point and hence (by the Hasse principle) a $\QQ$-point. By a $\GL_3\ZZ$-transformation, we set this point to $[1:0:0]$, the tangent line there to go through $[0:1:0]$, and then $A$ must take the form
    \[
    \begin{bmatrix}
      & & \pm 1 \\
      & -1 & a \\
      \pm 1 & a & b
    \end{bmatrix},
    \]
    which one easily sees is similar to $A_1$.
    \item If $\C$ has no $\RR$-points, then $A$ is positive (or negative) definite. It is well known that the only positive definite integral unimodular matrix of any rank is the identity. \qedhere
  \end{itemize}
\end{proof}

Hence $h^{\ntc}\big(C,(2)\big)$ decomposes into the pairs $(A,B)$ of ``type $A_1$'' and of ``type $I$'' according to the value of $A$ after an appropriate $\SL_3\ZZ$-transformation.

\subsubsection{Type \texorpdfstring{$A_1$}{A1}}
To understand pairs $(A_1,B)$, we capitalize on the fact that $A_1 \sim A_0$ over $\QQ$. Namely, the transformation
\[
T = \begin{bmatrix}
  1 & & \\
  & 1 & \\
  & & 2
\end{bmatrix}
\]
satisfies $TA_0T^{\top} = A_1$. Let
\[
(A_1,B) = \( A_1, \begin{bmatrix}
  a & b & c' \\
  b & c & d \\
  c' & d & e
\end{bmatrix} \).
\]
Then the pair
\[
\( T^{-1}A_1\( T^{-1} \)^{\top} , T^{-1} B \( T^{-1} \)^{\top} \) = (A_0, B')
\]
is determined up to $\SO(\QQ, A_0)$ by $(A_1, B)$, and the form
\[
f(x,y) = a x^4 + 2b x^3 y + (c + c') x^2 y^2 + d x y^3 + \frac{1}{4} e y^4
\]
is determined up to $\PSL_2(\QQ)$. This is a form of a peculiar shape, the \emph{$(1,2,1,1,\frac{1}{4})$-forms}. Since the Wood embedding is resolvent-preserving, we see that the resolvent of a $(1,2,1,1,\frac{1}{4})$-form is actually integral, which can also be deduced directly from the formula for the resolvent of a binary quartic.

The $(1,2,1,1,\frac{1}{4})$-forms do not naturally have an action by $\PGL_2\ZZ$, but rather by a group that we can reveal as $\SO(\ZZ, A_1)$:

\begin{lem}
  We have $\SO(\QQ, A_1) \isom \PGL_2(\QQ)$ via the isomorphism
  \[
  \epsilon_{A_1}\colon
  M = 
  \begin{bmatrix}
    a & b \\
    c & d
  \end{bmatrix} \mapsto T\epsilon_{A_0}(M)T^{-1}
  = \frac{1}{ad - bc}\begin{bmatrix}
    a^2 & a b & \frac{1}{2} b^2 \\
    2 a c & a d + b c & b d \\
    2 c^2 & 2 c d & d^2.
  \end{bmatrix}
  \]
  Under this map, the subgroup corresponding to $\SO(\ZZ, A_1)$ is $G = \GGamma^0(2) \sqcup \tau\GGamma^0(2)$, where
  \[
  \GGamma^0(2) = \left\{\begin{bmatrix}
    a & b \\
    c & d
  \end{bmatrix} \in \PGL_2(\ZZ) : b \equiv 0 \mod 2\right\}
  \]
  is a congruence subgroup, and
  \[
  \tau = \begin{bmatrix}
    & 2 \\
    1 &
  \end{bmatrix},
  \]
\end{lem}
\begin{proof}
  The first statement follows easily from considering the action of an element of $\SO(\QQ, A_1)$ on the locus of isotropic points for $A_1$, a conic in $\PP^2$ that is rationally isomorphic to $\PP^1$. Note that $\epsilon_{A_1}$ is compatible with the map $\epsilon_{A_0}$ found earlier: 
  \begin{equation}
    \epsilon_{A_1} = T \epsilon_{A_0} T^{-1}.
  \end{equation}
  As for the second statement, if
  \[
  M = \begin{bmatrix}
    a & b \\
    c & d
  \end{bmatrix} \in \PGL_2(\QQ)
  \]
  is given such that $\epsilon_{A_1}(M)$ is integral, we first multiply by $\tau$ if needed to make $v_2(\det M)$ even, and then scale $M$ so that $a,b,c,d$ are coprime integers. Then we argue that if a prime $p$ were to divide $ad-bc$, it must divide each of $a$,$b$,$c$,$d$ by the integrality of $\epsilon_{A_1}(M)$, which is a contradiction.
\end{proof}

We have now mapped each $\SL_3$-orbit of pairs of integral symmetric matrices of $A_1$-type to a $G$-orbit of binary quartic $(1,2,1,1,\frac{1}{4})$-forms; indeed, it is not hard to show that $G$ is in fact the subgroup of $\PGL_2(\QQ)$ that preserves the lattice of forms of this shape, and we have, by an argument similar to Lemma \ref{lem:bq aut},
\[
h^{A_1\text{-type}}\big(C,(2)\big) = \sum_{\substack{\text{quartic rings } Q \\ \text{with resolvent } C \\ \text{of $A_1$ type,} \\ \text{up to $C$-isom}}} \frac{1}{\Size{\Aut_C Q}}
= \sum_{\substack{(1,2,1,1,\frac{1}{4})\text{-forms } f(x,y) \\ \text{with resolvent $g$,} \\ \text{up to $G$}}} \frac{1}{\Size{\Stab_G f}}.
\]
\subsubsection{Type \texorpdfstring{$I$}{I}}
Following the same method, we can write
\[
h^{I\text{-type}}\big(C,(2)\big) = \sum_{\substack{\text{quartic rings } Q \\ \text{with resolvent } C \\ \text{of $A$ type,} \\ \text{up to $C$-isom}}} \frac{1}{\Size{\Aut_C Q}}
= \sum_{\substack{\text{$(I,B)$ with resolvent $g$} \\ \text{up to $\SO(\ZZ,I)$}}} \frac{1}{\Size{\Stab_{\SO(\ZZ,I)} B}}.
\]
At this point we make two striking observations:
\begin{itemize}
  \item The resolvent condition $\det(xI - yB) = g$ is equivalent to $B$ having \emph{characteristic polynomial} $g(x,1)$, so we have connected counting quartic rings to another classical problem, namely counting symmetric matrices of given characteristic polynomial;
  \item Since $\SO(\ZZ,I)$ is a finite group, isomorphic to $\S_4$ (in its representation as the group of rotations of a cube), there is no need to count \emph{orbits} of symmetric matrices; the matrices themselves will be finite in number.
\end{itemize}
Thus
\[
h^{I\text{-type}}\big(C,(2)\big) = \frac{1}{24} \Size{ \{B \in \Mat^{3\times 3} \ZZ : \charpoly(B) = g(x,1)\}}.
\]

\subsubsection{Conditions at \texorpdfstring{$\infty$}{infinity}}
The interpretations of class numbers of quartic rings that we have here developed can be modified to take into account local conditions or weightings at a prime. Here we only consider the prime at infinity.

Over $\RR$ there are only two nondegenerate cubic algebras, $\RR \cross \CC$ and $\RR \cross \RR \cross \RR$. If $C_g \tensor \RR \isom \RR \cross \CC$ (that is, $g$ has only one real root), then there is only one quartic algebra with resolvent $C_g$ up to $C_g$-isomorphism, so it does not make sense to impose local conditions at the infinite place. If, on the other hand, $C_g \tensor \RR \isom \RR \cross \RR \cross \RR$, then the three factors of $\RR$ are non-interchangeable, being labeled by the three real roots of $g$, and there are four non-$C_g$-isomorphic quartic algebras with resolvent $C_g$ (one isomorphic to $\RR \cross \RR \cross \RR \cross \RR$ and three to $\CC \cross \CC$), parametrized by the four Kummer elements $\delta \in C_g^{N=1}/(C_g^{N=1})^2.$ The following is not hard to verify:
\begin{lem}
  Let $g$ be a monic binary cubic form over $\RR$ whose dehomogenization has three real roots $\xi_1 < \xi_2 < \xi_3$. Identify the corresponding $\RR$-algebra $C_g$ with $\RR \cross \RR \cross \RR$ with the coordinates ordered so that $\xi \mapsto (\xi_1,\xi_2,\xi_3)$. Let $L/\RR$ be a quartic algebra with resolvent $C_g$. The corresponding pair $(A,B)$ of real symmetric matrices is related to the sign of the corresponding Kummer element $\delta \in C_g^{N=1}/(C_g^{N=1})^2$ in the following way:
  \begin{itemize}
    \item If $\sgn \delta = (+,+,+)$, then $(A,B)$ is of type $A_1$ and yields an indefinite binary quartic form with four real roots.
    \item If $\sgn \delta = (+,-,-)$, then $(A,B)$ is of type $A_1$ and yields a positive definite binary quartic form.
    \item If $\sgn \delta = (-,-,+)$, then $(A,B)$ is of type $A_1$ and yields a negative definite binary quartic form.
    \item If $\sgn \delta = (-,+,-)$, then $(A,B)$ is of type $I$.
  \end{itemize}
\end{lem}

\subsubsection{Statements of results}
We leave it to the reader to furnish the modifications of the condition at $\infty$ in the proof of Theorem* \ref{thm*:O-N quartic} to yield the following identities. Because $\ZZ$ is unramified at $2$, we can prove them unconditionally, but only for tamely ramified resolvent at present.
\begin{thm}[\textbf{Quartic O-N for binary quartic forms}]\label{thm:BQ}
  Let $g$ be a monic integral binary cubic form whose splitting field is unramified at $2$. Denote by $h(g)$ the number of integral binary quartic forms of resolvent $g$, up to $\PGL_2(\ZZ)$-equivalence and weighted by inverse of $\PGL_2(\ZZ)$-stabilizer. Denote by $h_4(g)$ the number of binary quartic $(1,2,1,1,\frac{1}{4})$-forms of resolvent $g$, up to $G$-equivalence and weighted by inverse of $G$-stabilizer. Denote by $s(g)$ the number of integral $3\times 3$ symmetric matrices of characteristic polynomial $g$. Then:
  \begin{itemize}
    \item If $\disc g < 0$, then
    \[
    2 h(g) = h_4(g).
    \]
    \item If $\disc g > 0$, then
    \begin{align*}
      h(g) &= 2 h_4^{\text{indef}}(g) \\
      h^{\text{indef or pos def}}(g) &= h_4^{\text{indef or pos def}}(g) \\
      h^{\text{indef or neg def}}(g) &= h_4^{\text{indef or neg def}}(g) \\
      24\big(h^{\text{indef}}(g) - h^{\text{def}}(g)\big) &= s(g)
    \end{align*}
    where the superscripts instruct one to count only forms satisfying the indicated condition at infinity, with the same weighting.
  \end{itemize}
\end{thm}
\begin{cor}\label{cor:BQ}
  Let $g$ be a monic integral binary cubic form with three real roots whose splitting field is unramified at $2$. Among integral binary quartics with resolvent $g$, at least half are indefinite when we weight by inverse size of $\PGL_2(\ZZ)$-stabilizer, with equality exactly when $g$ is not the characteristic polynomial of an integral $3\times 3$ symmetric matrix.
\end{cor}
We state these unconditionally because they apply only to the number field $K = \QQ$. We do not attempt to generalize to other number fields. While the quartic rings of types corresponding to soluble conics ($A_0$ and $A_1$ in our notation) continue to be connected to binary quartic forms, the number of insoluble types grows with the degree of $K$.

\subsection{The conductor property of the resolvent ring}
We conclude this part of the paper with a family of results which at first do not look at all like reflection theorems.

Let $R$ be an \'etale algebra over a local field $K$, and let $E/R$ be an abelian extension whose Artin map $\psi_{E/R} : R^\cross \to \Gal(E/R)$ vanishes on the base $K^\cross$. Call an order $\OO \subseteq R$ an \emph{admissible ring} for $E$ if $\psi(\OO^\cross) = 0$. Such rings exist (e.g{.} $\OO = \OO_K + \ff$ where $\ff \subseteq R$ is the conductor ideal) and are stable under passage to suborders. If there is a unique maximal admissible ring, we call it the \emph{conductor ring} of the extension $E/R$.

In like manner, we define admissible and conductor rings for an abelian extension $E/R$ over a global field $K$, if the Artin map $\psi_{E/R} : I(R, \mm) \to \Gal(E/R)$ vanishes on the ideals $I(K, \mm)$ of the base. By Lemma \ref{lem:ring_cl_fld}, if $E/R/K$ are fields, $\OO$ is an admissible ring of $E$ if and only if $E$ is contained in the ring class field of $\OO$.

In general, an extension $E/R$ can have multiple maximal admissible rings, and there is no reason for a conductor ring to exist.

\begin{examp}\label{ex:no cdr ring}
  Let $p \geq 5$ be a prime. Take $K = \QQ_p$, $R = \QQ_p^4$ and let $\chi : K^\cross \to \FF_p^\cross$ be any multiplicative homomorphism extending the natural projection from $\ZZ_p^\cross$. Define $\psi : R^\cross \to \FF_p^\cross$ by
  \[
  \psi(a;b;c;d) = \chi\left(\frac{a b}{c d}\right).
  \]
  This is the Artin map of a certain $\FF_p^\cross \isom \ZZ / (p-1)\ZZ$-torsor $E/R$. By construction, $\psi$ vanishes on $\QQ_p^\cross$, and the orders
  \begin{align*}
    \OO_1 &= \{(a;b;c;d) \in \ZZ_p^4 : a \equiv c, b \equiv d \mod p\} \\
    \OO_2 &= \{(a;b;c;d) \in \ZZ_p^4 : a \equiv d, b \equiv c \mod p\}
  \end{align*}
  are admissible rings for $\psi$. However, no ring strictly containing either $\OO_1$ or $\OO_2$ (of which there are very few) is an admissible ring for $\psi$; in particular, $\OO_1 \union \OO_2$ generates the whole $\OO_R$, which is certainly not an admissible ring for $\psi$. Thus $\psi$ has no conductor ring.
\end{examp}

However, in two special cases the conductor ring not only exists but has a striking characterization: it is the \emph{resolvent ring} of a certain maximal order. These cases are those of general cubic and quartic algebras.

\begin{prop}
  Let $L/K$ be a cubic \'etale algebra over a global or local field. Let $T$ be its quadratic resolvent torsor and $E = L T$ its $\S_3$-closure. Then the quadratic resolvent ring $\OO \subseteq T$ of $L$ is the conductor ring of $E/T$.
\end{prop}
\begin{proof}
  The global case reduces immediately to the local one. The Artin map $\phi_{E/T}$ vanishes on $K^\cross$ by the Galois symmetry of the situation (the same argument is carried out in a global context in Nakagawa \cite[p.~110]{Nakagawa}), so $L/K$ has admissible rings. Now the orders in $T$ are totally ordered: they are simply of the form $\OO = \OO_K + \pi^i \OO_T$ for $i \geq 0$. It is evident that the conductor ring of $E/T$ must be $\OO_K + \ff_{E/T} \OO_T$, which has discriminant
  \[
  \disc T \cdot \ff_{E/T}^2 = \disc T \cdot \disc(E/T).
  \]
  The proposition is now reduced to the identity
  
  \[
  \disc L = \disc T \cdot \disc(E/T).
  \]
  This is a form of the ``Brauer relation'' between the absolute discriminants of $L$, $T$, and $E$ and follows quickly from an Artin-conductor argument: see \cite{CohON}, equation (2.7).
\end{proof}

The above proof is not very deep and does not use reflection theorems at all. Let it be noted that over $K = \QQ$, a very similar result was proved, if not stated, by Hasse (\cite{Hasse}, table on p.~568) and forms a foundation to Nakagawa's proof of Ohno-Nakagawa \cite[Lemma 1.3]{Nakagawa}.

However, the quartic analogue of this statement, which we state in an identical way, is much deeper. Note that, in addition to the cubic resolvent $R$, a quartic \'etale algebra $L/K$ has a natural sextic resolvent $S$ coming from the map $\S_4 \to \S_6$ that sends a permutation of $\{1,2,3,4\}$ to the corresponding permutation of its $2$-element subsets. $S$ is naturally a quadratic \'etale extension of $R$ with the same Kummer element $\delta \in R^{N=1}/\(R^{N=1}\)^2$ that parametrized $L$ in \ref{thm:Kummer_new}.

\begin{thmstar}[\textbf{Conductor rings}]\label{thm*:cond_ring}
  Assume Theorem* \ref{thm*:quartic_local}.
  Let $L/K$ be a quartic \'etale algebra over a global number field or a $p$-adic field. Let $R$ and $S$ be its cubic and sextic resolvent algebras, respectively. Then the cubic resolvent ring $S_0 \subseteq R$ of $\OO_L$ is the conductor ring of $S/R$.
\end{thmstar}
\begin{proof}
  The global case reduces immediately to the local one. To see the vanishing of $\phi_{S/R}$ on $K^\cross$, let $S = R(\sqrt{\delta})$ where $S_{R/K}(\delta) = 1$. Then for $a \in K^\cross$,
  \[
  \phi_{S/R}(a) = \< \delta, a\> = \<a, \delta \> = \phi_{R(\sqrt{a})/R}(\delta) = \phi_{K(\sqrt{a})/K}\( S_{R/K}(\delta)\)
  = 1.
  \]
  
  The conjecture now has two parts:
  \begin{enumerate}[(a)]
    \item \label{it:cyes} $S_0$ is an admissible ring for $S/R$;
    \item \label{it:cno} Any admissible ring for $S/R$ is contained in $S_0$.
  \end{enumerate}
  As mentioned, this result does not on the surface look like a reflection theorem. But we will prove both \ref{it:cyes} and \ref{it:cno} using Theorem* \ref{thm*:quartic_local}.
  
  Let $g(L, S, \tt)$ denote the number of $\tt$-traced orders in $L$ with reduced resolvent $S$. Then Theorem* \ref{thm*:quartic_local} states that
  \begin{equation} \label{eq:xz}
    \hat g(L, S, (2)) = c \cdot g(L, S, (1))
  \end{equation}
  for the appropriate positive constant $c = \size{\OO_L / 2\OO_L}$. Now $g(L, S, (1))$ is the number of orders in $L$ of resolvent $S$. In particular, it is $1$ if $S = S_0$ and $0$ if $S \not\subseteq S_0$. On the other hand, $g(L', S, (2))$ can be interpreted as the number of ideals $I \subseteq R$ such that $(S, I, \delta')$ is balanced, and overall
  \begin{align}
    \hat g(L, S, (2)) &= \frac{1}{\size{H^0(M_R)}} \sum_{\delta' \in H^1(M_R)} \< \delta, \delta' \> g(L_{\delta'}, S, (2)) \nonumber \\
    &= \frac{1}{\size{H^0(M_R)}} \sum_{\substack{\delta' \in H^1(M_R) \\ (S, I, \delta') \text{ balanced}}} \< \delta, \delta' \> \nonumber \\
    &= \frac{1}{\size{H^0(M_R)}} \sum_{\substack{\delta' \in H^1(M_R) \\ (S, I, \delta') \text{ balanced}}} \phi_{S/R}(\delta'). \label{eq:xsum}\\
  \end{align}
  Assume for the sake of contradiction that $S = S_0$ is not an admissible ring for $S/R$. Then there exists $\epsilon \in S^\cross$ such that $\phi_{S/R}(\epsilon) = -1$. The rearrangement of terms $(S, I, \delta') \mapsto (S, I, \epsilon \delta')$ flips the sign of the sum, so $\hat g(L, S, (2)) = 0$, a contradiction, since $g(L, S, (1)) = 1$. This proves \ref{it:cyes}.
  
  Now assume for the sake of contradiction that there is an admissible ring $S \not\subseteq S_0$. Choose such an $S$ maximal for this property. Then divide the summands of \eqref{eq:xsum} into two cases:
  \begin{itemize}
    \item If $I$ is invertible in $S$, then $I = \alpha S$ for some $\alpha$, and $\delta = \alpha^2 \epsilon$ for some $\epsilon \in S^\cross$. These terms contribute
    \[
    \phi_{S/R}(\delta) = \phi_{S/R}(\epsilon) = 1,
    \]
    since $S$ is an admissible ring. There is at least one term of this type, namely $\delta = 1$, $I = S$.
    \item If $I$ is not invertible in $S$, then $\End I^2 = S' \supsetneq S$. (If we had $\End I^2 = S$, then by Lemma \ref{lem:sqrs are inv}, $I^2$ would be invertible in $S$ and then $I$ would also.) By maximality, $S'$ is \emph{not} an admissible ring for $S/R$ and there is an $\epsilon \in S'^\cross$ such that $\phi_{S/R}(\epsilon) = -1$. The rearrangement $(S, I, \delta') \mapsto (S, I, \epsilon \delta')$ permutes the terms with the same $S'$ and flips their signs. So the terms of this type contribute nil.
  \end{itemize}
  Overall, we get $\hat g(L, S, (2)) > 0$, a contradiction, since $g(L, S, (1)) = 0$. This proves \ref{it:cno}.
\end{proof}

\section{Tame quartic rings with non-split resolvent, by multijection}

In this section, we will adapt the methods of Section \ref{sec:cubic_tame} to prove Theorem* \ref{thm*:quartic_local} in the case that $K$ is tame (not $2$-adic) and $R \not\isom K \cross K \cross K$.

\subsection{Invertibility of ideals in orders}
We begin with a technical inquiry that has interest in its own right. It is well known that every $\ZZ$-lattice $\aa$ in a quadratic field is invertible with respect to \emph{some} order, namely its endomorphism ring $\End \aa$. In a cubic or higher-degree field this is not so. However, the following two lemmas will help us understand the structure of orders and ideals in such a setting.

\begin{lem} \label{lem:order is prod}
Let $K$ be a local field, and let $\OO$ be an order in a finite-rank \'etale algebra $L$ over $K$. Then there is a decomposition
\[
  L = L_1 \cross \cdots \cross L_s,
\]
each $L_i$ being the product of some field factors of $L$, with the following properties:
\[
  \OO = \OO_1 \cross \cdots \cross \OO_s
\]
is the product of orders in the $L_i$, and each $\OO_i$ has only a single prime $\pp_i$ above the valuation ideal $\pp$, so that every element of $\OO_i$ not lying in $\pp_i$ is a unit.
\end{lem}
\begin{proof}
Let $L = K_1 \cross \cdots \cross K_r$ be the field factor decomposition of $L$. Each $K_i$ is a local field; let $\pp_i$ be the pullback to $L$ of the valuation ideal of $K_i$. Define an equivalence relation on the $\pp_i$ by
\[
  \pp_i \sim \pp_j \iff \pp_i \intsec \OO = \pp_j \intsec \OO.
\]
Thus if $\pp_i \nsim \pp_j$, then there is an $\alpha \in \OO$ such that either
\begin{equation} \label{eq:nsim}
  \pp_i \nmid \alpha, \pp_j | \alpha \textor \pp_i | \alpha, \pp_j \nmid \alpha.
\end{equation}
Suppose for the moment that it is the first. We first claim that we can take $\alpha \equiv 1$ mod $\pp_i$. Note that $\OO_L/\pp_i$ is a finite field extension of $\OO_K/\pp$, so $\alpha$ satisfies a polynomial congruence
\[
  \alpha^m + u_{m-1}\alpha^{m-1} + \cdots + u_0 \equiv 0 \mod \pp_i, \quad u_i \in \OO_K, \quad u_0 \in \OO_K^\cross.
\]
Then 
\[
  \alpha' = -u_0^{-1}(\alpha^m + u_{m-1}\alpha^{m-1} + \cdots + u_1 \alpha)
\]
is $1$ mod $\pp_i$ and $0$ mod $\pp_j$. Also, note that we can switch $\alpha'$ with $1 - \alpha'$ to satisfy these congruences for any $i$ and $j$, regardless of which condition in \eqref{eq:nsim} held to begin with.

Fix $i$ and multiply the resulting values of $\alpha'$, which are $1$ mod $\pp_i$ (and hence $1$ mod any $\pp_j \sim \pp_i$) but $0$ mod $\pp_j$ for any chosen $\pp_j \nsim \pp_i$. We get a single $\alpha''$ such that
\[
  \alpha'' \equiv \begin{cases}
    1 \mod \pp_j, & \pp_j \sim \pp_i \\
    0 \mod \pp_j, & \pp_j \nsim \pp_i.
  \end{cases}
\]
As a final step, we can iterate the polynomial
\[
  f(\x) = \x^2(3-2x),
\]
which takes $\pp_i^m$ to $\pp_i^{2m}$ and $1 + \pp_i^m$ to $1 + \pp_i^{2m}$, and take the limit to derive that the idempotent $e = (e_j)_j \in \OO_L$, defined by
\[
  e_j = \begin{cases}
    1, & \pp_j \sim \pp_i \\
    0, & \pp_j \nsim \pp_i,
  \end{cases}
\]
lies in $\OO$ ($\OO$ is closed in the $\pp$-adic topology). These $e_j$'s form a set of orthogonal idempotents decomposing $\OO$ into a direct product of $\OO_i$, one for each $\sim$-equivalence class, that have the properties we seek.
\end{proof}
\begin{lem}\label{lem:inv=pri}
  Let $K$ be a local field, and let $\OO$ be an order in a finite-rank \'etale algebra $L$ over $K$. Then a fractional ideal of $\OO$ is invertible if and only if it is principal.
\end{lem}
\begin{rem}
  This implies that a if $\OO$ is an order in a finite-rank algebra over a Dedekind domain, then a fractional ideal of $\OO$ is invertible if and only if it is locally principal, where here ``locally'' denotes localization at each prime of $K$. Thus our statement and proof differ slightly from the corresponding statement in Neukirch \cite{Neukirch} (Theorem I.12.4), which is built by localization at the primes \emph{of $\OO$} (and also assumes that $\OO$ is a domain).
\end{rem}
\begin{proof}
By the preceding lemma, an ideal of a product $\OO_1 \cross \cdots \cross \OO_s$ is just a product $\aa_1 \cross \cdots \cross \aa_s$, which is principal (resp{.} invertible) if and only if every $\aa_i$ is: hence we can assume that $s = 1$.

The reverse direction is trivial (principal fractional ideals are invertible), so let $\aa$ be an ideal of $\OO$ with inverse $\aa^{-1}$, $\aa\aa^{-1} = \OO$. Since all ideals of $\OO_L$ are principal, we can assume that $\aa\OO_L = \OO_L$ and hence $\aa^{-1}\OO_L = \OO_L$ as well. We can express
\[
  1 = \sum_i \alpha_i \beta_i, \quad \alpha_i \in \aa, \quad \beta_i \in \aa^{-1}.
\]
Let $\pp_1,\ldots,\pp_r$ be the valuation ideals coming from the field factors of $L$. Since $1 \notin \pp_1$, some term $\alpha_i\beta_i$, say $\alpha_1\beta_1$, is nonzero mod $\pp_1$. But because $\pp_i \intsec \OO = \pp_1 \intsec \OO$ is the unique maximal ideal of $\OO$ for all $i$, we have $\alpha_1\beta_1 \notin \pp_i$ for all $i$. This implies that $\alpha_1$ and $\beta_1$ are units, whose product lies in $\OO^\cross$. Now since
\[
  \aa \supset \alpha_1\OO, \quad \bb \supset \beta_1\OO
\]
and $\aa\bb = \OO$, equality must hold.
\end{proof}
\begin{lem}\label{lem:sqrs are inv}
If $\cc \subseteq R$ is a lattice in a cubic algebra over a local field $K$, then $\cc^2$ is invertible in its endomorphism ring $\End(\cc^2)$.
\end{lem}
\begin{proof}
First, $\cc\OO_{R}$ is an invertible $\OO_{R}$-ideal, which, since $\OO_{R}$ is a product of PID's, we can scale to be $\OO_{R}$.

We first claim that $\cc$ contains a unit, or else has a special form for which $\cc^2 = \OO_{R}$ is clearly invertible. Let $p$ be a uniformizer for $\OO_K$ and $k = \OO_K/p\OO_K$ the residue field. The units of $\OO_{R}$ are those elements whose projections to the cubic $k$-algebra $\ba{C_1} = C_1/pC_1$ are non-units (that is, zero divisors). The non-units of $\ba{C_1}$ are the union of at most three proper subspaces (the projections of the valuation ideals of each field factor). The projection $\bar\cc$ of $\cc$ down to $\ba{C_1}$ cannot lie in any of these subspaces since $\cc\OO_{R} = C_1$. An easy theorem in linear algebra is that a vector space over a field $k$ cannot be the union of fewer than $\size{k} + 1$ proper subspaces. We conclude that $\cc$ contains a unit except if $\size{k} = 2$ and ${R} = K \cross K \cross K$ has three field factors. In this case, the only $\bar\cc \subseteq \ba{C_1}$ instantiating this case is
\[
  \bar\cc = \{(a;b;c) \in \FF_2^3 : a + b + c = 0\}.
\]
It is evident that $\bar\cc^2$ is the whole of $\FF_2^3$, whence by Nakayama's lemma, $\cc^2$ is the whole of $\OO_{R}$.

Now we can assume that $\cc$ contains a unit, which we scale to equal $1$. We claim that $\cc^3 = \cc^2$. By the theory of modules over a PID, we can find a basis $\{1, \alpha, \beta\}$ for $\OO_K$ such that $\{1, p^i\alpha, p^j\beta\}$ is a basis for $\cc$ for some integers $i, j \geq 0$. Now by translation, we can assume that $\alpha\beta = t \in \OO_K$. We then have
\[
  \cc^3 = \cc^2 + \<\alpha^3, \alpha^2\beta, \alpha\beta^2, \beta^3\> = \cc^2 + \<\alpha^3, t \alpha, t \beta, \beta^3\>
\]
The elements $t\alpha$ and $t\beta$ are certainly already in $\cc \subseteq \cc^2$. As for $\alpha^3$, since $\alpha$ is an integral element of $R$, its characteristic polynomial expresses $\alpha^3$ as an $\OO_K$-linear combination of $\alpha^2$, $\alpha$, and $1$, all of which lie in $\cc^2$. So $\cc^3 = \cc^2$.

We conclude that $\cc^4 = \cc^3 = \cc^2$, so $\cc^2$ is closed under multiplication and hence is an order. In particular, it coincides with its endomorphism ring and in particular is invertible.
\end{proof}
\begin{rem}
Although Lemma \ref{lem:sqrs are inv} is simple to state, we have not found it anywhere in the literature. In general, we suspect that if $\cc$ is a lattice in an algebra $L$ of rank $n$, then $\cc^{n-1}$ and all higher powers of $\cc$ are invertible in their common endomorphism ring. This is not hard to prove if $\ch k_K > n$. That the exponent $n-1$ is sharp is seen from the cute example
\[
  \cc = \ZZ_p + (0;1;\ldots;n-1)\ZZ_p + p\ZZ_p^n \subseteq \ZZ_p^n.
\]
The power $\cc^i$ consists of all sequences $(a_0;\ldots;a_{n-1}) \in \ZZ_p^n$ that are congruent modulo $p$ to the values
\[
  (f(0); f(1); \ldots; f(n-1))
\]
of a polynomial $f$ of degree at most $i$ with coefficients in $k$. If $p > n$, then this power stabilizes to the whole of $\ZZ_p^n$ only for $i \geq n - 1$. 
\end{rem}


\subsection{Self-duality of the count of quartic orders}
\begin{thm}[\textbf{Local quartic O-N in the tame, not totally split case}] \label{thm:O-N_quartic_local_tame}
Assume $K$ is a local field of residue characteristic not $2$. Let $C \subseteq R$ be a cubic \'etale order that is \emph{not} totally split. Then the assignment $f_C$ to each $L$ of the number of orders $\OO \subseteq L$ with resolvent $C$ is self-dual.
\end{thm}
\begin{proof}
As in the cubic case, the proof proceeds by reduction to the zero case (i.e{.} that $f(0) = \hat{f}(0)$).

Note that $M \cong M'$ as Galois modules (one can even make this canonical, using the unique alternating bilinear form on $M$). The fixity of this Galois module is easy to compute:
\[
  \size{H^0} = \size{H^0(K, R)} = \begin{cases}
    1 & \text{if $R$ is a field} \\
    2 & \text{if $R \cong K \cross K_2$ for some quadratic field $K_2$} \\
    4 & \text{if $R \cong K \cross K \cross K$.}
  \end{cases}
\]
This can be written concisely as
\[
  \size{H^0} = \frac{\size{R^\cross[2]}}{2}
\]
(the latter formula will work especially well for our case).

Moreover, since the unramified cohomology is self-orthogonal, we have
\[
  \size{H^1} = \size{H^0}^2 = \frac{\size{R^\cross[2]}^2}{4}.
\]
Since we are excluding the case $R \isom K \cross K \cross K$, there are just two possibilities:
\begin{itemize}
  \item If $R$ is a field, then $\size{H^1} = 1$, and there is nothing to prove, as any function on $H^1$ is self-dual.
  \item If $R$ is the product of two fields, then $\size{H^1} = 4$. Pick an $\FF_2$-basis $\<\sigma_1,\sigma_2\>$. The Tate pairing is given by the unique alternating pairing on $H^1$. The space of functions on $H^1$ is four-dimensional, and a basis is
  \[
    \1_{\<\sigma_1\>}, \1_{\<\sigma_2\>}, \1_{\<\sigma_1 + \sigma_2\>}, \1_{\{0\}}.
  \]
  Note that the first three basis elements are self-dual, while the fourth differs from its dual even at $0$. This proves that if $f$ is a function on $H^1$ with $f(0) = \hat f(0)$, then $f = \hat f$.
\end{itemize}
So we have reduced local O-N to the following lemma:
\end{proof}
\begin{lem}\label{lem:tame0}
Assume $K$ is a local field of residue characteristic not $2$. Let $C \subseteq R$ be an order in an \'etale algebra. Then the assignment $g_C$ to each $L$ of the number of orders $\OO \subseteq L$ with resolvent $C$ satisfies self-duality at $0$:
\[
  \hat g_C(0) = g_C(0).
\]
\end{lem}
\begin{proof}
As in the cubic case, the proof is by explicit multijection.

On the one hand, 
\[
  \size{H^0} \cdot \hat g_C(0) = \sum_\sigma g_C(\sigma)
\]
counts all quartic orders with cubic resolvent $C$, and using Theorem \ref{thm:hcl_quartic_sbi}, these can be parametrized by self-balanced ideals $(C, \cc, \delta)$, where $\delta$ ranges over a set of representatives for $R^{N=1}/\big(R^{N=1}\big)^2$. On the other hand, $g_C(0)$ is the number of orders with resolvent $C$ in $K \cross R$. Write such an order as $\OO = \OO_K + 0 \cross \aa$, where $\aa \subseteq R$ is a lattice. The condition that $\OO$ be a ring is (by Theorem \ref{thm:hcl_quartic}) subsumed by the resolvent conditions, namely that
\begin{enumerate}[(i)]
  \item $N_C(\aa) = 1$,
  \item $\Phi_{4,3}(0; \alpha) = \alpha'\alpha'' \in C$ for all $\alpha \in \aa$. 
\end{enumerate}
Our aim is to associate $\size{H^0}$ values of $(\cc,\delta)$ to each value of $\cc$.

The multijection is as follows. First, $\cc$ may not be invertible. Let $C_1 = \End \cc^2$ and $\cc_1 = \cc C_1$, an invertible and thus a principal $C_1$-ideal. Let
\begin{equation} \label{eq:multij quartic}
  \aa_1 = \frac{[C_1 : C]}{\delta [\cc_1 : \cc]} \( \cc_1^2\) ^{-1}.
\end{equation}
Finally, since $\cc_1$ and $\aa_1$ are both principal and thus scalar multiples of each other, we can take $\aa$ to be an ideal that sits inside $\aa_1$ as $\cc$ sits inside $\cc_1$: that is, if $\cc_1 = \gamma_0 C_1$ and $\aa_1 = \alpha_0 C_1$, then
\[
  \aa = \frac{\alpha_0}{\gamma_0} \cc.
\]
Before checking that this $\aa$ yields a valid ring, we check how many-to-one our multijection is. First note that $\aa$ determines $C_1 = \End \aa^2$ and $\aa_1 = \aa C_1$, and in particular the index $[\cc_1 : \cc] = [\aa_1 : \aa]$. Then, by \eqref{eq:multij quartic}, the ``shadow'' $\bb = \delta \cc^2 = \delta \cc_1^2$ is determined. Note that $\bb$ is an invertible $C_1$-ideal of norm $[C_1:C]^2/[\aa_1 : \aa]^2$, a square. 
The pairs $(\cc_1, \delta)$ satisfying $\bb = \delta \cc_1$, where $\cc_1$ is an invertible $C_1$-ideal and $\delta$ is one of the representatives for $R^{N=1}/\big(R^{N=1}\big)^2$, are found to be $\size{H^0}$ in number by an argument identical to the cubic case. Finally, locating $\cc$ within $\cc_1$ involves the same choice as locating $\aa$ within $\aa_1$. So we have a string of many-to-one correspondences
\begin{equation} \label{eq:multij quartic dgm}
  (\cc, \delta) \xrightarrow{\text{$n$ to 1}}
  (\cc_1, \delta) \xrightarrow{\text{$\size{H^0}$ to 1}}
  \aa_1 \xleftarrow{\text{$n$ to 1}}
  \aa,
\end{equation}
and thus overall there are $\size{H^0}$ times as many $(\cc, \delta)$ as $\aa$.

It remains to prove that the correspondence \eqref{eq:multij quartic dgm} preserves the resolvent and balancing conditions. As for the first condition, regarding the discriminant of the ring, we leave it to the reader to verify that
\begin{equation} \label{eq:xnorm}
  N_C(\aa) = \frac{1}{N(\delta) N_C(\cc)^2}.
\end{equation}
Now we may assume that both sides of \eqref{eq:xnorm} are $1$. Since $\cc_1$ is principal, we may assume that $\cc_1 = C_1$, adjusting $\delta$ by a square if necessary. Then by the conditions of Theorem \ref{thm:hcl_quartic_sbi}, $N(\delta)$ is a square $t^2$ and
\[
  \frac{1}{t} = N_C(\cc) = \frac{[C_1 : C]}{[\cc_1 : \cc]};
\]
thus
\[
  \aa_1 = \frac{t}{\delta} C_1.
\]

We first prove that if $\cc$ satisfies its resolvent condition, so does $\aa$. Any $\alpha \in \aa \subseteq \aa_1$ has the form $\alpha = \frac{t}{\delta} \beta$, $\beta \in C_1$, and then
\[
  \alpha' \alpha'' = \frac{t^2}{\delta' \delta''} \beta' \beta'' = \delta \beta' \beta'';
\]
and we note that if $\beta \in C_1$, then $\beta' \beta'' \in C_1$ as well, by the relation
\begin{equation} \label{eq:resolvent trick}
  \beta' \beta'' = \underbrace{\beta \beta' + \beta \beta'' + \beta' \beta''}_{{} \in \OO_K}
  - \beta(\underbrace{\beta + \beta' + \beta''}_{{} \in \OO_K} - \beta).
\end{equation}
Thus if $\cc$ satisfies the resolvent condition ($\delta \cc^2 \subseteq C$), then $\aa$ satisfies the resolvent condition ($\alpha' \alpha'' \in C$ for all $\alpha \in \aa$). To prove the converse, it suffices to show that
\[
  B = \{\beta' \beta'' : \beta \in \cc\}
\]
spans $C_1$ over $\OO_K$, where $\cc \subseteq C_1$ is a sublattice with $\cc^2 = C_1$.

We first claim that $\cc$ contains a unit. Let $p$ be a uniformizer for $\OO_K$ and $k = \OO_K/p\OO_K$ the residue field. In the cubic $k$-algebra $\ba{C_1} = C_1/pC_1$, the non-units are the union of at most three subspaces. The projection $\bar\cc$ of $\cc$ down to $\ba{C_1}$ cannot lie in any of these subspaces since $\cc^2 = C_1$, so, since $\size{k} \geq 3$, $\bar{\cc}$ must contain a unit, which lifts to a unit in $\cc$. There is no harm in rescaling $\cc$ so that $1 \in \cc$.

Now $1 \in B$. Also, for each $\beta \in \cc$,
\[
  \<B\> \ni \beta'\beta'' - (\beta' + 1)(\beta'' + 1) + \tr \beta - 1 = \beta,
\]
and thus
\[
  \<B\> \ni \beta' \beta'' - (\beta \beta' + \beta \beta'' + \beta' \beta'') + \beta(\beta + \beta' + \beta'') = \beta^2.
\]
But since $\ch k \neq 2$, the elements $\beta^2$, for $\beta \in \cc$, generate $\cc^2 = C_1$, completing the proof.
\end{proof}

For the totally split case, the method of proof of Theorem \ref{thm:O-N_quartic_local_tame} fails, because $\dim H^1 = 4$ and $\hat f(0) = f(0)$ is no longer enough to imply $\hat f = f$. But when we count by discriminant instead of resolvent, it can be rescued, due to the following symmetry argument.

\begin{thm}
Fix a cubic algebra $R$ over a local field $K$, and let $g_D : H^1 \to \ZZ$ count the number of orders in a quartic algebra $L \in H^1$ with discriminant $D$. That is, $g_D$ is the sum of all the $g_C$'s in Theorem \ref{thm:O-N_quartic_local_tame} over $C \subseteq R$ of discriminant $D$. Then $g_D$ is self-dual.
\end{thm}
\begin{proof}
All cases are covered by Theorem \ref{thm:O-N_quartic_local_tame} except for the totally split case $R \cong K \cross K \cross K$, where $\dim H^1 = 4$. We can write $H^1 = \<\sigma_1, \tau_1, \sigma_2, \tau_2\>$, where the $\sigma$'s and $\tau$'s correspond to Kummer elements
  \[
    \sigma_1 : (u,u,1), \quad \tau_1 : (p,p,1), \quad \sigma_2 : (1,u,u), \quad \tau_2 : (1, p, p),
  \]
  where $p \in \OO_K$ is a uniformizer and $u$ is a non-square unit. The Tate pairing is given, by Theorem \ref{thm:Tate_pairing}, by the $\mu_2$-valued pairing
  \[
    \<\sigma_1,\tau_2\> = \<\sigma_1, \tau_2\> = -1, \quad \<\sigma_1,\sigma_2\> = \<\tau_1,\tau_2\> = \<\sigma_1,\tau_1\> = \<\sigma_2,\tau_2\> = 1.
  \]
  The group $\Aut R \cong \S_3$ acts on $H^1$, permuting $\sigma_1, \sigma_2, \sigma_1 + \sigma_2$ and $\tau_1, \tau_2, \tau_1 + \tau_2$ in the permutation manner.  There are five orbits, represented by $0$, $\sigma_1$, $\tau_1$, $\sigma_1 + \tau_1$, and $\sigma_1 + \tau_2$. Note that $g_D$ must be constant on each orbit, because its definition is $S_3$-invariant. The functions
  \[
    \1_{\<\sigma_1,\sigma_2\>}, \1_{\<\tau_1,\tau_2\>}, \1_{\<\sigma_1 + \tau_1, \sigma_2 + \tau_2\>}, \sum_{\pi \in \Aut R} \1_{\pi\( \<\sigma_2, \sigma_1 + \tau_2\>\) }, \1_{\{0\}}
  \]
  form a basis for the $5$-dimensional space of $\S_3$-invariant functions on $H^1$. The first four are self-dual, while the last differs from its dual even at $0$; so, since $g_D(0) = \hat g_D(0)$ by Lemma \ref{lem:tame0}), the coefficient of the last basis element must be $0$ and $g_D$ is self-dual.
\end{proof}


\part{Counting quartic rings with prescribed resolvent}
\label{part:quartic_count}

\section{Introduction}

Here end the cases in which a conceptual, bijective argument has been found to suffice for proving local reflection for quartic rings. To win the remaining cases, we attack a problem that has interest in its own right: counting orders $\OO$ in a quartic algebra $L$ over a local field $K$ whose cubic resolvent ring $C \subseteq R$ is fixed.

The index $[\OO_L : \OO]$ being fixed by the condition $\Disc \OO = \Disc C$, we must analyze the resolvent condition $\Phi_{4,3}(\OO/\OO_K) \subseteq C/\OO_K$, where $\Phi_{4,3} : L/K \to R/K$ is the resolvent map. Recall that, with respect to bases of $\OO$ and $C$, $\Phi_{4,3}$ is given by a pair
\[
  (\M, \N) = \( [M_{ij}], [N_{ij}] \)
\]
of symmetric $3 \times 3$ matrices, and the resolvent condition can be viewed as the $\OO_K$-integrality of the entries (properly scaled to account for the tracedness condition). By suitably choosing coordinates, we can ensure that only the integrality of the entries
\[
  M_{11}, \quad N_{11}, \quad M_{12}, \textand M_{22}
\]
is in doubt. Of these, the condition on $M_{11}$ is the most challenging. It amounts to a quadratic condition on the first basis vector $\xi_1$ of $C$, that is, a conic on some pixel (determined by the $N_{11}$-condition) in $\PP^2(\OO_K)$. The solubility of this conic over $K$ is governed by the \emph{Hilbert symbol,} which we analyze. It is very hard in general to tell if any $K$-points of the conic lie in the requisite pixel, but if there is even one such $K$-point, then, using the rational parametrization of a conic with a basepoint, the volume of points in the pixel is easy to determine. Accordingly, our approach to solving the $M_{11}$-condition is a three-step one:
\begin{itemize}
  \item Determine the sum of the solution volumes for $\xi_1$ over all quartic algebras $L$.
  \item Find restrictions on what $L$ can yield a nonzero volume and what that volume can be, providing an upper bound (the \emph{bounding step}).
  \item If the sums of these upper bounds agree, deducing that the bound is attained everywhere (the \emph{summing step}).
\end{itemize}
The $M_{12}$ and $M_{22}$ conditions are essentially linear. We use the computer program LattE to sum the ring totals over all possible values of the discrete data and verify the local reflection theorem.

The case of wildly ramified resolvent (splitting type $1^21$) is still in progress. Except for brief remarks, it has been omitted from this edition. Also omitted are  the adaptations to be made when $\ch k_K > 2$, where, in view of Theorem \ref{thm:O-N_quartic_local_tame}, only splitting type $(111)$ need be considered. It involves only the black, brown, beige, and white zones; the conics are all very easy to solve and yield the same answers as the wild case upon substituting $e = 0$. So in the sequel, \textbf{$R$ is a tamely ramified \'etale algebra over a $2$-adic local field $K$.}

Corresponding to $R$, there is a Galois module $M_R$ whose underlying group is $\C_2 \cross \C_2$. We will work extensively with $H^1(K, M_R$), which we abbreviate to $H^1$.

\section{The group \texorpdfstring{$H^1$}{H1} of quartic algebras with given resolvent}
We fix a local field $K$ and an separable closure $\bar{K}$. Let $\OO_{\bar{K}}$ be the ring of integers in $\bar{K}$. (If the reader is uncomfortable with non-Noetherian rings, he can take $\bar{K}$ to be instead the compositum of all extensions of $K$ of degree at most $4$; the Galois cohomology and all proofs will be unaffected.)

By Theorem \ref{thm:Kummer_new}\ref{it:Kum_quartic}, we can identify $H^1$ naturally with
\[
  R^{N=1}/\(R^{N=1}\)^2.
\]
Now there is a natural isomorphism
\[
  R^\cross/(R^\cross)^2 \isom K^\cross/(K^\cross)^2 \cross  R^{N=1}/\(R^{N=1}\)^2.
\]
Thus for any $\alpha \in R^\cross$, we can talk about the class $[\alpha]$ of $\alpha$ in $H^1$, that is, the class of $\alpha^3/N(\alpha)$.

Hence the structure of $H^1$ can be uncovered by taking a suitable Shafarevich basis of $R^\cross/(R^\cross)^2$ and removing a basis of $K^\cross /(K^\cross)^2$, which, by Lemma \ref{lem:dist_base_change}, maps in isometrically:
\begin{lem} \label{lem:levels_quartic}
If $R/K$ is a cubic \'etale extension, then $H^1$ is an $\FF_2$-vector space of dimension $2 \dim_{\FF_2} H^0 + 2ef$. It has a basis of $\dim_{\FF_2} H^0$ nonunits, $\dim_{\FF_2} H^0$ intimate units, and $2ef$ generic units; the generic units can be chosen as follows:
\begin{enumerate}[$($a$)$]
  \item If $R$ is unramified, we take $2ef$ units of the form $1 + x\pi^{2i + 1}$, where $0 \leq i < e$ and $x$ ranges over $f$ elements whose reductions mod $\pi$ form an $\FF_2$-basis of $k_R^{\tr = 0}$ for each $i$.
  \item If $R$ is totally tamely ramified, we take $2ef$ units of the form $1 + x\pi^{2i+1}$, where $0 \leq i < 3e$ but $3 \nmid i$, and $x$ ranges over an $\FF_2$-basis of $k_K$ for each $i$.
\begin{wild}
  \item If $R \isom K \cross Q$ is partially wildly ramified, we take $2ef$ units of the form $(1 + y\pi^{i+1}; 1 + x\pi_Q^{2i+1})$ where $0 \leq i < 2e$, where $x$ ranges over an $\FF_2$-basis of $k_K$ for each $i$, and where $y \in \OO_K$ is adjusted for each $(i,x)$.
\end{wild}
\end{enumerate}
\end{lem}
In the unramified case, we define the \emph{level space} 
\[
  \L_{i} = \begin{cases}
    H^1 & i = -1 \\
    \{[\alpha] \in H^1 : \alpha \equiv 1 \mod \pi^{2i}\} & 0 \leq i \leq e \\
    \{[1]\} & i = e+1,
  \end{cases}
\]
noting that $\size{\L_i} = \size{H^0}q^{2(e-i)}$ for $0 \leq i \leq e$ and that $\L_i^\perp = \L_{e-i}$ for all $i$ by Lemma \ref{lem:Hilb_prod_size}.

Likewise, in the ramified case, we define the \emph{level space}
\[
  \L_{i} = \begin{cases}
  H^1 & i = -1 \\
  \{[\alpha] \in H^1 : \alpha \equiv 1 \mod \pi^{i}\} & 0 \leq i \leq 2e \\
  \{[1]\} & i = e+1,
  \end{cases}
\]
noting that $\size{\L_i} = \size{H^0}q^{2e-i}$ for $0 \leq i \leq e$ and that $\L_i^\perp = \L_{2e-i}$ for all $i$ by Lemma \ref{lem:Hilb_prod_size}.

We will occasionally let
\[
  e' = \begin{cases}
    e, & R \text{ unramified} \\
    2e, & R \text{ ramified}
  \end{cases}
\]
to shorten lemma statements.

We define the \emph{level} $\ell(\alpha)$ of an element $[\alpha] \in H^1$ as the largest $i \geq 0$ for which $[\alpha] \in \L_i$. We have $\ell(1) = e'+1$. By convention, if $[\alpha] \notin \L_0$, we set
\[
  \ell(\alpha) = -1/2
\]
to shorten some future statements.

\section{Reduced bases}

Define a valuation on $\bar K^n$ by
\[
  v(x_1;\ldots;x_n) = \min\{v(x_1),\ldots, v(x_n)\}.
\]

Let $R$ be a rank-$n$ \'etale algebra over a local field $K$. We can Minkowski-embed $R$ into $\bar{K}^n$. 
\begin{defn}\label{defn:red_basis}
Let $I$ be an $\OO_K$-lattice in $R$, and let $\omega \in \bar{K}^n$ be a multiplier with the following property:
\begin{enumerate}[$(*)$]
  \item\label{iota:*} If $\iota$, $\iota'$ are two coordinates of the same field factor of $R$, then $\omega^{(\iota)}$ and $\omega^{(\iota')}$ have the same valuation.
\end{enumerate} 
A basis $(\rho_1,\ldots,\rho_n)$ for $\omega I$ is called \emph{reduced} if
\begin{enumerate}[$($a$)$]
  \item $v(\rho_1) \leq \cdots \leq v(\rho_n)$;
  \item If $\rho \in \omega R$ is decomposed as
  \[
  \rho = \sum_i c_i \rho_i, \quad c_i \in K,
  \]
  then for each $i$,
  \[
  v(c_i \rho_i) \geq v(\rho).
  \]
\end{enumerate}
\end{defn}

This notion has the following properties:
\begin{prop}\label{prop:rb} Let $\omega I$ be as above.
\begin{enumerate}[$($a$)$]
\item\label{rb:exists} There exists a reduced basis $(\rho_1, \ldots, \rho_n)$ for $\omega I$.
\item\label{rb:max} If $(\rho_1', \ldots, \rho_n')$ is any other basis for $\omega I$, sorted so that $v(\rho_1') \leq \cdots \leq v(\rho_n')$, then for each $k$,
\[
  v(\rho_k') \leq v(\rho_k).
\]
In particular, if both bases are reduced, equality holds.
\item\label{rb:unique} If $(\rho_1', \ldots, \rho_n')$ is another reduced basis for $\omega I$, then
\[
  \rho_i' = \sum_j c_{ij} \rho_j
\]
for some change-of-basis matrix
\begin{equation} \label{eq:change_basis}
  \left[c_{ij}\right] \in \GL_n \OO_K, \quad v\(c_{ij}\) \geq v(\rho_i) - v(\rho_j).
\end{equation}
Conversely, any matrix $[c_{ij}]$ satisfying \eqref{eq:change_basis} yields a new reduced basis $(\rho_i')_i$.
\item\label{rb:span} As $\OO_{\bar K}$-modules,
\[
  \<\pi^{-v(\rho_i)} \rho_i : 1 \leq i \leq n\>
  = \<\pi^{-v(\rho)} \rho : \rho \in \omega R^\cross\>.
\]
\end{enumerate}
\end{prop}

\begin{proof}
\begin{enumerate}[$($a$)$]
  \item Choose a basis $(\rho_1, \ldots, \rho_n)$ such that the sum of the valuations $v(\rho_1) + \cdots + v(\rho_n)$ is maximal. This can be done because there are only finitely many possible valuations of primitive vectors in $\omega I$. Sort the $\rho_i$ in increasing order of valuation. We claim $(\rho_i)_i$ is reduced.
    
  Let $\rho = \sum_i c_i \rho_i$ be given. Let $a$ be the minimal valuation $v(c_i \rho_i)$ of a term, and suppose that $a < v(\rho)$. Then we have a linear dependency
  \[
    \sum_{v(c_i\rho_i) = a} c_i \rho_i \equiv 0 \mod \pi^a \mm_{\bar K}.
  \]
  Choose $j$ such that $v(c_j\rho_j) = a$ and $v(c_j)$ is minimal. Then
  \[
    \rho_j' = \sum_{v(c_i\rho_i) = a} \frac{c_i}{c_j} \rho_j
  \]
  is an element of $\omega I$ whose valuation exceeds $v(\rho_j)$. Since the coefficient of $\rho_j$ in $\rho_j'$ is $1$, replacing $\rho_j$ by $\rho_j'$ does not change the span $\omega I$ but increases the valuation sum $\sum_i c_i \rho_i$, contradicting the choice of basis $(\rho_i)_i$.
  \item Since $(\rho_i)_i$ and $(\rho_i')_i$ are bases for the same module $\omega I$, such a $[c_{ij}] \in \GL_n(\OO_K)$ certainly exists. Applying the reducedness property to each decomposition
  \[
    \rho_i' = \sum_{j} c_{ij} \rho_j
  \]
  yields a bound
  \[
    v(c_{ij}) \geq v(\rho_i') - v(\rho_j).
  \]
  Suppose that $v(\rho_k) < v(\rho_k')$ for some $k$. Then for $j \leq k \leq i$,
  \[
  v(\rho_j) \leq v(\rho_k) < v(\rho_k') \leq v(\rho_i'),
  \]
  so $c_{ij}$ has positive valuation. Thus, when the matrix $[c_{ij}]$ is reduced modulo $\pi$, it has an $(n - k + 1) \cross k$ block of $0$'s, large enough to make the determinant vanish, which is a contradiction.
  
  \item By the preceding part, $v(\rho'_k) = v(\rho_k)$. So the associated matrix $[c_{ij}]$ must satisfy
  \begin{equation}
  v(c_{ij}) \geq v(\rho_i') - v(\rho_j) = v(\rho_i) - v(\rho_j).
  \end{equation}
  Conversely, if $[c_{ij}]$ is an invertible matrix satisfying this inequality, we get a new basis $\rho_i'$ with $v(\rho_i') \geq v(\rho_i)$. Equality must hold, and now since $\sum_k v(\rho_k')$ achieves the maximal value, $(\rho_k')_k$ is reduced by the proof of part \ref{rb:exists}.
  \item The $\subseteq$ direction is obvious. For the $\supseteq$ direction, let $\rho \in \omega R^\cross$ be given. Since $(\rho_i)_i$ is reduced,
  \[
    \rho = \sum_i c_i \rho_i, \quad v(c_i \rho_i) \geq v(\rho),
  \]
  so
  \[
    \pi^{-v(\rho)} \rho = \sum_i \(\pi^{v(\rho_i) - v(\rho)}c_i\) \pi^{-v(\rho_i)} \rho_i,
  \]
  and the parenthesized coefficients belong to $\OO_{\bar K}$, as desired. \qedhere
\end{enumerate}
\end{proof}

We can find reduced bases with added structure.
\begin{defn}
Fix an ordering $R = R_1 \cross \cdots \cross R_{r}$ of the field factors of $R$. For $\rho \in \omega R$, let $k$ be the minimal index such that $v(\rho^{(k)}) = v(\rho)$. We say that $\rho$ is \emph{$R^{(k)}$-led,} and we define the \emph{leader} of $\rho$ to be the normalization
\[
  \ldr(\rho) = \(\frac{\rho}{\pi^u \omega}\)^{(i)},
\]
where $u$ is the unique integer for which $\ldr(\rho)$ is a primitive vector in $\OO_R$. We say that a reduced basis $(\rho_i)_i$ is \emph{well-led} if the leaders $\ldr(\rho_i)$ consist of a reduced basis for $\OO_R^{(k)}$ for each $k$.
\end{defn}

\begin{prop} \label{prop:well_led}
Every $\omega I$, as above, admits a well-led basis.
\end{prop}
\begin{proof}
Consider the element
\[
  \eta_\epsilon = (\pi^\epsilon ; \pi^{2 \epsilon} ; \ldots ; \pi^{r \epsilon})
\]
where $\epsilon$ is a positive rational number, small enough that if $a_1 < a_2$ are two valuations of elements in $\omega I$, then $r\epsilon < a_2 - a_1$. By Proposition \ref{prop:rb}\ref{rb:exists}, there is a reduced basis $\rho_1\eta_\epsilon, \ldots, \rho_n\eta_\epsilon$ for $\omega\eta_\epsilon I$. Each basis element $\rho_i \eta_\epsilon$ has some valuation $u + k\epsilon$, $u \in \QQ$, $k \in \{1,\ldots,r\}$, indicating that $\rho_i$ is $R^{(k)}$-led. Since replacing $\epsilon$ by $0$ preserves non-strict inequalities among valuations in $\eta_\epsilon \omega I$, the $\rho_i$ form a reduced basis for $\omega I$, which we claim is well-led.

Given $\alpha \in \OO_{R^{(k)}}$ primitive, decompose $\omega\alpha = \sum_j c_j \rho_j$ as an element of $\omega R$. We have
\[
  v(c_j\rho_j \eta_\epsilon) \geq v(\alpha \eta_\epsilon) = v(\alpha) + i\epsilon,
\]
and for a nonempty subset of $j$, equality must hold and, in particular, $\rho_j$ must be $R^{(k)}$-led. Let $L_k$ be the set of indices $j$ for which $\rho_j$ is $R^{(k)}$-led, and let $B_k = \{\ldr(\rho_j) : j \in L_k\}$. Now we have
\[
  \omega\alpha = \sum_{j \in L_k} c_j \rho_j^{(k)} + \omega\alpha'
\]
with each term of valuation at least $v(\alpha)$, and where $\alpha' \in \OO_R$ with $v(\alpha') > v(\alpha)$. We can rewrite this as
\[
  \alpha = \sum_{j \in L_k} c_j' \ldr(\rho_j) + \alpha'
\]
where $c_j' = c_j\pi^u \in \OO_K$. We can iteratively decompose $\alpha'$ the same way, and as its valuation goes to infinity, we get a decomposition
\begin{equation} \label{eq:c_j''}
  \alpha = \sum_{j \in L_k} c_j'' \ldr(\rho_j).
\end{equation}
So $B_k$ generates $\OO_{R^{(k)}}$, and in particular, $\size{L_k} \geq [R : R^{(k)}]$. However,
\[
  \sum_{k} [R : R^{(k)}] = n = \sum_k \size{L_k}.
\]
So equality holds and each $B_k$ is a basis for $\OO_R$. Since the decomposition \eqref{eq:c_j''} has every term of valuation at least $v(\alpha)$, and $\alpha$ was any primitive vector, $B_k$ is in fact reduced.
\end{proof}
It is evident that reduced indices take a limited number of values modulo $1$. Indeed, we have the following:
\begin{cor}\label{cor:idxs_mod_1}
  If $I \subseteq R$ is a lattice, then the multiset of valuations $\{a_i\}_i = \{v(\rho_i)\}_i$ mod $1$ of reduced basis elements for $\omega I$ depends only on $R$ and $\omega$, not on $I$. It consists of $f_{R_i/K}$ copies of
  \[
    v(\omega^{(R_i)}) + \frac{j}{e_{R_i/K}}, \quad j = 0, 1, \ldots, e_{R_i/K} - 1,
  \]
  where $R_i$ ranges over the field factors of $R = R_1 \cross \cdots \cross R_r$.
\end{cor}
\begin{proof}
Taking a well-led basis and passing to the leaders, we reduce to the case that $R = R_i$ is a field. Then $\omega$ has equal valuations in all coordinates, and we may assume that $\omega = 1$. Let $\rho_i = \pi_R^{e_{R/K}a_i} \xi_i$, where the $\xi_i \in \OO_R^\cross$ differ only by a unit from the normalizations used before. For each congruence class of $a_i$ modulo $1$, note that if the set of corresponding $\xi_i$ is linearly dependent modulo $\pi_R$, then one of the $\rho_i$ could be increased by an $\OO_K$-linear combination of the others to increase its valuation, contradicting the hypothesis that our basis is reduced. So the $\xi_i$ corresponding to each congruence class of $a_i$ modulo $1$ are linearly independent, and in fact must form a basis for $k_R$ over $k_K$ in order for there to be the full number $e_{R/K}f_{R/K}$ of $\xi_i$. This establishes the claimed multiset.
\end{proof}

A reduced basis for $\omega I$ does not always remain reduced when we extend the ground field $K$. To study this, we make the following definition.
\begin{defn}
An \emph{extender basis} for $\omega I$ is a basis $(\rho_1, \ldots, \rho_n)$ for $\omega I \tensor_{\OO_K} \OO_{\bar{K}}$ such that the vectors
\[
  \xi_i = \pi^{-v(\rho_i)} \rho_i
\]
form an $\OO_{\bar{K}}$-basis for $\OO_{\bar{K}}^n$. The valuations $a_i = v(\rho_i)$ are called the \emph{extender indices} of the basis, and the $\xi_i$ are called the \emph{extender vectors}.
\end{defn}

If it consists of $\rho_i \in \omega I$, an extender basis is easily seen to be reduced. Fortunately, in the cases of tamely ramified resolvent, this always holds:

\begin{prop}\label{prop:ext_basis_tame}
If $R/K$ is tamely ramified, then any reduced basis of $\omega I$ is an extender basis.
\end{prop}
\begin{proof}
In view of Proposition \ref{prop:rb}\ref{rb:span}, it is enough to show that
\[
  \OO_{\bar K} \<\pi^{-v(\rho)} \rho : \rho \in \omega R^\cross\> = \OO_{\bar K}^n.
\]
We immediately reduce to the case that $R$ is a field, and then we may assume $\omega = 1$. It suffices to prove that, for some reduced basis $\rho_1, \ldots, \rho_n$ of $R$,
\[
  \det \left[\pi^{-v(\rho_i)}\rho_i^{(j)}\right]_{i,j} \sim 1.
\]
But that follows from the familiar formula for the discriminant of a tamely ramified extension.
\end{proof}

\begin{wild}
\subsection{Splitting type \texorpdfstring{$1^21$}{1²1}}
When $R$ is wildly ramified, it is no longer possible for a reduced basis to be an extender basis (again by discriminant considerations). Here the case that concerns us is that $R = K \cross Q$ is a partially ramified cubic extension. Let $(\rho_1, \rho_2, \rho_3)$ with valuations $(a_1, a_2, a_3)$ be a reduced basis. Since $R \tensor Q \isom Q \cross Q \cross Q$ is an unramified extension of $Q$, there exists a well-led extender basis $(\bar{\rho}_1, \bar{\rho}_2, \bar{\rho}_3)$ of elements of $\omega I \tensor_K Q$. The extender indices $\bar a_i$ satisfy $\bar a_i \geq a_i$ by Proposition \ref{prop:rb}\ref{rb:max}. Our task in this section is to pin them down precisely.
\begin{lem} \label{lem:indices_wedge}
If $(\bar{\rho}_1, \ldots, \bar{\rho}_n)$ is an extender basis of $\omega I$, then the elementary indices $\bar a_i = v\(\bar\rho_i\)$ can be computed from
\[
  \bar a_1 + \cdots + \bar a_k = \min \left\{v(\eta) : \eta \in \Lambda^k(\omega I)\right\},
\]
the wedge power being taken in $\Lambda^k(\bar K^n) \isom \bar K^{\binom{n}{k}}$.
\end{lem}
\begin{proof}
Obvious.
\end{proof}
So the first index $\bar a_1 = a_1$ is unchanged by extension (this is rather obvious), while the sum of all the indices computes a discriminant, so
\begin{equation} \label{eq:disc_defect}
  (\bar a_1 + \bar a_2 + \bar a_3) - (a_1 + a_2 + a_3) = \frac{d_0 - 1}{2}.
\end{equation}

It remains to compute $\bar a_2$:

\begin{lem} \label{lem:bar_a_2}
Let $(\rho_i)_i$ be a well-led reduced basis of $\omega I$.
\begin{enumerate}[$($a$)$] 
  \item If $\rho_1$ or $\rho_2$ is $K$-led, then $\bar a_2 = a_2$, and an extender basis for $\omega I$ is
  \[
    \(\bar\rho_1, \bar\rho_2, \bar\rho_3\) = \(\rho_1, \rho_2, \pi^{a_3 + \frac{d_0 - 1}{2}} (0; 0; 1)\).
  \]
  \item If $\rho_3$ is $K$-led, then
  \begin{equation}
    \label{eq:bar_a_2}
    \bar a_2 = \min \left\{a_2 + v^{(K)}(\xi_1),\; a_2 + v^{(K)}(\xi_2),\; a_2 + \frac{d_0 - 1}{2},\; a_3\right\},
  \end{equation}
  and there is an extender basis $(\bar\rho_1, \bar\rho_2, \bar\rho_3)$ for $\omega I$ of one of the following shapes:  \begin{enumerate}[$($i$)$]
    \item If $\bar a_2$ is one of the first three arguments to the minimum in \ref{eq:bar_a_2}, then we can take
    \[
      \bar\rho_1 = \rho_1, \quad \bar\rho_2 = \rho_2 + u \rho_1
    \]
    for some $u \in \OO_{\bar K}^\cross$.
    \item If $\bar a_2 = a_3$, then we can take
    \[
      \bar\rho_1 = \rho_1, \quad \bar\rho_2 = \rho_3.
    \]
  \end{enumerate}
\end{enumerate}
\end{lem}

\begin{proof}
We denote the three places of $\bar K^3$ by $(K)$, $(Q1)$, and $(Q2)$ according to their relation to $R$. By Lemma \ref{lem:indices_wedge}, $\bar a_1 + \bar a_2$ is the minimum of the valuations
\[
  v(\rho_i^{(k)} \rho_j^{(\ell)} - \rho_j^{(\ell)} \rho_i^{(k)})
\]
of the nine $2 \times 2$ minors of the basis matrix
\[
  \begin{bmatrix}
    \rho_1^{(K)} & \rho_1^{(Q1)} & \rho_1^{(Q2)} \\
    \rho_2^{(K)} & \rho_2^{(Q1)} & \rho_2^{(Q2)} \\
    \rho_3^{(K)} & \rho_3^{(Q1)} & \rho_3^{(Q2)}
  \end{bmatrix}.
\]
If $\rho_1$ is $K$-led, then $\rho_2$ is $Q$-led and, in view of the valuation
\[
  v(\rho_1^{(K)} \rho_2^{(Q1)} - \rho_1^{(Q1)} \rho_2^{(K)})
  = v(\rho_1^{(K)} \rho_2^{(Q1)}) = a_1 + a_2
\]
of the upper left minor, we get $\bar a_2 \leq a_2$. But we know that $\bar a_2 \geq a_2$, so equality holds. In particular, $\rho_1 \wedge \rho_2$ is primitive in $\Lambda^2(\omega I) \tensor_{\OO_K} \OO_{\bar K}$, so we can keep $\rho_1$ and $\rho_2$ as our first two vectors in the extender basis. As for the last basis vector $\pi^{\bar a_3} \xi_3$, we already know that $\bar a_3 = a_3 + (d_0 - 1)/2$, and we can pick any $\xi_3 \in \OO_{\bar K}^3$ coprimitive with $\xi_1$ and $\xi_2$ without changing the span
\[
  \<\rho_1, \rho_2, \bar\rho_3\> = \<\rho_1, \rho_2, \pi^{\bar a_3}\OO_{\bar K}^3 \>.
\]
Since $\xi_1$ and $\xi_2$ necessarily lie in the sublattice
\[
  \left\{\(\xi^{(K)}, \xi^{(Q1)}, \xi^{(Q2)}\) \in \OO_{\bar K}^3 : \quad \xi^{(Q1)} \equiv \xi^{(Q2)} \mod \pi^{\frac{d_0 - 1}{2}} \right\},
\]
we get that $\xi_3 = (0; 0; 1)$ is a suitable coprimitive vector.

A similar argument works if $\rho_2$ is $K$-led.

If $\rho_3$ is $K$-led, we look at each of the minors in turn:
\begin{itemize}
  \item For $\rho_1^{(K)} \rho_2^{(Qi)} - \rho_1^{(Qi)} \rho_2^{(K)}$, there cannot be any cancellation, because the $K$-valuations differ by an integer (if both are finite) while the $Q$-valuations must differ by an integer plus $1/2$ to make the leaders a reduced basis of $Q$. So
  \begin{align*}
    v\(\rho_1^{(K)} \rho_2^{(Qi)} - \rho_1^{(Qi)} \rho_2^{(K)}\)
    &= \min\left\{v\(\rho_1^{(K)} \rho_2^{(Qi)}\), v\(\rho_1^{(Qi)} \rho_2^{(K)}\)\right\} \\
    &= \min\left\{a_1 + a_2 + v^{(K)}(\xi_1), a_1 + a_2 + v^{(K)}(\xi_2)\right\},
  \end{align*}
  accounting for the first two arguments to the minimum.
  \item The minor $\rho_1^{(Q1)} \rho_2^{(Q2)} - \rho_1^{(Q2)} \rho_2^{(Q1)}$ measures the discriminant of $\<\rho_1^{(Q)}, \rho_2^{(Q)}\>$ and hence has valuation
  \[
    a_1 + a_2 + \frac{d_0 - 1}{2},
  \]
  accounting for the third argument to the minimum. The other minors involving the two $Q$-valuations are at least as large and so can be ignored.
  \item Finally, the minor $\rho_1^{(K)} \rho_3^{(Qi)} - \rho_1^{(Qi)} \rho_3^{(K)}$ has valuation
  \[
    v(\rho_1^{(K)} \rho_3^{(Qi)} - \rho_1^{(Qi)} \rho_3^{(K)})
    = v(\rho_1^{(K)} \rho_3^{(Qi)}) = a_1 + a_3,
  \]
  accounting for the last argument to the minimum. The minor $\rho_2^{(K)} \rho_3^{(Qi)} - \rho_2^{(Qi)} \rho_3^{(K)}$ has valuation $a_2 + a_3 > a_1 + a_3$ and can be ignored.
\end{itemize}
This completes the determination of $\bar a_i$. We now go about constructing a suitable extender basis $\(\bar \rho_i\)$. We take $\bar \rho_1 = \rho_1$. The heart of the matter lies in finding a suitable $\bar \rho_2$.

If $\rho_1 \wedge \rho_2$ achieves the minimal valuation $\bar a_1 + \bar a_2$, that is, if the minimum corresponds to one of the first three arguments in \eqref{eq:bar_a_2}, then dividing by $\pi^{a_1 + a_2}$, we get
\[
  \xi_1 \wedge \xi_2 \equiv 0 \mod \pi^{\bar a_2 - a_2}.
\]
Since the $\xi_i$ are primitive, we get that for some $u \in \OO_{\bar K}^\cross$,
\begin{equation} \label{eq:x_bar_a_2}
  \xi_2 \equiv u \xi_1 \mod \pi^{\bar a_2 - a_2}.
\end{equation}
Then $\bar \rho_2 = \rho_2 - u \rho_1$ satisfies the needed divisibility
\[
  v\(\bar \rho_2\) \geq \bar a_2,
\]
and since $v\(\bar \rho_1 \wedge \bar \rho_2\) = \bar a_1 + \bar a_2$, equality must hold, that is, the primitive normalizations
\[
  \bar \xi_1 = \frac{\bar \rho_1}{\pi^{\bar a_1}} \textand
  \bar \xi_2 = \frac{\bar \rho_2}{\pi^{\bar a_2}}
\]
are linearly independent modulo $\mm_{\bar K}$. Then we can take any $\bar \xi_3$ independent from the previous two to complete the basis, as $\bar \rho_3 = \pi^{\bar a_3} \xi_3$ will automatically lie in $\omega \bar I$.

If $a_3$ achieves the minimum in \eqref{eq:bar_a_2}, we use the same method, with $\rho_3$ in place of $\rho_2$. The exponent $\bar a_2 - a_3$ appearing in \eqref{eq:x_bar_a_2} vanishes, allowing us to take $u = 0$ and $\bar\rho_2 = \rho_3$.
\end{proof}
\end{wild}

\subsection{The extender basis of a cubic resolvent ring}
Let $C$ be a candidate resolvent for $\tt$-traced quartic rings, and let $C_{\tt}$ be the corresponding reduced resolvent, that is, the unique ring such that $C = \OO_K + \tt^2 C_{\tt}$. First, look at the reduced basis of $C_{\tt} \subset R$ as a lattice; and look at its extender basis, a basis of $\bar C_{\tt} = C_{\tt} \tensor_{\OO_K} \OO_{\bar{K}}^3$. Because $1 \in C_{\tt}$ is an element of minimal valuation, there are not so many cases:
\begin{itemize}
  \item If $R$ is tamely ramified, then the reduced basis
  \[
    C_{\tt} = \< 1, \pi^{b_1} \theta_1, \pi^{b_2} \theta_2 \>, \quad b_1 \leq b_2
  \]
  is also an extender basis. If $R$ is unramified, the $b_i$ are of course integers; if $R$ is tamely ramified, then by Proposition \ref{prop:well_led}, we have
  \begin{equation} \label{eq:b_i_1^3}
    \left\{b_1, b_2\right\} \equiv \left\{\frac{1}{3}, \frac{2}{3}\right\} \mod \ZZ.
  \end{equation}
\begin{wild}
    \item If $R \isom K \cross Q$ is partially wildly ramified, then by Lemma \ref{lem:bar_a_2} we have a reduced basis
  \[
    C_{\tt} = \< 1, \pi^{b_1} \theta_1, \pi^{b_2} \theta_2 \>
  \]
  and an extender basis
  \[
    \bar{C}' = \< 1, \pi^{b_1} \theta_1, \pi^{\bar b_2} \bar \theta_2 \>
  \]
  with
  \[
    \bar b_2 = b_2 + \frac{d_0 - 1}{2} \textand \bar \theta_2 = (0;0;1).
  \]
  We have $b_1, b_2, \bar b_2 \in \frac{1}{2} \ZZ$. By Proposition \ref{prop:well_led}, the two $b_i$ are of different classes modulo $\ZZ$, so
  \begin{equation} \label{eq:b_i_1^21}
    \bar b_2 \equiv b_1 + \frac{d_0}{2} \mod \ZZ.
  \end{equation}
\end{wild}
\end{itemize}

Because $\theta_1$ is coprimitive to $1$, at most one pair of its three coordinates can be congruent modulo $\mm_{\bar K}$. We let $\bar s = \bar s_C$, the \emph{idempotency index} of $C$. Note that $\bar s$ is infinite only when two coordinates of $\theta_1$ are exactly equal. Since the $\theta_i$ are determined only up to finite precision, we can, and will, assume that $\bar s$ is finite.

\begin{lem}\label{lem:s}
  If finite, the value of $\bar s$ is constrained as follows:
  \begin{itemize}
    \item If $R$ is unramified, then $\bar s$ is an integer. For simplicity we let $s = \bar s$.
    \item If $R$ has splitting type $1^3$ (residue characteristic $2$), then $\bar s = 0$.
\begin{wild}
      \item If $R = K \cross Q$ has splitting type $1^21$, then either
    \begin{itemize}[\ensuremath{\blacktriangleright}]
      \item $b_1 \in \ZZ$ and $\bar s = \frac{d_0}{2} + s'$ for some $s' \in \ZZ_{\geq 0}$
      \item $b_1 \in \ZZ + \frac{1}{2}$ and $\bar s = \frac{d_0 - 1}{2}$.
    \end{itemize}
    In particular, $b_1 + \bar s \equiv \frac{d_0}{2} \mod \ZZ$.
\end{wild}
  \end{itemize}
\end{lem}
\begin{proof}
  In the tame splitting types this is immediate, knowing that $b_0 = 1$ is an extender vector for $C$ of minimal valuation.
\begin{wild}
    Assume that $R = K \cross Q$ is wildly ramified. Since 
  \[
  \theta_1 = \(\theta_1^{(K)}, \theta_1^{(Q1)}, \theta_1^{(Q2)}\)
  \]
  is coprimitive to $1$ and the two $Q$-components are already congruent mod $\pi^{\frac{d_0 - 1}{2}}$, we must have
  \[
  v\(\theta_1^{(K)} - \theta_1^{(Q1)}\) = 0 \textand \bar s = v\(\theta_1^{(Q1)} - \theta_1^{(Q2)}\) \geq \frac{d_0 - 1}{2}.
  \]
  If $b_1 \in \ZZ$, then $\theta_1 \in \OO_R$. Let $\I$ be the linear functional on $Q$ defined by
  \[
  \I(\xi) = \frac{\xi^{(Q1)} - \xi^{(Q2)}}{\pi_Q^{(Q1)} - \pi_Q^{(Q2)}},
  \]
  so that
  \[
  \I(1) = 0, \quad \I(\pi_Q) = 1.
  \]
  Then
  \[
  \bar s = \frac{d_0}{2} + v\(\I(\theta_1)\) = \frac{d_0}{2} + s',
  \]
  where $s' \geq 0$. If $b_1 \in \ZZ + 1/2$, then $\theta_1 \in \sqrt{\pi} \OO_R$, and in order for $\theta_1$ to be primitive, we must have $\sqrt{\pi}\theta_1^{(Q)} \sim \pi_Q$, so
  \[
  \bar s = \frac{d_0 - 1}{2} + v\(\I(\sqrt{\pi}\theta_1^{(Q)})\) = \frac{d_0 - 1}{2},
  \]
  as desired.
\end{wild}
\end{proof}

\begin{nota}\label{nota:coords}
If $\bar s > 0$, then there is a unique coordinate of $\bar K^3 \isom R \tensor \bar K$ at which $\omega_C$ has positive valuation. This defines a splitting $R \isom K \cross Q$ into a linear and a quadratic (possibly split) factor. We denote the three coordinates of $R$ by $(K)$, $(Q1)$, $(Q2)$, where $(K)$ is the distinguished one; thus we can write an element $\xi \in \bar K^3$ as
\[
  \xi = \(\xi^{(K)}; \xi^{(Q1)}; \xi^{(Q2)}\) = \(\xi^{(K)}; \xi^{(Q)}\)
\]
where $\xi^{(Q)} = \(\xi^{(Q1)}; \xi^{(Q2)}\) \in \bar K^2 = Q \tensor_K \bar K$.
\end{nota}

A common tool in understanding nonmaximal orders is their duals under the trace pairing. Hence it is fitting that we should understand the element $\omega_C \in \OO_{\bar K}^3$, unique up to scaling, that satisfies the relations
\[
  \tr \omega_C = \tr (\theta_1 \omega_C) = 0;
\]
that is, $\bar K \<\omega \> = \bar K \<1, \theta_1\>^\perp$ under the trace pairing.

One explicit choice of $\omega_C$ is as follows: If $\tilde\theta = \pi^{b_1}\theta_1 = (\tilde\theta^{(1)}, \tilde\theta^{(2)}, \tilde\theta^{(3)}) \in C$ is the second reduced basis vector, then
\[
  \hat{\omega}_C =
  \(\tilde\theta^{(2)} - \tilde\theta^{(3)}\) \(\tilde\theta^{(3)} - \tilde\theta^{(1)}\) \(\tilde\theta^{(1)} - \tilde\theta^{(2)}\) \(\tilde\theta^{(2)} - \tilde\theta^{(3)}; \tilde\theta^{(3)} - \tilde\theta^{(1)}; \tilde\theta^{(1)} - \tilde\theta^{(2)}\).
\]
The symmetry ensures that $\hat{\omega}_C \in R$. Note that $N\(\hat{\omega}_C\) = (\disc \tilde\theta)^2$ is a square in $K^\cross$. Note also that 
\[
  \tr\(\omega_C\) = \tr\(\theta_1\omega_C\) = 0.
\]
Then 
\[
\vec{v}(\hat{\omega}_C) = (4b_1+2\bar s, 4b_1+\bar s, 4b_1+\bar s).
\]
Two other rescalings of $\hat{\omega}_C$, chosen for primitivity rather than the property of lying in $R$, will also be used:
\begin{align*}
\omega_C &= \pi^{-4b_1-\bar s} \omega_C, & \vec{v}(\omega_C) &= (\bar s, 0, 0) \\
\bar\omega_C &= \pi^{-4b_1-2\bar s} \omega_C, & \vec{v}(\bar{\omega}_C) &= (0, -\bar s, -\bar s).
\end{align*}
They have the properties that $\omega_C$ and $\bar\omega_C^{-1}$ are primitive (that is, have valuation $0$) in $\OO_{\bar K}^{\oplus 3}$.

\section{Resolvent conditions}
Let $C,C_{\tt} \subseteq R$ be a resolvent cubic ring and its corresponding reduced resolvent, whose reduced bases are related by
\[
C = \<1, \pi^{b_1 + 2t} \theta_1, \pi^{b_2 + 2t} \theta_2\> \textand
\bar{C}_\tt = \< 1, \pi^{b_1} \theta_1, \pi^{b_2} \theta_2 \>, \quad 0 \leq b_1 \leq b_2,
\]
where $t = v_K(\tt)$. Let $L$ be a quartic algebra with resolvent $R$. As we noted in the proof of Theorem \ref{thm:hcl_quartic}, an order $\OO \subseteq L$ is completely determined by the lattice $I \subseteq R$ such that $\OO/\OO_K = \kappa(I)$, where
\begin{align*}
  \kappa : K &\to \bar K^4 \\
  \xi &\mapsto \(\tr_{\bar K^3/K} \xi \omega \sqrt{\delta}\)_\omega
\end{align*}
is the map in Theorem \ref{thm:Kummer_new}\ref{it:Kum_quartic}. For reasons that will become clear below, we take the reduced and extender bases, not of $I$ itself, but of
\[
\sqrt{\pi^{\bar s} \delta {\omega}_C} I.
\]
Let the reduced basis be
\[
\sqrt{\pi^{\bar s} \delta {\omega}_C} \cdot I = \<\pi^{a_1} \xi_1, \pi^{a_2} \xi_2, \pi^{a_3} \xi_3 \>.
\]
and let the extender basis be
\[
 \sqrt{\pi^{\bar s} \delta {\omega}_C} \cdot \bar{I} = \<\pi^{\bar a_1} \bar \xi_1, \pi^{\bar a_2} \bar \xi_2, \pi^{\bar a_3} \bar \xi_3 \>.
\]
Both $\delta$ and ${\omega}_C$ satisfy property \ref{iota:*} in Definition \ref{defn:red_basis}, so the foregoing theory applies. (In this edition, since $R$ will always be tamely ramified, the overbars can be ignored.)
\begin{wild}
In the tamely ramified splitting types, the reduced and extender bases coincide by Lemma \ref{prop:ext_basis_tame} and the bars can be dropped. On the other hand, in splitting type $1^21$, only the first basis vector $\pi^{a_1} \xi_1 = \pi_{\bar a_1} \bar \xi_1$ is held in common in general, and we will need both bases.
\end{wild}
We will take
\[
  a_i' = a_i + 2b_1 \textand \bar a'_i = \bar a_i + 2b_1,
\]
the reduced and extender indices of $\sqrt{\delta \hat\omega_C} I$.
Then we can write the resolvent conditions as follows:
\begin{lem} \label{lem:rsv}
  With respect to the above setup, a lattice $I$ yields a ring $\OO$ with a $\tt$-traced resolvent to $\C$ if and only if the following conditions hold:
  \begin{itemize}
    \item Discriminant $(\Theta)$ condition:
    \begin{equation}
      \bar a_{1} + \bar a_{2} + \bar a_{3} = \bar b_{1} + \bar b_{2} + 2\bar s + 4t - 4e \label{eq:disc_cond_quartic}
    \end{equation}
    \item Resolvent $(\Phi)$ conditions for $\bar\theta_2$-coefficients:
    \[
      \M_{ij} \colon \quad  \tr(\bar\xi_i \bar\xi_j) \equiv 0 \mod \pi^{\bar m_{ij}}
    \]
    where
    \[
      \bar m_{ij} = \begin{cases}
        \bar b_2 + 2t - 2e - 2\bar a_i + \bar s, & i = j \\
        \bar b_2 + 3t - 3e - \bar a_i - \bar a_j + \bar s, & i \neq j.
      \end{cases}
    \]
    \item Resolvent $(\Phi)$ conditions for $\bar\theta_1$-coefficients:
    \[
      \N_{ij} \colon \quad \text{All coordinates of} \quad \bar\omega_C^{-1} \cdot \bar\xi_i \bar\xi_j \quad \text{are congruent} \mod \pi^{\bar n_{ij}},
    \]
    where
    \[
      \bar n_{ij} = \begin{cases}
      \bar b_1 + 2t - 2e - 2\bar a_i + 2\bar s, & i = j \\
      \bar b_1 + 3t - 3e - \bar a_i - \bar a_j + 2\bar s, & i \neq j.
      \end{cases}
    \]
  \end{itemize}
\end{lem}
\begin{proof}
The conditions that $C$ is a $\tt$-traced resolvent for $\OO$ are that all coefficients in the coordinate representations of $\Theta$, $\Theta^{-1}$, and $\Phi$ have nonnegative valuation. In particular, it is equivalent to study when $\bar{C} = C \tensor_{\OO_K} \OO_{\bar{K}}$ is a resolvent of $\bar{\OO} = \OO \tensor_{\OO_K} \OO_{\bar{K}}$.

We have
\[
\bar{C} = \< 1, \pi^{2t + \bar b_1} \bar\theta_1, \pi^{2t + \bar b_2} \bar\theta_2 \>,
\]
so
\[
  \Lambda^2 (\bar{C}/\OO_{\bar{K}}) = \< \pi^{4t + \bar b_1 + \bar b_2} \bar\theta_1 \wedge \bar\theta_2 \>
  = \pi^{4t + \bar b_1 + \bar b_2} \Lambda^2 (\OO_{\bar{K}}^3 / \OO_{\bar{K}})
\]
and, by the formula for $\Theta$ in Proposition \ref{prop:Kummer_resolvent_quartic},
\[
  \Theta(\Lambda^2 (\bar{C}/\OO_{\bar{K}})) = \frac{1}{16\sqrt{N(\delta)}} \pi^{4t + \bar b_1 + \bar b_2} \cdot \kappa(\OO_{\bar{K}}^3)
\]
Meanwhile,
\begin{equation} \label{eq:I_basis_copy}
  \bar{I} = \frac{1}{\sqrt{\pi^{\bar s} \delta \omega_C}}\<\pi^{\bar a_1} \bar\xi_1, \pi^{\bar a_2} \bar\xi_2, \pi^{\bar a_3} \bar\xi_3 \>
\end{equation}
so
\begin{align*}
  \Lambda^3 \bar{I} &= \frac{1}{\sqrt{\pi^{3 \bar s} N(\delta \omega_C)}} \pi^{\bar a_1 + \bar a_2 + \bar a_3} \<\bar\xi_1 \wedge \bar\xi_2 \wedge \bar\xi_3\> \\
  &= \frac{1}{\sqrt{\pi^{3\bar s} N(\delta \omega_C)}} \pi^{\bar a_1 + \bar a_2 + \bar a_3} \Lambda^3(\OO_{\bar{K}}^4/ \OO_K).
\end{align*}
So the condition for $\Theta$ to define an isomorphism between $\Lambda^2(\bar{C}/\OO_{\bar{K}})$ and $\Lambda^3(\bar{\OO}/\OO_{\bar{K}})$ is that
\begin{align*}
4t + \bar b_1 + \bar b_2 - 4e - \frac{1}{2} v(N(\delta)) = \bar a_1 + \bar a_2 + \bar a_3 - \frac{3}{2}\bar s - \frac{1}{2} v(N(\delta\omega_C)),
\end{align*}
or, since $v(N(\omega_C)) = \bar s$,
\[
  4t - 4e + \bar b_1 + \bar b_2 + 2\bar s = \bar a_1 + \bar a_2 + \bar a_3,
\]
as desired.

Likewise, we use the formula
\[
  \Phi(\kappa(\alpha)) = 4\delta\alpha^2
\]
from Proposition \ref{prop:Kummer_resolvent_quartic} to transform the $\Phi$-condition to the following:
\begin{enumerate}[(i)]
  \item\label{xon} For every $\xi \in \bar{I}$, we have $4 \delta \xi^2 \in \bar{C} + \bar{K}$
  \item\label{xoff} For every $\xi,\eta \in \bar{I}$, we have $8 \pi^{-t} \delta \xi \eta \in \bar{C} + \bar{K}$.
\end{enumerate}
In terms of the basis \eqref{eq:I_basis_copy} for $\bar{I}$, this is to say that the diagonal entries of the matrix representing $\Phi$ belong to $\bar{C} + \bar{K}$ and the off-diagonal entries to $\tt\bar{C} + \bar{K}$. Hence it suffices to consider \ref{xon} for $\xi$ a basis element and \ref{xoff} for $\xi,\eta$ distinct basis elements.

To test whether an $\alpha \in \bar{K}^3$ lies in
\[
  \bar{C} = \< 1, \pi^{2t + \bar b_1} \bar\theta_1, \pi^{2t + \bar b_2} \bar\theta_2 \>,
\]
we can pair it with a basis of the dual lattice $\bar{C}^\vee$ with respect to the trace pairing. Let $(\lambda_0, \lambda_1, \lambda_2)$ be the dual basis to $(\bar\theta_0, \bar\theta_1, \bar\theta_2)$ (that is, $\tr(\bar\theta_i \lambda_j) = \1_{i = j}$). Then $(\lambda_1,\lambda_2)$ is a basis for $(\OO_K^3)^{\tr = 0}$, and we have already met $\lambda_2$: it is $\omega_C$, up to a unit. Hence
\begin{align*}
  \bar{C}^\vee &= \<\lambda_0, \pi^{-2t - \bar b_1} \lambda_1, \pi^{-2t - \bar b_2} \lambda_2 \> \\
  &= \<\lambda_0\> + \pi^{-2t - \bar b_1}(\OO_{\bar{K}}^3)^{\tr = 0} + \pi^{-2t - \bar b_2} \<{\omega}_C\>
\end{align*}
We actually wish to test not whether $\alpha \in \bar{C}$, but the weaker condition $\alpha \in \bar{C} + \bar{K}$, so (due to the natural duality between $\bar{C}/\OO_{\bar{K}}$ and $(\bar{C}^\vee)^{\tr = 0}$) we pair only with elements of
\[
  (\bar{C}^\vee)^{\tr = 0} = \pi^{-2t - \bar b_1}(\OO_{\bar{K}}^3)^{\tr = 0} + \pi^{-2t - \bar b_2} \<\hat{\omega}_C\>.
\]
Hence $\alpha \in \bar{C} + \bar{K}$ if and only if
\begin{itemize}
  \item $\pi^{-2t - \bar b_2}\tr({\omega}_C\alpha)$ is integral, and
  \item $\pi^{-2t - \bar b_1}\tr(\alpha \kappa)$ is integral for $\kappa \in \{(1;-1;0), (0; 1; -1)\}$; that is, all coordinates of $\alpha$ are congruent modulo $\pi^{2t + \bar b_1}$.
\end{itemize}
This is the origin of the $\M$- and $\N$-conditions respectively. Applying this to the values $\alpha = 4\delta\xi^2, 8\pi^{-t}\delta\xi\eta$ derived from the basis above yields the desired form of all the $\Phi$-conditions.
\end{proof}

We say that the condition $\M_{ij}$ or $\N_{ij}$ is \emph{active} if its corresponding modulus $\bar m_{ij}$ or $\bar n_{ij}$ is positive. An inactive condition is automatically satisfied (noting that $\bar\omega_C^{-1}$ has nonnegative valuations).

Because the $a$'s and $b$'s have been sorted in increasing order, and because $0 \leq t \leq e$, we have the following implications among the activity of the $\M_{ij}$ and $\N_{ij}$:
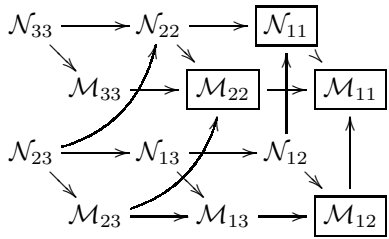
\[
\xymatrix@!0{
  \N_{33} \ar[rr] \ar[rd] & & \N_{22} \ar[rr] \ar[rd] & & \boxed{\N_{11}} \ar[rd] \\
  & \M_{33} \ar[rr] & & \boxed{\M_{22}} \ar[rr] & & \boxed{\M_{11}} \\
  \N_{23} \ar[rr] \ar[rd] \ar@/_3ex/[urur] & & \N_{13} \ar[rr] \ar[rd] & & \N_{12} \ar[rd] \ar[uu] \\
  & \M_{23} \ar[rr] \ar@/_3ex/[urur] & & \M_{13} \ar[rr] & & \boxed{\M_{12}} \ar[uu]
}
\]
The next lemma limits our concern to the four boxed conditions:
\begin{lem} \label{lem:inactives}
    \begin{enumerate}[$($a$)$]
        \item\label{inact:quart} Suppose that the $\bar a_i$ and $\bar\xi_i$ come from the extender decomposition of a quartic ring and a resolvent thereof. Then:
        \begin{itemize}
          \item No conditions are active except $\M_{11}$, $\M_{12}$, $\M_{22}$, $\N_{11}$, and $\N_{12}$.
          \item $\M_{12}$ and $\M_{22}$ are not both active.
          \item $\N_{12}$ is \emph{very weakly active,} that is, $\bar n_{12} \leq \bar s/2$.
        \end{itemize}
        \item\label{inact:weak} Suppose that the $\bar a_i$ and $\bar\xi_i$ come from the extender decomposition of \underline{some} lattice $I \subseteq R$. Suppose that the inactivity restrictions from part \ref{inact:quart} hold and that conditions $\M_{11}$, $\M_{12}$, $\M_{22}$, and $\N_{11}$ are satisfied. Then $\N_{12}$ is satisfied, and the $\bar a_i$ and $\bar\xi_i$ actually come from a quartic ring. That is, ``$\N_{12}$ is automatic if it is very weakly active.''
    \end{enumerate}
\end{lem}
\begin{proof}
\begin{enumerate}[$($a$)$]
  \item
Suppose that the $\bar a_i$ and $\bar\xi_i$ come from a quartic ring.

If $\M_{13}$ is active, so are $\M_{12}$ and $\M_{11}$. Since the $\bar\xi_i$ are supposed to form an $\OO_{\bar{K}}$-basis of $\OO_{\bar K}^3$, we obtain for all $\bar\xi \in \OO_{\bar{K}}^3$,
\[
  \tr(\bar\xi_1\bar\xi) \equiv 0 \mod \mm_{\bar K}.  
\]
Since $\bar\xi_1$ is primitive, this is a contradiction.

If $\M_{12}$ and $\M_{22}$ are active, then so is $\M_{11}$. We have a $2$-dimensional subspace $V = \<\bar\xi_1, \bar\xi_2\>$ of the $3$-dimensional space $k_{\bar K}^3$ that is isotropic for the trace pairing. But the trace pairing is nondegenerate, so this is a contradiction.

If $\M_{33}$ is active, so are $\M_{22}$ and $\M_{11}$. If $\ch k_K \neq 2$, then $t = e = 0$ so $\M_{12}$ is also active, and we have a contradiction as above. If $\ch k_K = 2$, we use that squaring is a linear operation mod $2$ to obtain that for all $\bar\xi \in \OO_{\bar{K}}^3$,
\[
  \tr(\bar\xi^2) \equiv 0 \mod \mm_{\bar K},
\]
which is a contradiction.

If $\N_{22}$ is active, note that $\ch k_K = 2$ since $\M_{22}$ is active. There are two cases. If $\bar s = 0$, then $\bar{\omega}_C^{-1}$ is a unit, so the condition
\[
  \N_{ii} \colon \quad \text{All coordinates of} \quad \bar\omega_C^{-1} \cdot \bar\xi_i^2 \quad \text{are congruent} \mod \pi^{\bar n_{ii}}
\]
determines $\bar\xi_i^2$ mod $\mm_{\bar K}$ up to scaling. But since we are in characteristic $2$, square roots are unique, and $\bar\xi_1$ and $\bar\xi_2$ are scalar multiples mod $\mm_{\bar K}$, a contradiction. Now assume $\bar s > 0$, so $\bar\omega_C^{-1} \equiv (u;0;0)$ for some unit $u \in \OO_{\bar{K}}^\cross$. Now $\N_{ii}$ gives $\bar\xi_i^{(K)} \equiv 0 \mod \mm_{\bar K}$. But now $\M_{ii}$ gives that $\bar\xi_i$ is a unit multiple of $(0;1;1)$ modulo $\mm_{\bar K}$, a contradiction.

Finally, assume that $\N_{12}$ is active and not very weakly active: that is, $\bar n_{12} > \bar s/2$. Note that $\bar n_{11} > \bar s$ since otherwise $\N_{22}$ would be active. If $\bar s = 0$, then $\N_{11}$ implies that $\bar\xi_1^2 \equiv \omega_C \mod \mm_C$, up to scaling, and then $\N_{12}$ implies that $\bar\xi_1\bar\xi_2 \equiv \omega_C \mod \mm_C$, up to scaling. Since $\omega_C$ is a unit, this is a contradiction. So assume $\bar s > 0$. Note that $\omega_C$ is a unit multiple of $(0;1;-1)$ modulo $\pi^{\bar s}$ and that $\bar{\omega}_C^{-1}$ is a unit multiple of $(1;0;0)$ modulo $\pi^{\bar s}$. Now $\N_{11}$ implies that
\[
  \bar\xi_1^2 \equiv \omega_C \mod \mm_{\bar K}\omega_C = \mm_{\bar K}(\pi^{\bar s};1;1).
\]
So $v_K(\bar\xi_1^{(K)}) = \bar s/2$ exactly (recalling the notion of $\xi^{(K)}$ from Notation \ref{nota:coords}). If $v_K(\bar\xi_2^{(K)}) > 0$ also, we get $\bar\xi_1 \equiv \bar\xi_2 \mod \mm_{\bar K}$ by the same argument as when $N_{22}$ is active. So $\bar\xi_2^{(K)}$ is a unit, and the $\N_{12}$ condition
\[
  \N_{12} \colon \quad \text{All coordinates of} \quad \bar\omega_C^{-1} \cdot \bar\xi_1 \bar\xi_2 \quad \text{are congruent} \mod \pi^{\bar n_{12}}
\]
is unsatisfied, because the $K$-coordinate has valuation $\bar s/2$ and the others have higher valuation.

\item
Note that $\bar n_{12} \leq \frac{1}{2} \bar n_{11}$ because otherwise $\N_{22}$ would be active. We have $\bar{\omega}_C^{-1}$ a unit multiple of $(1;0;0)$ modulo $\pi^{2\bar n_{12}}$, so $\N_{11}$ implies that $\pi^{2\bar n_{12}} | (\bar\xi_1^2)^{(K)}$, that is, $\pi^{\bar n_{12}} | \bar\xi_1^{(K)}$. Now $\pi^{\bar n_{12}} \mid \bar \omega_C^{-1} \bar\xi_1$, so $\N_{12}$ is satisfied. \qedhere
\end{enumerate}
\end{proof}

Based on this, we will count quartic rings with fixed resolvent. Since a lot will happen with various things being fixed and others varying, it is worthwhile to lay down the following:

\begin{conv}\label{conv:quartic}
We fix variables in the following order:
\begin{itemize}
  \item First, we fix the \emph{resolvent data,} which comprise
  \begin{itemize}
    \item a resolvent $C \subseteq R$;
    \item an extender decomposition $C = \<1, \pi^{\bar b_1} \bar\theta_1, \pi^{\bar b_2} \bar\theta_2\>$, which fixes $\hat\omega_C$ and $\bar s$. We can, and do, assume that $\bar s$ is finite;
    \item a tracedness parameter $t$, $0 \leq t \leq e$, which defines a reduced resolvent $C_{\tt}$.
  \end{itemize} 
  \item Then we fix the \emph{discrete data} of a quartic ring, which comprises
  \begin{itemize}
    \item a choice of \emph{coarse coset} $[\delta] \in H^1/\L_0$. There are $\size{H^0}$ cosets $\delta_0 \L_0$. Then $\delta = \delta_0 \tau$, where $\tau \in \OO_R^\cross$ can vary;
    \item its extender indices $\bar a_i$,
\begin{wild}
      together with the regularity or irregularity of each extender vector,
\end{wild}
which are constrained by the integrality needed for a sublattice of $\sqrt{\delta \hat\omega_C} R$ and the inactivity inequalities of Lemma \ref{lem:inactives}.
  \end{itemize}
  \item Then we choose $\delta$ and $\bar\xi_1$, which are constrained by the $\M_{11}$ and $\N_{11}$ conditions.
  \item Then we choose $\bar\xi_2$, which is constrained by its coprimitivity with $\bar\xi_1$ and by the $\M_{12}$ and $\M_{22}$ conditions.
  \item Lastly, we choose $\bar\xi_3$, which is constrained by its coprimitivity with $\bar\xi_1$ and $\bar\xi_2$.
\end{itemize}
Whenever we speak about possibilities for any of the items on this list, it will be implicitly assumed (if not stated) that all the previous items have been fixed in conformity with their respective restrictions.
\end{conv}

Since $\bar\xi_1 \in \pi^{-\bar a_1'} \sqrt{\delta \hat\omega_C} R$,
conditions $\M_{11}$ and $\N_{11}$ can be viewed in another way, which will be simpler for some purposes:
\begin{lem} \label{lem:to_box}
A $\bar\xi_1 \in \sqrt{\delta \hat\omega_C} R \intsec \OO_{\bar{K}}^3$ satisfies the $\M_{11}$ and $\N_{11}$ resolvent conditions if and only if the quotient
\[
  \beta = \frac{\bar\xi_1^2}{\omega_C} \in \pi^{-2\bar a_1' + \bar s}R
\]
is a linear combination of the reduced basis vectors of $C$ of the form
\[
    \beta = x + y\pi^{n_{11} - \bar s} \bar\theta_1 + z\pi^{m_{11}} \bar\theta_2, \quad x \in \pi^{-2\bar a_1 + \bar s}K, \quad y,z \in \OO_K.
\]
\end{lem}
\begin{proof}
Since $\frac{\pi^{\bar a_1'} \bar\xi_1}{\sqrt{\delta\hat\omega_C}}$ is the first basis element for $\bar{I}$, the $\M_{11}$ and $\N_{11}$ conditions are equivalent to
\begin{equation*}
\bar{K} + \bar{C} \ni 4\delta\(\frac{\pi^{\bar a_1'} \bar\xi_1}{\sqrt{\delta\hat\omega_C}}\)^2 =
\frac{4\pi^{2\bar a_1' }\bar\xi_1^2}{\omega_C} = 4\pi^{2\bar a_1 + 4\bar b_1 - 4\bar b_1 - \bar s} \beta = 4\pi^{2\bar a_1 - \bar s} \beta,
\end{equation*}
that is,
\[
  \beta \in \bar{K} + \pi^{\bar s - 2\bar a_1 - 2 e} \bar{C} = \bar{K} + \pi^{n_{11} - \bar s} \OO_{\bar{K}} \bar\theta_1 + \pi^{m_{11}} \OO_{\bar{K}} \bar\theta_2.
\]
Since $(1, \pi^{\bar b_1}\bar\theta_1, \pi^{\bar b_2}\bar\theta_2)$ form an $\OO_K$-basis for $C$, we can always find $x \in \pi^{\bar s - 2\bar a_1} K$ and $y, z \in K$ such that
\[
  \beta = x + y\pi^{n_{11} - \bar s} \bar\theta_1 + z\pi^{m_{11}} \bar\theta_2.
\]
Then the resolvent conditions simplify to $y, z \in \OO_K$.
\end{proof}

\subsection{Transformation, and ring volumes in the white zone}
We will proceed to compute the volumes of the solution sets in which the reduced vectors $\xi_i$ lie. (We use reduced vectors $\xi_i$, not extender vectors $\bar\xi_i$, because the latter do not lie in a controllable $\OO_K$-lattice.) For simplicity, we will transform everything to $R$ itself, which we normalize so that $\OO_R$ has volume $1$, and to its projectivization $\PP(\OO_R)$, to which we give a volume of $1 + 1/q + 1/q^2$, so that a distinguished affine open has volume $1$.

The simplest way to do this is as follows.
\begin{lem}\label{lem:gamma_white}
Fix the discrete data. In particular, $\delta = \tau \delta_0$ lies in a fixed coarse coset.\begin{wild}
  If $\xi_i$ is regular, 
\end{wild}
There is a $\gamma = \gamma_{i} \in \bar{K}^3$ with the properties that, letting
\[
  \gamma_i = \frac{\gamma_{i,0}}{\sqrt{\tau}},
\]
we have that $\xi_i' = \gamma_{i}^{-1} \xi_i$ is a primitive vector in $\OO_R$. 
\end{lem}
\begin{proof}
Note that $\xi_i$ must lie in
\begin{equation} \label{eq:J_for_gamma}
  J = \pi^{-a_i'} \sqrt{\delta \hat\omega_C} R \intsec \OO_{\bar{K}}^3 = \sqrt{\tau} J_0,
\end{equation}
where
\[
  J_0 = \pi^{-a_i'} \sqrt{\delta_0 \hat\omega_C} R \intsec \OO_{\bar{K}}^3
\]
is an $\OO_R$-lattice of dimension $1$. As $\OO_R$ is a principal ideal ring (it's a product of DVR's), we obtain that $J_0 = \gamma_i \OO_R$ for some $\gamma_i$, clearly not a zero-divisor. Since $\xi_i \in J \setminus \pi J$, we get $\xi_i' = \gamma_i^{-1} \xi_i \in \OO_R \setminus \pi\OO_R$, as desired.
\end{proof}

The valuations of $\gamma$ may be computed by observing the smallest nonnegative valuation of an element of $J$ at each place.
The $\xi_i$ form a set of reduced vectors for a sublattice of $\sqrt{\delta \hat\omega_C}I$ if and only if the $\xi_i'$ lie in certain explicit subsets of $\OO_R$, computed below. 

\paragraph{Unramified.}
\begin{itemize}
  \item If $[\delta \hat\omega_C] \in \L_0$, then  all $a_i$ are integers, and $\gamma_i = \gamma$ is a unit. The three $\xi'$ must form a basis of $\OO_R$. If they are found successively, their ring volumes are respectively $\boxed{1 + 1/q + 1/q^2}$, $\boxed{1 + 1/q}$, and $\boxed{1}$.
  \item If $[\delta \hat\omega_C] \in (1; \pi; \pi) \L_0$, then  there is one $a_{j_0} \in \ZZ$; there $\gamma_{j_0} \sim (1; \sqrt{\pi}; \sqrt{\pi})$, and
  \[
    \xi'_{j_0} \in \OO_K^\cross \cross \OO_Q,
  \]
  a subset whose projectivization has volume $\boxed{1}$. Meanwhile, two $a_{j_1}, a_{j_2}$ lie in $\ZZ + 1/2$; there $\gamma_{j_k} \sim (\sqrt{\pi}; 1; 1)$, and
  \[
    \xi'_{j_k} \in \OO_K \cross \OO_Q
  \]
  with their $\OO_Q$-coordinates forming a basis of $\OO_Q$; thus $\xi'_{j_1}$ has volume $\boxed{1 + 1/q}$ and $\xi'_{j_2}$ has volume $\boxed{1}$.
\end{itemize}
\paragraph{Splitting type $1^3$.}
\begin{itemize}
  \item Here $[\delta \hat\omega_C] \in \L_0$. The $a_i$ fill out the classes in $\frac{1}{3}\ZZ / \ZZ$, but all $\gamma_i$ are units and all $\xi'_i$ lie in $\OO_R^\cross$, a subset whose projectivization has volume $\boxed{1}$.
\end{itemize}
\begin{wild}
  \paragraph{Splitting type $1^21$.} 
\begin{itemize}
  \item If $[\delta \hat\omega_C] \in \L_0$, then by Corollary \ref{cor:idxs_mod_1}, the $a_i'$ form the multiset $\{0, 0, 1/2\}$ modulo $1$. For the two $a_i \in \ZZ$, we have $\gamma_i \sim 1$, and the value of $\xi'_i$ must lie off the $1$-pixel of $(0; \pi_R$), yielding a ring volume of $\boxed{1 + 1/q}$ for the first such $i$ and $\boxed{1}$ for the second (one $1$-pixel being excluded by linear independence). For the one $a_i \in \ZZ + 1/2$, we have $\gamma_i \sim (\sqrt{\pi} ; 1)$, and $\xi'_i$ must have a unit $Q$-component, yielding a ring volume of $\boxed{1}$. 
  \item If $[\delta \hat\omega_C] \in (\pi; \pi_R)\L_0$, then by Corollary \ref{cor:idxs_mod_1}, the $a_i'$ form the multiset $\{1/4, 1/2, 3/4\}$ modulo $1$. The $\gamma_i$ in the three cases are associate to
  \[
    (\pi^{1/4}; 1; 1), \quad (1; \pi^{1/4}; \pi^{1/4}), \textand (\pi^{3/4}; 1; 1)
  \]  
  No matter what $a_i'$ is, it determines whether $\xi'_i$ must have a unit $K$-component or unit $Q$-component, yielding a ring volume of $\boxed{1}$ in all cases.
\end{itemize}
For future reference, it is worth noting that, in splitting type $1^21$,
\[
  v\big(N(\gamma))\big) = \{a_i'\}.
\]

Hence we get ring volumes in the white zone.
\end{wild}
\begin{lem}\label{lem:white}
In the white (i.e.~free) zone where no $\M_{ij}$ or $\N_{ij}$ is active, we take all $\gamma_i$ as in Lemma \ref{lem:gamma_white}. The ring volume of triples $(\xi'_1, \xi'_2, \xi'_3)$ is given in terms of the discrete data as follows:
\begin{enumerate}
  \item If $R$ is unramified and $[\delta \hat\omega_C] \in \L_0$, the ring volume is $(1 + 1/q + 1/q^2)(1 + 1/q)$.
  \item If $R$ is unramified and $[\delta \hat\omega_C] \notin \L_0$, the ring volume is $1 + 1/q$.
  \item If $R$ is totally ramified, the ring volume is $1$.
\begin{wild}
    \item If $R$ is partially ramified and $[\delta \hat\omega_C] \in \L_0$, the ring volume is $1 + 1/q$.
  \item If $R$ is partially ramified and $[\delta \hat\omega_C] \notin \L_0$, the ring volume is $1$.
\end{wild}
\end{enumerate}
\end{lem}

\subsection{From ring volumes to ring counts}

\begin{lem} \label{lem:white_vol}
Let $\OO$ be a quartic ring. The set $V$ of triples $(\xi'_1, \xi'_2, \xi'_3)$ in $\PP^2(\OO_R)^3$ whose associated ring is $\OO$ has a volume determined by the discrete data of $\OO$ alone. It is given by
\[
  \mu(V) = q^{\ds -\ceil{a_2 - a_1} - \ceil{a_3 - a_2} - \ceil{a_3 - a_1} + \frac{3d_0^\mathrm{tame}}{2} + v_K\big(N(\gamma_1\gamma_2\gamma_3)\big)} \cdot c,
\]
where
\[
d_0^{\mathrm{tame}} = \sum_{R_i} \(e_{R_i/K} - 1\)
\]
is the standard lower bound for the discriminant valuation, attained for tame extensions, and
\[
  c = \begin{cases}
  1 & \text{if } a_1 = a_2 = a_3 \\
  \left(1 + \frac{1}{q} + \frac{1}{q^2}\right) & \text{if } a_1 = a_2 < a_3 \text{ or } a_1 < a_2 = a_3 \\
  \left(1 + \frac{1}{q}\right)\left(1 + \frac{1}{q} + \frac{1}{q^2}\right) & \text{if } a_1 < a_2 < a_3. \\
  \end{cases}
\]
\end{lem}
\begin{proof}
We have
\[
  \xi_i, \gamma_i \in \pi^{-a_i'} \sqrt{\delta\hat\omega_C} R, \quad \xi_i' \in \OO_R.
\]

To determine whether $(\Xi_1', \Xi_2', \Xi_3') \in V$, there are two conditions: firstly, the generators $\pi^{a_i} \Xi_1$, where $\Xi_i = \Xi_i'\gamma_i$ are the associated reduced basis vectors, belong to the correct lattice
\begin{equation} \label{eq:Xi_in}
  \pi^{a_i} \Xi_i \in \sqrt{\pi^{\bar s} \delta \omega_C} I = \<\pi^{a_1} \xi_1, \pi^{a_2} \xi_2, \pi^{a_3} \xi_3\>;
\end{equation}
and secondly, they generate the whole of $\sqrt{\pi^{\bar s} \delta \omega_C} I$.
Since the $\pi^{a_i} \xi_i$ form a $K$-basis for $\sqrt{\delta\hat\omega_C} R$, we can write 
\begin{equation}\label{eq:c_to_Xi}
  \Xi_i = \sum_{j = 1}^3 \pi^{a_j - a_i} c_{ij} \xi_i
\end{equation}
for some coefficients $c_{ij} \in K$. Condition \eqref{eq:Xi_in} is then equivalent to
\[
  v\(c_{ij}\) \geq \max\{a_i - a_j, 0\},
\]
while the condition that the $\pi^{a_i} \Xi_i$ generate the whole of $\sqrt{\pi^{s'} \delta \omega_C} I$ is equivalent to the change of basis being invertible:
\[
  v\(\det\left[c_{ij}\right]\) = 0.
\]
Thus we have parametrized $V$ by the group
\[
  \Gamma = \left\{[c_{ij}] \in \GL_3(\OO_K) : v\(c_{ij}\) \geq a_i - a_j\right\}.
\]
More precisely, $V$ is in continuous bijection with the cosets $T \backslash \Gamma$, where $T = \(\OO_K^\cross\)^3$ is the subgroup of diagonal matrices, because the $\Xi_i \in \PP(\OO_R)$ are defined only up to scaling.

Without the invertibility condition, the volume of matrices in $\Mat_3(\OO_K)$ satisfying the valuation restrictions defining $\Gamma$ is
\[
  q^{\ds -\sum_{1 \leq i,j \leq 3} \max\{0, \ceil{a_{i} - a_{j}}\}} = q^{\ds -\ceil{a_2 - a_1} - \ceil{a_3 - a_2} - \ceil{a_3 - a_1}}.
\]
The invertibility depends only on the $c_{ij}$ modulo $\pi$, and the fraction of matrices over $k_K$ of the shapes
\[
  \begin{bmatrix}
    * & * & * \\
    * & * & * \\
    * & * & *
  \end{bmatrix},
  \begin{bmatrix}
    * & * & * \\
    0 & * & * \\
    0 & * & *
  \end{bmatrix},
  \textand
  \begin{bmatrix}
    * & * & * \\
    0 & * & * \\
    0 & 0 & *
  \end{bmatrix}
\]
that are nondegenerate is seen to be
\[
  \(1 - \frac{1}{q}\)^3 \cdot c,
\]
accounting for the three cases in the definition of $c$. Projectivizing, $T\backslash \Gamma$ is a subset of $(\PP^3 \OO_K)^3$ of volume
\[
  q^{\ds -\ceil{a_2 - a_1} - \ceil{a_3 - a_2} - \ceil{a_3 - a_1}} \cdot c.
\]

It remains to compute how the volume transforms under the bijection $\Psi : T\backslash \Gamma \isom V$ that we have constructed. This map is $K$-linear and is a product of the three maps
\begin{align*}
  \Psi_i : K \cross K \cross K &\to R \\
  \(c_{i1}, c_{i2}, c_{i3}\) &\mapsto \sum_{j=1}^3 c_{ij} \pi^{a_j - a_i} \frac{\gamma_j}{\gamma_i} \xi_j'.
\end{align*}
On the domain where it sends primitive vectors to primitive vectors, $\Psi_i$ scales volumes by $q^{-n_i}$, where $n_i$ is the determinant valuation, i.e.
\begin{equation} \label{eq:x_det_rho}
  \Psi_i\(\Lambda^3 \OO_K^3\) = \pi^{n_i} \Lambda^3 \OO_R.
\end{equation}
Extending scalars to $\OO_{\bar K}$, the left side of \eqref{eq:x_det_rho} becomes
\begin{align*}
  \Psi_i\(\Lambda^3 \OO_{\bar K}^3\) &= \prod_{j} \pi^{a_j - a_i} \bigwedge_j \frac{\gamma_j \xi_j'}{\gamma_i} \\
  &= \frac{\pi^{\sum_j a_j - 3 a_i}}{N(\gamma_i)} \bigwedge_j \xi_j
\end{align*}
When $R$ is tamely ramified, the wedge product of the $\xi_j$ generates the whole of $\Lambda^3 \OO_{\bar K}$, because the $\xi_j$ are an extender basis.
\begin{wild}
  When $R = K\cross Q$ is partially wildly ramified, the wedge product of the $\xi_j$ generate only $\pi^{(d_0 - 1)/2} \Lambda^3 \OO_{\bar K}$, as we saw in \eqref{eq:disc_defect}.
\end{wild}
Indeed, in all cases, if we let
\[
  d_0^{\mathrm{wild}} = d_0 - d_0^{\mathrm{tame}},
\]
then
\[
  \<\xi_1 \wedge \xi_2 \wedge \xi_3 \> = \pi^{d_0^{\mathrm{wild}}/2} \Lambda^3 \OO_{\bar K}^3.
\]
Accordingly, we get 
\begin{align*}
  \Psi_i\(\Lambda^3 \OO_{\bar K}\) &= \frac{\pi^{\sum_j a_j - 3 a_i + \frac{d_0^{\mathrm{wild}}}{2}}}{N(\gamma_i)} \Lambda^3 \OO_{\bar K}.
\end{align*}
Meanwhile, the right side of \eqref{eq:x_det_rho} is
\begin{align*}
  \pi^{n_i} \Lambda^3 \(\OO_R \tensor_{\OO_K} \OO_{\bar K}\) &= \pi^{n_i + \frac{d_0}{2}} \Lambda^3 \OO_{\bar K}.
\end{align*}
Hence
\[
  n_i = \sum_{j} a_j - 3a_i - v_K\(N(\gamma_i)\) - \frac{d_0^{\mathrm{tame}}}{2}
\]
so
\[
  \sum_i n_i = -v_K\(N(\gamma_1\gamma_2\gamma_3)\) - \frac{3d_0^{\mathrm{tame}}}{2}
\]
and
\[
  \mu(V) = \mu\(T\backslash\Gamma\) \cdot q^{-\sum_i n_i}
  = q^{\ds -\ceil{a_2 - a_1} - \ceil{a_3 - a_2} - \ceil{a_3 - a_1} + \frac{3d_0^{\mathrm{tame}}}{2} + v_K\(N(\gamma_1\gamma_2\gamma_3)\)} \cdot c,
\]
as desired.
\end{proof}

Consequently, we can compute the number of rings with any given discrete data by finding the volume of permissible $(\xi'_1, \xi'_2, \xi'_3)$, and dividing by $\mu(V)$. We carry out the computation of this volume in the succeeding sections. Observe that $\mu(V)$ equals $q^{2a_1 - 2a_3}$ times a correction that depends only on the $a_i$ mod $1$ and whether any $a_i$ are equal. This will simplify the entry of the ring volumes into Sage at the end of the proof.

\section{The conic over \texorpdfstring{$ \OO_K $}{OK}}

For each $\alpha \in R^\cross$, the equation
\[
  \tr(\alpha\xi^2) = 0
\]
defines a conic on the projectivization of $R$. Its determinant is $D_0 N(\alpha)$, up to squares of units, with respect to any $\OO_K$-basis of $\OO_R$, where $D_0$ is the discriminant of $\OO_R$. As we will find, it is preferable to transform the conic so that its discriminant has as low valuation as possible:

\begin{defn}
  Let $K$ be a local field, $\ch K \neq 2$. By a \emph{conic} over $\OO_K$ we mean a lattice $V$ of dimension $3$ over $\OO_K$ equipped with an integral bilinear form $\C$, or equivalently an integer-matrix quadratic form $\C : V \to \OO_K$, up to scaling by $\OO_K^\cross$. We say that $\C$ is
  \begin{itemize}
    \item \emph{unimodular} if $\det \C \sim 1$ (note that $\det \C$ is uniquely defined up to squares of units);
    \item \emph{tiny} if $\det \C \sim \pi$ and there exists a $v \in V$ such that $\C(v)$ is a unit;
    \item \emph{relevant} if it is either unimodular or tiny.
  \end{itemize}
\end{defn}
\begin{rem}
  Between changing basis and rescaling the whole form $\C$, we can scale $\det \C$ by any unit: hence we will sometimes assume that $\det \C$ is exactly $1$ or $\pi$.
\end{rem}

Let $\diamondsuit$ be a generator of the different ideal $\dd_{R/K}$. For instance, we can take
\begin{equation} \label{eq:diamond_choice}
\diamondsuit = \begin{cases}
  1 & R \text{ unramified} \\
  \pi_R^2 & R \text{ totally tamely ramified}
\begin{wild}
  \\
    (1; \bar{\zeta}_2\sqrt{D_0}), & R \isom K \cross Q \text{ partially wildly ramified.}
\end{wild}
\end{cases}
\end{equation}
Then by the definition of the different, the formula
\[
  \lambda^\diamondsuit(\xi) = \tr \frac{\xi}{\diamondsuit}
\]
defines a linear functional $\lambda^\diamondsuit : \OO_R \to \OO_K$ that is \emph{perfect,} that is, $\lambda^\diamondsuit$ generates the dual $\OO_R^\vee$ as an $\OO_R$-module, and hence the pairing
\[
  (x, y) \mapsto \lambda^\diamondsuit(xy)
\]
is a perfect $\OO_K$-linear pairing on $\OO_R$. If $\alpha \in \OO_R$, then the conic
\[
  \lambda^\diamondsuit(\alpha \xi^2)
\]
is $\OO_K$-integral on $\OO_R$ (because the corresponding bilinear form $\lambda^\diamondsuit(\alpha\xi\eta)$ is integral) with determinant $N(\alpha)$. We will put the conic defined by the $\M_{11}$- and $\M_{22}$-conditions in this form.

The entities involved in transformation will be marked by the symbol $\odot$ (``odot''). This is the symbol for a circle in Euclidean geometry, and it is chosen to reflect a particular simplifying fact: after the transformation, the conic is \emph{self-congruent}, that is, any two points on it can be taken to one another by an \emph{isometry} of $\PP^2(\OO_K)$ preserving the conic. This will follow from the independence of basepoint in Lemmas \ref{lem:conic_1} and \ref{lem:conic_pi}. The general conic over a $p$-adic field is not self-congruent.

\begin{lem}\label{lem:tfm_conic}
  Let $i \in \{1,2\}$. Fix the discrete data such a way that $\M_{ii}$ is active. Recall that $\delta = \delta_0\tau$ is in a fixed coarse coset. Then there is a multiplier $\gamma^{\odot}_i \in R^\cross$ with the following property:
  \[
    \xi_i^\odot = \frac{\xi_i'}{\gamma^{\odot}_i} = \frac{\xi_i}{\gamma_i \gamma^{\odot}_i \sqrt{\tau}}
  \]
  is a primitive vector in $\OO_R$, and the $\M_{11}$ condition $\tr(\xi_1^2) \equiv 0 \mod \pi^{m_{11}}$ is equivalent to
\begin{wild}
    either a linear condition on $\xi_i^\odot$ or
\end{wild}
a condition of the form
  \[
    \lambda^\diamondsuit\(\delta^\odot {\xi^\odot}^2\) \equiv 0 \mod \pi^{m_{11}^\odot},
  \]
  where $m_{11}^\odot = m_{11} - p^\odot$ is an integer and $\delta^\odot = \delta_0^{\odot}\tau \in \OO_R$, where $\delta_0^\odot$ depends on the discrete data alone and satisfies
  \[
    [\delta^\odot_0] = [\delta_0 \hat\omega_C \diamondsuit] \in H^1 \textand v_K(N(\delta^\odot_0)) \in \{0,1\},
  \]
  so that the conic $\M^\odot(\xi^\odot) = \lambda^\diamondsuit\(\delta^\odot {\xi^\odot}^2\)$ is relevant. Moreover, $\M^\odot$ is unimodular exactly when
  \[
    [\delta_0 \hat\omega_C \diamondsuit] \in \L_0.
  \]
\end{lem}

\begin{proof}
We have $\xi_1 \in \pi^{-a_1 - 2b_1} \sqrt{\delta \hat\omega_C} \cdot R$. 

Note that whatever $\gamma^{\odot}_i$ we pick, the conic takes the form
\[
  \lambda^\diamondsuit \(\diamondsuit {\gamma_i^\odot}^2 \gamma_i^2 \tau {\xi_i^\odot}^2\) \equiv 0 \mod \pi^{\bar m_{ii}},
\]
or, for any $p^\odot$,
\[
  \lambda^\diamondsuit\( \frac{\diamondsuit {\gamma_i^\odot}^2 \gamma_i^2}{\pi^{p^\odot}} \tau {\xi_i^\odot}^2 \) \equiv 0 \mod \pi^{\bar m_{ii} - p^\odot}.
\]
So we seek to pick $\gamma_i^{\odot}$ and $p^\odot$ so that
\[
  \delta^\odot_0 = \frac{\diamondsuit {\gamma_i^\odot}^2 \gamma_i^2 \delta_0}{\pi^{p^\odot}}
\]
lies in $\OO_R$ with norm of valuation $0$ or $1$, and
\[
  m_{ii}^\odot = \bar m_{ii} - p^\odot
\]
is an integer. The second condition is easily seen to follow from the first.

\paragraph{Unramified.}
If $[\delta_0 \hat\omega_C] \in \L_0$, then all extender indices are in $\ZZ$ and $\gamma_i$ is a unit, so $\delta^\odot_0$ is a unit as well, choosing $\gamma^{\odot}_i = 1$ and $p^\odot = 0$. Thus we get 

If $[\delta_0 \hat\omega_C] \in (1;\pi;\pi)\L_0$ for some ordering of the coordinates, then $[\delta] = [(1;\pi;\pi)\delta']$ for some unit $\delta'$. By Corollary \ref{cor:idxs_mod_1}, either 
\begin{itemize}
  \item $a_i' \in \ZZ$ and $\xi_i \in \sqrt{\delta'}(\OO_K^\cross \cross \sqrt{\pi} \OO_Q)$, or
  \item $a_i' \in \ZZ + \frac{1}{2}$ and $\xi_i \in \sqrt{\delta'}(\sqrt{\pi} \OO_K \cross \OO_Q)$.
\end{itemize}
In the first case, $\M_{11}$ is unsatisfiable if active, because $\alpha\xi_i^2$ has exactly one coordinate of zero valuation. So we have the second case. Observe that $\vec{v}(\gamma_i) = (1/2,0,0)$, so choosing $\gamma^{\odot}_i = 1$ and $p^\odot = 0$, we get $\vec{v}(\delta^\odot) = (1,0,0)$: the conic has determinant $\sim \pi$. Note that $\M^\odot \nequiv 0 \mod \pi$ as a quadratic form: after passing to an unramified extension we may assume that $L = K \cross K \cross K$, and then $\M^\odot$ is diagonal with two of the three coefficients units. So $\M^\odot$ is tiny.

For compatibility with the other splitting types, we let
\[
  h_i = \begin{cases}
    0, & [\delta_0 \hat\omega_C] \in \L_0 \\
    1, & \text{otherwise.}
  \end{cases}
\]

\paragraph{Splitting type $(1^3)$.} Here $\OO_R$ is generated by a uniformizer $\pi_R$ with $\pi_R^3 = \pi$ (for a suitably chosen uniformizer $\pi$). Under the Minkowski embedding, $\pi_R = \sqrt[3]{\pi} \cdot \bar{\zeta}_3$, where
\[
  \bar{\zeta}_3 = (1; \zeta_3; \zeta_3^2).
\]
We have $\diamondsuit = \pi_R^2$.

Here $[\delta_0 \hat\omega_C] \in \L_0$ always. The extender indices are in $\frac{1}{3}\ZZ$, and $\gamma_{i}$ is always a unit. Let $h = h_i \in \{0,1,-1\}$ be the integer such that
\[
  a_i' \in \ZZ - \frac{h}{3}.
\]
Then
\[
  \xi_1 \in \pi^{h/3} \sqrt{\delta^\odot} R = \bar{\zeta}_{3}^{-h} \sqrt{\delta^\odot} R;
\]
indeed, since $\xi_1$ is primitive in $\OO_{\bar{K}}^3$,
\[
  \xi_1 \in \bar{\zeta}_{3}^{-h} \sqrt{\delta^\odot} \OO_R^\cross.
\]

If $h = 0$, then $\M_{11}$ is unsatisfiable because the trace of a unit in $\OO_R$ is always a unit.

If $h = 1$, the choice $\gamma_i^{\odot} = 1$, $p^\odot = 2/3$ works, making $\delta^\odot$ a unit.

If $h = -1$, we can no longer take $\gamma_i^{\odot} = 1$, because the maximal possible value for $p^\odot$ is $1/3$ and the corresponding conic has determinant $\sim \pi^2$. Instead, take $\gamma_i^{\odot} = \pi_R^{-2}$. Then the corresponding values of $\xi^\odot$, instead of being units, have valuation $2/3$ and thus are still primitive in $\OO_R$. Take $p^\odot = -2/3$ and observe that $\delta^\odot$ is again a unit.

\begin{wild}
  \paragraph{Splitting type $(1^21)$.} Here things are more delicate because
\[
\diamondsuit = \bar{\zeta}_2\sqrt{D_0} \sim (1; \pi_R^{d_0})
\]
is significantly far from being a unit.

First we isolate the cases where $\M_{11}$ turns into a linear condition.
\begin{lem} \label{lem:conic_linear}
If
\[
  m_{11} \leq \frac{d_0 - 1}{2};
\]
slightly more generally, if (TODOWILD see if we need to be so precise)
\[
  a_1' \in \frac{1}{2} \ZZ \textand m_{11} \leq \floor{\frac{d_0}{2}},
\]
or if
\[
  a_1' \in \frac{1}{4} + \frac{1}{2} \ZZ \textand m_{11} \leq \ceil{\frac{d_0}{2}} - \frac{1}{2},
\]
then $\M_{11}$ is equivalent to the linear condition
\[
  \xi_1^{(K)} \equiv 0 \mod \pi^{m_{11}/2}.
\]
\end{lem}
\begin{proof}
The $\M_{11}$-condition is given as
\[
  \tr(\xi_1^2) \equiv 0 \mod \pi^{m_{11}}.
\]
We have
\[
  \tr(\xi_1^2) = \lambda^\diamondsuit\(\diamondsuit \xi^2\) = \(\diamondsuit \xi^2\)^{(K)} + \I\(\(\diamondsuit \xi^2\)^{(Q)}\).
\]
Our strategy will be to prove that the second term vanishes modulo $\pi^{m_{11}}$, so that $\M_{11}$ simplifies to
\begin{align*}
  \diamondsuit^{(K)} {\xi_1^{(K)}}^2 &\equiv 0 \mod \pi^{m_{11}} \\
  {\xi_1^{(K)}}^2 &\equiv 0 \mod \pi^{m_{11}} \\
  \xi_1^{(K)} &\equiv 0 \mod \pi^{m_{11}/2}.
\end{align*}
To this end, we analyze
\[
  \I\(\(\diamondsuit \xi^2\)^{(Q)}\) = \I\(\bar\zeta_2\sqrt{D_0} {\xi^{(Q)}}^2\)
\]
in the two cases:
\begin{itemize}
  \item If $a_1' \in \frac{1}{2} \ZZ$, then ${\xi^{(Q)}}^2 \in \OO_Q$, so the argument to $\I$ belongs to
  \[
    \bar\zeta_2 \sqrt{D_0} \OO_Q = \pi_Q^{d_0} \OO_Q = \OO_K \< \pi^{\ceil{d_0/2}}, \pi^{\floor{d_0/2}} \pi_Q \>.
  \]
  The $\I$-values of the last two generators are $0$ and $\pi^{\floor{d_0/2}}$ respectively, yielding the claimed congruence modulo $\pi^{\floor{d_0/2}}$.
  \item If $a_1' \in \frac{1}{4} + \frac{1}{2} \ZZ$, then ${\xi^{(Q)}}^2 \in (\pi_Q/\sqrt{\pi})\OO_Q$, so the argument to $\I$ belongs to
  \[
    \bar\zeta_2 \sqrt{D_0} \frac{\pi_Q}{\sqrt{\pi}} \OO_Q = \frac{\pi_Q^{d_0 + 1}}{\sqrt{\pi}} \OO_Q = \OO_K \< \pi^{\floor{d_0/2} + 1/2}, \pi^{\ceil{d_0/2} - 1/2} \pi_Q\>.
  \]
  The $\I$-values of the last two generators are $0$ and $\pi^{\ceil{d_0/2} - 1/2}$ respectively, yielding the claimed congruence modulo $\pi^{\ceil{d_0/2} - 1/2}$. \qedhere
\end{itemize}
\end{proof}

So assume that $m_{11} \geq d_0/2$. Already, the preceding lemma ensures that $v_K(\xi_1^{(K)}) \geq (d_0 - 1)/4$. In particular, $\rho_1$ is $Q$-led, so either
\begin{itemize}
  \item $[\delta\hat\omega_C] \in \L_0$ and $a_1' \in \frac{1}{2}\ZZ$, or
  \item $[\delta\hat\omega_C] \in (\pi ; \pi_Q) \L_0$ and $a_1' \in \frac{1}{4} + \frac{1}{2}\ZZ$.
\end{itemize}
Let $h_1 \in \{0,1,2,3\}$ be the integer such that
\[
  a_1' \in \ZZ + \frac{h_1 - d_0}{4}.
\]
Then in $H^1$,
\[
  [\delta\hat\omega_C] =  (\pi ; \pi_Q)^{h_1 - d_0} \cdot [\upsilon], \quad \upsilon \in \OO_R^\cross,
\]
and
\[
  \xi_1 \in \pi^{\frac{d_0 - h_1}{4}} \cdot (\pi ; \pi_Q)^{\frac{h_1 - d_0}{2}} \sqrt{\upsilon} R.
\]
We take $\gamma^{\odot}$ so that
\[
  \gamma\gamma^{\odot} = \pi^{\frac{d_0 - h_1}{4}} \cdot \(1 ; \pi_Q^{\frac{-d_0 + \1_{2\nmid h_1}}{2}}\) \sqrt{\upsilon}.
\]
It is evident that $\xi_1^\odot = \gamma^{-1}{\gamma^\odot}^{-1}\xi_1$ belongs to $R$; we need that it is a primitive vector of $\OO_R$. Since $\xi_1^{(Q)}$ is a unit,
\[
  v_Q({\xi_1^\odot}^{(Q)}) = -v_Q\((\gamma\gamma^{\odot})^Q\) = 2\(\frac{- d_0 + h_1}{4}\) + \frac{d_0 - \1_{2\nmid h_1}}{2}
  = \begin{cases}
    0, & h_1 = 0, 1 \\
    1, & h_1 = 2, 3
  \end{cases}
\]
while, since $v_K(\xi_1^{(K)}) \geq \frac{d_0-1}{4}$ as noted above,
\[
  v_K({\xi_1^\odot}^{(K)}) \geq \frac{d_0 - 1}{4} + \frac{-d_0 + h_1}{4} = \frac{h_1 - 1}{4} \geq -\frac{1}{4}.
\]
Since ${\xi_1^\odot}^{(K)} \in K$, its valuation is in fact nonnegative, and we have established that $\xi^\odot$ is primitive. Pick
\[
  p^\odot = \frac{d_0 - h_1}{2} \equiv -2a_1' \equiv m_{11} \mod \ZZ,
\]
thus getting
\[
  \delta^\odot = \frac{\diamondsuit {\gamma_i^\odot}^2 \gamma_i^2}{\pi^{p^\odot}} \sim \frac{(1; \pi^{d_0/2}) \cdot \pi^{\frac{d_0 - h_1}{2}}\(1; \pi^{\frac{-d_0 + \1_{2\nmid h_1}}{2}}\)}{\pi^{\frac{d_0 - h_1}{2}}} \sim (1; \pi_Q^{\1_{2\nmid h_1}}).
\]
So the conic has determinant $\sim 1$ or $\sim \pi$ according as $h_1$ is even or odd, that is, as $[\delta\hat\omega_C\diamondsuit]$ belongs to $\L_0$ or $(\pi; \pi_Q)\L_0$. Note that $\M \nequiv 0 \mod \pi$ as a quadratic form because $\M(1;0) \sim \lambda^\diamondsuit(1; 0) \sim 1$. So the conic is either unimodular or tiny in the two respective cases, as desired.

\end{wild}
\end{proof}
To summarize, the salient data of the transformation is shown here:
\begin{equation} \label{tab:tfm_conic}
\begin{tabular}{cccc|ccc|ccc}
 spl.t. & $a_1' \bmod \ZZ$ & $h_1$ & $[\delta\hat\omega_C\diamondsuit] \in$ & $v_K\(\gamma^{(K)} {\gamma^\odot}^{(K)}\)$ & $v_K\(\gamma^Q{\gamma^\odot}^Q\)$ & $p^\odot$ & $\delta^\odot \sim$ & $\xi' \sim$ \\ \hline
 ur & $0$ & $0$ & $\L_0$ & $0$ & $0$ & $0$ & $1$ & ? \\
 ur & $1/2$ & $1$ & $(1;\pi;\pi)\L_0$ & $1/2$ & $0$ & $0$ & $(\pi ; 1 ; 1)$ & $(? ; 1 ; 1)\vphantom{\dfrac{}{0}}$ \\
 $1^3$ & $-1/3$ & $1$ & $\L_0$ & $0$ & $0$ & $2/3$ & $1$ & $1$ \\
 $1^3$ & $1/3$ & $-1$ & $\L_0$ & $-2/3$ & $-2/3$ & $-2/3$ & $1$ & $\pi_R^2 \vphantom{\dfrac{}{0}}$
\begin{wild}
   \\
 $1^21$ & & $0$ & $\L_0$ & $d_0/4$ & $0$ & $d_0/2$ & $1$ & $(? ; 1)$ \\
 $1^21$ & & $1$ & $(\pi; \pi_Q) \L_0$ & $(d_0 - 1)/4$ & $0$ & $(d_0 - 1)/2$ & $(1; \pi_Q)$ & $(? ; 1)$ \\
 $1^21$ & & $2$ & $\L_0$ & $(d_0 - 2)/4$ & $-1/2$ & $(d_0 - 2)/2$ & $1$ & $(? ; \pi_Q)$ \\
 $1^21$ & & $3$ & $(\pi; \pi_Q) \L_0$ & $(d_0 - 3)/4$ & $-1/2$ & $(d_0 - 3)/2$ & $(1; \pi_Q)$ & $(? ; \pi_Q)$
\end{wild}
\end{tabular}
\end{equation}

The advantage of making the conic's determinant associate to either $1$ or $\pi$ is that we have to solve very few isomorphism types of conics. Although we do not prove the following classification, it animates the choice of what invariants we compute:
\begin{conj}\label{conj:conic_classfn}
  Let $\C$ be a relevant conic over the ring of integers $\OO_K$ of a local field $K$.
  \begin{enumerate}[$($a$)$]
    \item If $\C$ is tiny, it is determined up to isomorphism by its \emph{Brauer class} $\epsilon(\C) \in \{\pm 1\}$, the single bit telling whether $\C(\vec{x}) = 0$ has a nonzero solution over $K$.
    \item If $\C$ is unimodular, it is determined up to isomorphism by its Brauer class $\epsilon(\C)$ and its \emph{squareness level} $\ell(\C)$, the largest $\ell \in \ZZ, 0 \leq \ell \leq e/2$, such that
    \[
      \C \equiv c\lambda^2 \mod \pi^{\min\{2\ell + 1, e\}}
    \]
    as a quadratic form, for some constant $c \in \OO_K^\cross$ and linear form $\lambda$. Moreover, all combinations of values $ (\epsilon, \ell) $ occur, except that for $e$ even, $\ell(\C) = e/2$ implies $\epsilon(\C) = 1$ by Proposition \ref{prop:conic_1_mod_2} below. Thus there are exactly $e + 1$ isomorphism classes of conics of determinant $1$.
  \end{enumerate}
\end{conj}

Note that if $\C$ is unimodular and $2\ell + 1 \leq e$, then $\C$ is congruent to a $c\lambda^2$ modulo $\pi^{2\ell}$ if and only if modulo $\pi^{2\ell + 1}$, as the $x^2$, $y^2$, $z^2$ coefficients have square ratios modulo $\pi^{2\ell + 1}$ and the cross-terms are multiples of $2$ anyway. So the squareness level carries the same amount of information as the \emph{squareness}
\[
  \square(\C) = \max \left\{i : \C \equiv c\lambda^2 \mod \pi^i \text{ as a quadratic form}\right\} \in \left\{0,2,4, \ldots, 2\floor{\frac{e-1}{2}}\right\} \union \{e\}.
\]
(We cannot have $\square(\C) > e$, or the determinant would vanish modulo $\pi$.)

When we use coordinates, we will generally use one of two explicit types of conics: the \emph{diagonal conic}
\[
  aX^2 + bY^2 + cZ^2 = 0
\]
and the \emph{basepoint conic}
\[
  2XZ - Y^2 + aZ^2 = 0,
\]
so called because it passes through the basepoint $[1:0:0]$ and is tangent to the line $Z = 0$ there. We begin with results concerning the diagonal conic.

\subsection{Diagonal conics}

\begin{lem}\label{lem:conic_diag}  
Any relevant conic $\C$ is diagonalizable, that is, there exists a basis $(v_1,\ldots,v_3)$ for the given lattice $V$ such that
\[
  \C(x_1 v_1 + x_2 v_2 + x_3 v_3) = a_1 x_1^2 + a_2 x_2^2 + a_3 x_3^2.
\]
\end{lem}
\begin{proof}
When $\ch k_K \neq 2$, we have that \emph{any} conic is diagonalizable by an easy Gram-Schmidt procedure (in fact more is true: see O'Meara \cite{OMeara}, 92:1). So we assume that $\ch k_K = 2$. Here a quadratic space is not diagonalizable in general, and we must use the restrictions given on $\C$.

Write the matrix of $\C$, with respect to any basis $(\xi_1, \xi_2, \xi_3)$, as
\[
  \C = \begin{bmatrix}
    a & f & e \\
    f & b & d \\
    e & d & c
  \end{bmatrix}.
\]
In the case that $\C$ is unimodular, we see from
\[
  \det \C = abc + 2def - ad^2 - be^2 - cf^2 \sim 1
\]
that at least one of the diagonal entries---say $a$---is a unit. Then we can use $a$ to eliminate $f$ and $e$ (that, is, add multiples of $\xi_1$ to $\xi_2$ and $\xi_3$). Now if $b$ (or, symmetrically, $c$), is nonzero modulo $\pi$, we use it to eliminate $d$, and we are done, as we have found the requisite diagonal form. However, it is possible that
\[
  \C \equiv \begin{bmatrix}
    a &   & \\
      &   & d \\
      & d &
  \end{bmatrix} \mod \pi
\]
for units $a$ and $d$. Rescaling $\xi_3$, we can assume that
\[
  \C \equiv a \cdot \begin{bmatrix}
    1 &   & \\
      &   & 1 \\
      & 1 &
  \end{bmatrix} \mod \pi.
\]
At first we are doubtful, because the unimodular form
\[
\C_2 = \begin{bmatrix}
  0 & 1 \\
  1 & 0
\end{bmatrix}
\]
on $\OO_K^2$ is not diagonalizable. However, we can use the identity
\[
  \begin{bmatrix}
    1 & 1 & 1 \\
    1 & 1 &   \\
    1 &   & 1
  \end{bmatrix}
  \begin{bmatrix}
    1 &   &   \\
      &   & 1 \\
      & 1
  \end{bmatrix}
  \begin{bmatrix}
    1 & 1 & 1 \\
    1 & 1 &   \\
    1 &   & 1
  \end{bmatrix}
  \equiv
  \begin{bmatrix}
    1 &   &   \\
      & 1 &   \\
      &   & 1
  \end{bmatrix}
  \mod 2
\]
to change to a basis in which
\[
  \C \equiv a \cdot
  \begin{bmatrix}
    1 &   &   \\
      & 1 &   \\
      &   & 1
  \end{bmatrix} \mod \pi.
\]
Then the diagonalization proceeds without a hitch.

If $\C$ is tiny, we proceed similarly. Taking $a$ a unit (since we are given $\C(\xi_1) \sim 1$ for some $\xi_1$), we can eliminate $f$ and $e$. Then since
\[
  \det \C = a(bc - d^2) \sim \pi,
\]
we must have at least one of $b$ and $c$ a unit, as otherwise $\det \C$ would be either a unit (if $d$ is a unit) or a multiple of $\pi^2$ (if $\pi \mid d$). So we can eliminate $d$ and again get the desired diagonalization.
\end{proof}

The following lifting lemma for solutions modulo $4\pi$ will be essential for us.
\begin{lem}\label{lem:conic_lift}
Let $\C$ be a diagonalized conic on an $\OO_K$-lattice $V$, and assume that $v(\det \C) \leq 1$. Let $v \in V$ be a primitive vector with
\[
  \C(v) \equiv 0 \mod \pi^m, \quad m > 2e.
\]
Then there exists a $v' \in V$ such that
\[
  v' \equiv v \mod \pi^{m-e} \textand \C(v') = 0.
\]
\end{lem}
\begin{proof}
We may write the conic in diagonal form
\[
  \C(x_1v_1 + \cdots + x_nv_n) = a_1 x_1^2 + \cdots + a_n x_n^2.
\]
Let $v = x_1v_1 + \cdots + x_nv_n$. Since $v$ is primitive, not all the $x_i$ are zero modulo $\pi$. We claim that there is an $i$ with
\begin{equation}\label{eq:x_move_coord}
  \pi \nmid a_i \textand \pi \nmid x_i.
\end{equation}
If not, then $\det \C \sim \pi$, and without loss of generality, $a_1,\ldots,a_{n-1}$ are units while $a_n \sim \pi$; and $x_n$ is a unit while $x_1,\ldots,x_{n-1}$ are multiples of $\pi$. Summing, we find that
\[
  \C(v) \equiv a_n x_n^2 \nequiv 0 \mod \pi^2,
\]
a contradiction.

Choose $i$ satisfying \eqref{eq:x_move_coord}. We will construct $v'$ by changing only the $x_i$ coordinate of $v$ to a different value $x_i'$. The desired condition $\C(v') = 0$ takes the form
\[
  x_i'^2 = y
\]
for some $y \equiv x_i^2$ modulo $\pi^{m}$. Since $m > 2e$, we have that $y$ is also a square and, indeed, has a (unique) square root $x_i'$ satisfying $x_i' \equiv x_i$ mod $\pi^{m-e}$. This constructs the desired $v'$.
\end{proof}

Here are two easy corollaries.
\begin{prop}\label{prop:conic_1_mod_2}
If $e$ is even and $\C$ is a unimodular conic of squareness $e$ (the maximal possible value), then $\epsilon(\C) = 1$, that is, $\C$ has a rational point.
\end{prop}
\begin{proof}
We may assume that $\C$ is diagonal:
\[
  \C(x_1 \xi_1 + x_2 \xi_2 + x_3 \xi_3) = a_1 x_1^2 + a_2 x_2^2 + a_3 x_3^2.
\]
Then $\epsilon(\C)$ is a Hilbert symbol,
\[
  \epsilon(\C) = \< \frac{-a_2}{a_1}, \frac{-a_3}{a_1} \>.
\]
Both arguments are squares of units modulo $2$. But since $e$ is even, they are actually squares modulo $2\pi$, so the Hilbert symbol is $1$ by Lemma \ref{lem:Hilb_prod_size}.
\end{proof}

\begin{prop}\label{prop:conic_perturb}
Let $\C$ and $\C'$ be conics of determinant $1$ and squareness level $\ell$. Suppose that the associated bilinear forms of $\C$ and $\C'$ are congruent modulo $\pi^{2e - 2\ell}$. Then $\epsilon(C) = \epsilon(C')$.
\end{prop}
\begin{proof}
We may assume that $\C$ is diagonal:
\[
\C(x_1 \xi_1 + x_2 \xi_2 + x_3 \xi_3) = a_1 x_1^2 + a_2 x_2^2 + a_3 x_3^2.
\]
Although $\C'$ need not be diagonal with respect to the same basis, the orthogonalization procedure furnished by the proof of Lemma \ref{lem:conic_diag} yields a basis $(\xi_1', \xi_2', \xi_3')$ with $\xi_j' \equiv \xi_j \mod \pi^{2e - 2\ell}$ such that $\C'$ is diagonal with respect to it,
\[
\C'(x_1 \xi'_1 + x_2 \xi'_2 + x_3 \xi'_3) = a'_1 x_1^2 + a'_2 x_2^2 + a'_3 x_3^2
\]
with $a_j \equiv a'_j \mod \pi^{2e - 2\ell}$. New compare
\[
  \epsilon(\C) = \< \frac{-a_2}{a_1}, \frac{-a_3}{a_1} \> \textand \epsilon(\C') = \< \frac{-a_2'}{a_1'}, \frac{-a_3'}{a_1'} \>
\]
The arguments to the Hilbert symbols are squares modulo $\pi^{2\ell + 1}$, so the value of the Hilbert symbol is unchanged under multiplying by units that are squares modulo $\pi^{\max\{2\ell+1, 2e - 2\ell\}}$ by Lemma \ref{lem:Hilb_prod_size}, and that is exactly what we have done.
\end{proof}

\subsection{The solution volume of the conic}

We now use the basepoint form to determine volumes of conics.

\begin{lem}[\textbf{Igusa zeta function of a conic}]\label{lem:conic_1}
  Let $\C$ be a conic of determinant $1$ on the projectivization $\PP(V)$ of a $3$-dimensional vector space $V$. Suppose that $\C$ has Brauer class $\epsilon(\C) = 1$, that is, it admits a basepoint $v_0$ such that $\C(v_0) = 0$.
  
  Let $U_{m^\odot, n^\odot}$ be the volume of $v \in \PP(V)$ (counting the whole $\PP(V)$ to have volume $1 + q^{-1} + q^{-2}$) such that
  \begin{align}
    v &\equiv v_0 \mod \pi^{n^\odot} \label{eq:conic_1_n} \\
    \C(v) &\equiv 0 \mod \pi^{m^\odot}. \label{eq:conic_1_m}
  \end{align}
  (Note that these $m^\odot$ and $n^\odot$ correspond to the $m_{11}^\odot$ and $n^\odot$ of Lemmas \ref{lem:tfm_conic} and \ref{lem:N11}.)
  
  Then for $m^\odot$ and $n^\odot$ integers with $m^\odot > 2e$ and $m^\odot \geq 2n^\odot$, the volume $U_{m^\odot,n^\odot}$ depends only on $m^\odot$, $n^\odot$, and the squareness level $\ell = \ell(\C)$. It is given by
  \[
  U_{m^\odot,n^\odot} = q^{2e-m^\odot-n^\odot}, \quad n^\odot \geq e
  \]
  and the recurrence
  \begin{align*}
    \frac{U_{m^\odot,n^\odot}}{U_{m^\odot,n^\odot+1}} &= \begin{cases}
      2 & n^\odot = e - 2\ell - 1 \geq 0 \\
      q & n^\odot \equiv e \mod 2, n^\odot > e - 2\ell - 1, n^\odot > 0 \\
      q+1 & n^\odot = 0, \ell = \frac{e}{2} \\
      1 & otherwise.
    \end{cases}
  \end{align*}
  Explicitly,
  \[
  {U_{m^\odot,n^\odot} = }
  \left\{\begin{tabular}{lll}
    $ q^{e-m^\odot-n^\odot} $ & $ n^\odot \geq e $ & (a \emph{black conic}) \\
    $ q^{-m^\odot + \floor{\frac{e - n^\odot}{2}}} $ & $ e \geq n^\odot > e - 2\ell - 1, n^\odot > 0 $ & (a \emph{blue conic}) \\
    $ 2 q^{-m^\odot + \ell} $ & $  0 \leq n^\odot \leq e - 2\ell - 1 $ & (a \emph{green conic}) \\
    $ \( 1 + \dfrac{1}{q}\) q^{-m^\odot + e/2}, $ & $ n^\odot = 0, \ell = \dfrac{e}{2}, \text{$ e $ even}  $ & (a \emph{large conic}) 
  \end{tabular}
  \right.
  \]
\end{lem}

\begin{proof}
  For the black-conic case, we diagonalize the conic to
  \[
  \C(X,Y,Z) = aX^2 + bY^2 + cZ^2.
  \]
  Let the basepoint be $v_0 = [X_0 : Y_0 : Z_0]$. Note that two coordinates of $v_0$, say $X_0$ and $Y_0$, are nonzero modulo $\pi$. We may scale so that $Y_0 = 1$ and so that all solutions we seek have $Y = 1$. This eliminates the issue of scaling ambiguity.
  
  The $n^\odot$-pixel of $[X:1:Z]$ satisfying \eqref{eq:conic_1_n} has volume $q^{-2n^\odot}$. For fixed $Z$, with $Z \equiv Z_0 \mod \pi^{n^\odot}$ the condition \eqref{eq:conic_1_m} simplifies to $X^2 \equiv u \mod \pi^{m^\odot}$, where $u$ is a unit with $u \equiv X_0^2 \mod \pi^{n^\odot + e}$. Hence its solutions form a congruence class mod $\pi^{m^\odot-e}$, and overall, the solution volume is $q^{n^\odot} \cdot q^{e-m^\odot} = q^{e-m^\odot-n^\odot}.$
  
  We use this as the base case to prove the recursive formula (and hence also the explicit formula) for $U_{m^\odot,n^\odot}$ by downward induction on $n^\odot$. Our aim is to determine the number $r$ of $(n^\odot+1)$-pixels within the $n^\odot$-pixel of $v_0$ that contain a solution to $\C(v) = 0$, or equivalently, to $\C(v) \equiv 0 \mod \pi^{m^\odot}$. Then, by induction, there is a volume $U_{m^\odot,n^\odot+1}$ of solutions in each of those, so $U_{m^\odot,n^\odot} = rU_{m^\odot,n^\odot+1}$ as desired. It remains to compute
  \[
    r = \frac{U_{m^\odot,n^\odot}}{U_{m^\odot,n^\odot+1}}.
  \]
  
  We now abandon the diagonalized form and choose coordinates such that the basepoint is $v_0 = [1:0:0]$ and the tangent line there is $Z = 0$. Then the conic has the form
  \[
  \C(X,Y,Z) = 2gXY + 2dXZ - cY^2 + 2bYZ + aZ^2.
  \]
  Note that $\pi$ does not divide both $d$ and $g$, for then the conic's determinant
  \[
  \begin{vmatrix}
    0 & g & d \\
    g & -c & b \\
    d & b & a
  \end{vmatrix} = 2bdg + cd^2 - ag^2
  \]
  would be divisible by $\pi^2$. So, by symmetry, we may assume $d \sim 1$. We scale the conic so that $d = 1$, and then the transformation $Z \mapsto Z - gY$ makes $g = 0$. Also, the transformation $X \mapsto X - bY$ makes $b = 0$. Now $c = 1$ to make the determinant $1$. Thus the conic takes the basepoint form
  \[
  \C(X,Y,Z) = 2XZ - Y^2 + aZ^2,
  \]
  where $a \in \OO_K$ is the only undetermined coefficient. By definition of squareness level, we know that $a$ is a square modulo $\pi^{\min\{2\ell + 1, e\}}$, but not modulo $\pi^{2\ell + 3}$ if $\ell < \floor{e/2}$. For any $a' \in \OO_K$, the transformation $Y \mapsto Y + a'Z$ can be used to increment $a$ by the square $a'^2$, followed by another $X \mapsto X - b'Y$ to remove the $YZ$ term. Picking $a'$ appropriately, we can arrange so that $v(a)$ \emph{reveals the squareness:} either
  \begin{itemize}
    \item $\ell < \floor{e/2}$ and $v(a) = 2\ell + 1$, or
    \item $\ell = \floor{e/2}$ and $2|a$. In this case, indeed, the transformation $X \mapsto X - a/2$ makes $a = 0$. (So we have proved one case of Conjecture \ref{conj:conic_classfn}: for $\ell = \floor{e/2}$ and $\epsilon = 1$, the conic takes the fixed form $2XZ - Y^2$.)
  \end{itemize}
  
  To parametrize $ \C(\OO_K) $, we use the age-old trick of \emph{stereographic projection,} that is, drawing lines of varying slope through the known basepoint $[1:0:0]$. An easy calculation shows that the second intersection of the line $sZ - tY = 0$ with the conic $\C$ is $[s^2 - at^2 : 2st : 2t^2]$, yielding an isomorphism
  \begin{align*}
    \PP^1(K) &\isom \C(K) \\
    [s:t] &\mapsto [s^2 - at^2 : 2st : 2t^2].
  \end{align*}
  If $ [s:t] $ is in lowest terms over $ \OO_K $, then $ [s^2 - at^2 : 2st : 2t^2] $ need not be in lowest terms over $ \OO_K $, but will have cancellation by $\pi^j$, where $ j = \min\{v_K(s^2 - at^2), e + v_K(st), e + v_K(t^2)\} $. Note that $ j \leq e $ because $ s $ and $ t $ are coprime, so
  \[
  j = \min\{v_K(s^2 - at^2), e\}.
  \] Note also that the resulting point
  \[
  [X : Y : Z] = \left[ \frac{s^2 - at^2}{\pi^j} : \frac{2st}{\pi^j} : \frac{2t^2}{\pi^j} \right]
  \]
  lies in the $ (e - j) $-pixel of $ [1:0:0] $, but not in the $ (e - j + 1) $-pixel if $ j > 0 $. Hence the points we are interested in, namely in the $ n^\odot $-pixel of the basepoint but outside the $ (n^\odot + 1) $-pixel, correspond exactly to values of $ [s : t] $ for which $ j = e - n^\odot $. That is, the valuation $ v(s^2 - at^2) $ must be exactly $ e-n^\odot $ (if $ n^\odot > 0 $) or at least $ e $ (if $ j = 0 $).
  
  Observe that when $n^\odot < e - 2\ell - 1$, there are no solutions. Also, when $n^\odot > e - 2\ell - 1$ is of the same parity as $e$, there are no solutions, because 
  \[
  s^2 - at^2 \equiv s^2 \mod \pi^{\min\{2\ell + 1, e\}}
  \]
  has even valuation if nonzero mod $\pi^{\min\{2\ell + 1, e\}}$. The ratio $ {U_{m^\odot,n^\odot}} / {U_{m^\odot, n^\odot+1}} $ is thus $ 1 $ in these cases, as claimed.
  
  \begin{itemize}
    \item Suppose that $ n^\odot = e - 2\ell - 1 $. If $ n^\odot > 0 $, we seek the $[s:t]$ such that $s^2 - a t^2$ attains its maximal valuation $2\ell+1$: this happens when $\pi^{\ell+1} | s$. If $n^\odot = 0$, we seek $\pi^{e}|s^2$, which is still equivalent to $\pi^{\ell+1} | s$. Hence we are looking at the $[s:t]$ with $ t = 1 $, $ s = \pi^{\ell + 1} s' $, where $ s' \in \OO_K$. Then the resulting point on $ \C(\OO_K) $ is
    \begin{align*}
      [X : Y : Z] &= \left[ \frac{s^2 - at^2}{\pi^{2\ell+1}} : \frac{2st}{\pi^{2\ell+1}} : \frac{2t^2}{\pi^{2\ell+1}} \right] \\
      &= \left[ \frac{\pi^{2\ell+2}s'^2 - a}{\pi^{2\ell+1}} : \frac{2\pi^{\ell+1}s'}{\pi^{2\ell+1}} : \frac{2t^2}{\pi^{2\ell+1}} \right] \\
      &= \left[ \pi s'^2 - \frac{a}{\pi^{2\ell + 1}} : \frac{2s'}{\pi^{\ell}} : \frac{2}{\pi^{2\ell + 1}} \right] \\
      &\equiv \left[ \pi s'^2 - \frac{a}{\pi^{2\ell + 1}} : 0 : \frac{2}{\pi^{2\ell + 1}} \right] \mod \pi^{n^\odot+1} = \pi^{e-2\ell}.
    \end{align*}
    We claim that this is actually the same point modulo $ \pi^{n^\odot+1} $ regardless of $ s' $, that is, the $ \pi s'^2 $ term contributes nothing. If $ 2\ell + 1 = e $, this is clear because $ n^\odot = 0 $. Otherwise, $ a' = a/\pi^{2\ell + 1} $ is a unit, and if we multiply all three coordinates by 
    \[
    \frac{a'}{a' + \pi s'} \equiv 1 \mod \pi,
    \]
    the last two coordinates do not change mod $ \pi^{n^\odot+1} $ because they are $ 0 $ mod $ \pi^{n^\odot} $. Thus all points $ [X : Y : Z] $ obtained lie in a single $ (n^\odot + 1) $-pixel, and hence the ratio $ {U_{m^\odot,n^\odot}} / {U_{m^\odot, n^\odot+1}} $ is $ 2 $.
    \item Suppose that $ n^\odot \equiv e \mod 2$, $n^\odot > e - 2\ell - 1 $, and $ n^\odot > 0 $. Then $ e - n^\odot = j = 2j' $ is even, with $ j' \leq \ell. $ The pairs $ [s:t] $ yielding $ s^2 - at^2 \sim \pi^j $ are exactly those with $ s \sim \pi^{j'} $. Write $ t = 1 $, $ s = \pi^{j'} s' $, where $ s' \in \OO_K^\cross $. Then the resulting point on $ \C(\OO_K) $ is
    \begin{align*}
      [X : Y : Z] &= \left[ \frac{s^2 - at^2}{\pi^j} : \frac{2st}{\pi^j} : \frac{2t^2}{\pi^j} \right] \\
      &= \left[ - \frac{a}{\pi^{2j'}} + s'^2 : \frac{2s'}{\pi^{j'}} : \frac{2}{\pi^{2j'}} \right] \\
      &\equiv \left[ s'^2 : 0 : \frac{2}{\pi^{2j'}} \right] \quad \mod \pi^{n^\odot + 1},
    \end{align*}
    where at the last step we multiplied all three coordinates by the unit
    \[
    \frac{s'^2}{s'^2 - \frac{a}{\pi^{2j'}}} \equiv 1 \mod \pi.
    \]
    We get $ q - 1 $ different $ (n^\odot + 1) $-pixels, one for each value of $ s' $ mod $ \pi $. Hence the ratio $ {U_{m^\odot,n^\odot}} / {U_{m^\odot, n^\odot+1}} $ is $ q $.
    \item Finally, suppose that $ n^\odot = 0 $ and $ \ell = \frac{e}{2} $. Write $ t = 1 $, $ s = \pi^{e/2} s' $, where $ s' \in \OO_K $. We get
    \begin{align*}
      [X : Y : Z] &= \left[ \frac{s^2}{\pi^e} : \frac{2st}{\pi^e} : \frac{2t^2}{\pi^e} \right] \\
      &= \left[ \frac{2s'}{\pi^{e/2}} + s'^2 : \frac{2s'}{\pi^{e/2}} : \frac{2}{\pi^{e}} \right] \\
      &\equiv \left[ s'^2 : 0 : \frac{2}{\pi^{e}} \right] \quad \mod \pi = \pi^{n^\odot+1}.
    \end{align*}
    We get $ q $ different $ 1 $-pixels, one for each value of $ s' $ mod $ \pi $. Hence the ratio $ {U_{m^\odot,n^\odot}} / {U_{m^\odot, n^\odot+1}} $ is $ q + 1 $. \qedhere
  \end{itemize}
  
\end{proof}

For determinant $\pi$, we use the same method. Fortunately, everything comes out much simpler.
\begin{lem}[\textbf{Igusa zeta function of a conic}]\label{lem:conic_pi}
  Let $\C$ be an integer-matrix conic over $\OO_K$ of determinant $\pi$. Suppose that $\C$ has Brauer class $\epsilon(\C) = 1$, that is, it admits a basepoint $v_0$ such that $\C(v_0) = 0$.
  
  Let $U_{m^\odot,n^\odot}$ be the volume of $v \in \PP(V)$ such that
  \begin{align}
    v &\equiv v_0 \mod \pi^{n^\odot} \label{eq:conic_pi_n} \\
    \C(v) &\equiv 0 \mod \pi^{m^\odot}. \label{eq:conic_pi_m}
  \end{align}
  
  Then for $m^\odot > 2e$ and $m^\odot \geq 2n^\odot$, the volume $U_{m^\odot,n^\odot}$ depends only on $m^\odot$ and $n^\odot$:
  \[
  {U_{m^\odot,n^\odot} = }
  \left\{\begin{tabular}{lll}
    $ q^{e-m^\odot-n^\odot} $ & $ n^\odot > e $ & (a \emph{black conic}) \\
    $ 2 q^{-m^\odot} $ & $ 0 \leq n^\odot \leq e $ & (a \emph{green conic})
  \end{tabular}
  \right.
  \]
\end{lem}
\begin{proof}
  Diagonalize the conic to the form
  \[
  \C(X,Y,Z) = a X^2 + b Y^2 + c \pi Z^2 = 0,
  \]
  where $a,b,c \in \OO_K^\cross$. Observe that $X_0$ and $Y_0$ are units, and scale so that $Y = Y_0 = 1$.
  
  If $n^\odot > e$, then for each $Z \equiv Z_0 \mod \pi^{n^\odot}$, the condition $\C(X,Y,Z) \equiv 0 \mod \pi^{m^\odot}$ simplifies to $X^2 \equiv u \mod \pi^{m^\odot}$, where $u \equiv X_0^2$ mod $\pi^{n^\odot+e+1}$, and hence the square roots $X$ with $X \equiv X_0$ mod $\pi^{n^\odot}$ form a single congruence class mod $\pi^{m^\odot-e}$. So the volume is $q^{e-m^\odot-n^\odot}$.
  
  If $n^\odot = e$, then we use the same method, but now the equation $X^2 \equiv u \mod \pi^{m^\odot}$, where $u \equiv X_0^2$ mod $\pi^{2e+1}$, has as solution set two classes mod $\pi^{m^\odot-e}$, each the negative of the other.
  
  We claim that these are all the solutions mod $\pi^{m^\odot}$; that is, that the whole conic lies within an $e$-pixel. Suppose there is such an $[X : 1 : Z] \nequiv [X_0 : 1 : Z_0] \mod 2$, and let $v_K(X - X_0) = i$, $v_K(Z - Z_0) = j$. Then
  \begin{align*}
    a (X^2 - X_0^2) = \pi c (Z^2 - Z_0^2).
  \end{align*}
  If $i < e$, then the left side has even valuation $2i$ which cannot be matched by the right side. If $j < e$, then the right side has odd valuation $2j+1$ which cannot be matched by the left side. This completes the proof.
  
\end{proof}

\subsection{The Brauer class}

\begin{normalsize}
In this section we understand the Brauer class of conics of the form
\[
\A : \quad   \tr (\alpha x^2) = 0.
\]
This is a conic in the $ \PP^2(K) $ of possible values of $ x $. Over $ K $, there are just two types of conic, one with points and one without. Our first task will be to understand which case occurs for each $ \alpha $.

Our main result will be the following.

\begin{lem} \label{lem:H_form}
  If $ \alpha \in R^{N=1} $, define
  \[
  \epsilon(\alpha) = \begin{cases}
  1 & \text{if $ \tr(\alpha x^2) = 0 $ for some nonzero $ x \in L $} \\
  -1 & \text{otherwise.}
  \end{cases}
  \]
  Then the map of $ \FF_2 $-vector spaces
  \begin{align*}
  H^1 &\to \mu_2 \\\
  \alpha &\mapsto \epsilon(\alpha) / \epsilon(1)
  \end{align*}
  is a nondegenerate quadratic form whose associated bilinear form is none other than the Hilbert pairing on $ H^1 $. That is,
  \[
  \epsilon(\alpha\beta) = \epsilon(1) \cdot \epsilon(\alpha) \cdot \epsilon(\beta) \cdot \<\alpha,\beta\>_{\epsilon}.
  \]
\end{lem}
The proof is not especially difficult, but it uses different tools than the rest of the paper and so will be deferred. See Appendix \ref{sec:GW}.

\begin{rem}
  $\epsilon$ comes up, in a related context, in the work of Bhargava and Gross (\cite{AIT}, \textsection 7.2), where it is stated to be a quadratic form, at least in the tamely ramified case.
\end{rem}

\begin{old}
For the reader's interest, we also include the following fact.
\begin{lem}\label{lem:H_extra}
    $ \epsilon(\omega) = 1 $ if and only if the associated quartic algebra $ L_\omega $ has an element $ \xi \neq 0 $ whose characteristic polynomial has the doubly depressed form
    \[
    \xi^4 - a\xi - b = 0.
    \]
    (Note that $ L_\omega = L_{\omega_0\delta} $ is in general distinct from the algebra $ L_\delta $ in which we have been trying to find orders.)
\end{lem}
\begin{proof}
  Note that $ \xi \in L_\omega $ has a doubly depressed characteristic polynomial if and only if
  \[
  \tr \xi = \tr \xi^2 = 0.
  \]
  Now the traceless elements $ \xi \in L_\omega $ have the form
  \[
  \kappa(x) = \tr_{RL_\omega/L_\omega} (x \sqrt{\omega}).
  \]
  For such $ \xi $, we compute
  \[
  \tr(\xi)^2 = \tr_{R/K} (\omega x^2),
  \]
  as desired.
\end{proof}
\end{old}

\begin{rem}
  Over fields of characteristic not $ 2 $, a quadratic form is uniquely determined by its associated bilinear form. However, over $ \FF_2 $, the local Hilbert pairing $ \<\bullet, \bullet\>_\epsilon $ on $ H^1 $ lifts to $ \size{(H^1)^*} = \size{H^1} $ quadratic forms, thanks to the ambiguity by adding a linear functional. It is not hard to show that these quadratic forms are exactly
  \[
  \alpha \mapsto \frac{\epsilon(\omega\alpha)}{\epsilon(\omega)},
  \]
  for each $ \omega \in H^1 $.
\end{rem}

\end{normalsize}
\subsection{The squareness}

\begin{normalsize}
For the cases in Lemma \ref{lem:tfm_conic} in which the transformed conic $\M$ is unimodular, we need also to compute its squareness.

If $\C$ is a unimodular conic, the maximal 

\begin{lem}\label{lem:squareness}
Let $[\heartsuit] = [1] \in H^1$. (The reason for this strange definition is that, in splitting type $1^2 1$, we will need a $[\heartsuit] \neq [1]$ in general.)

\begin{wild}
  Let $[\heartsuit] \in H^1$ be defined as follows:
\begin{itemize}
  \item If $R = K \cross Q$ is partially ramified, let $\pi_Q$ be any uniformizer that does \emph{not} have trace zero, and let
  \[
    [\heartsuit] = [(1; \tr \pi_Q)].
  \]
  \item In all other cases, we may take
  \[
    [\heartsuit] = 1.
  \]
\end{itemize}

Then
\end{wild}
The conic 
\[
  \M(\xi^\odot) = \lambda^\diamondsuit \(\delta^\odot {\xi^\odot}^2\), \quad \delta^\odot \in \OO_R^\cross
\]
has squareness
\begin{align} \label{eq:sqness}
  \square(\M) &= \max \left\{ \ell : [\delta^\odot\heartsuit] \equiv 1 \mod \pi^{\ell}; \quad \ell = e, \text{ or } \ell < e \text{ and } \ell \text{ is odd}\right\} \\
  &= \begin{cases}
    \min\left\{2\ell(\delta\hat\omega_C\diamondsuit\heartsuit) + 1, e\right\}, & \text{$R$ unramified} \\
    \min\left\{2\floor{\dfrac{\ell(\delta\hat\omega_C\diamondsuit\heartsuit)}{2} }+ 1, e\right\}, & \text{$R$ ramified}.
  \end{cases}
\end{align}
\end{lem}
\begin{rem}
  As stated, the lemma only requires $\heartsuit$ to be defined modulo $\L_{\floor{e/2}}$ (unramified types) resp{.} $\L_{2\floor{e/2}}$ (ramified types). We will mostly use $\heartsuit$ in this way, but when we do the brown zone, we will need a finer definition and will mention this.
\end{rem}
\begin{proof}
In unramified splitting type, we first claim that going up to an \emph{unramified} extension $K'/K$ does not change either the left or the right side of \eqref{eq:sqness}. The right-hand side is less than $\floor{e/2}$ only if $[\delta \hat{\omega}_C \heartsuit]$ is represented by a generic unit $\delta^\odot = 1 + \alpha \pi^{2i+1}$, $2i + 1 < e$, and this generic unit remains generic in $R' = K' \tensor_K R$. As to the left side, we can diagonalize the conic $\M$ to have the form
\[
  \M(X,Y,Z) = aX^2 + bY^2 + cZ^2, \quad abc = 1.
\]
Then
\[
  \square(\M) = \min\left\{2\ell(b/a) + 1, 2\ell(c/b) + 1, e\right\},
\]
and this remains invariant over $K'$.

Therefore, we may assume that $R \isom K \cross K \cross K$ is totally split. Let $\delta^\odot = (a; b; c)$. Then $\M$ is diagonal and
\begin{align*}
  \square(\M) &= \min\left\{2\ell(b/a) + 1, 2\ell(c/b) + 1, e\right\} \\
  &= \min\left\{2\ell(a) + 1, 2\ell(b) + 1, 2\ell(c) + 1, e\right\} \\
  &= \min\left\{2\ell(\delta^\odot) + 1, e\right\},
\end{align*}
as desired.

In splitting type $1^3$, we can scale $\delta^\odot$ by $(\OO_R^\cross)^2$ so that its level is manifest:
\[
  \delta^\odot = 1 + \alpha \cdot \pi_R^{2j+1},
\]
where $\alpha \in \OO_R^\cross$ and where
\[
  j \in \{1, 3, 4, 6, \ldots, 3e - 2, 3e, \infty \}
\]
controls $\ell(\delta^\odot)$. Let $\delta^\odot = a + b\pi_R + c\pi_R^2$ and note that, in the basis $(1, \pi_R, \pi_R^2)$, the conic
\[
  \M(\xi') = \tr(\pi_R^{-2} \delta^\odot \xi'^2)
\]
has matrix
\[
  \begin{bmatrix}
    c & b & a \\
    b & a & \pi c \\
    a & \pi c & \pi b
  \end{bmatrix}.
\]
Since $a$ is a unit, we get
\begin{align*}
  \ell(\M) &= \max\left\{ i \leq e : \frac{c}{a} \text{ and } \frac{\pi b}{a} \text{ are squares mod } \pi^{i} \right\} \\
  &= \min\left\{2\floor{\frac{j - 1}{3}}, e\right\} \\
  &= \min\left\{\floor{\frac{\ell(\delta^\odot)}{2}},e\right\},
\end{align*}
as desired.

\begin{wild}
  
In splitting type $1^21$, we first examine the peculiar element $\heartsuit$, defined by 
\[
  [\heartsuit] = (1; t), \quad t = \tr \pi_Q.
\]
The minimal polynomial of $\pi_Q$ is Eisenstein and can be written as
\[
  \pi_Q^2 = t\pi_Q + u\pi, \quad u \in \OO_K^\cross.
\]
Let $d_0' = v_K(t)$. It is worth noting that
\[
  \begin{cases}
    d_0' = \frac{d_0}{2} & \text{ if } d_0 \leq 2e \\
    d_0' \geq e & \text{ if } d_0 = 2e + 1.
  \end{cases}
\]
We may scale
\[
  \heartsuit = \(1; \frac{\pi_Q^{2 d_0'}}{t}\)
\]
to be a unit.

We then rescale $\delta^\odot$ such that $\delta^{\odot(K)} = 1$; this is seen to be an optimal scaling for making the conic approximate a square, since $\lambda^\diamondsuit((1;0)) = 1$.  We then scale $\delta^{\odot(Q)}$ by $(\OO_Q^\cross)^2$ so that the level of $\delta^\odot/\heartsuit$ is manifest:
\[
  \delta^\odot = \heartsuit (1; (1 + b\pi^{j+1}) + a\pi^j\pi_Q + 2\beta) = \(1; \pi_Q^{2d_0'} \cdot \frac{1+b\pi^{j+1} + a\pi^j\pi_Q + 2\beta}{t}\),
\]
where $j \in \{0, 1, 2,\ldots, e-1, \infty\}$ controls $\ell(\delta^\odot\heartsuit)$, and where $a \in \OO_K^\cross$, $b \in \OO_K$, $\beta \in \OO_R$.

To compute the squareness of $\M$, it suffices to check the squareness of $\M(\xi^\odot)$ where $\xi^\odot$ ranges over a basis of $\OO_R$. We compute
\begin{align*}
  \M(1;0) &= \lambda^\diamondsuit(1;0) \\
  &= 1 \\
  \M\(0;\pi_Q^{-d_0'}\pi^{\ceil{\frac{d_0'}{2}}}\) &= \pi^{2\ceil{d_0'/2}}\lambda^\diamondsuit\(0; \frac{1 + b\pi^{j+1} + a\pi^j\pi_Q + 2\beta}{t}\) \\
  &\equiv a \cdot \frac{\pi^{2\ceil{\frac{d_0'}{2}}}}{t} \cdot \pi^j \mod 2 \\
  \M\(0;\pi_Q^{1-d_0'}\pi^{\floor{\frac{d_0'}{2}}}\) &= \pi^{2\floor{d_0'/2}}\lambda^\diamondsuit\(0; \pi_Q^2 \cdot \frac{1 + b\pi^{j+1} + a\pi^j\pi_Q + 2\beta}{t}\) \\
  &\equiv (1 + \pi^{j+1} b)\pi^{2\floor{\frac{d_0'}{2}}} - au \cdot \frac{\pi^{2\floor{\frac{d_0'}{2}} + 1}}{t} \cdot \pi^j + a\pi^{j + 2 \floor{\frac{d_0'}{2}}} \mod 2,
\end{align*}
noting that the three arguments to $\M$ form a basis of $\OO_R$ (being associates to $(1;0)$, $(0;1)$, and $(0; \pi_Q)$ in an order depending on the parity of $d_0'$). We then observe that:
\begin{itemize}
  \item All values are squares modulo $\pi^j$, and the first value assures us that the conic is scaled as squarely as possible.
  \item If $j = \infty$, then all values are squares modulo $2$ so the squareness of the conic is $e$. For the remaining bullet-points we assume that $\floor{j/2} < \floor{e/2}$:
  \item If $d_0'$ is even and $j$ is odd, the second value is not a square modulo $\pi^j$.
  \item If $d_0'$ is even and $j$ is even, the third value is not a square modulo $\pi^{j+1}$.
  \item If $d_0'$ is odd and $j$ is odd, the third value is not a square modulo $\pi^j$.
  \item If $d_0'$ is odd and $j$ is even, the second value is not a square modulo $\pi^{j+1}$.
\end{itemize}
Thus in all cases, the squareness of the conic is $\min\{2\floor{j/2} + 1, e\}$, as desired.
\end{wild}
\end{proof}
\end{normalsize}

We have the following corollary:
\begin{lem}\label{lem:Brauer_const}
The Brauer class $\epsilon(\delta)$ takes the same value for all $\delta$ in the coset $\hat\omega_C \heartsuit \diamondsuit \L_{\ceil{e'/2}}$.
\end{lem}
\begin{proof}
Let $[\delta] = [\kappa \hat\omega_C \heartsuit \diamondsuit]$, where
\[
  \kappa \equiv 1 \mod \begin{cases}
    \pi^{2\ceil{e/2} + 1}, & \text{$R$ unramified} \\
    \pi^{2\ceil{e/2}}, & \text{$R$ ramified.}
  \end{cases}
\]
Then $[\delta^\odot \heartsuit] \in \L_{\ceil{e'/2}}$, so the conic $\M = \M_\delta$ has maximal squareness $e$. If $e$ is even, we know that $\epsilon(\delta) = 1$ by Proposition \ref{prop:conic_1_mod_2}. So we may assume that $e$ is odd.

 Now $\kappa \equiv 1 \mod 2\pi$, so the associated bilinear forms
 \[
   \M(\xi,\eta) = \lambda^\diamondsuit(\delta^\odot\xi\eta) \textand \M'(\xi,\eta) = \lambda^\diamondsuit(\kappa\delta^\odot\xi\eta)
 \]
 are congruent modulo $2\pi$. So by Proposition \ref{prop:conic_perturb}, the two conics have the same Brauer class.
\end{proof}

\subsection{\texorpdfstring{$\N_{11}$}{N11}}
In this section, we will transform the $\N_{11}$-condition, which says that all coordinates of $\bar\omega_C^{-1} \cdot \xi_1^2$ are congruent modulo $\pi^{n_{11}}$, into a more manageable form.

We will sometimes need to make some subtle reductions, and thus we make the following definition:

\begin{defn}
  A \emph{first vector problem} $\P$ consists of a choice of resolvent algebra $R$ and as much of the discrete data as is needed to make $\M_{11}$ and $\N_{11}$ meaningful: the coarse coset $\delta_0\L_0$, the resolvent extender vector $\theta_1$ (which determines $\hat\omega_C$ and $\bar s$), and the moduli $m_{11}$ and $n_{11}$. These are required to satisfy the requisite integrality properties, which essentially say that
  \[
  B_{\theta_1}(m_{11},n_{11}-\bar s) = \pi^{ m_{11}} \OO_K {\theta_1} + \pi^{n_{11}-\bar s} \OO_K \theta_2
  \]
  is a subset of $R$, but are otherwise untethered from a cubic or a quartic ring. The \emph{answer} to a first vector problem is the weighting
  \[
  W_{\P} = W_{{\theta_1},m_{11},n_{11}} : \delta_0 \L_0 \to \QQ_{\geq 0}
  \]
  that attaches to each quartic algebra $\delta \in \delta_0 \L_0$ the volume of $\xi'_1 \in \PP(\OO_R)$ such that the corresponding $\xi_1 = \xi'_1\gamma_{1}$ satisfies the resolvent conditions
  \begin{alignat*}{2}
    \M_{11} &:& \tr(\xi_1^2) \equiv 0 &\mod \pi^{ m_{11}} \\
    \N_{11} &:& \quad \text{All coordinates of } \bar\omega_C^{-1} \cdot \xi_1^2 \text{ are congruent} &\mod \pi^{n_{11}}.
  \end{alignat*}
  We normalize volumes so that
  \[
  \mu(R) = 1 \textand \mu(\PP(R)) = 1 + \frac{1}{q} + \frac{1}{q^2}.
  \]
\end{defn}

We write $W^{\odot}$ instead of $W$ when we wish to normalize instead by the vector $\xi_1^\odot = \xi_1'/\gamma^{\odot}$ in Lemma \ref{lem:tfm_conic}. Thus
\[
W^{\odot}_{\P} = q^{v\(N_{R/K}(\gamma^{\odot})\)}W_{\P}.
\]

First vector problems will be sorted into \emph{zones}, given by linear inequalities on $m_{11}$ and $n_{11}$, and having the properties that within each zone, the answer has a uniform description. Zones will be named by colors in such a way that a brightening of the color correlates with a lowering of $m_{11}$ and/or $n_{11}$ and an increase in the answer. Brightening is governed by the following poset:

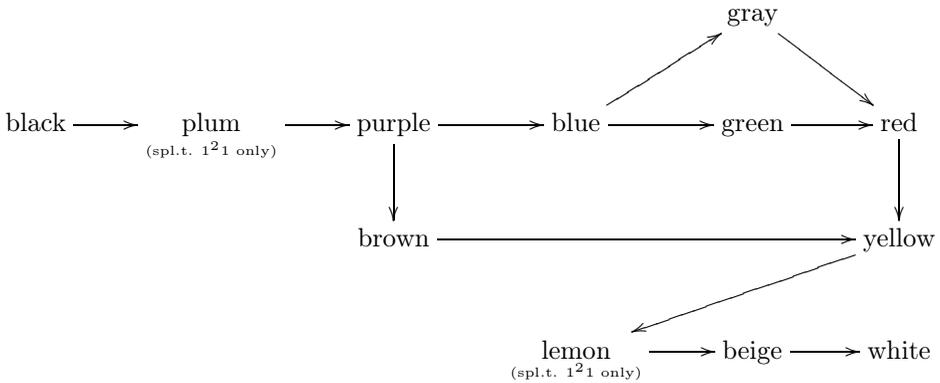
\[
\xymatrix{
  &&&& \text{gray} \ar[rd] \\
  \text{black} \ar[r] &
  \underset{\text{\tiny (spl.t. $1^21$ only)}}{\text{plum}}  \ar[r] &
  \text{purple} \ar[r] \ar[d] &
  \text{blue} \ar[r] \ar[ru] &
  \text{green} \ar[r] &
  \text{red} \ar[d] \\
  && \text{brown} \ar[rrr]
  &&& \text{yellow} \ar[lld] \\
  &&& \underset{\text{\tiny (spl.t. $1^21$ only)}}{\text{lemon}}  \ar[r] &
  \text{beige} \ar[r] &
  \text{white}
}
\]

\begin{lem} \label{lem:N11}
  Fix the data of a first vector problem $\P$ in such a way that $\N_{11}$ is active with
  \begin{equation}\label{eq:n11_leq_2e}
    0 < n_{11} \leq 2e < m_{11}^{\odot}
  \end{equation}
  and there is a solution $\xi_0$ to $\P$. Let $\xi^{\odot}_0 = \xi_0 / \gamma_1^\odot$ be its transform. Then there is an $n^\odot \in \ZZ_{\geq 0}$ such that, for any $\xi_1$ satisfying $\M_{11}$,
  \[
    \xi_1 \text{ satisfies } \N_{11} \iff \xi_1^\odot \equiv \xi_0^\odot \mod \pi^{n^\odot}.
  \]
  The value of $n^\odot$ is given as follows:
  \begin{itemize}
    \item In unramified splitting types,
    \[
      n^\odot = \ceil{\frac{n_{11}}{2}}.
    \]
    \item In splitting type $(1^3)$,
    \[
      n^\odot = \ceil{\frac{n_{11}}{2} - \frac{h_1}{3}}.
    \]
\begin{wild}
      \item In splitting type $(1^21)$,
    \[
      n^\odot = \max\left\{\ceil{\frac{n_{11} - d_0}{2}} + \frac{h_1}{2}, 0\right\}.
    \]
\end{wild}
  \end{itemize}
\end{lem}
\begin{rem}
  The condition $n_{11} \leq 2e$ (which, as we will see, restricts us to the blue, green, red, yellow, and lemon zones) can be removed, but then our conclusion must be that there is a family $\{\xi_{0(1)}, \ldots, \xi_{0(r)}\}$ of basic solutions, $r \in \{1, 2, 4\}$. The formula for $n^\odot$ becomes more complicated, and we will be able to solve these zones by other means.
\end{rem}
\begin{proof}[Proof of Lemma \ref{lem:N11}]
  In view of Lemma \ref{lem:conic_lift}, we may assume $m_{11} = \infty$, replacing $\xi_1^\odot$ by a value in the same $e$-pixel that satisfies $\M^\odot(\xi_1^\odot) = 0$ exactly.
  
  Let $\xi_{0}$ be a fixed solution to $\P$, and let $\xi_1$ be any solution to $\M_{11}$. Observe that $\xi_1^2$ and $\xi_{0}^2$ are both traceless, so their wedge product $\xi_1^2 \wedge \xi_{0}^2$ is a scalar multiple of $(1;1;1) \in \OO_{\bar K}^3$. (Here we identify $\Lambda^2 \OO_{\bar K}^3$ with $\OO_{\bar K}^3$ via the trace pairing and standard orientation, so that the wedge product is given by the same formula as the cross product on $\RR^3$.) Let $\{\xi_{0}^2, \alpha\}$ be an $\OO_{\bar K}$-basis for the traceless plane in $\OO_{\bar K}^3$. Write
  \[
  \xi_1^2 = c_0 \xi_{0}^2 + c_1 \alpha.
  \]
  The coefficient $c_1$ controls how far $\xi_1$ deviates from $\xi_{0}$ and thus the satisfaction of $\N_{11}$:
  \begin{align}
    \N_{11} &\iff \text{All coordinates of} \quad \bar\omega_C^{-1} \cdot \xi_1^2 \quad \text{are congruent} \mod \pi^{n_{11}} \nonumber \\
    &\iff \text{All coordinates of} \quad c_1 \bar\omega_C^{-1} \alpha \quad \text{are congruent} \mod \pi^{n_{11}} \nonumber
  \end{align}
  We claim that the element $\bar\omega_C^{-1} \alpha \in \OO_{\bar K}^3$ does \emph{not} have all coordinates congruent mod $\mm_{\bar K}$:
  \begin{itemize}
    \item If $\bar s = 0$, then $\bar\omega_C$ is a unit so this is equivalent to $\alpha$ and $\xi_0^2$ being linearly independent modulo $\mm_{\bar K}$;
    \item If $\bar s > 0$, then $\bar\omega_C \sim (1; \pi^{\bar s}; \pi^{\bar s})$, so $\xi_0^{(K)}$ has positive valuation. Hence $\alpha^{(K)}$ and $(\bar\omega_C \alpha)^{(K)}$ are units, while $(\bar\omega_C \alpha)^{(Q)}$ is not.
  \end{itemize}
  
  Consequently
  \begin{align}
    \N_{11} &\iff c_1 \equiv 0 \mod \pi^{n_{11}} \label{eq:x_c1_tricky} \\
    &\iff \xi_1^2 \wedge \xi_{0}^2 \equiv 0 \mod \pi^{n_{11}}. \label{eq:N11_wedge^2}
  \end{align}
  
  Now \eqref{eq:N11_wedge^2} is advantageous, because the wedge product $\xi_1^2 \wedge \xi_{0}^2$ has all its coordinates equal, so we can test $\N_{11}$ by looking at any one of them. We have (coordinate indices mod $3$)
  \begin{align*}
    \(\xi_1^2 \wedge \xi_{0}^2\)^{(i)} &= {\xi_1^{(i+1)}}^2 {\xi_{0}^{(i-1)}}^2 - {\xi_1^{(i-1)}}^2 {\xi_{0}^{(i+1)}}^2 \\
    &= \(\xi_1^{(i+1)} \xi_{0}^{(i-1)} - \xi_1^{(i-1)} \xi_{0}^{(i+1)}\)
    \(\xi_1^{(i+1)} \xi_{0}^{(i-1)} + \xi_1^{(i-1)} \xi_{0}^{(i+1)}\)
  \end{align*}
  and the two factors are congruent modulo $2$, so, since $n_{11} \leq 2e$,
  \begin{equation}
    \N_{11} \iff \(\xi_1 \wedge \xi_{0}\)^{(i)} \equiv 0 \mod \pi^{n_{11}/2}. \label{eq:N11_wedge}
  \end{equation}
  We now examine this for each coordinate $i$ in turn, and for each splitting type.

\paragraph{Unramified.}
We first dispose of the case $h_1 = 1$. Here the conic is tiny, and by Lemma \ref{lem:conic_pi}, all solutions $\xi_1$ satisfy
\begin{align*}
  \xi_1^\odot &\equiv \xi_0^\odot \mod \pi^e \\
  \xi_1 &\equiv \xi_0 \mod \(\pi^{e + 1/2}; \pi^e; \pi^e\) \\
  \xi_1^2 &\equiv \xi_0^2 \mod \(\pi^{2e + 1}; \pi^{2e}; \pi^{2e}\).
\end{align*}
Since we are assuming $n_{11} \leq 2e$, we find that $\N_{11}$ is automatic, and $n^\odot \leq e$ may be chosen at will.

Now assume that $h_1 = 0$. Here $\delta^\odot \in \OO_R^\cross$ and $\xi_0, \xi_1 \in \sqrt{\delta^\odot}\OO_R$ are $\OO_{\bar K}$-primitive. We scale $\xi_1$ by $\OO_K^\cross$ to be as close to $\xi_0$ as possible. Then $k = v(\xi_1 - \xi_0)$ is an integer, and
\[
  \xi_0 \textand \frac{\xi_1 - \xi_0}{\pi^k}
\]
are linearly independent elements of $\sqrt{\delta^\odot}\OO_R$. In particular, their wedge product is primitive, so
\[
  v(\xi_0 \wedge \xi_1) = k.
\]
Hence 
\begin{align*}
  \N_{11} &\iff k \geq \frac{n_{11}}{2} \\
  &\iff k \geq \ceil{\frac{n_{11}}{2}} \\
  &\iff \xi_1 \equiv \xi_0 \mod \pi^{\ceil{n_{11}/2}} \\
  &\iff \xi_1^\odot \equiv \xi_0^\odot \mod \pi^{\ceil{n_{11}/2}},
\end{align*}
as desired.
\paragraph{Splitting type $(1^3)$.} Scale $\xi_1 \in \bar\zeta_3^{-h_1} \sqrt{\delta^\odot} \OO_R^\cross$ to be as close to $\xi_0$ as possible, and consider the valuation $k = v(\xi_1 - \xi_0) \in \frac{1}{3}\ZZ$. If $k \geq e$, then both $\N_{11}$ and its claimed transformation are easily seen to hold, so assume that $k < e$. Then:
\begin{itemize}
  \item We cannot have $k \in \ZZ$, for then rescaling $\xi_1$ would bring it closer to $\xi_0$.
  \item If $k \in \ZZ + h_1/3$, then $v(\xi_1^2 - \xi_0^2) = 2k \in \ZZ - h_1/3$, and $\xi_1^2 - \xi_0^2$ has its first-order term a multiple of
  \begin{align*}
    &\(\bar\zeta_3^{-h_1}\sqrt{\delta^\odot}\)^2 \cdot \pi_R^{2k} \\
    &= \bar\zeta_3^{h_1} \delta^\odot \cdot \pi^{2k} \bar\zeta_3^{-h_1} \\
    &= \delta^\odot \pi^{2k},
  \end{align*}
  which has trace $\sim \pi^{2k}$, contradicting the constraint that both $\xi_1^2$ and $\xi_0^2$ are traceless.
\end{itemize}
Hence $k \in \ZZ - h_1/3$. Note that
\[
  v(\xi_1 \wedge \xi_0) = v\big((\xi_1 - \xi_0) \wedge \xi_0\big) = k,
\]
since the leading terms of $\xi_0$ and $\xi_1 - \xi_0$ are multiples of different powers of $\zeta_3$. So
\begin{align*}
  \N_{11} &\iff k \geq n_{11}/2 \\
  &\iff k \geq -\frac{h_1}{3} + \ceil{\frac{n_{11}}{2} - \frac{h_1}{3}} - r, \quad \text{any $r$ in the interval $[0,1)$} \\
  &\iff \xi_1 \equiv \xi_0 \mod \pi^{-\frac{h_1}{3} + \ceil{\frac{n_{11}}{2} - \frac{h_1}{3}} - r} \\
  &\iff \xi_1^\odot \equiv \xi_0^\odot \mod \pi^{\frac{1}{3} + \ceil{\frac{n_{11}}{2} - \frac{h_1}{3}} - r}.
\end{align*}
We take $r = 1/3$ to get the claimed
\[
  n^\odot = \ceil{\frac{n_{11}}{2} - \frac{h_1}{3}}.
\]

\begin{wild}
  \paragraph{Splitting type $(1^21)$.}
When $h_1 = 1$, the conic is tiny and $n^\odot$ is immaterial for similar reasons to those in unramified splitting type, $h_1 = 1$.

When $h_1 = 3$, we claim that $\M_{11}$ has no solutions. Observe that $n_{11} > 0$, so $m_{11} > d_0/2$ and $m_{11}^\odot = m_{11} - p^\odot = m_{11} - (d_0/3)/2 \geq 2$. Then write
\begin{align*}
  \M^\odot(\xi_1^\odot) &= \(\delta^\odot {\xi_1^\odot}^2\)^{(K)} - \I\(\(\delta^\odot {\xi_1^\odot}^2\)^{(Q)}\).
\end{align*}
The first term has even valuation since $\delta^{\odot(K)}$ is a unit. But
\[
  \delta^\odot {\xi_1^\odot}^2 \sim \pi_Q^3
\]
has $\I$-value $\sim \pi$, so $\M^\odot$ is unsatisfiable mod $\pi^2$. Therefore we may assume $h_1 \in \{0, 2\}$.

We have
\begin{align*}
  \N_{11} &\iff \xi_1^{(Q1)} \xi_{0}^{(Q2)} - \xi_1^{(Q2)} \xi_{0}^{(Q1)} \equiv 0 \mod \pi^{n_{11}/2} \\
  &\iff \frac{\xi_1^{(Q1)}}{\xi_{0}^{(Q1)}} - \frac{\xi_1^{(Q2)}}{\xi_{0}^{(Q2)}} \equiv 0 \mod \pi^{n_{11}/2} \\
  &\iff \I\(\frac{\xi_1^{(Q)}}{\xi_{0}^{(Q)}}\) \equiv 0 \mod \pi^{\frac{n_{11} - d_0}{2}},
\end{align*}
by the definition of $\I$. Note that the left-hand side belongs to $\OO_K$, so
\[
  \N_{11} \iff \I\(\frac{\xi_1^{(Q)}}{\xi_{0}^{(Q)}}\) \equiv 0 \mod \pi^{\ceil{\frac{n_{11} - d_0}{2}} - r}
\]
for any real $r$, $0 \leq r < 1$. We will pick $r$ below. Now, transforming to the coordinate $\xi_1^\odot$,
\begin{align}
  \N_{11} &\iff \I\(\frac{\xi_1^{\odot(Q)}}{\xi_{0}^{\odot(Q)}}\) \equiv 0 \mod \pi^{\ceil{\frac{n_{11} - d_0}{2}} - r} \nonumber \\
  &\iff \frac{\xi_1^{\odot(Q)}}{\xi_{0}^{\odot(Q)}} \equiv c \mod \pi^{\ceil{\frac{n_{11} - d_0}{2}} + \frac{1}{2} - r} \\
  &\iff \xi_1^{\odot(Q)} \equiv c \cdot \xi_{0}^{\odot(Q)} \mod \pi^{\ceil{\frac{n_{11} - d_0}{2}} + \frac{h_1}{4} + \frac{1}{2} - r}, \quad \text{some } c \in \OO_K^\cross \label{eq:N11_Q}
\end{align}
We may scale so that $c = 1$. We take $r = \frac{2 - h_1}{4} \in \{0, 1/2\}$, ensuring that
\[
  n^\odot = \ceil{\frac{n_{11} - d_0}{2}} + \frac{h_1}{2}
\]
is an integer. We claim it works overall, that is, that we can drop the superscript $(Q)$'s in \eqref{eq:N11_Q}. To do this, we look at the $Q1$-coordinate. We have
\begin{align*}
  \N_{11} &\iff \xi_1^{(K)} \xi_{0}^{(Q1)} - \xi_1^{(Q1)} \xi_{0}^{(K)} \equiv 0 \mod \pi^{n_{11}/2} \\
  &\iff \xi_1^{(K)} \equiv \xi_{0}^{(K)} \cdot \frac{\xi_1^{(Q1)}}{\xi_{0}^{(Q1)}} \mod \pi^{n_{11}/2} \\
  &\iff \xi_1^{\odot(K)} \equiv \xi_{0}^{\odot(K)} \cdot \frac{\xi_1^{\odot(Q1)}}{\xi_{0}^{\odot(Q1)}} \mod \pi^{\frac{2n_{11} - d_0 + h_1}{4}} \\
  &\implies \xi_1^{\odot(K)} \equiv \xi_{0}^{\odot(K)} \cdot \frac{\xi_1^{\odot(Q1)}}{\xi_{0}^{\odot(Q1)}} \mod \pi^{n^\odot} \\
  &\iff \xi_1^{\odot(K)} \equiv \xi_{0}^{\odot(K)} \mod \pi^{n^\odot},
\end{align*}
as desired, where the lone non-reversible step employs the inequalities $d_0 \geq 2 \geq h_1$ to get
\[
  \frac{n_{11}}{2} + \frac{-d_0 + h_1}{4} \leq \frac{n_{11} - d_0}{2} \leq n^\odot.
\]
If $n^\odot < 0$, we can obviously set $n^\odot = 0$ without changing the conclusion.

\end{wild}
\end{proof}

\section{Boxgroups}
If $\N_{11}$ is strongly active, then
\[
  \beta = \frac{\xi_1^2}{\omega_C}
\]
is a unit. By Lemma \ref{lem:to_box}, the solutions to $\P$ arise from the $\beta$ that lie in the box
\[
  x + y\pi^{n_{11} - s} \theta_1 + z\pi^{m_{11}} \theta_2, \quad x \in \pi^{-2a_1 + s}K, \quad y,z \in \OO_K.
\]
 Necessarily $-2a_1 + s \in \ZZ$ and $x \in \OO_K^\cross$. Also, $[\delta] = [\beta] \in H^1$. So the support of $\delta$ is bound up with the $H^1$-classes of units in various boxes. Certain boxes have pride of place: those for which the corresponding subset of $H^1$ is a group, which we will call a \emph{boxgroup.}

In this section, our aim is to define certain subgroups of $H^1$. We fix the resolvent data. We do not fix the discrete data, but we will reference the transformation $\gamma_{1,0}$ of the conic that occurs in Lemma \ref{lem:tfm_conic} when $L \cong K \cross R$ is the algebra for $\delta \in \L_0$ and $m_{11}$ is large enough.

\subsection{Signatures}
Recall that in Lemma \ref{lem:levels_quartic}, we filtered $H^1$ by level spaces $\L_0 \supset \L_1 \supset \cdots \supset \L_{e'}$, where
\[
  e' = \begin{cases}
    e, & \text{$R$ unramified} \\
    2e, & \text{$R$ ramified.}
  \end{cases}
\]
We would like to define some additional subgroups of $H^1$. We use the following notion.

\begin{defn}
If $S \subseteq H^1$ is a subgroup, define the \emph{signature} of $S$ to be the sequence of $e' + 2$ subgroups
\[
  S_i = \frac{S \intsec \L_i}{S \intsec \L_{i+1}} \subseteq \L_i/\L_{i+1}, \quad -1 \leq i \leq e.
\]
\end{defn}
The following subgroups $S_i \subseteq \L_i/\L_{i+1}$ will occur frequently and will be given names:
\begin{itemize}
  \item $\0$ denotes the zero subgroup $\L_{i+1}/\L_{i+1}$;
  \item $\*$ denotes the entire group $\L_{i}/\L_{i+1}$;
  \item $\tee$, in unramified resolvent for $0 \leq i < e$,  denotes the order-$q$ subgroup
  \[
    \tee = \{[1 + a\pi^{2i+1}\theta_1 : a \in \OO_K]\};
  \]
\end{itemize}
Thus, for instance, $\L_i$ for $0 \leq i \leq e'$ has signature $\0 .\0^i \*^{e'-i} . \*$. We separate the first and last terms of a signature by periods, because they carry less information than the other elements in general. In splitting types $3$ and $1^3$, we can omit these terms.

Moreover, $\L_i$ is the only subgroup with its signature. In general, however, the signature does not uniquely determine the subgroup, though it does determine the \emph{size} of the subgroup, since
\[
  \size{S} = \prod_i \size{S_i}.
\]
Note also that if $S$ has signature $S_{-1}S_0\ldots S_{e'-1}S_e$, then $S^\perp$ has signature $S_{e}^\perp S_{e'-1}^\perp \ldots S_0^\perp S_{-1}^\perp$, since the Tate pairing on $H^1$ induces a perfect pairing between $\L_{i}/\L_{i+1}$ and $\L_{e'-i-1}/\L_{e'-i}$.

In this section, our aim will be to define a family of \emph{boxgroups} in terms of which the ring totals will be written. These boxgroups will depend on the resolvent data alone. In the unramified splitting types they will have signature
\[
  \0 .\0^{\ell_0} \tee^{\ell_1} \*^{\ell_2} . \*, \quad \sum_i \ell_i = e
\]
and will be denoted by $T(\ell_0, \ell_1, \ell_2)$. In splitting type $1^3$ they will have signatures
\[
  \0 .(\0\0)^{\ell_0} (\0\*)^{\ell_1} (\*\*)^{\ell_2} . \* \textand \0 .(\0\0)^{\ell_0} (\*\0)^{\ell_1} (\*\*)^{\ell_2} . \*, \quad \sum_i \ell_i = e
\]
and will be denoted by $T_{-1}(\ell_0, \ell_1, \ell_2)$ and $T_1(\ell_0, \ell_1, \ell_2)$ respectively.

\begin{wild}
  \begin{lem} \label{lem:iota}
  Let $Q = K[\sqrt{D_0}]$ be a quadratic extension, and let $a \in K^\cross$. The level $\ell_Q(a)$ considered as an element of $Q^\cross/(Q^\cross)^2$ is determined by the levels of $a$ and $D_0 a$ in $K^\cross/(K^\cross)^2$ in the following way:
  \begin{enumerate}[$($a$)$]
    \item If $Q/K$ is unramified, then
    \[
    \ell_Q(a) = \begin{cases}
      \ell_K(a), & \ell_K(a) < e \\
      e+1, & \ell_K(a) \geq e.
    \end{cases}
    \]
    \item If $Q/K$ is ramified, then
    \[
    \ell_Q(a) = \begin{cases}
      e + \ell_K(a), & \ell_K(a) \geq e - d_0/2 \\
      e + \ell_K(D_0 a), & \ell_K(D_0 a) \geq e - d_0/2 \\
      2\ell_K(a) + d_0/2, & -1/2 \leq \ell_K(a) < e - d_0/2.
    \end{cases}
    \]
    (Here we put $\ell_K(a) = -1/2$ if $a \in \pi\OO_K^\cross$.)
    \item\label{iota:custom} In particular,
    \[
      \ell_Q(a) \geq d_0 + \floor{\ell_K(a)},
    \]
    and for any integer $m$, $0 \leq m \leq e - d_0/2$,
    \[
      \ell_K(a) \geq m \implies \ell_Q(a) \geq 2m + \frac{d_0}{2}.
    \]
  \end{enumerate}
\end{lem}
\begin{proof}
  Although a direct proof is not difficult, it is more illuminating to use the connection between levels and the discriminants of the corresponding quadratic extensions. Observe that, by Theorem \ref{thm:disc_Kummer_aff},
  \[
  (\Disc K[\sqrt{a}]) = \begin{cases}
    (\pi^{2e - 2\ell(a)}), & -1/2 \leq \ell(a) \leq e \\
    (1), & \ell(a) \geq e.
  \end{cases}
  \]
  In particular,
  \[
    \ell_K(D_0) = e - d_0/2.
  \]
  Thus, except for distinguishing levels $e$ and $e+1$, finding $\ell_Q(a)$ is equivalent to determining the relative discriminant $(\disc_Q Q[\sqrt{a}])$. Since
  \begin{equation*}
    N_{Q/K}(\disc_Q Q[\sqrt{a}]) = \frac{\disc_K Q[\sqrt{a}]}{\disc_K Q},
  \end{equation*}
  it is enough (excluding the trivial case $Q = K \cross K$) to compute the absolute discriminant $\disc_K Q[\sqrt{a}] = \disc_K K[\sqrt{D_0}, \sqrt{a}]$.
  
  Now, since the regular representation of $\ZZ/2\ZZ \cross \ZZ/2\ZZ$ is the direct sum of $1$-dimensional representations, we have by an Artin-conductor argument the identity
  \begin{align*}
    \disc_K K[\sqrt{D_0}, \sqrt{a}] = \disc_K K[\sqrt{D_0}] \cdot \disc_K K[\sqrt{a}] \cdot \disc_K K[\sqrt{aD_0}].
  \end{align*}
  This gives the required result in all cases except when $\ell(a) \geq e$ or $\ell(D_0 a) \geq e$, where it gives the incomplete conclusion $\ell_Q(a) \geq e$. Without loss of generality, $\ell(a) \geq e$. If $a$ is actually a square, it remains a square in $Q$, so assume $a = 1 + 4u$ is an intimate unit, $\tr_{k_K/\FF_2} (u \mod \pi) = 1$. If $Q$ is unramified, then $a$ is a square in $Q$: indeed we could have taken $D_0 = a$. But if $Q$ is ramified, then $a$ remains a nontrivial intimate unit in $Q$, so $\ell_Q(a) = 2e$, as claimed.
\end{proof}
\end{wild}
\subsection{Boxgroups in unramified splitting type}

\begin{lem}\label{lem:boxes_basic_ur} Let $ m, n \in \NN^+ \union \{\infty\} $, $ m \geq n > 0 $. Let $ B_{\theta_1}(m,n) $ be the box
  \[
  B_{\theta_1}(m,n) = \{\pi c_0 + \pi^{n}c_1{\theta_1} + \pi^m c_2 \theta_2 : c_i \in \OO_K \}.
  \]
  \begin{enumerate}[$($a$)$]
    \item\label{boxes_basic_ur:recenter} For every $ \xi \in 1 + B_{{\theta_1}(m,n)} $,
    \begin{align}
    B_{\theta_1}(m,n) &= \xi B_{\xi^{-1} {\theta_1}}(m,n) \label{eq:recenter_1} \\
    1 + B_{\theta_1}(m,n) &= \xi (1 + B_{\xi^{-1} {\theta_1}}(m,n)) \label{eq:recenter_2}.
    \end{align}
    \item If
    \begin{equation*}
    m \leq 2n + s,
    \end{equation*}
    then $ B_{\theta_1}(m,n) $ is closed under multiplication and the translate $ 1 + B_{\theta_1}(m,n) $ is a group under multiplication.
  \end{enumerate}
\end{lem}
\begin{proof}
  \begin{enumerate}[$($a$)$]
    \item Since $ \xi \equiv 1 $ mod $ \pi $, we can take $ \xi^{-1}\theta_2 $ in the role of $ \theta_2 $ for defining $ \xi B_{\xi^{-1} {\theta_1}}(m,n) $, which is a lattice with basis $ [\pi^n\xi, \pi^n{\theta_1}, \pi^m\theta_2] $, all of which are contained in $ B_{{\theta_1}}(m,n) $. This proves the reverse inclusion of \eqref{eq:recenter_1}, and equality follows by comparing volumes. To get \eqref{eq:recenter_2}, we add $ \xi $ to both sides and use $ \xi - 1 \in B_{\theta_1}(m,n) $ to simplify the left-hand side.
    \item Observe that $ B_{\theta_1}(m,n) $ is a lattice with basis $ [\pi^n, \pi^n{\theta_1}, \pi^m\theta_2] $. Since $ m \geq n $, the only product that does not clearly lie in the lattice is $ (\pi^{n} {\theta_1})^2 $, whose $ 1 $- and $ {\theta_1} $-components are divisible by $ \pi^{2n} $, and whose $ \theta_2 $-component is divisible by $ \pi^{2n + s} $. Since $ m \leq 2n + s $, this product lies in the lattice.
    
    Thus $ B_{\theta_1}(m,n) $ is closed under multiplication and so is $ 1 + B_{\theta_1}(m,n) $. To show the existence of inverses, simply note that
    \[
    \frac{1}{1 + \xi} = 1 - \xi + \xi^2 - \cdots
    \]
    converges to an element of $ 1 + B_{\theta_1}(m,n) $ for every $ \xi \in B_{\theta_1}(m,n) $.
\end{enumerate}
\end{proof}
\begin{lem} \label{lem:eta_ur}
Let
\[
  \square_C = \min \left\{ 2\ell(\hat\omega_C) + 1, e \right\}
\]
be the squareness of the conic for $\delta = 1$; put $\square_C = 0$ if $[\hat\omega_C] \notin \L_0$ (i.e{.} $s$ is odd).

If $\theta_1$ is translated by a suitable element of $\OO_K$ and scaled by a suitable element of $\OO_K^\cross$ (neither of which change the associated resolvent $C$), then there is an $\eta \in \OO_R$ such that
\begin{equation} \label{eq:eta_sqrt}
  \eta^2 \equiv \theta_1 \mod \pi^{s + \square_C}
\end{equation}
and
\begin{equation} \label{eq:eta_in}
  \eta \in \<1, \theta_1, \pi^{\ceil{\frac{s}{2}}} \theta_2, 2\theta_2\>.
\end{equation}
\end{lem}
\begin{proof}
If $s$ is odd, then we can scale and translate $\theta_1$ so that
\[
  \theta_1 \equiv (1; 0; 0) \mod \pi^s.
\]
Then $\eta = \theta_1$ satisfies the desired conditions.

If $s$ is even, then by definition of $\square_C$, there is a linear form $\lambda$ such that
\[
  \tr(1^\odot{\xi^\odot}^2) \equiv c \cdot \lambda(\xi^\odot)^2 \mod \pi^{\square_C}
\]
as functions of $\xi^\odot \in \OO_R$. Here $1^\odot$, the transform of $\delta = 1$ under Lemma \ref{lem:tfm_conic}, is a unit whose class in $H^1$ is $[\hat\omega_C]$; for concreteness, we may take
\[
  1^\odot = \frac{\hat\omega_C}{(\pi^s; 1; 1)}.
\]
Since $\hat\omega_C$ is traceless, the conic has a distinguished basepoint, namely $\xi'_0 = [\pi^{s/2}; 1; 1]$. Pick a $\xi' = \xi'_1$ in the kernel of $\lambda$ that does \emph{not} lie in the same $1$-pixel as the basepoint. We claim that the choice
\[
  \eta = (1; \pi^{s/2}; \pi^{s/2}) \xi'_1
\]
fulfills the conditions.

The $\theta_2$-coefficient of $\eta'^2$ is given by
\[
  \tr(\omega_C \eta^2) = \pi^s \tr(1^\odot \xi_1'^2) \equiv \pi^s \lambda(\xi'_1) = 0 \mod \pi^{s + \square_C}.
\]
Hence there are $a, b \in \OO_K$ such that
\[
  \eta^2 \equiv a + b \theta_1 \mod \pi^{s + \square_C}.
\]
We claim that $\pi \nmid b$, which makes it possible to replace $\theta_1$ by $a + b \theta_1$. Suppose not. If $s = 0$, we get $\xi'_1 \equiv \xi'_0$ modulo $\pi$, contrary to hypothesis. If $s > 0$, we get $\pi | a$ so $\pi | (\xi'_1)^{(K)}$; we also know that $\pi | (\xi'_0)^{(K)}$. But $1^\odot$ is a unit, so $\lambda'$ is a perfect linear functional, and its kernel in $\PP^2(k_K)$ intersects the line $(\xi')^{(K)} \equiv 0$ in only one $1$-pixel, a contradiction.

As for \eqref{eq:eta_in}, it can be rewritten as 
\[
  v_K(\eta^{(2)} - \eta^{(3)}) \geq \min\left\{\frac{s}{2}, e\right\}.
\]
To prove this, observe that
\[
  \(\eta^{(2)} - \eta^{(3)}\)\(\eta^{(2)} + \eta^{(3)}\) = (\eta^{(2)})^2 - (\eta^{(3)})^2 \equiv \theta_1^{(2)} - \theta_1^{(3)} \equiv 0 \mod \pi^s,
\]
with the two factors on the left-hand side congruent modulo $2 \sim \pi^e$.
\end{proof}

\begin{lem} \label{lem:boxgps_ur} Let $0 < n \leq m \leq 2e$ be integers such that
\begin{align}
m &\leq 2n + s && (\text{the gray-red inequality})\label{eq:box_red_ur} \\
m &\leq n + \square_C + s && (\text{the gray-green inequality})\label{eq:box_green_ur} \\
m &\leq e + \frac{n+s+1}{2} && (\text{the gray-blue inequality})\label{eq:box_blue_ur}.
\end{align}
Then the projection $[1 + B_{\theta_1}(m,n)]$ of $ 1 + B_{\theta_1}(m,n) $ onto $H^1$ is a subgroup of signature
\[
  \0 .\0^{\ell_0} \tee^{\ell_1} \*^{\ell_2} . \*
\]
where
\begin{align*}
\ell_0 &= \floor{\frac{n}{2}} \\
\ell_1 &= \floor{\frac{m}{2}} - \floor{\frac{n}{2}} \\
\ell_2 &= e - \floor{\frac{m}{2}}.
\end{align*}
\end{lem}
\begin{proof}
  
The gray-red inequality \eqref{eq:box_red_ur} ensures that the projection $T = [1 + B_{\theta_1}(m,n)]$ is a subgroup. It is clear that
  \[
    \L_{\floor{\frac{n}{2}}} \subseteq T \subseteq \L_{\floor{\frac{m}{2}}},
  \]
  so the signature of $T$ has the shape $  \0 .\0^{\ell_0} ?^{\ell_1} \*^{\ell_2} . \* $. Those middle $\ell_1$ components of the signature are at least $\tee$, because for $i \geq \floor{n/2}$ and for all $a \in K$, we have
  \[
    [1 + a \pi^{2i + 1} \theta_1] \in T.
  \]
  Thus the signature is at least the one claimed.
  
  To prove that equality occurs, we fix $ m \leq 2e $ and proceed by downward induction on $ n $. The base case $ n = m $ is clear since $ T = \L_{\floor{\frac{m}{2}}} $. When moving from $ n + 1 $ to $ n $, note that $ \size{T} $ can grow by at most a factor of
    \[
      [1 + B(m,n) : (1 + B(m, n+1))(1 + \pi^n \OO_K)] = q.
    \]
    If $ n $ is odd, there is nothing to prove, as we claim that $ \size{T} $ actually grows by a factor of $ q $. If $ n $ is even, we claim that $ T $ does not change. It suffices to prove that each of the $ q $ cosets in 
    \[
      \frac{(1 + B(m,n))}{(1 + B(m,n+1))(1 + \pi^n \OO_K)}
    \]
    contains a square. For $ c \in \OO_K $, consider
    \begin{align*}
    (1 + \pi^{\frac{n}{2}}c\eta)^2 = 1 + 2\pi^{\frac{n}{2}} c \eta + \pi^n c^2 \eta^2.
    \end{align*}
    The last term is $ \pi^n c^2 {\theta_1} $ up to an error in $ \pi^{n + s + \square_C}\OO_R$, which is in $\pi^m \OO_R$ by the gray-green inequality. We claim that the cross term $2\pi^{n/2} c \eta$ lies in $B(m,n+1)$ also. If $s$ is odd, this is trivial since we took $\eta = \theta_1$. Otherwise, we have
    \[
      2\pi^{n/2} \eta \in 2\pi^{n/2} \<1, {\theta_1}, \pi^{s/2}\theta_2, 2\theta_2\>
      = B_{\theta_1}\( e + \frac{n}{2} + \min\left\{\frac{s}{2}, e\right\}, e + \frac{n}{2} \)
    \]
    We get the needed inequality
    \[
      e + \frac{n}{2} + \frac{s}{2} \geq m
    \]
    from the gray-blue inequality, the difference of whose sides lies in $\ZZ + 1/2$ by parity considerations. So we have found a square in the coset $ 1 + \pi^n c^2 {\theta_1} + B(m,n+1) = (1 + \pi^n c^2 {\theta_1})(1 + B(m,n+1)) $, as desired.
\end{proof}
As a corollary, we have:
\begin{lem}\label{lem:boxgpS_ur}
    For every triple $ (\ell_0, \ell_1, \ell_2) $ of nonnegative integers satisfying
    \begin{align}
    \ell_0 + \ell_1 + \ell_2 &= e \\
    \ell_1 &\leq \ell_0 + \frac{s}{2} + 1 && (\text{the gray-red inequality})\label{eq:boxgp_red} \\
    \ell_1 &\leq \frac{s + \square_C + 1}{2} && (\text{the gray-green inequality}) \label{eq:boxgp_green} \\ 
    \ell_1 &\leq \ell_2 + \frac{s}{2} + 1, && (\text{the gray-blue inequality})\label{eq:boxgp_blue}
    \end{align}
    there is a \emph{boxgroup} $T(\ell_0,\ell_1,\ell_2) \subseteq H^1$ of signature $\0 .\0^{\ell_0} \tee^{\ell_1} \*^{\ell_2} . \*$. such that, if $m$, $n$ are integers satisfying the conditions of Lemma \ref{lem:boxgps_ur}, then
    \[
    [1 + B_{\theta_1}(m,n)] = T \( \floor{\frac{n}{2}}, \floor{\frac{m}{2}} - \floor{\frac{n}{2}}, e - \floor{\frac{m}{2}} \).
    \]
\end{lem}
\begin{proof}
    If $\ell_1 = 0$, take $T(\ell_0, 0, \ell_2) = \L_{\ell_0}$, the unique subgroup with the correct signature.
    
    Otherwise, let $ m = 2\ell_0 + 2\ell_1 $, $ n = 2\ell_0 + 1 $ in the preceding lemma. The transformation of the gray-red, gray-green, and gray-blue conditions is routine.
    
    For the last claim, note that decreasing $m$ or increasing $n$ can only make the conditions of Lemma \ref{lem:boxgps_ur} truer, with the exception of the condition $m \geq n$. If $\floor{m/2} = \floor{n/2}$, then clearly $[1 + B_{\theta_1}(m,n)] = \L_{\floor{m/2}}$, so we can assume that
    \[
      n \leq 2 \floor{\frac{n}{2}} + 1 < 2 \floor{\frac{m}{2}} \leq m.
    \]
    Clearly
    \[
    [1 + B_{\theta_1}(m,n)] \supseteq \left[1 + B_{\theta_1}\(2\floor{\frac{m}{2}}, n\)\right],
    \]
    but both sides have the same signature, so equality holds. Likewise,
    \[
    \left[1 + B_{\theta_1}\(2\floor{\frac{m}{2}}, n\)\right] \subseteq \left[1 + B_{\theta_1}\(2\floor{\frac{m}{2}}, 2\floor{\frac{n}{2}} + 1\)\right] = T\( \floor{\frac{n}{2}}, \floor{\frac{m}{2}} - \floor{\frac{n}{2}}, e - \floor{\frac{m}{2}} \),
    \]
    but both sides have the same signature, so equality holds.
\end{proof}

\subsubsection{Supplementary boxgroups}
As thus defined, all boxgroups $T$ satisfy $\L_e \subseteq T \subseteq \L_0$. Groups not satisfying these inclusions occur will be denoted as follows.

If $s > 0$, so that a distinguished splitting $R \isom K \cross Q$ exists, consider the image $\iota(K^\cross)$ of the map
\begin{align*}
  \iota\colon K^\cross/(K^\cross)^2 &\to H^1 \\
  a &\mapsto [(a;1)] = [(1;a)].
\end{align*}
We find that $\iota(K^\cross)$ has signature
\[
  \tee.\tee^{e}.\tee
\]
where the middle $e$-many $\tee$'s denote the usual subgroups
\[
  \tee_i = \{1 + \pi^{2i+1} a \theta_1 \} = \{1 + \pi^{2i+1} a (1; 0)\} \subseteq \L_i/\L_{i+1}
\]
and where the initial $\tee_{-1} = \{[1], [(1;\pi)]\} \subseteq H^1/\L_0$ and the final $\top_e = \{1 + 4a : a \in \OO_K\} \subseteq \L_e$ have size $2$ and $\size{H^0}/2$ respectively. In particular,
\[
  \size{\iota(K^\cross)} = \size{H^0} q^{e} = \sqrt{\size{H^1}}.
\]
From the explicit description in terms of the Hilbert pairing, we find that $\iota(K^\cross)$ is isotropic and hence maximally isotropic for $\<\bullet, \bullet\>$.

The group $\iota(K^\cross)$ is always important, but it does not behave well with respect to boxgroups unless $s > 2e$, in which case we give it the name $T(\emptyset, e, \emptyset)$.

If $s > 2\ell_1 $, we let
\begin{align*}
T(e - \ell_1, \ell_1, \emptyset) &= \iota(K^\cross) \intsec \L_{\ell_{e - \ell_1}} = \{[(a ; 1)] \in H^1 : a \equiv 1 \mod \pi^{2e - 2\ell_1 + 1}\} \\
T(\emptyset, \ell_1, e - \ell_1) &= \iota(K^\cross) \cdot \L_{\ell_1} = \{[(a ; \alpha)] \in H^1 : \alpha \equiv 1 \mod \pi^{2\ell_1+1}\}.
\end{align*}
\begin{wild}
  (In the first case, the equivalence of the two definitions is established using Lemma \ref{lem:iota}.)
\end{wild}
Their signatures are, respectively, $\0.\0^{\ell_0}\tee^{e-\ell_0}.\top$ and $\tee.\tee^{\ell_1}\x^{e - \ell_1}.\x$.

The restrictions on $s$ ensure that these boxgroups satisfy such natural relations as
\begin{align*}
  T(e - \ell_1, \ell_1, \emptyset) \cdot \L_e &= T(e - \ell_1, \ell_1, 0) \\
  T(\emptyset, \ell_1, e - \ell_1) \intsec \L_0 &= T(0, \ell_1, e - \ell_1),
\end{align*}
which we will often use without comment.

Finally, in all cases, we let
\begin{align*}
T(e, \emptyset, \emptyset) &= \{1\} \\
T(\emptyset, \emptyset, e) &= H^1.
\end{align*}

It will turn out that $T(\ell_0,\ell_1,\ell_2)$ and $T(\ell_2,\ell_1,\ell_0)$ are orthogonal complements whenever both are defined (Lemma \ref{lem:orth}). Actually, this is simple to prove in the case that one of $\ell_0,\ell_1,\ell_2$ is the symbol $\emptyset$. The sizes of these groups follow immediately from their signatures:
\begin{lem}\label{lem:111_T_size}
  If $T(\ell_0,\ell_1,\ell_2)$ is defined and $\ell_1 \neq \emptyset$, then
  \[
  \size{T(\ell_0,\ell_1,\ell_2)} = \size{H^0}q^{e + \ell_2 - \ell_0},
  \]
  where if $\emptyset$ occurs as either $\ell_0$ or $\ell_2$, it must be replaced by $-1/[k_K : \FF_2] = \log_q(1/2)$.
\end{lem}


\subsection{Boxgroups in splitting type \texorpdfstring{$1^3$}{1³}}
Let $h \in \{1, -1\}$ be the integer such that
\[
b_1 \in \ZZ + \frac{h}{3}, \quad \theta_1 \in \bar\zeta_3^h \OO_R^\cross, \quad \theta_2 \in \bar\zeta_3^{-h} \OO_R^\cross.
\]
Note the tight connection with the $h_i$ of Lemma \ref{lem:tfm_conic}. Namely, if $\N_{11}$ is strongly active, then $\beta$ is a unit in Lemma \ref{lem:to_box}, from which we get $m_{11} \in \ZZ$, $a_1 \in \ZZ$, and $h_1 = h$.

Let
\[
\square_C = \min \left\{ 2\floor{\frac{\ell(\hat\omega_C)}{2}} + 1, e \right\}
\]
be the squareness of the conic
\[
  \M(\xi^\odot) = \lambda^\diamondsuit(1^\odot{\xi^\odot}^2) 
\]
that occurs for $\delta = 1$. Note that there is just one conic, with $[1^\odot] = [\hat\omega_C]$, although we turn our attention to the part where $\xi^\odot \sim 1$ resp{.} $\xi^\odot \sim \pi_R^2$ according as $h_i = 1$ resp{.} $-1$. Then

\begin{lem}\label{lem:eta_1^3}
Fix a resolvent $C$. If $\theta_1$ is translated and scaled appropriately, there is an $\eta \in \bar\zeta_3^{-h} \OO_R^\cross$ such that
\[
  \eta^2 \equiv \theta_1 \mod \pi^{\square_C - \frac{2 h}{3}}.
\]
\end{lem}
\begin{proof}
By Lemma \ref{lem:squareness}, $\M$ has squareness $\square_C$, which means that there is a linear form $\lambda$ and a scalar $c \in \OO_K^\cross$ such that
\[
  \M(\xi^\odot) \equiv c \cdot \lambda(\xi^\odot)^2 \mod \pi^{\square_C}
\]
The zero locus of $\M$ modulo $\pi$, or equivalently of $\lambda$ modulo $\pi$, consists of $(q+1)$-many $1$-pixels, $q$ of which consist of units and the remaining one of elements $\xi^\odot \sim \pi_R^2$. If we were searching for a $\xi^\odot = \xi^\odot_1$ with $h_1 = h$, we would have $v(\xi^\odot) = \frac{1 - h}{3}$. We pick a $\xi^\odot$ of the \emph{other} valuation $v(\xi^\odot) = \frac{1 + h}{3}$ which lies in the kernel of $\lambda$, ensuring that
\[
  \M(\xi^\odot) \equiv 0 \mod \pi^{\square_C}.
\]

Then take
\[
  \eta = \pi^{2b_1 - \frac{h}{3}} \sqrt{\frac{1^\odot}{\hat\omega_C \diamondsuit}} \xi^\odot \in \bar{\zeta}_3^{-h} R.
\]
(See pp.~XI.304--05 for motivation and details.)  Note that
\begin{align*}
  \coef_{\theta_2}(\eta^2) &= \pi^{-4b_{1} - s + 2(2b_1 - \frac{h}{3})} \M\({\xi^\odot}^2\) \\
  &\equiv 0 \mod \pi^{\square_C - \frac{2h}{3}}
\end{align*}
Hence there are
\[
  a \in \begin{cases}
    \sqrt[3]{\pi^2}\OO_K, & h = 1 \\
    \sqrt[3]{\pi}\OO_K, & h = -1
  \end{cases}
\]
and $b \in \OO_K$ such that 
\[
  \eta^2 \equiv a + b {\theta_1} \mod \pi^{\square_C - \frac{2h}{3}}.
\]
Since $\eta$ is a unit, we must have $\pi \nmid b$, so $\theta_1$ can be replaced by $a + b\theta_1$.
\end{proof}

The following lemma is proved just like Lemma \ref{lem:boxes_basic_ur}.

\begin{lem}\label{lem:boxes_basic_1^3} If $m$ and $n$ satisfy
  \[
    0 < n < m < 2e, \quad m \in \ZZ - \frac{h}{3}, \quad n \in \ZZ + \frac{h}{3},
  \]
  Let $ B_{\theta_1}(m,n) $ be the box
  \[
  B_{\theta_1}(m,n) = \{\pi c_0 + \pi^{n} c_1 {\theta_1} + \pi^m c_2 \theta_2 : c_0,c_1,c_2 \in \OO_K \} \subseteq \OO_R.
  \]
  \begin{enumerate}[$($a$)$]
    \item For every $ \xi \in 1 + B_{{\theta_1}(m,n)} $,
    \begin{align}
    B_{\theta_1}(m,n) &= \xi B_{\xi^{-1} {\theta_1}}(m,n) \\
    1 + B_{\theta_1}(m,n) &= \xi (1 + B_{\xi^{-1} {\theta_1}}(m,n)).
    \end{align}
    \item If
    \begin{equation*}
    m \leq 2n,
    \end{equation*}
    then $ B_{\theta_1}(m,n) $ is closed under multiplication and the translate $ 1 + B_{\theta_1}(m,n) $ is a group under multiplication.
  \end{enumerate}
\end{lem}
Our goal is to study the projection of $1 + B_{\theta_1}(m,n)$ onto $H^1$. The following yields the conditions under which a useful group is formed thereby:
\begin{lem} \label{lem:boxgps_1^3} Let $0 < n < m \leq 2e$ be rational numbers with $ m \in \ZZ - \frac{h}{3}, n \in \ZZ + \frac{h}{3}$. Assume that
  \begin{align}
  m &\leq 2n && (\text{the gray-red inequality}) \label{eq:box_red_1^3} \\
  m &\leq n + \square_C - \frac{2h}{3} && (\text{the gray-green inequality}) \label{eq:box_green_1^3} \\
  m &\leq \frac{n + 1}{2} + e && (\text{the gray-blue inequality})\label{eq:box_blue_1^3}.
  \end{align}
  Write
  \[
  \Dot m = m - \frac{2h}{3} \in \ZZ, \quad \Dot n = n + \frac{2h}{3} \in \ZZ.
  \]
  Then the projection $T = [1 + B_{\theta_1}(m,n)]$ of $ 1 + B_{\theta_1}(m,n) $ onto $H^1$ is a subgroup of signature
  \[
    . (\0\0)^{\ell_0} (\*\0)^{\ell_1} (\*\*)^{\ell_2}.
  \]
  or
  \[
    . (\0\0)^{\ell_0} (\0\*)^{\ell_1} (\*\*)^{\ell_2}.
  \]
  for $h = 1$ and $h = -1$ respectively, where
  \begin{align*}
  \ell_0 &= \floor{\frac{\Dot n}{2}} \\
  \ell_1 &= \floor{\frac{\Dot m}{2}} - \floor{\frac{\Dot n}{2}}  \\
  \ell_2 &= e - \floor{\frac{\Dot m}{2}} 
  \end{align*}
  except for the case $h = 1$, $\Dot m = 2i + 1$, $\Dot n = 2i + 2$ ($i \in \ZZ$), where $T = \L_{2i+1}$.
\end{lem}
\begin{proof}

Suppose $h = 1$. For $a \in \OO_K^\cross$, there are elements in $T$ of the form
\begin{align*}
[1 + a \pi^{2i + \frac{1}{3}}\theta_1], \quad i &\geq \ceil{\frac{3n - 1}{6}} = \floor{\frac{\Dot n}{2}} \\
[1 + a \pi^{2i + \frac{5}{3}}\theta_2], \quad i &\geq \ceil{\frac{3m - 5}{6}} = \floor{\frac{\Dot m}{2}}.
\end{align*}
Likewise, in the case $h = -1$, there are elements in $T$ of the form
\begin{align*}
[1 + a \pi^{2i + \frac{5}{3}}\theta_1], \quad i &\geq \ceil{\frac{3n - 5}{6}} = \floor{\frac{\Dot n}{2}} \\
[1 + a \pi^{2i + \frac{1}{3}}\theta_2], \quad i &\geq \ceil{\frac{3m - 1}{6}} = \floor{\frac{\Dot m}{2}}.
\end{align*}
This shows that the signature of $T$ is at least as large as claimed. 

To show equality, we fix $m < 2e$ and proceed by downward induction on $n$. The base case $n = m - 1/3$ (for $h = 1$) or $n = m - 2/3$ (for $h = -1$) is clear since $T = \L_{\Dot m - 1}$ or $\L_{2\floor{\Dot m/2} + 1}$ respectively. When moving from $ n + 1 $ to $ n $, note that $ \size{T} $ can grow by at most a factor of
\[
[1 + B(m,n) : (1 + B(m, n+1))(1 + \pi^{\ceil{n}} \OO_K)] = q.
\]
If $ \Dot n$ is odd, there is nothing to prove, as we claim that $ \size{T} $ actually grows by a factor of $ q $. If $ \Dot n $ is even, we are claiming that $ T $ does not change. It suffices to prove that each of the $ q $ cosets in 
\[
\frac{(1 + B(m,n))}{(1 + B(m,n+1))(1 + \pi^{\ceil{n}} \OO_K)}
\]
contains a square. Recall the approximate square root $\eta \in \bar\zeta_3^{-h} \OO_R^\cross$ from Lemma \ref{lem:eta_1^3}, which satisfies
\[
\eta^2 \equiv \theta_1 \mod \pi^{\square_C - \frac{2 h}{3}}.
\]
For $ c \in \OO_K $, consider
\begin{equation} \label{eq:x_square_1^3}
(1 + \pi^{\frac{n}{2}}c\eta)^2 = 1 + 2\pi^{\frac{n}{2}} c \eta + \pi^n c^2 \eta^2.
\end{equation}
The last term is $ \pi^n c^2 {\theta_1} $ up to an error in $ \pi^{n + \square_C - \frac{2h}{3}}\OO_R$. To say that this is in $\pi^m \OO_R$, we need the inequalities
\[
  m \leq n + \square_C - \frac{2h}{3} \textand m \leq n + e - \frac{2h}{3}.
\]
The first of these is \eqref{eq:box_green_1^3}, and the second follows easily from \eqref{eq:box_red_1^3} and \eqref{eq:box_blue_1^3}. We claim that the middle term of \eqref{eq:x_square_1^3} lies in $ \pi^m \OO_R $ also, that is,
\[
  m \leq \frac{n}{2} + e.
\]
This follows from \eqref{eq:box_blue_1^3} and the fact that $m - n/2 \in \ZZ$. So we have found a square in the coset $ 1 + \pi^n c^2 {\theta_1} + B(m,n+1) = (1 + \pi^n c^2 {\theta_1})(1 + B(m,n+1)) $, as desired.
\end{proof}

As a corollary, just like Lemma \ref{lem:boxgpS_ur}, we get the following:
\begin{lem}\label{lem:boxgpS_1^3}
  For every triple $ (\ell_0, \ell_1, \ell_2) $ of nonnegative integers satisfying
  \begin{align}
  \ell_0 + \ell_1 + \ell_2 &= e \\
  \ell_1 &\leq \ell_0 && (\text{the gray-red inequality})\label{eq:boxgp_red_1^3} \\
  \ell_1 &\leq \frac{\square_C + 1}{2} - h && (\text{the gray-green inequality})\label{eq:boxgp_green_1^3} \\ 
  \ell_1 &\leq \ell_2, && (\text{the gray-blue inequality})\label{eq:boxgp_blue_1^3}
  \end{align}
  there is a \emph{boxgroup} $T_h(\ell_0,\ell_1,\ell_2) \subseteq H^1$ of signature
    \[
  . (\0\0)^{\ell_0} (\*\0)^{\ell_1} (\*\*)^{\ell_2}.
  \]
  or
  \[
  . (\0\0)^{\ell_0} (\0\*)^{\ell_1} (\*\*)^{\ell_2}.
  \]
  for $h = 1$ and $h = -1$ respectively, such that, if $m$, $n$ are rational numbers satisfying the conditions of Lemma \ref{lem:boxgps_1^3}, then
  \[
  [1 + B_{\theta_1}(m,n)] = T \( \ell_0', \ell_1', \ell_2' \)
  \]
  where $\ell_0', \ell_1', \ell_2'$ are the numbers $\ell_0, \ell_1, \ell_2$ defined in Lemma \ref{lem:boxgps_1^3}. 
\end{lem}
\begin{proof}
If $\ell_1 = 0$, take $T(\ell_0, 0, \ell_2) = \L_{\ell_0}$. Otherwise, take 
\[
  \Dot n = 2\ell_0 + 1 \in 2\ZZ + 1, \quad \Dot m = 2e - 2\ell_2 \in 2\ZZ,
\]
In other words,
\[
  n = 2\ell_0 + 1 - \frac{2h}{3}, \quad m = 2e - 2\ell_2 + \frac{2h}{3}.
\]
Conditions \eqref{eq:boxgp_red_1^3}--\eqref{eq:boxgp_blue_1^3} immediately imply \eqref{eq:box_red_1^3}--\eqref{eq:box_blue_1^3}. Just as in Lemma \ref{lem:boxgpS_ur}, we then argue that increasing $m$ by $1$ (resp{.} decreasing $n$ by $1$), if it does not violate \eqref{eq:box_red_1^3}--\eqref{eq:box_blue_1^3}, yields a boxgroup of the same signature that is contained in (resp{.} contains) $T_h(\ell_0, \ell_1, \ell_2)$ and thus must equal $T_h(\ell_0, \ell_1, \ell_2)$.
\end{proof}
The subscript ``$h$'' in $T_h(\ell_0, \ell_1, \ell_2)$ is logically superfluous, because ${\theta_1}$ is fixed. But it allows the following manipulation. Define
\[
  T_{-1}(\ell_0, \ell_1, \ell_2) = T_1 (\ell_0 + \ell_1, -\ell_1, \ell_1 + \ell_2)
\]
for all $(\ell_0, \ell_1, \ell_2)$ for which either side has been defined. Note that $T_{1}(\ell_0, 0, e - \ell_0) = T_{-1}(\ell_0, 0, e - \ell_0) = \L_{2\ell_0}$ already fulfills this relation, while allowing
\[
  T_1(\ell_0, -1, e+1 - \ell_0) = T_{-1}(\ell_0 - 1, 1, e - \ell_0) = \L_{2\ell_0+1}
\]
saves us the trouble of excluding the case $h = 1$, $\Dot m = 2i + 1$, $\Dot n = 2i + 2$ from Lemma \ref{lem:boxgps_1^3}. We do not use any other boxgroups with negative $\ell_1$ within this paper, but in the code we do, converting everything to a $T_{-1}$.

\begin{wild}
\subsection{Boxgroups in splitting type \texorpdfstring{$1^21$}{1²1}}
In splitting type $1^21$, the interaction between $d_0$ and the other parameters of the resolvent leads to a complexity that we must address first.

\begin{lem}\label{lem:types_1^21}
The values of $d_0$ and $s'$ constrain the value of $\ell(\hat\omega_C \diamondsuit\heartsuit)$, and hence of $\square_C$, as follows:
\begin{enumerate}[(A)]
  \item\label{type:A} If $s' = -1/2$, then $\ell(\hat\omega_C \diamondsuit\heartsuit) = -1/2$, so $\square_C = 0$.
  \item\label{type:B} If $s' < d_0/2 - 1$, then $\ell(\hat\omega_C \diamondsuit\heartsuit) = s'$, so $\square_C = 2\floor{s'/2}+1$.
  \item\label{type:C} If $s' = d_0/2 - 1$, then $\hat\omega_C \diamondsuit\heartsuit \in \L_{d_0/2 - 1}$, so $\square_C \geq 2\floor{(d_0-2)/4} + 1$.
  \item\label{type:D} If $s' \geq d_0/2$ and either $s' \equiv d_0/2$ mod $2$ or $d_0 = 2e + 1$, then $\hat\omega_C \diamondsuit\heartsuit \in \L_{\floor{d_0/2}}$, so $\square_C \geq 2\floor{d_0/4} + 1$.
  \item\label{type:E} If $s' \geq d_0/2$ and $s' \equiv d_0/2 + 1$ mod $2$, then $\hat\omega_C \diamondsuit\heartsuit \in \spadesuit \L_{d_0/2}$, so $\square_C = 2\floor{(d_0 - 2)/4} + 1$. Here $\spadesuit$ is the class in $H^1$ of $\iota(\pi) = (1; \pi; \pi)$, well defined up to $\iota(\OO_K^\cross) \subseteq \L_{d_0/2}$.
\end{enumerate}
\end{lem}
\begin{proof}
Case \ref{type:A} follows immediately from Lemma \ref{lem:tfm_conic}. In the remaining cases, $b_1 \in \ZZ$ and $\theta_1 \in \OO_R$. By translation we may assume that $\theta_1 = (0; c_0 + c_1\pi_Q)$ for some $c_0 \in \OO_K^\cross$, $c_1 \in \OO_K$. We have
\[
  s' = v(c_1).
\]
Let $\iota_Q : Q \to Q$ denote the Galois involution. We compute
\begin{align*}
  \hat\omega_C &= \pi^{4b_1} \(\theta_1^{(2)} - \theta_1^{(3)}\) \(\theta_1^{(3)} - \theta_1^{(1)}\) \(\theta_1^{(1)} - \theta_1^{(2)}\) \(\theta_1^{(2)} - \theta_1^{(3)}; \theta_1^{(3)} - \theta_1^{(1)}; \theta_1^{(1)} - \theta_1^{(2)}\) \\
  &= \pi^{4b_1} \(\theta_1^{(Q)} - \iota_Q\(\theta_1^{(Q)}\)\) N_{Q/K}\(\theta_1^{(Q)}\) \(\theta_1^{(Q)} - \iota_Q\(\theta_1^{(Q)}\); \iota_Q(\theta_1^{(Q)}); -\theta_1^{(Q)}\) \\
  &= \pi^{4b_1} \cdot N_{Q/K}\(\theta_1^{(Q)}\) \cdot c_1\sqrt{D_0} \(c_1\sqrt{D_0}; \iota_Q(\theta_1^{(Q)}); -\theta_1^{(Q)}\),
\end{align*}
which is seen to have the same class in $H^1$ as
\[
  \(1; c_1 \sqrt{D_0} \iota_Q\(\theta_1^{(Q)}\); -c_1 \sqrt{D_0} \theta_1^{(Q)}\).
\]
The last two coordinates are $\iota_Q$-conjugate, as they should be.
Multiplying by
\[
  \diamondsuit = \(1; \sqrt{D_0}; -\sqrt{D_0}\) \textand \heartsuit = \big(1; \tr(\pi_Q); \tr(\pi_Q)\big),
\]
we get
\[
  [\hat\omega_C \diamondsuit \heartsuit] = \big(c_1 \tr \pi_Q; \iota_Q(\theta_1); \theta_1\big).
\]
So
\begin{align*}
  \ell(\hat\omega_C \diamondsuit \heartsuit) &= \ell_R(b \tr \pi_Q; \bar\theta_1; \theta_1) \\
  &= \ell_Q\(c_1 \tr \pi_Q \cdot \theta_1\) \\
  &= \ell_Q\(c_1 \cdot \tr \pi_Q \cdot c_0 \cdot \(1 + \frac{c_1}{c_0} \pi_Q\)\).
\end{align*}
The last factor has exact level $s'$ if $s' < 2e$; otherwise it is a square. The remaining factors belong to $K$. By Lemma \ref{lem:iota}, their product is in $\L_e$ if $d_0 = 2e+1$; otherwise it is in either $\L_{d_0/2}$ or $\spadesuit \L_{d_0/2}$ according to the parity of
\[
  v\(c_1 \cdot \tr \pi_Q \cdot c_0\) = s' + \frac{d_0}{2}.
\]
It is now easy to work out the level possibilities for each value of $s'$.
\end{proof}

In cases \ref{type:A} and \ref{type:B}, as well as in the cases when $d_0 = 2e + 1$, there are no boxgroups worthy of discussion. Therefore in this subsection we will make the following assumptions:
\begin{itemize}
  \item $d_0 \leq 2e$ is even. We let $d_0 = 2d_0'$.
  \item $s' \geq d_0' - 1$, implying that $\hat\omega_C \diamondsuit\heartsuit \in \L_{d_0' - 1}$ and $\square_C \geq 2\floor{(d_0'-1)/2} + 1$.
\end{itemize}
First we need an $\eta$-lemma, as in the other splitting types:
\begin{lem}\label{lem:eta_1^21}
Let
\[
\square_C = \min \left\{ 2\floor{\frac{\ell(\hat\omega_C)}{2}} + 1, e \right\}
\]
be the squareness of the conic for $\delta = 1$. Assume we are in type \ref{type:C}, \ref{type:D}, or \ref{type:E} of Lemma \ref{lem:types_1^21}, and let
\[
  h_\eta = \begin{cases}
    1 & \text{ in type \ref{type:C}} \\
    0 & \text{ in types \ref{type:D} and \ref{type:E}}.
  \end{cases}
\]
Then, if $\theta_1$ is translated by a suitable element of $\OO_K$ and scaled by a suitable element of $\OO_K^\cross$ (neither of which change the associated resolvent $C$), there is an $\eta \in \OO_R$ satisfying the congruence
\begin{equation} \label{eq:eta_main}
  \eta^2 \equiv \pi^{h_\eta} \mod \pi^{m_\eta + h_\eta}
\end{equation}
where
\[
  m_\eta = 2\ceil{\frac{s'}{2}} + \square_C - d_0'.
\]
Also,
\begin{equation} \label{eq:eta_in_1^21}
  \eta \in \<1, \theta_1, \pi^{\ceil{\frac{s' - d_0'}{2}} + \frac{1}{2}}\theta_2, \pi^{e - d_0' + \frac{1}{2}}\theta_2\>.
\end{equation}
\end{lem}
\begin{proof}
In type \ref{type:E}, we simply scale $\theta_1$ so that
\[
  \theta_1 \equiv (1;0) \mod \pi^{s'}
\]
and take $\eta = \theta_1$. This achieves the desired modulus $m_\eta = s'$, and \eqref{eq:eta_in_1^21} is also trivially satisfied. We may therefore assume that $\theta_1$ is of type \ref{type:C} or \ref{type:D}.

Similar to the other splitting types, we view $\coef_{\theta_2} (\eta^2)$ as being a form of the conic $\M$ and carry out the transformation of the conic, as detailed in Lemmas \ref{lem:gamma_white} and \ref{lem:tfm_conic}:
\begin{align*}
  \frac{1}{\sqrt{\pi}}\coef_{\theta_2}(\eta^2) &\sim \pi^{- d_0'} \coef_{\bar\theta_2}(\eta^2) \\
  &\sim \pi^{- d_0'} \tr(\omega_C \eta^2) \\
  &= \pi^{- d_0'} \lambda^\diamondsuit(\diamondsuit \omega_C \eta^2) \\
  &= \pi^{\bar s - d_0'} \lambda^\diamondsuit(\diamondsuit \bar\omega_C \eta^2) \\
  &= \pi^{s'} \lambda^\diamondsuit(\diamondsuit \bar\omega_C \eta^2) \\
  &= \pi^{- \frac{h_1}{2}} \lambda^\diamondsuit(1^\odot {\eta^\odot}^2)\\
  &= \pi^{- \frac{h_1}{2}} \M(\eta^{\odot}).
\end{align*}
Here $1^\odot \in\OO_R^\cross$ is the coefficient arising in Lemma \ref{lem:tfm_conic}, with
\[
  \frac{\bar\omega_C \diamondsuit}{1^\odot} \in (R^\cross)^2,
\]
and
\[
  \eta^{\odot} = \pi^{\ceil{s'/2}} \sqrt{\frac{\bar\omega_C \diamondsuit}{1^\odot}} \cdot \eta.
\]
Thus points $\eta^{\odot}$ with $\M(\eta^{\odot}) \approx 0$ yield values $\eta$ such that $\coef_{\theta_2}(\eta^2) \approx 0$. The conic has a notable basepoint
\[
  \eta^{\odot}_0 = \pi^{\ceil{s'/2}} \sqrt{\frac{\bar\omega_C \diamondsuit}{1^\odot}} \sim \(\pi^{\ceil{s'/2}}; \pi_Q^{h_1/2}\),
\]
which corresponds to $\eta = 1$. It is evident that this is \emph{not} the $\eta$ we seek. However, it allows us to write the coordinate change succinctly as
\[
  \eta = \frac{\eta^{\odot}}{\eta^{\odot}_0}.
\]

By definition of $\square_C$, there is a linear form $\lambda$ such that
\[
\lambda^\diamondsuit(1^\odot{\eta^\odot}^2) \equiv c \cdot \lambda(\eta^{\odot})^2 \mod \pi^{\square_C}
\]
as functions of $\eta^{\odot} \in \OO_R$. Since $\M(\eta^{\odot}_0) = 0$, there is no harm in picking $\lambda$ such that $\lambda(\eta^{\odot}_0) = 0$ also. 
Our approach will be to complete to a basis $\{ \eta^{\odot}_0, \psi\}$ for $\ker \lambda$ and take $\eta^{\odot} = \pi^i \psi$, where $i$ is chosen such that $\eta = \eta^{\odot}/\eta^{\odot}_0$ is primitive in $\OO_R$.

In the proof of Squareness Lemma \ref{lem:squareness}, we rescaled $1^\odot$ by a scalar in $\OO_K^\cross$ and a square in $\(\OO_R^\cross\)^2$ so that the level $\ell(1^\odot \heartsuit)$ is manifest:
\[
1^\odot = \heartsuit (1; (1 + b\pi^{j+1}) + a\pi^j\pi_Q + 2\beta) = \(1; \pi_Q^{2d_0'} \cdot \frac{1+b\pi^{j+1} + a\pi^j\pi_Q + 2\beta}{t}\),
\]
where $j \in \{0, 1, 2,\ldots, e-1, \infty\}$ controls $\ell(1^\odot\heartsuit)$, and where $a \in \OO_K^\cross$, $b \in \OO_K$, $\beta \in \OO_R$.

We then looked at $\M(\xi^\odot)$ where $\xi^\odot$ ranges over the basis of $\OO_R$
\[
  (1;0), \quad (0; \alpha_0) = \(0;\pi_Q^{-d_0'}\pi^{\ceil{\frac{d_0'}{2}}}\), \quad (0; \alpha_1) = \(0;\pi_Q^{1-d_0'}\pi^{\floor{\frac{d_0'}{2}}}\)
\]
and got
\begin{align}
  \M(1;0) &= \lambda^\diamondsuit(1;0) \nonumber\\
  &= 1 \nonumber\\
  \M\(0;\alpha_0\) &= \pi^{2\ceil{d_0'/2}}\lambda^\diamondsuit\(0; \frac{1 + b\pi^{j+1} + a\pi^j\pi_Q + 2\beta}{t}\) \nonumber\\
  &\equiv a \cdot \frac{\pi^{2\ceil{\frac{d_0'}{2}}}}{t} \cdot \pi^j \mod 2 \label{eq:alpha0} \\
  \M\(0;\alpha_1\) &= \pi^{2\floor{d_0'/2}}\lambda^\diamondsuit\(0; \pi_Q^2 \cdot \frac{1 + b\pi^{j+1} + a\pi^j\pi_Q + 2\beta}{t}\) \nonumber \\
  &\equiv (1 + \pi^{j+1} b)\pi^{2\floor{\frac{d_0'}{2}}} - au \cdot \frac{\pi^{2\floor{\frac{d_0'}{2}} + 1}}{t} \cdot \pi^j + a\pi^{j + 2 \floor{\frac{d_0'}{2}}} \mod 2. \label{eq:bas_alpha1}
\end{align}
We may scale $\lambda$ so that $\lambda(1;0) = 1$. Then
\[
  \ker \lambda = \<\(-\lambda(0; \alpha_0); \alpha_0\), \(-\lambda(0; \alpha_1); \alpha_1\)\>;
\]
we choose one of these elements for $\psi$. The details are slightly different in each of types \ref{type:C}, \ref{type:D} and \ref{type:E}:
\begin{itemize}
  \item In type \ref{type:C}, we have $h_1/2 \equiv s' = d_0' - 1$ mod $2$, so
  \[
    \alpha_0 \sim \pi_Q^{1 - h_1/2} \textand \alpha_1 \sim \pi_Q^{h_1/2}.
  \]
  We can therefore take
  \[
    \psi = \(-\lambda(0; \alpha_0); \alpha_0\).
  \]
  Observe that 
  \begin{align*}
    v(\M(0; \alpha_0)) &\geq \min\left\{2\ceil{\frac{d_0'}{2}} - d_0' + j, e\right\} \\
    &= \min\left\{1 - \frac{h_1}{2} + j, e\right\} \\
    &\geq d_0' - \frac{h_1}{2} \\
    &= 2\floor{\frac{d_0' - 1}{2}} + 1.
  \end{align*}
  This lower bound is an odd integer, so we get
  \begin{align*}
    v(\lambda(0; \alpha_0)) &\geq \min\left\{\floor{\frac{d_0' - 1}{2}} + 1, \frac{\square_C + 1}{2}\right\} \\
    &= \floor{\frac{d_0' - 1}{2}} + 1 \\
    &= \floor{\frac{s'}{2}} + 1.
  \end{align*}
  This ensures that
  \[
  \eta^{\odot} = \pi^{h_1/2} \psi
  \]
  yields an $\eta = \eta^{\odot}/\eta^{\odot}_0$ with
  \[
    \pi | \eta_K \textand \eta_Q \sim \pi_Q.
  \]
  We get 
  \begin{alignat*}{2}
    \frac{1}{\sqrt{\pi}}\coef_{\theta_2}(\eta^2) &= \pi^{h_1/2} \M(\psi) \\
    &\equiv 0 &\mod &\pi^{h_1/2 + \square_C} \\
    &&={}&\pi^{ 2\ceil{s'/2} - s' + \square_C}\\
    &&={}&\pi^{m_\eta + h_\eta}.
  \end{alignat*}
  
  Hence $\eta^2$ is congruent modulo $\pi^{m_\eta + h_\eta + 1/2}$ to a value of the form $a + b\theta_1$, where $a,b \in \OO_K$. We have $v^{(K)}(\eta^2) \geq 2$, $v^{(Q)}(\eta^2) = 1$, so $\pi | a, \pi \sim b$. Thus $\eta^2 \equiv \pi\theta_1$ for suitably renormalized $\theta_1$. Note that \eqref{eq:eta_in_1^21} is trivial, as the right-hand side is all of $\OO_R$. This completes the proof in type \ref{type:C}.
  
  \item In type \ref{type:D}, we have $h_1/2 \equiv s' = d_0'$ mod $2$, so
  \[
  \alpha_0 \sim \pi_Q^{h_1/2} \textand \alpha_1 \sim \pi_Q^{1 - h_1/2}.
  \]
  We can therefore take
  \[
  \psi = \(-\lambda(0; \alpha_1); \alpha_1\).
  \]
  We have
  \[
  v\(\M(0;\alpha_1)\) = 2 \floor{\frac{d_0'}{2}}:
  \]
  the first term of \eqref{eq:bas_alpha1} dominates due to the bound on $j$ in type \ref{type:D}. Since $\square_C \geq 2\floor{d_0'/2} + 1$, we get
  \[
    v\(\lambda(0;\alpha)\) = \floor{\frac{d_0'}{2}}.
  \]
  Then we get a
  \[
    \psi = (-\lambda(\alpha_0); \alpha_0) \sim \(\pi^{\floor{d_0'/2}}; \pi_Q^{1-h_1/2}\)
  \]
  and take
  \[
    \eta^{\odot} = \pi^{\ceil{s'/2} - \floor{d_0'/2}} \psi
  \]
  for an $\eta = \eta^{\odot}/\eta^{\odot}_0$ whose $K$-component is a unit. As to the $Q$-component:
  \begin{align*}
    v(\eta^{(Q)}) &= \ceil{\frac{s'}{2}} - \floor{\frac{d_0'}{2}} + v(\psi^{(Q)}) \\
    &= \frac{2s' + h_1}{4} - \frac{2d_0' - h_1}{4} + \frac{2-h_1}{4}\\
    &\geq \frac{2 + h_1}{4}\\
    &\geq \frac{1}{2}.
  \end{align*}
  So $\eta$ is integral. We get
  \begin{alignat*}{2}
    \frac{1}{\sqrt{\pi}}\coef_{\theta_2}(\eta^2) &= \pi^{- \frac{h_1}{2} + 2\(\ceil{s'/2} - \floor{d_0'/2}\)} \M(\psi) \\
    &\equiv 0 &\mod &\pi^{- \frac{h_1}{2} + 2\(\ceil{s'/2} - \floor{d_0'/2}\) + \square_C} \\
    &&={}&\pi^{ 2\ceil{s'/2} - d_0' + \square_C}\\
    &&={}&\pi^{m_\eta + h_\eta}.
  \end{alignat*}
  Hence $\eta^2$ is congruent modulo $\pi^{m_\eta + h_\eta + 1/2}$ to a value of the form $a + b\theta_1$, where $a,b \in \OO_K$. We have $\pi \nmid b$ because the $K$-component of $\eta^2$ has zero valuation, the $Q$-component positive valuation. So $\eta^2 \equiv \theta_1 \mod \pi^{m_\eta + h_\eta + 1/2}$ for suitably renormalized $\theta_1$, as desired.
  
  To prove \eqref{eq:eta_in_1^21}, we observe that, since $\square_C > 0$,
  \[
    \(\eta^{(Q1)}\)^2 - \(\eta^{(Q2)}\)^2 \equiv \theta_1^{(Q1)} - \theta_2^{(Q1)} \equiv 0 \mod \pi^{\bar s}.
  \]
  Thus
  \[
    \eta^{(Q1)} - \eta^{(Q2)} \equiv 0 \mod \pi^{\min\{\bar s/2, e\}}
  \]
  (recalling that $\bar s = s' + d_0'$), that is,
  \begin{align*}
    v(\I(\eta^{(Q)})) &\geq \min\left\{\frac{\bar s}{2}, e\right\} - d_0' \\
    &= \min \left\{\frac{s' - d_0'}{2}, e - d_0'\right\}.
  \end{align*}
  Since the left-hand side is an integer, we can take the ceiling and get
  \[
    \eta \in \<1, (1;0), \pi^{\ceil{\frac{s' - d_0'}{2}}}(0; \pi_Q), \pi^{e - d_0'}(0; \pi_Q)\>.
  \]
  We can replace $(0; \pi_Q)$ by $\sqrt{\pi}\theta_2$, and since
  \[
    \ceil{\frac{s' - d_0'}{2}} \leq s',
  \]
  we can replace $(1;0)$ by $\theta_1$ to get the desired relation
  \[
    \eta \in \<1, \theta_1, \pi^{\ceil{\frac{s' - d_0'}{2}} + \frac{1}{2}}\theta_2, \pi^{e - d_0' + \frac{1}{2}}\theta_2\>. \qedhere
  \]
\end{itemize}
\end{proof}

We are now ready to construct boxgroups.
\begin{lem}\label{lem:boxes_basic_1^21} Assume that $b_1, s' \in \ZZ$, so $\theta_1 \in \OO_R$ is in one of types \ref{type:B}--\ref{type:E} of Lemma \ref{lem:types_1^21}. Let $ m_c, n_c \in \NN^+ \union \{\infty\} $, $ m_c \geq n_c > 0 $. Let $ B_{\theta_1}(m_c,n_c) $ be the box
  \[
  B_{\theta_1}(m_c,n_c) = \{\pi c_0 + \pi^{n_c}c_1\theta_1 + \pi^{m_c} c_2 \theta_2 : c_i \in \OO_K \}.
  \]
  \begin{enumerate}[$($a$)$]
    \item\label{boxes_basic_1^21:recenter} For every $ \xi \in 1 + B_{{\theta_1}(m_c,n_c)} $,
    \begin{align}
      B_{\theta_1}(m_c,n_c) &= \xi B_{\xi^{-1} {\theta_1}}(m_c,n_c) \label{eq:1^21_recenter_1} \\
      1 + B_{\theta_1}(m_c,n_c) &= \xi (1 + B_{\xi^{-1} {\theta_1}}(m_c,n_c)) \label{eq:1^21_recenter_2}.
    \end{align}
    \item\label{boxes_basic_1^21:group} If
    \begin{equation*}
      m_c \leq 2n_c + s',
    \end{equation*}
    then $ B_{\theta_1}(m_c,n_c) $ is closed under multiplication and the translate $ 1 + B_{\theta_1}(m_c,n_c) $ is a group under multiplication.
  \end{enumerate}
\end{lem}
\begin{rem}
  We will eventually apply this result with
\[
  m_c = m_{11}, \quad n_c = n_{11} - \bar s.
\]
\end{rem}
\begin{proof}
The proof proceeds just as in the unramified splitting types (Lemma \ref{lem:boxes_basic_ur}). For \ref{boxes_basic_1^21:group}, the key is to note that, in the basis $\{1, \theta_1, \sqrt{\pi} \theta_2\}$ for $\OO_R$,
\[
  v\(\coef_{\sqrt{\pi}\theta_2}(\theta_1^2)\) = s',
\]
a useful reinterpretation of the (renormalized) idempotency index $s'$.
\end{proof}
\begin{lem}\label{lem:boxgps_1^21}
Assume that $d_0$ is even and $\theta_1$ is of one of types \ref{type:C}--\ref{type:E} in Lemma \ref{lem:types_1^21}. Let $0 \leq n_c \leq m_c \leq 2e - d_0'$ be integers such that
\begin{align}
m_c &\geq n_c + d_0' - 1 \label{eq:gray_domino} \\
m_c &\leq 2n_c + s' &&(\text{the gray-red inequality}) \\
m_c &\leq n_c + 2\ceil{\frac{s'}{2}} - d_0' + \square_C &&(\text{the gray-green inequality}) \\
m_c &\leq e + \frac{n_c + s' - d_0' + 1}{2} &&(\text{the gray-blue inequality})
\end{align}
Then the projection $[1 + B_{\theta_1}(m_c,n_c)]$ of $1 + B_{\theta_1}(m_c,n_c)$ onto $H^1$ is a subgroup of signature
\begin{align*}
  &\0 .\0^{d_0'}.(\0\0)^{\ell_0} (\0\*)^{\ell_1} (\*\*)^{\ell_2} .\*^{d_0'}. \* \quad (\text{type \ref{type:C}}) \\
  &\0 .\0^{d_0'}.(\0\0)^{\ell_0} (\*\0)^{\ell_1} (\*\*)^{\ell_2} .\*^{d_0'}. \*\quad (\text{types \ref{type:D}--\ref{type:E}})
\end{align*}
where
\begin{align*}
  \ell_0 &= \floor{\frac{n_c - h_\eta}{2}} \\
  \ell_1 &= \floor{\frac{m_c - d_0' + h_\eta}{2}} - \floor{\frac{n_c - h_\eta}{2}} \\
  \ell_2 &= e - \floor{\frac{m_c + d_0' + h_\eta}{2}}
\end{align*}
\end{lem}
\begin{proof}
The gray-red inequality ensures that 
\[
  T = [1 + B_{\theta_1}(m_c, n_c)] = \{[1 + c_1 \pi^{n_c} \theta_1 + c_2 \pi^{m_c} (0; \pi_R)]\}
\]
is a subgroup. As in the other splitting types, we will establish its signature by downward induction on $n_c$, fixing $m_c$.

The base case is when \eqref{eq:gray_domino} is an equality. Here $\ell_1 \in \{0,1,-1\}$ is such that the signature simplifies, and we wish to prove that
\[
  T = \L_{m_c}.
\]
The $\supseteq$ direction is straightforward:
\[
  T \supseteq \{[\alpha] : \alpha \equiv 1 \mod \pi^{m_c + 1/2}\} = \L_{m_c}.
\]
For the $\subseteq$ direction, we simply need to show that the other generators $1 + c_1\pi^{n_c} \theta_1$ do not add any new elements to the box. Since $s' \geq d_0' - 1$ (we are in types \ref{type:C}--\ref{type:E}), we can scale and translate $\theta_1$ so that
\[
  \theta_1 \equiv (1; 0; 0) \mod \pi^{d_0' - 1/2}.
\]
Then
\[
  [1 + c_1\pi^{n_c} \theta_1] \equiv [(1 + c_1\pi^{n_c} ; 1 ; 1)] = \iota(1 + c_1\pi^{n_c}) \mod \L_{m_c}.
\]
The element $1 + c_1\pi^{n_c}$ has $K$-level at least $\floor{n_c/2}$ and, by Lemma \ref{lem:iota}\ref{iota:custom}, $Q$-level at least $n_c + d_0' - 1$, getting the needed bound $\floor{n_c/2} \leq e - d_0'$ from the bound $m_c \leq 2e - d_0'$. This completes the proof of the base case.

For the induction step, assume that $[1 + B_{\theta_1}(m_c, n_c + 1)]$ has the correct signature. Note that decreasing $n_c$ by $1$ enlarges the boxgroup at most by a factor of $q$. There are two main cases.

\paragraph{Moving case.} The first case is when
\begin{itemize}
  \item $n_c$ is odd in types \ref{type:D} and \ref{type:E}
  \item $n_c$ is even in type \ref{type:C}.
\end{itemize}
Here we claim that the boxgroup actually does grow by a factor of $q$, and in particular that the new generators $1 + c_1\pi^{n_c}$, $\pi \nmid c_1$, are all of level $n_c + d_0' - 1$. The method depends slightly on the type.

In types \ref{type:D} and \ref{type:E}, we have $s' \geq d_0'$ and thus can normalize $\theta_1$ so that
\[
  \theta_1 \equiv (1;0;0) \mod \pi^{d_0' + 1/2}.
\]
Then, for $n_c$ odd,
\[
[1 + c_1\pi^{n_c}\theta_1] \equiv [(1 + c_1\pi^{n_c}; 1; 1)] = \iota(1 + c_1\pi^{n_c}) \mod \L_{n_c + d_0'}.
\]
The element $1 + c_1\pi^{n_c}$ has $K$-level $(n_c - 1)/2$. Plugging into Lemma \ref{lem:iota}, we use $n_c \leq 2e - 2d_0'$ (from $m_c \leq 2e - d_0'$ and the fact that we are beyond the base case) to nail down the $Q$-level
\[
  \ell\(\iota(1 + c_1\pi^{n_c})\) = 2\(\frac{n_c - 1}{2}\) + d_0' = n_c + d_0' - 1.
\]

In type \ref{type:C}, we have $s' = d_0'$. We may normalize so that
\[
  \theta_1 = \(1; \pi^{d_0' - 1}u\pi_Q\), \quad u \in \OO_K^\cross.
\]
Then, for $n_c$ even,
\begin{align*}
\left[1 + c_1\pi^{n_c}\right] &= \left[\(1 + c_1\pi^{n_c}; 1 + c_1\pi^{n_c + d_0' - 1}u\pi_Q)\)\right] \\
&= \iota(1 + c_1\pi^{n_c}) \cdot \left[\(1; 1 + c_1\pi^{n_c + d_0' - 1}u\pi_Q)\)\right].
\end{align*}
The second factor has the desired level $n_c + d_0' - 1$. As to the first factor, an element of $K$-level at least $n_c/2$ is fed to $\iota$, yielding $Q$-level at least $n_c + d_0'$, using the inequality $n_c \leq 2e - 2d_0'$ again.

As a consequence, the boxgroup grows by a factor of $q$, entirely in level $n_c + d_0' - 1$: so we get its signature by replacing the $\0$ in that slot by a $\*$.

\paragraph{Stationary case.} We now come to the second and subtler case of the induction step, when
\begin{itemize}
  \item $n_c$ is odd in types \ref{type:D} and \ref{type:E}
  \item $n_c$ is even in type \ref{type:C}.
\end{itemize}
Here our aim is to prove that the boxgroup is unchanged, that is, that the new generators $1 + c_1\pi^{n_c}\theta_1$ actually lie in the existing group $[1 + B_{\theta_1}(m_c, n_c + 1)]$. It suffices to show that there is a square in each of the $q$ cosets
\[
  (1 + c_1^2\pi^{n_c}\theta_1)\(1 + B_{\theta_1}(m_c, n_c + 1)\) = 1 + c_1^2\pi^{n_c}\theta_1 + B_{\theta_1}(m_c, n_c + 1), \quad c_1 \in \OO_K
\]
(only $c_1$ mod $\pi$ matters; the two forms of the coset are equal by the gray-red inequality). We produce these squares using the $\eta$-lemma (Lemma \ref{lem:eta_1^21}).

In types \ref{type:D} and \ref{type:E}, we have an $\eta \in \OO_R$ such that (renormalizing $\theta_1$ appropriately)
\[
  \eta^2 \equiv \theta_1 \mod \pi^{m_\eta + 1/2}.
\]
Given $n_c$ even, we claim that the square
\[
  \(1 + \pi^{n_c/2} c_1 \eta\)^2 = 1 + 2\pi^{n_c/2} c_1 \eta + \pi^{n_c} c_1^2 \eta^2
\]
lies in the proper coset, that is,
\[
  1 + 2\pi^{n_c/2} c_1 \eta + \pi^{n_c} c_1^2 \eta^2 \equiv 1 + \pi^{n_c} c_1^2\theta_1 \mod B_{\theta_1}(m_c, n_c + 1).
\]
This requires two considerations:
\begin{itemize}
  \item To replace $\eta^2$ by $\theta_1$, we need
  \[
    m_\eta + n_c \geq m_c,
  \]
  which is equivalent to the gray-green inequality.
  \item For the cross term, we use
  \[
    \eta \in B_{\theta_1}\(\ceil{\frac{s' - d_0'}{2}}, 0\)
  \]
  from Lemma \ref{lem:eta_1^21} to get
  \[
    2\pi^{n_c/2} c_1 \eta \in B_{\theta_1}\(e + \frac{n_c}{2} + \ceil{\frac{s' - d_0'}{2}}, e + \frac{n_c}{2}\) \subseteq B_{\theta_1}(m_c, n_c + 1)
  \]
  by the gray-blue inequality, the difference of whose sides is in $\ZZ + 1/2$ in type \ref{type:D}. (In type \ref{type:E}, since $\eta = \theta_1$, this step is even easier.)
\end{itemize}

Similarly, in type \ref{type:C}, we have an $\eta \in \OO_R$ such that (renormalizing $\theta_1$ appropriately)
\[
\eta^2 \equiv \pi\theta_1 \mod \pi^{m_\eta + 1 + 1/2}.
\]
Given $n_c$ odd, we claim that the square
\[
  \(1 + \pi^{\frac{n_c - 1}{2}} c_1 \eta\)^2 = 1 + 2\pi^{\frac{n_c - 1}{2}} c_1 \eta + \pi^{n_c - 1} c_1^2 \eta^2
\]
lies in the proper coset, that is,
\[
  1 + 2\pi^{\frac{n_c - 1}{2}} c_1 \eta + \pi^{n_c - 1} c_1^2 \eta^2 \equiv 1 + \pi^{n_c} c_1^2 \theta_1 \mod B_{\theta_1}(m_c, n_c + 1).
\]
This requires two considerations:
\begin{itemize}
  \item To replace $\eta^2$ by $\pi \theta_1$, we need
  \[
    m_\eta + n_c \geq m_c
  \]
  (the term $1 = h_\eta$ cancels), which is equivalent to the gray-green inequality.
  \item For the cross term, since $\eta \equiv 0 \mod \sqrt{\pi}$, we have
  \[
    2\pi^{(n_c-1)/2} c_1 \eta \equiv 0 \mod \pi^{e + \frac{n_c - 1}{2} + \frac{1}{2}},
  \]
  the exponent being at least $m_c + 1/2$ by the gray-blue inequality (the difference of whose sides is in $\ZZ + 1/2$). \qedhere
\end{itemize}
\end{proof}

\begin{lem}\label{lem:boxgpS_1^21}
  For every triple $ (\ell_0, \ell_1, \ell_2) $ of integers satisfying
  \begin{align}
    \ell_0 + \ell_1 + \ell_2 &= e - d_0' \\
    \ell_1 &\leq \ell_0 + \ceil{\frac{s' - d_0'}{2}} + 1 + h_\eta && (\text{the gray-red inequality})\label{eq:boxgp_red_1^21} \\
    \ell_1 &\leq \frac{\square_C + 1}{2} + \ceil{\frac{s'}{2}} - d_0' + h_\eta && (\text{the gray-green inequality}) \label{eq:boxgp_green_1^21} \\ 
    \ell_1 &\leq \ell_2 + \ceil{\frac{s' - d_0'}{2}} + 1 + h_\eta, && (\text{the gray-blue inequality})\label{eq:boxgp_blue_1^21} \\
    \ell_0 &\geq 0, \quad \ell_1 \geq -1 + h_\eta, \quad \ell_2 \geq 0, \nonumber
  \end{align}
  there is a \emph{boxgroup} $T(\ell_0,\ell_1,\ell_2) \subseteq H^1$ of signature 
  \begin{align*}
    &\0 .\0^{d_0'}.(\0\0)^{\ell_0} (\0\*)^{\ell_1} (\*\*)^{\ell_2} .\*^{d_0'}. \* \quad (\text{type \ref{type:C}}) \\
    &\0 .\0^{d_0'}.(\0\0)^{\ell_0} (\*\0)^{\ell_1} (\*\*)^{\ell_2} .\*^{d_0'}. \*\quad (\text{types \ref{type:D}--\ref{type:E}})
  \end{align*}
  such that, if $m_c$, $n_c$ are integers satisfying the conditions of Lemma \ref{lem:boxgps_1^21}, then
  \[
  [1 + B_{\theta_1}(m,n)] = T\( \floor{\frac{n_c - h_\eta}{2}}, \floor{\frac{m_c - d_0' + h_\eta}{2}} - \floor{\frac{n_c - h_\eta}{2}}, e - \floor{\frac{m_c + d_0' + h_\eta}{2}}
   \).
  \]
\end{lem}
\begin{proof}
  If $\ell_1 = -1 + h_\eta$, take the appropriate level space, the unique subgroup with the correct signature.
  
  Otherwise, let
  \[
    m_c = 2e - 2\ell_2 - d_0' - h_\eta, \quad n_c = 2\ell_0 + h_\eta + 1
  \]
  in the preceding lemma. When we transform the conditions on $m_c$ and $n_c$ to conditions on the $\ell_i$, they become a priori
  \begin{align*}
    \ell_1 &\leq \ell_0 + \(\frac{s' - d_0' + h_\eta}{2}\) + 1 + h_\eta && (\text{the gray-red inequality}) \\
    \ell_1 &\leq \frac{\square_C + 1}{2} + \ceil{\frac{s'}{2}} - d_0' + h_\eta && (\text{the gray-green inequality}) \\
    \ell_1 &\leq \ell_2 + \(\frac{s' - d_0' + h_\eta}{2}\) + 1 + h_\eta, && (\text{the gray-blue inequality}) \\
  \end{align*}
  
  The gray-red and gray-blue inequalities are then put into the desired form using the identity
  \[
    \floor{\frac{s' - d_0' + h_\eta}{2}} = \ceil{\frac{s' - d_0'}{2}},
  \]
  which is easily proved using the interaction between $h_\eta$ and the parity of $s' - d_0'$ in each of types \ref{type:C}, \ref{type:D}, and \ref{type:E}.
  
  If the $\ell_i$ are fixed, let $(m_0, n_0)$ be the pair
  \[
    m_0 = 2e - 2\ell_2 - d_0' - h_\eta, \quad n_0 = 2\ell_0 + h_\eta + 1
  \]
  used in the proof of the lemma.
  
  The other pairs $(m_c, n_c)$ yielding the same triple $(\ell_0, \ell_1, \ell_2)$ are
  \[
    (m_c, n_c) = (m_0 + 1, n_0), \quad (m_0, n_0 - 1), \textor (m_0 + 1, n_0 - 1).
  \]
  
  Except for the base-case inequality, which is satisfied as long as $\ell_1 > -1 + h_\eta$, the conditions on $m_c$ and $n_c$ only become truer upon decreasing $m_c$ or increasing $n_c$. Now
  \[
  [1 + B_{\theta_1}(m_c,n_c)] \supseteq \left[1 + B_{\theta_1}\(m_0, n_c\)\right],
  \]
  but both sides have the same signature, so equality holds. Likewise,
  \[
  [1 + B_{\theta_1}(m_0,n_c)] \subseteq \left[1 + B_{\theta_1}\(m_0, n_0\)\right],
  \]
  but both sides have the same signature, so equality holds.
\end{proof}

\subsubsection{Supplementary boxgroups}
The boxgroups constructed thus far fit in the region
\[
  \L_{d_0/2} \supseteq T \supseteq \L_{e - d_0/2}.
\]
Outside this region, there is a less rich variety of groups appearing in the ring totals. They can be constructed using the following:
\begin{lem} \label{lem:iota_image}
Consider the map 
\begin{align*}
  \iota\colon K^\cross/(K^\cross)^2 &\to H^1 \\
  a &\mapsto [(a;1)] = [(1;a)].
\end{align*}
\begin{enumerate}
  \item If $d_0 = 2e+1$, then $\iota(K^\cross) = \L_e$.
  \item If $d_0 = 2d_0'$ is even, then the image $\iota(K^\cross)$ is a group of signature
  \[
  \0.\0^{d_0'-1}\spadesuit.(\*\0)^{e - d_0'}.\spadesuit^{\perp}\*^{d_0'-1}.\*
  \]
  where $\spadesuit \subseteq \L_{d_0' - 1}/\L_{d_0'}$ is the $2$-element subgroup generated by $\iota(\pi)$, and $\spadesuit^{\perp} \subseteq \L_{e - d_0'}/\L_{e - d_0' + 1}$ is its orthogonal complement.
\end{enumerate}
\end{lem}
\begin{proof}
Use Lemma \ref{lem:iota} to determine the size of the signature component at each level. Since $\iota(K^\cross)$ is self-orthogonal, we get the term $\spadesuit^{\perp}$ at level $e - d_0'$. 
\end{proof}

Now assume that $d_0 = 2d_0'$ is even. The group $\iota(K^\cross)$ is always important, but it does not behave well with respect to boxgroups unless $s' > 2e - d_0$, in which case we give it the name $T(\emptyset, e - d_0', \emptyset)$.

If $0 \leq \ell_1 \leq e - d_0'$ and $s' \geq d_0' + 2\ell_1$ (so we are in type \ref{type:D} or \ref{type:E}), we let
\begin{align*}
  T(e - d_0' - \ell_1, \ell_1, \emptyset) &= \iota(K^\cross) \intsec \L_{2e - d_0' - 2\ell_1} = \iota(1 + \pi^{2(e - d_0' - \ell_1) + 1}\OO_K) \\
  T(\emptyset, \ell_1, e - d_0' - \ell_1) &= \iota(K^\cross) \cdot \L_{d_0' + 2\ell_1}.
\end{align*}
(In the first case, the equivalence of the two definitions is established using Lemma \ref{lem:iota}.) Their signatures are, respectively,
\begin{gather*}
  \0.\0^{d_0'}.(\0\0)^{e - d_0' - \ell_1}(\*\0)^{\ell_1}.\spadesuit^\perp\*^{d_0' - 1}.\*
  \intertext{and}
  \0.\0^{d_0' - 1}\spadesuit . (\*\0)^{\ell_1} (\*\*)^{e - d_0' - \ell_1} . \*^{d_0'}.\*.
\end{gather*}

The restrictions on $s'$ ensure that these boxgroups satisfy such natural relations as
\begin{align*}
  T(e - d_0' - \ell_1, \ell_1, \emptyset) \cdot \L_{2e - d_0'} &= T(e - \ell_1, \ell_1, 0) \\
  T(\emptyset, \ell_1, e - d_0' - \ell_1) \intsec \L_{d_0'} &= T(0, \ell_1, e - d_0' - \ell_1),
\end{align*}
which we will often use without comment.

It will turn out that $T(\ell_0,\ell_1,\ell_2)$ and $T(\ell_2,\ell_1,\ell_0)$ are orthogonal complements whenever both are defined. Actually, this is simple to prove in the case that one of $\ell_0,\ell_1,\ell_2$ is the symbol $\emptyset$. The sizes of these groups follow immediately from their signatures:
\begin{lem}\label{lem:1^21_T_size}
  If $T(\ell_0,\ell_1,\ell_2)$ is defined and $\ell_1 \neq \emptyset$, then
  \[
  \size{T(\ell_0,\ell_1,\ell_2)} = \size{H^0}q^{e + \ell_2 - \ell_0},
  \]
  where if $\emptyset$ occurs as either $\ell_0$ or $\ell_2$, it must be replaced by $-1/[k_K : \FF_2] = \log_q(1/2)$.
\end{lem}

\end{wild}
\subsection{The recentering lemma}
When ${\theta_1}$ is nearly a square, we will sometimes be able to assume that it \emph{is} a square, thanks to the following lemma, which we state separately for each splitting type.
\begin{lem}\label{lem:recentering_ur}
  If $C$ is unramified and $s$ is even, then there is an element $ \psi_C \in 1 + B_{\theta_1}(\infty, \square_C) $ and $a \in \OO_K$, $b \in \OO_K^\cross$ such that $ \theta_1' = a + b\psi_C^{-1}\theta_1 = \eta^2 $ is a square in $\OO_R$.
\end{lem}
\begin{proof}
  It is easy to see that tweaking $\theta_1$ by addends in $\OO_K$ or multipliers in $\OO_K^\cross$, which do not change the underlying resolvent ring $C$, do not affect the truth of the lemma either. We may therefore assume, by Lemma \ref{lem:eta_ur}, that there is an $\eta$ such that
  \[
  \eta^2 \equiv \theta_1 \mod \pi^{s + \square_C}.
  \] 
  In particular, $ \eta^2 \equiv \theta_1 $ mod $ \pi^{2\ell + 1} \OO_K[\theta_1] $. Let $ k = \square_C $. Notice that $ \theta_1' $ can take all values in the orbit
  \[
  \theta_1' = \frac{a_{11} \theta_1 + a_{12}}{a_{21} \theta_1 + a_{22}} = \frac{a_{12}}{a_{22}} + \frac{(a_{11}a_{22} - a_{12}a_{21})\theta_1}{a_{22}(a_{21}\theta_1 + a_{22})}
  \]
  of $ \theta_1 $ under the congruence subgroup
  \[
  \Gamma(k) = \left\{\gamma =
  \begin{bmatrix}
  a_{11} & a_{12} \\ a_{21} & a_{22}
  \end{bmatrix} \in \GL_2(\OO_K): \gamma \equiv
  \begin{bmatrix}
  1 & 0 \\ 0 & 1
  \end{bmatrix}
  \mod \pi^k \right\}.
  \]
  We claim that this orbit is precisely the pixel $ \theta_1 + \pi^k \OO_K[\theta_1] $. The orbit is clearly contained in this pixel and contains all elements of the form
  \[
  \gamma(\theta_1) = \frac{(1 + a_{11}'\pi^k) \theta_1 + a_{12}'\pi^k}{a_{21}' \pi^k \theta_1 + (1 + a_{22}'\pi^k)} \equiv \theta_1 + \pi^k(a_{12}' + (a_{11}' - a_{22}')\theta_1 + a_{21}'\theta_1^2) \mod \pi^{2 k} \OO_K[\theta_1].
  \]
  So at least the orbit contains a point in each congruence class mod $ \pi^{2 k} \OO_K[\theta_1] $ in the claimed pixel. But applying general elements of $ \Gamma(2 k) $ to each of those, we get a point in each congruence class mod $ \pi^{4 k} \OO_K[\theta_1] $, and so on. Hence, the orbit is dense in the pixel, and being compact, it coincides with the pixel, establishing the desired result.
\end{proof}
\begin{lem}\label{lem:recentering_1^3}
  If $R$ is of splitting type $1^3$, there is an element $ \psi_C \in 1 + B_{\theta_1}(\infty, \square_C - 2h/3) $ and $a \in \pi^{-h/3} K \intsec \OO_{\bar K}$, $b \in \OO_K^\cross$ such that $ \theta_1' = a + b\psi_C^{-1}\theta_1 = \eta^2 $ is a square, $\eta \in \bar{\zeta}_3^{-h} \OO_R$.
\end{lem}
\begin{proof}
  It is easy to see that tweaking $\theta_1$ by addends in $\OO_K$ or multipliers in $\OO_K^\cross$, which do not change the underlying resolvent ring $C$, do not affect the truth of the lemma either. We may therefore assume, by Lemma \ref{lem:eta_1^3}, that there is an $\eta$ such that
  \[
  \eta^2 \equiv \theta_1 \mod \pi^{k}
  \]
  where $k = \square_C - 2h/3$. Notice that $ \theta_1' $ can take all values in the orbit
  \[
  \theta_1' = \frac{a_{11} \theta_1 + a_{12}}{a_{21} \theta_1 + a_{22}} = \frac{a_{12}}{a_{22}} + \frac{(a_{11} a_{22} - a_{12} a_{21})\theta_1}{a_{22}(a_{21}\theta_1 + a_{22})}
  \]
  of $ \theta_1 $ under the congruence subgroup
  \[
  \Gamma(k) = \left\{\gamma =
  \begin{bmatrix}
  a_{11} & a_{12} \\ a_{21} & a_{22}
  \end{bmatrix} \in \GL_2(\OO_{\bar K}): a_{ij} \in \pi^{\frac{h}{3}(i - j)} K, \gamma \equiv
  \begin{bmatrix}
  1 & 0 \\ 0 & 1
  \end{bmatrix}
  \mod \pi^k \right\}.
  \]
  We claim that this orbit is precisely the pixel
  \[
  \{ \theta_1' \in \bar{\zeta}_3^h \OO_R^\cross : \theta_1' \equiv \theta_1 \mod \pi^k \}.
  \]
  The orbit is clearly contained in this pixel and contains all elements of the form
  \[
  \gamma(\theta_1) = \frac{(1 + a'_{11}\pi^k) \theta_1 + a'_{12}\pi^k}{a'_{21} \pi^k \theta_1 + (1 + a'_{22}\pi^k)} \equiv \theta_1 + \pi^k(a'_{12} + (a'_{11} - a'_{22})\theta_1 + a'_{21}\theta_1^2) \mod \pi^{2 k}
  \]
  where the $a'_{ij}$ are integral in the appropriate groups $\pi^{\frac{h}{3}(i-j) - k} K$.
  
  So at least the orbit contains a point in each congruence class mod $ \pi^{2 k} \OO_K[\theta_1] $ in the claimed pixel. But applying general elements of $ \Gamma(2 k) $ to each of those, we get a point in each congruence class mod $ \pi^{4 k} \OO_K[\theta_1] $, and so on. Hence, the orbit is dense in the pixel, and being compact, it coincides with the pixel, establishing the desired result.
\end{proof}
\begin{wild}
  \begin{lem}\label{lem:recentering_1^21}
  If $C$ has splitting type $1^21$ and type \ref{type:C} or \ref{type:D}, then there is an element $ \psi_C \in 1 + B_{\theta_1}(\infty, \square_C - d_0' + h_1/2) $ and $a \in \OO_K$, $b \in \OO_K^\cross$ such that $ \theta_1' = a + b\psi_C^{-1}\theta_1 $ has the property that $\pi^{h_\eta} \theta_1' = \eta^2 $ is a square in $\OO_R$.
\end{lem}
\begin{proof}
  It is easy to see that tweaking $\theta_1$ by addends in $\OO_K$ or multipliers in $\OO_K^\cross$, which do not change the underlying resolvent ring $C$, do not affect the truth of the lemma either. We may therefore assume, by Lemma \ref{lem:eta_ur}, that there is an $\eta$ such that
  \[
  \frac{\eta^2}{\pi^{h_\eta}} \equiv \theta_1 \mod \pi^{s' + \frac{1}{2} + k}
  \]
  where $k = m_\eta - s' = \square_C - d_0' + h_1/2$.
  In particular, $ \eta^2 \equiv \theta_1 $ mod $ \pi^{k} \OO_K[\theta_1] $. Notice that $ \theta_1' $ can take all values in the orbit
  \[
  \theta_1' = \frac{a_{11} \theta_1 + a_{12}}{a_{21} \theta_1 + a_{22}} = \frac{a_{12}}{a_{22}} + \frac{(a_{11}a_{22} - a_{12}a_{21})\theta_1}{a_{22}(a_{21}\theta_1 + a_{22})}
  \]
  of $ \theta_1 $ under the congruence subgroup
  \[
  \Gamma(k) = \left\{\gamma =
  \begin{bmatrix}
    a_{11} & a_{12} \\ a_{21} & a_{22}
  \end{bmatrix} \in \GL_2(\OO_K): \gamma \equiv
  \begin{bmatrix}
    1 & 0 \\ 0 & 1
  \end{bmatrix}
  \mod \pi^k \right\}.
  \]
  We claim that this orbit is precisely the pixel $ \theta_1 + \pi^k \OO_K[\theta_1] $. The orbit is clearly contained in this pixel and contains all elements of the form
  \[
  \gamma(\theta_1) = \frac{(1 + a_{11}'\pi^k) \theta_1 + a_{12}'\pi^k}{a_{21}' \pi^k \theta_1 + (1 + a_{22}'\pi^k)} \equiv \theta_1 + \pi^k(a_{12}' + (a_{11}' - a_{22}')\theta_1 + a_{21}'\theta_1^2) \mod \pi^{2 k} \OO_K[\theta_1].
  \]
  So at least the orbit contains a point in each congruence class mod $ \pi^{2 k} \OO_K[\theta_1] $ in the claimed pixel. But applying general elements of $ \Gamma(2 k) $ to each of those, we get a point in each congruence class mod $ \pi^{4 k} \OO_K[\theta_1] $, and so on. Hence, the orbit is dense in the pixel, and being compact, it coincides with the pixel, establishing the desired result.
\end{proof}

\end{wild}
\subsection{Charmed cosets}
\begin{lem} \label{lem:charm}
Let $H$ be a finite $2$-torsion group, and let $\epsilon : H \to \{\pm 1\}$ be a nondegenerate quadratic form over $\FF_2$. Let $V \subseteq H$ be a subspace that is \emph{coisotropic;} that is, $V^\perp$ is isotropic, or equivalently, $V$ contains a maximal isotropic subspace. Then:
\begin{enumerate}[$($a$)$]
  \item \label{charm:exist} There is exactly one coset $\alpha V \subseteq H$ such that
\[
  \sum_{\delta \in \alpha V} \epsilon(\alpha) \neq 0.
\]
Indeed, the sum is
\[
  \pm_{\epsilon} \sqrt{\size{H}},
\]
where $\pm_{\epsilon} \in \{\pm 1\}$ is an invariant of the quadratic space $(V, \epsilon)$. We call $\alpha V$ the \emph{charmed coset} of $V$.
\item \label{charm:constant} On any coset $\beta V^\perp$ inside the charmed coset $\alpha V$, $\epsilon$ is constant. By contrast, on any coset $\beta V^\perp$ outside the charmed coset, $\epsilon$ is equidistributed.
\end{enumerate}
\end{lem}
\begin{proof}
\begin{enumerate}[$($a$)$]
\item Assume first that $V$ is maximal isotropic. Then $\size{V} = \sqrt{\size{H}}$. (If an $\FF_2$-space admits a nondegenerate quadratic form, its dimension is even if finite.) On each coset $\alpha V$, $\epsilon$ looks like a linear form, that is, there is a $\lambda_\alpha \in V^*$ such that
\[
  \epsilon(\alpha \beta) = \epsilon(\alpha) \cdot \lambda_\alpha(\beta).
\]
Note that the linear form $\lambda_\alpha$ is independent of coset representative, so we have a mapping
\[
  \lambda : H/V \to V^*.
\]
If $\lambda$ takes the same value on two different cosets $\alpha_1 V, \alpha_2 V$, then we see that $\epsilon$ is linear on the union $\alpha_1 V \sqcup \alpha_2 V$. Then the associated bilinear form $\<,\>_\epsilon$ is isotropic on the space $V \oplus \{0, \alpha_2/\alpha_1\}$, which is too big to be isotropic. So $\lambda$ is injective. Comparing sizes, we see that $\lambda$ is surjective also. So there is one coset $\alpha V$ on which $\epsilon$ is identically $1$ or $-1$. This is the charmed coset. On the remaining cosets, the values of $\epsilon$ are those of a nontrivial linear functional on $V$ and hence are equidistributed between $1$ and $-1$.

For a general $V$, take $N \subseteq V$ maximal isotropic. Every $V$-coset decomposes into $N$-cosets, and only the one containing the charmed coset of $N$ will yield a nonzero sum for $\epsilon$, namely $\pm \size{V} = \pm \sqrt{\size{H}}$.

A priori the sign $\pm$ of the sum on the charmed coset depends on both $\epsilon$ and $V$. But if $N \subseteq V$ are coisotropic, then it is easy to see that $N$ and $V$ yield the same sign. Then, taking $V = H$, we obtain that one sign holds for all coisotropic subspaces.

\item Take $N \subseteq V$ maximal isotropic. Then the charmed coset $\alpha V$ of $V$ is the one containing a charmed coset $\alpha N$ of $N$ on which $\epsilon = \pm_{\epsilon} 1$ is constant.

If $\beta/\alpha \in V$, then for all $\gamma \in V^\perp$, we have $\<\alpha\beta, \gamma\> = 1$ so
\begin{align*}
  \epsilon(\beta \gamma) &= \epsilon(\beta) \epsilon(\alpha) \epsilon(\alpha \gamma) \<\alpha\beta, \gamma\> \\
  &= \epsilon(\beta),
\end{align*}
since $\alpha, \alpha \gamma \in \alpha N$. By contrast, if $\beta/\alpha \notin V$, there exists $\gamma \in V^\perp$ such that $\<\alpha\beta, \gamma\> = -1$. Then $\epsilon(\beta \gamma) = -\epsilon(\beta)$, so $\epsilon$ is nonconstant on $\beta V^\perp$. But $V^\perp$ is isotropic, so $\epsilon$ is a linear form (plus a constant) on $\beta V^\perp$ and is therefore equidistributed.
\end{enumerate}
\end{proof}

As you might expect, we apply this lemma to the space $H = H^1$ with its quadratic form $\epsilon$. We will eventually find that $\pm_{\epsilon} = 1$ (it is ``positively charmed,'' one might say) though this is not obvious.

Let $\epsilon_C$ be the following translation of $\epsilon$: for $[\alpha] \in H^1$,
\[
  \epsilon_C(\alpha) = \epsilon\(\hat\omega_C \alpha\) = \begin{cases}
    1, & \coef_{\bar\theta_2}(\alpha \xi^2) = 0 \text{ for some } \xi \in R^\cross \\
    -1, & \text{otherwise.}
  \end{cases}
\]
Note that $\epsilon_C$ is still a quadratic form on $H^1$ whose associated bilinear form is the Hilbert symbol $\<\bullet,\bullet\>_R$. Note that $\epsilon_C(1) = 1$. Lemma \ref{lem:Brauer_const} can be interpreted as saying that 
\[
  \heartsuit \diamondsuit \L_{\floor{e'/2}}
\]
is charmed for $\epsilon_C$.

If $V \subseteq H^1$ is a subspace, we let
\[
  F_V = \1_{V}
\]
be its characteristic function. If $V$ is coisotropic, we denote by $G_{\epsilon_a,V}$, resp{.} $G_{\epsilon_C,V}$ the characteristic function of its charmed coset with respect to one of the quadratic forms $\epsilon_a, \epsilon_C$ whose associated bilinear form is the Hilbert pairing. The $\epsilon_a$ will be omitted if clear. The following results will power the computation of Fourier transforms of ring totals, a necessary step in our desired reflection theorems.
\begin{lem} \label{lem:FT_charm}
  Let $V$ be a subspace of $H^1$.
  \begin{enumerate}[(a)]
    \item\label{FT:F}
    \[
      \widehat{F_V} = \frac{\size{V}}{\size{H^0}} F_{V^\perp}.
    \]
    \item\label{FT:G} If $V$ is coisotropic, then
    \[
      \widehat{G_{\epsilon_a,V}} = \frac{\size{V}}{\size{H^0}} \epsilon_a(1) \epsilon_a F_{V^\perp}.
    \]
    \item\label{FT:eG} If $V$ is coisotropic, then
    \[
      \widehat{\epsilon_a G_V} = \pm_{\epsilon_a} q^e \epsilon_a(1) \epsilon_a F_V.
    \]
  \end{enumerate}
\end{lem}
\begin{proof}
  Part \ref{FT:F} is a standard property of the Fourier transform. For parts \ref{FT:G} and \ref{FT:eG}, let $N \subseteq V$ be a maximal isotropic subspace, so
  \[
    V \supseteq N = N^\perp \supseteq V^\perp.
  \]
  Let $\beta N$ be the charmed coset of $N$, so $\beta V$ is the charmed coset of $V$.
  
  For \ref{FT:G}, we compute
  \begin{align*}
    \widehat{G_V}(\delta) &= \frac{1}{\size{H^0}} \sum_{\alpha \in V} \<\alpha \beta, \delta\> \\
    &= \frac{1}{\size{H^0}} \<\beta, \delta\> \sum_{\alpha \in V} \<\alpha, \delta\> \\
    &= \frac{\size{V}}{\size{H^0}} \<\beta, \delta\> F_{V^\perp} (\delta).
  \end{align*}
  Hence it remains to prove that, for $\delta \in V^\perp$,
  \[
    \<\beta, \delta\> = \epsilon_a(\delta).
  \]
  By the definition of the associated bilinear form,
  \[
    \<\beta, \delta\> = \epsilon_a(1) \epsilon_a(\delta) \epsilon_a(\beta) \epsilon_a(\delta \beta).
  \]
  But $\epsilon_a(\beta) = \epsilon_a(\delta \beta)$ since both arguments lie in the charmed isotropic coset $\beta N$. This establishes \ref{FT:G}.
  
  For \ref{FT:eG}, we compute
  \begin{align*}
    \widehat{\epsilon_a G_V} &= \frac{1}{\size{H^0}} \sum_{\alpha \in V} \epsilon_a(\alpha\beta) \<\alpha\beta, \delta\> \\
    &= \frac{1}{\size{H^0}} \sum_{\alpha \in V} \epsilon_a(1) \epsilon_a(\delta) \epsilon_a(\alpha\beta\delta) \\
    &= \frac{\epsilon_a(1)\epsilon_a(\delta)}{\size{H^0}} \sum_{\alpha \in V} \epsilon_a(\alpha\beta\delta).
  \end{align*}
  The last sum equals $\pm_{\epsilon_a} \sqrt{\size{H^1}}$ if $\beta\delta V$ is charmed, $0$ otherwise. But $\beta V$ is charmed, so the relevant condition is that $\delta \in V$. So
  \begin{align*}
    \widehat{\epsilon_a G_V} &= \frac{\epsilon_a(1)\epsilon_a(\delta)}{\size{H^0}} \cdot \pm_{\epsilon_a} \sqrt{\size{H^1}} F_V(\delta) \\
    &= \pm_{\epsilon_a} q^e \epsilon_a(1)\epsilon_a(\delta) F_V(\delta),
  \end{align*}
  as desired.
\end{proof}

\subsection{The projectors}

On the space of complex- (or even rational-) valued functions on $H^1$, we can define certain projectors that divide up the work to be done. First look at the cosets of $\L_0$. Let 
\[
  I = \begin{cases}
    \L_0 \union (1;\pi;\pi)\L_0, & \bar s > 0 \\
    \L_0, & \bar s = 0
  \end{cases}
\]
be the union $\L_0 \union \alpha\L_0$ of up to two cosets, using the distinguished splitting $R \isom K \cross Q$ if $\bar s > 0$. Let $I_0$, $I_1$, $I_2$ be the restriction operators that restrict the support of a function to $\L_0$, $I \setminus \L_0$, and $H^1 \setminus I$, respectively. They are orthogonal idempotents ($I_1$ and/or $I_2$ may vanish). Let $J_i$ be the conjugate of $I_i$ under the Fourier transform. Each $J_i$ is convolution by a certain function supported on $\L_{e'}$; $J_0$ is none other than the \emph{smear operator} $\Sm_{e'}$ which will occur below. Since $\L_{e'} \subseteq \L_0$, each $I_i$ commutes with each $J_j$, so the $I_i J_j$ form a system of nine orthogonal idempotents. For orderliness of presentation, we transform all ring answers to a sum of terms each in the image of one idempotent. The Fourier transform interchanges the images of $I_i J_j$ and $I_j J_i$. In splitting type $(111)$, all nine idempotents are nonzero, although $I_1 J_2$ and $I_2 J_1$ will be found to annihilate every ring total. In the remaining splitting types, some of the idempotents vanish, and correspondingly some terms of our answers can be ignored.

\begin{defn}\label{defn:x's_ur} In unramified splitting types, we define the use of a symbol $x$ as follows, where the $\ell_i$ are such that the relevant boxgroups are well defined:
  \begin{itemize}
    \item $xF(0, \ell_1, \ell_2) = I_1(F(\emptyset, \ell_1,\ell_2))$, so that
    \[
      F(\emptyset, \ell_1,\ell_2) = F(0, \ell_1, \ell_2) + xF(0, \ell_1, \ell_2).
    \]
    \item $xxF(0, 0, e) = I_2(F(\emptyset, \emptyset, e))$, so that
    \[
      F(\emptyset, \emptyset, e) = F(0, 0, e) + xF(0, 0, e) + xxF(0, 0, e).
    \]
    \item $Fx(\ell_0, \ell_1, 0) = 2 \cdot J_1(F(\ell_0, \ell_1, \emptyset))$, so that, if there is a distinguished coarse coset,
    \[
      2 F(\ell_0, \ell_1, \emptyset) = F(\ell_0, \ell_1, 0) + Fx(\ell_0, \ell_1, 0).
    \]
    \item $Fxx(e, 0, 0) = \size{H^0} \cdot J_2(F(e, \emptyset, \emptyset))$, so that
    \[
      \size{H^0} F(e, \emptyset, \emptyset) = F(e, 0, 0) + Fx(e, 0, 0) + Fxx(e, 0, 0).
    \]
    \item $xFx(0, e, 0) = 2 \cdot I_1 J_1(F(\emptyset, e, \emptyset))$ (note there must be a distinguished coarse coset for this to be meaningful), so that
    \begin{align*}
      2 F(\emptyset, e, \emptyset)) &= 2 F(0, e, \emptyset) + 2 xF(0, e, \emptyset) \\
      &= F(\emptyset, e, 0) + Fx(\emptyset, e, 0) \\
      &= F(0, e, 0) + xF(0,e,0) + Fx(0,e,0) + xFx(0,e,0).
    \end{align*}
    \item The same definitions with $G$ replacing $F$, as appropriate. We find that $xG$, $xxG$, and $xGx$ are applicable.
  \end{itemize}
\end{defn}
Observe that these definitions are crafted so that the following cute rule applies:
\begin{lem}\label{lem:FT_x}
Let $N$ be one of the symbols $Fx, xF$, $Fxx$, $xxF$, $xFx$, and let $N'$ be the symbol made by spelling $N$ backward. If the $\ell_i$ are integers such that $N(\ell_0, \ell_1, \ell_2)$ is meaningful, then
\[
  \widehat{N(\ell_0, \ell_1, \ell_2)} = q^{\ell_2 - \ell_0} N'(\ell_2, \ell_1, \ell_0).
\]
(This lemma also holds for $N = F$, though this will be proved later: see Lemma \ref{lem:orth}.)
\end{lem}
This will allow us to write the ring totals for all three unramified splitting types in a uniform way and verify reflection for them simultaneously.

For the unshifted quadratic form $\epsilon$, $\L_i$ is charmed for $i \geq e'/2$, an easy consequence of  for $i \geq e'/2$, from which we derive:

\begin{lem} \label{lem:e_projectors}
  The function $\epsilon_C = \epsilon_C F(\emptyset, \emptyset, e)$ has at most the following projections nonzero:
  \begin{enumerate}[$($a$)$]
    \item If $[\hat{\omega}_C] \in \L_0$, then $I_0 J_0$, $I_1 J_1$, and $I_2 J_2$.
    \item If $[\hat{\omega}_C] \notin \L_0$, then $I_0 J_1$, $I_1 J_0$, and $I_2 J_2$.
  \end{enumerate}
\end{lem}
\begin{proof}
It suffices to prove that
\[
  I_i \(\epsilon_C F(\emptyset, \emptyset, e)\) = (x^i F)(0, 0, e)
\]
is in the image of $J_j$, for each pair $(i,j)$ mentioned in the lemma. Now $\L_0$ is charmed for the unshifted quadratic form $\epsilon$, an easy consequence of Proposition \ref{prop:conic_perturb}. Hence $\hat{\omega}_C \L_0$ is charmed for $\epsilon_C$, from which we get
\[
  (x^i F)(0, 0, e) = (x^j G)(0, 0, e).
\]
Taking the Fourier transform by Lemma \ref{lem:FT_charm}\ref{FT:eG},
\[
  \widehat{(x^i F)(0, 0, e)} = \pm_{\epsilon_C} (x^j F)(0, 0, e),
\]
which lies in the image of $I_j$. Hence the original $(x^i F)(0, 0, e)$ lies in the image of $J_j$, as desired.
\end{proof}

\subsection{Notation}

If $T(\ell_0, \ell_1, \ell_2)$ is a boxgroup, we let $F(T) = F(\ell_0, \ell_1, \ell_2)$ be its characteristic function, and $G(T) = G(\ell_0, \ell_1, \ell_2)$ be the characteristic function of its charmed coset, if applicable. Let $T^\cross(\ell_0, \ell_1, \ell_2)$ denote the subset of elements of $T(\ell_0, \ell_1, \ell_2)$ having minimal level, assuming this level is less than $e'$:
\[
  T^\cross(\ell_0, \ell_1, \ell_2) = \begin{cases}
    T(\ell_0, \ell_1, \ell_2) \setminus T(\ell_0 + 1, \ell_1 - 1, \ell_2), & \ell_0 \geq 0, \ell_1 \geq 1 \\
    T(\ell_0, \ell_1, \ell_2) \setminus T(\ell_0 + 1, 0, \ell_2 - 1), & \ell_0 \geq 0, \ell_1 = 0 \\
    T(\emptyset, \ell_1, \ell_2) \setminus T(0, \ell_1, \ell_2), & \ell_0 = \emptyset
  \end{cases}
\]
Let $F^\cross(\ell_0, \ell_1, \ell_2)$ and $G^\cross(\ell_0, \ell_1, \ell_2)(\delta)$ be the characteristic functions of $T^\cross(\ell_0, \ell_1, \ell_2)$ and $\psi_C T^\cross(\ell_0, \ell_1, \ell_2)$, respectively. This will provide enough notation to write the ring totals in the succeeding sections.

\section{Ring volumes for \texorpdfstring{$\xi'_1$}{xi'1}}
\label{sec:xi1}
In this section, we will compute the volume of vectors $\xi'_1 \in \OO_R$ satisfying the $\M_{11}$ and $\N_{11}$ conditions. A sample of our answers are tabulated in Appendix \ref{sec:zone_examples}.

Because $\bar\xi_1 = \xi_1$, we freely omit the bar on $m_{11} = m_{11}$ and $n_{11} = n_{11}$.

\subsection{The smearing lemma}
\label{sec:smear}
We now prove a simple lemma that allows us to reduce to the case $m_{11}$ large.
\begin{lem} \label{lem:smear}
  For a first vector problem $\P$, define $u_1$ to be the unique value such that, when the conic is transformed to minimal discriminant in accordance with \ref{lem:tfm_conic}, $\beta^\odot = \delta^\odot {\xi_1^\odot}^2$ must lie in $\pi^{u_1} \OO_R \setminus \pi^{u_1 + 1} \OO_R$, to wit:
  \[
  u_1 = \begin{cases}
    1 & \text{if } (\sigma, h_1) \in \{(1^3, -1)
\begin{wild}
      , (1^21, 2), (1^21, 3)
\end{wild}
    \} \\
    0 & \text{otherwise.}
  \end{cases}
  \]
  Assume that the $m_{11}$ and $n_{11}$ of $\P$ satisfy
  \[
  m_{11}^\odot > u_1
  \]
  (so $\M_{11}$ is active even after transformation) and $m_{11}^\odot > n_{11}$. Let $m_{11}^\sharp > m_{11}$ lie in the same class mod $\ZZ$, and let $m_{11}^{\odot\sharp}$ be the corresponding value of $m'_{11}$. Let $\P^\sharp$ be the first vector problem with $m_{11} = m_{11}^\sharp$ and the rest of the data the same.
  
  Then the answer to $\P$ can be computed from that of $\P^\sharp$ by the following formula:
  \[
  W_{\P} = q^{m_{11}^\sharp - m_{11}} \Sm_{r} W_{\P^\sharp}
  \]
  Here the \emph{smear operator} $\Sm_{r}$ is defined by the following convolution:
  \[
  \Sm_{r} W(\delta) = \frac{1}{\size{\L_r}} \sum_{\alpha \in \L_r} W(\alpha\delta),
  \]
  and $r$ is the level for which
  \[
  \L_r = \left\{ [\eta] : \eta \equiv 1 \mod \pi^{m_{11}^\odot - u_1} \right\}.
  \]
\end{lem}
We call this the \emph{smearing lemma} because it states that the function $W_{\P}$ can be obtained from $W_{\P^\sharp}$ by averaging over the cosets of $\L_{r}$, like reducing the resolution of a picture by averaging over larger pixels. Here the symbol $\sharp$ is used to mark the ``sharper'' image given by the solutions of $\P^\sharp$ and should not be confused with the use of the same symbol in the context of tilting.
\begin{todofn}
(See section \ref{sec:tilt}.)
\end{todofn}

The level $r$ is given explicitly as follows:
\begin{itemize}
  \item In unramified type,
  \[
  r = \floor{\frac{m_{11}}{2}}, \quad 0 < m_{11} \leq 2e
  \]
  \item In splitting type $1^3$,
  \[
  r = \floor{m_{11}}, \quad 0 < m_{11} \leq 2e.
  \]
\begin{wild}
    \item In splitting type $1^21$,
    \[
    r = m_{11} - \frac{d_0}{2} + \left\{\frac{h_1}{2}\right\},
    \quad \frac{d_0 + 1}{2} \leq m_{11} \leq 2e + \frac{d_0}{2}.
    \]
\end{wild}
\end{itemize}

\begin{proof}
  We will prove the identity by computing in two ways the volume of the set
  \begin{equation}
    \S = \left\{
    (\eta, \xi'):
    \begin{aligned}
      \eta &\equiv 1 \mod \pi^{m_{11}^\odot - u_1}, \\
      \lambda^\diamondsuit(\eta \delta^\odot \xi'^2) &\equiv 0 \mod \pi^{m_{11}^{\odot\sharp}}, \\
      \bar{\omega}_C^{-1} \gamma^2 \xi'^2 &\equiv a \mod \pi^{n_{11}} \text{ for some } a \in \OO_K
    \end{aligned}
    \right\}
    \subseteq \OO_R^\cross \cross \PP(\OO_R).
  \end{equation}
  First, fix $\eta \in 1 + \pi^{m_{11}^\odot - u_1}\OO_R$. The conditions on $\xi'$ are seen to be the $\M_{11}$ and $\N_{11}$ conditions for $\delta$ replaced by $\eta \delta$; the omission of the $\eta$ factor in $\N_{11}$ makes no difference, since the left-hand side is a multiple of $\pi^{u_1}$ and the addition would have valuation at least
  \[
  (m_{11}^\odot - u_1) + u_1 = m_{11}^\odot \geq n_{11}.
  \]
  So the volume of $\xi'$ for fixed $\eta$ is $W_{\P^\sharp} (\delta\eta)$, and since $[\eta]$ takes all classes in $\L_{r}$ equally often while ranging in a pixel of volume $q^{-3(m_{11}^\odot - u_1)}$, we get
  \[
  \mu(\S) = \frac{q^{-3(m_{11}^\odot - u_1)}}{\size{\L_r}} \sum_{\eta \in \L_{r}} W_{\P^\sharp}(\delta\eta).
  \]
  On the other hand, a fixed $\xi' \in \OO_R$ has a chance of being the second coordinate of a pair in $\S$ only if
  \begin{itemize}
    \item it satisfies the $\N_{11}$ condition $\bar{\omega}_C^{-1} \gamma^2 \xi'^2 \equiv a \mod \pi^{n_{11}}$ for some $a \in \OO_K$, and 
    \item it satisfies the $\M_{11}$ condition for some $\eta = 1 + \pi^{m_{11}^\odot - u_1} \eta'$, $\eta' \in \OO_R$; in particular,
    \begin{align*}
      0 &\equiv \lambda^\diamondsuit(\eta \delta^\odot \xi'^2) \\
      &= \lambda^\diamondsuit(\delta \xi'^2) + \pi^{m_{11}^\odot - u_1} \lambda^\diamondsuit(\eta' \delta \xi'^2) \\
      &= \lambda^\diamondsuit(\delta \xi'^2) + \pi^{m_{11}^\odot} \lambda^\diamondsuit\(\pi^{-u_1} \delta \xi'^2 \) \\
      &\equiv \lambda^\diamondsuit(\delta \xi'^2) \mod \pi^{m_{11}^\odot}.
    \end{align*}
  \end{itemize}
  The volume of $\xi'$ satisfying these conditions is, by definition, none other than $W_{\P}$. For fixed $\xi'$, the value of $\eta'$ is constrained by $\M_{11}$ alone:
  \begin{align*}
    \pi^{m_{11}^\odot - u_1} \lambda^\diamondsuit(\eta' \delta^\odot \xi'^2) &\equiv 0 \mod \pi^{m_{11}^{\odot\sharp}} \\
    \lambda^\diamondsuit\(\eta' \cdot \pi^{-u_1}\delta^\odot\xi'^2\) &\equiv 0 \mod \pi^{m_{11}^{\odot\sharp} - m_{11}^\odot} = \pi^{m_{11}^\sharp - m_{11}}.
  \end{align*}
  Since $\lambda^\diamondsuit$ is a perfect linear functional and $\pi^{-u_1}\delta^\odot\xi'^2$ is a primitive vector in $\OO_R$, the volume of $\eta'$ satisfying this congruence is $q^{m_{11} - m_{11}^\sharp}$, which makes a volume of $q^{-3(m_{11}^\odot - u_1) + m_{11} - m_{11}^\sharp}$ for $\eta$. So
  \[
  \mu(\S) = \pi^{-3(m_{11}^\odot - u_1) + m_{11} - m_{11}^\sharp} W_{\P}(\delta).
  \]
  Comparing the two expressions for $\mu(\S)$, the result follows.
\end{proof}

\subsection{The zones when \texorpdfstring{$\N_{11}$}{N11} is strongly active (black, plum, purple, blue, green, and red)}\label{sec:strong}
In this section, we solve first vector problems in which $\N_{11}$ is strongly active, that is, $n_{11} > \bar s$. In view of the smearing lemma, we assume that $m_{11} > 2e$. Let $n_c$ (``n for the colorful zones'') be $n_{11} - \bar s$.

The following little symmetry will be occasionally useful:
\begin{lem}\label{lem:recenter}
Let $\psi = a + b \theta_1 \in \OO_R$ with $a \in \OO_K^\cross$, $b \in \pi^{n_{11}} \OO_K^\cross$, and let $\P'$ be the first vector problem derived from $\P$ by replacing the pertinent extender vector $\theta_1$ by $\theta_1' = \psi^{-1} \theta_1$. Then:
\begin{enumerate}[$($a$)$]
  \item \label{recenter:w_C} $\hat\omega_{\P'} = \psi^{-1} \hat\omega_\P$.
  \item \label{recenter:e_C} $\epsilon_{\P'}(\delta) = \epsilon_\P(\psi^{-1} \delta)$.
  \item \label{recenter:W} $W_{\theta_1',m_{11},n_{11}}(\delta) = W_{\theta_1,m_{11},n_{11}}(\psi^{-1}\delta)$.
\end{enumerate}
\end{lem}
\begin{proof}
The left-hand side is the volume of $\xi'$ for which
\begin{align*}
  \delta \(\frac{\xi'}{\xi'_0}\)^2 &\in \OO_K^\cross + B_{\psi^{-1}\theta_1}(m_{11},n_{11}) \\
  &= \psi^{-1} \cdot \(\OO_K^\cross + B_{\theta_1}(m_{11},n_{11})\)
\end{align*}
by Lemma \ref{lem:boxes_basic_ur}\ref{boxes_basic_ur:recenter}, since $\psi \in \OO_K^\cross + B_{\theta_1}(m_{11},n_{11})$. So the condition on $\xi'$ can be written as
\[
  \delta\psi \(\frac{\xi'}{\xi'_0}\)^2 \in \OO_K^\cross + B_{\theta_1}(m_{11},n_{11}),
\]
of which the solution volume is seen to be $W_{\theta_1,m_{11},n_{11}}(\psi\delta)$.

\end{proof}

The following formula for the sum of the values of $W_{m_{11},n_{11}}$ will be essential:
\begin{lem} \label{lem:sum_strong}
If the values of $m_{11}$ and $n_{11}$ make $\N_{11}$ strongly active, then $W_{m_{11},n_{11}} = 0$ unless the chosen coarse coset is $\L_0$. In this case
\[
  \sum_{\delta \in \L_0} W_{m_{11},n_{11}}(\delta) = \size{H^0} q^{2e - m_{11} - n_{11} + \frac{\bar s}{2} + \frac{d_0}{2} + v(N(\gamma))}
\]
where $d_0 = v_K(\Disc_K R)$; equivalently,
\[
  \sum_{\delta \in \L_0} W^{\odot}_{m_{11},n_{11}}(\delta) = \size{H^0} q^{2e - m_{11} - n_{11} + \frac{\bar s}{2} + \frac{d_0}{2} + v(N(\gamma\gamma^{\odot}))}.
\]
\end{lem}
\begin{proof}
By Lemma \ref{lem:to_box}, the support of $W_{m_{11},n_{11}}$ consists of the classes $[\delta] = [\beta]$ in $H^1$ of elements $\beta$ of the box
\[
  \OO_K^\cross + B_{m_{11},n_c} = \{x + y \pi^{n_c} \theta_1 + z \pi^{m_{11}} \theta_2, \quad x \in \OO_K^\cross, \quad y,z \in \OO_K\}.
\]
Our method is to show that $W_{m_{11},n_{11}}(\delta)$ can be interpreted as the volume of $\beta$ of class $[\delta]$ in the box, up to a scalar. Then since each $\beta$ belongs to just one square-class, the sum is known.

Since $m_{11},n_c > 0$, $\beta$ must be a unit, explaining why $[\delta] \in \L_0$. For fixed $\delta \in \OO_R^\cross$, as $\xi'$ ranges over the solution set of its transformed conditions, $\beta$ ranges over the elements of its box of class $[\delta]$, up to scaling. The correspondence is given by a relation of the form
\begin{equation}
  \beta = \delta \(\frac{\xi'}{\xi'_0}\)^2,
\end{equation}
where
\[
  \xi'_0 = \frac{1}{\gamma}\sqrt{\delta\omega_C}
\]
is an element of $\OO_R$ constructed from the $\gamma$ of Lemma \ref{lem:gamma_white} whose valuations represent a lower bound on the valuations of $\xi'$.

Now we compare the projective volumes of $\beta$ and $\xi'$ satisfying the $\M_{11}$ and $\N_{11}$ conditions. In the sequence
\[
  \xi' \longmapsto \frac{\xi'}{\xi'_0} \longmapsto \(\frac{\xi'}{\xi'_0}\)^2 \longmapsto \delta \(\frac{\xi'}{\xi'_0}\)^2 = \beta,
\]
each member is a primitive vector in $\OO_R$, so we can speak of projective volumes.

Dividing by $\xi'_0$ is a one-to-one operation that scales both affine and projective volumes by
\[
  q^{v_K(N_{R/K}(\xi'_0))} = q^{\bar s/2 - v(N(\gamma))}.
\]
So there is a volume $q^{\bar s/2 - v(N(\gamma))} W_{m_{11},n_{11}}(\delta)$ of $\xi/\xi'_0$.

Squaring, on units, multiplies small projective volumes by $q^{-2e}$ (since it takes the $i$-pixel about $1$ to the $(i+e)$-pixel for $i > e$). But it is $\size{H^0}$-to-one since there are $\size{H^0}$-many square roots of $1$ in $R$, up to scaling by $\pm 1$. Since the resolvent conditions are invariant under multiplying $\xi'$ by a square root of $1$, the volume of $(\xi/\xi'_0)$ is
\[
  \frac{q^{-2e + \bar s/2 - v(N(\gamma))}}{\size{H^0}} W_{m_{11},n_{11}}(\delta).
\]

Lastly, $\delta$ is a unit, so multiplying by it does not change volumes. Hence
\[
  \frac{q^{-2e + \bar s/2 - v(N(\gamma))}}{\size{H^0}} \sum_{\delta \in \L_0} W_{m_{11},n_{11}}(\delta)
\]
is the volume of the box $\PP(\OO_K^\cross + B_{m_{11},n_c})$. Thus it suffices to prove that
\[
  \mu(\PP(\OO_K^\cross + B_{m_{11},n_c})) = q^{-m_{11}-n_{11}+\bar s + \frac{d_0}{2}}.
\]
Converting to affine volumes, with $\mu(R) = 1$,
\begin{align*}
  \mu(\PP(\OO_K^\cross + B_{m_{11},n_c})) &= \frac{1}{1 - \frac{1}{q}} \mu(\OO_K^\cross + B_{m_{11},n_c}) \\
  &= \mu(\OO_K + B_{m_{11},n_c}) \\
  &= q^{i}
\end{align*}
where $i$ is the integer such that
\begin{align*}
  \Lambda^3(\OO_K + B_{m_{11},n_c}) &= \pi^i \Lambda^3 \OO_R \\
  \Lambda^3(\OO_{\bar{K}}\<1, \pi^{n_c} \theta_1, \pi^m_{11} \theta_2\>) &= \pi^i \Lambda^3 \OO_R \\
  n_c + m_{11} &= i + \frac{d_0}{2}.
\end{align*}
So
\[
  \sum_{\delta \in \L_0} W_{m_{11},n_c}(\delta) = \size{H^0} \cdot q^{2e - \bar s/2 + v(N(\gamma))} \cdot q^{-m_{11} - n_c + \frac{d_0}{2}} = \size{H^0} \cdot q^{2e - m_{11} - n_{11} + \bar s/2 + d_0' + v(N(\gamma))}
\]
as desired.
\end{proof}

We tabulate:
\begin{equation} \label{tab:sum_strong}
  \begin{tabular}{rc|cl}
    spl.t. & $h$ & $\xi'_0 \sim$ & $v(N(\gamma\gamma^{\odot}))$ \\ \hline
    ur & $0$ & $(\pi^{s/2};1;1)$, $s$ even & $0$ \\
    ur & $1$ & $(\pi^{(s-1)/2} ; 1 ; 1)$, $s$ odd & $1/2$ \\
    $1^3$ & $1$ & $1$ & $0$ \\
    $1^3$ & $-1$ & $\pi_R^2$ & $-2$ \\
      \\
  \end{tabular}
\end{equation}
As shown, we obtain as a by-product that $s$ determines $h$ in unramified splitting types.
\begin{wild}
  The last row is left blank because the integrality of the $K$-valuation yields $s' \in 1/2 + 2\ZZ$, contradicting the known fact $s' \in \{-1/2\} \union \ZZ_{\geq 0}$!
\end{wild}

We now come to our main lemma, which computes $ W_{m_{11},n_{11}} $ for large $m_{11}$.

\subsubsection{Unramified}
\begin{lem}\label{lem:111_strong_zones}
  Suppose $R$ is unramified over $K \supseteq \QQ_2$.
  Let $ m_{11} $ and $ n_c $ be integers, $ m_{11} > 2e \geq n_c > 0 $. Let $\square_C =  \min\left\{2\ell(\hat\omega_C) + 1, e\right\}$ be the squareness of the $\delta = 1$ conic if $s$ is even; let $\square_C = 0$ if $s$ is odd. If $s$ is even, let
  \[
  \tilde{n} = \floor{\frac{2e - s - n_c + 2}{4}} = \ceil{\frac{e - \frac{s}{2} - \ceil{\frac{n_{11}}{2}}}{2}}.
  \]
  Then $ W_{m_{11},n_{11}} $ is given as follows:
  \begin{enumerate}[$($a$)$]
    \item If $n_c > 2e$ (black zone), then
    \begin{align*}
    W_{m_{11},n_{11}} &= \size{H^0} q^{2e - m_{11} - n_c - \floor{\frac{s}{2}}} F\( e, \emptyset, \emptyset \) \\
    &= q^{2e - m_{11} - n_c - \floor{\frac{s}{2}}} [F(e,0,0) + Fx(e,0,0) + Fxx(e,0,0)]
    \end{align*}
    \item If $ 2e - s < n_c \leq 2e $ (purple zone), then
    \begin{align*}
    W_{m_{11},n_{11}} &= 2q^{e - m_{11} - \ceil{\frac{n_c}{2}} - \floor{\frac{s}{2}}} F\( \floor{\frac{n_c}{2}}, e - \floor{\frac{n_c}{2}}, \emptyset \) \\
     &= q^{e - m_{11} - \ceil{\frac{n_c}{2}} - \floor{\frac{s}{2}}} \left[F\(\floor{\frac{n_c}{2}}, e - \floor{\frac{n_c}{2}}, 0\) + Fx\(\floor{\frac{n_c}{2}}, e - \floor{\frac{n_c}{2}}, 0\)\right]
    \end{align*}
    \item If $ 2e - s - 2\square_C < n_c \leq 2e - s $ and $ n_c \geq \dfrac{2e - s}{3} $ (blue zone), then $s$ is even and
    \[
    W_{m_{11},n_{11}} = q^{-m_{11} + \floor{\frac{2e - s - n_c}{4}}} F\( \floor{\frac{n_c}{2}}, e - \tilde{n} - \floor{\frac{n_c}{2}}, \tilde{n} \)
    \]
    \item If $ \square_C < n_c \leq 2e - s - 2\square_C $ (green zone), then
    \[
    W_{m_{11},n_{11}} = q^{-m_{11} + \ceil{\ell_C}}(1 + \epsilon_C)F\( \floor{\frac{n_c}{2}}, \ell_C + \frac{s}{2} + \1_{2\nmid n_c}, e - \ell_C - \frac{s}{2} - \ceil{\frac{n_c}{2}}\)
    \]
    where $\ell_C = \frac{\square_C - 1}{2} \in \{-1/2\} \union \ZZ_{\geq 0}$.
    \item If $n_c < \dfrac{2e - s}{3}$ and $n_c \leq \square_C$ (red zone), then $s$ is even and
    \begin{align*}
    W_{m_{11},n_{11}} &= \sum_{\floor{\frac{n_c}{2}} \leq \ell < \tilde{n}} q^{-m_{11} + \ell}(1 + \epsilon_C) G^\cross\( \ell, \ceil{\frac{n_c}{2}} + \frac{s}{2}, e - \ceil{\frac{n_c}{2}} - \frac{s}{2} - \ell \)  + {}\\
    & \quad {} + q^{-m_{11}+ \floor{\frac{2e - s - n_c}{4}}} G(\tilde{n}, e - 2\tilde{n}, \tilde{n}).
    \end{align*}
  \end{enumerate}
\end{lem}
\begin{proof}
First note that if $s \geq 2e$, then the blue, green, and red zones are empty, and if $s < 2e$ is odd, then since $\ell_0 = -1/2$, the blue and red zones are empty. This ensures that the answers are at least well defined.

By Lemma \ref{lem:to_box}, the support of $W_{m_{11},n_{11}}$ consists of the classes $[\delta] = [\beta]$ in $H^1$ of elements $\beta$ of the box
\[
1 + B_{m_{11},n_c} = \{1 + y \pi^{n_c} \theta_1 + z \pi^{m_{11}} \theta_2, \quad x \in \OO_K^\cross, \quad y,z \in \OO_K\}.
\]
Since $m_{11} > 2e$, the $z$ term has no effect on $[\beta]$, and we ignore it. In particular, the conic
\[
  \coef_{\theta_2} (\beta \xi^2) = 0, \quad \text{that is,} \quad \tr(\hat\omega_C \beta \xi^2) = 0
\]
has a solution $\xi = 1$, so $\epsilon_C(\beta) = \epsilon(\hat\omega_C \beta) = 1$ for all $\beta$ in the box.
  
  \paragraph{Black zone.} In the black zone, we have $\beta \equiv 1 \mod 4\pi$, so $[\beta] = 1$. Hence only $\delta = 1$ yields a nonzero volume, which is, by Lemma \ref{lem:sum_strong},
  \[
    \size{H^0} q^{2e - m_{11} - n_c + \frac{d_0}{2} - (s/2 - v(N(\gamma)))} = \size{H^0} q^{2e - m_{11} - n_c - \floor{\frac{s}{2}}}.
  \]
  
  For the remaining zones, let $ W_{m_{11},n_{11}}' $ denote the claimed value of $W_{m_{11},n_{11}}$ in each case. Our proof method will consist of two steps:
  \begin{itemize}
    \item We prove that $W_{m_{11},n_{11}}(\delta) \leq W_{m_{11},n_{11}}'(\delta)$ for every $\delta$ (the \emph{bounding step}).
    \item We check that
    \[
    \sum_{\delta \in H^1} W'_{m_{11},n_{11}}(\delta) = \size{H^0} q^{2e - m_{11} - n_c} = \sum_{\delta \in H^1} W_{m_{11},n_{11}}(\delta)
    \]
    (the \emph{summing step}), implying that equality must hold for every $\delta$.
  \end{itemize}
  
  \paragraph{Purple zone.} In the purple zone, $s > 0$ defines a splitting $R = K \cross Q$. If we translate $\theta_1$ so that $\theta_1^Q \equiv 0 \mod \pi^s$, then we get $\beta^Q \equiv 1 \mod 4\pi$. Also, $\beta \equiv 1$ mod $\pi^{n_c}$, and indeed, $\beta^{(K)}$ can achieve any value $\equiv 1 \mod \pi^{n_c}$, each congruence class modulo $4\pi$ achieved equally often. So $[\delta] = [\beta]$ ranges uniformly over $T\(\floor{\frac{n_c}{2}}, e - \floor{\frac{n_c}{2}}, \emptyset\)$, and for each class $[\delta]$ that is attained,
  \[
    W_{m_{11},n_{11}}(\delta) = \frac{\size{H^0} q^{2e - m_{11} - n_c - \floor{\frac{s}{2}} }}{\size{T\(\floor{\frac{n_c}{2}}, e - \floor{\frac{n_c}{2}}, \emptyset\)}} = \frac{\size{H^0} q^{2e - m_{11} - n_c - \floor{\frac{s}{2}}}}{\frac{\size{H^0}}{2} q^{e - \floor{\frac{n_c}{2}}}} = 2q^{e - m_{11} - \ceil{\frac{n_c}{2}} - \floor{\frac{s}{2}}},
  \]
  as claimed.
  
  \paragraph{Blue zone.} Note that $s < 2e$ and that $s$ is even (as the bounds imply $\square_C > 0$). Our strategy is to note that
  \[
    \beta \in 1 + B_{m_{11},n_c} \subseteq 1 + B_{m', n_c}
  \]
  for some $m' \leq m_{11}$ for which $[1 + B_{m',n_c}]$ is a boxgroup. Here, we find that the gray-blue condition \eqref{eq:box_blue_ur} in Lemma \ref{lem:boxgps_ur} is the most stringent one, so we take
  \[
    m' = \ceil{\frac{n_c}{2}} + e + \frac{s}{2}
  \]
  and get
  \begin{align*}
  \beta &\in \left[1 + B_{\theta_1}\( \ceil{\frac{n_c}{2}} + e + \frac{s}{2}, n_c\)\right] \\
   &= T\( \floor{\frac{n_c}{2}}, \floor{\frac{\ceil{\frac{n_c}{2}} + e + \frac{s}{2}}{2}} - \floor{\frac{n_c}{2}}, e - \floor{\frac{\ceil{\frac{n_c}{2}} + e + \frac{s}{2}}{2}} \) \\
  &= T\( \floor{\frac{n_c}{2}}, e - \tilde{n} - \floor{\frac{n_c}{2}}, \tilde{n} \).
  \end{align*}
  For each such $\delta$, the value of
  \[
    W_{m_{11},n}(\delta)
  \]
  is controlled by the conic via Lemma \ref{lem:conic_1}, once we know the level
  \[
    \ell = \min\left\{\ell(\delta\hat\omega_C\diamondsuit\heartsuit), \floor{\dfrac{e}{2}}\right\}.
  \]
  We claim that all these conics are blue in the sense of Lemma \ref{lem:conic_1}; this requires
  \[
    2 \ell + 1 \stackrel{?}{\geq} e - n' = e - \ceil{\frac{n_{11}}{2}} = e - \ceil{\frac{n_c + s}{2}},
  \]
  that is,
  \[
    \delta \hat\omega_C \diamondsuit \heartsuit \stackrel{?}{\in} \L_{\ceil{\frac{e - \ceil{\frac{n_c + s}{2}}}{2}}} = \L_{\ceil{\frac{e - \floor{\frac{n_c+s+1}{2}}}{2}}} = \L_{\ceil{\frac{2e - n_c - s - 1}{4}}}.
  \]
  When $\delta = 1$, the required relation
  \[
    \hat\omega_C \diamondsuit \heartsuit \in \L_{\ceil{\frac{2e - n_c - s - 1}{4}}}
  \]
  follows from the given inequality $n_c > 2e - s - 2\square_C$. So it suffices to show that
  \[
    \beta \in \L_{\ceil{\frac{2e - n_c - s - 1}{4}}}.
  \]
  Since $\beta \equiv 1 \mod \pi^{n_c}$, it suffices to show that
  \[
    \floor{\frac{n_c}{2}} \stackrel{?}{\geq} \ceil{\frac{2e - n_c - s - 1}{4}}.
  \]
  But the given red-blue inequality $ n_c \geq \frac{2e - s}{3} $ gives
  \[
    \frac{n_c - 1}{2} \geq \frac{2e - n_c - s - 1}{4},
  \]
  from which the desired inequality follows by taking ceilings. So all conics are blue, and for every $\delta$, $W_{\theta_1,m_{11},n_{11}}(\delta) = q^{-m_{11} + \floor{\frac{2e - s - n_c}{4}}}$ if nonzero. The summing step is now straightforward:
  \begin{align*}
  \sum_{\delta \in H^1} W'_{m_{11},n_{11}}(\delta)
  &= q^{-m_{11} + \floor{\frac{2e - s - n_c}{4}}} \Size{T\( \floor{\frac{n_c}{2}}, e - \tilde{n} - \floor{\frac{n_c}{2}}, \tilde{n} \)} \\
  &= \size{H^0} \cdot q^{-m_{11} + \floor{\frac{2e - s - n_c}{4}} + e - \tilde{n} - \floor{\frac{n_c}{2}} + 2\tilde{n}} \\
  &= \size{H^0} \cdot q^{-m_{11} + e - \floor{\frac{n_c}{2}} + \floor{\frac{2e - s - n_c}{4}} + \floor{\frac{2e - s - n_c + 2}{4}}} \\
  &= \size{H^0} \cdot q^{-m_{11} + e - \floor{\frac{n_c}{2}} + \floor{\frac{2e - s - n_c}{2}}} \\
  &= \size{H^0} \cdot q^{-m_{11} + e - \floor{\frac{n_c}{2}} + e - \frac{s}{2} - n_c + \floor{\frac{n_c}{2}}} \\
  &= \size{H^0} \cdot q^{2e - m_{11} - n_c - \frac{s}{2}}.
  \end{align*}
  This completes the proof, and in particular shows that $\epsilon = 1$ identically on $\hat\omega_C \cdot T\( \floor{\frac{n_c}{2}}, e - \tilde{n} - \floor{\frac{n_c}{2}}, \tilde{n} \)$. This result will be important in proving the remaining zones.
  
  \paragraph{Green zone.} Again, we write
\[
\beta \in 1 + B_{m_{11},n_c} \subseteq 1 + B_{m', n_c}
\]
where $m'$ is as large as possible to make a boxgroup. This time, we find that the gray-green inequality \eqref{eq:box_green_ur} is the most stringent of the conditions in Lemma \ref{lem:boxgps_ur}, so we take $m' = 2 \ceil{\frac{n_c}{2}} + \square_C + s$ (noting that $m'$ is an odd integer) and find that the support of $W_{m_{11},n_{11}}$ is contained in
\[
\left[1 + B_{\theta_1}\( 2 \ceil{\frac{n_c}{2}} + \square_C + s, n_c\)\right] = T\( \floor{\frac{n_c}{2}}, \ell_C + \frac{s}{2} + \1_{2\nmid n_c}, e - \ell_C - \frac{s}{2} - \ceil{\frac{n_c}{2}}\).
\]
We claim all conics are green of the same squareness $\square_C$. The zone boundaries easily imply $\ell_C < \floor{e/2}$, so
\[
[\hat\omega_C \diamondsuit \heartsuit] \in \L_{\floor{\ell_C}} \setminus \L_{\floor{\ell_C} + 1}.
\]
We then note that $[\beta] \in \L_{\floor{\ell_C} + 1}$, because $\beta \equiv 1 \mod \pi^{n_c}$ and we have the inequality $n_c > \square_C$. So $[\beta\hat\omega_C \diamondsuit \heartsuit]$ is also of exact level $\ell_C$. Thus all conics are green of the same squareness, and by Lemmas \ref{lem:conic_1} and \ref{lem:conic_pi}, we have the bound
\begin{align*}
W_{m_{11},n_c} &\leq 2 q^{-m_{11} + \ceil{\ell_C}} T\( \floor{\frac{n_c}{2}}, \ell_C + \frac{s}{2} + \1_{2\nmid n_c}, e - \ell_C - \frac{s}{2} - \ceil{\frac{n_c}{2}}\),
\end{align*}
indeed
\begin{align*}
W_{m_{11},n_c} &\leq q^{-m_{11} + \ceil{\ell_C}} (1 + \epsilon_C) T\( \floor{\frac{n_c}{2}}, \ell_C + \frac{s}{2} + \1_{2\nmid n_c}, e - \ell_C - \frac{s}{2} - \ceil{\frac{n_c}{2}}\).
\end{align*}
where $\epsilon_C(\delta) = \epsilon(\hat\omega_C \delta)$. This completes the bounding step.

To perform the summing step, we need to compute the sum
\[
q^{-m_{11} + \ceil{\ell_C}} \sum_{\delta \in T\( \floor{\frac{n_c}{2}}, \ell_C + \frac{s}{2} + \1_{2\nmid n_c}, e - \ell_C - \frac{s}{2} - \ceil{\frac{n_c}{2}}\)} (1 + \epsilon_C(\delta)).
\]
The term $1$ is found to sum to the desired total $\size{H^0} \cdot q^{2e - m_{11} - n_c - \floor{\frac{s}{2}}}$. 
We claim that
\[
\sum_{\delta \in T\( \floor{\frac{n_c}{2}}, \ell_C + \frac{s}{2} + \1_{2\nmid n_c}, e - \ell_C - \frac{s}{2} - \ceil{\frac{n_c}{2}}\)} \epsilon_C(\delta) = 0,
\]
in other words that $\epsilon$ is \emph{equidistributed} between $1$ and $-1$ in this boxgroup. This follows from Lemma \ref{lem:charm}\ref{charm:constant}: because $\L_{\floor{\ell_C} + 1} \neq \hat\omega_C \L_{\floor{\ell_C} + 1}$ is an uncharmed coset, we have $\epsilon$ equidistributed on cosets of $\L_{e - \floor{\ell_C} - 1}$.

\paragraph{Red zone.}

We first recenter. Changing ${\theta_1}$ to the element $a + b\psi_C^{-1}{\theta_1}$ from Lemma \ref{lem:recentering_ur}, keeping the rest of the resolvent data fixed, gives us a new first vector problem $\P'$ whose associated cubic ring $C'$ has first extender vector $\theta_1' = \eta^2$ is a square. By Lemma \ref{lem:recenter}, the ring-count function $W_{\P}$ simply shifts by $\psi_C$. Note that all boxgroups $T(\ell_0, \ell_1, \ell_2)$ in claimed totals satisfy the gray-green inequality
\[
  \ell_1 \leq \frac{s + \square_C + 1}{2}
\]
Since $\theta_1' \equiv {\theta_1} \mod \pi^{s + \square_C}$, the replacement does not change any of the boxgroups, and their charmed cosets merely translate by $[\psi]$ along with the quadratic form $\epsilon_{C'}(\delta) = \epsilon_C(\psi\delta)$.
   
So it suffices to prove the result in the case that ${\theta_1}$ is a square. Note that $[\hat\omega_C] \in \L_{\floor{e/2}}$ since $\omega_C \equiv \theta_1 \mod 2$.
   
We will actually prove something stronger:
\begin{equation}\label{eq:red_F}
  W_{m_{11},n_{11}} = \sum_{\floor{\frac{n_c}{2}} \leq \ell < \tilde{n}_c} q^{-m_{11} + \ell}(1 + \epsilon_C) F^\cross\( \ell, \ceil{\frac{n_c}{2}} + \frac{s}{2}, e - \ceil{\frac{n_c}{2}} - \frac{s}{2} - \ell \) + q^{-m_{11} + \floor{\frac{2e + s - n_c}{4}}} F(\tilde{n}, e - 2\tilde{n}, \tilde{n}),
\end{equation}
in which we have replaced all $G$'s by $F$'s. We take a moment to realize why this is actually stronger. The claim that $T = T(\tilde{n}, e - 2\tilde{n}, \tilde{n})$ is in the support of $W_{m_{11},n_{11}}$ implies the following:
\begin{itemize}
  \item $\epsilon_C$ is identically $1$ on $T$, and hence
  \item $T$ is maximal isotropic for the Hilbert pairing, and also
  \item the identity coset $T$ is charmed, so $F_T = G_T$, and
  \item the quadratic form $\epsilon_C$ is positively charmed.
\end{itemize}
All the remaining terms use $F_V$ where $V \supseteq T$, so $F$ is interchangeable with $G$ there too. We now prove \eqref{eq:red_F}.

Begin with an arbitrary
\[
\beta = 1 + \pi^{n_c} b {\theta_1}.
\]
Assume first that $\pi^{n_c} b$ is \emph{not} a square modulo $\pi^{e - \frac{s}{2} - \ceil{\frac{n_c}{2}} + 1}$, and let $k$ be the largest integer such that $\pi^{n_c} b$ is a square modulo $\pi^{2k+1}$. Note that $\floor{\frac{n_c}{2}} \leq k < \tilde{n}$. We may write
\[
\pi^{n_c} b = \(\pi^{\ceil{\frac{n_c}{2}}}a\)^2 + \pi^{2k + 1} c, \quad \pi \nmid c.
\]
Let $\zeta = 1 + \pi^{\ceil{\frac{n_c}{2}}}a \eta$. We claim that 
\[
\frac{\beta}{\zeta^2} \in 1 + B_{\theta_1}\( 2\ceil{\frac{n_c}{2}} + s + 2k + 1, 2k + 1 \),
\]
implying that $[\beta] \in T\( k, \ceil{\frac{n_c}{2}} + \frac{s}{2}, e - \ceil{\frac{n_c}{2}} - \frac{s}{2} - k \)$ (compare the $k$th term of the sum). Write
\begin{align*}
\frac{\beta}{\zeta^2}
&= \frac{1 + \( \pi^{\ceil{\frac{n_c}{2}}}a \)^2{\theta_1} + \pi^{2k + 1} c {\theta_1}}{\zeta^2} \\
&= \frac{\left[ 1 + \( \pi^{\ceil{\frac{n_c}{2}}}a \)^2 {\theta_1}\right] \( 1 + \pi^{2k+1} c {\theta_1} \) - \pi^{2\ceil{\frac{n_c}{2}} + 2k + 1} a^2 c \theta_1^2}{\zeta^2} \\
&= \frac{\left[ \( 1 + \pi^{\ceil{\frac{n_c}{2}}}a \eta \)^2 - 2\pi^{\ceil{\frac{n_c}{2}}} a \eta\right] \( 1 + \pi^{2k+1} c {\theta_1} \) - \pi^{2\ceil{\frac{n_c}{2}} + 2k + 1} a^2 c {\theta_1}^2}{\zeta^2} \\
&= \frac{\( \zeta^2 - 2\pi^{\ceil{\frac{n_c}{2}}} a \eta \) \( 1 + \pi^{2k+1} c \theta_1 \) - \pi^{2\ceil{\frac{n_c}{2}} + 2k + 1} a^2 c \theta_1^2}{\zeta^2}.
\end{align*}

We first claim that the denominator $\zeta^2$ belongs to $\OO_K[\theta_1] = B_{\theta_1}(s,0)$. Since
\[
\zeta^2 = \( 1 + \pi^{\ceil{\frac{n_c}{2}}}a \eta\)^2 = 1 + \pi^{2\ceil{\frac{n_c}{2}}}a^2 \theta_1 + 2 \pi^{\ceil{\frac{n_c}{2}}}a \eta,
\]
only the last term is in question, and since $v_K\(\coef_{\theta_2} \eta\) \geq s/2$ (as we saw in Lemma \ref{lem:eta_ur}), the inequality needed is
\[
e + \ceil{\frac{n_c}{2}} + \frac{s}{2} \geq s,
\]
a consequence of $s < 2e$.

Therefore the last term of $\beta/\zeta^2$ is
\begin{align*}
\frac{-\pi^{2\ceil{\frac{n_c}{2}} + 2k + 1} a^2 c \theta_1^2}{\zeta^2} &\in \pi^{2\ceil{\frac{n_c}{2}} + 2k + 1} \OO_K[\theta_1] \\
&= B_{\theta_1}\(2\ceil{\frac{n_c}{2}} + s + 2k + 1, n_c + 2k + 1\) \\
&\subseteq B_{\theta_1}\(2\ceil{\frac{n_c}{2}} + s + 2k + 1, 2k + 1\).
\end{align*}
Thus it is enough to prove that 
\[
\frac{\( \zeta^2 - 2\pi^{\ceil{\frac{n_c}{2}}} a \eta \) \( 1 + \pi^{2k+1} c \theta_1 \) }{\zeta^2} \in 1 + B_{\theta_1}(n_c + s + 2k + 1, 2k + 1) \subseteq 1 + B_{\theta_1}\(2\ceil{\frac{n_c}{2}} + s + 2k + 1, 2k + 1\).
\]
Since the right-hand side is a group and contains $1 + \pi^{2k+1} c \theta_1$, it is enough to show that
\[
\frac{\zeta^2 - 2\pi^{\ceil{\frac{n_c}{2}}} a \eta}{\zeta^2} \in 1 + B_{\theta_1}(n_c + s + 2k + 1, 2k + 1),
\]
that is,
\[
\frac{2\pi^{\ceil{\frac{n_c}{2}}} a \eta}{\zeta^2} \in B_{\theta_1}(n_c + s + 2k + 1, 2k + 1).
\]
But, since $\eta$ belongs to the ring $B_{\theta_1}(s/2, 0)$ and $\zeta^2$ is a unit in that ring,
\begin{align*}
\frac{2\pi^{\ceil{\frac{n_c}{2}}} a \eta}{\zeta^2} &\in \pi^{e + \ceil{\frac{n_c}{2}}} B_{\theta_1}\( \frac{s}{2}, 0 \) \\
&= B_{\theta_1} \( e + \ceil{\frac{n_c}{2}} + \frac{s}{2}, e + \ceil{\frac{n_c}{2}} \) \\
&\subseteq B_{\theta_1}(n_c + s + 2k + 1, 2k + 1),
\end{align*}
where the last step uses $k < \tilde{n}$. This establishes the claim that
\[
[\beta] \in T\( k, \ceil{\frac{n_c}{2}} + \frac{s}{2}, e - \ceil{\frac{n_c}{2}} - \frac{s}{2} - k \).
\]
To replace the $T$ by $T^\cross$, note that
\[
\frac{\beta}{\zeta^2} \equiv 1 + \pi^{2k+1} c \theta_1 \mod \pi^{2k + 2},
\]
a generic unit of exact level $k$. Hence
\[
\beta \in T^\cross\( k, \ceil{\frac{n_c}{2}} + \frac{s}{2}, e - \ceil{\frac{n_c}{2}} - \frac{s}{2} - k \).
\]
Also $\ell(\beta\hat\omega_C) = k$, so at $\delta = \beta$ there is a green conic of squareness $k$ and
\[
W_{m_{11},n_{11}}(\beta) = 2q^{-m_{11} + k}.
\]

If it so happens that $b\pi^{n_c}$ \emph{is} a square modulo $\pi^{e - \frac{s}{2} - \ceil{\frac{n_c}{2}} + 1}$, then writing $\pi^{n_c} b = \(\pi^{\ceil{\frac{n_c}{2}}}a\)^2 + \pi^{2\tilde{n}} c$ and carrying out the above computations, \emph{mutatis mutandis,} shows that
\[
[\beta] \in T(\tilde{n}, e - 2\tilde{n}, \tilde{n})
\]
and
\[
\ell(\beta\hat\omega_C) \geq \tilde n,
\]
so there is a blue conic at $\beta$ and
\[
W_{m_{11},n_{11}} = q^{-m_{11} + \floor{\frac{2e - s - n_c}{4}}}.
\]
Overall,
\[
W_{m_{11},n_{11}} \leq \sum_{\floor{\frac{n_c}{2}} \leq \ell < \tilde{n}} q^{-m_{11} + \ell}(1 + \epsilon_C)F^\cross\( \ell, \ceil{\frac{n_c}{2}} + \frac{s}{2}, e - \ceil{\frac{n_c}{2}} - \frac{s}{2} - \ell \) + q^{-m_{11} + \floor{\frac{2e + s - n_c}{4}}} F(\tilde{n}, e - 2\tilde{n}, \tilde{n}) \\
\]
This completes the bounding step. For the summing step, we note that $T^\cross\( \ell, \ceil{\frac{n_c}{2}} + \frac{s}{2}, e - \ceil{\frac{n_c}{2}} - \frac{s}{2} - \ell \)$ is a union of cosets of $\L_{e - \ell - 1}$ that do not lie in the charmed coset $\L_{\ell + 1}$, so $\epsilon_C$ is equidistributed. The summation then proceeds routinely.
\end{proof}

\paragraph{Further remarks on the red zone.} For general $C$, we end up proving that the coset $\psi_C^{-1} T$ is charmed for each $T = T(\ell_0, \ell_1, \ell_2)$ appearing (either positively or negatively) in the sum. Now $\psi_C \in 1 + B_{\theta_1}(\infty, \square_C)$, which, if
\begin{equation} \label{eq:G_is_F}
  \ell_0 \leq \ell_C,
\end{equation}
is contained in the box defining $T$. Thus the identity coset is charmed and we can replace $G$ by $F$.

\subsubsection{Splitting type \texorpdfstring{$1^3$}{1³}}

\begin{lem}\label{lem:1^3_strong_zones}
  Suppose $R$ has splitting type $1^3$.
  Let $ m_{11} $ and $ n_{11} $ be rational numbers with
  \[
    m_{11} \in \ZZ - \frac{h}{3}, \quad n_{11} \in \ZZ + \frac{h}{3}, \quad m_{11} > 2e \geq n_{11} > 0.
  \]
  Write
  \[
    \Dot m = m_{11} - \frac{2h}{3} \in \ZZ, \quad \Dot n = n_{11} + \frac{2h}{3} \in \ZZ.
  \]
  Let $\square_C = \min\left\{2\floor{\ell(\hat\omega_C) / 2} + 1, e\right\}$ be the squareness of the $\delta = 1$ conic. Also let
  \[
  \tilde{n} = \floor{\frac{2e - \Dot n + 2 + 2h}{4}}.
  \]
  Then $ W_{m_{11},n_{11}} $ is given as follows:
  \begin{enumerate}[$($a$)$]
    \item If $n_{11} > 2e$ (black zone), then
    \[
      W_{m_{11},n_{11}} = q^{2e - m_{11} - n_{11} + 1} F\( e, 0, 0 \).
    \]
    \item If $ 2e - 2\square_C + 2 + \dfrac{4h}{3} \leq n_{11} \leq 2e $ and $ n_{11} \geq \dfrac{2e}{3} $ (blue zone), then
    \[
    W_{m_{11},n_{11}} = q^{-\Dot m + 1 - h + \floor{\frac{2e - \Dot n + 2h}{4}}} F_h\( \floor{\frac{\Dot n}{2}} , e - \floor{\frac{\Dot n}{2}} - \tilde{n}, \tilde{n} \)
    \]
    \item If $ \square_C + 1 - \dfrac{2h}{3} \leq n_{11} \leq 2e - 2\square_C + \dfrac{4h}{3} $ (green zone), then
    \[
    W_{m_{11},n_{11}} = q^{-\Dot m + \ell_C + 1 - h}(1 + \epsilon_C)F_h\( \floor{\frac{\Dot n}{2}}, \ell_C - h + \1_{2\nmid \Dot n}, e - \ceil{\frac{\Dot n}{2}} - \ell_C + h\)
    \]
    where $\ell_C = \frac{\square_C - 1}{2} \in \ZZ_{\geq 0}$.
    \item If $n_{11} < \dfrac{2e}{3}$ and $n_{11} \leq \square_C - \dfrac{2h}{3}$ (red zone), then
    \begin{align*}
    W_{m_{11},n_{11}} &= \sum_{\floor{\frac{\Dot n}{2}} \leq \ell < \tilde{n}} q^{-\Dot m + 1 - h + \ell}(1 + \epsilon_C)G_h^\cross\( \ell, \ceil{\frac{\Dot n}{2}} - h, e - \ceil{\frac{\Dot n}{2}} - \ell + h \) + \\
    &\quad{} + q^{-\Dot m + 1 - h + \floor{\frac{2e - \Dot n + 2h}{4}}} G_h(\tilde{n}, e - 2\tilde{n}, \tilde{n}) \\
    \end{align*}
  \end{enumerate}
\end{lem}
\begin{proof}
  
  By Lemma \ref{lem:to_box}, the support of $W_{m_{11},n_{11}}$ consists of the classes $[\delta] = [\beta]$ in $H^1$ of elements $\beta$ of the box
  \[
  1 + B_{m_{11},n_{11}} = \{1 + y \pi^{n_{11}} \theta_1 + z \pi^{m_{11}} \theta_2, \quad x \in \OO_K^\cross, \quad y,z \in \OO_K\}.
  \]
  Since $m_{11} > 2e$, the $z$ term has no effect on $[\beta]$, and we ignore it. In particular, the conic
  \[
  \coef_{\theta_2} (\beta \xi^2) = 0, \quad \text{that is,} \quad \tr(\hat\omega_C \beta \xi^2) = 0
  \]
  has a solution $\xi = 1$, so $\epsilon(\hat\omega_C \beta) = 1$ for all $\beta$ in the box.
  
  In this case, Lemma \ref{lem:sum_strong} yields
  \[
    \sum_{\delta \in \L_0} W_{m_{11},n_{11}}(\delta) = q^{2e - m_{11} - n_{11} + 1},
  \]
  in particular verifying the total in the black zone.
  
  Again, let $ W_{m_{11},n_{11}}' $ denote the claimed value of $W_{m_{11},n_{11}}$ in each case. Our proof method will consist of bounding and summing, as in the preceding splitting types.
  
  \paragraph{Blue zone.} Our strategy is to note that
  \[
  \beta \in 1 + B_{m_{11},n_{11}} \subseteq 1 + B_{m', n_{11}}
  \]
  for some $m' \leq m_{11}$ for which $[1 + B_{m',n_{11}}]$ is a boxgroup. Here, we find that the gray-blue condition \eqref{eq:box_blue_1^3} in Lemma \ref{lem:boxgps_1^3} is the most stringent one, so we take
  \[
    \Dot{m'}= \ceil{\frac{\Dot n + 1}{2}} + e - h \textand m' = \Dot{m'} + \frac{2h}{3}
  \]
  We easily verify that
  \[
    e - \floor{\frac{\Dot m}{2}} = \tilde{n},
  \]
  so $[1 + B_{\theta_1}(m', n_{11})]$ is exactly the support of $W'_{m_{11},n_{11}}$. We claim all conics are blue. This requires that the squareness
  \[
    \square = \min \left\{2\floor{\frac{\ell(\beta\hat\omega_C)}{2}}+1, e\right\}
  \]
  satisfy
  \[
    \square \stackrel{?}{\geq} e - n' = e - \ceil{\frac{\Dot n}{2}} + h,
  \]
  which simplifies to
  \[
    \beta \hat\omega_C \in \L_{\floor{\frac{2e - \Dot n + 2h}{4}}}.    
  \]
  When $\beta = 1$, we have $\ell = \ell_C$, and the required relation follows from the given blue-green inequality $n_{11} \geq 2e - 2\square_C + 1 + \dfrac{4h}{3}$. So it suffices to show that
  \[
  \beta \in \L_{\floor{\frac{2e - \Dot n + 2h}{4}}}.
  \]
  Using the known relation $\ell(\beta) \geq 2\floor{\Dot n / 2}$
  and the blue-red inequality $n_{11} \geq 2e/3$, this is not hard to prove. So all conics are blue, and the only possible nonzero value of $W_{m_{11},n_{11}}(\delta)$ is
  \[
    q^{-m_{11}^\odot + \floor{\frac{e - n_{11}'}{2}}} = q^{-\Dot m + \floor{\frac{2e - \Dot n + 2h}{4}}}.
  \]
  This completes the bounding step. The summing step is routine. 
  
  This completes the proof, and in particular shows that $\epsilon = 1$ identically on boxgroups of the shape in the lemma. This result will be important in proving the remaining zones.
  
  \paragraph{Green zone.} Again, we write
  \[
  \beta \in 1 + B_{m_{11},n_{11}} \subseteq 1 + B_{m', n_{11}}
  \]
  where $m'$ is as large as possible to make a boxgroup. This time, we find that the gray-green inequality \eqref{eq:box_green_1^3} is the most stringent of the conditions in Lemma \ref{lem:boxgps_1^3}, so we take $m' = n_{11} + \square_C - \frac{2h}{3}$, that is,
  \[
    \Dot{m'} = m_{11} - \frac{2h}{3} = \Dot n + \square_C - 2h \in \ZZ.
  \]
  and find that the support of $W_{m_{11},n_{11}}$ is contained in
  \[
  \left[1 + B_{\theta_1}\(m', n_{11}\)\right] = T_h\( \floor{\frac{\Dot n}{2}}, \ell_C - h + \1_{2\nmid \Dot n}, e - \ceil{\frac{\Dot n}{2}} - \ell_C + h\)
  \]
  We claim all conics are green of the same squareness $\square_C$. It is easy to prove that $\square_C < e - 1$ in this zone, so
  \[
  [\hat\omega_C \diamondsuit \heartsuit] \in \L_{2 \ell_C} \setminus \L_{2 \ell_C + 2}.
  \]
  We then note that $[\beta] \in \L_{2 \ell_C + 2}$, because $[\beta] \in \L_{2\floor{\Dot n/2}}$ and we have the inequality $\Dot n \geq 2\ell_C + 2$. So $[\beta\hat\omega_C \diamondsuit \heartsuit]$ is also of exact level $2\ell_C$ or $2\ell_C + 1$. Thus all conics are green of the same squareness, and by Lemmas \ref{lem:conic_1} and \ref{lem:conic_pi}, we have the bound
  \begin{align*}
  W_{m_{11},n_{11}} &\leq 2 q^{-\Dot m + 1 - h + \ell_C} F_h\( \floor{\frac{\Dot n}{2}}, \ell_C - h + \1_{2\nmid \Dot n}, e - \ceil{\frac{\Dot n}{2}} - \ell_C + h\),
  \end{align*}
  indeed
  \begin{align*}
  W_{m_{11},n_{11}} &\leq q^{-\Dot m + 1 - h + \ell_C} (1 + \epsilon_C) F_h\( \floor{\frac{\Dot n}{2}}, \ell_C - h + \1_{2\nmid \Dot n}, e - \ceil{\frac{\Dot n}{2}} - \ell_C + h\).
  \end{align*}
This completes the bounding step.
  
  To perform the summing step, we need to compute the sum of $q^{-\Dot m + \ell_C} (1 + \epsilon_C)$ over the stated boxgroup. The term $1$ is found to sum to the desired total $\size{H^0} \cdot q^{2e - m_{11} - n_{11} + h}$. We claim that
  \[
  \sum_{\delta \in T\( 2\floor{\frac{\Dot n}{2}} + \frac{1 - h}{2}, 2\ell_C - 1 - h + \1_{2\nmid \Dot n}, 2e - 2\ell_C - 2\ceil{\frac{\Dot n}{2}} + \frac{3h + 1}{2}\)} \epsilon_C(\delta) = 0,
  \]
  in other words that $\epsilon$ is equidistributed between $1$ and $-1$ in this boxgroup. This follows from Lemma \ref{lem:charm}\ref{charm:constant}: because $\L_{2\ell_C + 2} \neq \hat\omega_C \L_{2\ell_C + 2}$ is an uncharmed coset, we have $\epsilon$ equidistributed on cosets of $\L_{2e - 2\floor{\ell_C} - 2}$, of which the boxgroup in question is a union by the green-blue inequality. 
  
  \paragraph{Red zone.} Considerations of space prevent us from writing out the proof, which is like that in the unramified splitting types with the following changes:
  \begin{itemize}
    \item We reduce to the case that $\theta_1 = \eta^2$ is a square using Lemma \ref{lem:recentering_1^3}, and there we will prove the result with the $G$'s replaced by $F$'s. We begin with an arbitrary
    \[
    \beta = 1 + \pi^{n_{11}} b \theta_1 = 1 + \pi^{\Dot{n}}b \cdot \pi^{-2h/3}\theta_1, \quad b \in \OO_K.
    \]
    \item We assume first that $\pi^{\Dot n} b$ is \emph{not} a square modulo $\pi^{e - \ceil{\frac{\Dot n}{2}} + 1 + h}$, and let $k$ be the largest integer such that $\pi^{\Dot n} b$ is a square modulo $\pi^{2k+1}$. We find that
    \[
      [\beta] \in T^\cross\( \ell, \ceil{\frac{\Dot n}{2}} - h, e - \ceil{\frac{\Dot n}{2}} - \ell + h \)
    \]
    is in the support of the $k$th term of the claimed answer and that there is a green conic of level $k$ there.
    \item If $\pi^{\Dot n} b$ in fact \emph{is} a square modulo $\pi^{e  - \ceil{\frac{n_{11}}{2}} + 1 + h}$, then we find that $[\beta]$ is in the support of the last term of the claimed answer and that there is a blue conic there.
  \end{itemize}
The rest of the proof, including the summing step, is completely like the unramified splitting types.
\end{proof}

The $G$'s can be replaced by $F$'s when the index $k$ is at most $\ell_C$, for reasons just like those named above. Another corollary is the following.

\subsubsection{The level parity lemma}
Note that in ramified splitting type, conics of given squareness $k < \floor{e/2}$ occur for $[\delta^\odot]$ of exactly two levels: $2k$ and $2k+1$. The following lemma tells when each occurs, at least when $\delta = 1$.
\begin{lem}\label{lem:level_parity_1^3}
  Take $\delta = 1$. If $[\delta^\odot] \notin \L_e$, then either
  \begin{itemize}
    \item $\delta^\odot$ has even level and $b_1 \in \ZZ + 1/3$, or
    \item $\delta^\odot$ has odd level and $b_1 \in \ZZ - 1/3$.
  \end{itemize}
\end{lem}
\begin{proof}
Consider the first vector problem $\P$ with the given $\theta_1$, with $m_{11} > 2e$ and with $\N_{11}$ as minimally active as can be:
\[
  n_{11} = \begin{cases}
    1/3, & b_1 \in \ZZ + 1/3 \\
    2/3, & b_1 \in \ZZ - 1/3.
  \end{cases}
\]
This lies in the red zone, and we get an answer of the form (using ellipses to mark unimportant portions)
\[
  W_{m_{11}, n_{11}} = \sum_{0 \leq \ell < \ceil{e/2}} (\cdots) G_h^\cross \(\ell, \ldots\) + (\cdots) G_h\(\ceil{\frac{e}{2}}, -\1_{2\nmid e}, \ceil{\frac{e}{2}}\).
\]
Since $\ell(\delta^\odot) < e$, the black-red comparison shows that $\delta = 1$ lies in the support of the $k$th summand, $k = \ell_C =  \floor{\ell(\delta^\odot)/2}$. Now simply note that the $k$th summand consists entirely of elements of exact level $2k$ (for $h = 1$) or $2k + 1$ (for $h = -1$).
\end{proof}

\begin{wild}
  \subsubsection{Splitting type \texorpdfstring{$1^21$}{1²1}}
\begin{lem}\label{lem:1^21_strong_zones}
Suppose $R$ has splitting type $1^21$.
Let $ m_c $ and $ n_c $ be integers with $m_c \geq 2e \geq n_c > 0$. Let $\square_C = \min\left\{2\floor{\ell(\hat\omega_C\diamondsuit\heartsuit) / 2} + 1, e\right\}$ be the squareness of the $\delta = 1$ conic. Also, in types \ref{type:C} and \ref{type:D}, let
\begin{align*}
  \tilde{n} &= \floor{\frac{2e - n_{11} - 2h_\eta + 2}{4}} = \ceil{\frac{e - \ceil{\frac{n_{11}}{2}} - h_\eta}{2}} \\
  \tilde{n}^- &= \floor{\frac{2e - n_{11} - 2h_\eta}{4}} =
  \floor{\frac{e - \ceil{\frac{n_{11}}{2}} - h_\eta}{2}}.
\end{align*}
Then $ W_{m_c,n_c} $ is given as follows:
\begin{enumerate}[$($a$)$]
  \item If $n_c > 2e$ (black zone), then
  \[
  W^{\odot}_{m_{11},n_{11}} = 2q^{2e - m_{11}^\odot - n_c - \frac{s'}{2} + \frac{h_1}{4} - \floor{\frac{h_1}{2}}} L\( 2e + 1 \).
  \]
  \item If $2e - d_0 + 2 \leq n_c \leq 2e$ and $n_c \geq 2e - 2s'$ (dark purple zone; note $s' \geq 0$), then
  \[
  W^{\odot}_{m_{11},n_{11}} = q^{e - m_{11}^\odot - \ceil{\frac{n_c}{2}} - \ceil{\frac{s'}{2}}} L\(e + \floor{\frac{n_c}{2}}\).
  \]
  which, upon simplification, gives
  \[
    W_{m_{11},n_{11}} = q^{e - m_{11} + \frac{d_0 + o_{a_1} - 2s'}{4} - \ceil{\frac{n_c}{2}}} L\(e + \floor{\frac{n_c}{2}}\).
  \]
  \item If $2e - d_0/2 - s' + 1 \leq n_c \leq 2e - d_0 + 1$ (purple zone; note that we must be in type \ref{type:D} or \ref{type:E}), then
  \[
  W^{\odot}_{m_{11},n_{11}} = 2q^{e - m_{11}^\odot - \ceil{\frac{n_c}{2}} - \ceil{\frac{s'}{2}}} F\(\floor{\frac{n_c}{2}}, e - d_0' - \floor{\frac{n_c}{2}}, \emptyset\).
  \]
  which, upon simplification, gives
  \[
    W_{m_{11},n_{11}} = 2 q^{e - m_{11} + \frac{d_0 + o_{a_1} - 2s'}{4} - \ceil{\frac{n_c}{2}}}
    F\(\floor{\frac{n_c}{2}}, e - d_0' - \floor{\frac{n_c}{2}}, \emptyset\).
  \]
  \item In types \ref{type:C} and \ref{type:D}, if $ 2e + d_0' - s' - 2\square_C + 1 - h_1 \leq n_c \leq 2e - d_0' - s'$ (blue zone), then
  \[
  W^{\odot}_{m_{11},n_{11}} = q^{-m_{11}^\odot + \floor{\frac{d_0' + h_\eta}{2}} + \tilde n^-} F\(\floor{\frac{n_c - h_\eta}{2}}, e - d_0' - \floor{\frac{n_c - h_\eta}{2}} - \tilde n , \tilde n\)
  \]
  which, upon simplification, gives
  \[
  W_{m_{11},n_{11}} = q^{-m_{11} + d_0' + \frac{2h_\eta + o_{a_1}}{4} + \tilde n^-} F\(\floor{\frac{n_c - h_\eta}{2}}, e - d_0' - \floor{\frac{n_c - h_\eta}{2}} - \tilde n , \tilde n\).
  \]
  \item In type \ref{type:A} (green zone),
  \[
    W^{\odot} = (1 + \epsilon_C) q^{-m_{11}^\odot} L\(n_c - \frac{1}{2}\)
  \]
  so (as $h_1 = 1$)
  \[
    W = (1 + \epsilon_C) q^{-m_{11} + \frac{d_0 - 1 + o_{a_1}}{4}} L\(n_c - \frac{1}{2}\).
  \]
  \item In type \ref{type:B}, if $n_c \leq 2e - 2s' - 1$ (green zone),
  \[
    W^{\odot} = (1 + \epsilon_C) q^{-m_{11}^\odot + \floor{s'/2}}L\(n_c + s'\),
  \]
  so
  \[
    W = (1 + \epsilon_C) q^{-m_{11} + \frac{d_0 + o_{a_1} + 2 s'}{4}}L\(n_c + s'\).
  \]
  \item In types \ref{type:C}--\ref{type:E}, if $\square_C - d_0' + h_1/2 < n_c \leq 2e + d_0' - s' - 2\square_C - h_1$ (green zone),
  \[
    W^{\odot} = (1 + \epsilon_C) q^{-m_{11}^\odot + \ell_C} F\(\floor{\frac{n_c - h_\eta}{2}}, \ceil{\frac{s'}{2}} + \ell_C - d_0' + h_\eta + \1_{2\nmid n_c + h_\eta}, e - \ceil{\frac{n_c + h_\eta}{2}} - \ceil{\frac{s'}{2}} - \ell_C\).
  \]
  \item If $n_c < (2e - s' - d_0')/3$ and $n_c \leq \square_C - d_0' + h_1/2$ (red zone; note that we must be in type \ref{type:C} or \ref{type:D}), then
  \begin{align*}
    W_{m_{11},n_{11}}^{\odot} &= \sum_{\floor{\frac{n_c - h_\eta}{2}} \leq k < \tilde{n}} q^{-m_{11}^\odot + \floor{\frac{d_0' + h_\eta}{2}} + k}(1 + \epsilon_C)G_{h_\eta}^\cross\( k, \ceil{\frac{s' - d_0' + n_c}{2}} + h_\eta, e - \ceil{\frac{s' + d_0' + n_c}{2}} - k - h_\eta \) + \\
    &\quad{} + q^{-m_{11}^\odot + \floor{\frac{d_0' + h_\eta}{2}} + \tilde n^{-} } G_{h_\eta}(\tilde{n}, e - d_0' - 2\tilde{n}, \tilde{n})
  \end{align*}
\end{enumerate}
\end{lem}

\begin{proof}
Specializing the summation lemma \ref{lem:sum_strong} to this splitting type and using the value of $v\(N\(\gamma\gamma^{\odot}\)\)$ given in Lemma \ref{lem:tfm_conic} yields the sum
\[
  \sum_{\delta \in \L_0} W_{m_{11},n_{11}}^\odot(\delta) = 2q^{2e - m_{11}^\odot - n_c - \frac{s'}{2} + \frac{h_1}{4} - 2\floor{\frac{h_1}{2}}},
\]
which yields the desired total in the black zone. 

\paragraph{Dark purple zone.} Here the support is
\begin{equation*}
  1 + B_{\theta_1}(m_c, n_c) = \left\{[1 + c_1 \pi^{n_c} \theta_1]\right\}.
\end{equation*}
We scale so that $\theta_1 \equiv (1;0)$ mod $\pi^{s' + 1/2}$ and get
\begin{align*}
  1 + \pi^{n_c} \theta_1 &\in \(1 + \pi^{n_c} (1;0)\) \cdot \L_{\min\{n_c + s', 2e\}} \\
  &\subseteq \iota\(\L_{\floor{n_c/2}}\) \cdot \L_{\min\{n_c + s', 2e\}} \\
  &= \L_{e + \floor{n_c/2}} \cdot \L_{\min\{n_c + s', 2e\}} \\
  &= \L_{e + \floor{n_c/2}}.
\end{align*}
Also, the box $B = 1 + B_{\theta_1}(m_c, n_c)$ is a group, at least for $m_c = 2e$ which we can assume. We claim $[B]$ has full signature $\0.\0^{e+\floor{n_c/2}}.\*^{e - \floor{n_c/2}}.\*$ by downward induction on $n_c$. In the base case $n_c = 2e$, we only need to get the intimate unit, which we can get by putting an appropriate unit for $c_1$. The induction step only has content if $n_c$ is odd, and then as $c_1$ varies over $\OO_K/\pi\OO_K$, the resulting element $1 + c_1 \pi^{n_c} \theta_1$ ranges through the $q$-many cosets of $\L_{e + (n_c-1)/2}/\L_{e + (n_c+1)/2}$, as desired.

Again because $B$ is a group, the thickness is uniform and can be retrieved from the summation lemma.

\paragraph{Purple zone.} The purple zone is done by the exact same method; now the supplementary boxgroup
\[
  \iota\(\L_{\floor{n_c/2}}\) = \F\(\floor{\frac{n_c}{2}}, e - d_0' - \floor{\frac{n_c}{2}}, \emptyset\)
\]
appears. (See p.~XII.159 for a bit more detail.)

\paragraph{Blue zone.} As in the other splitting types, our strategy is to note that
\[
\beta \in 1 + B_{m_c,n_c} \subseteq 1 + B_{m', n_c}
\]
for some $m' \leq m_c$ for which $[1 + B_{m',n_c}]$ is a boxgroup. Here, we find that the gray-blue condition is the most stringent one, and we take
\[
  m' = \floor{e + \frac{n_c + s' - d_0' + 1}{2}}
\]
to get
\[
  [\beta] \in [1 + B_{m',n_c}] = T\(\floor{\frac{n_c - h_\eta}{2}}, e - d_0' - \floor{\frac{n_c - h_\eta}{2}} - \tilde n , \tilde n\).
\]
The thickness $W_{m_{11},n_{11}}(\delta)$ is found by noting that the conics are all blue. The summing step then proceeds analogously to the other splitting types.

\paragraph{Green zone.} We have labeled three cases ``green zone,'' because although the details vary by letter type, they all have in common that $\square_C$ is so low that we get green conics of fixed squareness:
\begin{itemize}
  \item Case \ref{type:A} is the simplest. We have $n_c \in \ZZ + 1/2$, and clearly
  \[
    [\beta] \in [\alpha : \alpha \equiv 1 \mod \pi^{n_c}] = \L_{n_c - 1/2}.
  \]
  The conics are all tiny and green. Since $\L_0$ is uncharmed, $\epsilon_C$ is equidistributed on this level space, rendering the summing step easy.
  \item In case \ref{type:B}, we may assume that $\theta_1 \equiv (1;0) \mod \pi^{s' + 1/2}$, and then
  \begin{align*}
    \beta &\in \(1 + \pi^{n_c}c_1(1;0)\) \L_{n_c + s'} \\
    &\subseteq \iota\(1 + \pi^{n_c} \OO_K\) \cdot \L_{n_c + s'}.
  \end{align*}
  The bound $s' < d_0/2 - 1$ defining type \ref{type:B} ensures that the second factor dominates, establishing the bounding step. Since
  \[
    n_c + s' \leq 2e - s' + 1
  \]
  and $\hat\omega_C\diamondsuit\heartsuit$ has level $s'$, we get equidistribution of $\epsilon_C$ on the level space, rendering the summing step easy.
  \item Finally, cases \ref{type:C}--\ref{type:E} parallel the green zone of the other splitting types. We take
  \[
    m' = n_c + 2 \ceil{\frac{s'}{2}} - d_0' + \square_C,
  \]
  making the gray-green inequality an equality, which is the most stringent of the boxgroup-defining inequalities, and get the desired bound. The equidistribution of $\epsilon_C$ holds by level considerations, and the summing step is routine.
\end{itemize}

\paragraph{Red zone.} Considerations of space prevent us from writing out the proof, which is like that in the unramified splitting types with the following changes:
\begin{itemize}
  \item We reduce to the case that $\pi^{h_\eta}\theta_1 = \eta^2$ is a square using Lemma \ref{lem:recentering_1^21}, and there we will prove the result with the $G$'s replaced by $F$'s. We begin with an arbitrary
  \[
  \beta = 1 + \pi^{n_c} b \theta_1 = 1 + \pi^{n_c - h_\eta} b \eta^2, \quad b \in \OO_K.
  \]
  \item Assume first that $\pi^{n_c - h_\eta} b$ is \emph{not} a square modulo $\pi^{2 \tilde n}$, and let $k$ be the largest integer such that $\pi^{n_c - h_\eta} b$ is a square modulo $\pi^{2k+1}$. Write
  \[
    \pi^{n_c - h_\eta} b = \(\pi^{\ceil{\frac{n_c - h_\eta}{2}}} a \)^2 + \pi^{2k + 1} c, \quad \pi \nmid c.
  \]
  \item Let $\zeta = 1 + \pi^{\ceil{\frac{n_c - h_\eta}{2}}} a \eta$. We will show that
  \[
    \frac{\beta}{\zeta^2} \in 1 + B, \quad B = B_{\theta_1}\(2\ceil{\frac{n_c + h_\eta}{2}} + 2k + 1 + s', 2k + 1 + h_\eta\),
  \]
  establishing that
  \begin{equation} \label{beta_in_boxgp_1^21}
    [\beta] \in T\(k, \ceil{\frac{s' - d_0' + n_c}{2}}, e - \ceil{\frac{s' + d_0' + n_c}{2}} - k \).
  \end{equation}
  \item Transform
  \begin{equation}
    \frac{\beta}{\zeta^2} = \(1 - \frac{2\pi^{\ceil{\frac{n_c - h_\eta}{2}}}a \eta}{\zeta^2}\)\(1 + \pi^{2k+1+h_\eta} c \theta_1\) - \frac{\pi^{2\ceil{\frac{n_c + h_\eta}{2} + 2k + 1}} a^2 c \theta_1^2}{\zeta^2}.
  \end{equation}
  \item Check that $\eta \in \OO_K[\theta_1] = B_{\theta_1}(s', 0)$, using the last part of Lemma \ref{lem:eta_1^21} and the relation $s' < 2e - d_0'$ needed for the red zone to be nonempty.
  \item Deduce (using the appropriate bounds) that the summands
  \begin{equation}\label{eq:3_summands_1^21}
    \frac{2\pi^{\ceil{\frac{n_c - h_\eta}{2}}}a \eta}{\zeta^2}, \quad \pi^{2k+1+h_\eta} c \theta_1, \textand \frac{\pi^{2\ceil{\frac{n_c + h_\eta}{2} + 2k + 1}} a^2 c \theta_1^2}{\zeta^2}
  \end{equation}
  all lie in the requisite box $B$, establishing \eqref{beta_in_boxgp_1^21}.
  \item To replace $T$ by $T^\cross$ in \eqref{beta_in_boxgp_1^21}, we note that, when reexamining the moving case of Lemma \ref{lem:boxgps_1^21}, we see that when the $n_c$-value of the box changes from $2k + 2 + h_\eta$ to $2k + 1 + h_\eta$, the new elements all have exact level $2k + h_\eta + d_0'$. (This proof relies on a normalization of $\theta_1$ inconsistent with the $\eta$-lemma, but the normalization does not affect the boxes or boxgroups involved.) So it suffices to show that
  \[
    \frac{\beta}{\zeta^2} \notin 1 + B', \quad B' = B_{\theta_1}\(2\ceil{\frac{n_c + h_\eta}{2}} + 2k + 1 + s', 2k + 2 + h_\eta\).
  \]
  Since first and third quantities in \eqref{eq:3_summands_1^21} lie in $B'$ but the second does not, this is straightforward to see.
  \item If $\pi^{n_c - h_\eta}$ in fact \emph{is} a square modulo $\pi^{2\tilde n}$, then we find by the same method that $[\beta]$ is in the support of the last term of the claimed ring total.
  \item Using bounds on levels, we find that the conics are green for the main sum and blue for the last term, and also that $\epsilon_C$ is equidistributed on each $T^\cross$ appearing in the main sum. The rest of the proof, including the summing step, is completely like the unramified splitting types. \qedhere
\end{itemize}
\end{proof}
\end{wild}

\subsection{The zones when \texorpdfstring{$\N_{11}$}{N11} is weakly active (brown and yellow)}\label{sec:weak}

We now turn our attention to first vector problems such that $N_{11}$ is weakly active and $m_{11} > 2e$. Here we use a significantly different framework. Note that $\bar s > 0$ induces a distinguished splitting $R = K \cross Q$. The resolvent conditions $\M_{11}$ and $\N_{11}$ simplify to
\begin{alignat*}{2}
  \M_{11}&\colon &\tr(\xi_1^2) &\equiv 0 \mod \pi^{m_{11}} \\
  \N_{11}&\colon &\xi_1^{(K)} &\equiv 0 \mod \pi^{n_{11}/2}
\end{alignat*}
For $\xi_1$ to satisfy this, its two $Q$-components must be units.
Under the transformation of Lemma \ref{lem:tfm_conic}, the conditions can also be written as
\begin{alignat*}{2}
  \M_{11}\colon &&\tr(\delta^\odot{\xi_1^\odot}^2) &\equiv 0 \mod \pi^{m_{11}^\odot} \\
  \N_{11}\colon &&\xi_1^{\odot(K)} &\equiv 0 \mod \pi^{\frac{n_{11}}{2} - v^{(K)}(\gamma\gamma^\odot)} \\
  \iff& &\xi_1^{\odot(K)} &\equiv 0 \mod \pi^{\ceil{\frac{n_{11}}{2} - v^{(K)}(\gamma\gamma^\odot)}}.
\end{alignat*}

For simplicity we let
\[
  n_{11}^\odot = \ceil{\frac{n_{11}}{2} - v^{(K)}(\gamma\gamma^\odot)},
\]
which sometimes differs slightly from the $n^\odot$ of Lemma \ref{lem:N11}.

In each of the three applicable splitting types, we will find a \emph{brown zone} where $n_{11}$ is so high that $\xi_1^{\odot(K)}$ can be taken to be $0$, so the only $\delta^\odot$ for which there is a solution are those where $\delta^{(Q)}$ is the class of an element on the traceless line of $Q$. Let $\clubsuit$ be such an element of the form
\[
  \clubsuit = (\alpha + \bar \alpha; \alpha - \bar\alpha; \bar\alpha - \alpha)
\]
where $\alpha$ is a generator for $\OO_Q$ as an $\OO_K$-module. Remarks are in order:
\begin{itemize}
  \item In splitting type $(111)$, we can take $\clubsuit = (1;1;-1)$.
  \item In splitting type $(12)$, we can take $\alpha = (1 + \sqrt{D_0})/2$, giving $\clubsuit = (1; \bar\zeta_2\sqrt{D_0})$.
\begin{wild}
    \item In splitting type $(1^21)$, taking $\alpha = \pi_Q$ gives $\clubsuit = (t; \bar\zeta_2\sqrt{D_0})$.
\end{wild}
\end{itemize}
Denote by $Y_{\P}(\delta^\clubsuit)$ the volume of $\xi'_1$ satisfying the $\M_{11}$ and $\N_{11}$ conditions when
\[
  [\delta^\clubsuit] = [\hat\omega_C \clubsuit \delta],
\]
that is,
\[
  Y_\P(\delta^\clubsuit) = W_{\P}\(\hat\omega_C \clubsuit \delta^\clubsuit\).
\]
Observe that the dependence on $\theta_1$ has been nullified and, if $m_{11}^\odot > 2e$, that $Y_{\P}$ is supported on the vanishing locus of the quadratic form
\[
  \epsilon^\clubsuit(\delta^\clubsuit) \coloneqq \epsilon(\clubsuit\delta^\clubsuit).
\]
Let $Y^\odot_{\P}$ be the corresponding volume of $\xi^\odot_1$.

As in the strong zones, we need a summation lemma. 
\begin{lem} \label{lem:sum_weak}
We have
\[
  \sum_{\delta^\clubsuit \in \delta_0^\clubsuit \L_0} Y_{\P}(\delta^\clubsuit) = \size{H^0} q^{2e - m_{11}^\odot - n_{11}^\odot + 2v_K\(\delta^{\odot(Q)}\)},
\]
where $m_{11}^\odot$, $\gamma$, and $\delta^\odot$ are given by Lemma \ref{lem:tfm_conic} and depend only on the class $\delta^\clubsuit\L_0 \in H^1/\L_0$.
\end{lem}
\begin{proof}
Since $\xi_1^{(Q)}$ must be a unit, $\xi_1^{\odot(Q)}$ has fixed valuation
\[
  v^{(Q))}\(\xi_1^\odot\) = -v^{(Q)}\(\gamma\gamma^\odot\) \eqqcolon r.
\]
Meanwhile, $v\(\xi_1^{\odot(K)}\) \geq n_{11}^\odot$ can vary. For $i \geq n^\odot_{11}$, let $Y_{\P}^{(i)}(\delta^\clubsuit)$ denote the volume of $\xi_1^\odot$ satisfying the first vector problem $\P$ and having $v_K(\xi^{\odot(K)}) = i$. We have
\[
  Y_{\P} = \sum_{i = i_0}^\infty Y_{\P}^{(i)},
\]
since only the measure-zero set where $\xi^{\odot(K)} = 0$ has been dropped. Let
\begin{equation} \label{eq:beta_weak}
  \beta = \delta^\odot {\xi^\odot}^2 \cdot \pi^{2r}
\end{equation}
Observe that $2r$ is an integer (by reference to \eqref{tab:tfm_conic}) so $\beta \in R$. Moreover, as we vary $[\delta^\clubsuit] \in \delta^\clubsuit \L_0$ and $\xi_1^\odot$ in the set whose volume is $Y_{m_{11}, n_{11}}^{(i)}$, we get that $\beta$ varies in the region of primitive members of $\OO_R$ such that
\begin{align*}
  v_K(\beta^{(K)}) &= v_K\(\delta^{\odot(K)}\) + 2i - 2r \\
  v_K(\beta^Q) &= 0 \\
  \lambda^\diamondsuit(\beta) &\equiv 0 \mod \pi^{m_{11}^\odot + 2r}.
\end{align*}
Conditions of the form $\beta^{(K)} \equiv 0 \mod \pi^{n_{11}}$ and $\lambda^\diamondsuit(\beta) \equiv 0 \mod \pi^{m_{11}}$ cut out a box of volume $q^{-m_{11}-n_{11}}$; so the volume of $\beta$ is
\begin{align*}
  &\left(1 - \frac{1}{q}\right)q^{-\(v_K\(\delta^{\odot(K)}\) + 2i + 2r\) - \(m_{11}^\odot - 2r\)} \\
  &= \left(1 - \frac{1}{q}\right) q^{-2i - m_{11}^\odot - v_K\(\delta^{\odot(K)}\) + 4r}.
\end{align*}
Consider the sequence
\[
  \xi' \longmapsto \frac{\xi'}{\(\pi^i; \pi_Q^{2r}\)}
  \longmapsto \frac{\xi'^2}{\(\pi^{2i}; \pi_Q^{4r}\)}
  \longmapsto \frac{\delta^\odot \xi'^2}{\pi^{2r}} = \beta.
\]
Every term is a primitive vector in $\OO_R$, so we can consider projective volumes. The linear map of dividing by $\(\pi^i; \pi_Q^{2r}\)$ scales volumes by $q^{i + 2r}$. Squaring by units scales volumes by $1/\(\size{H^0} \cdot q^{2e}\)$ on regions symmetric under multiplication by $R^\cross[2]$, as we noted above in the proof of Lemma \ref{lem:sum_strong}. Finally, multiplication by
\[
  \frac{\delta^\odot \cdot \(\pi^{2i}; \pi_Q^{4r}\)}{\pi^{2r}}
\]
scales volumes by
\[
  q^{-v_K\(\delta^{\odot(K)}\) - 2v_K\(\delta^{\odot(Q)}\) - 2i + 2r}.
\]
So overall, a volume $Y_{\P}^{(i)}(\delta)$ of $\xi'$ transforms to a volume of
\[
  \frac{1}{\size{H^0}} q^{-2e-i-v_K\(\delta^{\odot(K)}\) - 2v_K\(\delta^{\odot(Q)}\) + 4r} Y_{\P}^{(i)}(\delta).
\]
Summing over $[\delta] \in \delta_0\L_0$,
\[
  \frac{1}{\size{H^0}} q^{-2e-i-v_K\(\delta^{\odot(K)}\) - 2v_K\(\delta^{\odot(Q)}\) + 4r} \sum_{\delta \in \delta_0\L_0} Y_{\P}^{(i)}(\delta) = \left(1 - \frac{1}{q}\right) q^{-2i - m_{11}^\odot - v_K\(\delta^{\odot(K)}\) + 4r},
\]
that is,
\[
  \sum_{\delta \in \delta_0\L_0} Y_{\P}^{(i)}(\delta)
  = \left(1 - \frac{1}{q}\right) q^{2e - i - m_{11}^\odot + 2v_K\(\delta^{\odot(Q)}\)}.
\]
Summing over $i \geq n^\odot_{11}$, the right-hand side becomes a geometric series and we get
\[
  \sum_{\delta \in \delta_0\L_0} Y_{\P}(\delta)
  = q^{2e - m^\odot_{11} - n_{11}^\odot + 2v_K\(\delta^{\odot(Q)}\)},
\]
as desired.
\end{proof}

We now use this to power the summing step in each splitting type.
\subsubsection{Unramified}

Here $\clubsuit = (1 ; \bar\zeta_2 \sqrt{D_0})$, where $D_0$ is scaled so that $D_0 \equiv 1 \mod 4$. Observe that $\clubsuit^Q$ is traceless.

\begin{lem}\label{lem:12_alpha}
There is an $\alpha_0 \in \OO_Q$ such that
\[
  \ker\(\tr_{\OO_R/\OO_K}\) = \<\clubsuit (0;1), \clubsuit (1; \alpha_0^2)\>
\]
as $\OO_K$-modules.
\end{lem}
\begin{proof}
This is a notable example of a lemma of simple form that can be proved using the machinery we have got. (Incidentally, if $R$ is totally split, the choice $\alpha_0 = (0;1)$ works, so we're really only concerned about splitting type $(12)$: but we have no need to separate the splitting types here.)

The conic
\[
  \A(\xi) = \tr(\clubsuit \xi^2)
\]
has determinant $1$. Since $\clubsuit \equiv 1 \mod 2$, $\A$ has maximal squareness $\floor{e/2}$. It has Brauer class $\epsilon(\A) = 1$, since $\xi = (0;1)$ is a solution. Hence, by Lemma \ref{lem:conic_1}, not all the $\OO_K$-points of $\A$ lie in a single $1$-pixel. The reduction of $\A$ modulo $\pi$ consists of $(q+1)$-many $1$-pixels, only one of which has vanishing $K$-coordinate. Hence there is a solution $\xi = (a; \alpha)$ with $\pi \nmid a$. Rescaling, we can take $a = 1$.
\end{proof}

\begin{lem}\label{lem:111_weak_zones}
In splitting types $111$ and $12$. If $h_1 = 0$, let
\[
\tilde{n} = \floor{\frac{2e - n_{11} + 2}{4}} = \ceil{\frac{e - n_{11}^\odot}{2}}.
\]
Then $ Y_{\P} $ is given in terms of $n_{11} = n_{11}$ and as follows:
\begin{enumerate}[$($a$)$]
  \item If $n_{11} > 2e$ (brown zone), then
  \begin{equation*}
  Y_{\P} = \begin{cases}
    2q^{e - m_{11} - n_{11}^\odot}  F\( 0, e, \emptyset \) & [\delta_0^\clubsuit] \in \L_0 \\
    2q^{e - m_{11} - n_{11}^\odot} xF\( 0, e, \emptyset \) & [\delta_0^\clubsuit] \in (1;\pi;\pi)\L_0
  \end{cases}
  \end{equation*}
  \item If $ 0 < n_{11} \leq 2e $ (yellow zone), then
  \begin{equation*}
  Y_{\P} = \begin{cases}
     \ds \sum_{0 \leq \ell < \tilde n} q^{-m_{11} + \ell}(1 + \epsilon^\clubsuit)F^\cross\( \ell, n^\odot_{11}, e - n^\odot_{11} - \ell \) + q^{-m_{11} + \floor{\frac{2e - n_{11}}{4}}} F(\tilde{n}, e - 2\tilde{n}, \tilde{n}), & h_1 = 0 \\
     \ds q^{-m_{11}} (1 + \epsilon^\clubsuit) xF\(0, n^\odot_{11}, e - n^\odot_{11}\) & h_1 = 1.
  \end{cases}
  \end{equation*}
\end{enumerate}
\end{lem}

\begin{proof}
In the brown zone, the conditions imply that
\[
  (\delta^\odot{\xi^{\odot}}^2)^Q \equiv a \clubsuit \mod \pi^{2e + 1}  
\]
for some $a \in \OO_K$, necessarily in $\OO_K^\cross$. So $[\delta^\odot] \in \clubsuit\iota(K) = \clubsuit T(\emptyset, e, \emptyset)$. It's easy to see that all values occur, and the solution volume is constant within the appropriate coarse coset, because the $K$-coordinate of $\beta = \delta^\odot{\xi^{\odot}}^2/\pi^{h_1}$ can range over all of $\pi^{n} \OO_K^\cross$ (for $n \geq 2n^\odot_{11} - h_1$ of the correct parity) while $\beta^{(Q)} \equiv \clubsuit^{(Q)} \mod \pi^{2e + 1}$ remains of constant class. So
\[
  Y_\P(\delta^\clubsuit) = 2q^{e - m_{11} - n_{11}^\odot}
\]
for each $\delta^\clubsuit$ in the support, as desired.

In particular, $T(\emptyset, e, \emptyset)$ is charmed for $\epsilon^\clubsuit$.

In the yellow zone, the conditions imply that
\[
  \beta^Q \equiv a \clubsuit \mod \pi^{2n_{11}^\odot + h_1}
\]
for some $a$ in $\OO_K$, necessarily in $\OO_K^\cross$. Hence
\[
  [\delta^\clubsuit] \in \begin{cases}
    \epsilon^{\clubsuit-1}(1) \intsec T^\cross\(\emptyset, n^\odot_{11}, e - n^\odot_{11} \) & h_1 = 1 \\
    \epsilon^{\clubsuit-1}(1) \intsec T\(0, n^\odot_{11}, e - n^\odot_{11} \) & h_1 = 0.
  \end{cases} 
\]
Note that the boxgroups have $\ell_1 = n^\odot_{11} \leq \ceil{s/2}$, so they are well defined. Also, $F = G$ in the answer, because everything contains $ \iota(K^\cross) \cdot \L_e$, which is charmed for $\epsilon^\clubsuit$. 

In the case $h_1 = 1$, the conics are all tiny and green, and we get the bound
\begin{align*}
  Y_\P \leq q^{-m_{11}}(1 + \epsilon^\clubsuit) xF\(0, n^\odot_{11}, e - n^\odot_{11}\),
\end{align*}
which is exactly as desired. The summing step precedes routinely, noting that $\epsilon^\clubsuit$ is equidistributed because everything is contained in a non-charmed coarse coset.

In the case $h_1 = 0$, that is, $[\delta^\odot] \in \L_0$, some further analysis must be done to narrow the support. By Lemma \ref{lem:conic_lift}, we can assume that $\M_{11}$ is an equality on the nose and also that the $K$-component of $\xi'$ is not \emph{exactly} $0$ (to allow recovery of $[\delta^\clubsuit] = [\delta^\odot\diamondsuit\clubsuit\xi'^2]$). Then $\tr(\delta^\odot\xi'^2) = 0$, so there are $u \in \OO_K^\cross$, $b \in \OO_K$ such that
\[
  \delta^\clubsuit\xi'^2 = u(0;1) + b \pi^{2n'}(1; \alpha_0^2).
\]
Hence we are curious about the $H^1$-class of the right-hand side. Since scaling by $K^\cross$ preserves $H^1$-class, we may assume that $u = 1$. We have $b \neq 0$.

Assume first that $b$ is \emph{not} a square modulo $\pi^{v_K(b) + 2\tilde{n}}$. Write
\[
  b = a^2(1 + \pi^{2k+1}c), \quad \pi \nmid c, \quad 0 \leq k < \tilde n.
\]
We will show that $[\delta^\clubsuit] $ lies in the $k$th term
\[
  T^\cross\(k, n^\odot_{11}, e - n^\odot_{11} - k \).
\]
We compute:
\begin{align*}
  [\delta^\clubsuit] &= [(0;1) + b \pi^{2n^\odot_{11}}(1; \alpha_0^2)] \\
  &= [(b ; 1 + b\pi^{2n^\odot_{11}} \alpha_0^2)] \\
  &= [1 + c\pi^{2k+1} ; 1 + a^2\pi^{2n^\odot_{11}} \alpha_0^2 + a^2 c\pi^{2n^\odot_{11} + 2k + 1} \alpha_0^2] \\
  &\equiv [1 + c\pi^{2k+1} ; 1 + a^2\pi^{2n^\odot_{11}} \alpha_0^2 ] \mod \L_{k + n^\odot_{11}} \\
  &\equiv [1 + c\pi^{2k+1} ; (1 + a \pi^{n^\odot_{11}} \alpha_0)^2 ] \mod \L_{\ceil{(e + n^\odot_{11})/2}} \subseteq \L_{k + n^\odot_{11}} \\
  &= [1 + c\pi^{2k+1} ; 1] \\
  &\in T^\cross\(k, n^\odot_{11}, e - n^\odot_{11} - k \).
\end{align*}
Also, the conics here are green of squareness $k$. Likewise, if $b$ \emph{is} a square modulo $\pi^{v_K(b) + 2\tilde{n}}$, the same computation shows that $[\delta^\clubsuit] \in S(\tilde{n}, e - 2\tilde{n}, \tilde{n})$. Here the conics are blue, and we have the bounding step.

For the summing step, we note that $\epsilon$ is equidistributed on each support $T^\cross(\ell, n^\odot_{11}, e - n^\odot_{11} - \ell)$, since it is a union of cosets of $\L_{e - \ell - 1}$ inside a non-charmed coset of $\L_{\ell+1}$. The sum is then easy to compute and compare against the total of Lemma \ref{lem:sum_weak}.
\end{proof}

When we translate back to $W$'s counting the volumes by the value of $\delta$, the answers change but slightly. In the brown zone, the support is contained in $\hat\omega_C \clubsuit \iota(K)$, but using $s > 2e$ (or the black-brown comparison, when the chosen coarse coset is $\L_0$), this coset is the identity. Intersecting this with the two possible coarse cosets yields the answers
\[
W_{\P} = \left\{\begin{alignedat}{3}
  &2q^{e - m_{11} - \ceil{n_{11}/2}} F\( 0, e, \emptyset \) && [\delta] \in \L_0, \quad &&\text{$s$ even} \\
  &2q^{e - m_{11} - \ceil{n_{11}/2}} xF\( 0, e, \emptyset \) \quad && [\delta] \in (1;\pi;\pi)\L_0, \quad &&\text{$s$ odd} \\
  &2q^{e - m_{11} - \floor{n_{11}/2}} F\( 0, e, \emptyset \) && [\delta] \in \L_0, \quad &&\text{$s$ odd} \\
  &2q^{e - m_{11} - \floor{n_{11}/2}} xF\( 0, e, \emptyset \) && [\delta] \in (1;\pi;\pi)\L_0, \quad &&\text{$s$ even}
\end{alignedat}\right.
\]

\paragraph{Further remarks on the yellow zone.}
In the ``long'' yellow zone answer ($h_1 = 0$), the $F$'s can be changed to $G$'s, because the answer implies that $T(\tilde n, e - 2\tilde n, \tilde n)$ is charmed. When translating from $Y$ back to $W$, we keep these $G$'s, switching quadratic forms from $\epsilon^\clubsuit$ to $\epsilon_C$. However, some of these $G$'s and $xG$'s admit simplifications that are of importance to us.
\begin{itemize}
  \item If $h = 1$ and $s$ is odd, the fact that $\delta = 1$ yields a nonzero ring volume implies that the $xG$ simplifies to $F$:
  \[
    W_{m_{11},n_{11}} = q^{-m_{11}}(1 + \epsilon_C) F\(0, \floor{\frac{n_{11}}{2}}, e - \floor{\frac{n_{11}}{2}}\)
  \]
  \item Still assuming $h = 1$ and $s$ is odd, in the very special case that $n_{11} = 2e$ (on the yellow-brown border), the term
  \[
    \epsilon_C F(0,e,0)
  \]
  admits a curious simplification. Since the charmed coset of $T(0,e,0)$ is $T^\cross(\emptyset,e,0)$, we have
  \[
    G(0, e, 0) = xF(0,e,0).
  \]
  Taking Fourier transforms of both sides,
  \[
    \epsilon_C F(0,e,0) = Fx(0,e,0).
  \]
  This accounts for the ``Fx-yellow-special'' zone in the code.
  
  \item If $h = 1$ and $s$ is even, we claim that the $xG$ simplifies to $xF$. Note that $\hat\omega_C$ is constrained by $s$: we have $\theta_1 \equiv (1;0;0) \mod \pi^s$, so $\omega \equiv (0; \clubsuit^Q) \mod \pi^s$ and
  \[
  [\hat\omega_C\clubsuit] \in \iota(K^\cross) \cdot (1 + \pi^s \OO_R)
  \subseteq   \cdot (1 + \pi^{n_{11}} \OO_R) = T\(\emptyset, \floor{\frac{n_{11}}{2}}, e - \floor{\frac{n_{11}}{2}}\).
  \]
  So the regions in which $\delta^\clubsuit$ lay are also good for $\delta$.
  \item If $h = 0$ and $s$ is even, then $\delta = 1$ must also lie in some term of the sum. Which term it is can be determined using the levels of the conics. If 
  \[
    \ell_C \geq \tilde{n},
  \]  
  then $\delta = 1$ lies in the last term and all $G$'s can be made $F$'s. Otherwise, the terms for
  \[
    \ell \leq \ell_C
  \]
  so simplify. 
\end{itemize}

\begin{wild}
  \subsubsection{Splitting type \texorpdfstring{$1^21$}{1²1}}
As usual, the partially ramified splitting type causes further complications.

Fix a first vector problem with $\N_{11}$ weakly active. Note that
$m_{11} \geq n_{11} + d_0/2 > d_0/2$, which is good enough for $\xi_1^\odot$ to be meaningful. Also, note that
\begin{equation*}
  n^\odot_{11} = \ceil{\frac{n_{11}}{2} - v^{(K)}(\gamma\gamma^\odot)} = 
  \ceil{\frac{2n_{11} - d_0 + h_1}{4}}.
\end{equation*}
If
\[
  n_{11} \leq \frac{d_0 - h_1}{2},
\]
then $\N_{11}$ is automatic, so the answer $W_{m_{11},n_{11}} = W_{m_{11},0}$ is the same as in the beige zone, to be computed below. Assume we are not in this case. Then the cases $h_1 = 1, 3$ offer no solutions, as there
\[
  \(\delta^\odot {\xi_1^\odot}^2\)^{(Q)} \sim \pi_Q^{h_1}
\]
has nonzero $\I$-value to first order, which cannot be canceled by the $K$-term since $\pi \mid \xi_1^{\odot(K)}$. So we assume that $h_1 \in \{0, 2\}$.

The $\beta$ of Lemma \ref{lem:sum_weak} simplifies to
\[
  \beta = \delta^\odot {\xi_1^\odot}^2 \pi^{h_1/2},
\]
and
\[
  \N_{11} \iff \beta^{(K)} \equiv 0 \mod \pi^{n_\beta},
\]
where
\[
  n_\beta = 2n^\odot_{11} - \frac{h_1}{2}
\]
is an integer, either even (if $h_1 = 0$) or odd (if $h_1 = 2$). The support of $Y_\P$ consists of the classes $[\delta^\clubsuit] = [\beta \diamondsuit\clubsuit]$ where $\beta$ ranges over the elements of $\OO_R$ satisfying
\[
  \lambda^\diamondsuit(\beta) = 0, \quad \beta^{(Q)} \sim 1, \quad v(\beta^{(K)}) \geq n_\beta, \quad v(\beta^{(K)}) \equiv n_\beta \mod 2.
\]

\begin{lem} \label{lem:1^2_1_weak_zones}
Assume that $n_{11} > (d_0 - h_1)/2$, $h_1 \in \{0, 2\}$, and $m_{11} > 2e + (d_0 - h_1)/2$. The volume $W^\odot_{\P}$ of solutions to the first vector problem is as follows.
\begin{itemize}
  \item If $n_{\beta} \geq 2e - \ceil{d_0/2} + 1$ (brown zone), then
  \[
    W^\odot_{\P} = \begin{cases}
      2q^{e - m_{11}^\odot - n_{11}^\odot} F(0, e - d_0', \emptyset), & d_0 - h_1 \equiv 0 \mod 4 \\
      2q^{e - m_{11}^\odot - n_{11}^\odot} xF(0, e - d_0', \emptyset), & d_0 - h_1 \equiv 2 \mod 4 \\
      q^{e - m_{11}^\odot - n_{11}^\odot} L(e), & d_0 = 2e + 1
    \end{cases}
  \]
  \item If $d_0/2 \leq n_\beta \leq 2e - d_0/2$ (yellow zone), then, noting that $d_0 = 2d_0'$ is even, the answer is as follows:
  \begin{itemize}
    \item If $d_0 - h_1 \equiv 2 \mod 4$, that is, $a_1' \in \ZZ + 1/2$, so $n_\beta \equiv d_0' + 1$ mod $2$, then
    \[
      W^\odot_{\P} = \begin{cases}
        \ds (1 + \epsilon_C) q^{-m_{11}^\odot + \floor{\frac{d_0' - 1}{2}}} F\(0, \frac{n_\beta - d_0' - 1}{2}, e - \frac{n_\beta + d_0' - 1}{2}\) & h_1 = 0 \\
        \ds (1 + \epsilon_C) q^{-m_{11}^\odot + \floor{\frac{d_0' - 1}{2}}} xF\(0, \frac{n_\beta - d_0' - 1}{2}, e - \frac{n_\beta + d_0' - 1}{2}\) & h_1 = 2
      \end{cases}
    \]
    \item If $d_0 - h_1 \equiv 0 \mod 4$, that is, $a_1' \in \ZZ$, so $n_\beta \equiv d_0'$ mod $2$, then
    \begin{align*}
      W^\odot_{\P} &= \sum_{k=0}^{\tilde n - 1}(1 + \epsilon_C) q^{-m_{11}^\odot + k + \floor{\frac{d_0'}{2}}} G^\cross\(k, \frac{n_\beta - d_0'}{2}, e - \frac{n_\beta + d_0'}{2} - k\) + \\
      &\quad + q^{-m_{11}^\odot + \floor{d_0'/2} + \tilde{n}^-} G(\tilde n, e - d_0' - 2\tilde n, \tilde n)
    \end{align*}
    where
    \[
      \tilde n = \floor{\frac{2e - d_0' - n_\beta}{4}} \textand
      \tilde n^- = \floor{\frac{2e - d_0' - n_\beta - 2}{4}}.
    \]
  \end{itemize}
  \item If $0 < n_\beta \leq d_0/2$ (lemon zone), then
  \begin{align*}
    W^\odot_{\P} &= \sum_{k = 0}^{\floor{\frac{e - n_\beta - 1}{2}}} (1 + \epsilon_C) q^{-m_{11}^\odot + \floor{n_\beta/2} + k} GL^\cross(n_\beta + 2k) + q^{-m_{11}^\odot + \floor{e/2}} GL(e)
  \end{align*}
\end{itemize}
\end{lem}
\begin{proof} We consider each zone in turn.
\paragraph{Brown zone.} Let $n = n_\beta$. Here it suffices to prove that, as $c$ varies over $\OO_K^\cross$, the element
\[
  \beta = (\pi^{n} c ; 1 + \pi_Q \pi^{n} c)
\]
attains all classes in $\iota(\pi^{n} \OO_K^\cross)$ equally often. Now $[\beta] \in \iota(c) \L_{2e - \ceil{d_0/2} + 1}$, but the claim that all classes occur equally often is less immediate. Nevertheless, as promised, we will do it directly, without reference to an $\N_{11}$-lemma.

Since
\[
  \iota : \L_{K, e - \ceil{d_0/2} + 1} \to \L_{2e - \ceil{d_0/2} + 1}
\]
is an isomorphism of groups, we may consider its inverse $\iota^{-1}$. Also, let
\[
  \sigma : \OO_K^\cross/\(\OO_K^\cross\)^2 \to \OO_K^\cross
\]
be a section of the natural projection that has the minimal distance property
\[
  \size{\sigma([x_1]) - \sigma([\x_2])} \leq \size{x_1 - x_2}.
\]
Such a $\sigma$ can be made, for instance, by multiplying elements in the basis given in Proposition \ref{prop:Sh_basis}. Now consider the map
\begin{align*}
  j : \OO_K^\cross &\to \OO_K^\cross \\
  c &\mapsto c \cdot \sigma\big(\iota^{-1}(1 + \pi_Q\pi^{n} c)\big)
\end{align*}
We claim that $j$ is an \emph{isometry} in the $p$-adic metric.

Let $c_1, c_2 \in \OO_K^\cross$ be given with $\size{c_1 - c_2} = \size{\pi}^{k}$. Let $n = n_\beta$. It suffices to prove that
\[
  \sigma\big(\iota^{-1}(1 + \pi_Q\pi^{n} c_1)\big) \equiv \sigma\big(\iota^{-1}(1 + \pi_Q\pi^{n} c_2)\big) \mod \pi^{k+1}.
\]
However, we have
\begin{equation*}
  \ell\(\frac{1 + \pi_Q\pi^{n} c_1}{1 + \pi_Q\pi^{n} c_2}\) \geq \min\{n + k, 2e+1\}
\end{equation*}
and since $\iota$ increases levels by $e$ in the aforementioned region of invertibility,
\[
  \ell\(\iota^{-1}\(\frac{1 + \pi_Q\pi^{n} c_1}{1 + \pi_Q\pi^{n} c_2}\)\) \geq \min\{n - e + k, e + 1\}
\]
which ensures that
\[
  v\(\sigma\(\iota^{-1}\(\frac{1 + \pi_Q\pi^{n} c_1}{1 + \pi_Q\pi^{n} c_2}\)\) - 1\) \geq 2(n - e + k) > k,
\]
as desired. So $j$ is an isometry. In particular, $j$ is bijective and preserves volumes. So the $\beta_c$, which satisfy
\[
[\beta_c] = \big[\iota\big(\pi^{n} j(c)\big)\big],
\]
are equidistributed in $\iota(\pi^{n} \OO_K^\cross)$. Summing over $k = n + 2i$ and getting the scaling from the summation lemma \ref{lem:sum_weak} yields the claimed brown-zone answer.

\paragraph{Yellow zone (short answer).} When $n_\beta \equiv d_0' + 1 \mod 2$, the element $\beta^\odot = (c_1 n_\beta; 1 + c_1 n_\beta \pi_Q)$, $2 \mid v(c_1)$, is seen to lie in
\[
  \iota\(\pi^{h_1/2} \OO_K^\cross\) \cdot L(n_\beta) = \iota(\pi^{h_1/2}) \cdot T\(0, \frac{n_\beta - d_0' - 1}{2}, e - \frac{n_\beta + d_0' - 1}{2}\),
\]
which is exactly the support of the claimed ring total. The conics are all green of level $\floor{(d_0' - 1)/2}$, and $\epsilon^\clubsuit$ is equidistributed because we are in an uncharmed coarse coset, so the summing step is routine.

\paragraph{Yellow zone (long answer).} 
 When $n_\beta \equiv d_0' \mod 2$, more work is in order. By analogy with the other splitting types, our intuition is to try to write 
  $\beta^\odot/\heartsuit$ as an approximate sum of squares. Begin with the multiplication law
\[
  \pi_Q^2 = t\pi_Q + u\pi, \quad t \sim \pi^{d_0'}, \quad u \sim 1.
\]
Note that the elements
\[
  \alpha_0 = (0;1) \textand \alpha_1 = (t; u\pi + t\pi_Q)
\]
span $\ker \lambda^\diamondsuit$ over $K$, but over $\OO_K$, they span the index-$q^{d_0'}$ sublattice
\[
  \Lambda = \<\alpha \in \OO_R : \lambda^\diamondsuit(\alpha) = 0 \textand \alpha_0^{(K)} \equiv 0 \mod \pi^{d_0'}\>.
\]
Due to the inequality $n_\beta \geq d_0'$, $\beta^\odot$ lies in $\Lambda$ and can be written in the form
\[
  \beta^\odot = c_0\alpha_0 + c_1\alpha_1 = (c_1 t; c_0 + c_1 u\pi + c_1 t\pi_Q).
\]
Here $c_0$ is a unit, and we can scale so that $c_0 = 1$; while $c_1 = b\pi^{n_\beta - d_0'},$ where $b$ has even valuation.

If $b$ is \emph{not} a square modulo $\pi^{v(b) + 2\tilde n}$, write
\[
  b = a^2(1 + \pi^{2k+1} c), \quad \pi\nmid c, \quad 0 \leq k < \tilde n.
\]
We claim $[\beta^\odot]$ belongs to the support of the $k$th term in the claimed ring total:
\begin{align*}
  [\beta^\odot] &= \left[\(a^2\pi^{n_\beta - d_0'}(1 + \pi^{2k+1}c) t; 1 + a^2\pi^{n_\beta - d_0'} \pi_Q^2 + a^2 \pi^{n_\beta - d_0' + 2k + 1}\pi_Q^2\)\right] \\
  &= \heartsuit \left[\(1 + \pi^{2k+1}c; 1 + a^2\pi^{n_\beta - d_0'} \pi_Q^2 + a^2 \pi^{n_\beta - d_0' + 2k + 1}\pi_Q^2\)\right] \\
  &\equiv \heartsuit \left[\(1 + \pi^{2k+1}c; (1 + a^2\pi^{n_\beta - d_0'} \pi_Q^2)(1 + a^2 \pi^{n_\beta - d_0' + 2k + 1}\pi_Q^2)\)\right] \\
  &\equiv \heartsuit \left[\(1 + \pi^{2k+1}c; (1 + 2a\pi^{(n_\beta - d_0')/2} + a^2\pi^{n_\beta - d_0'} \pi_Q^2)(1 + a^2 \pi^{n_\beta - d_0' + 2k + 1}\pi_Q^2)\)\right]\\
  &= \heartsuit \left[\(1 + \pi^{2k+1}c; (1 + a\pi^{(n_\beta - d_0')/2})^2(1 + a^2 \pi^{n_\beta - d_0' + 2k + 1}\pi_Q^2)\)\right] \\
  &= \heartsuit \left[\(1 + \pi^{2k+1}c; (1 + a^2 \pi^{n_\beta - d_0' + 2k + 1}\pi_Q^2)\)\right] \\
  &= \heartsuit \left[\(1 + \pi^{2k+1}c; 1 + a^2 \pi^{n_\beta - d_0' + 2k + 1}t\pi_Q + a^2 \pi^{n_\beta - d_0' + 2k + 2} u\)\right] \\
  &\equiv \heartsuit \left[\(1 + \pi^{2k+1}c; 1 + a^2 \pi^{n_\beta - d_0' + 2k + 2} u\)\right] \\
  &= \heartsuit \cdot \iota\((1 + \pi^{2k+1}c)(1 + a^2 \pi^{n_\beta - d_0' + 2k + 2} u)\) \mod \L_{n_\beta + 2k + 1},
\end{align*}
the dropped terms are being of valuation at least $n_\beta + 2k + 1 + 1/2$ by various evident combinations of the known inequalities. We get that 
\[
  [\beta^\odot] \in \heartsuit \iota(1 + \pi^{2k+1}\OO_K)\L_{n_\beta + 2k + 1} = \heartsuit T\(k, \frac{n_\beta - d_0'}{2}, e - \frac{n_\beta + d_0'}{2} - k\),
\]
and evidently \emph{not} in $\heartsuit \iota(1 + \pi^{2k+3}\OO_K)\L_{n_\beta + 2k + 1}$, allowing us to refine $T$ to $T^\cross$.

If $b$ is a square modulo $\pi^{v(b) + 2\tilde n}$, a similar calculation shows that
\[
  [\beta^\odot] \in \heartsuit T\(\tilde n, e - d_0' - 2\tilde n, \tilde n\).
\]
Keeping track of the levels and colors of the conics, we get the bounding step, with the $G$'s specified to $F$'s. Then, using the equidistribution of $\epsilon^\clubsuit$, we get the summing step together with the extra information that $T\(\tilde n, e - d_0' - 2\tilde n, \tilde n\)$ is charmed for $\epsilon^\clubsuit$.

\paragraph{Lemon zone.} Here $n^\odot = 0$. The proof consists of
\begin{itemize}
  \item restricting the level of $\delta^\clubsuit$, using the level parity lemma \ref{lem:level_parity_1^2_1} below and the easy inequality that
  \[
    v(\I(\beta^\odot)) \geq k, \quad k < d_0/2 \implies \ell(\beta^\odot) \geq k.
  \]
  \item dealing with the large conics that appear when $\ell(\heartsuit \beta^\odot) \geq 2\floor{e/2}$. Here a conic for fixed $\delta^\odot$ can have points both in $q$-many \emph{generic} $1$-pixels where $\xi_1^{\odot(Q)} \sim 1$ and a single \emph{special} $1$-pixel where $\xi_1^{\odot(Q) \sim \pi_Q}$. But we do not count these all together, but only for $h_1 = 0$ and $h_1 = 2$ respectively. The computation is identical to that for the beige zone in splitting type $1^3$ (see Lemma \ref{lem:beige_1^3}). \qedhere
\end{itemize}
\end{proof}

We are left with the following level parity lemma:
\begin{lem}\label{lem:level_parity_1^2_1}
  In splitting type $1^21$, if a first vector problem $\P$ is solvable for large $m_{11}$ and $[\delta^\clubsuit] \notin \L_e$, then either
  \begin{itemize}
    \item $[\delta^\clubsuit]$ has even level and $h_1 = 0$, or
    \item $[\delta^\clubsuit]$ has odd level and $h_1 = 2$.
  \end{itemize}
\end{lem}
\begin{proof}
A proof of all cases except $s' = d_0' - 1, \ell_C = (d_0' - 1)/2$ is found on pp.~XII.180--182. TODO \textsc{Level Parity}
\end{proof}

In the brown zone, consideration of where $\hat\omega_C$ lies leads to the answer
\[
W^\odot_{\P} = \begin{cases}
  \ds (1 + \epsilon_C) q^{-m_{11}^\odot + \floor{\frac{d_0' - 1}{2}}} F\(0, \frac{n_\beta - d_0' - 1}{2}, e - \frac{n_\beta + d_0' - 1}{2}\) & h_1 = 0 \\
  \ds (1 + \epsilon_C) q^{-m_{11}^\odot + \floor{\frac{d_0' - 1}{2}}} xF\(0, \frac{n_\beta - d_0' - 1}{2}, e - \frac{n_\beta + d_0' - 1}{2}\) & h_1 = 2.
\end{cases}
\]
In the yellow 
 and lemon zones, the $F$'s can be changed to $G$'s, some of which become $F$'s when we go back from $Y$ to $W$, entirely analogous to the first splitting types.

\end{wild}
\subsection{The beige zone}\label{sec:beige}

\subsubsection{Unramified}
The following is immediate from Lemmas \ref{lem:conic_1}, \ref{lem:conic_pi}, and \ref{lem:squareness}.
\begin{lem} \label{lem:beige_ur}
In unramified splitting type, if $m_{11} > 2e$ but $\N_{11}$ is inactive (beige zone), then the ring volume for $\xi'_1$ is given by
  \begin{align*}
W_{m_{11},0} &= \begin{cases}
\ds q^{-m_{11}} \left[(1 + \epsilon_C) \sum_{0 \leq \ell \leq \floor{\frac{e}{2}} - 1}
q^{\ell}G^{\cross\cross}\(\ell, 0, e - \ell\)
+  (\text{core})\right] & [\delta^\odot] \in \L_0 \\
\ds q^{-m_{11}} (1 + \epsilon_C) \big(xG(0,0,e) + xxG(0,0,e) \big) & [\delta^\odot] \notin \L_0
\end{cases}
\end{align*}
where
\[
(\text{core}) = \begin{cases}
 q^{\frac{e}{2}} \displaystyle \( 1 + \frac{1}{q} \) G\(\frac{e}{2}, 0, \frac{e}{2} \) &
e \text{ even} \\
 q^{\frac{e-1}{2}} \displaystyle (1 + \epsilon_C) G\(\frac{e-1}{2}, 0, \frac{e+1}{2}\) &
e \text{ odd}
\end{cases}
\]
\end{lem}

\paragraph{Further remarks on the beige zone.} In the unramified $h = 1$ case, the answer would more strictly be written as a restriction to the particular coarse coset specified by the discrete data, but we do not do so, as all $(\size{H^0} - 1)$-many non-charmed coarse cosets admit the same extender indices and will be immediately summed. We know which coset of $\L_0$ is charmed, and hence:
\[
  xG(0,0,e) + xxG(0,0,e) = \begin{cases}
    xF(0,0,e) + xxF(0,0,e) & \text{$s$ even} \\
    F(0,0,e) + xxF(0,0,e) & \text{$s$ odd}.
  \end{cases}
\]

\subsubsection{Splitting type \texorpdfstring{$1^3$}{1³}}
Here, a little more care is required to deal with the restrictions on the valuation of $\xi^\odot$ that remain active in the beige zone.

We first use a summation lemma. (We could have proved a summation lemma in the unramified splitting type, but it was unnecessary for finding the answers.)

\begin{lem}\label{lem:sum_beige_1^3}
In splitting type $1^3$, let $m_{11} > 0$. Then
\[
  \sum_{\delta \in H^1} W_{m_{11}, 0} = q^{2e + \frac{h-1}{2} - \Dot m_{11}}.
\]
\end{lem}
\begin{proof}
Analogous to Lemma \ref{lem:sum_strong}, map each $\xi' \in \OO_R^\cross$ in the solution set to a corresponding $\beta = \delta^\odot \xi'^2$. The $\M_{11}$-condition restricts $\beta$ to a space of volume $q^{\frac{h-1}{2} - \Dot m_{11}}$, and the squaring multiplies projective volumes by $q^{-2e}$, establishing the result.
\end{proof}

\begin{lem}\label{lem:beige_1^3}
In splitting type $1^3$, a first vector problem with $m_{11} > 2e$ and $n_{11} \leq 0$ (beige zone) has the answer, for $h_1 = 1$,
\begin{align*}
  W_{m_{11}, 0} &= q^{-\Dot m_{11}}(1 + \epsilon_C) \sum_{k = 0}^{\ceil{e/2} - 1} q^k \big(G_1(k, 0, e - k) - G_1(k+1, -1, e - k)\big) + {} \\
  &\quad + q^{-\Dot m_{11} + \ceil{e/2}} G_1\(\ceil{\frac{e}{2}}, -\1_{2 \nmid e}, \ceil{\frac{e}{2}}\)
\end{align*}
and, for $h_1 = -1$,
\begin{align*}
  W_{m_{11}, 0} &= q^{-\Dot m_{11} + 2}(1 + \epsilon_C) \sum_{k = 0}^{\floor{e/2} - 1} q^k \big(G_{-1}(k, 1, e - k - 1) - G_{-1}(k+1, 0, e - k - 1)\big) + {} \\
  &\quad + q^{-\Dot m_{11} + 2 + \floor{e/2}} G_{-1}\(\floor{\frac{e}{2}}, \1_{2 \nmid e}, \floor{\frac{e}{2}}\).
\end{align*}
\end{lem}
\begin{proof}
Since $\L_e$ is charmed for $\epsilon$, all of the $G$'s can be viewed as selectors for the coset of the indicated boxgroup containing $[\hat\omega_C]$.

For each $[\delta^\odot] = [\delta\hat\omega_C]$ with $\epsilon(\delta^\odot) = 1$, the conic $\M_{\delta^\odot} : \lambda^\diamondsuit(\delta^\odot {\xi^\odot}^2)$ has \emph{some} rational point $\xi^\odot$. Let $h'$ be the value of $h_1$ to which $\xi^\odot$ contributes to the beige-zone answer:
\[
  h' = \begin{cases}
    1 & \xi^\odot \sim 1 \text{ (here $\xi^\odot$ lies in a \emph{generic} $1$-pixel)} \\
    -1 & \xi^\odot \sim \pi_R^2 \text{ (here $\xi^\odot$ lies in the \emph{special} $1$-pixel)}
  \end{cases}
\]
(A potential third case $\xi^\odot \sim \pi_R$ does not satisfy $\M_{\delta^\odot}$ even mod $\pi$.) Then let $\theta_1' \in \bar\zeta_3^{h'} R$ whose associated traceless $\omega_{\theta_1'}$ is the known
\[
  \omega' = \frac{\diamondsuit \delta^\odot {\xi^\odot}^2}{\pi^{4 + 2h}}
\]
 
When $[\delta^\odot] = [\delta\hat\omega_C] \not\in \L_{2\floor{e/2}}$, the corresponding conic $\M_{\delta^\odot}$ has squareness $k < \floor{e/2}$. The solutions to the conic lie in a single $1$-pixel, which must be either generic or special. By the level parity result in Lemma \ref{lem:level_parity_1^3} applied to $\theta_1'$, values with $\ell(\delta^\odot)$ even and odd give the generic and special pixels respectively. This establishes the claimed result when $[\delta^\odot] \notin \L_{2\floor{e/2}}$, which corresponds to the terms $k \leq \floor{e/2} - 1$ of each sum.

We now turn to the case that $[\delta^\odot] \in \L_{2\floor{e/2}}$, so the conic has squareness $\floor{e/2}$. Refer to Lemma \ref{lem:conic_1} for the analysis of conics. If $e$ is even, the conic has solutions in $(q+1)$-many $1$-pixels. These form the line in $\PP(\OO_R/\pi\OO_R)$ given by reducing the conic mod $\pi$; they are $q$ generic and $1$ special. Because the conic has equal volume in each $1$-pixel, the volume can be computed explicitly and contributes
\[
  q^{-\Dot m_{11} + e/2} G_1\(\frac{e}{2}, 0, \frac{e}{2}\)
\]
to the $h_1 = 1$ case and $q$ times as much to the $h_1 = -1$ case (the special pixel gets inflated by $q$ in the $\xi^\odot \mapsto \xi'$ transition), as desired.

If $e$ is odd, the conic has solutions in two $1$-pixels and we need to know whether one of them is special. For $\delta^\odot$ of \emph{exact} level $e-1$, we know, again by Lemma \ref{lem:level_parity_1^3}, that both pixels are generic. For $\delta^\odot \in \L_e$, we get \emph{at least} one generic pixel and \emph{at most} one special pixel, leading to the inequalities
\begin{align*}
  W_{m_{11},0}^{h_1 = 1} \big|_{\hat\omega_C \L_{e-1}} &\geq q^{-\Dot m_{11} + (e - 1)/2}(1 + \epsilon_C) \(G_{-1}\(\frac{e-1}{2}, 0, \frac{e+1}{2}\) - \frac{1}{2}G_{-1}\(\frac{e-1}{2}, 1, \frac{e-1}{2}\)\) \\  
  &= q^{-\Dot m_{11} + (e - 1)/2}(1 + \epsilon_C) \(G_{-1}\(\frac{e-1}{2}, 0, \frac{e+1}{2}\) - G_{-1}\(\frac{e-1}{2}, 1, \frac{e-1}{2}\)\) + {} \\
  &\quad + q^{-\Dot m_{11} + (e-1)/2} G_{-1}\(\frac{e-1}{2}, 1, \frac{e-1}{2}\)
\end{align*}
and 
\begin{align*}
  W_{m_{11},0}^{h_1 = -1} \big|_{\hat\omega_C \L_{e-1}} &\leq q^{-\Dot m_{11} + 2 + (e-1)/2} G_{-1}\(\frac{e-1}{2}, 1, \frac{e-1}{2}\).
\end{align*}
Then, summing and comparing against Lemma \ref{lem:sum_beige_1^3}, we find that equality must hold.
\end{proof}

\begin{wild}
  \subsubsection{Splitting type \texorpdfstring{$1^21$}{1²1}}
TODO. We have an answer, pretty much identical to spl.t{.} $1^3$, but are stuck on level parity.

The same proof works in splitting type $1^21$, we think; we simply state the result.
\begin{lem}
In the beige zone of splitting type $1^21$, if $m_{11}^\odot > 2e$, the total $W_{m_{11},0}^{\odot}$ is given as follows:
\begin{itemize}
  \item If $h_1 = 3$, there are no solutions.
  \item If $h_1 = 1$, then
  \begin{align*}
    W_{m_{11},0}^{\odot} &=
    (1 + \epsilon_C)q^{-m_{11}^\odot} GL^\cross(-1) \\
    &= \begin{cases}
      (1 + \epsilon_C)q^{-m_{11}^\odot} L(0) & \text{types \ref{type:A} ($ d_0 $ even); \ref{type:B} or \ref{type:D} ($d_0$ odd)}\\
      (1 + \epsilon_C)q^{-m_{11}^\odot} \(L(-1) - L(0)\) & \text{types \ref{type:A} ($d_0$ odd); \ref{type:B}--\ref{type:E} ($d_0$ even)}.
    \end{cases}
  \end{align*}
  \item If $h_1 = 0$ or $h_1 = 2$, then
  \[
    W^\odot_{m_{11},0} = \sum_{k = 0}^{\floor{\frac{e - h_1/2 - 1}{2}}} (1 + \epsilon_C) q^{-m_{11}^\odot + \floor{n_\beta/2} + k} GL^\cross(2k + h_1/2) + q^{-m_{11}^\odot + \floor{e/2}} GL(e).
  \]
\end{itemize}
\end{lem}

\end{wild}
\subsection{Orthogonality}

Before proceeding to first-vector problems with $m_{11} \leq 2e$, we prove the following result, which will enable us to compute Fourier transforms of ring totals.
\begin{lem}\label{lem:orth}
For any boxgroup $ T_{\theta_1}(\ell_0, \ell_1, \ell_2) $, its orthogonal complement is given by $ T_{\theta_1}(\ell_2, \ell_1, \ell_0) $.
\end{lem}
\begin{rem}
This is an example of an \emph{explicit reciprocity law}, that is, a formula for the Hilbert symbol in a certain region. There is a wide literature on explicit reciprocity laws, but we suspect that this one is new. In our proof, the only fact we use about the Hilbert pairing is that it is the associated bilinear form to $\epsilon$ (and $\epsilon_C$). This enables us to use various facts about $\epsilon$ gleaned in the preceding sections. We will be concocting various values of the resolvent datum $\theta_1$ and of the discrete datum $n_{11}$ to plug into the lemmas regarding the ring volumes.
\end{rem}
\begin{proof}[Proof of Lemma \ref{lem:orth}]
We carry out the proof in the unramified splitting type only, the proof in the other types being very similar.

We first note that if any of $\ell_0, \ell_1, \ell_2$ is the symbol $\emptyset$, the result follows easily from the self-orthogonality of $\iota(K^\cross)$ (if applicable), as mentioned above. So we can assume that the $\ell_i$ are integers. By definition, they must satisfy
 \begin{align}
\ell_0 + \ell_1 + \ell_2 &= e \\
\ell_1 &\leq \ell_0 + \frac{s}{2} + 1 && (\text{the gray-red inequality})\label{eq:orth_red} \\
\ell_1 &\leq \frac{s + \square_C + 1}{2} && (\text{the gray-green inequality}) \label{eq:orth_green} \\
\ell_1 &\leq \ell_2 + \frac{s}{2} + 1, && (\text{the gray-blue inequality})\label{eq:orth_blue}
\end{align}

We first reduce to the case $\square_C = e$. If $s$ is even, this is accomplished, as in the proof of the red zone, by replacing $\theta_1$ by $\theta_1' = \eta^2$ and noting that, by the gray-green inequality, the boxgroups $T_{\theta_1}(\cdots) = T_{\theta_1'}(\cdots)$ are unchanged. If $s$ is odd, we simply replace $\theta_1$ by a $\theta_1' \equiv (1;0;0) \mod \pi^{\max\{2e + 1, s\}}$ whose corresponding $\hat\omega_{C'}$ is in $\L_{\floor{e/2}}$. (For example, $\theta_1' = (1; 8 \pi^{2s} \bar\zeta_2 \sqrt{D_0}))$ is found to work.) Then since $\theta_1' \equiv \theta_1 \mod \pi^s$ and all boxgroups satisfy the gray-green inequality $\ell_1 \leq (s + 1)/2$, the boxgroups are unchanged. Incidentally, we can also assume as a result of this reduction that $s$ is even.

Now the gray-green inequality is subsumed by the gray-blue and gray-red ones.  Using $\L_i^\perp = \L_{e-i}$, the truth of the lemma for a triple $(\ell_0, \ell_1, \ell_2)$, $\ell_1 \geq 1$, implies its truth for the triples $(\ell_0 + 1, \ell_1 - 1, \ell_2)$ and $(\ell_0, \ell_1 - 1, \ell_2 + 1)$. Hence we can run these reductions backward, increasing $\ell_1$ until we reach an obstruction. This usually happens if either the gray-blue or the gray-red inequality becomes an equality, but it can also happen in two special cases, which we dispatch now:
\begin{itemize}
  \item $\ell_0 = \ell_2 = 0$. Here $s \geq 2e - 2$, and the self-orthogonality of $T(0, e, 0)$ follows from that of $T(\emptyset, e, \emptyset) = \iota(K^\cross)$, unless $e = 1$ and $s = 0$, in which case the gray-blue and gray-red inequalities are also equalities.
  \item Both the gray-blue and gray-red inequalities are $1$ away from equality, that is,
  \[
    (\ell_0, \ell_1, \ell_2) = \(\frac{2e - s}{6}, \frac{e + s}{3}, \frac{2e - s}{6}\).
  \]
  As we see, this is only possible if $2e \equiv s \mod 3$. This space shows up as the support in the blue zone for the first vector problem
  \[
    m_{11} = 2e + 1, \quad n_{11} - s = \frac{2e - s}{3},
  \]
  right on the blue-red boundary. It is therefore isotropic and, by virtue of its size, maximal isotropic.
\end{itemize}
So we are left with the case that, without loss of generality, the gray-blue inequality is an equality
\[
  \ell_1 = \ell_2 + \frac{s}{2} + 1.
\]
We thus have
\[
  (\ell_0, \ell_1, \ell_2) = \( e - 2\ell_2 - \frac{s}{2} - 1, \ell_2 + \frac{s}{2} + 1, \ell_2\).
\]
Let $V = T(\ell_0, \ell_1, \ell_2)$ and $W = T(\ell_2, \ell_1, \ell_0)$ be the claimed orthogonals. By the gray-red inequality, $\ell_0 \geq \ell_2$ and $V \subseteq W$.

If $\ell_0 = \ell_2$, we again have a unique group
\[
  (\ell_0, \ell_1, \ell_2) = \(\frac{2e - s - 2}{6}, \frac{e + s + 2}{3}, \frac{2e - s - 2}{6}\)
\]
which shows up as the support in the blue zone for the first vector problem
\[
  m_{11} = 2e + 1, \quad n_c = \frac{2e - s - 2}{3},
\]
right on the blue-red boundary. It is therefore isotropic and, by virtue of its size, maximal isotropic.

So we may assume that $\ell_0 > \ell_2$ and $V \subsetneq W$. The space $V$ shows up as the blue-zone support for the first vector problem
\[
  m_{11} = 2e + 1, \quad n_c = 2e - 4\ell_2 - 1.
\]
Hence $V$ is isotropic. Also, the first vector problem
\[
  m_{11} = 2e + 1, \quad n_c = 2\ell_2 + 1
\]
lies in the red zone. The first term of its answer is a positive multiple of
\[
  (1 + \epsilon_C) G^\cross(W),
\]
but since $\square_C = e$, the $G$ can be replaced by $F$. Our strategy is as follows. Since $W^\cross = T^\cross(\ell_2, \ell_1, \ell_0)$ generates $W$ as a group, it's enough to show that any $\alpha \in W^\cross$ and $\beta \in V$ are orthogonal. We may write
\begin{align*}
\<\alpha, \beta\> &= \epsilon_C(1) \epsilon_C(\beta) \epsilon_C(\alpha) \epsilon_C(\alpha\beta) \\
&= \epsilon_C(\alpha) \epsilon_C(\alpha\beta).
\end{align*}
So if $\alpha$ and $\beta$ are \emph{not} orthogonal, then, applying the transformation $\alpha \mapsto \alpha\beta$ if need be, we may assume
\[
  \epsilon_C(\alpha) = 1 \textand \epsilon_C(\alpha\beta) = -1.
\]
Since $\epsilon_C(\alpha) = 1$, $\alpha$ is in the support of the red-zone answer $W_{m_{11} = 2e + 1, n_c = 2\ell_2 + 1}$. So there is a $\psi$ in the box $1 + B_{\theta_1}(\infty, 2\ell_2 + 1)$ representing the class $[\alpha] \in H^1$. Recenter, using Lemma \ref{lem:recenter}, and consider the $\P'$ with $\theta_1' = \psi^{-1}\theta_1$, $m_{11}' = 2e+1$, and
\[
  n_c' = 2e - 4\ell_2 - 1.
\]
Since
\[
  \ell_{\P'} = \ell(\alpha) = \ell_2,
\]
this problem is still in the blue zone (right on the blue-green boundary), and we get that $\epsilon_{C'} = 1$ on the boxgroup
\[
  1 + B_{\theta_1'}(\ell_0, \ell_1, \ell_2).
\]
Since $\theta_1' \equiv \theta_1 \mod \pi^{2\ell_2 + 1} \OO_K[\theta_1]$, the subscript can be changed from $\theta_1'$ to $\theta_1$ without changing the boxgroup. So
\[
  1 = \epsilon_{\P'}(\beta) = \epsilon_\P(\alpha \beta),
\]
a contradiction. This completes the proof in the unramified splitting types.
\end{proof}
%

\subsubsection{\texorpdfstring{$E$}{E}-forms of the long answers}
In the red, yellow, and beige zones when $h = 0$, the answer, as announced in \ref{lem:111_strong_zones}, \ref{lem:111_weak_zones}, and \ref{lem:beige_ur}, is a sum of $G$ and $\epsilon_C G$ terms. It is capable of a simplification.

\begin{defn} \label{defn:E}
  If $ T(\ell_0,\ell_1,\ell_2) $ is defined, define
  \[
  E(\ell_0,\ell_1,\ell_2) = \widehat{\epsilon_C \widehat F}(\ell_0,\ell_1,\ell_2)
  \]
\end{defn}
Note that if the $\ell_i$ are integers, we have
\[
E(\ell_0,\ell_1,\ell_2)
= \begin{cases}
G(\ell_0,\ell_1,\ell_2), & \ell_0 \leq \ell_2 \\
q^{\ell_2 - \ell_0} \epsilon_C G(\ell_2,\ell_1,\ell_0), & \ell_2 \leq \ell_0.
\end{cases}
\]

\begin{lem} \label{lem:to_E}
  \begin{enumerate}[$($a$)$]
    \item \label{to_E:red_yellow} A sum of the form
    \[
    \sum_{\ell = a}^{\ceil{b/2} - 1} q^\ell (1 + \epsilon_C) G^\cross(\ell, e-b, b-\ell) + q^{-\floor{b/2}} G^\cross\(\ceil{\frac{b}{2}}, e - 2\ceil{\frac{b}{2}}, \ceil{\frac{b}{2}}\),
    \]
    where $b \geq 2a$, can be rewritten as
    \[
    \sum_{\ell = a}^{b - a} q^\ell E(\ell, e-b, b-\ell) - \sum_{\ell = a}^{b - a - 1} q^\ell E(\ell + 1, e - b - 1, b-\ell).
    \]
    \item \label{to_E:beige} In unramified splitting types, a sum of the form
    \[
    (1 + \epsilon_C) \sum_{a \leq \ell \leq \floor{\frac{e}{2}} - 1}
    q^{\ell}G^{\cross\cross}\(\ell, 0, e - \ell\)
    +  (\text{core}),
    \]
    where $(core)$ denotes the beige-zone core (see Lemma \ref{lem:beige_ur}) and $a < \floor{e/2}$, can be rewritten as
    \[
    \sum_{\ell = a}^{e - a} q^{\ell} E(\ell, 0, e-\ell) - \sum_{\ell = a+1}^{e-a-1} q^{\ell - 1} E(\ell, 0, e-\ell)
    \]
    where the second sum may be empty or may have to be interpreted according to the natural convention
    \[
    \sum_{\ell = 1}^{-1} x_i = -x_0.
    \]
\begin{wild}
      \item TODOWILD Write analogues, if any.
\end{wild}
  \end{enumerate}
\end{lem}

\begin{proof}
  The proof is straightforward, converting each $G$ and $\epsilon_C G$ into an $E$ and merging the ranges of summation.
\end{proof}

\begin{cor}\label{cor:E-forms}
For unramified splitting type, we get in the red zone: 
\begin{align*}
  W_{m_{11},n_{11}} &= \sum_{\ell = \floor{\frac{n_c}{2}}}^{e - n_{11} + \frac{s}{2}} q^{-m_{11} + \ell} E\( \ell, \ceil{\frac{n_{11}}{2}}, e - \ceil{\frac{n_{11}}{2}} - \ell \) \\
  & \quad {} - \sum_{\ell = \floor{\frac{n_c}{2}}}^{e - n_{11} + \frac{s}{2} - 1} q^{-m_{11} + \ell} E\( \ell + 1, \ceil{\frac{n_{11}}{2}} - 1, e - \ceil{\frac{n_{11}}{2}} - \ell\),
\end{align*}
and in the yellow zone:
\begin{align*}
W_{m_{11},n_{11}} &= 
\sum_{\ell = 0}^{e - \ceil{n_{11}/2}} q^{-m_{11} + \ell}E\( \ell, \ceil{\frac{n_{11}}{2}}, e - \ceil{\frac{n_{11}}{2}} - \ell \) \\
& \quad {}- \sum_{\ell = 0}^{e - \ceil{n_{11}/2} - 1} q^{-m_{11} + \ell}E\( \ell+1, \ceil{\frac{n_{11}}{2}} - 1, e - \ceil{\frac{n_{11}}{2}} - \ell \),
\end{align*}
and in the beige zone:
\[
  W_{m_{11}} = \sum_{\ell = 0}^e q^{-m_{11} + \ell} E\(k, 0, e-k\) - \sum_{\ell = 1}^{e-1} q^{-m_{11} + k - 1} E(k, 0, e-k).
\]

For splitting type $1^3$, we get in the red zone:
\begin{align*}
   W_{m_{11},n_{11}} &= \sum_{\ell = \floor{\frac{\Dot n}{2}}}^{e - \Dot n + h} q^{-\Dot m + 1 - h + \ell} E_h\( \ell, \ceil{\frac{\Dot n}{2}} - h, e - \ceil{\frac{\Dot n}{2}} - \ell + h \) \\
&\quad {} - \sum_{\ell = \floor{\frac{\Dot n}{2}}}^{e - \Dot n + h} q^{-\Dot m + 1 - h + \ell} E_h\( \ell + 1, \ceil{\frac{\Dot n}{2}} - h - 1, e - \ceil{\frac{\Dot n}{2}} - \ell + h \),
\end{align*}
and in the beige zone $(h_1 = 1)$:
\[
  W_{m_{11}} = \sum_{\ell = 0}^{e} q^{-\Dot m + \ell} E_{1}(k,0,e-k) - \sum_{\ell = 0}^{e-1} q^{-\Dot m + \ell} E_{1}(k+1,-1,e-k)
\]
and in the beige zone $(h_1 = -1)$:
\[
  W_{m_{11}} = \sum_{\ell = 0}^{e-1} q^{-\Dot m + 2 + \ell} E_{-1}(k,1,e-k-1) - \sum_{\ell = 0}^{e-2} q^{-\Dot m + 2 + \ell} E_{-1}(k + 1,0,e-k-1).
\]
\end{cor}

The following are to be kept in mind when manipulating terms $E(\ell_0, \ell_1, \ell_2)$:
\begin{itemize}
  \item When $\ell_0 \leq \ell_C$ and $\ell_0 \leq \ell_2$, the $E$ came from a $G(\ell_0, \ell_1, \ell_2)$ with $\ell_0 \leq \ell_C$. As we observed in the ``Further remarks'' sections following the red, yellow, and beige zones, such a $G$ is interconvertible with an $F$. Hence such an $E$ will be changed to $F$ if it appears in the final answer (after smearing and applying the $\xi_2$ restrictions: see below).
  \item When $\ell_2 \leq \ell_C$ and $\ell_2 < \ell_0$, the $E$ came from an $\epsilon_C G(\ell_2, \ell_1, \ell_0)$ with $\ell_2 \leq \ell_C$. Such an $\epsilon_C G$ is interconvertible with $\epsilon_C F$ and hence is its own Fourier transform, up to the inevitable factor of $q^e$. We annotate it as $E_{bal}$ (``bal'' for ``balanced''). The same can happen to $\epsilon_C Fx$ and $\epsilon_C Fxx$ terms, which we accordingly notate as $E_{bal} x$ and $E_{bal} xx$.
  \item When $\ell_C < \ell_0$ and $\ell_C < \ell_2$, the above transformations do not apply. We annotate the $E$ as $E_{side}$ and note that, for reflection to hold, either
  \begin{itemize}
    \item The $E_{side}$ pairs with its Fourier transform, an $\epsilon_C F$ from the green zone, or
    \item The $E_{side}$ cancels with a like term for a different value of the discrete data. Indeed, we notice that increasing $n_{11}$ by $2$ in the red or yellow zone causes most of the positive terms to reappear with a negative sign.
  \end{itemize}
\end{itemize}

\subsection{Smeared answers}
It is now necessary to compute $\Sm_{r}(W)$ to solve first vector problems with small $m_{11}$.

It is convenient to express as much as possible in terms of \emph{sparks} that vanish suddenly as the smear index $r$, or equivalently $m_{11}$, is decreased.

\begin{defn}
  A function $W : H^1 \to \CC$ is a \emph{spark of level $r_0$} if for all $r$, $0 \leq r \leq e' + 1$,
  \[
  \Sm_r(W) = \begin{cases}
    W & r \geq r_0 \\
    0 & r < r_0.
  \end{cases}
  \]
\end{defn}
\begin{lem}
  Let $0 \leq r_0 \leq e' + 1$. A function $W$ is a spark of level $r_0$ if and only if its Fourier transform $\widehat W$ is supported on the set $\L_{e' - r_0} \backslash \L_{e' - r_0 + 1}$ of elements of exact level $r_0$ (or $-1/2$, in the case $r_0 = e'+1$). 
\end{lem}
\begin{proof}
  Using the familiar Fourier duality between multiplication and convolution, we have the relation
  \begin{equation} \label{eq:Sm}
    \Sm_r(W) = \( \1_{\L_{e'-r}} \cdot \widehat{W} \) \afterhat\, .
  \end{equation}
  So $W$ is a spark of level $r_0$ if and only if
  \[
  \1_{\L_{\ell}} \cdot \widehat{W} = \begin{cases}
    \widehat{W} & \ell \leq e' - r_0 \\
    0 & \ell > e' - r_0.
  \end{cases}
  \]
  This evidently happens exactly when $\widehat{W}$ is supported on $\L_{e' - r_0} \backslash \L_{e' - r_0 + 1}$, as desired.
\end{proof}
Drawing on the repertory of Fourier transforms we computed in Lemma \ref{lem:FT_charm}, as well as the definition of $E$-functions, we get the following.
\begin{lem}\label{lem:sparks}
  Let $\fcr(\ell)$ be the minimal level of elements in a boxgroup denoted $T(\ell, \ell_1, \ell_2)$; to wit,
  \begin{itemize}
    \item $\fcr(\ell) = \ell$ in unramified splitting type
    \item $\fcr(\ell) = 2\ell$ in splitting type $1^3$ if $h = 1$
    \item $\fcr(\ell) = 2\ell+1$ in splitting type $1^3$ if $h = -1$
\begin{wild}
      \item $\fcr(\ell) = 2\ell + d_0'$ in splitting type $1^21$, letter types \ref{type:D} and \ref{type:E}
      \item $\fcr(\ell) = 2\ell + d_0'+1$ in splitting type $1^21$, letter type \ref{type:C}.
\end{wild}
  \end{itemize}
  The following functions are sparks of the indicated levels.
  \begin{enumerate}[$($a$)$]
    \item \label{sparks:green} If $ T(\ell_0,\ell_1,\ell_2) $ is defined and $ \ell_0, \ell_2 \geq \ell_C + 1 $, then
    \[
    W = \epsilon_C F(\ell_0,\ell_1,\ell_2)
    \]
    is a spark of level $e' - \fcr(\ell_C)$. This will be used in the green zone. 
    \item \label{sparks:red} If both terms are defined and $\ell_1 \geq 1$, then
    \[
    W = q \cdot E(\ell_0,\ell_1,\ell_2) - E(\ell_0, \ell_1 - 1, \ell_2 + 1)
    \]
    is a spark of level $\fcr(e - \ell_2)$. This will be used in the red and yellow zones, as well as the beige zone in ramified splitting types.
    \item \label{sparks:beige} In unramified splitting type, if both terms are defined, then
    \[
    W = q^2 \cdot E(\ell_0,0,\ell_2) - E(\ell_0 + 1, 0, \ell_2-1)
    \]
    is a spark of level $\fcr(\ell_0)$. This will be used in the beige zone.
    \item \label{sparks:Fxx} Expressions of the form
    \[
    W = Gx(), \quad xGx(), \quad \textand Gxx()
    \]
    are sparks of level $e' + 1$. 
  \end{enumerate}
\end{lem}

Although a single $F$ is not generally a spark, we do have the relation $\L_{e' - \fcr(\ell_2)} \subseteq T(\ell_0, \ell_1, \ell_2)$, from which $F(\ell_0, \ell_1, \ell_2)$ and $G(\ell_0, \ell_1, \ell_2)$ are stable under smears of levels $r \geq e' - \fcr(\ell_2)$.

We are now ready to compute explicit answers for the smear. Note that we do \emph{not} try to write $W_{m_{11},n_{11}}$ for each value of $m_{11}$ and $n_{11}$. Instead, we express $W_{m_{11},n_{11}}$ for $m_{11}$ large as a sum of sparks and stable terms whose appearance and disappearance can be coded simply.

One region that we do \emph{not} have to work out is the \emph{gray zone} where $m_{11}$ is so low as to satisfy all the conditions of Lemma \ref{lem:boxgps_ur}, resp{.} \ref{lem:boxgps_1^3},
\begin{wild}
  resp{.} \ref{lem:boxgps_1^21}
\end{wild}
for the defining of boxgroups. There the zone total is simply
\begin{equation}
  W_{m_{11},n_{11}} = \frac{\size{H^0} q^{2e - m_{11} - n_c + \frac{d_0}{2} - (s/2 - v(N(\gamma)))}}{\size{T(\ell_0, \ell_1, \ell_2)}} \cdot F(\ell_0, \ell_1, \ell_2)),
\end{equation}
where $T(\ell_0,\ell_1,\ell_2)$ is the corresponding boxgroup.

We now consider each zone in turn:
\begin{itemize}
  \item The black and brown zones need no smear, as $m_{11} \geq n_{11} > 2e$ automatically in them.
  \item The purple and blue zones have as answer a single $F$. It is stable as long as $m_{11}$ is above the gray zone, since we computed the support by relaxing $\M_{11}$ until we hit the gray zone.
  \item The green-zone answer is a sum of type $F + \epsilon_CF$. The $F$ is stable, as we got it by relaxing $\M_{11}$ until we hit the gray zone. The $\epsilon_CF$ is a spark by Lemma \ref{lem:sparks}\ref{sparks:green}.
  \item In the red, yellow, and beige zones in the charmed coarse coset, the answer is a difference of two series of $E$'s. The $E$'s pair up to form sparks of the types in \ref{lem:sparks}, leaving one singleton (two in the beige zone), a mostly stable $F$.
  In the yellow and beige zones, if the positive sum gets cut down to $1$ or $0$ terms respectively, the negative sum has $-1$ term and must be coded up in a \emph{special zone}. This happens in a few cases, as shown in the code.
\end{itemize}

\subsection{The average value of a quadratic character on a box}
The results in this subsection, coupled with the strong-zone answers in Section \ref{sec:strong}, yield a quick solution to a problem that, at first glance, is unrelated to the topic of this paper.

\begin{thm}\label{thm:char_box}
  Let $K \supseteq \QQ_2$ be a finite extension, and let $R \supset K$ be a tamely ramified \'etale extension of degree $3$. In other words, $R$ is one of the following:
  \begin{itemize}
    \item $K \cross K \cross K$
    \item $K \cross Q$, where $Q$ is the unramified quadratic extension field
    \item the unramified cubic extension field
    \item a totally ramified cubic extension field.
  \end{itemize}
  Let
  \[
  \chi : R^\cross \to \{\pm 1\}
  \]
  be a character, that is, a group homomorphism, such that $\chi(a) = 1$ for all $a \in K^\cross$. (All such characters can be put in the form
  \[
  \chi_\alpha(\xi) = \<\alpha \cdot N_{R/K}(\alpha), \xi\>,
  \]
  where $\alpha \in R^\cross$ and $\<\bullet, \bullet\>$ is the Hilbert symbol.) Let $B$ be an $\OO_K$-sublattice contained in the Jacobson radical of $\OO_R$. (That is, $B$ is a subgroup of $\OO_R$ of finite index closed under multiplication by $\OO_K$, and all elements of $B$ have positive valuation at every field factor of $R$.) Then the average value
  \[
  \frac{1}{\int_B 1} \int_B \chi(1 + \xi) \, \mathrm{d}\xi
  \]
  takes on one of the following values:
  \begin{itemize}
    \item $0$
    \item $1$
    \item $q^{-i}$ for some $i$, $1 \leq i \leq e$. (Here $q = \size{k_K}$, and $e = v_K(2)$ is the absolute ramification index.)
  \end{itemize}
\end{thm}
\begin{proof}
  We may assume that $B$ contains $\pi \OO_{K}$, as enlarging $1 + B$ to $(1 + \pi \OO_K)(1 + B)$ does not change the average of a character that vanishes on $K^\cross$. Now we can take a reduced basis
  \[
  B = \pi \OO_K + \pi^{n} \theta_1 + \pi^m \theta_2,
  \]
  and observe that $B$ is one of the boxes that came up in Lemma \ref{lem:to_box}. The strong-zone total $W_{m,n}(\delta)$ that we have computed in Section \ref{sec:strong} can also be interpreted (up to scaling) as the volume of $\beta \in B$ of class $\delta$. Hence the average in question is
  \[
  I(\chi) = \frac{\sum_\delta \chi(\delta) W_{m,n}(\delta)}{\sum_\delta W_{m,n}(\delta)}
  = \frac{\widehat W_{m,n}(\chi)}{\widehat W_{m,n}(1)}.
  \]
  We wish to understand the possible values of this as $\chi$ ranges over $H^1$.
  In view of the smearing lemma (Lemma \ref{lem:smear}), increasing $m$ only makes the theorem stronger, so we can assume that we are in the case of Lemma \ref{lem:111_strong_zones} or \ref{lem:1^3_strong_zones}.
  
  In the black, purple, and blue zones, $W_{m,n}(\delta) = c \cdot F(T)$ for some subgroup $T \subseteq H_1$, so $I$ takes the value $1$ or $0$ according as $\chi(T) = 1$ or not.
  
  In the green zone,
  \[
  W_{m,n}(\delta) = c (F(T) + \epsilon_C F(T))
  \]
  for some boxgroup $T = T(\ell_0, \ell_1, \ell_2)$ on which $\epsilon_C$ is equidistributed, so
  \[
  I = \begin{cases}
    F(T^\perp) + G(T^\perp), & T^\perp \supseteq T \\
    \ds F(T^\perp) + \epsilon_C \frac{\sqrt{\size{V}}}{\size{T}} G(T), & T^\perp \subsetneq T.
  \end{cases}
  \]
  In either case, the identity coset $T$ is uncharmed, so the two terms have disjoint supports. In the first case, $I$ is either $1$ or $0$. In the second case, we can additionally get a value of 
  \[
  \epsilon_C \frac{\sqrt{\size{V}}}{\size{T}} = \pm q^{\ell_0 - \ell_2}.
  \]
  The value of $i = \ell_2 - \ell_0$ evidently satisfies $1 \leq i \leq e$ by our setup of boxgroups.
  
  Finally, in the red zone, the Fourier transform is easier to compute using the $E$-form (Corollary \ref{cor:E-forms}), which is of the form
  \[
  W_{m,n} = \sum_{\ell = a}^{b - a} q^{c + \ell} E(\ell, e-b, b-\ell) - \sum_{\ell = a}^{b - a - 1} q^{c - \ell} E(\ell + 1, e - b - 1, b-\ell).
  \]
  The Fourier transform, by definition of $E$, is
  \begin{align*}
    \widehat{W}_{m,n} = c \epsilon_C \cdot \Bigg(\sum_{\ell=a}^{b-a} q^{b - \ell} F(b - \ell, e - b, \ell) - \sum_{\ell = a}^{b - a - 1} q^{b - \ell - 1} F(b - \ell, e - b - 1, \ell + 1) \Bigg).
  \end{align*}
  Scaling by $\widehat{W}_{m,n}(1) = c q^{b-a}$, we get 
  \[
  I = \epsilon_C \Bigg(\sum_{\ell=a}^{b-a} q^{a - \ell} F(b - \ell, e - b, \ell) - \sum_{\ell = a}^{b - a - 1} q^{a - \ell - 1} F(b - \ell, e - b - 1, \ell + 1) \Bigg).
  \]
  The boxgroups on which the terms are supported form a nested chain
  \[
  T(b-a,e-b,a) \subset T(b-a,e-b-1,a+1) \subset T(b-a-1,e-b,a+1) \subset \cdots \subset T(a, e-b,a)
  \]
  that appear alternately with positive and negative coefficients. So there are two types of behavior upon plugging in any individual $\chi$:
  \begin{itemize}
    \item $\chi$ lies in an even number of boxgroups in the chain, and they cancel in pairs to yield $I = 0$.
    \item $\chi$ lies in an odd number of boxgroups in the chain, and only the smallest one yields a contribution $I = \pm q^{a - \ell}$. The negative of the exponent satisfies
    \[
    0 \leq \ell - a \leq b - a \leq e.
    \]
  \end{itemize}
  We must exclude the possibility that $I = -1$. This can be done by noting that $I$ is the average value of a character that takes only the values $1$ and $-1$, and that in a small neighborhood of $1 \in 1 + B$, $\chi$ is identically $1$ (continuity of $\chi$ is automatic, because all values near $1$ are squares).
\end{proof}
\section{Ring volumes for \texorpdfstring{$\xi'_2$}{xi'2}}

Fix $\xi_1$ and all the data leading up to it. By Lemma \ref{lem:inactives}, there are only three possibilities for $\bar \xi_2$: either it is unrestricted, in which case the volume is given by the white-zone answer in Lemma \ref{lem:white}, or it is restricted by $\M_{12}$ or $\M_{22}$.

\begin{wild}
  \subsection{Considerations in splitting type \texorpdfstring{$1^21$}{1²1}}
Notice the extra complications that occur when $R = K \cross Q$ is partially wildly ramified. At first glance, it would seem that we have made everything harder by using the reduced indices $a_1, a_2, a_3$ for our discrete data because it is the extender indices $\bar a_1, \bar a_2, \bar a_3$ that govern whether the resolvent conditions $\M_{ij}$ and $\N_{ij}$ are active. Nevertheless, we shall get by, using the Lemma \ref{lem:bar_a_2} to recover $\bar a_2$ from the $a_i$:
\begin{equation}
  \label{eq:bar_a_2_copy}
  \bar a_2 = \min \left\{a_2 + v^{(K)}(\xi_1),\; a_2 + v^{(K)}(\xi_2),\; a_2 + \frac{d_0 - 1}{2},\; a_3\right\}.
\end{equation}
We first examine the summand $v^{(K)}(\xi_1)$. Because of the third argument of the minimum, we only need $\min\{v^{(K)}(\xi_1), (d_0 - 1)/2\}$, which can be controlled as follows:
\begin{lem}\label{lem:wild_v_xi_1} For any $\xi_1$ satisfying the $\M_{11}$ and $\N_{11}$ conditions, if $m_{11}$ is large enough that $m_{11} > 2e$ and $m_{11}^\odot > 2e$, the value of
  \[
    g_1 = \min\{v^{(K)}(\xi_1), (d_0 - 1)/2\}
  \]
  is constrained as follows: (TODO strengthen to $d_0/2$ as seems to be required by the code)
  \begin{enumerate}[$($a$)$]
    \item In the yellow and darker zones, it depends only on the value of $s'$:
    \begin{itemize}
      \item In type \ref{type:A} of Lemma \ref{lem:types_1^21}, we have
      \[
        g_1 = \frac{d_0 - 1}{4}.
      \]
      \item In type \ref{type:B},
      \[
        g_1 = \frac{2s' + d_0}{4}.
      \]
      \item In types \ref{type:C}--\ref{type:E},
      \[
        g_1 = \frac{d_0 - 1}{2}.
      \]
    \end{itemize}
    \item In the lemon and beige zones, it depends only on what term $GL(\ell)$ or $GL^\cross(\ell)$ we are in:
    \begin{itemize}
      \item On $GL^\cross(-1)$,
      \[
      g_1 = \frac{d_0 - 1}{4}.
      \]
      \item On $GL^\cross(\ell)$, $0 \leq \ell < d_0/2 - 1$,
      \[
      g_1 = \frac{2\ell + d_0}{4}.
      \]
      \item Within $GL(\ceil{d_0/2} - 1)$,
      \[
      g_1 = \frac{d_0 - 1}{2}.
      \]
    \end{itemize}
  \end{enumerate}
\end{lem}
\begin{proof}
The cases have been divided up by the value of $\ell = \ell\(\delta^\odot\) = \ell\(\delta_0 \hat\omega_C \diamondsuit\)$.

If $\ell = -1$, we must have $h_1 = 1$, $\delta^{\odot(Q)} \sim \pi_Q, \xi_1^{\odot(Q)} \sim 1$. Then
\[
  \I\(\delta^{\odot(Q)}{\xi_1^{\odot(Q)}}^2\) \sim 1,
\]
so to satisfy $\M_{11}$, we must have
\[
  v^{(K)}\(\xi_1^\odot\) = 0,
\]
so by Table \eqref{tab:tfm_conic},
\[
  v^{(K)}\(\xi_1\) = \frac{d_0 - h_1}{4} = \frac{d_0 - 1}{4},
\]
as desired.

If $0 \leq \ell < d_0/2 - 1$, then $h_1$ is either $0$ or $2$. We have
\[
  v\(\I\(\delta^{\odot(Q)}{\xi_1^{\odot(Q)}}^2\)\) = \ell + \frac{h_1}{2},
\]
because the square ${\xi_1^{\odot(Q)}}^2$ is too close to being in $K$ to cancel the main term of $\I\(\delta^{\odot(Q)}\)$, 
 so to satisfy $\M_{11}$,
\[
  v^{(K)}\(\delta^{\odot}{\xi_1^\odot}^2\) = \ell + \frac{h_1}{2}.
\]
Since $\delta^{\odot(K)}$ is a unit, this forces $h_1/2 \equiv \ell \mod 2$ and
\[
  v^{(K)}\(\xi_1^{\odot}\) = \frac{\ell + \frac{h_1}{2}}{2} = \frac{2\ell + h_1}{4}.
\]
Hence
\[
  v^{(K)}\(\xi_1\) = \frac{2\ell + h_1}{4} + \frac{d_0 - h_1}{4} = \frac{2\ell + d_0}{4},
\]
as desired.

Finally, if $\ell \geq d_0/2 - 1$, then $h_1$ is either $0$ or $2$. We have
\[
  v\(\I\(\delta^{\odot(Q)}{\xi_1^{\odot(Q)}}^2\)\) \geq \frac{d_0}{2} - 1 + \frac{h_1}{2},
\]
so to satisfy $\M_{11}$,
\[
v^{(K)}\(\delta^{\odot}{\xi_1^\odot}^2\) \geq \frac{d_0}{2} - 1 + \frac{h_1}{2}.
\]
Since $\delta^{\odot(K)}$ is a unit, the $\delta^{\odot}$ factor can be dropped, and
\[
  v^{(K)}\(\xi_1\) \geq \frac{\frac{d_0}{2} - 1 + \frac{h_1}{2}}{2} + \frac{d_0 - h_1}{4} = \frac{d_0 - 1}{2}.
\]
So $g_1 = (d_0 - 1)/2$, as desired.
\end{proof}

\end{wild}
\subsection{\texorpdfstring{$\M_{12}$}{M12}}
\begin{lem} \label{lem:M12_tame_tfm}
  Assume that $R$ is tamely ramified and $\M_{12}$ is active. Fix $\xi_1$ satisfying the $\M_{11}$, $\N_{11}$ conditions, and normalize $\gamma_2$ as in Lemma \ref{lem:gamma_white}. Then $\M_{12}$ is equivalent to a relation of the form
  \[
  \lambda^\diamondsuit(\alpha_1 \xi'_2) \equiv 0 \mod \pi^{\floor{m_{12}}}
  \]
  where $\alpha_1 \in \OO_R$ is primitive.
\end{lem}
\begin{proof}
  The $\M_{12}$ condition says that
  \[
  \tr(\xi_1\xi_2) \equiv 0 \mod \pi^{m_{12}}.
  \]
  We have
  \[
  \tr(\xi_1\xi_2) = \tr(\xi_1\gamma_2 \xi_2') = \lambda^\diamondsuit(\diamondsuit \xi_1 \gamma_2 \xi_2').
  \]
  Observe that
  \[
  \diamondsuit \xi_1 \gamma_2 \in \pi^{-a_1' - a_2'} R.
  \]
  Let $r_{12} \in a_1' + a_2' + \ZZ$ be the unique value such that
  \[
  \alpha_0 = \pi^{r_{12}} \diamondsuit \xi_1 \gamma_2
  \]
  is a primitive vector in $\OO_R$. Then 
  \[
  \M_{12} : \lambda^\diamondsuit(\alpha_1 \xi'_2) \equiv 0 \mod \pi^{m_{12} + r_{12}},
  \]
  and the exponent is seen to be an integer. To show that it is $\floor{m_{12}}$, it's enough to prove that
  \[
  -1 < r_{12} \leq 0.
  \]
  We examine the cases.
  \begin{itemize}
    \item If $R$ is unramified and $[\delta \hat\omega_C] \in \L_0$, then $\diamondsuit$, $\gamma_2$ are units and $\xi_1$ is primitive, so $r_{12} = 0$.
    \item If $R$ is unramified and $[\delta \hat\omega_C] \in (1; \pi; \pi) \L_0$, then $\diamondsuit$ is a unit. Since $\M_{11}$ is satisfied, we have $j_0 \neq 1$ and (with respect to the na\"ive choice of $\gamma_1$ from Lemma \ref{lem:gamma_white}, which differs from how we actually found the ring volume for $\xi_1$), $\gamma_1 \sim (\sqrt{\pi}; 1 ; 1)$, $(\xi_1')^Q$ is primitive. There are then two subcases, $j_0 = 2$ and $j_0 = 3$. In both cases we find that the scaling of $\alpha_0$ is controlled by the $Q$-components and $r_{12} = -1/2$ or $0$ respectively.
    \item If $R$ is totally ramified, then $\diamondsuit \sim \pi_R^2$, $\xi_1 \sim 1$, $\gamma_2 \sim 1$. According as the product $\diamondsuit \xi_1 \gamma_2$ lies in $R$, $\bar\zeta_3 R$, or $\bar\zeta_3^2 R$, we must take $r_{12} = 0$, $-2/3$, or $-1/3$ respectively.
  \end{itemize}
  Thus in all cases $-1 < r_{12} \leq 0$, as desired.
\end{proof}
This allows us to compute the ring volume for $\xi_2$:
\begin{lem}\label{lem:M12_tame}
  If $\M_{12}$ is active, the solution volume for $\xi'_2$ is
  \[
  q^{-\floor{m_{12}}},
  \]
  except when $m_{12} \in \ZZ$ in splitting type $1^3$, in which case $\M_{12}$ has no solutions. 
\end{lem}
\begin{proof}
  By Lemma \ref{lem:M12_tame_tfm}, the $\M_{12}$-condition is given by one of the form
  \[
  \lambda^\diamondsuit(\alpha_1 \xi'_2) \equiv 0 \mod \pi^{\floor{m_{12}}}.
  \]
  Here $\alpha_1$ is primitive, so we have a linear relation modulo $\pi^{\floor{m_{12}}}$ which yields a solution volume of $(1 + 1/q) q^{-\floor{m_{12}}}$. The solutions must be further whittled down using the restrictions on $\xi'_{2}$ in Lemma \ref{lem:gamma_white} as well as the condition that $\xi_1$ and $\xi_2$ be linearly independent mod $\mm_{\bar K}$. Here it's important to note that the solutions to $\M_{12}$ are distributed equally among $(q + 1)$-many $1$-pixels.
  
  \begin{itemize}
    \item If $R$ is unramified and $[\delta \hat\omega_C] \in \L_0$, then $a_1' \in \ZZ$. One of the $1$-pixels is that of $\xi_1$, which violates the linear independence, so we eliminate it.
    \item If $R$ is unramified, $[\delta \hat\omega_C] \in (1;\pi;\pi)\L_0$, and
    \[
      (a_1, a_2, a_3) \equiv \(\frac{1}{2}, 0, \frac{1}{2}\) \mod 1,
    \]
    then $m_{12} \in \ZZ + 1/2$, and the condition that $\pi \nmid (\xi'_2)^{(K)}$ eliminates one $1$-pixel (no more, because $(\xi'_1)^Q$ is primitive).
    \item If $R$ is unramified, $[\delta \hat\omega_C] \in (1;\pi;\pi)\L_0$, and
    \[
    (a_1, a_2, a_3) \equiv \(\frac{1}{2}, \frac{1}{2}, 0\) \mod 1,
    \]
    then $m_{12} \in \ZZ$ and $\gamma_1 = \gamma_2 \sim (\sqrt{\pi} ; 1 ; 1)$. When the $\M_{12}$-condition
    \[
      \tr(\gamma_1^2 \xi_1' \xi_2') \equiv 0 \mod \pi^{m_{12}}
    \]
    is looked at mod $\pi$, it uniquely determines $\xi_2'^Q \equiv \xi_1'^Q \mod \pi$ (by $\M_{11}$, the value $\xi_2' = \xi_1'$ is a solution). But then $\xi_2 \equiv \xi_1 \mod \sqrt{\pi}$, violating linear independence. So $\M_{12}$ is unsatisfiable in this case.
    \item If $R$ is totally ramified, then the condition that $\xi'_2$ be a unit eliminates one $1$-pixel, unless all its solutions are non-units. This happens exactly when $\alpha_1 \sim \pi_R^2$, which is seen to be equivalent to $m_{12} \in \ZZ$.
  \end{itemize}
  Thus, in all but the stated exceptional case, we eliminate one $1$-pixel, leaving a ring volume of $q^{-\floor{m_{12}}}$.
\end{proof}

\subsection{\texorpdfstring{$\M_{22}$}{M22}}
When $\M_{22}$ is active, of course $\M_{11}$ is also, and we normalize both $\gamma_1$ and $\gamma_2$ (that is, $\xi'_1$ and $\xi'_2$) according to Lemma \ref{lem:tfm_conic}. In tame splitting types, we find that the conic
\[
  \M : \lambda^\diamondsuit (\delta^\odot \xi_i'^2) \equiv 0 \mod m'_{ii}
\]
is actually the same conic, but that the $\xi'_i$ cannot even lie in the same $1$-pixel. The following two lemmas detail when this can happen.

\begin{lem}\label{lem:A_M22_1}
Let $\A$ be a conic of determinant $1$ on a lattice $\Lambda$ over $\OO_K$, and let $m$ be an integer, $1 \leq m \leq e$. Then there are coprimitive $\vec{x}_1, \vec{x}_2 \in \PP(\Lambda)$ satisfying
\begin{equation}\label{eq:A_M22_1}
  \A(\vec{x}_1) \equiv \A(\vec{x}_2) \equiv 0 \mod \pi^m
\end{equation}
if and only if the squareness $\square(\A)$ satisfies
\[
  \square(\A) \geq m.
\]
Moreover, if $\square(\A) + 1 \geq m$, then for fixed $\vec{x}_1$, the volume of $\vec{x}_2$ satisfying \eqref{eq:A_M22_1} and coprimitive to $\vec{x}_1$ is $q^{-\ceil{m/2}}$, split evenly among $q$-many $1$-pixels.
\end{lem}
\begin{proof}
If coprimitive $\vec{x}_1, \vec{x}_2$ satisfy \eqref{eq:A_M22_1}, we can complete them to a basis $(\vec{x}_1, \vec{x}_2, \vec{x}_3)$ of $\Lambda$. Note that $\pi \nmid \A(\vec{x}_3)$, or else the determinant could not be $1$. So we may scale so that $\A(\vec{x}_3) = 1$, and now we see that $\A$ is a square modulo $\pi^m$.

Conversely, suppose that $\A$ is a square $\lambda^2$ of a linear form modulo $\pi^m$. Then for $\vec{x} \in \PP(\Lambda)$,
\begin{align*}
  \A(\vec{x}) &\equiv 0 \mod \pi^m \\
  \iff \lambda(\vec{x})^2 &\equiv 0 \mod \pi^m \\
  \iff \lambda(\vec{x}) & \equiv 0 \mod \pi^{\ceil{m/2}}.
\end{align*}
Since $\lambda \nequiv 0 \mod \pi$, this has solution volume $(1 + 1/q)q^{-\ceil{m/2}}$, split evenly among $(q+1)$-many $1$-pixels. If $\vec{x}_1$ is given, then $\vec{x}_2$ can occupy any $1$-pixel except the one containing $\vec{x}_1$.
\end{proof}

The following lemma limits $m_{22}$:
\begin{lem}\label{lem:m22<=e} The conditions on $\xi'_2$ can be satisfied only if $m_{22} \leq e$ and
  $m_{22}^{\odot} \leq e$.
\end{lem}
\begin{proof}
For the first part, note that if $m_{22} > e$, then
\[
  m_{12} = \frac{m_{11} + m_{22}}{2} - e + t > \frac{e + e}{2} - e = 0,
\]
so $\M_{12}$ is active, contradicting Lemma \ref{lem:inactives}. In unramified splitting type, this is the entire content of the lemma, since $m_{22}^{\odot} = m_{22}$. In splitting type $1^3$, note that
\[
  m_{11}^\odot = m_{11} - \frac{2h_1}{3} \textand
  m_{22}^{\odot} = m_{22} - \frac{2h_2}{3}
\]
satisfy $m_{11}^\odot + m_{22}^{\odot} = m_{11} + m_{22}$, since $\{h_1, h_2\} = \{1, -1\}$. Also, $m_{11}^\odot > m_{22}^{\odot}$ since $m_{11} - m_{22} = 2(a_2 - a_1) \geq 2/3$. So if $m_{22}^{\odot} > e$, then $m_{11}^\odot > e$ and $\M_{12}$ is active as above.
\end{proof}

This enables us to compute the ring volume for $\xi'_2$.
\begin{lem}\label{lem:M22}
Let $R$ be tamely ramified.
Fix the discrete data of a quartic ring and a $\xi'_1$ satisfying its conditions. The $\M_{22}$-condition is solvable for $\xi'_2$ if and only if the following conditions are satisfied:
\begin{itemize}
  \item $[\delta \hat\omega_C] \in \L_0$;
  \item the value of $a_2$ mod $\ZZ$ allows for a $\gamma_2$ and $m_{22}^{\odot}$ according to Lemma \ref{lem:tfm_conic}, and
  \item $m_{22}^{\odot} \leq \min\{2 \ell(\M) + 1, e\}$.
\end{itemize}
In such cases, the volume of $\xi'_2$ is 
\begin{itemize}
  \item $q^{1 - \ceil{m_{22}^{\odot}/2}}$ if $R$ is totally ramified, $h_1 = 1$ and $h_2 = -1$, 
  \item $q^{-\ceil{m_{22}^{\odot}/2}}$ otherwise.
\end{itemize}
\end{lem}
\begin{proof}
The necessity of the restrictions on $m_{22}^{\odot}$ is shown by the foregoing lemmas. If they are satisfied, then $\M_{22}$ transforms to a linear condition with solution volume $(1 + 1/q)q^{-\ceil{m_{22}^{\odot}/2}}$, distributed equally among $(q+1)$-many $1$-pixels. We must check its solutions against the other restrictions on $\xi'_2$:
\begin{itemize}
  \item If $R$ is unramified and $[\delta \hat\omega_C] \in \L_0$, then $\xi'_1$, $\xi'_2$ are required to be coprimitive, eliminating one of the $1$-pixels.
  \item If $R$ is unramified and $[\delta \hat\omega_C] \in (1;\pi;\pi)\L_0$ then $j_0$ must be $3$ in order for both $\gamma_1$ and $\gamma_2$ to exist. Then since $m_{22}^{\odot} = m_{22} \geq 1$, we must have two solutions $\xi'_1$, $\xi'_2$ to the transformed conic $\M$ modulo $\pi$ whose $Q$-components are coprimitive. But $\M(\xi') \mod \pi$ only depends on the $Q$-component and has only a unique solution in $\PP(\OO_Q/\pi \OO_Q)$, so $\M_{22}$ is unsatisfiable if active.
  \item If $R$ has splitting type $1^3$, then for $\gamma_1, \gamma_2$ to exist, we must have $\{h_1, h_2\} = \{1, -1\}$. Here, coprimitivity between the $\xi_i$ is subsumed by the condition that each $\xi^\odot_i$ lie in its correct domain
  \[
    \xi^\odot_i \sim \pi_R^{1-h_i}.
  \]
  As we noted in the proof of Lemma \ref{lem:eta_1^3}, the solutions to $\M$ mod $\pi$ comprise $q$-many $1$-pixels of $\xi^\odot \sim 1$ and one $1$-pixel of $\xi^\odot \sim \pi_R^2$. Hence if $h_1 = -1$ and $h_2 = 1$, we retain $q$ of the $(q+1)$-many $1$-pixels, getting a volume $q^{-\ceil{m_{22}^{\odot}/2}}$. But if $h_1 = 1$ and $h_2 = -1$, then $\xi^\odot_2$ is restricted to one $1$-pixel. This gives a volume of $q^{-1-\ceil{m_{22}^{\odot}/2}}$, but we multiply back by $q^2$ since $\xi_2' = \pi_R^{-2} \xi_2^\odot$. \qedhere
\end{itemize}
\end{proof}
Of the conditions, only
\begin{equation}\label{eq:M22}
m_{22}^{\odot} \leq \min\{2 \ell(\M) + 1, e\} 
\end{equation}
is not trivial to verify. The following solves it:
\begin{lem}
Suppose that the discrete data is fixed in such a way that
\begin{itemize}
  \item $\M_{22}$ is active,
  \item $[\delta \hat\omega_C] \in \L_0$ (so the conic has determinant $1$),
  \item the value of $a_2$ mod $\ZZ$ allows for a $\gamma_2$ and $m_{22}^{\odot}$ according to Lemma \ref{lem:tfm_conic}.
\end{itemize}
Also suppose that $\xi_1$ is fixed, satisfying the conditions $\M_{11}$, $\N_{11}$ governing it. Then the remaining condition \eqref{eq:M22} can be checked as follows:
\begin{itemize}
  \item In the black, purple, and blue zones, it is automatic.
  \item In the green zone, it is equivalent to
  \[
    m_{22}^{\odot} \leq \square_C
  \]
  \item In the red, yellow, and beige zones, it restricts the sum to only use terms $G^\cross (\ell_0,\ell_1,\ell_2)$ with
  \[
    m_{22}^{\odot} \leq 2 \ell_0 + 1.
  \]
\end{itemize}
\end{lem}
\begin{proof}
When the conic $\M$ is green, its level was computed as part of the finding of the zone total for $\xi'_1$. So it remains to prove that if the conic is black or blue, \eqref{eq:M22} is satisfied. That the conic is black or blue implies that
\[
  n_{11} \geq 2e - 4\ell(\M) - 1.
\]
We already know $m_{22} \leq e$. Suppose that
\[
  m_{22} \geq 2 \ell(\M) + 2.
\]
Since $m_{12}$ is inactive,
\[
  m_{11} = 2 m_{12} - m_{22} + 2e - 2t \leq 2e - 2\ell(\M) - 2.
\]
But then
\[
  n_{22} = m_{22} - m_{11} + n_{11} \geq 3,
\]
so $\N_{22}$ is active and we have a contradiction.
\end{proof}
\begin{rem}
In the code, $o_2$ and $h_2$ are defined by
\[
  o_2 = 3(b_2 \bmod 1) \in \{0,1,2\}, \quad h_2 = 2o_2 - 3 \in \{-3, -1, 1\}.
\]
The condition
\[
  \Dot m_{22} \leq 2k + 1
\]
(where $k = \ell_C$ in the green zone, or $k$ is the index of summation in the red zone) is coded as
\[
  2k + \frac{1}{3} - m_{22} + \frac{2o_2}{3} \geq 0.
\]
One verifies that this is the same thing when $o_1 \in \{1, 2\}$, while when $o_1 = 0$, it reduces to $m_{22} \leq 0$, as desired, since $\M_{22}$ is unsatisfiable if active.
\end{rem}

\section{Further remarks on the code}
In the attached code, we use the computer programs SAGE and LattE to compute the generating function of rings. First we count ``zone tuples'' consisting of integer values of the following variables:
\begin{itemize}
  \item $e, t, b_1, b_2, s$.
  \item $\verb|a1f| = \floor{a_1}, \verb|a2f| = \floor{a_2}, o_1, o_2, o_3$. Here we've decomposed
  \[
    a_i = \floor{a_i} + \frac{o_i}{o},
  \]
  where 
  \[
    o = \begin{cases}
    1 & \text{in splitting types $(111)$, $(12)$, and $(3)$} \\
    3 & \text{in splitting type $(1^3)$}
\begin{wild}
     \\ 2 & \text{in splitting type $(1^21)$.}
\end{wild}
\end{cases}.
  \]
  The $o_i$ belong to one of a finite number of ``flavors'' coding the classes of the $a_i$ and $b_i$ mod $\ZZ$.
  
  Note that $\floor{a_3}$ is missing from the variable list, as its value is uniquely determined by the discriminant identity (Lemma \ref{lem:rsv})
  \item 
  \[
    \verb|lCf| = \begin{cases}
      -1, & \square_C = 0\\
      \ell_C, & 0 < \square_C < e \\
      \floor{\dfrac{e}{2}}, & \square_C = e.
    \end{cases}
  \]
  This unambiguously determines $\square_C$. We note that in our answers, $\square_C$ only appears in the answer and bounds of the green zone. When $\square_C = e$, there is no green zone and its bounds are far from being achieved, so we can equate $\floor{\ell_C}$ with $\verb|lCf|$ without changing anything.
  \item $k$, an index of summation needed to input the answers in the red, yellow, and beige zones.
  \item If $\floor{m_{11}/2}$ appears in a zone answer, we add a variable $\verb|m11fl|$ and impose one of the equations
  \[
    m_{11} = 2 \cdot \verb|m11fl| \textor m_{11} = 2 \cdot \verb|m11fl| + 1.
  \]
  Similarly we treat $\floor{n_{11}/2}$, $\floor{n_c/2}$, $\tilde{n}$, $\floor{m_{22}/2}$, $\floor{e/2}$, $\floor{s/2}$.
\end{itemize}
The result comes out as a rational function in $\verb|RRR| = \ZZ((E, T, B_1, B_2, S, \verb|A1F|,\ldots))$, Here we use the convention that the value of a lowercase variable appears as an exponent of the corresponding uppercase variable. For instance, the generating function of the three lattice points on the line segment $y = 2x$, $0 \leq x \leq 2$ is
\[
  1 + XY^2 + X^2Y^4.
\]

Because each variable is bounded below in terms of the preceding variables, the power series is formally convergent. Because each zone is delimited by finitely many linear inequalities with $\ZZ$-coefficients, the generating function is a rational function, computed by Barvinok's algorithm as implemented in LattE. We then encode ring totals as substitutions that land us in a common ring
\[
  \verb|RINGS_RING| = \ZZ[[\verb|E_|,\verb|T_|,\verb|LCF_|,\verb|SFL_|,\verb|B1_|,\verb|B2_|,q,\verb|L0|,\verb|L2|]].
\]
The trailing underscores are to prevent the computer from confounding certain elements of $\verb|RRR|$ and $\verb|RINGS_RING|$, although the reader can think of them as identified. Two new variables $\verb|L0|,\verb|L2|$, whose exponents are the $\ell_0$ and $\ell_2$ of a boxgroup, complete the description of a ring total, along with a string \verb|Ftype| that tells the kind of weighting ($F$, $Fxx$, $E_{bal}$, etc.).

\begin{old}
\subsection{[TODO Old] The program \texttt{latte.sage}}

The remainder of the theorem can be proved automatically. To carry out the computations, we use Sage and the interface to LattE \cite{LattE}. LattE is a specialized tool that, given a convex polyhedral region
\[
\R = \{\vec{v} \in \RR^n : A \vec{v} + \vec{b} \geq 0\}
\]
given by a finite list of linear inequalities with rational coefficients, can compute the $n$-dimensional generating function for the lattice points in $\R$:
\[
F(X_1, \ldots, X_n) = \sum_{(x_1,\ldots,x_n) \in \R \intsec \ZZ^n} X_1^{x_1} \cdots X_n^{x_n}.
\]
The generating function $F$ will be a rational function if it converges, which happens for instance if all coordinates $x_i$ are bounded below.

We apply LattE to the variables $i_1, i_2, i_3, \ell_1, \ell_2, r$ and the inequalities defining each zone. The following adaptations are to be noted:
\begin{itemize}
  \item Since LattE counts only integer points, we decompose $i_k = i_{k,\text{floor}} + i_{k,\text{type}}/2$, where $i_{k,\text{floor}}$ and $i_{k,\text{type}}$ are integers and $0 \leq i_{k,\text{type}} \leq 1$. The three variables $i_{k,\text{type}}$ thus encapsulate the type of any quartic order:
  \[
  \begin{tabular}{cccc}
  Type & $i_{1,\text{type}}$ & $i_{2,\text{type}}$ & $i_{3,\text{type}}$ \\ \hline
  ur & $0$ & $0$ & $0$ \\
  1  & $0$ & $1$ & $1$ \\
  2  & $1$ & $0$ & $1$ \\
  3  & $1$ & $1$ & $0$
  \end{tabular}
  \]
  In addition, we create three variables $\verb|twoi|_k$, $1 \leq k \leq 3$, defined by $\verb|twoi|_k = 2i_k = 2i_{k,\text{floor}} + i_{k,\text{type}}$.
  \item The ring totals in the zones other than IV and V are simply exponentials where $q$ is raised to a power linear in the $i_k$, $\ell_k$, and $r$. Thus they can be evaluated by plugging the appropriate powers of $q$ for the corresponding variables in the generating function. We define $\verb|rq| = \sqrt{q}$ to deal with fractional powers of $q$ appearing in intermediate steps of the calculation.
  \item The ring totals in zones IV and V are not quite pure exponentials: they depend on the parity of $b_{11}$. Since $b_{11} = \ell_{1} - 2i_1$, this is equivalent to the parity of $\ell_{1}$. We therefore let $\ell_{1} = 2\ell_{1,\text{floor}} + \ell_{1,\text{parity}}$ with $0 \leq \ell_{1,\text{parity}} \leq 1$.
  \item Likewise, we deal with the restriction on the parity of $r$ in Zones VI and VII.
  \item Only the zone supports for even $r$ are entered by hand, since the zone supports for odd $r$ are generated by flipping about a vertical axis.
\end{itemize}
The program finishes by printing two rows of zeros, indicating that Theorem \ref{thm:O-N_quartic_local_tame} has been verified for both even and odd $r$.

\end{old}

\part{Unanswered questions} \label{part:unans}
\begin{todofn}
\section{Extension to function fields}

We hope that there are extensions of the reflection theorem for $3$-torsion found in Lee \cite{LeeRefl}.

(Old: TODOfn needs review)

We carry out the analogous program to Tate's thesis where the additive group $K$ is replaced by a cohomology group $H^1(K, M)$. For instance, if $M \cong \mu_m$, we are looking at a quotient $H^1(K, M) \cong K^\cross/(K^\cross)^m$ of the multiplicative group. Fortunately, all the theoretical groundwork has been laid for us by \v Cesnavi\v cius \cite{Ces_Topo, Ces_Poitou}.

Let $K$ be a local field, and let $M$ be a commutative finite flat group scheme over $K$, such as a Galois module over $K$. The pairing needed for step \ref{it:Tate_local} is the \emph{Tate pairing}
\[
\<\bullet, \bullet\> : H^1(K,M) \cross H^1(K,M') \to \QQ/\ZZ,
\]
where $M' = \underline{\Hom}(M, \GG_m)$ is the arithmetic dual. If $\ch K \nmid \size{M}$, then $H^1(K,M)$ and $H^1(K,M')$ are finite abelian groups of the same order, but in general, they are Pontryagin duals as Hausdorff locally compact abelian topological groups (\cite{Ces_Topo}, \textsection 1.1). 

Next, let $K$ be a global field, $M$ a commutative finite flat group scheme over $K$, and $\M$ a \emph{model} for $M$; that is, a commutative finite flat group scheme over the ring $K_{\Sigma_0}$ of ${\Sigma_0}$-integral elements, for some finite set $\Sigma_0$ of places, equipped with an identification $\M \cross_{K_\Sigma} K \cong M$. Such a model exists by the general technique of ``spreading out.'' This is a technical step, and the reader is invited to think mainly of the case in which $M$ and its dual are Galois modules.

Then \v Cesnavi\v cius proves that the integral cohomology $H^1(\OO_v, \M)$ is a compact open subgroup of $H^1(K, M)$, which one may think of as generalizing the \emph{unramified} cohomology $H^1_\ur(K,M)$; and that the resulting adelic cohomology
\[
H^1(\AA_K, M) = \sideset{}{'}\prod_v H^1(K_v,M) = \bigcup_{\Sigma_0 \supseteq \Sigma} \Big( \prod_{v \in \Sigma} H^1(K_v, M) \cross \prod_{v \notin \Sigma} H^1(\OO_v, \M) \Big)
\]
is independent of $\M$. It is clear that $H^1(\AA_K, M)$ is locally compact. So, after recalling that the integral cohomology $H^1(\OO_v, \M)$ is dual to $H^1(\OO_v, \M')$ in the Tate pairing for $v \notin \Sigma_0$ (\cite{Ces_Poitou}, remarks in the proof of Proposition 4.10(b)), we can multiply the Tate pairings for all $v$ to get a perfect pairing
\[
\<\bullet, \bullet\> : H^1(\AA_K, M) \cross H^1(\AA_K, M') \to \QQ/\ZZ.
\]
Fix a scaling for the Haar measure on $H^1(\AA_K, M)$. Then every locally constant, compactly supported function (hereafter a \emph{step function}) $f : H^1(\AA_K, M) \to \CC$ has a Fourier transform
\[
\hat f : H^1(\AA_K, M') \to \CC, \quad \hat f(\beta) = \int_{\alpha \in H^1(\AA_K, M)} f(\alpha) e^{-2\pi i \<\alpha, \beta\>} \, d\alpha
\]
that is also a step function.

The Poisson summation formula will be applied to this class of functions.

\begin{thm}[\textbf{Poisson summation for adelic cohomology}] \label{thm:Poisson}
  Let $K$ be a global field, $M$ a commutative finite flat group scheme over $K$, and $\M$ a model for $M$. Fix a scaling for the Haar measure on $H^1(\AA_K, M)$.There exists a constant $c > 0$ such that for every step function $f : H^1(\AA_K, M) \to \CC$,
  \begin{equation} \label{eq:Poisson}
    \sum_{\alpha \in H^1(K, M)} f(\alpha) = c \sum_{\beta \in H^1(K, M')} \hat f(\beta).
  \end{equation}
\end{thm}
\begin{proof}
  The natural map $\loc_M : H^1(K, M) \to H^1(\AA_K, M)$, which we are implicitly using, has discrete image and finite kernel (\cite{Ces_Poitou}, Proposition 4.12 and Lemma 4.4(b)). So, up to changing the constant $c$, \eqref{eq:Poisson} is equivalent to
  \[
  \sum_{\alpha \in \im(\loc_M)} f(\alpha) = c \sum_{\beta \in \im(\loc_{M'})} \hat f(\beta).
  \]
  But since $\im(\loc_M)$ and $\im(\loc_{M'})$ are dual lattices, by 
  Tate duality (\cite{Ces_Poitou}, Proposition 4.10(c)), this literally follows from Poisson summation.
\end{proof}

\subsection{Tame function-field reflection}

Given a stock of local weightings and their Fourier transforms, we can apply Poisson summation (Theorem \ref{thm:Poisson}) to derive global reflection theorems. We refrain from stating results over $\QQ$ and other number fields, of which O-N is the prime example, until we have proved the appropriate results at the prime $3$ in the next section. However, we can now state an analogue of O-N for function fields.

TODOfn state it ($q \equiv 2 \mod 3$)

\subsection{The scale factor \texorpdfstring{$c$}{c}}

The constant $c$ is the covolume of the lattice of global cohomology elements, adjusted by the orders of the kernels of $\loc_M$ and $\loc_{M'}$. The method of proof of \cite{Ces_Poitou} does not lend itself easily to computing $c$.

To even ask what $c$ is, one must somehow write down a measure on $H^1(\AA_K, M)$, or on all the groups $H^1(K_v, M)$. For a general $M$, this seems quite hard. The easiest case occurs when $M \cong M'$ is self-dual (for instance, when $K$ has characteristic $p$ and $M = \alpha_p$, the kernel of Frobenius). Here, taking the unique self-dual measure on $H^1(\AA_K, M)$, or equivalently the product of the unique self-dual measure on each $H^1(K_v, M)$, we see from symmetry considerations that $c = 1$.

One other case in which something can be said is the following.

\begin{conj} \label{conj:Poisson_c}
  Suppose $M$ is \'etale, that is, $M$ is a Galois module over $K$. Give $H^1(K_v, M)$ (which is now discrete) the measure for which each single point has measure $1/\size{H^0(K_v,M)}$. Then with respect to the measure thereby induced on $H^1(\AA_K, M)$,
  \[
  c = \frac{\size{H^0(K, M)}}{\size{H^0(K, M')}}.
  \]
\end{conj}
We note that the product measure on $H^1(\AA_K, M)$ converges because, if $\M$ is a model for $M$, then for all but finitely many $v$, we have $\size{H^1(\OO_v, \M)} = \size{H^0(\OO_v, \M)} = \size{H^0(K_v, M)}$.

\begin{thm} \label{thm:Poisson_c} Conjecture \ref{conj:Poisson_c} holds when:
  \begin{enumerate} [(a)]
    \item \label{it:Poisson_tame} $\ch K \nmid \size{M}$; or
    \item \label{it:Poisson_split} $M$ is split, that is, $M$ is a module with trivial Galois action.
  \end{enumerate}
\end{thm}
\begin{proof}
  \begin{enumerate}[(a)]
    \item See for instance \cite{DDT}, Theorems 2.18--19(?), where the authors deduce from Poitou-Tate duality a theorem equivalent to Theorem \ref{thm:Poisson} for $f$ the characteristic function of a compact open box $\L = \prod_v \L_v \subseteq H^1(\AA_K, M)$ (a \emph{Selmer system} on $H^1(\AA_K, M)$). The constant factor is computed using the Euler-characteristic formula
    \[
    \frac{\size{H^0(K_v,M)} \size{H^2(K_v,M)}}{\size{H^1(K_v,M)}} = \Big\lvert\size{M}\Big\rvert_v.
    \]
    We expect the method to generalize readily to the function-field case.
    \item We can assume that $K$ is a function field, the number-field case being covered by \ref{it:Poisson_tame}. We can assume that $M = \ZZ/m\ZZ$ is cyclic. Note that since $M$ is split, it has a global model $\M$, namely $\ZZ/m\ZZ$ considered as a group scheme over the curve $C(K)$, with Cartier dual $\M' = \mu_m$. We take
    \[
    f = \1_{H^1(\OO_{\AA_K}, \ZZ/m\ZZ)} \textand \hat f = \1_{H^1(\OO_{\AA_K}, \mu_m)}
    \]
    for the purposes of computing
    \[
    c = \frac{\size{H^1(\OO_K, \ZZ/m\ZZ)}}{\size{H^1(\OO_K, \mu_m)}}.
    \]  
    Now
    \[
    H^1(K, \ZZ/m\ZZ) = \Hom(\Gal(K^{\mathrm{sep}}/K), \mu_m
    \]
    parametrizes degree-$m$ cyclic extensions of $K$, and the integral cohomology subgroup $H^1(\OO_K, \ZZ/m\ZZ)$ parametrizes those that are unramified. By global class field theory, these are in bijection with characters $\chi : \Pic(K) \to \mu_m$ of the Picard group. Now $\Pic(K) \cong \ZZ \oplus \Pic^0(K)$ with $\Pic^0(K)$ finite, so the number of such $\chi$ is
    \[
    m \cdot \size{\Pic^0(K) / (\Pic^0(K))^m}.
    \]
    On the other hand, $H^1(K, \mu_m)$ is parametrized by $K^\cross/(K^\cross)^m$, and the integral cohomology subgroup $H^1(\OO_K, \mu_m)$ is the subgroup of functions whose divisor is a multiple of $m$, that is, the kernel of the map $\operatorname{div} : K^\cross/(K^\cross)^m \to \Div^0(K)/m\Div^0(K)$. Applying the Snake Lemma to the diagram
    \[
    \xymatrix{
      0 \ar[r] & K^\cross/\mu_m(K) \ar[r]^{\bullet^m} \ar[d] & K^\cross \ar[r] \ar[d] &  K^\cross/(K^\cross)^m \ar[r] \ar[d] & 0 \\
      0 \ar[r] & \Div^0(K) \ar[r]^{m\bullet} & \Div^0(K) \ar[r] & \Div^0(K)/m\Div^0(K) \ar[r] & 0
    }
    \]
    shows that 
    \[
    \size{H^1(\OO_K, \mu_m)} = \size{\mu_m(K)} \cdot \size{\Pic^0(K) / (\Pic^0(K))^m}.
    \]
    Hence
    \[
    c = \frac{\size{H^1(\OO_K, \ZZ/m\ZZ)}}{\size{H^1(\OO_K, \mu_m)}} = \frac{m}{\size{\mu_m(K)}} = \frac{\size{H^0(K, M)}}{\size{H^0(K, M')}}.
    \]
  \end{enumerate} 
\end{proof}

The most familiar case of Theorem \ref{thm:Poisson} is when $\ch K \nmid \size{M}$ and the adelic step function $f$ is the characteristic function of an intersection of local subgroups $\L_v \subseteq H^1(K, M_v)$. Then Poisson summation recovers the celebrated \emph{Greenberg-Wiles formula} for the ratio of the sizes of the Selmer groups corresponding to $\{\L_v\}$ and its orthogonal complement $\{\L_v^\perp\}$.

\subsection{Reflection on orbit counts}
A large number of identities can be derived by applying Poisson summation to various functions $ f $. In this paper, we will mainly take $f$ to be the \emph{number of integral orbits} for a representation whose rational orbits are parametrized by the appropriate cohomology group. We first define the situations in which we want to work in rather ungainly generality.
\begin{defn}
  Let $ K $ be a number field. A \emph{reflectible space} over $ K $ consists of the following data:
  \begin{itemize}
    \item An algebraic group $ \Gamma $ over $ \OO_K $;
    \item A finite-dimensional representation $ V $ of $ \Gamma $ over $ \OO_K $;
    \item A (not necessarily complete) invariant $ I $ of $ V $, that is, a polynomial map $ I $ from $ V $ to an affine space $ \AA^r $ such that the identity $ I(v) = I(\gamma v) $ ($ v \in V $, $ \gamma \in \Gamma $) holds formally;
    \item A specific value $ I_0 \in \AA^r(\OO_K) $;
    \item A finite Galois module $ M $ over $ K $;
    \item For each extension $ K'/K $, a bijection
    \[
    u_{K'} : H^1(K',M) \to \Gamma(K')\backslash V(K')_{I=I_0}
    \]
    that identifies the orbits with invariant $ I_0 $ with a Galois cohomology group;
    \item For each place $ v $ of $ K $, a local weighting
    \[
    w_v : \Gamma(\OO_{K_v}) \backslash V(K_v)_{I=I_0} \to \CC
    \]
    on the integral orbits
  \end{itemize}
  with the following properties:
  \begin{enumerate}[$($a$)$]
    \item The parametrization $u$ respects base change, that is, for extensions $ K''/K' $ and $ \alpha \in H^1(K',M) $, $ u_{K'}(\alpha) $ and $ u_{K''}(\Res_{K'}^{K''}\alpha) $
    belong to the same $ \Gamma(K'') $-orbit;
    \item For all $ K'/K $ and all $ x \in V(K') $, the stabilizer has order
    \[
    \size{\Stab_{\Gamma(K')} x} = \size{H^0(K', M)};
    \]
    \item The local weighting $ w_v $ vanishes on all but finitely many integral orbits within each rational orbit, so that the local class number
    \begin{align*}
      g_{v} : H^1(K_v, M) &\to \CC \\
      \alpha &\mapsto \sum_{\gamma \in \Gamma(\OO_{K_v})\backslash \Gamma(K_v)} w_{v}(\gamma v_\alpha)
    \end{align*}
    is finite;
    \item For all but finitely many $ v $,
    \[
    g_v = \1_{H^1_\ur(K_v, M)}.
    \]
  \end{enumerate}
\end{defn}

\begin{defn}
  The \emph{global class number} of a reflectible system is
  \[
  h = \sum_{x \in \Gamma(\OO_K)\backslash V(K)} \left(\frac{1}{\size{\Stab_{\Gamma(\OO_K)}} x} \prod_v w_v(x)\right)
  \]
\end{defn}

Before stating our main theorem regarding these complicated objects, we review some standard theory of algebraic groups and make some observations. The \emph{class number} of $ \Gamma $ is (TODOfn what?).

\begin{thm} \label{thm:loc_glo}
  Let $ K $ be a number field. Let
  \[
  (\Gamma^1, V^1, I^1, I_0^1, M^1, \{u^1_{K'}\}_{K'}, \{w_v^1\}_v) \textand (\Gamma^2, V^2, I^2, I_0^2, M^2, \{u^2_{K'}\}_{K'}, \{w_v^2\}_v)
  \]
  be two reflectible spaces over $ K $. Suppose that
  \begin{enumerate}
    \item $ \Gamma^1 $ and $ \Gamma^2 $ have class number $ 1 $ (TODOfn explain);
    \item $ M^2 $ is identified with the Tate dual $ (M^1)' = \Hom(M^1,\mu) $;
    \item At each place $v$, the local class numbers $ g_v^i : H^1(K_v, M^i) $ satisfy \emph{local reflection} 
    \[
    \widehat{g_v^1} = g_v^2.
    \]
  \end{enumerate}
  Then the global class numbers $ h^1 $, $ h^2 $ of the reflectible systems are equal.
\end{thm}

\subsection{Function-field reflection theorems via tilting}
When doing tame function field case, don't forget the following remark from GGS \cite{GGS}:

P. Deligne has observed that the bijection of orbits and rings established in Proposition 4.2 holds over any base scheme $S$. There is an equivalence of categories between the following two kinds of objects, with morphisms being the isomorphisms:

(a) a vector bundle $V$ of rank $2$ with $p$ in $\Sym^3(V)\tensor \Lambda^2(V)^{-1}$;

(b) a vector bundle $A$ of rank $3$ with a (commutative) algebra structure.

\label{sec:tilt}
\emph{Tilting} is a technique for transporting knowledge about local fields from mixed to pure characteristic (or, more rarely, vice versa) using the intuition that a characteristic-$p$ local field is a ``limit'' of extensions of $\QQ_p$ whose ramification indices tend to infinity. Here our application will be to produce O-N-style reflection theorems over \emph{wild} function fields, that is, those $K$ for which $\ch K \mid \size{M}$.

As in the number field case, the proof of O-N has two parts: local and global. The global part will proceed in essentially identical fashion using Theorem \ref{thm:Poisson}. It is the local part that we concern ourselves with here. 

We begin with the following generalization of Theorem \ref{thm:O-N_traced_local}:
\begin{thm}[\textbf{Traced O-N in the wild function-field case}]\label{thm:O-N_traced_local_wildfn}
  Let $K$ be a local field of residue field $k_K = \FF_q$. Let $V_1$ be the space of binary cubic forms, a representation of $\Gamma_1 := \SL_2$ with invariant
  \[
  \disc(ax^3 + bx^2y + cxy^2 + dy^3) = b^2c^2 - 4ac^3 - 4b^3d - 27a^2d^2 + 18abcd.
  \]
  Parametrize the $K$-orbits $V_1(K)^{\disc = D}$, via the corresponding cubic algebras, by the cohomology group
  \[
  U_{1,D} := H^1(K, M_D) \cong \begin{cases}
    K[\sqrt{-3D}]^{N=1}/\text{cubes} & \ch K \neq 3 \\
    K\sqrt{D}/\wp(K\sqrt{D}), & \ch K = 3.
  \end{cases}
  \]
  
  Then let $V_2$ be the space of symmetric trilinear forms over $K$, a representation of $\Gamma_2 := \SL_2$ with invariant
  \[
  \disc\( ax_1x_2x_3 + b \sum_{\sym} x_1x_2y_3 + c \sum_{\sym} x_1y_2y_3 + dy_1y_2y_3\)
  = -3b^2c^2 + 4ac^3 + 4b^3d + a^2d^2 - 6abcd.
  \]
  Parametrize the $K$-orbits $V_2(K)^{\disc = D}$, via the corresponding cubic algebras, by the flat cohomology group
  \[
  U_{2,D} := H^1(K, M_D \tensor \mu_3) \cong K[\sqrt{D}]^{N=1}/\text{cubes}.
  \]
  Note that $U_{1,D}$ and $U_{2,D}$ are Pontryagin duals under the local Tate pairing. Scale the measures on them so that TODOfn.
  
  For $0 \leq t \leq e := v_K(3)$, define
  \begin{align*}
    \Lambda_{1,t} &:= \{ax^3 + bx^2y + cxy^2 + dy^3 \in V_1 : a,d \in \OO_K, b,c \in \pi^t \OO_K\} \\
    \Lambda_{2,t} &:= \left\{ ax_1x_2x_3 + b \sum_{\sym} x_1x_2y_3 + c \sum_{\sym} x_1y_2y_3 + dy_1y_2y_3 \in V_2 : b,c \in \OO_K, a,d \in \pi^t \OO_K\right \}.
  \end{align*}
  Let, as above,
  \[
  g_i(\alpha, \Lambda_{i,t}, D) = \size{\{c \in \Gamma_i(\OO_K) \backslash \Gamma_i(K) : c v_{D,\alpha} \in \Lambda_{i,t}\}}
  \]
  where $v_{D,\alpha} \in V_i(K)$ is an arbitrary vector with discriminant $D$ in the orbit parametrized by the class $\alpha \in U_{i,D}$.
  
  Then $g_1$ and $g_2$ are Fourier duals:
  \[
  \hat{g}_2(\alpha, \Lambda_{2,t}, D) = q^t g(\alpha, \Lambda_{1,t}, \pi^{2t}D).
  \]
\end{thm}
\begin{proof}
  In the case that $\ch K \neq 3$, this is Theorem \ref{thm:O-N_traced_local}, with slight adaptations of notation. But in the case that $\ch K = 3$, $U_{1,D}$ is a countable discrete group while $U_{2,D}$ is an infinite profinite compact group. Thus, if we want to say that $U_i$ is a ``limit'' of the $U_i$ occurring in the characteristic-0 case, it will have to be a direct limit for $U_1$ and an inverse limit for $U_2$.
  
  Let $K = \FF_q((\pi))$ be given, where $q$ is a power of $3$, and let
  \[
  K_e = \QQ_q[\sqrt[e]{-3}],
  \]
  a field with uniformizer $\pi_e = \sqrt[e]{-3}$. The negative sign ensures that, for $e$ even, $K_e$ contains $\mu_3$. We define a map
  \begin{align*}
    \sharp \colon K &\to K_e \\
    \sum_{i} a_i \pi^i &\mapsto \sum_{i} \tilde a_i \pi_e^i
  \end{align*}
  where $\tilde a$ denotes the Teichm\"uller lift. As is customary in this field, the image $\sharp(f)$ is denoted by $f^\sharp$, and the inverse map $\sharp^{-1}$ is denoted by $\flat$.
  
  Note that $\sharp$ is a homeomorphism (indeed an isometry) of the two valued fields, but not a ring isomorphism. It is easy to see, however, that $\sharp$ is an isomorphism ``mod $\pi^e$'' in the sense that, for all $a,b \in \OO_K$,
  \begin{equation}
    \label{eq:hom mod e}
    \left.
    \begin{aligned}
      (a + b)^\sharp &\equiv a^\sharp + b^\sharp \\
      (ab)^\sharp &\equiv a^\sharp b^\sharp
    \end{aligned}
    \right\}
    \mod \pi_e^e.
  \end{equation}
  Our strategy will be to let $e \to \infty$ in relation to other parameters in the problem and exploit \eqref{eq:hom mod e} to show that $K$ and $K_e$ behave identically in the needed senses. We use the notation $e \gg_M 1$ to avoid repeating ``if $e$ is large enough, for fixed $m$.'' Let $\sharp$ and $\flat$ act coordinatewise on vectors and matrices. Equip $V_i$ and $\Gamma_i$ with the $\ell^\infty$-norm.
  
  The algebraic group $\Gamma = \SL_2$ is a bit awkward to use in conjunction with tilting, inasmuch as the condition that a determinant be exactly $1$ is not preserved under $\sharp$. Therefore we introduce, for $N \in \NN^+$, a ``thickened'' group $\Lambda_N$, defined by
  \begin{equation*}
    \Lambda_N(K) = \{\gamma \in \GL_2 K : \det \gamma \equiv 1 \mod \pi^N\}.
  \end{equation*}
  Let $\Lambda_N(\OO_K) = \Lambda_N(K) \intsec \GL_2 \OO_K$, and define $\Lambda_N(K_e)$ and $\Lambda_N(\OO_{K_e})$ analogously. Note that $\Lambda_N(\OO_{K_e}) = \Lambda_N(\OO_{K})^\sharp$. Also note that $\Lambda_N(K)$ is the direct product of $\Gamma_K$ with $U_N := 1 + \pi^N \OO_K$, and likewise for $\Lambda_N(\OO_K)$, $\Lambda_N(K_e)$, $\Lambda_N(\OO_{K_e})$.
  
  \begin{lem}
    \label{lem:tilt1}
    For all $v \in V_1(\OO_K)$, if $N \gg_{\size{v}} 1$ and $e \gg_{\size{v}, N} 1$,
    \begin{equation} \label{eq:x tilt1}
      \Lambda_N\(\OO_{K_e}\) v^\sharp = \( \Lambda_N(\OO_K) v \)^\sharp.
    \end{equation}
  \end{lem}
  \begin{proof}
    Let $q^M \geq \size{v}$. Note that since, for $\gamma \in \Lambda_N(\OO_K)$,
    \[
    (\gamma v)^\sharp \equiv \gamma^\sharp v^\sharp \mod \pi_e^{e-M},
    \]
    each point in the right-hand side of \eqref{eq:x tilt1} is in the same $(e-M)$-pixel as a point of the left-hand side, and conversely. So it suffices to show that both sides of \eqref{eq:x tilt1} are unions of $(e-M)$-pixels. Let us show that the pixel $v + \pi^{e-M}$ lies in the orbit $\Lambda_N v$. Note that $\Lambda_N$ contains the $N$-pixel
    \[
    P := I + \pi^N \Mat_2 \OO_K.
    \]
    We consider the map $\bullet v : P \to V_1$. It is not hard to compute that the differential of $\bullet v$ has determinant of size  $\size{D}$ at all points of $P$, so by a version of the inverse function theorem (which we could just as well verify explicitly), $\bullet v$ surjects onto the pixel $v + \pi^{N - v_K(D)}$. This shows that the right-hand side of \eqref{eq:x tilt1} is a union of $(e - M)$-pixels for $e \gg_{M,N,D} 1$, and the same proof applies to the left-hand side.
  \end{proof}
  \begin{rem}
    The above lemma is stated only for $V_1$. Attempting to prove it for $V_2$ runs into difficulties (the map $\bullet v$ has everywhere singular differential), as it should: the orbits are no longer open, and there are infinitely many of them, even for fixed $D$.
  \end{rem}
  
  \begin{lem}
    \label{lem:tilt cosets}
    Let $M > 0$ be given. If $e \gg_M 1$, and if $\gamma_i$ ($1 \leq i \leq n$) are a set of representatives for the cosets $\gamma\Gamma(\OO_K) \subset \Gamma(K)$ for which $\size{\gamma} \leq q^M$, then the cosets $\delta \Gamma(\OO_{K_e}) \subset \Gamma(K_e)$ for which $\size{\delta} \leq q^M$ have a set of representatives $u_i \gamma_i^\sharp$, where the $u_i$ are units.
  \end{lem}
  \begin{proof}
    Note that we cannot simply take $\gamma_i^\sharp$ as a coset representative, as its determinant may differ ever so slightly from $1$. However, since $\det \gamma_i^\sharp \equiv 1$ mod $\pi_e^{e - M}$, we can find a unique unit $u_i$, $u_i \equiv 1$ mod $\pi_e^{e - M}$ such that $\det (u_i \gamma_i^\sharp) = 1$. Now let $\delta \in \Gamma(K_e)$ be given, $\size{\delta} \leq q^M$. Then $\size{\delta^\flat} \leq q^M$ and there exists a unique $i$ such that $u\delta^\flat \in g_i \Gamma(\OO_K)$. It is not hard to verify that $\delta$ belongs to $u_i g_i^\sharp \Gamma(\OO_{K_e})$ and to no other $u_j g_j^\sharp \Gamma(\OO_{K_e})$.
  \end{proof}
  
  \begin{lem}
    \label{lem:tilt2}
    Let $M > 0$ and $t \in \NN$ be given. If $e \gg_M 1$, then for bounded $v \in V_i(K)$ with $\size{v} < q^M$ and $v_K(\disc v) < M$, then for all $\gamma \in \Gamma_i(\OO_K)$,
    \begin{equation} \label{eq:tilt2}
      \gamma v \in \Lambda_{i,t}(\OO_K) \iff \gamma^\sharp v^\sharp \in \Lambda_{i,t}(\OO_{K_e}).
    \end{equation}
    Moreover, both sides of \eqref{eq:tilt2} are false if $\size{\gamma} \gg_M 1$.
  \end{lem}
  \begin{proof}
    Note that for any $N > 0$, if $\size{\gamma} < q^N$, then
    \begin{equation} \label{eq:sharp_and_act}
      \gamma^\sharp v^\sharp \equiv (\gamma v)^\sharp \mod \pi_e^{e - N - M}.
    \end{equation}
    Since $\Lambda_{i,t}(\OO_K)^\sharp = \Lambda_{i,t}(\OO_{K_e})$, this implies the desired conclusion as soon as $e > M + N + t$. We now claim that $\size{\gamma}$ is bounded for all $\gamma$ for which either side of \eqref{eq:tilt2} holds. Note that we may multiply $v$ by a sufficiently high power of $\pi$ so that it is integral, as this only weakens the conditions \eqref{eq:tilt2}. We do this using the algebraic interpretation of points of $V_i$:
    \begin{itemize}
      \item If $i = 1$, then $v$ corresponds to a nondegenerate cubic ring $C_1$ in a maximal ring $C_0$, and $\gamma v$, if it is integral, corresponds to a ring $C_2$ of the same discriminant. We can write $\gamma = \gamma_2 \gamma_1^{-1}$ where $\gamma_i$ is a change-of-basis matrix to pass from $C_0$ to $C_i$. Then $\gamma_i$ is integral since $C_i \subseteq C_0$, and we have
      \[
      \size{\gamma} = \size{\gamma_2 \gamma_{1}^{-1}} \leq \size{\gamma_1^{-1}} \leq [C_0 : C_1] \leq q^{v(\disc C_1)} \leq q^M.
      \]
      \item If $i = 2$, then $v$ corresponds to a self-balanced ideal $(\OO_D, I_1, \delta)$ with $\delta I_1^3 \subseteq \OO_D$. We may scale $\delta$ so that $\delta^{-1} \in \OO_{D_0}$, the maximal order in the field containing $\OO_D$. Then $I_1 \subseteq \OO_{D_0}$ because $\OO_{D_0}$ is integrally closed. Now if $\gamma v$ is integral, it corresponds to another self-balanced ideal $(\OO_D, I_2, \delta)$. We can write $\gamma = \gamma_2 \gamma_1^{-1}$ where $\gamma_i$ is a change-of-basis matrix to pass from $\OO_{D_0}$ to $I_i$. Then $\size{\gamma} \leq q^M$ similarly to the preceding case.
    \end{itemize}
    The proof over $K_e$ is identical.
  \end{proof}
  
  The Hilbert pairing is stable under tilting in the following sense: 
  \begin{lem} \label{lem:tilt Hilbert}
    Let $M > 0$ be given. If $e$ is even and sufficiently large compared to $M$, then for all fundamental discriminants $D_0 \in \OO_K$, all $\beta \in K[\sqrt{D_0}]^{N=1}$ and all $\alpha \in K\sqrt{D_0}/\wp(K\sqrt{D_0})$ with $\size{\alpha} < q^M$,
    \[
    (\beta, \alpha]_K = (\beta^\sharp, 1 - 3 \sqrt{-3} \alpha^\sharp)_{K_e}.
    \]
  \end{lem}
  \begin{proof}
    This follows easily from the theory of explicit reciprocity laws. See Sen (\cite{SenI}, Theorem 3) for a result easily strong enough to compute the right-hand side. The left-hand side, as is well known (\cite{SerreLF}, Prop.~XIV.15), is computable from the residue of a local differential
    \[
    (\beta, \alpha]_K = \tr_{k_K/\FF_3} \Res(\alpha\, d\beta/\beta)
    \]
  \end{proof}
  
  Now fix $\beta \in K[D]^{N=1}/\text{cubes}$. Given that for all $e$,
  \[
  g_2(\beta^\sharp, \Lambda_{2,t}, D^\sharp) = q^t \sum_{\alpha \in U_{1,D}(K_e)} (\alpha, \beta) \cdot g_1(\alpha, \Lambda_{1,t}, \pi_e^{2t}D^\sharp),
  \]
  we wish to deduce the corresponding identity over $K$. First we claim that
  \begin{equation} \label{eq:x tilt h2}
    g_2(\beta, \Lambda_{2,t}, D) = g_2(\beta^\sharp, \Lambda_{2,t}, D^\sharp)
  \end{equation}
  when $e \gg_D 1$. Note that the class $\beta$ has a representative
  \[
  v_\beta = x_1 x_2 x_3 - 3 \sum_{\sym} x_1 y_2 y_3 - (\tr \beta) y_1 y_2 y_3
  \]
  for which $\size{v_\beta}$ is bounded (indeed by $1/q$) and $v_\beta^\sharp$ also represents $\beta^\sharp$. The conclusion \eqref{eq:x tilt h2} now follows from Lemmas \ref{lem:tilt cosets} and \ref{lem:tilt2}.
  
  With $g_1$, we use a bit more indirection. First recall that, by (TODO the Artin-Schreyer analogue of Theorem \ref{thm:disc_Kummer_aff}), the cubic orders over $\OO_K$ of discriminant $D$ have corresponding Artin-Schreyer elements $s\sqrt{D} \in K\sqrt{D_0}$ with $\size{s} \leq q^M$, where $M$ depends only on $D$. A representative vector of this class is
  \begin{equation}
    v_{s} = D x^3 - x y^2 - s y^3;
  \end{equation}
  note that $\disc v_{s} = D$ and $\size{v_s} \leq q^M$. By Lemma \ref{lem:tilt1}, the orbits of $\Lambda_N(\OO_K)$ on vectors of approximate discriminant $D(1 + \pi^N\OO_K)$ are in bijection, via $\sharp$, with the orbits of $\Lambda_N(\OO_{K_e})$ on vectors of the same approximate discriminant. Moreover, by \eqref{eq:sharp_and_act}, this bijection respects the partition of these orbits among $\Lambda_N(K)$-orbits, respectively $\Lambda_N(K_e)$ orbits. 
  
  \begin{lem} \label{lem:tilt_AS}
    The tilted orbit representative
    \[
    v_s^\sharp = D^\sharp x^3 - x y^2 - s^\sharp y^3
    \]
    has discriminant $D^\sharp$, up to squares of units, and corresponding Kummer element
    \[
    \alpha = 1 + 3 s^\sharp \sqrt{-3 D^\sharp}.
    \]
  \end{lem}
  \begin{proof}
    The Kummer element for $v_{\alpha}^\sharp$ is the same as for
    \begin{equation*}
      x^3 - D^\sharp x y^2 - (D^{\sharp})^2 s^\sharp
    \end{equation*}
    which, by the cubic formula, is
    \begin{equation*}
      \kappa = \frac{s^\sharp(D^\sharp)^2}{2} + \sqrt{\frac{(s^\sharp)^2(D^\sharp)^4}{4} - \frac{(D^\sharp)^3}{27}}
      = \frac{s^\sharp(D^\sharp)^2}{2} +\frac{D^\sharp \sqrt{-3D^\sharp}}{9} \sqrt{1 - \frac{27(s^\sharp)^2D^\sharp}{4}},
    \end{equation*}
    where the last square root exists in $K_e$ because the radicand is close to $1$: specifically,
    \begin{equation*}
      1 - \frac{27(s^\sharp)^2D^\sharp}{4} = u^2, \quad \size{u - 1} = \size{27 s^2 D_0}.
    \end{equation*}
    So
    \begin{equation*}
      \kappa = \frac{s^\sharp(D^\sharp)^2}{2} + \frac{u^2 D^\sharp \sqrt{-3D^\sharp}}{9}
      = \( \frac{-\sqrt{-3 D^\sharp}}{3} \)^3 \( u^2 - \frac{3 s^\sharp \sqrt{-3 D^\sharp}}{2} \).
    \end{equation*}
    Now it suffices to show that the ratio
    \[
    \frac{u^2 - \dfrac{3 s^\sharp \sqrt{-3 D^\sharp}}{2}}{1 + 3 s^\sharp \sqrt{-3 D^\sharp}}
    \]
    is a cube. But
    \[
    \Size{\frac{u^2 - \dfrac{3 s^\sharp \sqrt{-3 D^\sharp}}{2}}{1 + 3 s^\sharp \sqrt{-3 D^\sharp}} - 1}
    = \Size{(u^2 - 1) - \frac{9 s \sqrt{-3 D^\sharp}}{2}}
    \leq \max\{\size{27 (s^\sharp)^2 D^\sharp}, \size{9 s^\sharp \sqrt{-3 D^\sharp}}\}
    < \size{3\sqrt{-3}},
    \]
    and as we determined before, any unit whose distance to $1$ is less than $\size{3\sqrt{-3}}$ is a cube.
  \end{proof}
  
  We derive that, for all $\alpha$ with $\size{\alpha} \leq q^M$,
  \[
  g_1(\alpha, \Lambda_{1,t}, \pi_e^{2t}D) = g_1(1 - 3 \sqrt{3} \cdot \alpha^\sharp, \Lambda_{1,t}, \pi^{2t}D).
  \]
  We then apply Lemma \ref{lem:tilt Hilbert} and note that, by Lemmas \ref{lem:tilt1} and \ref{lem:tilt_AS} the Kummer elements $1 + 3 \sqrt{-3D} s^\sharp$, $s \in K$, $\size{s} \leq q^M$, attain each $\alpha$ with nonzero $g_1(\alpha, \Lambda_{1,t}, (\pi^\sharp)^{2t}D)$ exactly once. We get
  \[
  \sum_{\alpha \in K_e[\sqrt{D^\sharp}]^{N=1}/\text{cubes}} (\beta, \alpha) \cdot g_1(\alpha, \Lambda_{1,t}, \pi_e^{2t}D^\sharp)
  = \sum_{\alpha \in K\sqrt{D}/\wp(K\sqrt{D})} (\beta, \alpha] \cdot g_1(\alpha, \Lambda_{1,t}, \pi^{2t}D)
  \]
  which finishes the proof of the theorem.
\end{proof}

TODOfn Write up the conjectured tilts of quadratic, cubic, and quartic O-N here.

\begin{rem}
  The version of cubic O-N for the group scheme $\alpha_3$, over function fields of curves in characteristic $3$, essentially recovers appropriate cases of the Riemann-Roch theorem. We are not the first to use the term ``arithmetic Riemann-Roch'' for various results provable by adelic Fourier analysis (see Mundy \cite{Mundy}).
\end{rem}
\end{todofn}

\section{Doubly traced quartic rings} \label{sec:doubly_traced}
It is an open question to classify, analogously to Theorem \ref{thm:Osborne}, the lattices in the space $V_K$ of pairs of ternary quadratic forms over $K$ invariant under $\GL_2(\OO_K) \cross \GL_3(\OO_K)$. This problem is not quite the most relevant to us because the factor $\GL_2(\OO_K)$, which changes the coordinates of the resolvent, is not relevant when we count rings with a \emph{fixed} resolvent, as we have done in this paper. We therefore look at lattices invariant under $\GL_3(\OO_K)$ alone. By inspection we find the following examples, which we conjecture exhaust all of them:
\begin{conj}
  Let $K$ be a number field, and let $\V(\OO_K)$ be the lattice of pairs of ternary quadratic forms over $\OO_K$.
  A primitive, $\GL_3(\OO_K)$-invariant lattice $L$ in $\V(\OO_K)$ is $\GL_2(\OO_K)$-equivalent to one whose completions $L_\pp$ are as follows:
  \begin{enumerate}[(a)]
    \item \label{it:intquartic} For $\pp|2$, the lattices $L_{\pp,t,s}$ of $(\pp^t,\pp^s)$-traced pairs of ternary quadratics
    \[
    L_{t,s} = L_{\pp,t,s} = \left\{ \left( \sum_{1 \leq i \leq j \leq 3} a_{ij} x_i x_j, \sum_{1 \leq i \leq j \leq 3} b_{ij} x_i x_j \right) :
    a_{ij} \equiv 0 \bmod \pp^t, b_{ij} \equiv 0 \bmod \pp^s\, (i < j)
    \right\}
    \]
    for $0 \leq t \leq s \leq e = v_p(2)$;
    \item For all other $\pp$, the maximal lattice $V_{\OO_{K,\pp}}$ only.
  \end{enumerate}
  In other words, they are of the form
  \[
  L_{\tt,\ss} = \bigcap_{\pp | 2} L_{\pp,\tt_\pp,\ss_\pp} = 
  \left\{ \left( \sum_{1 \leq i \leq j \leq 3} a_{ij} x_i x_j, \sum_{1 \leq i \leq j \leq 3} b_{ij} x_i x_j \right) :
  a_{ij} \equiv 0 \bmod \tt, b_{ij} \equiv 0 \bmod \ss\, (i < j)
  \right\}
  \]
  for ideals $(1) \supseteq \tt \supseteq \ss \supseteq (2)$.
\end{conj}
As in the cubic case, the same lattices $L_{\tt,\ss}$ necessarily appear in the analogues of $\V(\OO_K)$ with any Steinitz class $\aa$, appropriately adjusting the ideals that the $a_{ij}$, $b_{ij}$ must lie in.

For now, we look at the ($2$-adic) local case $K/\QQ_2$. An element of $L_{t,s}$ can be visualized as a pair of symmetric matrices
\[
(\A, \B) = 
\left(
\begin{bmatrix}
a_{11} & \frac{1}{2} a_{12} & \frac{1}{2} a_{13} \\
\frac{1}{2} a_{12} & a_{22} & \frac{1}{2} a_{23} \\
\frac{1}{2} a_{13} & \frac{1}{2} a_{23} & a_{33}
\end{bmatrix},
\begin{bmatrix}
b_{11} & \frac{1}{2} b_{12} & \frac{1}{2} b_{13} \\
\frac{1}{2} b_{12} & b_{22} & \frac{1}{2} b_{23} \\
\frac{1}{2} b_{13} & \frac{1}{2} b_{23} & b_{33}
\end{bmatrix}
\right)
\]
with diagonal entries $a_{ii}, b_{ii} \in \OO_K$ and off-diagonal entries $a_{ij} \in 2^{-1}\pi_K^{t}\OO_K$, $b_{ij} \in 2^{-1}\pi_K^{s}\OO_K$. It is easy to check that the cubic resolvent $f(\x,\y) = 4 \det (\A x - \B y)$ has the form
\[
f(\x,\y) = \pi^{2 s} a \x^3 + \pi^{2 s} b \x^2 \y + \pi^{2 t} c \x \y^2 + \pi^{2 t} d \y^3
\]
with discriminant $D \in \pi^{4 t + 4 s}$. Now $(A,B)$ parametrizes a quartic ring with a resolvent $\Phi : \OO \to C$, and $f(\x,\y)$ is the index form of $C$. The divisibility conditions on $f(\x,\y)$ can be interpreted as a non-maximality condition on $C$: since $t \leq s$,
\[
\frac{f(\pi^{-t}\x,\pi^{-t}\y)}{\pi^{-2 t}} = \pi^{2 s - 2 t} a \x^3 + \pi^{2 s - 2 t} b \x^2 \y + c \x \y^2 + d \y^3
\]
and
\[
\frac{f(\pi^{-t}\x,\pi^{-s}\y)}{\pi^{-t-s}} = a \x^3 + \pi^{s - t} b \x^2 \y + c \x \y^2 + \pi^{s - t} d \y^3
\]
are integral and thus are the index forms of certain overrings $C_1$, $C_2$ which we call the \emph{reduced resolvent} and the \emph{reduced coresolvent.} Appropriately lifting the basis $\<\bar\xi, \bar\eta\>$ of $C/\OO_K$ in which $(A,B)$ and hence $f$ are written, we have $C = \<1, \xi, \eta\>$, $C_1 = \<1, \pi^{-t}\xi, \pi^{-t}\eta\>$, and $C_2 = \<1, \pi^{-t}\xi, \pi^{-s}\eta\>$. Note that $C_1$ is a subring of $C_2$ of index $\pi^{s-t}$ and moreover is a \emph{unidirectional} subring in the sense that $C_1/C_2 \isom \OO_K/\pi^{s-t}$ is generated by one element. The integrality properties of $A$ and $B$ translate readily into relations between the resolvent and the rings $C_1$ and $C_2$.
\begin{defn}
  Let $\OO_K$ be a Dedekind domain, and let $L/K$ be a quartic algebra over its field of fractions with resolvent $\Phi : L \to R$. Let $C_1 \subseteq C_2 \subseteq R$ be a pair of subrings with $C_2/C_1 \isom \OO_K/\dd$ for some divisor $\dd$ of $2$ in $\OO_K$. Let $\tt$ be a divisor of $2\dd^{-1}$ and let $\ss = \tt \dd$.
  
  An order $\OO \subseteq L$ is \emph{$(\tt, \ss)$-traced with reduced resolvent $C_1$ and reduced coresolvent $C_2$} if
  \begin{enumerate}[(a)]
    \item $C = \OO_K + \tt^2 C_1$ is a resolvent for $\OO$, that is, $\disc \OO = \disc C$ and $\Phi(\OO) \subseteq C$;
    \item The associated bilinear form
    \[
    \tilde \Phi(\x,\y) = \Phi(\x + \y) - \Phi(\x) - \Phi(\y) = xy' + \x'\y + \x''\y''' + \x'''\y''
    \]
    maps $\OO \cross \OO$ into $\OO_K + \tt^3 \ss C_2$.
  \end{enumerate}
\end{defn}
If $C$ has a basis $\<1, \xi, \eta\>$ for which $C_1 = \<1, \tt^{-1}\xi, \tt^{-1}\eta\>$ and $C_2 = \<1, \tt^{-1}\xi, \ss^{-1}\eta\>$, as always happens when $\OO_K$ is a PID, then this is easily seen to be equivalent to the condition that $\Phi : \OO \to C$ is a resolvent whose matrix under this basis is in $L_{\tt,\ss}$.

The functional equation for the Shintani zeta functions on $V_\ZZ$ was stated and proved by Sato and Shintani (\cite{SatoShintani}). It relates pairs of integer-coefficient ternary quadratic forms (over $\ZZ$) with pairs of integer-matrix forms, that is, $L_{(1),(1)}$ to $L_{(2),(2)}$ in our notation. In contrast to the lattice of cubic forms, there is no $\SL_3(K)$-invariant inner product in $V_K$: as a representation of $\SL_3(K)$, $V_K$ is not isomorphic to its dual. However, $V_K$ \emph{is} isomorphic to its dual twisted by the automorphism of $\SL_3(K)$ given by inverse transpose, and under this duality, it is easy to see that $L_{\tt,\ss}$ is interchanged with $L_{2\ss^{-1},2\tt^{-1}}$. Therefore, it is a pretty conjecture that the corresponding composed varieties are naturally dual.

\begin{conj} \label{conj:O-N local quartic traced}
  Let $K$ be a $2$-adic local field and $0 \leq t \leq s \leq e = v_K(2)$ be integers. Let $C_1 \subseteq C_2$ be orders in a cubic $K$-algebra $R$ such that $C_2/C_1 \cong \OO_K/\pi^{s-t}$ is an $\OO_K$-module with one generator. For each quartic $K$-algebra $L$ with resolvent $R$, denote by $g(L, C_1, C_2, t, s)$ the number of $(t,s)$-traced quartic rings in $L$ with reduced resolvent and coresolvent $C_1$ and $C_2$, respectively. Then the dual of $g(L, C_1, C_2, t, s)$, considered as a function of $L$, is
  \begin{equation}\label{eq:O-N local quartic wild}
  \hat g(L, C_1, C_2, t, s) = q^{t+s} g(L, C_1, C_2, e-s, e-t)
  \end{equation}
  where $q = \Size{k_K}$ is the order of the residue field.
\end{conj}
Note that the quartic rings counted on either side of \eqref{eq:O-N local quartic wild} actually have resolvents $\OO_K + \pi^{t}C_1$ and $\OO_K + \pi^{e - s}C_1$ and discriminants $\pi^{2 t} \disc C_1$ and $\pi^{2 e - 2 s} \disc C_1$, respectively. If $t + s = e$, this conjecture asserts the self-duality (up to the correct scaling) of the indicated local weighting.

If this conjecture is true, then by our reflection engine, we immediately get the following corresponding global result.
\begin{conj} \label{conj:O-N quartic doubly traced}
  Let $K$ be a number field. Let $C_1 \subseteq C_2$ be orders in a cubic $K$-algebra $R$ such that $C_2/C_1 \cong \OO_K/\dd$ is an $\OO_K$-module with one generator. Let $\tt,\ss \subseteq \OO_K$ be ideals such that $\ss = \dd \tt \mid (2)$. Let $h(C_1, C_2, \tt, \ss)$ count the number of $(\tt, \ss)$-traced quartic rings with reduced resolvent and coresolvent $C_1$ and $C_2$, respectively, each weighted by the reciprocal of its number of resolvent-preserving automorphisms. Let $h^{\ntc}(C_1, C_2, \tt, \ss)$ count the subset of the foregoing that are ntc, weighted in the same way. Then
  \[
  h(C_1, C_2, \tt, \ss) = \frac{N(\tt\ss)}{2^{r_\infty}} \cdot h^{\ntc}(C_1, C_2, 2\ss^{-1}, 2\tt^{-1}),
  \]
  where $r_\infty$ is the number of real places of $K$ over which $R$ is not totally real plus twice the number of complex places of $K$.
\end{conj}

\section{Reflection for \texorpdfstring{$2 \times n \times n$}{2xnxn} boxes}

We close with a conjectural generalization to pairs of symmetric matrices $(\A,\B)$ of any odd order $n$. These play an important role in understanding $2$-torsion in $n$-ic rings, owing to a parametrization of Wood \cite{W2xnxn} as well as the Selmer groups of hyperelliptic curves of genus $(n - 1)/2$ (see \cite{bhargava2013hyperelliptic}). The connection to quartic rings, on the other hand, has no known analogue.

Let $K$ be a field, $\ch K \neq 2$ and $n$ be an odd integer. The group $\Gamma = \SL_n$ acts on the space $V$ of pairs of symmetric matrices $(\A,\B)$ over $K$, preserving the resolvent
\[
  f(x,y) = 2^{n-1}\det (\A x - \B y).
\]
The point stabilizer of an orbit with resolvent $f$ is isomorphic to $\C_2^{n-1}$, with Galois action permuting a hyperbasis as it permutes the roots of $f$, and the variety $(V_f, \Gamma)$ of pairs with fixed resolvent is a composed variety (closely related to the ``third representation'' of $\SO_n$ in \cite{AIT}). If $\OO_K \subseteq K$ is a Dedekind subring with fraction field $K$, there are integral forms $(\V_{\aa,\tt,\ss}, \G_\aa)$ of $V$ for each $[\aa] \in \Cl(\OO_K)$ and all ideals $(1) \supseteq \tt \supseteq \ss \supseteq (2)$, defined by a straightforward extension of the $2\times 3\times 3$ case. For simplicity we look only at the case when $[\aa] = [1]$, $\tt = (\tau) = \ss$ are principal:
\[
  \V_{\tau}(\OO_K) = \left\{[a_{ij}, b_{ij}] : a_{ij}, b_{ij} \in (2^{-1}\tau)^{\1_{i \neq j}}\OO_K \right\}.
\]
Note that the resolvent $f$ of such a box is divisible by $\tau^{n-1}$, since
\[
  f(x,y) = \tau^{n-1} \cdot 2^{-1}\tau \det\(2\tau^{-1} (\A x - \B y)\),
\]
and the argument to the determinant is an $\OO_K$-integral matrix that is skew-symmetric, hence singular (being of odd order), modulo $2\tau^{-1}$.

\begin{conj}[\textbf{Local O-N for $2\times n\times n$ boxes}]
If $K$ is a local field, $\ch K \neq 2$, then
\[
  \V_{\tau, \tau^{n-1}f} \textand 
  \V_{2\tau^{-1}, (2\tau^{-1})^{n-1}f}
\]
are naturally dual with duality constant $q^{(n-1)v_K(\tau)}$; in order words, for every binary $n$-ic form $f$, the associated local orbit counters
\[
g_{\tau, f} : H^1(K, M_f) \to \NN
\]
satisfy the local reflection theorem 
\[
\hat{g}_{\tau,f} = \Size{\OO_K/\tau\OO_K}^{n-1} \cdot  g_{2\tau^{-1},4f}.
\]
\end{conj}
If this conjecture is true, our local-to-global reflection engine yields the following.
\begin{conj}[\textbf{O-N for $2\times n\times n$ boxes}]
Let $K$ be a number field, and let $\tau$ be a divisor of $2$ in $\OO_K$. Denote by $h_\tau(f)$ the number of $\Gamma(\OO_K) = \SL_n(\OO_K)$-orbits of pairs $(\A,\B)$ of $n\times n$ symmetric matrices whose on- and off-diagonal elements belong to $\OO_K$ and $2^{-1} \tau \OO_K$ respectively, each $(\A,\B)$ weighted by the reciprocal of the order of its stabilizer in $\Gamma(\OO_K)$. Denote by $h_\tau^{\ntc}(f)$ the count (by the same weighting) of the subset of orbits whose corresponding self-balanced ideal $(R_f, I, \delta)$ has $\delta > 0$ at every real place of $K$. Then we have a global reflection theorem
\[
h_\tau(\tau^{n-1}f) = \frac{N(\tt)^2}{2^{r_\infty}} \cdot h^{\ntc}_{2\tau^{-1}}\((2\tau^{-1})^{n-1}\),
\]
where $r_\infty \in \ZZ$ depends only on the splitting type of $f$ at each of the real places of $K$.
\end{conj}
While some low-discriminant cases of this conjecture can be verified using the known structure of self-balanced ideals, the general case---which involves the Igusa zeta function of an intersection of two quadrics in $\PP^{n-1}$---is quite far from solvable using the techniques in this paper. If true, it furnishes a more satisfactory answer to the question addressed by Cohen--Rubinstein-Salzedo--Thorne \cite{CohON}, namely the production of a family of O-N-like reflection theorems for representations whose dimensions go to infinity.

\appendix
\part{Appendices}

\section{The Grothendieck-Witt ring and the proof of Lemma \ref{lem:H_form}}
\label{sec:GW}
In this section we prove Lemma \ref{lem:H_form}, which is of a different character than the other results in this paper.

To prove this lemma, we must recall some facts about the Grothendieck-Witt ring of a local field. The nondegenerate quadratic forms over a field $ K $ ($ \ch K \neq 2 $), up to isomorphism, form a semiring under the operations of orthogonal direct sum $ \perp $ and tensor product $ \tensor $; it is cancellative (the so-called \emph{Witt cancellation theorem}), and the ring obtained by adjoining formal additive inverses is called the \emph{Grothendieck-Witt ring} $ GW(K) $ of $ K $. In the case that $ K $ is a local field, a form is determined (see O'Meara \cite{OMeara}, Theorem 63:20) by three invariants: its dimension $ n \in \ZZ_{\geq 0} $, its determinant $ D \in K^\cross /(K^\cross)^2 $, and one other bit of information, the \emph{Hasse symbol} $ \epsilon \in \{ \pm 1\} $. For a diagonal form $ a_1 x_1^2 + \cdots + a_n x_n^2 $, the determinant is given by $ D = a_1 \cdots a_n $ and the Hasse symbol by
\[
\epsilon = \prod_{i<j} \<a_i, a_j\>_K.
\]
(The determinant is also called the \emph{discriminant}; our choice of terminology is influenced partly by the clash in sign with the discriminant of a binary quadratic form.)

Every combination of dimension, determinant, and Hasse symbol determines a unique element of $ GW(K) $. For dimension at least $3$, all elements are actually realized by a quadratic form; in dimension $ 3 $, the Hasse symbol carries the same information as the class in the Brauer group of the associated conic. The structure of $ GW(K) $ in terms of these invariants is easily computed by reducing to the case of diagonal forms; the formulas are here recorded, as they will be useful to us.
\begin{align}
  \dim (f \perp g) &= \dim f + \dim g \label{eq:dim_perp}\\
  \det (f \perp g) &= \det f \det g \\
  \epsilon(f \perp g) &= \epsilon(f) \epsilon(g) \<\det f, \det g\>_K \label{eq:Witt_add} \\
  \dim (f \tensor g) &= \dim f \dim g \\
  \det (f \tensor g) &= (\det f)^{\dim g} (\det g)^{\dim f} \\
  \epsilon(f \tensor g) &= \epsilon(f)^{\dim g} \epsilon(g)^{\dim f} \<\det f, -1\>_K^{\binom{\dim g}{2}} \<\det g, -1\>_K^{\binom{\dim f}{2}} \<\det f, \det g\>_K^{\dim f \dim g - 1}. \label{eq:Hasse_tensor}
\end{align}
We denote by $ \Ell_K $ the unique class in $ GW(K) $ of dimension $ 0 $, determinant $ 1 $, and Hasse symbol $ -1 $ (the ``elliptic class''). Note that if $ f $ and $ g $ are nonisomorphic quadratic forms over $ K $ with the same dimension and determinant, then $ [f] = [g] + \Ell_K $ in $ GW(K) $.

Suppose that $ L/K $ is a field extension. If $ q : W \to L $ is a quadratic form over $ L $, we can view $ W $ as a $ K $-vector space and postcompose with the trace $ \tr_{L/K} $ to get a quadratic form $ \tr_{L/K} q $. Since $ \tr_{L/K} $ respects orthogonal direct sums, it induces a group homomorphism (though not a ring homomorphism) from $ GW(L) $ to $ GW(K) $ (the opposite direction to the more familiar extension-of-scalars morphism). We easily compute that
\[
\dim \tr_{L/K} q = [L:K] \cdot \dim q \textand
\det \tr_{L/K} q = (\disc(L/K))^{\dim q} \cdot N_{L/K}(\det q).
\]
We wish to understand how $ \epsilon(\tr_{L/K} q) $ behaves. The following is the most important result needed.
\begin{lem}
  For any extension $ L/K $ of local fields not of characteristic $ 2 $,
  \[
  \tr_{L/K} \Ell_{L} = \Ell_K.
  \]
  In other words, $ \tr_{L/K} : GW(L) \to GW(K) $ preserves the Hasse symbol on classes of dimension $ 0 $ and determinant $ 1 $.
\end{lem}
\begin{proof}
  We may assume that $ L/K $ is a \emph{primitive} extension, that is, has no nontrivial intermediate extensions, since $ \tr_{E/K} \circ \tr_{L/E} = \tr_{L/K} $ for a tower $ L/E/K $.
  
  Since $ \tr_{L/K} \Ell_{L} $ is of dimension $ 0 $ and determinant $ 1 $, the only other possibility is that $ \tr_{L/K} \Ell_{L} = 0 $. We prove that this cannot hold.
  
  Let $ a \in K^\cross $ be an element that does not become a square in $ L $. Such an $ a $ exists because $ L/K $ is primitive; if not, then $ L/K $ would contain both an unramified and a ramified quadratic extension. Then choose $ \theta \in L $ such that $ \<a, \theta\>_{L} = -1 $. We also have $ \<a, N_{L/K}(\theta)\>_K = -1 $ by the standard relation $ \<a, \theta\>_{L} = \<a, N_{L/K}(\theta)\>_K $ ($ a \in K $, $ \theta \in L $). Consider the following quadratic forms over $ L $:
  \[
  f(x,y) = x^2 - \theta y^2, \quad g(x,y) = a x^2 - a \theta y^2.
  \]
  In other words, $ f = q_1 \perp q_{-\theta} $ and $ g = q_a + q_{-a\theta} $ where $ q_a(x) = a x^2 $. Both $ f $ and $ g $ have dimension $ 2 $ and discriminant $ -\theta $, but their Hasse symbols are $ 1 $ and $ \<a, -a \theta\> = -1 $, respectively. Hence $ [f] = [g] + \Ell_{L} $. But $ g = af $, so $ \tr_{L/K}(g) = a\tr_{L/K}(f) = q_a \tensor \tr_{L/K}(f) $. So by the formula \eqref{eq:Hasse_tensor} for the Hasse symbol of a tensor product,
  \begin{align*}
    \epsilon(\tr_{L/K} g) &= \epsilon(q_a \tensor \tr_{L/K} f) \\
    &= \epsilon(q_a)^{2[L:K]} \cdot \epsilon(\tr_{L/K} f)^{1} \cdot \<\det q_a, -1\>^{\binom{2[L:K]}{2}} \cdot \<\det \tr_{L/K} f, -1\>^{\binom{1}{2}} \cdot \<\det q_a, \det \tr_{L/K} f\>^{4 [L:K] - 1} \\
    &= \epsilon(\tr_{L/K} f) \cdot \<a, -1\>^{[L:K]} \cdot \<a, (\disc (L/K))^{2} \cdot N_{L/K}(\det f)\>^{4 [L:K] - 1} \\
    &= \epsilon(\tr_{L/K} f) \cdot \<a, -1\>^{[L:K]}\cdot \<a, N_{L/K}(-\theta)\> \\
    &= \epsilon(\tr_{L/K} f) \cdot \<a, N_{L/K}(\theta)\> \\
    &= -\epsilon(\tr_{L/K} f),
  \end{align*}
  So $ [\tr_{L/K} f] = [\tr_{L/K} g] + \Ell_{K}, $ yielding the desired conclusion.
\end{proof}

\begin{proof}[Proof of Lemma \ref{lem:H_form}]
  We now relate $ \epsilon(\alpha) $ to Hasse symbols. Denote by $ q_{R,\alpha} $ the quadratic form
  \[
  q_{R,\alpha}(\xi) = \alpha \xi^2
  \]
  over $ R $. Then $ \epsilon(\alpha) = 1 $ if and only if the form $ \tr_{R/K}(q_{R,\alpha}) $ is isotropic, where $ \tr_{R/K} $ is to be interpreted in the obvious way if $ R $ is not a field. Given $ \alpha,\beta \in R^\cross $ of norm $ 1 $, the forms
  \[
  f = q_{R,1} \perp q_{R,\alpha\beta} \textand
  g = q_{R,\alpha} \perp q_{R,\beta}
  \]
  have dimension $ 2 $ and determinant $ \alpha \beta $ over $ R $, and in $ GW(K), $
  \[
  \epsilon(1)\epsilon(\alpha)\epsilon(\beta)\epsilon(\alpha\beta) = \epsilon([\tr_{R/K} f] - [\tr_{R/K} g]).
  \]
  Decompose $ R = \prod_{i = 1}^r R_i $ into its field factors ($ 1 \leq r \leq 3 $), and let $ \alpha = (\alpha_1 ; \ldots ; \alpha_r) $ and $ \beta = (\beta_1; \ldots; \beta_r) $. Then
  \begin{align*}
    [\tr_{R/K} f] - [\tr_{R/K} g] = \sum_{i=1}^r \tr_{R_i/K} ([q_{R_i,\alpha_i\beta_i}] - [q_{R_i,\alpha_i}] - [q_{R_i,\beta_i}] + [q_{R_i, 1}]).
  \end{align*}
  Since the invariants of each $q_{R_i, \gamma}$ are known, the invariants of the class in parentheses can be computed by repeated application of \eqref{eq:dim_perp}--\eqref{eq:Witt_add}. We find that it has dimension $ 0 $, determinant $ 1 $ and Hasse symbol $ \<\alpha_{i}, \beta_{i}\>_{R_i} $. By the preceding lemma, its trace has the same invariants. Hence the whole sum has Hasse symbol
  \[
  \prod_i \<\alpha_{i}, \beta_{i}\>_{R_i} = \<\alpha, \beta\>_R,
  \]
  as desired.
\end{proof}

\section{Examples of zone totals}
\label{sec:zone_examples}
The following tables serve to illustrate some of the totals computed in Section \ref{sec:xi1}.

The following are valid in the unramified splitting types ($111$, $12$, and $3$): more specifically, the first table (with $s=0$) applies to all three, the remaining ones to splitting types $111$ and $12$.

For brevity, the following conventions have been observed:
\begin{itemize}
  \item In the red zone, the answer (which is independent of $\ell_C$) is the sum of the red-colored entries of the corresponding row.
  \item In the other color zones, the answer is the single entry corresponding to the appropriate values of $\ell_C$ and $n_{11}$.
  \item The invariable factor of $q^{-m_{11}}$ has been omitted.
\end{itemize}
%
%

\noindent\begin{minipage}{\textwidth}
For $e = 8, s = 0$ (similar results hold whenever $s = 0$):

\begin{sideways}
\begin{tabular}{r|ccccc}
$\ell_C\rightarrow$ & $0$ & $1$ & $2$ & $3$ & $4$ \\ \hline$n_{11} = 16$ & \textcolor{blue}{$q^{0} F(8,0,0)$} & \textcolor{blue}{$q^{0} F(8,0,0)$} & \textcolor{blue}{$q^{0} F(8,0,0)$} & \textcolor{blue}{$q^{0} F(8,0,0)$} & \textcolor{blue}{$q^{0} F(8,0,0)$} \\
$15$ & \textcolor{blue}{$q^{0} F(7,1,0)$} & \textcolor{blue}{$q^{0} F(7,1,0)$} & \textcolor{blue}{$q^{0} F(7,1,0)$} & \textcolor{blue}{$q^{0} F(7,1,0)$} & \textcolor{blue}{$q^{0} F(7,1,0)$} \\
$14$ & \textcolor{green}{$(1 + \epsilon_C)q^{0} F(7,0,1)$} & \textcolor{blue}{$q^{0} F(7,0,1)$} & \textcolor{blue}{$q^{0} F(7,0,1)$} & \textcolor{blue}{$q^{0} F(7,0,1)$} & \textcolor{blue}{$q^{0} F(7,0,1)$} \\
$13$ & \textcolor{green}{$(1 + \epsilon_C)q^{0} F(6,1,1)$} & \textcolor{blue}{$q^{0} F(6,1,1)$} & \textcolor{blue}{$q^{0} F(6,1,1)$} & \textcolor{blue}{$q^{0} F(6,1,1)$} & \textcolor{blue}{$q^{0} F(6,1,1)$} \\
$12$ & \textcolor{green}{$(1 + \epsilon_C)q^{0} F(6,0,2)$} & \textcolor{blue}{$q^{1} F(6,1,1)$} & \textcolor{blue}{$q^{1} F(6,1,1)$} & \textcolor{blue}{$q^{1} F(6,1,1)$} & \textcolor{blue}{$q^{1} F(6,1,1)$} \\
$11$ & \textcolor{green}{$(1 + \epsilon_C)q^{0} F(5,1,2)$} & \textcolor{blue}{$q^{1} F(5,2,1)$} & \textcolor{blue}{$q^{1} F(5,2,1)$} & \textcolor{blue}{$q^{1} F(5,2,1)$} & \textcolor{blue}{$q^{1} F(5,2,1)$} \\
$10$ & \textcolor{green}{$(1 + \epsilon_C)q^{0} F(5,0,3)$} & \textcolor{green}{$(1 + \epsilon_C)q^{1} F(5,1,2)$} & \textcolor{blue}{$q^{1} F(5,1,2)$} & \textcolor{blue}{$q^{1} F(5,1,2)$} & \textcolor{blue}{$q^{1} F(5,1,2)$} \\
$9$ & \textcolor{green}{$(1 + \epsilon_C)q^{0} F(4,1,3)$} & \textcolor{green}{$(1 + \epsilon_C)q^{1} F(4,2,2)$} & \textcolor{blue}{$q^{1} F(4,2,2)$} & \textcolor{blue}{$q^{1} F(4,2,2)$} & \textcolor{blue}{$q^{1} F(4,2,2)$} \\
$8$ & \textcolor{green}{$(1 + \epsilon_C)q^{0} F(4,0,4)$} & \textcolor{green}{$(1 + \epsilon_C)q^{1} F(4,1,3)$} & \textcolor{blue}{$q^{2} F(4,2,2)$} & \textcolor{blue}{$q^{2} F(4,2,2)$} & \textcolor{blue}{$q^{2} F(4,2,2)$} \\
$7$ & \textcolor{green}{$(1 + \epsilon_C)q^{0} F(3,1,4)$} & \textcolor{green}{$(1 + \epsilon_C)q^{1} F(3,2,3)$} & \textcolor{blue}{$q^{2} F(3,3,2)$} & \textcolor{blue}{$q^{2} F(3,3,2)$} & \textcolor{blue}{$q^{2} F(3,3,2)$} \\
$6$ & \textcolor{green}{$(1 + \epsilon_C)q^{0} F(3,0,5)$} & \textcolor{green}{$(1 + \epsilon_C)q^{1} F(3,1,4)$} & \textcolor{green}{$(1 + \epsilon_C)q^{2} F(3,2,3)$} & \textcolor{blue}{$q^{2} F(3,2,3)$} & \textcolor{blue}{$q^{2} F(3,2,3)$} \\
$5$ & \textcolor{green}{$(1 + \epsilon_C)q^{0} F(2,1,5)$} & \textcolor{green}{$(1 + \epsilon_C)q^{1} F(2,2,4)$} & \textcolor{red}{$(1 + \epsilon_C)q^{2} G^\cross(2,3,3)$} & \textcolor{red}{$q^{2} G(3,2,3)$} & \textcolor{red}{---} \\
$4$ & \textcolor{green}{$(1 + \epsilon_C)q^{0} F(2,0,6)$} & \textcolor{green}{$(1 + \epsilon_C)q^{1} F(2,1,5)$} & \textcolor{red}{$(1 + \epsilon_C)q^{2} G^\cross(2,2,4)$} & \textcolor{red}{$q^{3} G(3,2,3)$} & \textcolor{red}{---} \\
$3$ & \textcolor{green}{$(1 + \epsilon_C)q^{0} F(1,1,6)$} & \textcolor{red}{$(1 + \epsilon_C)q^{1} G^\cross(1,2,5)$} & \textcolor{red}{$(1 + \epsilon_C)q^{2} G^\cross(2,2,4)$} & \textcolor{red}{$q^{3} G(3,2,3)$} & \textcolor{red}{---} \\
$2$ & \textcolor{green}{$(1 + \epsilon_C)q^{0} F(1,0,7)$} & \textcolor{red}{$(1 + \epsilon_C)q^{1} G^\cross(1,1,6)$} & \textcolor{red}{$(1 + \epsilon_C)q^{2} G^\cross(2,1,5)$} & \textcolor{red}{$(1 + \epsilon_C)q^{3} G^\cross(3,1,4)$} & \textcolor{red}{$q^{3} G(4,0,4)$} \\
$1$ & \textcolor{red}{$(1 + \epsilon_C)q^{0} G^\cross(0,1,7)$} & \textcolor{red}{$(1 + \epsilon_C)q^{1} G^\cross(1,1,6)$} & \textcolor{red}{$(1 + \epsilon_C)q^{2} G^\cross(2,1,5)$} & \textcolor{red}{$(1 + \epsilon_C)q^{3} G^\cross(3,1,4)$} & \textcolor{red}{$q^{3} G(4,0,4)$} 
\end{tabular}
\end{sideways}
\end{minipage}
\vspace{0.1in}

\noindent\begin{minipage}{\textwidth}
For $e = 6, s = 14$ (similar results hold whenever $s \geq 2e$):

\begin{tabular}{r|c}
$\ell_C\rightarrow$ & $\text{any}$ \\ \hline$n_{11} = 12$ & \textcolor{violet}{$2q^{0} F(6,0, \emptyset)$} \\
$11$ & \textcolor{violet}{$2q^{0} F(5,1, \emptyset)$} \\
$10$ & \textcolor{violet}{$2q^{1} F(5,1, \emptyset)$} \\
$9$ & \textcolor{violet}{$2q^{1} F(4,2, \emptyset)$} \\
$8$ & \textcolor{violet}{$2q^{2} F(4,2, \emptyset)$} \\
$7$ & \textcolor{violet}{$2q^{2} F(3,3, \emptyset)$} \\
$6$ & \textcolor{violet}{$2q^{3} F(3,3, \emptyset)$} \\
$5$ & \textcolor{violet}{$2q^{3} F(2,4, \emptyset)$} \\
$4$ & \textcolor{violet}{$2q^{4} F(2,4, \emptyset)$} \\
$3$ & \textcolor{violet}{$2q^{4} F(1,5, \emptyset)$} \\
$2$ & \textcolor{violet}{$2q^{5} F(1,5, \emptyset)$} \\
$1$ & \textcolor{violet}{$2q^{5} F(0,6, \emptyset)$} 
\end{tabular}
\vspace{0.1in}
\end{minipage}

\noindent\begin{minipage}{\textwidth}
For $e = 11, s = 7$ (similar results hold whenever $s < 2e$ is odd):

\begin{tabular}{r|c}
$\ell_C\rightarrow$ & $3$ \\ \hline$n_{11} = 22$ & \textcolor{violet}{$2q^{0} F(11,0, \emptyset)$} \\
$21$ & \textcolor{violet}{$2q^{0} F(10,1, \emptyset)$} \\
$20$ & \textcolor{violet}{$2q^{1} F(10,1, \emptyset)$} \\
$19$ & \textcolor{violet}{$2q^{1} F(9,2, \emptyset)$} \\
$18$ & \textcolor{violet}{$2q^{2} F(9,2, \emptyset)$} \\
$17$ & \textcolor{violet}{$2q^{2} F(8,3, \emptyset)$} \\
$16$ & \textcolor{violet}{$2q^{3} F(8,3, \emptyset)$} \\
$15$ & \textcolor{green}{$(1 + \epsilon_C)q^{3} F(7,4,0)$} \\
$14$ & \textcolor{green}{$(1 + \epsilon_C)q^{3} F(7,3,1)$} \\
$13$ & \textcolor{green}{$(1 + \epsilon_C)q^{3} F(6,4,1)$} \\
$12$ & \textcolor{green}{$(1 + \epsilon_C)q^{3} F(6,3,2)$} \\
$11$ & \textcolor{green}{$(1 + \epsilon_C)q^{3} F(5,4,2)$} \\
$10$ & \textcolor{green}{$(1 + \epsilon_C)q^{3} F(5,3,3)$} \\
$9$ & \textcolor{green}{$(1 + \epsilon_C)q^{3} F(4,4,3)$} \\
$8$ & \textcolor{green}{$(1 + \epsilon_C)q^{3} F(4,3,4)$} \\
$7$ & \textcolor{green}{$(1 + \epsilon_C)q^{3} F(3,4,4)$} \\
$6$ & \textcolor{green}{$(1 + \epsilon_C)q^{3} F(3,3,5)$} \\
$5$ & \textcolor{green}{$(1 + \epsilon_C)q^{3} F(2,4,5)$} \\
$4$ & \textcolor{green}{$(1 + \epsilon_C)q^{3} F(2,3,6)$} \\
$3$ & \textcolor{green}{$(1 + \epsilon_C)q^{3} F(1,4,6)$} \\
$2$ & \textcolor{green}{$(1 + \epsilon_C)q^{3} F(1,3,7)$} \\
$1$ & \textcolor{green}{$(1 + \epsilon_C)q^{3} F(0,4,7)$} 
\end{tabular}
\vspace{0.1in}
\end{minipage}

\noindent\begin{minipage}{\textwidth}
For $e = 10, s = 4$ (similar results hold whenever $s < 2e$ is even):

\begin{sideways}
\begin{tabular}{r|ccccc}
$\ell_C\rightarrow$ & $0$ & $1$ & $2$ & $3$ & $\geq4$ \\ \hline$n_{11} = 20$ & \textcolor{violet}{$2q^{0} F(10,0, \emptyset)$} & \textcolor{violet}{$2q^{0} F(10,0, \emptyset)$} & \textcolor{violet}{$2q^{0} F(10,0, \emptyset)$} & \textcolor{violet}{$2q^{0} F(10,0, \emptyset)$} & \textcolor{violet}{$2q^{0} F(10,0, \emptyset)$} \\
$19$ & \textcolor{violet}{$2q^{0} F(9,1, \emptyset)$} & \textcolor{violet}{$2q^{0} F(9,1, \emptyset)$} & \textcolor{violet}{$2q^{0} F(9,1, \emptyset)$} & \textcolor{violet}{$2q^{0} F(9,1, \emptyset)$} & \textcolor{violet}{$2q^{0} F(9,1, \emptyset)$} \\
$18$ & \textcolor{violet}{$2q^{1} F(9,1, \emptyset)$} & \textcolor{violet}{$2q^{1} F(9,1, \emptyset)$} & \textcolor{violet}{$2q^{1} F(9,1, \emptyset)$} & \textcolor{violet}{$2q^{1} F(9,1, \emptyset)$} & \textcolor{violet}{$2q^{1} F(9,1, \emptyset)$} \\
$17$ & \textcolor{violet}{$2q^{1} F(8,2, \emptyset)$} & \textcolor{violet}{$2q^{1} F(8,2, \emptyset)$} & \textcolor{violet}{$2q^{1} F(8,2, \emptyset)$} & \textcolor{violet}{$2q^{1} F(8,2, \emptyset)$} & \textcolor{violet}{$2q^{1} F(8,2, \emptyset)$} \\
$16$ & \textcolor{blue}{$q^{2} F(8,2,0)$} & \textcolor{blue}{$q^{2} F(8,2,0)$} & \textcolor{blue}{$q^{2} F(8,2,0)$} & \textcolor{blue}{$q^{2} F(8,2,0)$} & \textcolor{blue}{$q^{2} F(8,2,0)$} \\
$15$ & \textcolor{blue}{$q^{2} F(7,3,0)$} & \textcolor{blue}{$q^{2} F(7,3,0)$} & \textcolor{blue}{$q^{2} F(7,3,0)$} & \textcolor{blue}{$q^{2} F(7,3,0)$} & \textcolor{blue}{$q^{2} F(7,3,0)$} \\
$14$ & \textcolor{green}{$(1 + \epsilon_C)q^{2} F(7,2,1)$} & \textcolor{blue}{$q^{2} F(7,2,1)$} & \textcolor{blue}{$q^{2} F(7,2,1)$} & \textcolor{blue}{$q^{2} F(7,2,1)$} & \textcolor{blue}{$q^{2} F(7,2,1)$} \\
$13$ & \textcolor{green}{$(1 + \epsilon_C)q^{2} F(6,3,1)$} & \textcolor{blue}{$q^{2} F(6,3,1)$} & \textcolor{blue}{$q^{2} F(6,3,1)$} & \textcolor{blue}{$q^{2} F(6,3,1)$} & \textcolor{blue}{$q^{2} F(6,3,1)$} \\
$12$ & \textcolor{green}{$(1 + \epsilon_C)q^{2} F(6,2,2)$} & \textcolor{blue}{$q^{3} F(6,3,1)$} & \textcolor{blue}{$q^{3} F(6,3,1)$} & \textcolor{blue}{$q^{3} F(6,3,1)$} & \textcolor{blue}{$q^{3} F(6,3,1)$} \\
$11$ & \textcolor{green}{$(1 + \epsilon_C)q^{2} F(5,3,2)$} & \textcolor{blue}{$q^{3} F(5,4,1)$} & \textcolor{blue}{$q^{3} F(5,4,1)$} & \textcolor{blue}{$q^{3} F(5,4,1)$} & \textcolor{blue}{$q^{3} F(5,4,1)$} \\
$10$ & \textcolor{green}{$(1 + \epsilon_C)q^{2} F(5,2,3)$} & \textcolor{green}{$(1 + \epsilon_C)q^{3} F(5,3,2)$} & \textcolor{blue}{$q^{3} F(5,3,2)$} & \textcolor{blue}{$q^{3} F(5,3,2)$} & \textcolor{blue}{$q^{3} F(5,3,2)$} \\
$9$ & \textcolor{green}{$(1 + \epsilon_C)q^{2} F(4,3,3)$} & \textcolor{green}{$(1 + \epsilon_C)q^{3} F(4,4,2)$} & \textcolor{blue}{$q^{3} F(4,4,2)$} & \textcolor{blue}{$q^{3} F(4,4,2)$} & \textcolor{blue}{$q^{3} F(4,4,2)$} \\
$8$ & \textcolor{green}{$(1 + \epsilon_C)q^{2} F(4,2,4)$} & \textcolor{green}{$(1 + \epsilon_C)q^{3} F(4,3,3)$} & \textcolor{blue}{$q^{4} F(4,4,2)$} & \textcolor{blue}{$q^{4} F(4,4,2)$} & \textcolor{blue}{$q^{4} F(4,4,2)$} \\
$7$ & \textcolor{green}{$(1 + \epsilon_C)q^{2} F(3,3,4)$} & \textcolor{green}{$(1 + \epsilon_C)q^{3} F(3,4,3)$} & \textcolor{blue}{$q^{4} F(3,5,2)$} & \textcolor{blue}{$q^{4} F(3,5,2)$} & \textcolor{blue}{$q^{4} F(3,5,2)$} \\
$6$ & \textcolor{green}{$(1 + \epsilon_C)q^{2} F(3,2,5)$} & \textcolor{green}{$(1 + \epsilon_C)q^{3} F(3,3,4)$} & \textcolor{green}{$(1 + \epsilon_C)q^{4} F(3,4,3)$} & \textcolor{blue}{$q^{4} F(3,4,3)$} & \textcolor{blue}{$q^{4} F(3,4,3)$} \\
$5$ & \textcolor{green}{$(1 + \epsilon_C)q^{2} F(2,3,5)$} & \textcolor{green}{$(1 + \epsilon_C)q^{3} F(2,4,4)$} & \textcolor{red}{$(1 + \epsilon_C)q^{4} G^\cross(2,5,3)$} & \textcolor{red}{$q^{4} G(3,4,3)$} & \textcolor{red}{---} \\
$4$ & \textcolor{green}{$(1 + \epsilon_C)q^{2} F(2,2,6)$} & \textcolor{green}{$(1 + \epsilon_C)q^{3} F(2,3,5)$} & \textcolor{red}{$(1 + \epsilon_C)q^{4} G^\cross(2,4,4)$} & \textcolor{red}{$q^{5} G(3,4,3)$} & \textcolor{red}{---} \\
$3$ & \textcolor{green}{$(1 + \epsilon_C)q^{2} F(1,3,6)$} & \textcolor{red}{$(1 + \epsilon_C)q^{3} G^\cross(1,4,5)$} & \textcolor{red}{$(1 + \epsilon_C)q^{4} G^\cross(2,4,4)$} & \textcolor{red}{$q^{5} G(3,4,3)$} & \textcolor{red}{---} \\
$2$ & \textcolor{green}{$(1 + \epsilon_C)q^{2} F(1,2,7)$} & \textcolor{red}{$(1 + \epsilon_C)q^{3} G^\cross(1,3,6)$} & \textcolor{red}{$(1 + \epsilon_C)q^{4} G^\cross(2,3,5)$} & \textcolor{red}{$(1 + \epsilon_C)q^{5} G^\cross(3,3,4)$} & \textcolor{red}{$q^{5} G(4,2,4)$} \\
$1$ & \textcolor{red}{$(1 + \epsilon_C)q^{2} G^\cross(0,3,7)$} & \textcolor{red}{$(1 + \epsilon_C)q^{3} G^\cross(1,3,6)$} & \textcolor{red}{$(1 + \epsilon_C)q^{4} G^\cross(2,3,5)$} & \textcolor{red}{$(1 + \epsilon_C)q^{5} G^\cross(3,3,4)$} & \textcolor{red}{$q^{5} G(4,2,4)$} 
\end{tabular}
\end{sideways}
\vspace{0.1in}
\end{minipage}



\bibliography{../Master}
\bibliographystyle{plain}

\end{document}